\documentclass[11pt]{report}

\usepackage{amssymb,amsmath,amsthm,enumitem,xy,mathabx,stmaryrd,dsfont,pstricks,pst-node,changepage,titlesec}
\xyoption{all}

\usepackage{imakeidx,blindtext}
\makeindex[name=gen,title={Index of Terms}]
\makeindex[name=not,title={Index of Notation}]
\makeindex
\def\gena#1{\index[gen]{#1}}
\def\nota#1{\index[not]{#1}}

\usepackage[cal=boondoxo]{mathalfa}

\usepackage{chngcntr}
\counterwithout{section}{chapter}

\renewcommand*\thesection{\arabic{section}}

\titleformat{\subsubsection}[runin]{\normalfont\bfseries}{\thesubsubsection.}{3pt}{}

\newcounter{temp}

\numberwithin{figure}{section}
\numberwithin{equation}{section}

\newtheorem{SpinPin}{SpinPin}
\newtheorem{RelSpinPin}{RelSpinPin}
\newtheorem{CROrient}{CROrient}

\newtheorem{thm}{Theorem}[section]
\newtheorem{prp}[thm]{Proposition}
\newtheorem{lmm}[thm]{Lemma}
\newtheorem{crl}[thm]{Corollary}

\theoremstyle{definition}
\newtheorem{dfn}[thm]{Definition}
\newtheorem{eg}[thm]{Example}
\newtheorem{rmk}[thm]{Remark}

\def\BE#1{\begin{equation}\label{#1}}
\def\EE{\end{equation}}
\def\eref#1{(\ref{#1})}

\def\ov#1{\overline{#1}}
\def\sf#1{\textsf{#1}}
\def\tn#1{\textnormal{#1}} 

\def\blr#1{\big\langle{#1}\big\rangle}
\def\lr#1{\langle{#1}\rangle}
\def\llrr#1{\langle\!\langle{#1}\rangle\!\rangle}
\def\bllrr#1{\big\langle\!\!\big\langle{#1}\big\rangle\!\!\big\rangle}

\def\wh#1{\widehat{#1}}
\def\wt#1{\widetilde{#1}}
\def\wch#1{\widecheck{#1}}
\def\mr#1{\mathring{#1}}
\def\un#1{\underline{#1}}
\def\sm#1{\begin{small}#1\end{small}}

\def\lra{\longrightarrow}

\def\Llra{\Longleftrightarrow}
\def\xlra#1{\xrightarrow{{#1}}}

\def\cA{\mathcal A}
\def\C{\mathbb C}
\def\cC{\mathcal C}
\def\bD{\mathbb D}
\def\cD{\mathcal D}
\def\fF{\mathfrak F}
\def\F{\mathbb F}
\def\G{\mathbb G}
\def\bH{\mathbb H}
\def\wcH{\widecheck{H}}
\def\cH{\mathcal H}

\def\cJ{\mathcal J}
\def\cK{\mathcal K}
\def\cL{\mathcal L}
\def\cM{\mathcal M}
\def\fM{\mathfrak M}
\def\cN{\mathcal N}
\def\cO{\mathcal O}
\def\fO{\mathfrak O}
\def\P{\mathbb P}
\def\cP{\mathcal P}
\def\cQ{\mathcal Q}

\def\R{\mathbb R}
\def\cR{\mathcal R}
\def\fR{\mathfrak R}
\def\cS{\mathcal S}

\def\cT{\mathcal T}
\def\cU{\mathcal U}
\def\Z{\mathbb Z}
\def\cZ{\mathcal Z}

\def\fI{\mathfrak i}
\def\fj{\mathfrak j}

\def\al{\alpha}
\def\be{\beta}
\def\de{\delta}
\def\ep{\epsilon}
\def\ga{\gamma}
\def\io{\iota}
\def\ka{\kappa}
\def\la{\lambda}
\def\na{\nabla}
\def\om{\omega}
\def\th{\theta}
\def\si{\sigma}
\def\ve{\varepsilon}
\def\vp{\varpi}
\def\vph{\varphi}
\def\ze{\zeta}

\def\De{\Delta}
\def\Ga{\Gamma}
\def\La{\Lambda}
\def\Om{\Omega}
\def\Si{\Sigma}
\def\Th{\Theta}

\def\fc{\mathfrak c}
\def\ce{\mathcal e}
\def\nd{\tn{d}}
\def\fd{\mathfrak d}
\def\ne{\tn{e}}
\def\ff{\mathfrak f}
\def\bh{\mathbf h}
\def\bI{\mathbb I}
\def\fI{\mathfrak i}
\def\fj{\mathfrak j}

\def\fo{\mathfrak o}
\def\os{\mathfrak{os}}
\def\fp{\mathfrak p}
\def\fs{\mathfrak s}
\def\ft{\mathfrak t}

\def\Abel{\tn{Abel}}
\def\ad{\tn{ad}}
\def\Ad{\tn{Ad}}
\def\Aut{\tn{Aut}}
\def\codim{\tn{codim}}
\def\Co{\tn{Co}}
\def\cok{\tn{cok}}
\def\Cov{\tn{Cov}}
\def\End{\tn{End}}
\def\ev{\tn{ev}}
\def\exp{\tn{exp}}
\def\Ext{\tn{Ext}}
\def\FS{\tn{FS}}
\def\gr{\tn{gr}}
\def\Hom{\tn{Hom}}
\def\id{\tn{id}}
\def\Im{\tn{Im}}
\def\ind{\tn{ind}}
\def\LB{\tn{LB}}
\def\nod{\tn{nd}}
\def\O{\tn{O}}
\def\OSp{\mathcal{OSp}}

\def\OVB{\tn{OVB}}
\def\PD{\tn{PD}}
\def\Pin{\tn{Pin}}
\def\PSL{\tn{PSL}}

\def\Prin{\tn{Prin}}
\def\pr{\tn{pr}}
\def\Re{\tn{Re}}
\def\Res{\tn{Res}}
\def\rk{\tn{rk}}
\def\ws{\tn{ws}}
\def\SG{\tn{SG}}
\def\SO{\tn{SO}}
\def\OSpin{\tn{OSpin}}
\def\Sp{\mathcal{Sp}}

\def\sp{\tn{sp}}
\def\Spin{\tn{Spin}}

\def\SU{\tn{SU}}
\def\St{\tn{St}}
\def\top{\tn{top}}
\def\Triv{\tn{Triv}}

\def\bu{\bullet}
\def\hb{\hbar}
\def\i{\infty}
\def\eset{\emptyset}
\def\1{\mathds 1}
\def\bfa{\mathbf a}
\def\bfc{\mathbf c}
\def\bfd{\mathbf d}

\def\p{\mathbf p}
\def\t{\mathbf t}
\def\u{\mathbf u}
\def\x{\mathbf x}
\def\prt{\partial}

\def\st{\star}
\def\bigst{\bigstar}

\def\bp{\bar\partial}
\def\dag{\dagger}

\def\pt{\tn{pt}}

\begin{document}

\title{Spin/Pin-Structures\\ and\\ Real Enumerative Geometry \vspace{1in}}

\vspace{2in}

\author{Xujia Chen and Aleksey Zinger}

\date{}

\maketitle

\thispagestyle{empty}
\tableofcontents
\thispagestyle{empty}
\addtocontents{toc}{\protect\thispagestyle{empty}}

\clearpage

\pagestyle{empty}
\chapter*{Preface}
\thispagestyle{empty}

Spin- and Pin-structures on vector bundles have long played an important role in differential geometry.
Since the mid-1990s, they have also been central to the open and real sectors
of Gromov-Witten theory, mirror symmetry, and enumerative geometry.
Spin- and Pin-structures in the classical perspective are (equivalence classes of)
principal Spin- and Pin-bundles that doubly cover 
the oriented orthonormal frame bundles of oriented vector bundles
and the orthonormal frame bundles of vector bundles, respectively, 
over paracompact (and usually locally contractible) topological spaces; 
this perspective motivates the terminology.
In the standard modern perspective of symplectic topology,
Spin- and Pin-structures are homotopy classes of trivializations of 
vector bundles over 2-skeleta of CW complexes.
In the perspective more recently introduced in~\cite{XCapsSigns},
Spin- and Pin-structures are collections of homotopy classes of trivializations of 
vector bundles over loops that respect cobordisms between the loops.
This last perspective applies over any topological space, is completely intrinsic, 
and connects directly with the usage of Spin- and Pin-structures 
in symplectic topology.
The present monograph provides an accessible introduction to Spin- and Pin-structures in general,
demonstrates their role in the orientability considerations in symplectic topology,
and exhibits their applications in enumerative geometry.\\

Chapter~\ref{topology_ch} contains a systematic treatment of Spin- and Pin-structures
in all three perspectives.
We in particular verify that these structures satisfy a collection of succinctly formulated properties
and describe natural correspondences between the three perspectives.
We also recall the notions of relative Spin- and Pin-structures originating
in the early versions of~\cite{FOOO} in a CW perspective and 
introduce an alternative perspective on them in the spirit of the perspective
of~\cite{XCapsSigns} on Spin- and Pin-structures.
In the latter perspective, relative Spin- and Pin-structures
are collections of homotopy classes of trivializations over boundaries
of compact bordered surfaces that respect certain cobordisms.
We conclude Chapter~\ref{topology_ch} by verifying that the relative Spin- and Pin-structures
in both perspectives satisfy a collection of succinctly formulated properties
and describe a natural correspondence between the two perspectives.\\

Real Cauchy-Riemann operators are first-order elliptic operators that arise
in symplectic topology as linearizations of the pseudoholomorphic map equation;
the latter has been central to the field since Gromov's seminal work~\cite{Gr}.
Spin- and Pin-structures and their relative counterparts induce orientations
of the determinants of such operators on vector bundle pairs over bordered surfaces
and on real bundle pairs over symmetric surfaces.
We detail constructions of these orientations,
which are central to the open and real sectors of Gromov-Witten theory,
in Chapter~\ref{CROrient_ch}. 
We also verify that these orientations satisfy a collection of succinctly formulated properties,
many of which are associated with the properties of relative Spin- and Pin-structures
established in Chapter~\ref{topology_ch}.
The intrinsic perspective on these structures introduced in Chapter~\ref{topology_ch}
fits perfectly with all considerations in Chapter~\ref{CROrient_ch}.\\

Chapter~\ref{Geom_chap} applies some of the key results of Chapters~\ref{topology_ch} 
and~\ref{CROrient_ch} to real enumerative geometry.
Welschinger's invariants~\cite{Wel4,Wel6,Wel6b} are signed counts of 
real (involution-invariant) holomorphic curves on K\"ahler surfaces 
and threefolds with anti-holomorphic involutions and more generally of 
pseudoholomorphic curves in symplectic fourfolds and sixfolds
with anti-symplectic involutions.
The interpretation of these invariants, up to a topological sign, in~\cite{Sol}
as signed counts of pseudoholomorphic maps from disks brought the potential of
applying the machinery of (complex) Gromov-Witten theory to study Welschinger's invariants.
In Chapter~\ref{Geom_chap}, we describe a further re-interpretation of these
invariants in~\cite{RealWDVV,RealWDVV3} as the degrees of relatively orientable pseudocycles 
from moduli spaces of real pseudoholomorphic maps with signed marked points introduced
in~\cite{Ge2}.
We also determine the exact topological signs relating different versions of 
Welschinger's invariants, using 
the classical perspective on the Spin- and Pin-structures 
in addition to the properties of orientations of real Cauchy-Riemann operators
established in Chapter~\ref{CROrient_ch}.
The examples of Sections~\ref{Wel4eg_subs} and~\ref{Wel6eg_subs} precisely match 
computations of these invariants in various perspectives.  
The first author's pseudocycle re-interpretation of Welschinger's invariants underpins
her proof in~\cite{RealWDVV} 
of Solomon's relations~\cite{Sol2} for these invariants envisioned almost 12~years earlier
and has led in~\cite{RealWDVV3} to similar relations for Welschinger's invariants of some real symplectic sixfolds.
As illustrated in~\cite{RealWDVVapp}, 
these relations provide an effective way of computing the invariants of~\cite{Wel4,Wel6,Wel6b}
for many real symplectic fourfolds.\\

We include two appendices containing standard statements to simplify 
the presentation in the main body of this monograph.
Appendix~\ref{CechH_app} recalls some relations between
the \v{C}ech cohomology of sheaves over paracompact spaces,
singular cohomology, vector bundles, and their characteristic classes.
Appendix~\ref{LG_app} relates the algebra and topology of covering projections
that are Lie group homomorphisms.\\

Many statements appearing in this monograph are well-known,
either available somewhere in the literature in some form or believed to be~true.
Some of them are gathered among the properties of Spin- and Pin-structures,
relative Spin- and Pin-structures, and
the orientations of the determinants of real Cauchy-Riemann operators in
Sections~\ref{SpinPinProp_subs}, \ref{RelSpinPinProp_subs}, \ref{OrientPrp_subs1}, 
and~\ref{OrientPrp_subs2}, ready for immediate use and fully established 
in later sections.
The more delicate and technical statements are new
and motivated by modern developments in symplectic topology and real enumerative geometry;
they come with complete proofs as~well. 
We hope this monograph overall will facilitate access~to and will further progress in 
these interrelated fields.\\ 

The authors were partially supported by the NSF grants DMS~1500875
and~1901979
and the Simons Collaboration grant~587036.\\

\vspace{.2in}

\noindent
{\it Department of Mathematics, Stony Brook University, Stony Brook, NY 11794\\
xujia@math.stonybrook.edu, azinger@math.stonybrook.edu}\\
\date{\today}

\clearpage
\setcounter{page}{1}
\pagestyle{plain}

\chapter{Spin- and Pin-Structures}
\label{topology_ch}

The purpose of this chapter is to provide a comprehensive introduction
to Spin- and Pin-structures in the classical perspective,
based on the Lie groups $\Spin(n)$ and~$\Pin^{\pm}(n)$,
and in the two modern perspectives of symplectic topology,
based on trivializations of vector bundles.
Section~\ref{intro_sec} formally presents Spin- and Pin-structures in
the three perspectives, collects their properties in a ready-to-use format,
and gives several explicit examples of these structures that are 
used in Chapters~\ref{CROrient_ch} and~\ref{Geom_chap}.
This section also contains one of the two theorems of the chapter,
Theorem~\ref{SpinStrEquiv_thm}, identifying the three perspectives on
Spin- and Pin-structures.\\

We review the definitions and key properties of the Lie groups $\Spin(n)$
and $\Pin^{\pm}(n)$ in Section~\ref{LG_sec} from purely topological considerations.
Explicit constructions of these Lie groups in terms of Clifford algebras
appear in \cite[\S3]{ABS}, in the first part of \cite[Section~1]{KirbyTaylor}, 
and in \cite[Sections~I.1,I.2]{LawsonMichel};
in the notation of \cite{LawsonMichel}, $\Pin^+(n)$ and $\Pin^-(n)$ are $\Pin_{n,0}$
and~$\Pin_{0,n}$, respectively.
However, these explicit constructions are not necessary for many purposes.
Section~\ref{LG_sec} is needed for the classical perspective on Spin- and Pin-structures,
but is irrelevant for the other two perspectives of Section~\ref{intro_sec}
and for formulating key properties of these structures.\\

Section~\ref{SpinPin_sec1} details the classical perspective on Spin- and Pin-structures,
provides further examples of these structures, and relates them
to \v{C}ech cohomology.
It also establishes the properties of Spin- and Pin-structures listed in Section~\ref{intro_sec}
in the classical perspective.
Section~\ref{SpinPin_sec2} establishes the same properties in the two modern perspectives 
of symplectic topology.
In Section~\ref{thmpf_sec}, we show that the three perspectives on Spin- and Pin-structures
are equivalent, when restricted to the appropriate categories of topological spaces,
by constructing natural correspondences between these structures in the classical perspective
and in each of the two modern perspectives;
these correspondences respect the properties listed in Section~\ref{intro_sec}.
Section~\ref{thmpf_sec} thus establishes Theorem~\ref{SpinStrEquiv_thm}.\\ 

Relative Spin- and Pin-structures trace their origins to the early versions of~\cite{FOOO},
where they appeared  as natural topological data for orienting 
the determinants of real Cauchy-Riemann operators on
vector bundle pairs over the unit disk $\bD^2_+\!\subset\!\C$.
In Section~\ref{RelSpinPin_sec}, we recall the original perspective on 
relative Spin- and Pin-structures, which involves an auxiliary choice
of a CW structure on the underlying pair of topological spaces $Y\!\subset\!X$,
introduce completely intrinsic perspective on~them, 
which fits directly with how these structures are used to orient real Cauchy-Riemann operators,
and collect important properties of relative Spin- and Pin-structures
in a ready-to-use format.
We also verify these properties in both perspectives and describe 
a natural correspondence between the two perspectives which respects 
the stated properties of relative Spin- and Pin-structures.
This establishes the other theorem of this chapter, Theorem~\ref{RelSpinStrEquiv_thm},
that the two perspectives on relative Spin- and Pin-structures are equivalent when 
restricted to CW complexes (or topological spaces homeomorphic to a CW~complex).\\

In retrospect, the CW perspectives on Spin- and Pin-structures and their relative
counterparts are rather artificial, in terms of both their formulations and
the intended applications.
The former is immediately apparent from the need to choose a CW structure on a topological space
and is further reflected in the verification of some of the properties of
these structures.
The latter is demonstrated by the need in
the construction of orientations of the determinants of real Cauchy-Riemann operators
in \cite[Section~8.1]{FOOO} to first choose a CW structure on 
the topological spaces $Y\!\subset\!X$ involved,
a vector bundle over the 3-skeleton $X_3\!\subset\!X$, and
a homotopy between a given continuous map from $(\bD^2_+,S^1)$ to $(X,Y)$
and a continuous map to $(X_3,Y_2)$ and then to show that the induced orientation does
not depend on all these choices.
We bypass these complications in Chapter~\ref{CROrient_ch} by using
the intrinsic perspective on relative Spin- and Pin-structures introduced
in Section~\ref{RelSpinPin_sec}.

\section{Introduction}
\label{intro_sec}

We define Spin- and Pin-structures in three perspectives and state 
the first theorem of this chapter, that these perspectives are essentially equivalent,
in Section~\ref{DfnThm_subs}.
The properties of these structures described in Section~\ref{SpinPinProp_subs}
include the obstructions to the existence of these structures and 
compatibility with short exact sequences of vector bundles.
Basic examples of Spin- and Pin-structures appear in Section~\ref{EgOutline_subs};
additional examples with more nuanced considerations appear in  Section~\ref{EgOutline_subs2}.

\subsection{Definitions and main theorem}
\label{DfnThm_subs}

For a topological space $Y$, let 
$$\tau_Y\equiv Y\!\times\!\R\lra Y \nota{tauy@$\tau_Y$}$$
be the trivial line bundle over~$Y$ and $\fo_Y$\nota{oy@$\fo_Y$}
be its standard orientation.
For a vector bundle~$V$ over~$Y$, we denote~by
$$\la(V)\nota{lav@$\la(V)$}\!\equiv\!\La_{\R}^{\top}V\lra Y$$
its top exterior power and by $\fO(V)$\nota{ov@$\fO(V)$} the set of the orientations on~$V$.
For each $\fo\!\in\!\fO(V)$, we denote by $\la(V,\fo)$\nota{lavo@$\la(V,\fo)$} 
the real vector bundle~$\la(V)$
with the orientation~$\la(\fo)$\nota{lao@$\la(\fo)$} induced by~$\fo$ and~by
\BE{Stdfn_e}\St(V,\fo)\nota{stvo@$\St(V,\fo)$}\equiv\big(\tau_Y\!\oplus\!V,\St_V(\fo)\big)\EE
the real vector bundle $\St(V)\!\equiv\!\tau_Y\!\oplus\!V$ with the induced orientation.
We make the canonical identifications
\BE{laStiden_e}\la(\tau_Y\!\oplus\!V)\!\equiv\!\tau_Y\!\otimes\!\la(V)=\la(V),
\quad
\la\big(\St(V,\fo)\big)\!\equiv\!
(\tau_Y,\fo_Y)\!\otimes\!\la(V,\fo)=\la(V,\fo)\,.\EE

\vspace{.15in}

For $n\!\in\!\Z^+$, the Lie groups $\Spin(n)$ and $\Pin^{\pm}(n)$ are double covers
of the \sf{$n$-th special orthogonal group}~$\SO(n)$\nota{son@$\SO(n)$} and 
the \sf{$n$-th orthogonal group}~$\O(n)\nota{on@$\O(n)$}$,
\BE{qnSpinPin_e}q_n\!:\Spin(n)\lra\SO(n)\nota{spinn@$\Spin(n)$}\gena{Spin group $\Spin(n)$} 
\qquad\hbox{and}\qquad \nota{qn@$q_n$}\nota{qnpm@$q_n^{\pm}$}
q_n^{\pm}\!:\Pin^{\pm}(n)\lra\O(n)\nota{pinn@$\Pin^{\pm}(n)$}
\gena{Pin group $\Pin^{\pm}(n)$};\EE
see Section~\ref{LG_sec}.
For $n\!\ge\!2$, $q_n$ is the unique connected double cover of~$\SO(n)$.
The double covers~$q_n^{\pm}$ restrict to~$q_n$ over~$\SO(n)$;
they  
are the same topologically, but have different group structures.
The preimages of an order~2 element of~$\O(n)$  
under $q_n^{\pm}$ generate a subgroup of~$\Pin^{\pm}(n)$ of order~4.
The subgroup generated by the preimages of an order~2 element of~$\O(n)$
with precisely one $(-1)$-eigenvalue is $\Z_2^{\,2}$ in the case of $\Pin^+(n)$ and 
$\Z_4$ in the case of $\Pin^-(n)$.\\

A rank~$n$ vector bundle $V$ over a paracompact space~$Y$
determines a principal $\O(n)$-bundle $\O(V)$\nota{ov@$\O(V)$} of orthonormal frames;
see Section~\ref{SpinPindfn_subs}.
Each orientation~$\fo$ on~$V$, determines 
a principal $\SO(n)$-bundle $\SO(V,\fo)$\nota{sovo@$\SO(V,\fo)$} of oriented orthonormal frames.
In the classical perspective, $\Pin^{\pm}$- and $\Spin$-structures are
liftings of these frame bundles over the projection maps~$q_n^{\pm}$ 
and~$q_n$. 

\begin{dfn}\label{PinSpin_dfn}
Let $Y$ be a paracompact space and $V$ be a  rank~$n$ real vector bundle over~$Y$.
\begin{enumerate}[label=(\alph*),leftmargin=*]

\item\label{PinStrDfn_it}
A \sf{$\Pin^{\pm}$-structure}~$\fp$\gena{Pin-structure!classical} 
on~$V$ 
is a principal $\Pin^{\pm}(n)$-bundle $\Pin^{\pm}(V)$ over~$Y$ with a $2\!:\!1$ 
covering~map
$$q_V\!:\Pin^{\pm}(V)\lra\O(V) $$ 
which commutes with the projections to~$Y$ and is equivariant with respect
to the group homomorphism~$q_n^{\pm}$. 

\item\label{SpinStrDfn_it} If $\fo\!\in\!\fO(V)$, 
a \sf{$\Spin$-structure}~$\fs$\gena{Spin-structure!classical} 
on~$(V,\fo)$
is a principal $\Spin(n)$-bundle $\Spin(V,\fo)$ over~$Y$ with a $2\!:\!1$ covering~map
$$q_V\!:\Spin(V,\fo)\lra\SO(V,\fo) $$ 
which commutes with the projections to~$Y$ and is equivariant with respect
to the group homomorphism~$q_n$. 

\end{enumerate}
\end{dfn}

\vspace{.1in}

\noindent
We call a pair $\os\!\equiv\!(\fo,\fs)$ consisting of an orientation~$\fo$ on a vector bundle~$V$ 
and a $\Spin$-structure~$\fs$ on~$(V,\fo)$ an 
\textsf{$\OSpin$-structure}\gena{OSpin-structure} on~$V$.
Two $\Spin$-structures 
\BE{2SpinStr_e}q_V\!:\Spin(V,\fo)\lra\SO(V,\fo) \qquad\hbox{and}\qquad 
q_V'\!:\Spin'(V,\fo)\lra\SO(V,\fo) \EE
are \sf{equivalent}\gena{Spin-structure!equivalent}
if there exists a $\Spin(n)$-equivariant isomorphism
$$\wt\Psi\!: \Spin(V,\fo)\lra\Spin'(V,\fo) \qquad\hbox{s.t.}\quad q_V=q_V'\!\circ\!\wt\Psi.$$
The analogous notions of equivalence\gena{Pin-structure!equivalent}\gena{OSpin-structure!equivalent} 
apply to $\Pin^{\pm}$ and 
$\OSpin$-structures.
When there is no ambiguity, we will not distinguish between the Spin- and Pin-structures
of Definition~\ref{PinSpin_dfn} and their equivalence classes.\\

\noindent
An $\OSpin$-structure on the trivial rank~$n$ vector bundle 
$$n\tau_Y\equiv Y\!\times\!\R^n\lra Y$$ 
with its canonical orientation is given
\BE{osVfodfn_e}
q_V\!=\!\id_Y\!\times\!q_n\!:\Spin\big(n\tau_Y,n\fo_Y\big)\!\equiv\!Y\!\times\!\Spin(n)
\lra Y\!\times\!\Spin(n)\!\equiv\!\SO(n\tau_Y,n\fo_Y\big)\,.\EE
An orientation~$\fo$ on a line bundle $V$ over a paracompact space~$Y$
determines a homotopy class of isomorphisms $(V,\fo)\!\approx\!(\tau_Y,\fo_Y)$.
A split of an oriented vector bundle~$(V,\fo)$ into $n$~oriented line bundles thus determines
a homotopy class of trivializations of the principal $\SO(n)$-bundle $\SO(V,\fo)$
and thus an $\OSpin$-structure~$\os_0(V,\fo)$\nota{os0vo@$\os_0(V,\fo)$} on~$V$.
In general, $\os_0(V,\fo)$  depends on the choice of the split.
Once its summands are chosen, the induced orientation on~$V$ depends 
on their ordering and on the orientations of the individual summands.
However, the $\Spin$-structures induced by 
orderings giving rise to different orientations are identified under 
the natural correspondence~\eref{OSpinRev_e}.\\

\noindent
For typical applications in symplectic topology,
it is more convenient to view these structures in terms
of trivializations of vector bundles.
The standard variation of this perspective, captured by Definition~\ref{PinSpin_dfn2} below,
goes back at least to the mid-1990s.
The second variation, captured by Definition~\ref{PinSpin_dfn3}, appears in
\cite[Section~5.1]{XCapsSigns}.
It is completely intrinsic and connects directly with its usage in symplectic topology.

\begin{dfn}\label{PinSpin_dfn2}
Let $(V,\fo)$ be an oriented vector bundle over a CW~complex~$Y$ with $\rk_{\R}V\!\ge\!3$.
A \sf{$\Spin$-structure}~$\fs$\gena{Spin-structure!CW trivializations}  
on~$(V,\fo)$ is a homotopy class 
of trivializations of~$(V,\fo)$ over the 2-skeleton~$Y_2$ of~$Y$.
 \end{dfn}

\noindent
A \sf{loop}\gena{loop} in a topological space~$Y$ is a continuous map $\al\!:S^1\!\lra\!Y$.
We denote the collection of all loops in~$Y$ by~$\cL(Y)$\nota{ly@$\cL(Y)$}. 
We call a compact two-dimensional topological manifold~$\Si$ with the boundary~$\prt\Si$ 
possibly nonempty a \sf{bordered surface}\gena{bordered surface}.
A \sf{closed surface}\gena{closed surface} is a bordered surface~$\Si$ with \hbox{$\prt\Si\!=\!\eset$}.

\begin{dfn}\label{PinSpin_dfn3}
Let $(V,\fo)$ be an oriented vector bundle over a topological space~$Y$
with \hbox{$\rk_{\R}V\!\ge\!3$}.
A \sf{$\Spin$-structure}~$\fs$\gena{Spin-structure!loop trivializations} 
on~$(V,\fo)$ is a collection $(\fs_{\al})_{\al\in\cL(Y)}$
of homotopy classes~$\fs_{\al}$ of trivializations of $\al^*(V,\fo)$
such that for every continuous map $F\!:\Si\!\lra\!Y$ 
from a  bordered surface 
the vector bundle $F^*\!(V,\fo)$ over~$\Si$ admits a trivialization whose restriction
to each component~$\prt_r\Si$ of~$\prt\Si$ lies in $\fs_{u|_{\prt_r\Si}}$
under any identification of~$\prt_r\Si$ with~$S^1$.
\end{dfn}

The existence of a $\Spin$-structure $\fs$ in the sense of  Definition~\ref{PinSpin_dfn3}
on an oriented vector bundle $(V,\fo)$ with \hbox{$\rk_{\R}V\!\ge\!3$}  
explicitly requires the bundle~$F^*V$ to be trivializable for 
every continuous map $F\!:\Si\!\lra\!Y$ from a closed surface.\\

In both variations of the trivializations perspective, a 
\sf{Spin structure}~$\fs$\gena{Spin-structure!CW trivializations}
\gena{Spin-structure!loop trivializations}  
on an oriented vector bundle~$(V,\fo)$ over~$Y$
with $\rk_{\R}V\!=\!1,2$ is a $\Spin$-structure on the vector bundle
$2\tau_Y\!\oplus\!V$ with the induced orientation in the first case and
on $\tau_Y\!\oplus\!V$ in the second.
A \sf{$\Pin^{\pm}$-structure}~$\fp$
\gena{Pin-structure!CW trivializations}\gena{Pin-structure!loop trivializations}  
on a real vector bundle~$V$ over~$Y$
is a $\Spin$-structure on the canonically oriented vector bundle 
\BE{Vpmdfn_e}V_{\pm}\equiv V\!\oplus\!(2\!\pm\!1)\la(V)\, \nota{Vpm@$V_{\pm}$};\EE
see the paragraph containing~\eref{Vpmdfn_e2}.
An \textsf{$\OSpin$-structure} on~$V$ is again a pair \hbox{$\os\!\equiv\!(\fo,\fs)$} 
consisting of an orientation~$\fo$ on~$V$ and a $\Spin$-structure~$\fs$ on~$(V,\fo)$.
Analogously to the perspective of Definition~\ref{PinSpin_dfn},
a split of an oriented vector bundle~$(V,\fo)$ into $n$~oriented line bundles 
determines an $\OSpin$-structure~$\os_0(V,\fo)$\nota{os0vo@$\os_0(V,\fo)$} on~$V$.\\

For a vector bundle~$V$ over a topological space~$Y$ satisfying the appropriate conditions,
we denote by $\cP^{\pm}(V)$\nota{pv@$\cP^{\pm}(V)$} and~$\OSp(V)$\nota{ospv@$\OSp(V)$}
the sets of $\Pin^{\pm}$-structures and $\OSpin$-structures, respectively, on~$V$
in any given of the three perspectives 
(up to equivalence in the perspective of Definition~\ref{PinSpin_dfn}).
For $\fo\!\in\!\fO(V)$,  we denote by $\Sp(V,\fo)$\nota{spvo@$\Sp(V,\fo)$} 
the set of $\Spin$-structures on~$(V,\fo)$.
We identity $\Sp(V,\fo)$ with a subset of~$\OSp(V)$ in the obvious way.\\

By Theorem~\ref{SpinStrEquiv_thm} below, there is no fundamental ambiguity in 
the definitions of the sets $\cP^{\pm}(V)$, $\OSp(V)$, and $\OSp(V,\fo)$
when $Y$ is a CW complex (and so the perspectives of 
Definitions~\ref{PinSpin_dfn}, \ref{PinSpin_dfn2}, and~\ref{PinSpin_dfn3} apply) or 
more generally $Y$ is a paracompact locally contractible space 
(and the perspectives of Definitions~\ref{PinSpin_dfn} and~\ref{PinSpin_dfn3} apply).
The local contractability restriction can in fact be weakened to 
the local $H^2(\cdot;\Z_2)$-simplicity of Definition~\ref{locH1simp_dfn}.

\begin{thm}\label{SpinStrEquiv_thm}
Let $Y$ be a topological space and $V$ be a vector bundle over~$Y$. 
\begin{enumerate}[label=(\arabic*),leftmargin=*]

\item\label{PinSpinPrp_it} The OSpin- and Pin-structures on~$V$
in the perspectives of Definitions~\ref{PinSpin_dfn} 
if $Y$ is paracompact and locally contractible,
\ref{PinSpin_dfn2} if $Y$ is a CW complex, and~\ref{PinSpin_dfn3} 
satisfy all properties of Section~\ref{SpinPinProp_subs}.

\item\label{SpinStrEquiv1vs3_it} If $Y$ is a paracompact locally contractible space,
there are canonical identifications of the sets $\OSp(V)$ in the perspectives of 
Definitions~\ref{PinSpin_dfn} and~\ref{PinSpin_dfn3} and
of the sets $\cP^{\pm}(V)$ in the two perspectives for every vector bundle~$V$ over~$Y$.
These identifications associate the distinguished elements~$\fs_0(V,\fo)$
in the two perspectives with each other for all oriented vector bundles~$(V,\fo)$
split into oriented line bundles and respect all structures 
and correspondences of Section~\ref{SpinPinProp_subs}.

\item\label{SpinStrEquiv1vs2_it} If $Y$ is a CW complex,
the same statements apply to the OSpin- and Pin-structures on~$V$
in the perspectives of Definitions~\ref{PinSpin_dfn} and~\ref{PinSpin_dfn2}.

\end{enumerate}
\end{thm}

\vspace{.1in}

\noindent 
This theorem is fundamentally a consequence of
\begin{gather}
\label{pi012SO_e}
\big|\pi_0\big(\SO(n)\big)\big|,\big|\pi_2\big(\SO(n)\big)\big|=1,
\quad\big|\pi_1\big(\SO(n)\big)\big|=2,\\
\label{pi012Spin_e}
\big|\pi_0\big(\Spin(n)\big)\big|,\big|\pi_1\big(\Spin(n)\big)\big|,
\big|\pi_2\big(\Spin(n)\big)\big|=1
\end{gather}
for $n\!\ge\!3$ for the following reasons.
Suppose $Y$ is a CW complex of dimension at most~2 
and $V$ is an oriented vector bundle over~$Y$ of rank $n\!\ge\!3$.
By~\eref{pi012Spin_e}, every bundle $\Spin(V,\fo)$
as in Definition~\ref{PinSpin_dfn}\ref{SpinStrDfn_it}
is trivializable as a principal bundle and thus admits a section~$\wt{s}$;
any two such sections are homotopic.
The section $q_V\!\circ\!\wt{s}$ of $\SO(V,\fo)$ then determines a trivialization
of the vector bundle~$V$ over~$Y$;
the trivializations of~$V$ determined by different sections of~$\Spin(V,\fo)$ are homotopic.\\

The above implies that a $\Spin$-structure in the sense of 
Definition~\ref{PinSpin_dfn} determines
a $\Spin$-structure in the sense of Definition~\ref{PinSpin_dfn2} and 
a $\Spin$-structure in the sense of Definition~\ref{PinSpin_dfn3},
under appropriate topological conditions on the underlying topological space~$Y$.
It is immediate that equivalent $\Spin$-structures in the sense of 
Definition~\ref{PinSpin_dfn} determine the same $\Spin$-structures
in the perspectives of  Definitions~\ref{PinSpin_dfn2} and~\ref{PinSpin_dfn3}
in the above construction.
The next observation ensures that this leads to a correspondence between the three notions
of $\Spin$-structure.

\begin{lmm}\label{SpinStrEquiv_lmm}
Let $(V,\fo)$ be an oriented vector bundle over~$S^1$ with $\rk\,V\!\ge\!3$.
The homotopy classes of trivializations of~$V$ determined by 
different equivalence classes of $\Spin$-structures on~$(V,\fo)$
in the sense of Definition~\ref{PinSpin_dfn} are different.
\end{lmm}

Theorem~\ref{SpinStrEquiv_thm} and Lemma~\ref{SpinStrEquiv_lmm} are 
proved in Section~\ref{thmpf_sec}.

\subsection{Properties of \Spin- and \Pin-structures}
\label{SpinPinProp_subs} 

Let $Y$ be a topological space. An isomorphism $\Psi\!:V'\!\lra\!V$ of vector bundles
over~$Y$ induces bijections 
\BE{SPinPullback_e}\Psi^*\!: \OSp(V)\lra\OSp(V') \qquad\hbox{and}\qquad 
\Psi^*\!: \cP^{\pm}(V)\lra\cP^{\pm}(V')\EE
between the $\OSpin$-structures on~$V$ and~$V'$ and the $\Pin^{\pm}$-structures on~$V$ and~$V'$
in the perspectives of Definitions~\ref{PinSpin_dfn} if $Y$ is paracompact,
\ref{PinSpin_dfn2} if $Y$ is a CW complex, and~\ref{PinSpin_dfn3}.
If the ranks of~$V$ and~$V'$ are at least~3,
the $\OSpin$- and $\Pin^{\pm}$-structures on~$V'$ in the last two perspectives
are obtained from the same types of structures on~$V$ by 
pre-composing the relevant trivializations with~$\Psi$.
If $Y$ is paracompact, $\Psi$ induces an isomorphism
$$\O(\Psi)\!:\O(V')\lra \O(V)$$
of principal $\O(\rk\,V)$-bundles over~$Y$.
If $q_V$ is a $\Pin^{\pm}$-structure on~$V$ in the sense of 
Definition~\ref{PinSpin_dfn}\ref{PinStrDfn_it},  then 
$$\O(\Psi)^*q_V\!:\O(\Psi)^*\Pin(V)\lra\O(V'), \qquad
\big(p',\wt{p}\big)\lra p',$$
is a $\Pin^{\pm}$-structure on~$V'$.
If in addition $\fo\!\in\!\fO(V)$ and $\Psi^*\fo$ is the orientation on~$V'$ induced by~$\fo$
via~$\Psi$, then~$\Psi$ also induces an isomorphism
$$\SO(\Psi,\fo)\!:\SO\big(V',\Psi^*\fo\big)\lra \SO(V,\fo)$$
of principal $\SO(\rk\,V)$-bundles over~$Y$.
If $q_V$ is a $\Spin$-structure on~$(V,\fo)$ in the sense of 
Definition~\ref{PinSpin_dfn}\ref{SpinStrDfn_it},  then 
$$\SO(\Psi,\fo)^*q_V\!:\SO(\Psi,\fo)^*\Spin(V,\fo)\lra\SO\big(V',\Psi^*\fo\big), \qquad
\big(p',\wt{p}\big)\lra p',$$
is a $\Spin$-structure on~$(V',\Psi^*\fo)$.\\

Let $V$ be a vector bundle over a topological space~$Y$ and $Y'$ be another topological space.
A continuous map $f\!:Y'\!\lra\!Y$ induces maps
\BE{fSPinPullback_e}f^*\!: \OSp(V)\lra\OSp(f^*V) \qquad\hbox{and}\qquad 
f^*\!: \cP^{\pm}(V)\lra\cP^{\pm}(f^*V)\EE
in the perspectives of Definitions~\ref{PinSpin_dfn} if $Y$ and~$Y'$ are paracompact,
\ref{PinSpin_dfn2} if $f$ is a map of CW complexes, and~\ref{PinSpin_dfn3}.
In the last two perspectives, the $\OSpin$- and $\Pin^{\pm}$-structures on~$f^*V$ 
are obtained from the same types of structures on~$V$ by 
pre-composing the relevant trivializations with the~map
\BE{fVpullPsidfn_e}\wt{f}_V\!:f^*V\lra V, \qquad (y',v)\lra v,\EE
covering~$f$.
If $Y$ is paracompact, this map induces an $\O(\rk\,V)$-equivariant map
$$\O(\wt{f}_V)\!:\O(f^*V)\lra \O(V)$$
covering~$f$.
If $q_V$ is a $\Pin^{\pm}$-structure on~$V$ in the sense of 
Definition~\ref{PinSpin_dfn}\ref{PinStrDfn_it},  then 
$$\O(\wt{f}_V)^*q_V\!:\O(\wt{f}_V)^*\Pin(V)\lra\O(f^*V), \quad
\big(p',\wt{p}\big)\lra p',$$
is a $\Pin^{\pm}$-structure on~$f^*V$.
If in addition $\fo\!\in\!\fO(V)$ and $f^*\fo$ is the orientation on~$f^*V$ induced by~$\fo$
via~$\wt{f}_V$, then~$\wt{f}_V$ also induces an $\SO(\rk\,V)$-equivariant map
$$\SO(\wt{f}_V,\fo)\!:\SO\big(f^*V,f^*\fo\big)\lra \SO(V,\fo)$$
covering~$f$.
If $q_V$ is a $\Spin$-structure on~$(V,\fo)$ in the sense of 
Definition~\ref{PinSpin_dfn}\ref{SpinStrDfn_it},  then 
$$\SO\big(\wt{f}_V,\fo\big)^*q_V\!:\SO\big(\wt{f}_V,\fo\big)^*\Spin(V,\fo)
\lra\SO\big(f^*V,f^*\fo\big), \qquad
\big(p',\wt{p}\big)\lra p',$$
is a $\Spin$-structure on~$(f^*V',f^*\fo)$.\\

If an oriented vector bundle~$(V,\fo)$ is split as a direct sum of oriented line bundles,
we take the induced splitting of $(V,\ov\fo)$ to be the splitting 
obtained from the splitting of~$(V,\fo)$ by negating the last component.
We take the induced splitting of $\St_V(V,\fo)$ to be the splitting obtained by combining
the canonical trivialization of~$\tau_Y$ as the first component
with the splitting of~$(V,\fo)$.\\

The SpinPin properties below apply in any of the three perspectives of 
Definitions~\ref{PinSpin_dfn}-\ref{PinSpin_dfn3} on \hbox{$\Pin^{\pm}$-,} $\Spin$-, and $\OSpin$-structures,
provided the topological space~$Y$ appearing in these statements satisfies
the appropriate restrictions:
\begin{enumerate}[label=$\bullet$,leftmargin=*]

\item $Y$ is a paracompact locally contractible space for the perspective of  
Definition~\ref{PinSpin_dfn};

\item $Y$ is a CW complex for the perspective of Definition~\ref{PinSpin_dfn2}.

\end{enumerate}
The naturality properties of the group actions and correspondences below
refer to the commutativity with the pullbacks~\eref{SPinPullback_e}
induced by isomorphisms of vector bundles over~$Y$ and 
the pullbacks~\eref{fSPinPullback_e} induced by the admissible continuous maps,
as appropriate.

\begin{SpinPin}[Obstruction to Existence]\label{SpinPinObs_prop}
\begin{enumerate}[label=(\alph*),leftmargin=*]

\item\label{PinObs_it} A real vector bundle $V$ over $Y$ admits a $\Pin^-$-structure 
(resp.~$\Pin^+$-structure) if and only if $w_2(V)\!=\!w_1^2(V)$
(resp.~$w_2(V)\!=\!0$).

\item\label{SpinObs_it} An oriented vector bundle $(V,\fo)$ over $Y$ admits 
a $\Spin$-structure if and only if $w_2(V)\!=\!0$.

\end{enumerate}
\end{SpinPin}

\begin{SpinPin}[Affine Structure]\label{SpinPinStr_prop}
Let $V$ be a real vector bundle over~$Y$.
\begin{enumerate}[label=(\alph*),leftmargin=*]

\item\label{PinStr_it}  If\, $V$ admits a $\Pin^{\pm}$-structure, then
the group $H^1(Y;\Z_2)$ acts naturally, freely, and transitively on the set~$\cP^{\pm}(V)$.

\item\label{SpinStr_it}  If $\fo\!\in\!\fO(V)$ and $(V,\fo)$
admits a $\Spin$-structure, then
the group $H^1(Y;\Z_2)$ acts naturally, freely, and transitively on the set~$\Sp(V,\fo)$.

\end{enumerate}
If $V'$ and $V''$ are vector bundles over~$Y$, 
with at least one of them of positive rank, 
then the action of the automorphism
\BE{SpinPinStr_e0} \Psi\!:V'\!\oplus\!V''\lra V'\!\oplus\!V'', \qquad \Psi(v',v'')=(-v',v''),\EE
on $\cP^{\pm}(V'\!\oplus\!V'')$ is given~by 
\BE{SpinPinStr_e}  \Psi^*\fp= 
\big((\rk\,V'\!-\!1)w_1(V')\!+\!(\rk\,V')w_1(V'')\big)\!\cdot\!\fp
\qquad\forall\,\fp\!\in\!\cP^{\pm}(V'\!\oplus\!V'').\EE
\end{SpinPin}

\begin{SpinPin}[Orientation Reversal]\label{OSpinRev_prop}
Let $V$ be a real vector bundle over~$Y$.
There is a natural $H^1(Y;\Z_2)$-equivariant involution
\BE{OSpinRev_e}\OSp(V)\lra \OSp(V), \qquad \os\lra\ov\os,\EE
which maps $\Sp(V,\fo)$ bijectively onto $\Sp(V,\ov\fo)$ for every $\fo\!\in\!\fO(V)$
and sends $\os_0(V,\fo)$ to $\os_0(V,\ov\fo)$ for every oriented vector bundle~$(V,\fo)$
split as a direct sum of oriented line bundles.
\end{SpinPin}

\begin{SpinPin}[Reduction]\label{Pin2SpinRed_prop}
Let $V$ be a real vector bundle over~$Y$.
For every $\fo\!\in\!\fO(V)$, there are natural 
$H^1(Y;\Z_2)$-equivariant bijections
\BE{Pin2SpinRed_e} \fR_{\fo}^{\pm}\!:\cP^{\pm}(V)\lra \Sp(V,\fo)\EE
so that $\fR_{\ov\fo}^{\pm}(\cdot)\!=\!\ov{\fR_{\fo}^{\pm}(\cdot)}$.
\end{SpinPin}

\begin{SpinPin}[Stability]\label{SpinPinStab_prop}
Let $V$ be a real vector bundle over~$Y$.
There are natural $H^1(Y;\Z_2)$-equivariant bijections
\BE{SpinPinStab_e}
\St_V\!: \OSp(V)\lra \OSp\big(\tau_Y\!\oplus\!V\big), \quad
\St_V^{\pm}\!:\cP^{\pm}(V)\lra \cP^{\pm}\big(\tau_Y\!\oplus\!V\big)\EE
so that $\St_V(\os_0(V,\fo))\!=\!\os_0(\St_V(V,\fo))$ for every oriented vector bundle~$(V,\fo)$
split as a direct sum of oriented line bundles,
\BE{SpinPinStab_e2}
\St_V(\ov\os)=\ov{\St_V(\os)}~~\forall\,\os\!\in\!\OSp(V), \quad
\St_V\!\circ\!\fR_{\fo}^{\pm}\!=\!\fR_{\St_V(\fo)}^{\pm}\!\circ\!\St_V^{\pm}
~~\forall\,\fo\!\in\!\fO(V).\EE
\end{SpinPin}

\vspace{.15in}

If $V$ is a real vector bundle over~$Y$, then the real vector bundles~\eref{Vpmdfn_e} 
over~$Y$ have canonical orientations, which we denote by~$\fo_V^-$ and~$\fo_V^+$.
They are described as follows.
The orientations on these bundles 
canonically correspond to the orientations~on the real line bundles
\BE{Vpmdfn_e2}\la(V_-)=\la(V)\!\otimes\!\la(V) \quad\hbox{and}\quad
\la(V_+)=\la(V)\!\otimes\!\la(V)
\!\otimes\!\la(V)\!\otimes\!\la(V),\EE
respectively.
We orient the fiber of the first line bundle 
 over a point $y\!\in\!Y$ by choosing
an orientation on the first factor $\la(V_y)$ and using the same orientation
on the second factor.
The resulting orientation on the product does not depend on the choice of 
the orientation on the first factor.
We thus obtain an orientation on~$\la(V_-)$.
The same reasoning applies to~$\la(V_+)$ as well.
The identifications~\eref{laStiden_e} induce identifications
\BE{SpinPinCorr_e0}\big((\tau_Y\!\oplus\!V)_{\pm},\fo_{\tau_Y\oplus V}^{\,\pm}\big)
=\big(\tau_Y\!\oplus\!V_{\pm},\St_{V_{\pm}}(\fo_V^{\pm})\big)
\equiv\St\big(V_{\pm},\fo_V^{\pm}\big).\EE

\begin{SpinPin}[Correspondences]\label{SpinPinCorr_prop}
Let $V$ be a real vector bundle over~$Y$.
There are natural $H^1(Y;\Z_2)$-equivariant bijections
\BE{SpinPinCorr_e}\Co_V^{\pm}\!:\cP^{\pm}(V)\lra\Sp\!\big(V_{\pm},\fo_V^{\pm}\big)\EE
so that 
$\Co_{\tau_Y\oplus V}^{\pm}\!\circ\!\St_V^{\pm}\!=\!\St_{V_{\pm}}\!\circ\!\Co_V^{\pm}$.
\end{SpinPin}

Suppose 
\BE{SpinPinSES_e0}0\lra V' \stackrel{\io}{\lra} V \stackrel{\fj}{\lra} V'' \lra 0\EE
is a short exact sequence of vector bundles over~$Y$.
Orientations~$\fo'$ on~$V'$ and $\fo''$ on~$V''$ determine an 
orientation~$\fo'_{\ce}\fo''$\nota{oeo@$\fo'_{\ce}\fo''$} on~$V$
as follows.
If $y\!\in\!Y$ and $v_1',\ldots,v_m'$ is an oriented basis for~$(V_y',\fo')$,
then 
$$\io(v_1'),\ldots,\io(v_m'),v_{m+1},\ldots,v_n\in V_y$$
is an oriented basis for $(V_y,\fo'_{\ce}\fo'')$ if and only if 
$\fj(v_{m+1}),\ldots,\fj(v_n)$ is an oriented basis for~$(V_y'',\fo'')$.
Whenever the vector bundle~$V'$ and~$V''$ are specified, we denote 
the canonical, direct short exact sequence
\BE{SpinPinDS_e0}0\lra V' \stackrel{\io}{\lra} V'\!\oplus\!V'' 
\stackrel{\fj}{\lra} V'' \lra 0\EE
by~$\oplus$.
If $\fo'$ and $\fo''$ are orientations on~$V'$ and $\fo''$ on~$V''$, respectively,
we denote by $\fo'\fo''$\nota{oo@$\fo'\fo''$} the induced orientations on~$V$.\\

If $\ce$ is the short exact sequence~\eref{SpinPinSES_e0} and~$\fo'$ is an orientation on~$V'$,
we denote by~$\ce_{\fo'}^{\pm}$ any short exact sequence
\BE{SpinPinSES_e0b}0\lra V' \stackrel{(\io,0)}{\lra} V_{\pm} 
\stackrel{\fj_{\fo'}^{\pm}}{\lra} V''_{\pm} \lra 0\EE
obtained by combining~\eref{SpinPinSES_e0} with the isomorphism
$$\la(V)\approx \la(V')\!\otimes\!\la(V'')\approx\la(V'')$$
induced by an orientation-preserving trivialization of~$(V',\fo')$.
If $Y$ is paracompact (as in the perspectives of 
Definitions~\ref{PinSpin_dfn} and~\ref{PinSpin_dfn2}), 
$\fo'$ determines a homotopy class of trivializations of~$V'$ and thus a homotopy
class of short exact sequences~$\ce_{\fo'}^{\pm}$ of vector bundles~\eref{SpinPinSES_e0b} over~$Y$.
For the purposes of Definition~\ref{PinSpin_dfn3}, it is most natural to define 
a \sf{short exact sequence} of vector bundles $V',V,V''$ over~$Y$ as a collection 
of homotopy classes of short exact sequences
$$0\lra \al^*V' \stackrel{\io_{\al}}{\lra} \al^*V \stackrel{\fj_{\al}}{\lra} \al^*V'' \lra 0$$
over~$S^1$ in the usual sense,
one homotopy class for each loop $\al$ in~$Y$, which extend to short exact sequences 
$$0\lra F^*V' \stackrel{\io_F}{\lra} F^*V \stackrel{\fj_F}{\lra} F^*V'' \lra 0$$
over bordered
surfaces~$\Si$ in the usual sense for each continuous map \hbox{$F\!:\Si\!\lra\!Y$}
as in Definition~\ref{PinSpin_dfn3}.
An orientation~$\fo'$ on~$V'$ then still determines a short exact sequence~$\ce_{\fo'}^{\pm}$
of vector bundles~\eref{SpinPinSES_e0b}.

\begin{SpinPin}[Short Exact Sequences]\label{SpinPinSES_prop}
Every short exact sequence~$\ce$ of vector bundles over~$Y$ as in~\eref{SpinPinSES_e0}
determines natural $H^1(Y;\Z_2)$-biequivariant maps
\BE{SpinPinSESdfn_e0}\begin{split}
\llrr{\cdot,\cdot}_{\ce}\!:\OSp(V')\!\times\!\OSp(V'')&\lra\OSp(V), \\
\llrr{\cdot,\cdot}_{\ce}\!:\OSp(V')\!\times\!\cP^{\pm}(V'')&\lra\cP^{\pm}(V)
\end{split}\EE
so that the following properties hold. 
\begin{enumerate}[label=(ses\arabic*),leftmargin=*]

\item\label{DSsplit_it} If oriented vector bundles $(V',\fo')$ and $(V'',\fo'')$ over~$Y$ 
are split as direct sums of oriented line bundles, then
$$\llrr{\os_0(V',\fo'),\os_0(V'',\fo'')}_{\oplus}
=\os_0\big((V',\fo')\!\oplus\!(V'',\fo'')\!\big).$$

\item\label{DSorient_it} If  $\ce$ is as in~\eref{SpinPinSES_e0},
$\fo'\!\in\!\fO(V')$, and $\fo''\!\in\!\fO(V'')$, then
\begin{alignat}{2}
\label{DSorient1_e}
\llrr{\os',\os''}_{\ce}&\in\Sp(V,\fo'_{\ce}\fo'')
&\quad&\forall~\os'\!\in\!\Sp(V',\fo'),\,\os''\!\in\!\Sp(V'',\fo''),\\
\label{DSorient2_e}
\fR_{\fo'_{\ce}\fo''}^{\pm}\big(\llrr{\os',\fp''}_{\ce}\big)
&=\bllrr{\os',\fR_{\fo''}^{\pm}(\fp'')}_{\ce}
&\quad&\forall~\os'\!\in\!\Sp(V',\fo'),\,\fp''\!\in\!\cP^{\pm}(V'').
\end{alignat}

\item\label{DSEquivSum_it} 
If $\ce$ is as in~\eref{SpinPinSES_e0}, $\fo'\!\in\!\fo(V')$, 
$\os'\!\in\!\Sp(V',\fo')$, and $\fp''\!\in\!\cP^{\pm}(V'')$, then
$$\llrr{\ov\os',\fp''}_{\ce}=w_1(V'')\!\cdot\!\llrr{\os',\fp''}_{\ce}, \quad
\Co_V^{\pm}\big(\llrr{\os',\fp''}_{\ce}\big)=
\bllrr{\os',\Co_{V''}^{\pm}(\fp'')}_{\ce_{\fo'}^{\pm}}\,.$$

\item\label{DSassoc_it} 
If $V_1',V_2',V''$ are vector bundles over~$Y$, $\os_1'\!\in\!\OSp(V_1')$, $\os_2'\!\in\!\OSp(V_2')$,  
and $\fp''\!\in\!\cP^{\pm}(V'')$, then
$$\bllrr{\os_1',\llrr{\os_2',\fp''}_{\oplus}\!}_{\oplus}
=\bllrr{\!\llrr{\os_1',\os_2'}_{\oplus},\fp''}_{\oplus}\,.$$

\item\label{DSstab_it} 
If $V$ is a vector bundle over~$Y$ and $\fp\!\in\!\cP^{\pm}(V)$, then
$$\St_V^{\pm}(\fp) =\llrr{\os_0(\tau_Y,\fo_Y),\fp}_{\oplus} .$$

\item\label{DSPin2Spin_it} 
If $V$ is a vector bundle over~$Y$, $\fo\!\in\!\fO(V)$, 
and $\fp\!\in\!\cP^{\pm}(V)$, then
$$\Co_V^{\pm}(\fp)=
\bllrr{\fR_{\fo}^{\pm}(\fp),\os_0((2\!\pm\!1)\la(V,\fo))}_{\oplus}\,.$$
\end{enumerate}

\end{SpinPin}

\vspace{.1in}

By SpinPin~\ref{Pin2SpinRed_prop}, \eref{DSorient2_e},
and the first statement in SpinPin~\ref{SpinPinSES_prop}\ref{DSEquivSum_it},
\BE{OspinDSEquivSum_e}  
\llrr{\ov\os',\os''}_{\ce}=\llrr{\os',\ov\os''}_{\ce}
=\ov{\llrr{\os',\os''}_{\ce}}\EE
for every short exact sequence~$\ce$ of vector bundles as in~\eref{SpinPinSES_e0}, 
$\os'\!\in\!\OSp(V')$, and $\os''\!\in\!\OSp(V'')$.
By SpinPin~\ref{Pin2SpinRed_prop} and~\ref{SpinPinSES_prop}\ref{DSassoc_it} 
and~\eref{DSorient2_e},
\BE{OspinDSassoc_e}  
\bllrr{\os_1',\llrr{\os_2',\os''}_{\oplus}\!}_{\oplus}
=\bllrr{\!\llrr{\os_1',\os_2'}_{\oplus},\os''}_{\oplus}\EE
for all vector bundles $V_1',V_2',V''$ over~$Y$, 
$\os_1'\!\in\!\OSp(V_1')$, $\os_2'\!\in\!\OSp(V_2')$,  and $\os''\!\in\!\OSp(V'')$.
By SpinPin~\ref{Pin2SpinRed_prop}, \ref{SpinPinStab_prop}, and~\ref{SpinPinSES_prop}\ref{DSstab_it}
and \eref{DSorient2_e}, 
\BE{OspinDSstab_e}\St_V(\os) =\llrr{\os_0(\tau_Y,\fo_Y),\os}_{\oplus}\EE
for every vector bundle~$V$ over~$Y$ and $\os\!\in\!\OSp(V)$.\\

Combining the second map in~\eref{SpinPinSESdfn_e0} with the canonical isomorphism
of $V'\!\oplus\!V''$ with $V''\!\oplus\!V'$, we obtain
a natural $H^1(Y;\Z_2)$-biequivariant map
\BE{SpinPinSESdfn_e0b}
\llrr{\cdot,\cdot}_{\ce}\!:\cP^{\pm}(V')\!\times\!\OSp(V'')\lra\cP^{\pm}(V)\,.\EE
By the SpinPin~\ref{SpinPinSES_prop} property, this map satisfies the obvious
analogues of \eref{DSorient2_e}, the first statement in SpinPin~\ref{DSEquivSum_it}, 
and SpinPin~\ref{DSassoc_it} and~\ref{DSPin2Spin_it}.

\begin{rmk}\label{ShExSeqCorr1_rmk}
The first statement of the SpinPin~\ref{SpinPinSES_prop}\ref{DSEquivSum_it} property
corrects \cite[Lemma~8.1]{Sol},
which suggests that the induced $\Pin^{\pm}$-structure on $V'\!\oplus\!V''$ 
does not depend on the orientation of~$V'$.
\end{rmk}

\subsection{Basic examples}
\label{EgOutline_subs}

The first two cases of the first homomorphism~\eref{qnSpinPin_e} are given~by
\BE{Spin1and2_e}
q_1\!:\Spin(1)\!\equiv\!\Z_2\lra\SO(1)\!=\!\{\1\}, \quad
q_2\!:\Spin(2)\!=\!S^1\lra\SO(2)\!=\!S^1,~~u\lra u^2\,.
\gena{Spin$(1)$, Spin$(2)$}\nota{q1a@$q_1,q_2$}\EE
The groups $\Pin^+(1)$ and $\Pin^-(1)$ are $\Z_2^{\,2}$ and $\Z_4$,
respectively, and 
\BE{q1minhm_e} q_1^-\!:\Pin^-(1)\!=\!\Z_4\lra 
\O(1)\!=\!\{\pm1\}\approx\!\Z_2\!\equiv\!\Z/2\Z\gena{Pin$^{\pm}(1)$}\nota{q1b@$q_1^{\pm}$}\EE
is the unique surjective homomorphism.
For concreteness, we take 
\BE{Pin1iden_e}
\wt\bI_{1;1}\!=\!1_{\Z_4}\!\equiv\!1\!+\!4\Z\in\Pin^-(1), \quad
\wt\bI_{1;1}\!=\!\big(0,1_{\Z_2}\big)\in\Pin^+(1), \nota{i11@$\wt\bI_{1;1}$}\quad 
\wh\bI_1\!=\!\big(1_{\Z_2},0\big)\in\Pin^+(1)\,.\nota{i1@$\wh\bI_1$}\EE
The last choice implies that 
\BE{q1pluhm_e} q_1^+\!\!=\!\pr_2\!:\Pin^+(1)\!=\!\Z_2\!\oplus\!\Z_2\lra 
\O(1)\!=\!\{\pm1\}\approx\!\Z_2\!\equiv\!\Z/2\Z
\nota{q1b@$q_1^{\pm}$}\gena{Pin$^{\pm}(1)$}\EE
is the projection to the second component.
We next explicitly describe $\Spin$- and $\Pin^{\pm}$-structures on 
some rank~1 and~2 vector bundles.

\begin{eg}\label{LBPin_eg}
Let $V$ be a real line bundle over a connected space~$Y$.
Every automorphism~$\Psi$ of~$V$ is homotopic to either the identity or 
the multiplication by~$-1$.
An automorphism~$\Psi$ of~$V$ homotopic to the identity clearly 
acts trivially on the sets $\cP^{\pm}(V)$ of $\Pin^{\pm}-$ 
and $\OSp(V)$ of $\OSpin$-structures on~$V$ in either of the three perspectives
of Definitions~\ref{PinSpin_dfn}-\ref{PinSpin_dfn3}.
By~\eref{SpinPinStr_e} with $V''$ of rank~0,
an automorphism~$\Psi$ of~$V$ homotopic to the multiplication by~$-1$
also acts trivially on~$\cP^{\pm}(V)$.
Thus, every automorphism $\Psi$ of~$V$ acts trivially on~$\cP^{\pm}(V)$.
If $\fo\!\in\!\fO(V)$, then every orientation-preserving automorphism of~$(V,\fo)$
acts trivially on the set~$\Sp(V,\fo)$ of $\Spin$-structures on~$(V,\fo)$.
\end{eg}

\begin{eg}\label{LBSpinStr_eg}
The  $\Spin$-structures on an oriented line bundle~$(V,\fo)$ over a path-connected 
paracompact space~$Y$
in the perspective of Definition~\ref{PinSpin_dfn}
are the topological double covers \hbox{$q_V\!:\wt{Y}\!\lra\!Y$}.
The distinguished element $\fs_0(V,\fo)$\nota{s0V@$\fs_0(V,\fo)$} 
of~$\Sp(V,\fo)$ is 
the disconnected double cover of~$Y$.
\end{eg}

Examples~\ref{Pin1pStr_eg}, \ref{Pin1mStr_eg}, and~\ref{Spin2Str_eg} below 
are the base inputs in the proof of the SpinPin~\ref{SpinPinObs_prop} property
for the perspective of Definition~\ref{PinSpin_dfn}.
This property for the perspective of Definition~\ref{PinSpin_dfn2} is 
a direct consequence of its validity in the perspective of Definition~\ref{PinSpin_dfn}.
The proof of the SpinPin~\ref{SpinPinObs_prop} property for the perspective of 
Definition~\ref{PinSpin_dfn3} is fundamentally different.

\begin{eg}\label{Pin1pStr_eg}
The real tautological line bundle $\ga_{\R}\lra\!\R\P^{\i}$\nota{gar@$\ga_{\R}$} admits two 
$\Pin^+$-structures in the sense of Definition~\ref{PinSpin_dfn}, up to equivalence.
Since $\O(\ga_{\R})\!=\!S^{\i}$ is simply connected, the domain $\Pin^+(V)$ of~$q_V$ 
is $\Z_2\!\times\!S^{\i}$ and
$$q_V\!:\Pin^+(V)\!=\!\Z_2\!\times\!S^{\i}\lra \O(\ga_{\R})\!=\!S^{\i}$$
is the projection to the second factor.
This map is equivariant with respect to the homomorphism~\eref{q1pluhm_e}
if the first $\Z_2$ factor  in $\Pin^+(1)$ acts by addition 
on the first factor of~$\Pin^+(V)$, leaving the second factor of~$\Pin^+(V)$ fixed, 
and the non-trivial element $1_{\Z_2}$ of the second $\Z_2$ factor in $\Pin^+(1)$
acts by the antipodal map on the second factor of~$\Pin^+(V)$,
either leaving the first factor fixed or interchanging its two elements.
Since any automorphism of the bundle 
$$\Pin^+(V)\!=\!\Z_2\!\times\!S^{\i}\lra \R\P^{\i}$$
is equivariant with respect to either $\Pin^+(1)$-action, 
the two $\Pin^+$-structures are not equivalent. 
\end{eg}

\begin{eg}\label{Pin1mStr_eg}
The line bundle $\ga_{\R}\lra\!\R\P^{\i}$ does not admit a $\Pin^-$-structure
in the sense of Definition~\ref{PinSpin_dfn}.
Since $\O(\ga_{\R})\!=\!S^{\i}$ is simply connected, the domain $\Pin^-(V)$ of~$q_V$ 
would have to be $\Z_2\!\times\!S^{\i}$.
The restriction of~$q_V$ to each connected component of $\Pin^-(V)$ 
would then have to be a homeomorphism onto~$S^{\i}$ and equivariant
with respect to the $\Z_2$ subgroup of $\Pin^-(1)\!\approx\!\Z_4$ 
preserving each connected component of~$\Pin^-(V)$.
This is impossible because this subgroup lies in the kernel of 
the homomorphism~\eref{q1minhm_e}.
\end{eg}

\begin{eg}\label{Pin1mStr1_eg}
The real tautological line bundle $\ga_{\R;1}\lra\!\R\P^1$\nota{gar1@$\ga_{\R;1}$} admits two
$\Pin^-$-structures  in the sense of Definition~\ref{PinSpin_dfn}.
By the reasoning in Example~\ref{Pin1mStr_eg}, the projection~$q_V$ 
is the connected double~cover
$$q_V\!:\Pin^-\big(\ga_{\R;1}\big)\!=\!S^1\lra \O(\ga_{\R;1})\!=\!S^1, \quad
q_V(w)=w^2~~\forall\,w\!\in\!S^1\!\subset\!\C.$$
This map is equivariant with respect to the homomorphism~\eref{q1minhm_e}
if $1_{\Z_4}\!\in\!\Z_4$ acts by the multiplication by either $\ne^{\fI\pi/2}$
or~$\ne^{-\fI\pi/2}$.
Since any automorphism of the bundle 
$$\Pin^-\big(\ga_{\R;1}\big)\!=\!S^1\lra \R\P^1$$
is a multiplication by a power of $\ne^{\fI\pi/2}$,
the two $\Pin^-$-structures are not equivalent.
We denote by~$\fp_0^-(\ga_{\R;1})$\nota{p01mgaR1@$\fp_0^-(\ga_{\R;1}),\fp_1^-(\ga_{\R;1})$} 
the $\Pin^-$-structure on~$\ga_{\R;1}$
in which $1_{\Z_4}\!\in\!\Z_4$ acts on $\Pin^-(\ga_{\R;1})$ by the multiplication by~$\ne^{\fI\pi/2}$ 
and by~$\fp_1^-(\ga_{\R;1})$ the other $\Pin^-$-structure on~$\ga_{\R;1}$. 
\end{eg}

\begin{eg}\label{Spin2Str_eg}
The complex tautological line bundle $\ga_{\C}\lra\!\C\P^{\i}$\nota{gac@$\ga_{\C}$} 
does not admit a $\Spin$-structure  in the sense of Definition~\ref{PinSpin_dfn}.
Since $\Spin(2)\!=\!S^1$ and $\SO(\ga_{\C})\!=\!S^{\i}$ are path-connected,
the domain $\Spin(V)$ of~$q_V$ would have to be path-connected.
However, this is impossible because $\SO(\ga_{\C})$ is simply connected and
thus does not admit any path-connected double covers.
\end{eg}

The perspectives of Definitions~\ref{PinSpin_dfn2} and~\ref{PinSpin_dfn3} 
on $\Spin$- and $\Pin^{\pm}$-structures
make it feasible to specify them explicitly in some cases central to
the open and real sectors of Gromov-Witten theory, mirror symmetry, and enumerative geometry.
This is indicated by the usage of the examples below later in this monograph.
Specifications of topological data on a target manifold analogous to
the description of Example~\ref{RP2pin_eg} below lead to explicit determinations
of orientations of moduli spaces of open and real maps and signed counts of curves
in \cite[Sections~5,6]{Teh} and \cite[Section~5]{RealGWsIII};
such computations are notoriously difficult to carry out correctly.

\begin{eg}\label{CanSpin_eg}
Let $V$ be a vector bundle over a topological space~$Y$.
\begin{enumerate}[label=(\alph*),leftmargin=*]

\item\label{CanSpin_it1} The vector bundle $V\!\oplus\!V$ has a canonical orientation.
Its restriction to the fiber $V_y\!\oplus\!V_y$ of $V\!\oplus\!V$
over a point $y\!\in\!Y$ is obtained by choosing any orientation in the first copy
of~$V_y$ and then taking the same orientation in the second copy.

\item\label{CanSpin_it2} If $V$ is oriented, then the oriented vector bundle $V\!\oplus\!V$
has a canonical $\Spin$-structure in the sense of Definition~\ref{PinSpin_dfn3}.
The canonical homotopy class of trivializations of $V\!\oplus\!V$ over a loop~$\al$
in~$Y$ is obtained by choosing any trivialization of the first copy of~$\al^*V$ and 
then taking the same trivialization of the second copy.
\end{enumerate}
\end{eg}

\vspace{.1in}

The standard metric on $\R^{n+1}$ determines an identification of the real tautological line 
bundle
$$\ga_{\R;n}\equiv 
\big\{\big(\ell,v\big)\!\in\!\R\P^n\!\times\!\R^{n+1}\!:\,v\!\in\!\ell\!\subset\!\R^{n+1}\big\}
\nota{garn@$\ga_{\R;n}$} $$ 
over~$\R\P^n$ with its dual~$\ga_{\R;n}^*$.
This identification and Euler's exact sequence,
\BE{RPnEES_e3}0\lra \tau_{\R\P^n} \stackrel{f}{\lra} (n\!+\!1)\ga_{\R;n}^*
\stackrel{g}{\lra} T(\R\P^n) \lra 0,\EE
of vector bundles over $\R\P^n$ determines a canonical homotopy class  of isomorphisms
\BE{RPnEES_e3b}\tau_{\R\P^n}\oplus T(\R\P^n) \approx (n\!+\!1)\ga_{\R;n}\EE
of vector bundles over~$\R\P^n$.
The bundle homomorphisms in~\eref{RPnEES_e3} can be described explicitly 
as in the proof of \cite[Lemma~2.1]{RealGWsIII}.

\begin{eg}\label{ga1Rpin_eg0}
We denote by $\prt_{\th}$ the vector field on $\R\P^1\!\equiv\!S^1/\Z_2$ induced 
by the vector field $(-y,x)$ on~$\R^2$ and by $\fo_{T\R\P^1}$ 
the induced orientation of~$\R\P^1$.
Under the standard identification of~$\ga_{\R;1}$ with~$\ga_{\R;1}^*$,
the homomorphisms in the $n\!=\!1$ case of~\eref{RPnEES_e3} are given~by
\begin{gather*}
f\big([x,y],a\big)=\bigg(\![x,y],a\frac{x(x,y)}{x^2\!+\!y^2},a\frac{y(x,y)}{x^2\!+\!y^2}\bigg)
\in \R\P^1\!\times\!\R^2\!\times\!\R^2
\quad\forall\,\big([x,y],a\big)\!\in\!\R\P^1\!\times\!\R,\\
g\big([x,y],(x_0,y_0),(x_1,y_1)\!\big)=
\frac{x(xx_1\!+\!yy_1)\!-\!y(xx_0\!+\!yy_0)}{x^2\!+\!y^2}\prt_{\th}\big|_{[x,y]}\,.
\end{gather*}
In particular,
$$f\big([1,0],1\big)=\big([1,0],(1,0),(0,0)\!\big) \quad\hbox{and}\quad
g\big([1,0],(0,0),(1,0)\!\big)=\prt_{\th}\big|_{[1,0]}\,.$$
Thus, the orientation on $2\ga_{\R;1}$ induced by the orientations
$\fo_{\R\P^1}$ of~$\tau_{\R\P^1}$ and $\fo_{T\R\P^1}$ of~$\R\P^1$ 
is the canonical orientation~$\fo_{\ga_{\R;1}}^-$ 
provided by 
Example~\ref{CanSpin_eg}\ref{CanSpin_it1}.
\end{eg}

\begin{eg}\label{ga1Rpin_eg}
Under the standard identification $\R^2\!=\!\C$, the bundle homomorphism
\BE{ga1Rpin_e} \Phi_0\!:2\ga_{\R;1}\lra \R\P^1\!\times\!\C, \quad
\big(\ell,c_1,c_2\big)\lra \big(\ell,c_1\!+\!\fI c_2\big), \EE
is a trivialization of the vector bundle $2\ga_{\R;1}$ over~$\R\P^1$;
this is the trivialization \cite[(5.11)]{RealGWsI}.
The homotopy class of
the trivialization $\id\!\oplus\!\Phi_0$ of $\tau_{\R\P^1}\!\oplus\!2\ga_{\R;1}$
is a $\Spin$-structure $\os_0(2\ga_{\R;1},\fo_{\ga_{\R;1}}^-)$
 on $(2\ga_{\R;1},\fo_{\ga_{\R;1}}^-)$ 
in the perspective of Definition~\ref{PinSpin_dfn3}
and a $\Pin^-$-structure on~$\ga_{\R;1}$ in this perspective.
By the SpinPin~\ref{SpinPinStr_prop}\ref{SpinStr_it} property,
the group $H^1(\R\P^1;\Z_2)$ acts freely and transitively 
on the set $\Sp(2\ga_{\R;1},\fo_{\ga_{\R;1}}^-)$ of $\Spin$-structures on 
$(2\ga_{\R;1},\fo_{\ga_{\R;1}}^-)$. 
There is thus one other $\Spin$-structure $\os_1(2\ga_{\R;1},\fo_{\ga_{\R;1}}^-)$
on this oriented vector bundle.
Example~\ref{SpinDfn1to3_eg} identifies $\os_0(2\ga_{\R;1},\fo_{\ga_{\R;1}}^-)$
and $\os_1(2\ga_{\R;1},\fo_{\ga_{\R;1}}^-)$ 
with the $\Pin^-$-structures 
$\fp_0^-(\ga_{\R;1})$\nota{p01mgaR1@$\fp_0^-(\ga_{\R;1}),\fp_1^-(\ga_{\R;1})$}
and $\fp_1^-(\ga_{\R;1})$, respectively, of Example~\ref{Pin1mStr1_eg}.
\end{eg}

\begin{eg}\label{RP2pin_eg}
By the SpinPin~\ref{SpinPinStab_prop} and~\ref{SpinPinCorr_prop} properties
and $n\!=\!2$ case of~\eref{RPnEES_e3b}, 
the $\Pin^-$-structures on (the tangent bundle of) $\R\P^2$ 
correspond to the $\Spin$-structures on the canonically oriented vector~bundle
$$\big(\St(T(\R\P^2)\!)\!\big)_-
\!\equiv\! \tau_{\R\P^2}\!\oplus\!T(\R\P^2)
\oplus \la\big(\tau_{\R\P^2}\!\oplus\!T(\R\P^2)\big)
\approx 3\ga_{\R;2}\oplus\ga_{\R;2}^{\otimes3}\,.$$
Since the line bundle $\ga_{\R;2}^{\otimes2}$ has a canonical trivialization,
there is thus a canonical homotopy class of isomorphisms 
$$\tau_{\R\P^2}\!\oplus\!T(\R\P^2)
\oplus \la\big(\tau_{\R\P^2}\!\oplus\!T(\R\P^2)\big)
\approx 4\ga_{\R;2}\,.$$
By the conclusion of Example~\ref{ga1Rpin_eg0}, the orientation on~$4\ga_{\R;2}$ 
induced via this homotopy class of isomorphism is the orientation~$\fo_{\ga_{\R;2}}^+$
provided by  Example~\ref{CanSpin_eg}\ref{CanSpin_it1}.
Example~\ref{CanSpin_eg}\ref{CanSpin_it2} endows this oriented vector bundle
with a canonical $\OSpin$-structure $\os_0(4\ga_{\R;2},\fo_{\ga_{\R;2}}^+)$
in the sense of Definition~\ref{PinSpin_dfn3}.
By the SpinPin~\ref{SpinPinStr_prop} property,
$(4\ga_{\R;2},\fo_{\ga_{\R;2}}^+)$ admits precisely one other $\Spin$-structure 
$\os_1(4\ga_{\R;2},\fo_{\ga_{\R;2}}^+)$.
The associated homotopy class of trivializations over a loop~$\al$ in~$\R\P^2$
is obtained by combining trivializations of the two copies of $\al^*(2\ga_{\R;2})$
in $\al^*(4\ga_{\R;2})$ that differ by one rotation.
By Corollary~\ref{RP2pin_crl}, the $\Pin^-$-structures 
$\fp_0^-(\R\P^2)$\nota{p01rp2@$\fp_0^-(\R\P^2),\fp_1^-(\R\P^2)$} 
and $\fp_1^-(\R\P^2)$ on $T\R\P^2$ associated to 
$\os_0(4\ga_{\R;2},\fo_{\ga_{\R;2}}^+)$ and $\os_1(4\ga_{\R;2},\fo_{\ga_{\R;2}}^+)$,
respectively, via~\eref{SpinPinCorr_e}
 are fixed by every automorphism of~$\R\P^2$.
\end{eg}

\begin{eg}\label{RP3spin_eg}
We denote by $\os_0(4\ga_{\R;3},\fo_{\ga_{\R;3}}^+)$ the canonical $\OSpin$-structure
on the vector bundle~$4\ga_{\R;3}$ over~$\R\P^3$ in the sense of Definition~\ref{PinSpin_dfn3}
provided by Example~\ref{CanSpin_eg}\ref{CanSpin_it2}.
By the SpinPin~\ref{SpinPinStr_prop} property,
$(4\ga_{\R;3},\fo_{\ga_{\R;3}}^+)$ admits precisely one other $\Spin$-structure 
$\os_1(4\ga_{\R;3},\fo_{\ga_{\R;3}}^+)$.
By the $n\!=\!3$ case of~\eref{RPnEES_e3b}, 
these two $\OSpin$-structures correspond to $\OSpin$-structures on~$\R\P^3$;
we denote the latter by $\os_0(\R\P^3)$ and~$\os_1(\R\P^3)$, respectively.
\end{eg}

\subsection{Further examples}
\label{EgOutline_subs2}

We now combine Examples~\ref{LBPin_eg}-\ref{Spin2Str_eg} with 
the properties of Spin- and Pin-structures collected in Section~\ref{SpinPinProp_subs}
to explore these structures on real line bundles more systematically.\\

As described in Section~\ref{DfnThm_subs}, an orientation $\fo\!\in\!\fO(V)$
on a line bundle~$V$ over a topological space~$Y$ determines an element 
$\fs_0(V,\fo)\!\in\!\Sp(V,\fo)$.
We denote its preimage under~\eref{Pin2SpinRed_e} by 
\hbox{$\fp_0^{\pm}(V)\!\in\!\cP^{\pm}(V)$}\nota{p0Vpm@$\fp_0^{\pm}(V)$}.
By the identity after~\eref{Pin2SpinRed_e}, $\fp_0^{\pm}(V)$ does not depend on the choice of~$\fo$.
By Example~\ref{Pin1pStr_eg}, the real tautological line bundle  
$\ga_{\R}\lra\!\R\P^{\i}$ admits two $\Pin^+$-structures. 
By \cite[Theorem~5.6]{MiSt}, for every real line bundle~$V$ over a paracompact space~$Y$
there exists a continuous map $f\!:Y\!\lra\!\R\P^{\i}$ such that 
$V\!\approx\!f^*\ga_{\R}$; this map is unique up to homotopy.
By~\eref{SpinPinStr_e}, the $\Pin^+$-structure on~$V$
induced from a $\Pin^+$-structure on~$\ga_{\R}$ does not depend on the choice of this isomorphism.

\begin{lmm}\label{PinpCan_lmm}
Let $V$ be a line bundle over a topological space~$Y$ and $\fo\!\in\!\fO(V)$.
If $f\!:Y\!\lra\!\R\P^{\i}$ is a continuous map such that 
$V\!\approx\!f^*\ga_{\R}$, then
the $\Pin^+$-structure~$\fp_0^+(V)$ on~$V$ corresponding to 
the canonical $\Spin$-structure $\fs_0(V,\fo)$ 
on~$(V,\fo)$ is the pullback by~$f$ of either $\Pin^+$-structure on~$\ga_{\R}$.
\end{lmm}

\begin{proof} Let $q\!:S^{\i}\!\lra\!\R\P^{\i}$ be the standard quotient projection.
Since $S^{\i}$ is simply connected,
$$\tau_{S^{\i}}\!\equiv\!S^{\i}\!\times\!\R\!=\!q^*\ga_{\R}\lra S^{\i}\,.$$
Since $V$ is orientable, there exists a continuous map $\wt{f}\!: Y\!\lra\!S^{\i}$
so that $q\!\circ\!\wt{f}$ is homotopic to~$f$ and
$$(V,\fo)\approx\wt{f}^*\big(\tau_{S^{\i}},\fo_{S^{\i}}\big)\lra Y$$
as oriented line bundles over~$Y$.
The canonical $\Spin$-structure~$\fs_0(V,\fo)$ on~$V$ is the pullback by~$\wt{f}$
of~$\fs_0(\tau_{S^{\i}},\fo_{S^{\i}})$.
By the SpinPin~\ref{SpinPinStr_prop}\ref{SpinStr_it} property (or Example~\ref{LBSpinStr_eg}), 
there are no other $\Spin$-structures on~$\tau_{S^{\i}}$.
By the  SpinPin~\ref{Pin2SpinRed_prop} property,
the pullback by~$q$ of either $\Pin^+$-structure on~$\ga_{\R}$ is thus
 the $\Pin^+$-structure~$\fp_0^+(\tau_{\R})$ on~$\tau_{S^{\i}}$ 
corresponding to~$\fs_0(\tau_{S^{\i}},\fo_{S^{\i}})$. 
Since the pullback by~$f$ of a $\Pin^+$-structure~$\fp$ on~$\ga_{\R}$
is the pullback by~$\wt{f}$ of the pullback of~$q$ of~$\fp$,
the claim now follows from the naturality of 
the correspondence of the SpinPin~\ref{Pin2SpinRed_prop} property with respect 
to continuous maps.
\end{proof}

\begin{lmm}\label{Pin1_lmm}
A real line bundle~$V$ over a topological space~$Y$ carries a canonical 
$\Pin^+$-structure $\fp_0^+(V)$.\nota{p0Vpm@$\fp_0^{\pm}(V)$}
For every unorientable line bundle~$V$ over an oriented circle~$Y$,
there is a natural bijection
\BE{Pin1Unor_e} \cP^-(V)\lra\Z_2;\EE
reversing the orientation of~$Y$ flips this bijection.
\end{lmm}

\begin{proof}
By the SpinPin~\ref{SpinPinObs_prop}\ref{PinObs_it} and~\ref{SpinPinStr_prop}\ref{PinStr_it} 
properties  (or Example~\ref{Pin1pStr_eg}),
there are two $\Pin^+$-structure
on the real tautological line bundle $\ga_{\R}$ over~$\R\P^{\i}$.
In the perspective of Definition~\ref{PinSpin_dfn} and Example~\ref{Pin1pStr_eg},
we take the one in which the non-trivial element~$1_{\Z_2}$ of the second $\Z_2$ factor 
in $\Pin^+(1)$ fixes the topological components of~$\Pin^+(V)$ 
to be the canonical $\Pin^+$-structure~$\fp_0^+(\ga_{\R})$ on~$\ga_{\R}$.
Given a real line bundle~$V$ over a paracompact space~$Y$,
let  $f\!:Y\!\lra\!\R\P^{\i}$ be a continuous map such that 
$V\!\approx\!f^*\ga_{\R}$; this map is unique up to homotopy.
We take the canonical equivalence class~$\fp_0^+(V)$ of $\Pin^+$-structures on~$V$
to be the pullback by~$f$ of~$\fp_0^+(\ga_{\R})$ in the perspective of  Definition~\ref{PinSpin_dfn}.
By~\eref{SpinPinStr_e}, $\fp_0^+(V)$ does not depend on the choice of 
isomorphism of $V$ with~$f^*\ga_{\R}$.
This also specifies~$\fp_0^+(\ga_{\R})$ in the perspective of  Definition~\ref{PinSpin_dfn3}
for an arbitrary topological space~$Y$.\\
We view $\R\P^1$ as the quotient of the unit circle $S^1\!\subset\!\C$ 
by the antipodal~map.
By the SpinPin~\ref{SpinPinObs_prop}\ref{PinObs_it} and~\ref{SpinPinStr_prop}\ref{PinStr_it} 
properties (or Example~\ref{Pin1mStr1_eg}), there are two $\Pin^-$-structures
on the real tautological line bundle $\ga_{\R;1}$ over~$\R\P^1$.
An oriented circle~$Y$ can be identified with~$\R\P^1$ as oriented manifold;
such an identification is unique up to homotopy.
An unorientable line bundle~$V$ over~$\R\P^1$ is isomorphic to~$\ga_{\R;1}$;
every two such isomorphisms differ by an automorphism of~$\ga_{\R;1}$.
By~\eref{SpinPinStr_e}, 
it thus remains to show that the orientation-reversing map
$$f\!:\R\P^1\lra\R\P^1, \qquad f\big([z]\big)=\big[\ov{z}\big]~~\forall\,z\!\in\!S^1,$$
interchanges the two $\Pin^-$-structures on~$\ga_{\R;1}$.\\

Under the standard identifications $\R\P^1\!=\!S^1/\Z_2$, $S^1\!\subset\!\C$, and
$\ga_{\R;1}\!\subset\!\R\P^1\!\times\!\C$,
\begin{gather*}
f^*\ga_{\R;1}=
\big\{([z],v)\!\in\!\R\P^1\!\!\times\!\C\!:v\!\in\!\R\ov{z}\big\},~
\O(f^*\ga_{\R;1})=f^*\O(\ga_{\R;1})=
\big\{([z],v)\!\in\!\R\P^1\!\!\times\!\C\!:v\!\in\!\{\pm\ov{z}\}\!\big\},\\
\pi_2^*\Pin^-(\ga_{\R;1})=
\big\{([z],v,w)\!\in\!\R\P^1\!\!\times\!\C\!\!\times\!S^1\!:v\!\in\!\{\pm\ov{z}\},\,w^2\!=\!v\big\},
\end{gather*}
where $\pi_2\!:f^*\O(\ga_{\R;1})\!\lra\!\O(\ga_{\R;1})$ is the projection to
the second component.
The map 
$$F\!:f^*\ga_{\R;1}\lra \ga_{\R;1}, \qquad F\big([z],v\big)\!=\!\big([z],\ov{v}\big),$$
is an isomorphism of line bundles over~$\R\P^1$  and
$$\big\{F|_{\O(f^*\ga_{\R;1})}\big\}^{\!*}\Pin^-\big(\ga_{\R;1}\big)
\!\equiv\!\big\{([z],v,w)\!\in\!\R\P^1\!\!\times\!\C\!\times\!S^1\!\!:
v\!\in\!\{\pm\ov{z}\},\,w^2\!=\!\ov{v}\big\}.$$
The~map 
$$\Psi\!: \pi_2^*\Pin^-(\ga_{\R;1})\lra 
\big\{F|_{\O(f^*\ga_{\R;1})}\big\}^{\!*}\Pin^-\big(\ga_{\R;1}\big), \quad
\Psi\big([z],v,w\big)=\big([z],v,\ov{w}\big),$$
is an isomorphism of $\Z_4$-bundles over~$\R\P^1$.
This map is $\Z_4$-equivariant if $1_{\Z_4}\!\in\!\Z_4$ acts on $\Pin^-(\ga_{\R;1})$
in the domain of~$\Psi$ by the multiplication by~$\ne^{\fI\pi/2}$ 
and on $\Pin^-(\ga_{\R;1})$
in the target of~$\Psi$ by the multiplication by~$\ne^{-\fI\pi/2}$.
Thus, the pullback by~$f$ interchanges the two $\Pin^-$-structures on~$\ga_{\R;1}$.
\end{proof}

\begin{crl}\label{Spin2Pin_crl2}
Let $V$ be an unorientable line bundle over~$S^1$.
If $f\!:S^1\!\lra\!S^1$ is a continuous map of even degree,
then the pullbacks by~$f$ of the two $\Pin^-$-structures on~$V$ are the same.
If the degree of~$f$ is divisible by~4, then these pullbacks are
the $\Pin^-$-structure~$\fp_0^-(f^*V)$ corresponding to~$\fs_0(f^*V,\fo)$
for any $\fo\!\in\!\fO(f^*V)$.
If the degree of~$f$ is even, but not divisible by~4, then
this is not the~case.
\end{crl}

\begin{proof} Let $f\!:S^1\!\lra\!S^1$ be a continuous map of even degree.
We can assume that it is given~by
$$f\!:S^1\lra S^1, \qquad z\lra z^{2d},$$
for some $d\!\in\!\Z$.
By Example~\ref{LBPin_eg}, we can assume that $V\!=\!\ga_{\R;1}$.
By Examples~\ref{Pin1mStr1_eg} and~\ref{LBSpinStr_eg}, 
the pullback by~$f$ of a $\Pin^-$-structure on~$V$ is~$\fp_0^-(f^*V)$ if and only~if
the total space of the~cover 
$$f^*\Pin^-\big(\ga_{\R;1}\big)\!\equiv\!\big\{(z,w)\!\in\!S^1\!\times\!S^1\!:
z^{2d}\!=\!w^4\big\}\lra S^1, \qquad (z,w)\lra z,$$
has 4 connected components. This is the case if and only if $d$ is even.
\end{proof}

\begin{rmk}\label{Pin1_rmk0} 
By the first statement of Lemma~\ref{Pin1_lmm} and the definition above Lemma~\ref{PinpCan_lmm},
the canonical equivalence classes of $\Pin^+$-structures on real line bundles
and of $\Pin^-$ and $\Spin$-structures on orientable line bundles
are preserved by the pullbacks by continuous maps.
By the second statement of Lemma~\ref{Pin1_lmm} and Corollary~\ref{Spin2Pin_crl2}, 
a continuous map~$f$ from~$S^1$ to itself need not pull back 
the canonical equivalence class~$\fp_0^-(V)$ of $\Pin^-$-structures
on an unorientable line bundle~$V$ (i.e.~the class corresponding to $0\!\in\!\Z_2$)
to~$\fp_0^-(f^*V)$.
\end{rmk}

\begin{rmk}\label{Pin1_rmk}
By the second statement of Lemma~\ref{Pin1_lmm}, the second claim of \cite[Lemma~8.2]{Sol}
is not correct in the $\Pin^-$-case:
there is no canonical way to choose a $\Pin^-$-structure on an unorientable real
line bundle~$V$ over a paracompact space~$Y$ that admits such a structure.
The existence of such a way would have implied that every real line bundle admits 
a $\Pin^-$-structure; this is not the case by Example~\ref{Pin1mStr_eg}.
However, there is a canonical way to choose a $\Pin^-$-structure on an unorientable real
line bundle~$V$ over an {\it oriented} circle~$Y$.
This is analogous to the phenomenon described by \cite[Corollary~5.6]{RealGWsI}.
\end{rmk}

\begin{eg}\label{Pin1pMB_eg}
The projection
$$\pi_V\!:V\!=\!\big([0,1]\!\times\!\R\big)\!\big/\!\!\sim\lra \R\P^1\!=\!S^1/\Z_2, \quad
(0,-a)\!\sim\!(1,a),~~\pi_V\big([t,a]\big)=\big[\ne^{\pi\fI t}\big],$$ 
defines a real line bundle over~$\R\P^1$.
The associated orthogonal frame bundle~is 
$$\O(V)\!=\!\big([0,1]\!\times\!\Z_2\big)\!\big/\!\!\sim\lra \R\P^1, \quad
\big(0,a\!+\!1_{\Z_2}\big)\sim(1,a), \quad [t,a]\lra\big[\ne^{\pi\fI t}\big].$$
The two $\Pin^+$-structures on~$V$ in the perspective of Definition~\ref{PinSpin_dfn}
are 
\begin{alignat*}{3}
q_{V;0}^+\!:\Pin_0^+(V)&\!\equiv\!\big([0,1]\!\times\!\Z_2^{\,2}\big)\!\big/\!\!\sim\lra\O(V),
&~\big(0,(a,b\!+\!1)\!\big)&\!\sim\!\big(1,(a,b)\!\big), &~
q_{V;0}^+\big([t,(a,b)]\big)&=[t,b],\\
q_{V;1}^+\!:\Pin_1^+(V)&\!\equiv\!\big([0,1]\!\times\!\Z_2^{\,2}\big)\!\big/\!\!\sim\lra\O(V),
&~\big(0,(a\!+\!1,b\!+\!1)\!\big)&\!\sim\!\big(1,(a,b)\!\big), &~ 
q_{V;1}^+\big([t,(a,b)]\big)&=[t,b].
\end{alignat*}
The group~$\Pin^+(1)$ acts on the two bundles via the addition on the $\Z_2^{\,2}$-factor,
making  the projections $q_{V;0}^+$ and $q_{V;1}^+$ above equivariant with respect to
the homomorphism~\eref{q1pluhm_e}.
Since the action of $\wt\bI_{1;1}\!\equiv\!(0,1_{\Z_2})$ preserves the connected components
of $\Pin_0^+(V)$, this bundle provides the canonical $\Pin^+$-structure~$\fp_0^+(V)$
on~$V$ of Lemma~\ref{Pin1_lmm}.
\end{eg}

\begin{eg}\label{Pin1mMB_eg}
Let $V$ be as in Example~\ref{Pin1pMB_eg}.
The two $\Pin^-$-structures on~$V$ in the perspective of Definition~\ref{PinSpin_dfn}
are 
\BE{Pin1mMB_e}\begin{aligned}
q_{V;0}^-\!:\Pin_0^-(V)&\!\equiv\!\big([0,1]\!\times\!\Z_4\big)\!\big/\!\!\sim\lra\O(V),
&~\big(0,a\!+\!1_{\Z_4}\!\big)&\!\sim\!(1,a), &~
q_{V;0}^-\big([t,a]\big)&=[t,a],\\
q_{V;1}^-\!:\Pin_1^-(V)&\!\equiv\!\big([0,1]\!\times\!\Z_4\big)\!\big/\!\!\sim\lra\O(V),
&~\big(0,a\!-\!1_{\Z_4}\!\big)&\!\sim\!(1,a), &~
q_{V;1}^-\big([t,a]\big)&=[t,a].
\end{aligned}\EE
The group~$\Pin^-(1)$ acts on $\Pin^-(V)$ via the addition on the $\Z_4$-factor,
making the projections $q_{V;0}^-$ and $q_{V;1}^-$ above equivariant with respect to
the homomorphism~\eref{q1minhm_e}.
The $\Pin^-$-structures on the real tautological line bundle~$\ga_{\R;1}$ over~$\R\P^1$
are described in Example~\ref{Pin1mStr1_eg}.
The~map
\BE{Pin1mMBiden_e}F\!:V\lra\ga_{\R;1}, \qquad 
F\big([t,a]\big)=\big([\ne^{\pi\fI t}],a\ne^{\pi\fI t}\big),\EE
is an isomorphism of line bundles over~$\R\P^1$.
The induced isomorphisms between the orthogonal frame bundles
and $\Pin^-$-structures of~$V$ and~$\ga_{\R;1}$  are given~by 
\begin{alignat*}{3}
\O(V)&\lra\O\big(\ga_{\R;1}\big),&\qquad 
\Pin_0^-(V)&\lra \Pin_0^-\big(\ga_{\R;1}\big),  &\qquad
\Pin_1^-(V)&\lra \Pin_1^-\big(\ga_{\R;1}\big),  \\
[t,a]&\lra (-1)^a\ne^{\pi\fI t}, &\qquad 
[t,a]&\lra \ne^{\pi\fI a/2}\ne^{\pi\fI t/2}, &\qquad 
[t,a]&\lra \ne^{-\pi\fI a/2}\ne^{\pi\fI t/2}\,.
\end{alignat*}
The middle and last maps above are equivariant with respect to the $\Z_4$-actions 
on the right-hand sides with $\wt\bI_{1;1}\!\equiv\!1_{\Z_4}$ 
acting by the multiplication by~$\ne^{\fI\pi/2}$ and~$\ne^{-\fI\pi/2}$, respectively. 
Thus, the first $\Pin^-$-structure on~$V$ in~\eref{Pin1mMB_e} corresponds to $0\!\in\!\Z_2$
in the sense of Lemma~\ref{Pin1_lmm}.
\end{eg}

\begin{eg}\label{PinpMB_eg}
Let $Y$ be the infinite Mobius band, i.e.~the total space of 
the tautological real line bundle~$\ga_{\R;1}$ over~$\R\P^1$,
and $S\ga_{\R;1}\!\subset\!Y$ be its unit circle bundle.
The bundle projection $\pi\!:Y\!\lra\!\R\P^1$ induces an exact sequence
\BE{ShExSeqImm0b_e1}0\lra \pi^*\ga_{\R;1}\lra TY \stackrel{\nd\pi}\lra \pi^*T\R\P^1\lra 0\EE
of vector bundles over~$Y$ and lifts an orientation~$\fo_{T\R\P^1}$
on~$\R\P^1$ to an orientation~$\fo_{TS\ga_{\R;1}}$ of~$S\ga_{\R;1}$. 
The inclusion of~$S\ga_{\R;1}$ into~$Y$ induces an exact sequence
\BE{ShExSeqImm0b_e3}0\lra T(S\ga_{\R;1})\lra TY|_{S\ga_{\R;1}}\lra 
\pi^*\ga_{\R;1}\big|_{S\ga_{\R;1}}\lra 0\,.\EE
We denote by $\os_0(\R\P^1,\fo_{T\R\P^1})$ the distinguished OSpin-structure
on~$T\R\P^1$ determined by~$\fo_{T\R\P^1}$ as in Example~\ref{LBSpinStr_eg} 
and by $\os_0(S\ga_{\R;1},\fo_{TS\ga_{\R;1}})$ the distinguished OSpin-structure
on~$T(S\ga_{\R;1})$ determined by~$\fo_{TS\ga_{\R;1}}$.
Let 
$$\fp_0^+\equiv\bllrr{\pi^*\fp_0^+(\ga_{\R;1}),
\pi^*\os_0(\R\P^1,\fo_{T\R\P^1})}_{\eref{ShExSeqImm0b_e1}}$$
be the $\Pin^+$-structure on $TY$ induced by $\fp_0^+(\ga_{\R;1})$ and 
$\os_0(\R\P^1,\fo_{T\R\P^1})$ via the exact sequence~\eref{ShExSeqImm0b_e1}.
By the naturality of the Spin-structure~$\fs_0(V,\fo)$ of Example~\ref{LBSpinStr_eg}
and the $\Pin^+$-structure~$\fp_0^+(V)$ of Lemma~\ref{Pin1_lmm},
$$\os_0\big(S\ga_{\R;1},\fo_{TS\ga_{\R;1}}\big)=
\big(\pi^*\os_0(\R\P^1,\fo_{T\R\P^1})\!\big)\!\big|_{S\ga_{\R;1}}
\quad\hbox{and}\quad
\fp_0^+\big(\pi^*\ga_{\R;1}|_{S\ga_{\R;1}}\big)
=\big(\pi^*\fp_0^+(\ga_{\R;1})\!\big)\!\big|_{S\ga_{\R;1}}\,.$$
Along with the definition and naturality of the map~\eref{SpinPinSESdfn_e0b},
this gives
\BE{ShExSeqImm0b_e5}
\fp_0^+\big|_{S\ga_{\R;1}}
=\bllrr{\os_0(S\ga_{\R;1},\fo_{S\ga_{\R;1}}),\fp_0^+\big(\pi^*\ga_{\R;1}
|_{S\ga_{\R;1}}\big)}_{\eref{ShExSeqImm0b_e3}}\,.\EE 
\end{eg}

\section{The Lie groups $\Spin(n)$ and $\Pin^{\pm}(n)$}
\label{LG_sec}

In this section, we describe the Lie groups $\Spin(n)$ and $\Pin^{\pm}(n)$ 
and key homomorphisms between them from purely topological considerations. 
These homomorphisms are used to establish the properties of Section~\ref{SpinPinProp_subs}
for the perspective of Definition~\ref{PinSpin_dfn} in Section~\ref{SpinPin_sec1}.\\

Let $m,n\!\in\!\Z^+$ with $m\!<\!n$. 
The standard identification of $\R^m\!\times\!\R^{n-m}$ with~$\R^n$ 
induces inclusions
\BE{iodfn_e}
\io_{n;m}\!:\SO(m)\!\times\!\SO(n\!-\!m)\lra \SO(n) \quad\hbox{and}\quad
\io_{n;m}\!:\SO(m)\!\times\!\O(n\!-\!m)\lra \O(n) \nota{iznm@$\io_{n;m},\io_{n;m}',\io_{n;m}''$}\EE
so that the diagram
$$\xymatrix{ \{\1\} \ar[r]& \SO(m)\!\times\!\SO(n\!-\!m) \ar[r]\ar[d]^{\io_{n;m}}& 
\SO(m)\!\times\!\O(n\!-\!m)\ar[r]\ar[d]^{\io_{n;m}}& \Z_2\ar[r]\ar[d]^{\id}&\{\1\}\\
\{\1\}\ar[r]& \SO(n) \ar[r]& \O(n)\ar[r]& \Z_2\ar[r]&\{\1\}}$$
commutes;
the homomorphisms to $\Z_2\!\equiv\!\Z/2\Z$ above are given by
the sign of the determinant.
We denote by
$$\io_{n;m}'\!:\SO(m)\lra \SO(n), \quad
\io_{n;n-m}''\!:\SO(n\!-\!m)\lra \SO(n), \quad
\io_{n;n-m}''\!:\O(n\!-\!m)\lra \O(n), \nota{iznm@$\io_{n;m},\io_{n;m}',\io_{n;m}''$}$$
the compositions of the maps in~\eref{iodfn_e}
with the canonical inclusions of~$\SO(m)$, 
$\SO(n\!-\!m)$, and $\O(n\!-\!m)$ in the domains of these maps.\\

We note that
\begin{gather}\notag
\SO(1)=\{1\}, \quad \O(1)=\{\pm1\}, \quad \SO(2)=S^1,
\quad \SO(3)\approx\R\P^3,\\
\label{SOsumm_e}
\big|\pi_0\big(\SO(n)\big)\big|,\big|\pi_2\big(\SO(n)\big)\big|=1~~\forall\,n\!\ge\!1,
\qquad \pi_1\big(\SO(n)\big)=\begin{cases}\Z,&\hbox{if}~n\!=\!2;\\
\Z_2,&\hbox{if}~n\!\ge\!3.
\end{cases}\end{gather}
The $\pi_2$ statement above and the $n\!\ge\!3$ case of the $\pi_1$ statement
are obtained by induction from the homotopy exact sequence
for the fibration
\BE{SnFB_e}\SO(n)\xlra{\io_{n+1;n}'} \SO(n\!+\!1)\lra S^n\,.\EE
The induced homomorphism
\BE{SnFB_e2}\io_{n+1;n\,*}'\!:\pi_1\big(\SO(n)\big)\lra \pi_1\big(\SO(n\!+\!1)\big)\EE
is an isomorphism for $n\!\ge\!3$ and is a surjection for $n\!=\!2$.

\subsection{The groups $\SO(n)$ and $\Spin(n)$}
\label{SO2Spin_subs}

\noindent
By the last statement in~\eref{SOsumm_e}, for every $n\!\ge\!2$ there exists 
a unique $2\!:\!1$ covering projection
\BE{Spindfn_e}q_n\!:\Spin(n)\lra \SO(n) \nota{qn@$q_n$|textbf}
\gena{Spin group $\Spin(n)$|textbf}\EE
from a connected Lie group. 
This {\it defines} the Lie group $\Spin(n)$ and determines an exact sequence
\BE{Spindfn_e2} \{\1\}\lra \Z_2 \lra \Spin(n) \stackrel{q_n}{\lra} \SO(n)\lra \{\1\}\EE
of Lie groups; see Lemma~\ref{CovExt_lmm}\ref{CovExtConn_it1}.
The preimage of the identity $\bI_n\!\in\!\SO(n)$\nota{in@$\bI_n,\wt\bI_n,\wh\bI_n$} 
under~\eref{Spindfn_e} consists
of the identity \hbox{$\wt\bI_n\!\in\!\Spin(n)$} and an element 
$\wh\bI_n\!\in\!\Spin(n)$\nota{in@$\bI_n,\wt\bI_n,\wh\bI_n$}
of order~2.
Since the inclusion of~$\Z_2$ in~\eref{Spindfn_e2} sends its non-trivial element~$1_{\Z_2}$
to $\wh\bI_n\!\in\!\Spin(n)$,
\BE{Spindfn_e4}\wh\bI_n\wt{A}=\wt{A}\wh\bI_n \qquad\forall\,\wt{A}\!\in\!\Spin(n).\EE

\begin{eg}\label{SU2_eg}
Let $\SU(n)$ denote the $n$-th \sf{special unitary group}.
We show~that 
\BE{Spin3_e}
q_3\!:\Spin(3)\!=\!\SU(2)\!\approx\!S^3\lra \SO(3)\!\approx\!\R\P^3.
\nota{q3@$q_3$}\gena{Spin$(3)$}\EE
The subspace $S^3\!\subset\!\C^2$ inherits a Lie group structure
from the quaternion multiplication on \hbox{$\bH\!=\!\C\!\oplus\!\fj\C$}.
The~map 
$$\SU(2) \lra S^3, \qquad
\left(\!\!\begin{array}{cc}a & -\ov{b} \\ b & ~\ov{a} \end{array}\!\!\right) \lra a+\fj b,$$
is a Lie group isomorphism; it intertwines the standard action of $\SU(2)$ on~$\C^2$
with the action of~$S^3$ on~$\bH$ by the multiplication on the left.
The~map
$$\bH\!\times\!\bH \lra \bH, \qquad 
(u,x) \lra  \Ad_u(x)\!\equiv\!ux\ov{u},$$
is linear in the second input~$x$ and restricts to an action of 
$S^3\!\subset\!\bH$ by isometries~on the subspace
$$\Im\,\bH\equiv\big\{x\!\in\!\bH\!: \Re\,x\!=\!0\big\}\approx\R^3$$
of~$\bH$.
The action of the differential of this restriction at $u\!=\!\1\!\in\!S^3$ is the homomorphism
$$\ad\!:\Im\,\bH\lra\End_{\R}\big(\Im\,\bH\big), \qquad
v\lra\ad_v,~~\ad_v(x)=2\,\Im(vx).$$
Since this homomorphism is injective, the Lie group homomorphism
$$\Ad\!: S^3\lra \SO(3)$$
induced by the action of $S^3$ on $\Im\,\bH$ is a covering projection.
\end{eg}

Let $m\!\in\!\Z^+$ with $m\!<\!n$.
By the last statement in~\eref{SOsumm_e} and 
the surjectivity of the homomorphism~\eref{SnFB_e2} for $n\!\ge\!2$,
the composition
$$\Spin(m)\!\times\!\Spin(n\!-\!m)\xlra{q_m\times q_{n-m}} \SO(m)\!\times\!\SO(n\!-\!m)
\xlra{\io_{n;m}} \SO(n)$$
induces the trivial homomorphism on the fundamental groups.
If $m,n\!-\!m\!\ge\!2$, Lemma~\ref{CovExt_lmm}\ref{CovExtConn_it2}
implies that the first embedding in~\eref{iodfn_e} lifts uniquely
to a Lie group 
homomorphism~$\wt\io_{n;m}$\nota{iznm2@$\wt\io_{n;m},\wt\io_{n;m}',\wt\io_{n;m}''$} 
so that the diagram 
\BE{wtioSpin_e}\begin{split}
\xymatrix{ \{\1\}\ar[r]& \Z_2\!\times\!\Z_2 \ar[r]\ar[d]^+& 
\Spin(m)\!\times\!\Spin(n\!-\!m) \ar[r]\ar[d]^{\wt\io_{n;m}}& 
\SO(m)\!\times\!\SO(n\!-\!m)  \ar[d]^{\io_{n;m}} \ar[r]& \{\1\}\\ 
\{\1\}\ar[r]& \Z_2 \ar[r]& \Spin(n) \ar[r]^{q_n}& \SO(n)\ar[r]& \{\1\}}
\end{split}\EE
commutes.
Since the homomorphisms on the fundamental groups induced by $\io_{n;m}'$ and 
$\io_{n;n-m}''$ surject onto $\pi_1(\SO(n))\!\approx\!\Z_2$ in this case,
\BE{wtioSpin_e2} 
\wt\io_{n;m}\big(\wt\bI_m,\wh\bI_{n-m}\big),
\wt\io_{n;m}\big(\wh\bI_m,\wt\bI_{n-m}\big)=\wh\bI_n\,.\EE
If $m\!=\!1$ or $n\!-\!m\!=\!1$, Lemma~\ref{CovExt_lmm}\ref{CovExtConn_it2} implies the existence
of a lift~$\wt\io_{n;m}$ of~$\io_{n;m}$ to the identity component 
of $\Spin(m)\!\times\!\Spin(n\!-\!m)$.
The condition~\eref{wtioSpin_e2} and the property~\eref{Spindfn_e4} 
extend it over the remaining component(s).
We denote by
$$\wt\io_{n;m}'\!:\Spin(m)\lra \Spin(n) \qquad\hbox{and}\qquad
\wt\io_{n;n-m}''\!:\Spin(n\!-\!m)\lra \Spin(n)
\nota{iznm2@$\wt\io_{n;m},\wt\io_{n;m}',\wt\io_{n;m}''$}$$
the compositions of~$\wt\io_{n;m}$ with the canonical inclusions of~$\Spin(m)$
and $\Spin(n\!-\!m)$ in the domain of~$\wt\io_{n;m}$.
By~\eref{wtioSpin_e2}, these maps are embeddings.\\

Let $\bI_{n;2}\!\in\!\SO(n)$\nota{in2@$\bI_{n;2},\wt\bI_{n;2},\wh\bI_{n;2}$} 
for $n\!\ge\!2$ and $\bI_{n;4}\!\in\!\SO(n)$\nota{in4@$\bI_{n;4},\wt\bI_{n;4},\wh\bI_{n;4}$} 
for $n\!\ge\!4$
be the diagonal matrices satisfying
\BE{In24dfn_e}\big(\bI_{n;2}\big)_{ii}=\begin{cases}1,&\hbox{if}~i\!\le\!n\!-\!2;\\
-1,&\hbox{otherwise};\end{cases}
\qquad
\big(\bI_{n;4}\big)_{ii}=\begin{cases}1,&\hbox{if}~i\!\le\!n\!-\!4;\\
-1,&\hbox{otherwise}.\end{cases}\EE
We denote the two elements of the preimage of~$\bI_{n;2}$ 
under~\eref{Spindfn_e} by $\wt\bI_{n;2}$ 
and~$\wh\bI_{n;2}$\nota{in2@$\bI_{n;2},\wt\bI_{n;2},\wh\bI_{n;2}$} 
and the two elements of the preimage of~$\bI_{n;4}$ by $\wt\bI_{n;4}$ 
and~$\wh\bI_{n;4}$.\nota{in4@$\bI_{n;4},\wt\bI_{n;4},\wh\bI_{n;4}$}
If $m\!\in\!\Z^+$ and $m\!<\!n$,
\BE{I2244_e} 
\wt\io_{n;m}''\big(\{\wt\bI_{m;2},\wh\bI_{m;2}\}\!\big)=\big\{\wt\bI_{n;2},\wh\bI_{n;2}\big\}
~\hbox{if}~m\!\ge\!2,
\quad
\wt\io_{n;m}''\big(\{\wt\bI_{m;4},\wh\bI_{m;4}\}\!\big)=\big\{\wt\bI_{n;4},\wh\bI_{n;4}\big\}
~\hbox{if}~m\!\ge\!4.\EE
Since the multiplication by $\wh\bI_n$ is the deck transformation of
the double cover~\eref{Spindfn_e},
\begin{equation*}
\begin{aligned}
\wh\bI_{n;2}&=\wh\bI_n\wt\bI_{n;2}=\wt\bI_{n;2}\wh\bI_n,
&\qquad
\wt\bI_{n;2}^2&=\!\wh\bI_{n;2}^2\in\big\{\wt\bI_n,\wh\bI_n\big\},\\
\wh\bI_{n;4}&=\wh\bI_n\wt\bI_{n;4}=\wt\bI_{n;4}\wh\bI_n,
&\qquad
\wt\bI_{n;4}^2&=\wh\bI_{n;4}^2\in\big\{\wt\bI_n,\wh\bI_n\big\}.
\end{aligned}\end{equation*}
By the second statement in~\eref{Spin1and2_e}, $\wt\bI_{2;2}^2\!=\!\wh\bI_2$.
Along with~\eref{I2244_e} and~\eref{wtioSpin_e2}, this implies~that
$$\wt\bI_{n;2}^2=\wt\io_{n;2}''\big(\wt\bI_{2;2}\big)^2
=\wt\io_{n;2}''\big(\wt\bI_{2;2}^2\big)=\wt\io_{n;2}''(\wh\bI_2)\!=\!\wh\bI_n
\quad\forall~n\!\ge\!2\,.$$
The next lemma gives a direct proof of this statement which readily adapts
to show that $\wt\bI_{n;4}^2\!=\!\wt\bI_n$ for all $n\!\ge\!4$.

\begin{lmm}\label{SO2Spin_lmm}
With the notation as above,
\BE{SO2Spin_e}\wt\bI_{n;2}^2=\wh\bI_{n;2}^2=\wh\bI_n \quad\forall~n\!\ge\!2
\qquad\hbox{and}\qquad  \wt\bI_{n;4}^2=\wh\bI_{n;4}^2=\wt\bI_n \quad\forall~n\!\ge\!4.\EE
\end{lmm}

\begin{proof} Let $n\!\ge\!2$, $n_-\!=\!n$, and $n_+\!=\!n\!+\!2$.
We define a path~$\ga_-$ in~$\SO(n)$ from~$\bI_n$ to~$\bI_{n;2}$ and 
a~path~$\ga_+$ in~$\SO(n\!+\!2)$ from~$\bI_{n+2}$ to~$\bI_{n+2;4}$ by
\begin{alignat*}{2}
&\ga_-\!:[0,1]\lra\SO(n), &\qquad &\ga_+\!:[0,1]\lra\SO(n\!+\!2),\\
&\ga_-(t)=\left(\!\!\begin{array}{ccc} \bI_{n-2} && \\
& \cos(\pi t) &   -\sin(\pi t)\\
& \sin(\pi t) & \cos(\pi t) \end{array}\!\!\!\right), &\qquad &
\ga_+(t)=\left(\!\!\begin{array}{ccc} \ga_-(t) && \\
& \cos(\pi t) &   -\sin(\pi t)\\
& \sin(\pi t) & \cos(\pi t) \end{array}\!\!\!\right).
\end{alignat*}
The endpoints of the lifts 
$$\wt\ga_-\!:[0,1]\lra\Spin(n) \qquad\hbox{and}\qquad \wt\ga_+\!:[0,1]\lra\Spin(n\!+\!2)$$
of these paths such that $\wt\ga_{\pm}(0)\!=\!\wt\bI_{n_{\pm}}$
lie in $\{\wt\bI_{n;2},\wh\bI_{n;2}\}$ and 
$\{\wt\bI_{n+2;4},\wh\bI_{n+2;4}\}$, respectively.
Let 
$$\al_{\pm}\!\equiv\!
\ga_{\pm}\!*\!\big(\ga_{\pm}(1)\ga_{\pm}\big)\!:[0,1]\lra\SO(n_{\pm})  \quad\hbox{and}\quad
\wt\al_{\pm}\!\equiv\!\wt\ga_{\pm}\!*\!\big(\wt\ga_{\pm}(1)\wt\ga_{\pm}\big)\!:
[0,1]\lra\Spin(n_{\pm})$$
be the concatenations (products) of the paths $\ga_{\pm}$ and $\wt\ga_{\pm}$
with the paths
\begin{alignat*}{2}
\ga_{\pm}(1)\ga_{\pm}\!:[0,1]&\lra\SO(n_{\pm}), &\quad
\big\{\ga_{\pm}(1)\ga_{\pm}\big\}(t)&=\ga_{\pm}(1)\ga_{\pm}(t), \qquad\hbox{and}\\
\wt\ga_{\pm}(1)\wt\ga_{\pm}\!:[0,1]&\lra\Spin(n_{\pm}), &\quad
\big\{\wt\ga_{\pm}(1)\wt\ga_{\pm}\big\}(t)&=\wt\ga_{\pm}(1)\wt\ga_{\pm}(t),
\end{alignat*}
respectively.
In particular,
\BE{SO2Spin_e5}
\wt\bI_{n;2}^2=\!\wh\bI_{n;2}^2=\wt\al_-(1),\qquad
\wt\bI_{n+2;4}^2=\!\wh\bI_{n+2;4}^2=\wt\al_+(1)\,.\EE

\vspace{.2in}

By the homotopy exact sequences for the fibrations~\eref{SnFB_e},
the~path
$$\al_-\!:[0,1]\lra\SO(n), \quad
\al_-(t)
=\left(\!\!\begin{array}{ccc} \bI_{n-2} && \\
& \cos(2\pi t) &   -\sin(2\pi t)\\
& \sin(2\pi t) & \cos(2\pi t) \end{array}\!\!\!\right)\,,$$
is a loop generating $\pi_1(\SO(n))$.
Thus, its lift~$\wt\al_-$ is not a loop and so $\wt\al_-(1)\!\neq\!\wt\bI_n$.
Along with the first equation in~\eref{SO2Spin_e5}, 
this establishes the first claim in~\eref{SO2Spin_e}.\\

By the homotopy exact sequences for the fibrations~\eref{SnFB_e}, the~path
$$\al_+\!:[0,1]\lra\SO(n\!+\!2), \quad
\al_+(t)
= \left(\!\!\begin{array}{ccc} \al_-(t) && \\
& \cos(2\pi t) &   -\sin(2\pi t)\\
& \sin(2\pi t) & \cos(2\pi t) \end{array}\!\!\!\right),$$
is a loop representing twice the generator of $\pi_1(\SO(n\!+\!2))\!\approx\!\Z_2$.
Thus, its lift~$\wt\al_+$ is a loop and so $\wt\al_+(1)\!=\!\wt\bI_{n+2}$.
Along with the second equation in~\eref{SO2Spin_e5}, 
this establishes the second claim in~\eref{SO2Spin_e}.
\end{proof}

For concreteness, we take 
\BE{I24dfn_e}
\wt\bI_{n;2}\!\equiv\!\wt\ga_-(1)\in\Spin(n),
\nota{in2@$\bI_{n;2},\wt\bI_{n;2},\wh\bI_{n;2}$}
\qquad \wt\bI_{n+2;4}\!\equiv\!\wt\ga_+(1)\in\Spin(n\!+\!2)
\nota{in4@$\bI_{n;4},\wt\bI_{n;4},\wh\bI_{n;4}$},\EE
with the notation as in the proof of Lemma~\ref{SO2Spin_lmm}.
With these choices,
\BE{I2244_e2}
\wt\io_{n;m}''\big(\wt\bI_{m;2}\big)=\wt\bI_{n;2}~~\hbox{if}~2\!\le\!m\!\le\!n
\qquad\hbox{and}\qquad
\wt\io_{n;m}''(\wt\bI_{m;4})=\wt\bI_{n;4}~~\hbox{if}~4\!\le\!m\!\le\!n.\EE

\vspace{.15in}

Let $m_1,m_2\!\in\!\Z^+$ with $m_1,m_2\!\le\!n$ and $m_1\!\neq\!m_2$.
We denote by 
$\bI_{n;(m_1,m_2)}\!\in\!\SO(n)$\nota{inm1m2@$\bI_{n;(m_1,m_2)},\wt\bI_{n;(m_1,m_2)},\wh\bI_{n;(m_1,m_2)}$} 
the diagonal matrix satisfying
\BE{Inmdfn_e1}\big(\bI_{n;(m_1,m_2)}\big)_{ii}=\begin{cases}-1,&\hbox{if}~i\!=\!m_1,m_2;\\
1,&\hbox{otherwise};\end{cases}\EE
and by $\wt\bI_{n;(m_1,m_2)},\wh\bI_{n;(m_1,m_2)}\!\in\!\Spin(n)$
\nota{inm1m2@$\bI_{n;(m_1,m_2)},\wt\bI_{n;(m_1,m_2)},\wh\bI_{n;(m_1,m_2)}$} 
its two preimages
under~\eref{Pindfn_e}.
By symmetry and Lemma~\ref{SO2Spin_lmm} (or directly by its proof),
\BE{SO2Spin_e2} \wt\bI_{n;(m_1,m_2)}^2,\wh\bI_{n;(m_1,m_2)}^2=\wh\bI_n\,.\EE
If in addition $m\!\in\!\Z^+$ with $m_1,m_2\!\le\!m\!<\!n$, then
\BE{I22m_e} 
\wt\io_{n;m}'\big(\{\wt\bI_{m;(m_1,m_2)},\wh\bI_{m;(m_1,m_2)}\}\!\big)
=\big\{\wt\bI_{n;(m_1,m_2)},\wh\bI_{n;(m_1,m_2)}\big\}.\EE

\vspace{.15in}

For $B\!\in\!\O(n)$, let 
$$\fc(B)\!:\SO(n)\lra\SO(n), \qquad \fc(B)A=BAB^{-1}, \nota{cb@$\fc(B),\wt\fc(B)$}$$
be the conjugation action of $B$ on $\SO(n)$.
By Lemma~\ref{CovExt_lmm}\ref{CovExtConn_it2}, it lifts uniquely to an isomorphism
$$ \wt\fc(B)\!:\Spin(n)\lra\Spin(n). $$
In particular,
\BE{laminbIn_e}\wt\fc(-\bI_n)\!=\!\id\!: \Spin(n)\lra\Spin(n)\,.\EE
For $m\!\in\!\Z^+$ with $m\!\le\!n$, 
we denote by $\bI_{n;(m)}\!\in\!\O(n)$\nota{inm@$\bI_{n;(m)},\wt\bI_{n;(m)},\wh\bI_{n;(m)}$} 
the diagonal matrix satisfying
\BE{Inmdfn_e}\big(\bI_{n;(m)}\big)_{ii}=\begin{cases}-1,&\hbox{if}~i\!=\!m;\\
1,&\hbox{otherwise}.\end{cases}\EE
Let $\bI_{n;1}\!=\!\bI_{n;(n)}$.\nota{in1@$\bI_{n;1},\wt\bI_{n;1},\wh\bI_{n;1}$}

\begin{crl}\label{SO2Spin_crl}
If $m,m_1,m_2\!\in\!\Z^+$ with $m,m_1,m_2\!\le\!n$ and $m_1\!\neq\!m_2$, then
\BE{SO2Spin2_e}\wt\fc(\bI_{n;(m)})\wt\bI_{n;(m_1,m_2)}=\wt\bI_{n;(m_1,m_2)}\cdot
\begin{cases}\wt\bI_n,&\hbox{if}~m\!\not\in\!\{m_1,m_2\};\\
\wh\bI_n,&\hbox{if}~m\!\in\!\{m_1,m_2\}.
\end{cases}\EE
\end{crl}

\begin{proof} By symmetry, we can assume that $m_1\!=\!m_2\!-\!1$ and $m_2\!\le\!m\!=\!n$.
If $m_2\!<\!n$, the group homomorphisms
$$\Spin(m_2)\lra \Spin(n), \qquad 
\wt{A}\lra\wt\io_{n;m_2}'\big(\wt{A}\big),
\wt\fc(\bI_{n;1})\wt\io_{n;m_2}'\big(\wt{A}\big),$$
lift $\io_{n;m_2}'$.
Lemma~\ref{CovExt_lmm}\ref{CovExtConn_it2} thus implies that
\BE{SO2Spin2_e3} \wt\fc(\bI_{n;1})\wt\io_{n;m_2}'\big(\wt{A}\big)=\wt\io_{n;m_2}'\big(\wt{A}\big)
\qquad\forall~\wt{A}\!\in\!\Spin(m_2)\,.\EE
Along with~\eref{I22m_e}, this implies that 
$$\wt\fc(\bI_{n;1})\wt\bI_{n;(m_1,m_2)}
=\wt\fc(\bI_{n;1})\wt\io_{n;m_2}'\big(\wt\bI_{m_2;(m_1,m_2)}\big)
=\wt\io_{n;m_2}'\big(\wt\bI_{m_2;(m_1,m_2)}\big)
=\wt\bI_{m_2;(m_1,m_2)}\,.$$
This establishes the first case in~\eref{SO2Spin2_e}.\\

\noindent
Let $\ga_-\!:[0,1]\!\lra\!\SO(n)$ be as in the proof of Lemma~\ref{SO2Spin_lmm}.
The~path
$$\fc(\bI_{n;1})\ga_-\!:[0,1]\lra\SO(n), \qquad
\fc(\bI_{n;1})\ga_-(t)=\left(\!\!\begin{array}{ccc} \bI_{n-2} && \\
& \cos(\pi t) &   \sin(\pi t)\\
& -\sin(\pi t) & \cos(\pi t) \end{array}\!\!\!\right),$$
is the reverse of the path $\ga_-(1)\ga_-$.
Since the lift~$\wt\al_-$ of the loop $\al_-\!\equiv\!\ga_-\!*\!(\ga_-(1)\ga_-)$
to $\Spin(n)$ is not a loop, the endpoints of the lifts 
$$\wt\ga_-,\wt\fc(\bI_{n;1})\wt\ga_-\!:\big([0,1],0\big)\lra\big(\Spin(n),\wt\bI_n\big)$$
of $\ga_-$ and $\fc(\bI_{n;1})\ga_-$ are different points of the preimage
$\{\wt\bI_{n;2},\wh\bI_{n;2}\}$ of~$\bI_{n;2}$ under~\eref{Spindfn_e}.
This establishes the second case in~\eref{SO2Spin2_e}.
\end{proof}

Let $m\!\in\!\Z^+$ with $m\!<\!n$.
Since $\Spin(m)$ is connected if $m\!\ge\!2$ and the group homomorphisms
$$\Spin(m)\lra \Spin(n), \qquad 
\wt{A}\lra \wt\io_{n;m}'\big(\wt\fc(\bI_{m;1})\wt{A}\big),
\wt\fc(\bI_{n;(m)})\wt\io_{n;m}'\big(\wt{A}\big),$$
lift $\fc(\bI_{n;(m)})\io_{n;m}'$, they agree; Lemma~\ref{CovExt_lmm}\ref{CovExtConn_it2}. 
Thus, 
\BE{SO2Spin2_e7} 
\wt\io_{n;m}'\big(\wt\fc(\bI_{m;1})\wt{A}\big)=
\wt\fc(\bI_{n;(m)})\wt\io_{n;m}'\big(\wt{A}\big)
\qquad\forall~\wt{A}\!\in\!\Spin(m)\,.\EE
Since $\Spin(m)$ is connected for $m\!\ge\!2$ and the group homomorphisms
$$\Spin(m)\lra \Spin(n), \qquad 
\wt{A}\lra\wt\io_{n;m}''\big(\wt\fc(\bI_{m;1})\wt{A}\big),
\wt\fc(\bI_{n;1})\wt\io_{n;m}''\big(\wt{A}\big),$$
lift $\fc(\bI_{n;1})\io_{n;m}''$, the two homomorphisms agree. 
Since the group homomorphisms
$$\Spin(m)\lra \Spin(n), \qquad 
\wt{A}\lra\wt\io_{n;m}''\big(\wt{A}\big),
\wt\fc(\bI_{n;(m)})\wt\io_{n;m}''\big(\wt{A}\big),$$
lift $\io_{n;m}''$, they also agree.
Thus,
\BE{SO2Spin2_e9}
\wt\io_{n;m}''\big(\wt\fc(\bI_{m;1})\wt{A}\big)=
\wt\fc(\bI_{n;1})\wt\io_{n;m}''\big(\wt{A}\big), ~~
\wt\io_{n;m}''\big(\wt{A}\big)=\wt\fc(\bI_{n;(m)})\wt\io_{n;m}''\big(\wt{A}\big)
\quad\forall~\wt{A}\!\in\!\Spin(m)\,.\EE
The $m\!=\!1$ cases of~\eref{SO2Spin2_e7} and~\eref{SO2Spin2_e9} follow
from~\eref{wtioSpin_e2}.

\subsection{The groups $\O(n)$ and $\Pin^{\pm}(n)$}
\label{O2Pin_subs}

The group $\O(n)$ is isomorphic to the semi-direct product $\SO(n)\!\rtimes\!\Z_2$ with
the action of the non-trivial element~$1_{\Z_2}$ of~$\Z_2$ on~$\SO(n)$ given 
by the conjugation by~$\bI_{n;(m)}$, for any fixed choice of~$m$, or 
by any other order~2 element of the non-identity component of~$\O(n)$.

\begin{eg}\label{O2_eg}
As a pointed smooth manifold, $(\O(2),\bI_2)$ is diffeomorphic to 
$(S^1\!\times\!\Z_2,1\!\times\!0)$.
However, the group structure of~$\O(2)$ is given~by
$$\big(\R/2\pi\Z\!\times\!\Z_2\big)\!\times\!\big(\R/2\pi\Z\!\times\!\Z_2\big)
\lra\big(\R/2\pi\Z,\Z_2\big), ~~
(\th_1,k_1)\!\cdot\!(\th_2,k_2)=\big(\th_1\!+\!(-1)^{k_1}\th_2,k_1\!+\!k_2\big).$$
\end{eg}

\vspace{.2in}

By the last statement in~\eref{SOsumm_e}, for every $n\!\ge\!2$ there exists 
a {\it topologically} unique double cover
\BE{Pindfn_e}q_n^{\pm}\!:\Pin^{\pm}(n)\lra \O(n)  \nota{qnpm@$q_n^{\pm}$|textbf}\EE
with two path components
that restricts to~\eref{Spindfn_e} over $\SO(n)\!\subset\!\O(n)$.
However, there are two Lie group structures on $\Pin(n)$,
denoted by $\Pin^{\pm}(n)$, so that the diagram 
$$\xymatrix{ \{\1\}\ar[r]& \Z_2 \ar[r]\ar[d]^{\id}& 
\Spin(n) \ar[r]^{q_n}\ar[d]& \SO(n)  \ar[d] \ar[r]& \{\1\}\\ 
\{\1\}\ar[r]& \Z_2 \ar[r]& \Pin^{\pm}(n) \ar[r]^{q_n^{\pm}}& \O(n)\ar[r]& \{\1\}}$$
of Lie group homomorphisms commutes;
the inclusions of~$\Z_2$ above send its non-trivial element~$1_{\Z_2}$
to $\wh\bI_n\!\in\!\Spin(n)\!\subset\!\Pin^{\pm}(n)$.
The two group structures are described below.\\

We denote by \hbox{$\wt\bI_{n;n},\wh\bI_{n;n}\!\in\!\Pin^{\pm}(n)$}
\nota{inn@$\bI_{n;n},\wt\bI_{n;n},\wh\bI_{n;n}$} the two preimages
of \hbox{$\bI_{n;n}\!\equiv\!-\bI_n$}\nota{inn@$\bI_{n;n},\wt\bI_{n;n},\wh\bI_{n;n}$}
under~\eref{Pindfn_e}.
For $m\!\in\!\Z^+$ with $m\!\le\!n$, let
\hbox{$\wt\bI_{n;(m)},\wh\bI_{n;(m)}\!\in\!\Pin^{\pm}(n)$}
\nota{inm@$\bI_{n;(m)},\wt\bI_{n;(m)},\wh\bI_{n;(m)}$} 
be the two preimages
of~$\bI_{n;(m)}$. 
In particular,
$$\wt\bI_{n;n}^2,\wh\bI_{n;n}^2,
\wt\bI_{n;(m)}^2,\wh\bI_{n;(m)}^2\in
q_n^{-1}(\bI_n)\!\equiv\!\big\{\wt\bI_n,\wh\bI_n\big\}.$$
If $m_1,m_2\!\in\!\Z^+$ with $m_1,m_2\!\le\!n$ and $m_1\!\neq\!m_2$, then
\BE{m1m2prod_e} 
\wt\bI_{n;(m_1)}\!\cdot\!\big\{\wt\bI_{n;(m_1,m_2)},\wh\bI_{n;(m_1,m_2)}\big\},
\big\{\wt\bI_{n;(m_1,m_2)},\wh\bI_{n;(m_1,m_2)}\big\}\!\cdot\!\wt\bI_{n;(m_1)}\in
\big\{\wt\bI_{n;(m_2)},\wh\bI_{n;(m_2)}\big\}\,.\EE
Let $\wt\bI_{n;1}\!=\!\wt\bI_{n;(n)}$
and $\wh\bI_{n;1}\!=\!\wh\bI_{n;(n)}$\nota{in1@$\bI_{n;1},\wt\bI_{n;1},\wh\bI_{n;1}$}.\\

\noindent
Since the multiplication by $\wh\bI_n$ is the deck transformation of
the double cover~\eref{Pindfn_e},
\BE{Pindfn_e2} \wh\bI_n\wt{A}=\wt{A}\wh\bI_n\quad\forall\,\wt{A}\!\in\!\Pin^{\pm}(n).\EE
Along with~\eref{m1m2prod_e}, \eref{SO2Spin2_e}, and~\eref{SO2Spin_e2},
this implies that 
\BE{wtInmsq_e}\begin{split}
\wt\bI_{n;(m)}^2,\wh\bI_{n;(m)}^2=\big(\wt\bI_{n;1}\wt\bI_{n;(m,n)}\big)^2
&=\big(\wt\fc(\bI_{n;1})\wt\bI_{n;(m,n)}\big)\wt\bI_{n;1}^2\wt\bI_{n;(m,n)}\\
&=\wt\bI_{n;(m,n)}\wh\bI_n\wt\bI_{n;1}^2\wt\bI_{n;(m,n)}
=\wt\bI_{n;1}^2,\wh\bI_{n;1}^2
\end{split}\EE
for all $m\!\in\!\Z^+$ with $m\!<\!n$.
Along with another application of~\eref{SO2Spin_e2}, this in turn implies~that 
\BE{wtInmsq_e2}
\wt\bI_{n;1}\wt\bI_{n;(m)}\wt\bI_{n;1}^{-1}\wt\bI_{n;(m)}^{-1}=
\big(\wt\bI_{n;1}\wt\bI_{n;(m)}\big)^{\!2\,}\wt\bI_{n;1}^{-2}\wt\bI_{n;(m)}^{-2}
=\wh\bI_n, \quad
\wt\bI_{n;(m)}\wt\bI_{n;n}=\wt\bI_{n;n}\wt\bI_{n;(m)}\wh\bI_n^{\,n-1}\EE
for all $m\!\in\!\Z^+$ with $m\!<\!n$.\\

\noindent
Every element of $\Pin^{\pm}(n)\!-\!\Spin(n)$ can be written uniquely as $\wt{A}\wt\bI_{n;1}$ 
with $\wt{A}\!\in\!\Spin(n)$ and
$$\wt\bI_{n;1}\wt{A}=\big(\wt\bI_{n;1}\wt{A}\wt\bI_{n;1}^{\,-1}\big)\wt\bI_{n;1}
=\big(\wt\fc(\bI_{n;1})\wt{A}\big)\wt\bI_{n;1}
\qquad\forall\,\wt{A}\!\in\!\Spin(n).$$
The element $\wt\fc(\bI_{n;1})\wt{A}\!\in\!\Spin(n)$ is determined
by the adjoint action~$\fc(\bI_{n;1})$ of~$\bI_{n;1}$ on~$\SO(n)$.
The group multiplication in $\Pin^{\pm}(n)$ is thus characterized by whether
$\wt\bI_{n;1}^2$ equals $\wt\bI_n$ or~$\wh\bI_n$ or equivalently whether 
the subgroup of $\Pin^{\pm}(n)$ generated by~$\wt\bI_{n;1}$ and~$\wh\bI_{n;1}$
is $\Z_2^{\,2}$ or~$\Z_4$.
The Lie group $\Pin^+(n)$\gena{Pin group $\Pin^{\pm}(n)$|textbf} 
is {\it defined} to be the version of~$\Pin^{\pm}(n)$
in the first case
and $\Pin^-(n)$\gena{Pin group $\Pin^{\pm}(n)$|textbf} 
is {\it defined} to be the version of~$\Pin^{\pm}(n)$  in the second case.
Since the subspace of order~2 elements of~$\O(n)$ with precisely one $(-1)$-eigenvalue
is connected, the distinction between $\Pin^+(n)$ and $\Pin^-(n)$ 
can be equivalently formulated in terms of any such element.
Example~\ref{PinPin_eg} provides a more algebraic reason for this.
The group $\Pin^+(n)$ is isomorphic to the semi-direct product $\Spin(n)\!\rtimes\!\Z_2$ with
the action of the non-trivial element~$1_{\Z_2}$ of~$\Z_2$ on~$\Spin(n)$ 
given by the conjugation by~$\wt\bI_{n;(m)}$ for any $m\!\in\!\Z^+$ with $m\!\le\!n$.

\begin{eg}\label{Pin1Pin2_eg}
With the identification of $\O(2)$ given by Example~\ref{O2_eg},
the $n\!=\!2$ case of the cover~\eref{Pindfn_e} can be written~as
$$q_2^{\pm}\!:\Pin^{\pm}(2)\!=\!\R/2\pi\Z\!\times\!\Z_2\lra
\O(2)\!=\!\R/2\pi\Z\!\times\!\Z_2, \quad q_2^{\pm}(\th,k)=(2\th,k),
\nota{q2pm@$q_2^{\pm}$}\gena{Pin$^{\pm}(2)$}$$
in the smooth manifold category (ignoring the group structures).
The group structures are given~by
\begin{alignat*}{2}
\Pin^+(2)\!\times\!\Pin^+(2)&\lra\Pin^+(2), &\quad 
(\th_1,k_1)\!\cdot\!(\th_2,k_2)&=\big(\th_1\!+\!(-1)^{k_1}\th_2,k_1\!+\!k_2\big),\\
\Pin^-(2)\!\times\!\Pin^-(2)&\lra\Pin^-(2), &\quad 
(\th_1,k_1)\!\cdot\!(\th_2,k_2)&=
\big(\th_1\!+\!(-1)^{k_1}\th_2\!+\!k_1k_2\pi,k_1\!+\!k_2\big).
\end{alignat*}
In particular, $\Pin^+(2)\!\approx\!\O(2)$.
Since the square of every element in the non-identity component of $\Pin^+(2)$ 
(resp.~$\Pin^-(2)$) is~$\wt\bI_2$ (resp.~$\wh\bI_2$),
$\Pin^+(2)$ and $\Pin^-(2)$ are not isomorphic as groups
(even ignoring the projections to~$\O(2)$).
\end{eg}

Let $m\!\in\!\Z^+$ with $m\!<\!n$.
By~\eref{SO2Spin2_e3} and the first statement in~\eref{SO2Spin2_e9},
\BE{wtiocomm_e}\begin{aligned}
\wt\bI_{n;1}\wt\io_{n;m}'\big(\wt{A}\big)
&=\wt\io_{n;m}'\big(\wt{A}\big)\wt\bI_{n;1}
&\quad&\forall\,\wt{A}\!\in\!\Spin(m),\\
\wt\io_{n;n-m}''\big(\wt\bI_{n-m;1}\wt{A}\wt\bI_{n-m;1}^{-1}\big)
&=\wt\bI_{n;1}\wt\io_{n;n-m}''\big(\wt{A}\big)\wt\bI_{n;1}^{-1}
&\quad&\forall\,\wt{A}\!\in\!\Spin(n\!-\!m)\,.
\end{aligned}\EE
Thus, the homomorphism~$\wt\io_{n;m}$ in~\eref{wtioSpin_e} 
can be extended to a Lie group 
homomorphism~$\wt\io_{n;m}$\nota{iznm2@$\wt\io_{n;m},\wt\io_{n;m}',\wt\io_{n;m}''$} 
lifting 
the second embedding in~\eref{iodfn_e} so that the diagram
\BE{wtioPin_e}\begin{split}
\xymatrix{ \{\1\}\ar[r]& \Z_2\!\times\!\Z_2 \ar[r]\ar[d]^+& 
\Spin(m)\!\times\!\Pin^{\pm}(n\!-\!m) \ar[r]\ar[d]^{\wt\io_{n;m}}& 
\SO(m)\!\times\!\O(n\!-\!m)  \ar[d]^{\io_{n;m}} \ar[r]& \{\1\}\\ 
\{\1\}\ar[r]& \Z_2 \ar[r]& \Pin^{\pm}(n) \ar[r]^{q_n}& \O(n)\ar[r]& \{\1\}}
\end{split}\EE
commutes by fixing 
\BE{Pinliftcond_e}
\wt\io_{n;m}\big(\wt\bI_m,\wt\bI_{n-m;1}\big)\in\big\{\wt\bI_{n;1},\wh\bI_{n;1}\big\}.\EE
By~\eref{SO2Spin2_e7}, the second statement in~\eref{SO2Spin2_e9}, and
the first statement in~\eref{wtInmsq_e2},
\BE{wtiocomm_e4}\begin{split}
&\wt\io_{n;m}
\big(\wt\bI_{m;1}\wt{A}'\wt\bI_{m;1}^{-1},\wt{A}''\big)=
\wt\bI_{n;(m)}\wt\io_{n;m}'\big(\wt{A}'\big)\wt\bI_{n;(m)}^{-1}
\wt\io_{n;n-m}''\big(\wt{A}''\big)\\
&=\wt\bI_{n;(m)}\wt\io_{n;m}\big(\wt{A}',\wt{A}''\big)\wt\bI_{n;(m)}^{-1}
\!\times\!
\begin{cases}
\wt\bI_n,& \hbox{if}~(\wt{A}',\wt{A}'')\!\in\!\Spin(m)\!\times\!\Spin(n\!-\!m);\\
\wh\bI_n,&
\hbox{if}~(\wt{A}',\wt{A}'')\!\in\!\Spin(m)\!\times\!(\Pin^{\pm}(n\!-\!m)\!-\!\Spin(n\!-\!m)).
\end{cases}
\end{split}\EE

\vspace{.15in}

There are two possible lifts $\wt\io_{n;m}$ in~\eref{wtioPin_e}; 
they differ by the composition with the automorphism
\BE{rhodfn_e}
\rho\!:\Pin^{\pm}(n)\lra \Pin^{\pm}(n)
\qquad\hbox{s.t.}\quad 
\rho(\wt{A})=\begin{cases}
\wt{A},& \hbox{if}~\wt{A}\!\in\!\Spin(n);\\
\wt{A}\wh\bI_n,&\hbox{if}~\wt{A}\!\in\!\Pin^{\pm}(n)\!-\!\Spin(n).
\end{cases}\EE
We fix one lift by requiring~that 
\BE{wtioPindfn_e} \wt\io_{n;m}\big(\wt\bI_m,\wt\bI_{n-m;1}\big)=\wt\bI_{n;1}\,.\EE
We denote by
$$\wt\io_{n;n-m}''\!:\Pin^{\pm}(n\!-\!m)\lra \Pin^{\pm}(n) 
\nota{iznm2@$\wt\io_{n;m},\wt\io_{n;m}',\wt\io_{n;m}''$}$$
the composition of~$\wt\io_{n;m}$ with the canonical inclusion of~$\Pin^{\pm}(n\!-\!m)$ 
in the domain of~$\wt\io_{n;m}$.\\

Suppose $m_1,m_2\!\in\!\Z^+$ with $m_1\!+\!m_2\!<\!n$.
By~\eref{wtioPindfn_e}, 
\BE{wtioassoci_e0b}
\wt\io_{n;m_1}\big(\wt\bI_{m_1},\wt\io_{n-m_1;m_2}(\wt\bI_{m_2},\wt\bI_{n-m_1-m_2;1})\!\big)
=\wt\io_{n;m_1+m_2}\big(\wt\io_{m_1+m_2;m_1}(\wt\bI_{m_1},\wt\bI_{m_2}),
\wt\bI_{n-m_1-m_2;1}\!\big).\EE
By Lemma~\ref{CovExt_lmm}\ref{CovExtConn_it2}, 
\BE{wtioassoci_e0a}\begin{split}
\wt\io_{n;m_1}\!\circ\!\big\{\id_{\Spin(m_1)}\!\times\!\wt\io_{n-m_1;m_2}\big\}
&=\wt\io_{n;m_1+m_2}\!\circ\!\big\{\wt\io_{m_1+m_2;m_1}\!\times\!\id_{\Spin(n-m_1-m_2)}\big\}
\!:\\
&\quad \Spin(m_1)\!\times\!\Spin(m_2)\!\times\!\Spin(n\!-\!m_1\!-\!m_2)\lra\Spin(n),
\end{split}\EE
since these two homomorphisms lift the homomorphism
\begin{equation*}\begin{split}
\io_{n;m_1}\!\circ\!\big\{\id_{\SO(m_1)}\!\times\!\io_{n-m_1;m_2}\big\}
&=\io_{n;m_1+m_2}\!\circ\!\big\{\io_{m_1+m_2;m_1}\!\times\!\id_{\SO(n-m_1-m_2)}\big\}
\!:\\
&\quad \SO(m_1)\!\times\!\SO(m_2)\!\times\!\SO(n\!-\!m_1\!-\!m_2)\lra\SO(n).
\end{split}\end{equation*}
By~\eref{wtioassoci_e0a} and~\eref{wtioassoci_e0b},
\BE{wtioassoci_e}\begin{split}
\wt\io_{n;m_1}\!\circ\!\big\{\id_{\Spin(m_1)}\!\times\!\wt\io_{n-m_1;m_2}\big\}
&=\wt\io_{n;m_1+m_2}\!\circ\!\big\{\wt\io_{m_1+m_2;m_1}\!\times\!\id_{\Pin^{\pm}(n-m_1-m_2)}\big\}
\!:\\
&\quad \Spin(m_1)\!\times\!\Spin(m_2)\!\times\!\Pin^{\pm}(n\!-\!m_1\!-\!m_2)\lra\Pin^{\pm}(n).
\end{split}\EE

\begin{rmk}\label{PinPin_rmk}
The natural embedding of $\O(m)\!\times\!\O(n\!-\!m)$ into~$\O(n)$
does not lift to a Lie group homomorphism
$$\wt\io\!: \Pin^{\pm}(m)\!\times\!\Pin^{\pm}(n\!-\!m)\lra \Pin^{\pm}(n)\,.$$
If such lift existed, it would satisfy 
$$\wt\io\big(\wt\bI_m,\wt\bI_{n-m;1}\big)\in\big\{\wt\bI_{n;1},\wh\bI_{n;1}\big\},\qquad
\wt\io\big(\wt\bI_{m;1},\wt\bI_{n-m}\big)\in\big\{\wt\bI_{n;(m)},\wh\bI_{n;(m)}\big\}.$$
Since $(\wt\bI_m,\wt\bI_{n-m;1})$ and $(\wt\bI_{m;1},\wt\bI_{n-m})$ commute 
in the domain of~$\wt\io$,
this would contradict the first statement in~\eref{wtInmsq_e2}.
\end{rmk}

The embeddings 
\begin{alignat*}{2}
&\io_n^-\!:\O(n)\lra \SO(n\!+\!1),  &\qquad &\io_n^+\!:\O(n)\lra \SO(n\!+\!3),\\
&\io_n^-(A)=\left(\!\!\begin{array}{cc} A & \\ & \det A \end{array}\!\!\!\right), 
&\qquad  &\io_n^+(A)=\left(\!\!\begin{array}{cccc} A & & & \\
 & \det A & & \\ & & \det A& \\ & & & \det A \end{array}\!\!\!\right), 
\nota{izn@$\io_n^{\pm},\wt\io_n^{\pm}$}
\end{alignat*}
also lift to Lie group homomorphisms.
 
\begin{crl}\label{Pin2Spin_crl}
The embeddings $\io_n^-$ and~$\io_n^+$ lift to injective Lie group 
homomorphisms~$\wt\io_n^{\,-}$ and~$\wt\io_n^{\,+}$\nota{izn@$\io_n^{\pm},\wt\io_n^{\pm}$} 
so that the diagrams 
$$\xymatrix{\Spin(n) \ar@/^1.8pc/[rr]_{\wt\io_{n+1;n}'}
\ar[r]\ar[d]&  \Pin^-(n) \ar[r]^>>>>>{\wt\io_n^{\,-}}\ar[d]& \Spin(n\!+\!1)\ar[d] &
\Spin(n)\ar@/^1.8pc/[rr]_{\wt\io_{n+3;n}'}
\ar[r]\ar[d]&\Pin^+(n) \ar[r]^>>>>>{\wt\io_n^{\,+}}\ar[d]& \Spin(n\!+\!3)\ar[d]\\
\SO(n) \ar@/_1.8pc/[rr]^{\io_{n+1;n}'} \ar[r]&\O(n) \ar[r]^>>>>>>>{\io_n^-}& \SO(n\!+\!1) & 
\SO(n) \ar@/_1.8pc/[rr]^{\io_{n+3;n}'} \ar[r]&\O(n) \ar[r]^>>>>>>>{\io_n^+}& \SO(n\!+\!3)}$$
commute.
These homomorphisms satisfy
\BE{Pin2Spin_e}\begin{aligned}
\wt\io_n^{\,-}\big(\wh\bI_n\big)&=\wh\bI_{n+1},&\quad
\wt\io_n^{\,-}\big(\big\{\wt\bI_{n;1},\wh\bI_{n;1}\big\}\big)
&=\big\{\wt\bI_{n+1;2},\wh\bI_{n+1;2}\big\},\\
\wt\io_n^{\,+}\big(\wh\bI_n\big)&=\wh\bI_{n+3},
&\quad
\wt\io_n^{\,+}\big(\big\{\wt\bI_{n;1},\wh\bI_{n;1}\big\}\big)
&=\big\{\wt\bI_{n+3;4},\wh\bI_{n+3;4}\big\}.
\end{aligned}\EE
\end{crl}

\begin{proof} 
Let $q_n$ be the topological double cover as in~\eref{Spindfn_e}.
The double covers
\BE{Pin2Spin_e2}\io_n^{\,-\,*}q_{n+1}\!:\io_n^{\,-\,*}\Spin(n\!+\!1)\lra\O(n)
\quad\hbox{and}\quad
\io_n^{\,+\,*}q_{n+3}\!:\io_n^{\,+\,*}\Spin(n\!+\!3)\lra\O(n)\EE
are Lie group homomorphisms. 
The embeddings
$$\wt\io_{n+1;n}'\!:\Spin(n)\lra\Spin(n\!+\!1) \qquad\hbox{and}\qquad
\wt\io_{n+3;n}'\!:\Spin(n)\lra\Spin(n\!+\!3)$$
induce Lie groups isomorphisms
\begin{alignat*}{2}
\Spin(n)&\lra\io_n^{\,-\,*}\Spin(n\!+\!1)\big|_{\SO(n)}
\subset\SO(n)\!\times\!\Spin(n\!+\!1), &\quad
\wt{A}&\lra\big(q_n(\wt{A}),\wt\io_{n+1;n}'(\wt{A})\big),\\
\Spin(n)&\lra\io_n^{\,+\,*}\Spin(n\!+\!3)\big|_{\SO(n)}
\subset\SO(n)\!\times\!\Spin(n\!+\!3), &\quad
\wt{A}&\lra\big(q_n(\wt{A}),\wt\io_{n+3;n}'(\wt{A})\big),
\end{alignat*}
commuting with the projections to~$\SO(n)$.
By the discussion after~\eref{Pindfn_e}, 
each of the two projections in~\eref{Pin2Spin_e2} is thus isomorphic 
to either~$q_n^-$ or~$q_n^+$.\\

\noindent
By definition, $\wt\bI_{n;1}$ is an element of $\Pin^-(n)$ of order~4 so that
$$\io_n^-\big(q_n^-(\wt\bI_{n;1})\!\big)=\io_n^-\big(\bI_{n;1}\big)
=\bI_{n+1;2}=q_{n+1}\big(\big\{\wt\bI_{n+1;2},\wh\bI_{n+1;2}\big\}\big).$$
By the first statement of Lemma~\ref{SO2Spin_lmm} with $n$ replaced by~$n\!+\!1$, 
$\wt\bI_{n+1;2}$ and $\wh\bI_{n+1;2}$ are elements of~$\Spin(n\!+\!1)$ of order~4.
Thus, the first projection in~\eref{Pin2Spin_e2} is isomorphic to~$q_n^-$;
this is equivalent to the $\Pin^-$-case of the corollary.\\

\noindent
By definition, $\wt\bI_{n;1}$ is an element of $\Pin^+(n)$ of order~2 so that
$$\io_n^{\,+}\big(q_n^+(\wt\bI_{n;1})\!\big)
=\io_n^{\,+}\big(\bI_{n;1}\big)
=\bI_{n+3;4}=q_{n+3}\big(\big\{\wt\bI_{n+3;4},\wh\bI_{n+3;4}\big\}\big).$$
By the second statement of Lemma~\ref{SO2Spin_lmm} with $n$ replaced by~$n\!+\!3$, 
$\wt\bI_{n+3;4}$ and $\wh\bI_{n+3;4}$ are elements of $\Spin(n\!+\!3)$ of order~2.
Thus, the second projection in~\eref{Pin2Spin_e2} is isomorphic to~$q_n^+$;
this is equivalent to the $\Pin^+$-case of the corollary.
\end{proof}

The extensions $\wt\io_n^{\,-}$ of $\wt\io_{n+1;n}'$ and $\wt\io_n^{\,+}$ of $\wt\io_{n+3;n}'$
provided by Corollary~\ref{Pin2Spin_crl} can be defined explicitly as follows.
Let
$$n_-=n\!+\!1  \qquad\hbox{and}\qquad n_+=n\!+\!3.$$
By~\eref{SO2Spin2_e7} with $(n,m)$ replaced by $(n_{\pm},n)$, \eref{SO2Spin2_e3}
with $(n,m_2)$ replaced by~$(n_{\pm},n)$, and three consecutive applications of
these equations in the~$+$ case,
\begin{gather*}
\wt\io_{n_{\pm};n}'\big(\wt\bI_{n;1}\wt{A}\wt\bI_{n;1}^{-1}\big)=
\wt\bI_{n_{\pm};(n)}\wt\io_{n_{\pm};n}'\big(\wt{A}\big)\wt\bI_{n_{\pm};(n)}^{-1}, \quad
\wt\io_{n_-;n}'\big(\wt{A}\big)=
\wt\bI_{n_-;1}\wt\io_{n_-;n}'\big(\wt{A}\big)\wt\bI_{n_-;1}^{-1},\\
\wt\io_{n_+;n}'\big(\wt{A}\big)=
\wt\bI_{n_+;(n_+-2)}\wt\bI_{n_+;(n_+-1)}\wt\bI_{n_+;(n_+)}\wt\io_{n_+;n}'\big(\wt{A}\big)
\wt\bI_{n_-;(n_+)}^{-1}\wt\bI_{n_-;(n_+-1)}^{-1}\wt\bI_{n_-;(n_+-2)}^{-1}
\end{gather*}
for all $\wt{A}\!\in\!\Spin(n)$.
Since
$$\wt\bI_{n_-;(n)}\wt\bI_{n_-;1}\in\big\{\wt\bI_{n_-;2},\wt\bI_{n_-;2}\wh\bI_{n_-}\big\}
\quad\hbox{and}\quad
\wt\bI_{n_+;(n)}\wt\bI_{n_+;(n_+-2)}\wt\bI_{n_+;(n_+-1)}\wt\bI_{n_+;(n_+)}
\in\big\{\wt\bI_{n_+;2},\wt\bI_{n_+;2}\wh\bI_{n_+}\big\},$$
the three equations above and~\eref{Pindfn_e2} imply that
\BE{wtiocomm_e2}\begin{aligned}
\wt\io_{n+1;n}'\big(\wt\bI_{n;1}\wt{A}\wt\bI_{n;1}^{-1}\big)
&=\wt\bI_{n+1;2}\wt\io_{n+1;n}'\big(\wt{A}\big)\wt\bI_{n+1;2}^{-1}
&\quad&\forall\,\wt{A}\!\in\!\Spin(n),\\
\wt\io_{n+3;n}'\big(\wt\bI_{n;1}\wt{A}\wt\bI_{n;1}^{-1}\big)
&=\wt\bI_{n+3;4}\wt\io_{n+3;n}'\big(\wt{A}\big)\wt\bI_{n+3;4}^{-1}
&\quad&\forall\,\wt{A}\!\in\!\Spin(n)\,.
\end{aligned}\EE
The two possible extensions $\wt\io_n^-$ of $\wt\io_{n+1;n}'$ are given~by
\begin{gather*}
\Pin^-(n)\!-\!\Spin(n)\lra \Spin(n\!+\!1), \\
\wt{A}\wt\bI_{n;1}\lra \wt\io_{n+1;n}'(\wt{A})\wt\bI_{n+1;2},
\wt\io_{n+1;n}'(\wt{A})\wh\bI_{n+1;2}\quad\forall\,\wt{A}\!\in\!\Spin(n).
\end{gather*}
Since  $\wt\bI_{n;1}$ is an element of $\Pin^-(n)$ of order~4 and
$\wt\bI_{n+1;2}$ and $\wh\bI_{n+1;2}$ are elements of $\Spin(n\!+\!1)$ of order~4,
the first equation in~\eref{wtiocomm_e2} and the second equation in~\eref{wtioSpin_e2}
with $(n,m)$ replaced by~$(n_-,n)$,
ensure that these extensions are Lie group homomorphisms.
The two possible extensions $\wt\io_n^+$ of $\wt\io_{n+3;n}'$ can similarly be
defined~by
\begin{gather*}
\Pin^+(n)\!-\!\Spin(n)\lra \Spin(n\!+\!3), \\
\wt{A}\wt\bI_{n;1}\lra \wt\io_{n+3;n}(\wt{A})\wt\bI_{n+3;4},
\wt\io_{n+3;n}(\wt{A})\wh\bI_{n+3;4} \quad\forall\,\wt{A}\!\in\!\Spin(n).
\end{gather*}
Since  $\wt\bI_{n;1}$ is an element of $\Pin^+(n)$ of order~2 and
$\wt\bI_{n+3;4}$ and $\wh\bI_{n+3;4}$ are elements of $\Spin(n\!+\!3)$ of order~2,
the second equation in~\eref{wtiocomm_e2} ensures that these extensions 
are Lie group homomorphisms.\\

\noindent
There are thus two lifts~$\wt\io_n^{\,-}$ of~$\io_n^-$ and 
two lifts $\wt\io_n^{\,+}$ of~$\io_n^+$;
they differ by the composition with the involution~$\rho$ in~\eref{rhodfn_e}
on the right.
For concreteness, we choose them so~that 
\BE{Pin2Spindfn_e}  \wt\io_n^{\pm}\big(\wt\bI_{n;1}\big)=
\wt\bI_{n_{\pm};3\pm1}\in\Spin(n_{\pm})\,.\EE
Along with~\eref{wtioPindfn_e} and~\eref{I2244_e2}, this implies that 
\BE{StCorrComm_e0b}
\wt\io_n^{\pm}\big(\wt\io_{n;m}(\wt\bI_m,\wt\bI_{n-m;1})\big)
=\wt\io_{n_{\pm};m}\big(\wt\bI_m,\wt\io_{n-m}^{\pm}(\wt\bI_{n-m;1})\big)\EE
for all $m,n\!\in\!\Z^+$ with $m\!<\!n$.
By Lemma~\ref{CovExt_lmm}\ref{CovExtConn_it2}, 
\BE{StCorrComm_e0a}
\wt\io_n^{\pm}\!\circ\!\wt\io_{n;m}
\!=\!\wt\io_{n_{\pm};m}\!\circ\!\big\{\id_{\Spin(m)}\!\times\!\wt\io_{n-m}^{\pm}\big\}\!:
\Spin(m)\!\times\!\Spin(n\!-\!m)\lra\Spin(n_{\pm}),\EE 
since these two homomorphisms lift the homomorphism
$$\io_{n_{\pm};n}'\!\circ\!\io_{n;m}
\!=\!\io_{n_{\pm};m}\!\circ\!\big\{\id_{\SO(m)}\!\times\!\io_{n_{\pm}-m;n-m}'\big\}\!:
\SO(m)\!\times\!\SO(n\!-\!m)\lra\SO(n_{\pm}).$$
By~\eref{StCorrComm_e0a} and~\eref{StCorrComm_e0b},
\BE{StCorrComm_e} 
\wt\io_n^{\pm}\!\circ\!\wt\io_{n;m}
\!=\!\wt\io_{n_{\pm};m}\!\circ\!\big\{\id_{\Spin(m)}\!\times\!\wt\io_{n-m}^{\pm}\big\}\!:
\Spin(m)\!\times\!\Pin^{\pm}(n\!-\!m)\lra\Spin(n_{\pm}).\EE

\section{Proof of Theorem~\ref{SpinStrEquiv_thm}\ref{PinSpinPrp_it}: classical perspective}
\label{SpinPin_sec1}

In this section, which relies heavily on the notation and notions introduced in Section~\ref{LG_sec},
we establish the statements of Section~\ref{SpinPinProp_subs} for the notions of 
$\Spin$-structure and $\Pin^{\pm}$-structure arising from  Definition~\ref{PinSpin_dfn}
and give additional examples of these structures.
Throughout this section, the terms $\Spin$-structure and $\Pin^{\pm}$-structure refer to
these notions.
In Section~\ref{SpinPindfn_subs}, we clarify the terminology 
around Definition~\ref{PinSpin_dfn} and then relate $\Spin$- and $\Pin^{\pm}$-structures
to \v{C}ech cohomology with values in Lie groups. 
In Section~\ref{EquivClass_subs}, we use \v{C}ech cohomology to
establish statements~\ref{PinStr_it} and~\ref{SpinStr_it} of 
the SpinPin~\ref{SpinPinStr_prop} property
and the SpinPin~\ref{OSpinRev_prop} and~\ref{Pin2SpinRed_prop} properties.
We obtain the SpinPin~\ref{SpinPinObs_prop}, \ref{SpinPinStab_prop}, and
\ref{SpinPinCorr_prop} properties in Section~\ref{SpinPinprp_subs}
and the SpinPin~\ref{SpinPinSES_prop} property and the last statement of 
the SpinPin~\ref{SpinPinStr_prop} property in Section~\ref{ShExSeq_subs}.

\subsection{Spin/Pin-structures and \v{C}ech cohomology}
\label{SpinPindfn_subs}

Every rank~$n$ real vector bundle~$V$ over a paracompact space~$Y$ 
admits a metric~$g_V$.
A metric~$g_V$ on~$V$ determines an \sf{orthonormal frame 
bundle~$\O(V)$}\nota{ov@$\O(V)$|textbf}
over~$Y$;
the fiber of $\O(V)$ over a point $y\!\in\!Y$ consists of ordered tuples
$(v_1,\ldots,v_n)$ so that $\{v_1,\ldots,v_n\}$ is an orthonormal basis
for the fiber~$V_y$ of~$V$ over~$y$.
This is a principal $\O(n)$-bundle over~$Y$ so~that 
$$V=\O(V)\!\times_{\O(n)}\!\R^n\equiv
\big(\O(V)\!\times\!\R^n\big)\big/\!\!\sim, \quad
(p,c)\sim(pA^{-1},Ac)~~\forall\,A\!\in\!\O(n).$$
If $\fo\!\in\!\fO(V)$ is an orientation on~$V$,
then $g_V$ also determines an \sf{oriented orthonormal frame 
bundle~$\SO(V,\fo)$}\nota{sovo@$\SO(V,\fo)$|textbf} 
over~$Y$;
it is the subspace of $\O(V)$ consisting of the tuples $(v_1,\ldots,v_n)$
that are oriented bases for the fibers of~$V$ with respect to~$\fo$.
This is a principal $\SO(n)$-bundle over~$Y$ so~that 
$$V=\SO(V,\fo)\!\times_{\SO(n)}\!\R^n\equiv
\big(\SO(V,\fo)\!\times\!\R^n\big)\big/\!\!\sim, \quad
(p,c)\sim(pA^{-1},Ac)~~\forall\,A\!\in\!\SO(n).$$
The inclusion $\SO(V,\fo)\!\lra\!\O(V)$ is equivariant with respect to 
the inclusion $\SO(n)\!\lra\!\O(n)$.
Since the space of metrics~$g_V$ on a vector bundle~$V$ over a paracompact space~$Y$ 
is contractible,
the isomorphism classes of the frame bundles~$\O(V)$ and $\SO(V,\fo)$
are independent of the choice of~$g_V$.
Thus, the notions of $\Pin^{\pm}$-structure and $\Spin$-structure of 
Definition~\ref{PinSpin_dfn} are also independent of the choice of~$g_V$.\\

The next example illustrates the notion of equivalence of Spin-structures 
formulated after Definition~\ref{PinSpin_dfn}. 
Along with the SpinPin~\ref{Pin2SpinRed_prop} property,
this example also implies the last statement of the SpinPin~\ref{SpinPinStr_prop} property
with $V'$ orientable of even~rank and $\rk\,V''\!=\!0$.

\begin{eg}\label{negaaut_eg0}
Let $V$ be a vector bundle over a paracompact space~$Y$
with \hbox{$n\!\equiv\!\rk\,V$} even and
\BE{negaaut_e}\Psi\!:V\lra V, \qquad v\lra -v.\EE
Suppose $\fo\!\in\!\fO(V)$ and  $\fs\!\equiv\!(\Spin(V,\fo),q_V)$ is a Spin-structure 
on~$(V,\fo)$ as in Definition~\ref{PinSpin_dfn}.
Since $n$ is even, $\wt\bI_{n;n}\!\in\!\Spin(n)$ and the~map
$$\Psi^*\Spin(V,\fo)\!\equiv\!\big\{(p,\wt{p})\!:p\!=\!-q_V(\wt{p})\big\}\lra \Spin(V,\fo),
\quad (p,\wt{p})\lra\wt{p}\!\cdot\!\wt\bI_{n;n},$$
is an isomorphism of $\Spin(n)$-fiber bundles over~$Y$.
By~\eref{laminbIn_e}, this map is $\Spin(n)$-equivariant and is thus 
an equivalence of $\Spin$-structures on~$(V,\fo)$.
We conclude~that $[\Psi^*\fs]\!=\![\fs]$ in~$\Sp(V,\fo)$.
\end{eg}

Let $Y$ be a topological space, $P$ be a principal $\SO(n)$-bundle over~$Y$,
and $\wt{P}$ be a principal $\Spin(n)$-bundle over~$Y$.
Equivariant trivializations~$\Phi$ of~$P$ and~$\wt\Phi$ of $\wt{P}$ 
over $U\!\subset\!Y$ correspond to sections
$$s\in\Ga\big(U;P\big) \qquad\hbox{and}\qquad
\wt{s}\in\Ga\big(U;\wt{P}\big)$$
of the two bundles over~$U$:
\BE{sbsPhidfn_e}s\big(\pi_1\big(\Phi(p)\!\big)\!\big)\!\cdot\!\pi_2\big(\Phi(p)\!\big)=p
~~\forall\,p\!\in\!P|_U, \quad
\wt{s}\big(\pi_1\big(\wt\Phi(\wt{p})\!\big)\!\big)\!\cdot\!\pi_2\big(\wt\Phi(\wt{p})\!\big)
=\wt{p}
~~\forall\,\wt{p}\!\in\!\wt{P}|_U,\EE
where
$$\pi_1,\pi_2\!:U\!\times\!\SO(n)\lra U,\SO(n) \quad\hbox{and}\quad
\pi_1,\pi_2\!:U\!\times\!\Spin(n)\lra U,\Spin(n)$$
are the projection maps.
If $q\!:\wt{P}\!\lra\!P$ is a double cover satisfying the conditions
of Definition~\ref{PinSpin_dfn}\ref{SpinStrDfn_it},
a $\Spin(n)$-equivariant trivialization~$\wt\Phi$ of~$\wt{P}$ over~$U$ thus induces 
an $\SO(n)$-equivariant trivialization~$\Phi$ of~$P$
over~$U$ so that the~diagram
\BE{SpinVSOV_e}\begin{split}
\xymatrix{ \wt{P}|_U \ar[r]^>>>>>{\wt\Phi}\ar[d]_q&  
U\!\times\!\Spin(n) \ar[d]^{\id\times q_n}\\
P|_U \ar[r]^>>>>>{\Phi}& U\!\times\!\SO(n)}
\end{split}\EE
commutes.\\

Let $\{U_{\al}\}_{\al\in\cA}$ be an open cover of~$Y$ and $q\!:\wt{P}\!\lra\!P$ be as above.
For $\al_0,\al_1,\ldots,\al_p\!\in\!\cA$, we~define
$$U_{\al_0\al_1\ldots\al_p}=U_{\al_0}\!\cap\!U_{\al_1}\!\cap\!\ldots\!\cap\!U_{\al_p}
\subset Y\,.$$
A collection 
\BE{SpinVPsi_e}\big\{\wt\Phi_{\al}\!:\wt{P}|_{U_{\al}}
\lra U_{\al}\!\times\!\Spin(n)\big\}_{\al\in\cA}\EE
of $\Spin(n)$-equivariant trivializations determines transition data
\begin{gather}\label{SpinVg_e}
\big\{\wt{g}_{\al\be}\!: U_{\al\be}\!\lra\!\Spin(n)\big\}_{\al,\be\in\cA}
\qquad\hbox{s.t.}\\
\notag
\wt\Phi_{\al}(\wt{p})=\big(\pi_{\be;1}(\wt\Phi_{\be}(\wt{p})),
\wt{g}_{\al\be}\big(\pi_{\be;1}(\wt\Phi_{\be}(\wt{p}))\big)\pi_{\be;2}(\wt\Phi_{\be}(\wt{p}))\big)
~~\forall\,\wt{p}\!\in\!\wt{P}|_{U_{\al\be}},\,\al,\be\!\in\!\cA,
\end{gather}
where 
$$\pi_{\be;1},\pi_{\be;2}\!:U_{\be}\!\times\!\Spin(n)\lra U_{\be},\Spin(n)$$
are the projection maps.
The collection~\eref{SpinVg_e} satisfies
\BE{SpinVg_e2}\wt{g}_{\be\ga}\big|_{U_{\al\be\ga}}
\wt{g}_{\al\ga}^{\,-1}\big|_{U_{\al\be\ga}}
\wt{g}_{\al\be}\big|_{U_{\al\be\ga}}=\wt\bI_n
\qquad\forall\,\al,\be,\ga\!\in\!\cA.\EE
The collection
\BE{SOVg_e}\big\{g_{\al\be}\!=\!q_n(\wt{g}_{\al\be})\!:
U_{\al\be}\!\lra\!\SO(n)\big\}_{\al,\be\in\cA}\EE
then consists of the transition data for the principal $\SO(n)$-bundle $P$
induced by the collection 
\BE{SOVg_e0}\big\{\Phi_{\al}\!:P\big|_{U_{\al}}\lra U_{\al}\!\times\!\SO(n)\big\}_{\al\in\cA}\EE
of its trivializations determined by the collection~\eref{SpinVPsi_e}.\\

Every collection $\{\wt{g}_{\al\be}\}_{\al,\be\in\cA}$
as in~\eref{SpinVg_e} satisfying~\eref{SpinVg_e2} conversely determines 
a principal $\Spin(n)$-bundle 
\begin{gather*}
\wt{P}\equiv\bigg(\bigsqcup_{\al\in\cA}\!U_{\al}\!\times\!\Spin(n)\!\!\bigg)\!\!\Big/\!\!\!\sim
\,\lra Y, \quad \big[y,\wt{A}\big]\lra y,\\
U_{\al}\!\times\!\Spin(n)\!\ni\!\big(y,\wt{g}_{\al\be}(y)\wt{A}\big)\sim
\big(y,\wt{A}\big) \!\in\! U_{\be}\!\times\!\Spin(n)
\quad\forall\,\big(y,\wt{A}\big)\!\in\!U_{\al\be}\!\times\!\Spin(n),
~\al,\be\!\in\!\cA.
\end{gather*}
Every collection $\{g_{\al\be}\}_{\al,\be\in\cA}$
as in~\eref{SOVg_e} satisfying~\eref{SpinVg_e2} with all tildes dropped
similarly determines a principal $\SO(n)$-bundle 
\begin{gather*}
P\equiv\bigg(\bigsqcup_{\al\in\cA}\!U_{\al}\!\times\!\SO(n)\!\!\bigg)\!\!\Big/\!\!\!\sim
\,\lra Y, \quad \big[y,A\big]\lra y,\\
U_{\al}\!\times\!\SO(n)\!\ni\!\big(y,g_{\al\be}(y)A\big)\sim
\big(y,A\big) \!\in\! U_{\be}\!\times\!\SO(n)
\quad\forall\,
\big(y,A\big) \!\in\! U_{\al\be}\!\times\!\SO(n),~\al,\be\!\in\!\cA.
\end{gather*}
Thus, every collection as in~\eref{SpinVg_e} satisfying~\eref{SpinVg_e2} determines
a double cover
\BE{qVdfn_e}q\!:\wt{P}\lra P, \quad
q\big([y,\wt{A}]\big)=\big[y,q_n(\wt{A})\big]~~
\forall\,\big(y,\wt{A}\big)\!\in\!U_{\al}\!\times\!\Spin(n),\,\al\!\in\!\cA,\EE
satisfying the conditions of Definition~\ref{PinSpin_dfn}\ref{SpinStrDfn_it}.\\

If $q'\!:\wt{P}\!\lra\!P$ is another double cover
satisfying the conditions of Definition~\ref{PinSpin_dfn}\ref{SpinStrDfn_it},
there exists a collection
\begin{gather}\notag
\{f_{\al}\!:U_{\al}\!\lra\!\SO(n)\big\}_{\al\in\cA}
\qquad\hbox{s.t.}\\
\label{qVpsial_e}
q'\big([y,\wt{A}]\big)=\big[y,f_{\al}(y)q_n(\wt{A})\big] \quad
\forall\,\big(y,\wt{A}\big)\!\in\!U_{\al}\!\times\!\Spin(n),\,\al\!\in\!\cA.
\end{gather}
This collection satisfies
\BE{SOVaut_e}g_{\al\be}=
f_{\al}\big|_{U_{\al\be}}g_{\al\be}f_{\be}^{-1}\big|_{U_{\al\be}}
\quad\forall\,\al,\be\!\in\!\cA.\EE
The $\Spin$-structures $q$ and $q'$ are equivalent if
there exists a collection
\begin{gather*}
\big\{\wt{f}_{\al}\!:U_{\al}\!\lra\!\Spin(n)\big\}_{\al\in\cA}
\qquad\hbox{s.t.}\\
\wt{g}_{\al\be}=
\wt{f}_{\al}\big|_{U_{\al\be}}\wt{g}_{\al\be}\wt{f}_{\be}^{-1}\big|_{U_{\al\be}}
\quad\forall\,\al,\be\!\in\!\cA, \qquad
q_n\!\circ\!\wt{f}_{\al}=f_{\al}~~\forall\,\al\!\in\!\cA;
\end{gather*}
such a collection $\{\wt{f}_{\al}\}$ induces a $\Spin(n)$-equivariant automorphism
$$\wt\Psi\!:\wt{P}\lra\wt{P}, \qquad 
\wt\Psi\big([y,\wt{A}]\big)=\big[x,\wt{f}_{\al}(y)\big]  \quad
\forall\,\big(y,\wt{A}\big)\!\in\!U_{\al}\!\times\!\Spin(n),\,\al\!\in\!\cA,$$
satisfying $q'\!=\!q\!\circ\!\wt\Psi$.\\

Combining the last three paragraphs with Section~\ref{CechPB_subs}, we obtain a~map
$$\wcH^1\big(Y;\Spin(n)\big)\lra \wcH^1\big(Y;\SO(n)\big)$$
from the set of equivalence classes of principal $\Spin(n)$-bundles
to the  set of equivalence classes of principal $\SO(n)$-bundles
so that for every principal $\Spin(n)$-bundle $\wt{P}$ 
in the preimage of a principal $\SO(n)$-bundle $P$ 
there exists a double cover $q\!:\wt{P}\!\lra\!P$ satisfying the conditions
of Definition~\ref{PinSpin_dfn}\ref{SpinStrDfn_it}.
However, two such double covers need not be equivalent;
see the last part of Section~\ref{EquivClass_subs}.
The analogous considerations apply in the setting of 
Definition~\ref{PinSpin_dfn}\ref{PinStrDfn_it} with the~map
$$\wcH^1\big(Y;\Pin^{\pm}(n)\big)\lra \wcH^1\big(Y;\O(n)\big)$$
from the set of equivalence classes of principal $\Pin^{\pm}(n)$-bundles
to the  set of equivalence classes of principal $\O(n)$-bundles
induced by the Lie group homomorphism~\eref{Pindfn_e}.

\subsection{The sets $\cP^{\pm}(V)$ and $\Sp(V,\fo)$}
\label{EquivClass_subs}

We now apply the \v{C}ech cohomology perspective of the previous section
to describe the structure of the sets $\cP^{\pm}(V)$ 
of $\Pin^{\pm}$-structures on~$V$ and 
$\Sp(V,\fo)$ of $\Spin$-structures on~$(V,\fo)$, up to equivalence.
By Proposition~\ref{PinSpin_prp1} below, 
the statements~\ref{PinStr_it} and~\ref{SpinStr_it} of 
the SpinPin~\ref{SpinPinStr_prop} property hold if $Y$ is paracompact and satisfies
the $k\!=\!1$ case of the mild topological condition in the following definition.

\begin{dfn}\label{locH1simp_dfn}
Let $k\!\in\!\Z^{\ge0}$.
A topological space~$Y$ is \sf{locally $H^k(\cdot;\Z_2)$-simple} 
if it is locally path-connected and
\BE{locH1simp_e}\varinjlim_{U\ni y}H^p(U;\Z_2)=0 \qquad\forall~
y\!\in\!Y,~p\!=\!1,\ldots,k,\EE
where the direct limit is taken with respect to the cohomology homomorphisms
induced by the inclusions of neighborhoods~$U$ of~$y$.
\end{dfn}

The direct limit in~\eref{locH1simp_e} vanishes if and only if 
for every neighborhood $U\!\subset\!Y$ of~$y$ and every $\eta\!\in\!H^p(U;\Z_2)$ 
there exists a neighborhood $U'\!\subset\!U$ 
of~$y$ such that $\eta|_{U'}\!=\!0$.
If $Y$ is locally $H^1$-simple, then the sheaf homomorphisms
$$\cS_Y\big(\Spin(n)\big)\lra \cS_Y\big(\SO(n)\big)
\qquad\hbox{and}\qquad 
\cS_Y\big(\Pin^{\pm}(n)\big)\lra \cS_Y\big(\O(n)\big)$$
between the sheaves of germs of continuous functions on~$Y$ with values in
$\Spin(n)$, $\SO(n)$, $\Pin^{\pm}(n)$, and~$\O(n)$ induced by
the Lie group homomorphisms~\eref{Spindfn_e} and~\eref{Pindfn_e} are surjective;
see the paragraph after Definition~\ref{locLGsimp_dfn}.
All CW-complexes are locally $H^k$-simple for every $k\!\in\!\Z^+$ and paracompact.

\begin{prp}\label{PinSpin_prp1}
Let $V$ be a vector bundle over a paracompact locally $H^1$-simple space~$Y$.
\begin{enumerate}[label=(\alph*),leftmargin=*]

\item\label{Pin1_it}  If\, $V$ admits a $\Pin^{\pm}$-structure, then
the group $H^1(Y;\Z_2)$ acts naturally, freely, and transitively on the set~$\cP^{\pm}(V)$.

\item\label{Spin1_it}  If $\fo\!\in\!\fO(V)$ and $(V,\fo)$
admits a $\Spin$-structure, then
the group $H^1(Y;\Z_2)$ acts naturally, freely, and transitively on the set~$\Sp(V,\fo)$.

\end{enumerate}
\end{prp}

\begin{proof}
By Proposition~\ref{CechSing_prp},
it is sufficient to establish the two claims with $H^1(Y;\Z_2)$
replaced by~$\wcH^1(Y;\Z_2)$.
Let $n\!=\!\rk\,V$,  $\fo\!\in\!\fO(V)$, 
and $\fs\!\equiv\!(\Spin(V,\fo),q_V)$ be a Spin-structure on~$(V,\fo)$.\\

Suppose $\eta\!\in\!\wcH^1(Y;\Z_2)$. 
Choose a collection $\{\wt{g}_{\al\be}\}$ of transition data 
for the principal $\Spin(n)$-bundle 
$$\wt{P}\equiv\Spin(V,\fo)$$ 
as in~\eref{SpinVg_e} and~\eref{SpinVg_e2} 
so that the collection~\eref{SOVg_e} consists of transition data for
the principal $\SO(n)$-bundle \hbox{$P\!\equiv\!\SO(V,\fo)$} and $q_V$ is given
by~\eref{qVdfn_e}.
Refining the open cover $\{U_{\al}\}$ if necessary, 
we can assume that~$\eta$ is represented by a collection
\BE{PinSpin1_e3}
\big\{\eta_{\al\be}\!:U_{\al\be}\!\lra\!\Z_2\big\}_{\al,\be\in\cA}
\quad\hbox{s.t.}\quad
\eta_{\be\ga}\big|_{U_{\al\be\ga}}\!-\!\eta_{\al\ga}\big|_{U_{\al\be\ga}}
\!+\!\eta_{\al\be}\big|_{U_{\al\be\ga}}=0~~\forall\,\al,\be,\ga\!\in\!\cA.\EE
We take 
\BE{PinSpin1_e4} \eta\!\cdot\!\fs \equiv \fs'\!\equiv\!\big(\Spin'(V,\fo),q_V'\big)\EE
to be the $\Spin$-structure on~$(V,\fo)$ determined by the transition~data
\BE{PinSpin1_e5} \big\{\wt{g}_{\al\be}'\!\equiv\!\wt{g}_{\al\be}\wh\bI_n^{\,\eta_{\al\be}}\!:
U_{\al\be}\!\lra\!\Spin(n)\big\}_{\al,\be\in\cA}\,.\EE
This collection satisfies~\eref{SpinVg_e2} with~$\wt{g}_{\al\be}$
replaced by~$\wt{g}_{\al\be}'$ because $\wh\bI_n$ is in the center of~$\Spin(n)$.
For the same reason, the equivalence class of~$\fs'$ depends
only on the equivalence classes of~$\fs$ and on~$\eta$.
It is immediate that this defines a group action of $\wcH^1(Y;\Z_2)$
on~$\Sp(V,\fo)$.\\
 
Suppose the $\Spin$-structure~$\fs'$ in~\eref{PinSpin1_e4}
is equivalent to~$\fs$.
Refining the open cover $\{U_{\al}\}$ if necessary, 
we can assume that there exists a collection
\BE{PinSpin1_e6}\big\{\wt{f}_{\al}\!:U_{\al}\!\lra\!\Spin(n)\big\}_{\al\in\cA}\EE
of continuous maps such~that 
\BE{PinSpin1_e7}\begin{aligned}
\wt{g}_{\al\be}\wh\bI_n^{\,\eta_{\al\be}}&=\wt{f}_{\al}\big|_{U_{\al\be}}\wt{g}_{\al\be}
\wt{f}_{\be}^{-1}\big|_{U_{\al\be}}
&\quad&\forall~\al,\be\!\in\!\cA, \\
q_n\big(\wt{A}\big)&=q_n\big(\wt{f}_{\al}(x)\wt{A}\big)
&\quad&
\forall~\big(x,\wt{A}\big)\!\in\!U_{\al}\!\times\!\Spin(n),\,\al\!\in\!\cA.
\end{aligned}\EE
The first condition above means that the collection~\eref{PinSpin1_e6} 
determines an equivalence between the principal $\Spin(n)$-bundles
determined by the collections~\eref{SpinVg_e} and~\eref{PinSpin1_e5}. 
The second condition means that this equivalence commutes with 
the corresponding projections~$q_V$ and~$q_V'$ to
the principal $\SO(n)$-bundle $\SO(V,\fo)$.
It implies that for every $\al\!\in\!\cA$ there exists a continuous function
$$\mu_{\al}\!:U_{\al}\lra \Z_2\qquad\hbox{s.t.}\quad \wt{f}_{\al}=\wh\bI_n^{\,\mu_{\al}}\,.$$
Since $\wh\bI_n$ is in the center of~$\Spin(n)$, this 
and the first condition in~\eref{PinSpin1_e7} imply~that 
$$\eta_{\al\be}=\mu_{\al}\big|_{U_{\al\be}}-\mu_{\be}\big|_{U_{\al\be}}
\qquad\forall~\al,\be\!\in\!\cA,$$
i.e.~$\eta$ represents the trivial element of $\wcH^1(Y;\Z_2)$.  
Thus, $\wcH^1(Y;\Z_2)$ acts freely on~$\Sp(V,\fo)$.\\ 

Suppose $\fs''\!\equiv\!(\Spin''(V,\fo),q_V'')$ is another $\Spin$-structure on~$(V,\fo)$.
Refining the open cover $\{U_{\al}\}$ if necessary, 
we can assume that this pair is determined by the collections 
$$\big\{\wt{g}_{\al\be}''\!:U_{\al\be}\!\lra\!\Spin(n)\big\}
\qquad\hbox{and}\quad
\big\{f_{\al}\!:U_{\al}\!\lra\!\SO(n)\big\}_{\al\in\cA}$$
satisfying~\eref{SpinVg_e2} with~$\wt{g}$ replaced by~$\wt{g}''$,
\eref{qVpsial_e} with $q'$ replaced by~$q_V''$, \eref{SOVaut_e}, and 
\BE{PinSpin1_e15}q_n(\wt{g}_{\al\be}'')=g_{\al\be} \qquad\forall\,\al,\be\!\in\!\cA.\EE
If $n\!\ge\!2$, let $\eta_n\!\in\!H^1(\SO(n);\Z_2)$ denote the generator.
Since $Y$ is locally $H^1$-simple, we can assume~that
$$f_{\al}^*\eta_n=0\in H^1\big(U_{\al};\Z_2\big)
\qquad\forall	\,\al\!\in\!\cA$$
and thus there exists a collection
\BE{PinSpin1_e16}\big\{\wt{f}_{\al}\!:U_{\al}\!\lra\!\Spin(n)\big\}_{\al\in\cA}
\qquad\hbox{s.t.}\quad
q_n\!\circ\!\wt{f}_{\al}=f_{\al}~~\forall\,\al\!\in\!\cA\,.\EE
If $n\!=\!1$, the existence of such a collection is immediate. 
By~\eref{SOVaut_e} and~\eref{PinSpin1_e15}, the existence of the collection~\eref{PinSpin1_e16}
implies that 
there exists a collection as in~\eref{PinSpin1_e3} such that 
\BE{PinSpin1_e19}\wt{g}_{\al\be}''=\wt{f}_{\al}\big|_{U_{\al\be}}
\big(\wt{g}_{\al\be}\wh\bI_n^{\,\eta_{\al\be}}\big)\wt{f}_{\be}^{\,-1}\big|_{U_{\al\be}}
\qquad\forall\,\al,\be\!\in\!\cA.\EE
Since $\{\wt{g}_{\al\be}\}_{\al\be}$ satisfies~\eref{SpinVg_e2},
$\{\wt{g}_{\al\be}''\}_{\al\be}$ satisfies~\eref{SpinVg_e2}
with~$\wt{g}$ replaced by~$\wt{g}''$, and $\wh\bI_n$ is in the center of $\Spin(n)$,
the collection $\{\eta_{\al\be}\}$ satisfies the (cocycle) condition
in~\eref{PinSpin1_e3} and thus determines an element~$\eta$ of $\wcH^1(Y;\Z_2)$.
By~\eref{PinSpin1_e19} and the condition in~\eref{PinSpin1_e16}, 
the collection in~\eref{PinSpin1_e16} determines an equivalence between
the $\Spin$-structures~$\eta\!\cdot\!\fs$ and~$\fs''$.	
Thus, $\wcH^1(Y;\Z_2)$ acts transitively on~$\Sp(V,\fo)$.\\

\noindent
This concludes the proof of~\ref{Spin1_it}.
The same reasoning, with $\Spin$ and $\SO(n)$ replaced by~$\Pin^{\pm}$ and~$\O(n)$
everywhere, applies to~\ref{Pin1_it}.
\end{proof}

\begin{eg}\label{negaautPin_eg}
Suppose $V$ is a vector bundle over a paracompact locally $H^1$-simple space~$Y$,
\hbox{$n\!\equiv\!\rk\,V$}, $\Psi$ is as in~\eref{negaaut_e}, 
and $\fp\!\equiv\!(\Pin^{\pm}(V),q_V)$ is a $\Pin^{\pm}$-structure on~$V$.
Let 
$$\big\{\wt{g}_{\al\be}\!:
U_{\al\be}\!\lra\!\Pin^{\pm}(n)\big\}_{\al,\be\in\cA}
\qquad\hbox{and}\qquad 
\big\{g_{\al\be}\!=\!q_n^{\pm}(\wt{g}_{\al\be})\!:
U_{\al\be}\!\lra\!\O(n)\big\}_{\al,\be\in\cA}$$
be a collection of transition data for the principal $\Pin^{\pm}(n)$-bundle 
$\Pin^{\pm}(V)$ and the associated collection of transition data for 
the principal $\O(n)$-bundle~$\O(V)$.
Define 
$$\big\{\eta_{\al\be}\!:U_{\al\be}\!\lra\!\Z_2\big\}_{\al,\be\in\cA}
\qquad\hbox{by}\quad (-1)^{\eta_{\al\be}}=\big(\det g_{\al\be}\big)^{n-1}\,.$$
By Corollary~\ref{w1_crl}, the \v{C}ech cohomology class~$[\eta]$ 
represented by $\{\eta_{\al\be}\}$ is $(n\!-\!1)w_1(V)$.
The $\Pin^{\pm}$-structure 
$$\eta\!\cdot\!\fp \equiv \fp'\!\equiv\!\big(\Pin'^{\pm}(V),q_V'\big)$$
on~$V$ is determined by the transition~data
$$\big\{\wt{g}_{\al\be}'\!\equiv\!\wt{g}_{\al\be}\wh\bI_n^{\,\eta_{\al\be}}\!:
U_{\al\be}\!\lra\!\Pin^{\pm}(n)\big\}_{\al,\be\in\cA}\,.$$
By~\eref{laminbIn_e} and the second equation in~\eref{wtInmsq_e2},  
$\wt\bI_{n;n}\wt{g}_{\al\be}\!=\!\wt{g}_{\al\be}'\wt\bI_{n;n}$.
Thus, the~map
\begin{gather*}
\wt\Psi\!:
\Psi^*\Pin^{\pm}(V)\!\equiv\!\big\{(p,\wt{p})\!:p\!=\!-q_V(\wt{p})\big\}\lra \Pin'^{\pm}(V),\\
\wt\Psi\big([x,-q_n(\wt{A})],[x,\wt{A}]\big)= \big[x,\wt\bI_{n;n}\wt{A}\big]
\quad\forall~\big(x,\wt{A}\big)\!\in\!U_{\al}\!\times\!\Pin^{\pm}(n),~\al\!\in\!\cA,
\end{gather*}
is well-defined and $\Pin^{\pm}(n)$-equivariant and
$\Psi^*q_V\!=\!q_V'\!\circ\!\wt\Psi$.
We conclude~that 
$$[\Psi^*\fp]=\big((n\!-\!1)w_1(V)\big)\!\cdot\![\fp],$$
establishing the $\rk\,V''\!=\!0$ case of~\eref{SpinPinStr_e}.
\end{eg}

\begin{proof}[{\bf{\emph{Proof of SpinPin~\ref{OSpinRev_prop} property}}}]
Let $V$ be a rank~$n$ vector bundle over a paracompact space~$Y$,
$\fo\!\in\!\fO(V)$, and  $\fs\!\equiv\!(\Spin(V,\fo),q_V)$ be a Spin-structure 
on~$(V,\fo)$.
Let 
$$\big\{\wt{s}_{\al}\!:U_{\al}\!\lra\!\Spin(V,\fo)\big\}_{\al\in\cA}$$ 
be a collection of local sections so that  $\{U_{\al}\}_{\al\in\cA}$
is an open cover of~$Y$, 
\eref{SpinVg_e} be the transition data for the principal $\Spin(n)$-bundle $\Spin(V,\fo)$
determined by the associated collection~\eref{SpinVPsi_e} of trivializations, 
and \eref{SOVg_e} be the transition data for the principal $\SO(n)$-bundle 
$$\SO(V,\fo)\subset \O(V)$$
determined by the associated collection~\eref{SOVg_e0} of trivializations.
Let $s_{\al}\!=\!q_V\!\circ\!\wt{s}_{\al}$.\\

The sections \hbox{$\ov{s}_{\al}\!\equiv\!s_{\al}\!\cdot\!\bI_{n;1}$} 
determine trivializations 
$\bI_{n;1}^{-1}\!\cdot\!\Phi_{\al}$ of the oriented orthonormal frame bundle 
$$\SO(V,\ov\fo)\subset \O(V)$$
for the opposite orientation~$\ov\fo$ on~$V$ with the transition data
$$\big\{\bI_{n;1}^{-1}g_{\al\be}\bI_{n;1}\!:
U_{\al\be}\!\lra\!\SO(n)\big\}_{\al,\be\in\cA}\,.$$
The collection 
\BE{conjSpin_e}\big\{\wh{g}_{\al\be}\!\equiv\!\wt\bI_{n;1}^{-1}\wt{g}_{\al\be}\wt\bI_{n;1}\!:
U_{\al\be}\!\lra\!\Spin(n)\big\}_{\al,\be\in\cA}\EE
then determines a $\Spin$-structure~$\ov\fs$ on $(V,\ov\fo)$,
\begin{gather*}
\ov{q}_V\!:
\Spin(V,\ov\fo)\equiv\bigg(\bigsqcup_{\al\in\cA}\!U_{\al}\!\times\!\Spin(n)\!\!\bigg)\!\!\Big/\!\!\!\sim
\,\lra  \SO(V,\ov\fo),
\quad \ov{q}_V\big([y,\wt{A}]\big)=s_{\al}(y)\!\cdot\!\bI_{n;1}q_n(\wt{A}),\\
U_{\al}\!\times\!\Spin(n)\!\ni\!\big(y,\wh{g}_{\al\be}(y)\wt{A}\big)\sim
\big(y,\wt{A}\big) \!\in\! U_{\be}\!\times\!\Spin(n)
\quad\forall\,\big(y,\wt{A}\big)\!\in\!U_{\al\be}\!\times\!\Spin(n),
~\al,\be\!\in\!\cA.
\end{gather*}
The equivalence class of~$\ov\fs$ depends only on the equivalence class of~$\fs$.\\

The above construction induces a bijection from $\Sp(V,\fo)$ to $\Sp(V,\ov\fo)$.
This bijection is natural with respect to continuous maps
and isomorphisms of vector bundles.
Since $\wh\bI_n$ lies in the center of $\Spin(n)$,
this bijection is equivariant with respect to the $H^1(Y;\Z_2)$-actions
of Proposition~\ref{PinSpin_prp1}.
It is immediate that this bijection satisfies the last condition of 
the SpinPin~\ref{OSpinRev_prop} property.
\end{proof}

\begin{rmk}\label{OSpinRev_rmk}
Since the space $\O(n)\!-\!\SO(n)$ is path-connected,
the equivalence class of $\ov\fs$ above would be the same if it were constructed using
any element~$\wt{B}$ of $\Pin^{\pm}(n)\!-\!\Spin(n)$ instead of~$\bI_{n;1}$.
An explicit equivalence between the two resulting Spin-structures would be induced by 
the~maps
$$U_{\al}\!\times\!\Spin(n)\lra U_{\al}\!\times\!\Spin(n), \qquad
\big(y,\wt{A}\big)\lra \big(y,\wt{B}^{-1}\wt\bI_{n;1}\wt{A}\big).$$
The bijections~\eref{OSpinRev_e} constructed in the proofs of the SpinPin~\ref{OSpinRev_prop} 
property in the perspectives of Definitions~\ref{PinSpin_dfn2} and~\ref{PinSpin_dfn3}
in Sections~\ref{SpinPindfn2_subs} and~\ref{SpinPindfn2_subs}, respectively,
likewise would not change if they were constructed using any element of $\O(n)\!-\!\SO(n)$
instead of~$\bI_{n;1}$. 
\end{rmk}

\begin{eg}\label{negaaut_eg}
Suppose $V$ is a vector bundle over a paracompact locally $H^1$-simple space~$Y$
with \hbox{$n\!\equiv\!\rk\,V$} odd, $\Psi$ is as in~\eref{negaaut_e},
$\fo\!\in\!\fO(V)$, and  $\fs\!\equiv\!(\Spin(V,\fo),q_V)$ is a Spin-structure 
on~$(V,\fo)$.
Let 
$$\big\{\wt{s}_{\al}\!:U_{\al}\!\lra\!\Spin(V,\fo)\big\}_{\al\in\cA}\qquad\hbox{and}\qquad 
\ov\fs\equiv(\Spin(V,\ov\fo),\ov{q}_V)$$ 
be as in the proof of the SpinPin~\ref{OSpinRev_prop} property above.
By~\eref{laminbIn_e},  
$$\wt\bI_{n;1}^{-1}\wt\bI_{n;n}\wt{g}_{\al\be}=\wh{g}_{\al\be}\wt\bI_{n;1}^{-1}\wt\bI_{n;n}.$$
Thus, the~map
\begin{gather*}
\wt\Psi\!:
\Psi^*\Spin(V,\fo)\!\equiv\!\big\{(p,\wt{p})\!:p\!=\!-q_V(\wt{p})\big\}\lra\Spin(V,\ov\fo),\\
\wt\Psi\big(s_{\al}(y)\!\cdot\!(-q_n(\wt{A})),\wt{s}_{\al}(y)\!\cdot\!\wt{A}\big)
= \big[y,\wt\bI_{n;1}^{-1}\wt\bI_{n;n}\wt{A}\big]
\quad\forall~\big(y,\wt{A}\big)\!\in\!U_{\al}\!\times\!\Spin(n),~\al\!\in\!\cA,
\end{gather*}
is well-defined and $\Spin(n)$-equivariant and
$\Psi^*q_V\!=\!\ov{q}_V\!\circ\!\wt\Psi$.
We conclude~that $[\Psi^*\fs]\!=\![\ov\fs]$ in~$\Sp(V,\ov\fo)$.
\end{eg}

\begin{proof}[{\bf{\emph{Proof of SpinPin~\ref{Pin2SpinRed_prop} property}}}]
Let $V$ be a rank~$n$ vector bundle over a paracompact space~$Y$ and 
$\fo\!\in\!\fO(V)$. 
Thus, 
\BE{Spin2Pin_e3}\O(V)\approx\SO(V,\fo)\!\times\!\Z_2\EE
as topological spaces.
By the SpinPin~\ref{SpinPinObs_prop} property,
$V$ admits a $\Pin^{\pm}$-structure if and only if
$(V,\fo)$ admits a $\Spin$-structure.
We can thus assume that there exists a $\Pin^{\pm}$-structure 
$$\fp\equiv\big(\Pin^{\pm}(V),q_V\big)$$ 
on~$V$ as in 
Definition~\ref{PinSpin_dfn}\ref{PinStrDfn_it}.
The restriction
$$q_V\big|_{q_V^{-1}(\SO(V,\fo)\times\{0\})}\!:
\Spin(V,\fo)\!\equiv\!q_V^{-1}\big(\SO(V,\fo)\!\times\!\{0\}\big)\lra Y$$
is then a $\Spin$-structure on~$(V,\fo)$.
This induces a well-defined natural $H^1(Y;\Z_2)$-equivariant
map~$\fR_{\fo}^{\pm}$ from $\cP^{\pm}(V)$ to~$\Sp(V,\fo)$.
Since the $H^1(Y;\Z_2)$-actions of Proposition~\ref{PinSpin_prp1} are free and transitive,
this map is bijective.\\

We now establish the last claim of the SpinPin~\ref{Pin2SpinRed_prop} property.
With $\fp$ and $\Spin(V,\fo)$ as above, let 
$\wt{s}_{\al}$, $\wt{g}_{\al\be}$, and $\ov{q}_V$ be as in the proof 
of the SpinPin~\ref{OSpinRev_prop} property.
In particular, $\wt{s}_{\al}\wt{g}_{\al\be}\!=\!\wt{s}_{\be}$ on~$U_{\al\be}$.
The~map
\begin{gather*}
\wt\Psi\!: \Spin(V,\ov\fo)\lra 
q_V^{-1}\big(\SO(V,\ov\fo)\!\times\!\{0\}\big)\subset\Pin^{\pm}(V),\\
\wt\Psi\big([y,\wt{A}]\big)=\wt{s}_{\al}(y)\!\cdot\!\wt\bI_{n;1}\wt{A}
\quad\forall\,\big(y,\wt{A}\big)\!\in\!U_{\al}\!\times\!\Spin(n),~\al\!\in\!\cA,
\end{gather*}
is then well-defined and $\Spin(n)$-equivariant and
$\ov{q}_V\!=\!q_V\!\circ\!\wt\Psi$.
Thus, the $\Spin$-structure $\ov\fs$ on~$(V,\ov\fo)$ constructed in  the proof 
of the SpinPin~\ref{OSpinRev_prop} property
is equivalent 
to the $\Spin$-structure obtained by restricting~$q_V$ to 
$q_V^{-1}(\SO(V,\ov\fo)\!\times\!\{0\})$.
\end{proof}

\noindent
For a topological space $Y$ and a Lie group~$G$, let 
$C(Y;G)$ denote the group of continuous maps from $Y$ to a Lie group~$G$.
If $Y$ is paracompact locally $H^1$-simple, the short exact sequence
\BE{SpinSOseq_e} \{\1\}\lra \Z_2\lra\Spin(n)\stackrel{q_n}{\lra} \SO(n)\lra \{\1\}\EE
of Lie groups induces an exact sequence
\BE{SpinLES_e}\begin{split}
\xymatrix{ \{\1\}\ar[r]& C(Y;\Z_2) \ar[r]& C\big(Y;\Spin(n)\big)\ar[r]^{\wch{q}_0}& 
C\big(Y;\SO(n)\big)\ar[lld]_{\wch\de_0}& \\
&H^1(Y;\Z_2) \ar[r]&
\wcH^1\big(Y;\Spin(n)\big)\ar[r]^{\wch{q}_1}& 
\wcH^1\big(Y;\SO(n)\big) \ar[r]^>>>>>>{\wch\de_1}& H^2(Y;\Z_2)}
\end{split}\EE
of based sets; see Proposition~\ref{SnakeLG_prp}.
The first three non-trivial maps in this sequence are group homomorphisms.\\

Suppose $n\!\ge\!2$.
By the naturality of~\eref{SpinLES_e}, there is a commutative diagram
$$ \xymatrix{ 
C\big(\SO(n);\Spin(n)\big) \ar[r]^{\wch{q}_0}\ar[d]^{f^*}&   
C\big(\SO(n);\SO(n)\big)\ar[r]^{\wch{\de}_0}\ar[d]^{f^*}& H^1(\SO(n);\Z_2) \ar[d]^{f^*}\\  
C\big(Y;\Spin(n)\big) \ar[r]^{\wch{q}_0}&  C\big(Y;\SO(n)\big)\ar[r]^{\wch{\de}_0}&
H^1(Y;\Z_2)}$$  
of group homomorphisms for every continuous map $f\!:Y\!\lra\!\SO(n)$.
Such a map lies in the image of $\wch{q}_0$ in~\eref{SpinLES_e} if and only if
$$f^*\eta_n=0\in H^1(Y;\Z_2),$$
where $\eta_n\!\in\!H^1(\SO(n);\Z_2)$ is the generator as before.
Thus,
\begin{gather*}
\wch\de_0\big(\id_{\SO(n)}\big)=\eta_n\in H^1\big(\SO(n);\Z_2\big), \quad
\wch\de_0(f)=\wch\de_0\big(f^*\id_{\SO(n)}\big)
=f^*\eta_n \in H^1(Y;\Z_2),\\
\Im\,\wch\de_0=\big\{f^*\eta_n\!:
f\!\in\!C^0\!\big(Y;\SO(n)\!\big)\big\}\subset H^1(Y;\Z_2) \,.
\end{gather*}
By Proposition~\ref{PinSpin_prp1}, $H^1(Y;\Z_2)$ acts freely and transitively on
the set~$\Sp(V,\fo)$ of $\Spin$-structures on an oriented vector bundle~$(V,\fo)$ 
over~$Y$ that admits a $\Spin$-structure.
By the exactness of~\eref{SpinLES_e}, 
$H^1(Y;\Z_2)/\Im\,\wch\de_0$ acts freely and transitively on the set of principal
$\Spin(n)$-bundles that doubly cover the trivial $\SO(n)$-bundle over~$Y$.
The latter action is induced by the former.
Thus, $\Im\,\wch\de_0$ in~~\eref{SpinLES_e} acts freely and transitively 
on the set of equivalence classes of double covers
$$q_V\!: \Spin(V,\fo_Y)\lra \SO(V,\fo_Y)\!=\!Y\!\times\!\SO(n)$$
as in Definition~\ref{PinSpin_dfn} with fixed domain and target.
The same considerations apply to $\Pin^{\pm}$-structures on a real vector bundle~$V$.\\

\noindent
If $n\!=\!2$ and $Y$ is a CW complex, then all groups in~\eref{SpinSOseq_e} are abelian
and $\Im\,\wch\de_0$ in~\eref{SpinLES_e} is the entire space $H^1(Y;\Z_2)$.
Thus, for every rank~2 oriented vector bundle~$(V,\fo)$ over~$Y$
there is at most one principal $\Spin(2)$-bundle $\Spin(V,\fo)$ that doubly covers 
the principal $\SO(2)$-bundle~$\SO(V,\fo)$ and 
$H^1(Y;\Z_2)$ acts freely and transitively on the set of equivalence classes
of associated projections~$q_V$ with the domain and target fixed.

\subsection{Correspondences and obstructions to existence}
\label{SpinPinprp_subs}

\noindent
We next establish the SpinPin~\ref{SpinPinStab_prop} property for the
$\Spin$- and $\Pin^{\pm}$-structures of Definition~\ref{PinSpin_dfn}.
Combining it with Examples~\ref{Pin1pStr_eg}-\ref{Spin2Str_eg},
we then obtain the SpinPin~\ref{SpinPinObs_prop} property;
it holds as long as $Y$ is a paracompact locally $H^2$-simple space
in the sense of Definition~\ref{locH1simp_dfn}.
We conclude this section by establishing the SpinPin~\ref{SpinPinCorr_prop} property.

\begin{proof}[{\bf{\emph{Proof of SpinPin~\ref{SpinPinStab_prop} property}}}]
Let $n\!=\!\rk\,V$, $W\!=\!\tau_Y\!\oplus\!V$, and
\BE{PinSpin3_e3}\begin{split}
\xymatrix{ \Pin^{\pm}(n) \ar[rr]^{\wt\io_{n+1;n}''} \ar[d]_{q_n^{\pm}} 
&& \Pin^{\pm}(n\!+\!1)\ar[d]^{q_{n+1}^{\pm}}  \\
\O(n) \ar[rr]^{\io_{n+1;n}''} && \O(n\!+\!1) }
\end{split}\EE
be as in Section~\ref{O2Pin_subs}.
We identify $\O(n)$ with a subspace of $\O(n\!+\!1)$ via the embedding~$\io_{n+1;n}''$
and $\Pin^{\pm}(n)$ with the subspace
$$\Pin^{\pm}(n\!+\!1)\big|_{\O(n)}\subset\Pin^{\pm}(n\!+\!1)$$
via the embedding~$\wt\io_{n+1;n}''$.
The natural embedding of $V$ into $W$
identifies $\O(V)$ with a subspace of $\O(W)$ preserved by the action of
$\O(n)\!\subset\!\O(n\!+\!1)$.\\

If $q_W\!:\Pin^{\pm}(W)\lra\O(W)$ 
is a $\Pin^{\pm}$-structure on~$W$, then the subspace
$$\Pin^{\pm}(V)\equiv \Pin^{\pm}(W)\big|_{\O(V)}\equiv q_W^{-1}\big(\O(V))$$
is preserved by the action of $\Pin^{\pm}(n)\!\subset\!\Pin^{\pm}(n\!+\!1)$.
Thus, the restriction
$$q_V\!: \Pin^{\pm}(V)\lra\O(V)$$
of $q_W$ determines a $\Pin^{\pm}$-structure on~$V$.\\

Let $\fp\!\equiv\!(\Pin^{\pm}(V),q_V)$ be a $\Pin^{\pm}$-structure on~$V$.
Transition data
$$\big\{\wt{g}_{\al\be}\!:U_{\al\be}\!\lra\!\Pin^{\pm}(n)\big\}_{\al,\be\in\cA}$$
for $\Pin^{\pm}(V)$ determines the transition data
\begin{gather*}
\big\{g_{\al\be}\!\equiv\!q_n^{\pm}(\wt{g}_{\al\be})\!:
U_{\al\be}\!\lra\!\O(n)\big\}_{\al,\be\in\cA}, \quad
\big\{\io_{n+1;n}''(g_{\al\be})\!:
U_{\al\be}\!\lra\!\O(n\!+\!1)\big\}_{\al,\be\in\cA},\\
\hbox{and}\qquad
\big\{\wt\io_{n+1;n}''(\wt{g}_{\al\be})\!:
U_{\al\be}\!\lra\!\Pin^{\pm}(n\!+\!1)\big\}_{\al,\be\in\cA}
\end{gather*}
for $\O(V)$, $\O(W)$, and a principal $\Pin^{\pm}(n\!+\!1)$-bundle $\Pin^{\pm}(W)$.
By the commutativity of~\eref{PinSpin3_e3}, there is a commutative diagram
$$\xymatrix{  \Pin^{\pm}(V) \ar[d]_{q_V}\ar[r] & \Pin^{\pm}(W)\ar[d]^{q_W}\\
\O(V) \ar[r]& \O(W)}$$
so that the top and bottom arrows are equivariant with respect to
the actions of $\Pin^{\pm}(n)$ and~$\O(n)$, respectively.
The right vertical arrow above determines a $\Pin^{\pm}$-structure $\St_V(\fp)$ on~$W$.\\

The two constructions above are mutual inverses and descend
to bijections between the collections $\cP^{\pm}(V)$ and~$\cP^{\pm}(W)$.
This establishes the existence of the second map in~\eref{SpinPinStab_e}.
The construction of the first map in~\eref{SpinPinStab_e} is analogous, 
with $\O$ and $\Pin^{\pm}$ replaced by $\SO$ and $\Spin$ throughout.
It is immediate that it satisfies the first condition after~\eref{SpinPinStab_e}.
The equivariance of the two bijections in~\eref{SpinPinStab_e} 
follows from~\eref{wtioSpin_e2}.\\

It remains to establish~\eref{SpinPinStab_e2}.
Suppose $\fo\!\in\!\fO(V)$, \hbox{$\fs\!\equiv\!(\Spin(V,\fo),q_V)$} is 
a Spin-structure on~$(V,\fo)$, and 
$$\St_V(\fs)\equiv\big(\Spin(W,\St_V(\fo)),q_W\big)$$ 
is the $\Spin$-structure on $(W,\St_V(\fo))$ constructed from a collection 
\BE{SpinPinStab_e15}
\big\{\wt{g}_{\al\be}\!:U_{\al\be}\!\lra\!\Spin(n)\big\}_{\al,\be\in\cA}\EE
of transition data for $\Spin(V,\fo)$ as above.
If $\fp$ is a $\Pin^{\pm}$-structure as before and $\fs$ is determined 
by~$\fp$ as in the proof of the SpinPin~\ref{Pin2SpinRed_prop} property, then
the $\Pin^{\pm}$-structure $\St_V(\fp)$ on~$W$ can  be constructed
using the same collection~$\{\wt{g}_{\al\be}\}_{\al,\be\in\cA}$.
The restriction of the projection
$$q_W\!:\Pin(W)\lra \O(W)$$
to the preimage of $\SO(W,\St_V(\fo))\!\subset\!\O(W)$ is 
then the projection for~$\St_V(\fs)$.
Thus, $\St_V(\fs)$ is the $\Spin$-structure on  $(W,\St_V(\fo))$ determined
by the $\Pin^{\pm}$-structure $\St_V(\fp)$ on~$W$ as 
 in the proof of the SpinPin~\ref{Pin2SpinRed_prop} property.
This establishes the last equality in~\eref{SpinPinStab_e2}.\\

Let $\{\wt{s}_{\al}\}_{\al\in\cA}$ be a collection of local sections of 
$\Spin(V,\fo)\!\subset\!\Spin(W,\St_V(\fo))$ 
determining a collection~\eref{SpinPinStab_e15} of transition data for $\Spin(V,\fo)$.
The collection 
\BE{SpinPinStab_e17}\big\{\wt\io_{n+1;n}''(\wt{g}_{\al\be})\!:
U_{\al\be}\!\lra\!\Spin(n\!+\!1)\big\}_{\al,\be\in\cA}\EE
then consists of the transition data for $\Spin(W,\St_V(\fo))$  
determined by $\{\wt{s}_{\al}\}_{\al\in\cA}$.
Let 
$$\ov{q}_V\!: \Spin(V,\ov\fo)\lra\SO(V,\ov\fo) ~~\hbox{and}~~
\ov{q}_W\!: 
\Spin\big(W,\ov{\St_V\fo}\big)\lra\SO\big(W,\ov{\St_V\fo}\big)\!=\!\SO\big(W,\St_V\ov\fo\big)$$
be the $\Spin$-structures determined by the collections~\eref{SpinPinStab_e15}
and~\eref{SpinPinStab_e17}
as in the proof of the SpinPin~\ref{OSpinRev_prop} property and
$$\ov{q}_W'\!: \Spin\big(W,\St_V\ov\fo\big)\lra\SO\big(W,\St_V\ov\fo\big)$$
be the $\Spin$-structure on $(W,\St_V\ov\fo)$ 
determined by the collection~\eref{conjSpin_e} as above.
The last two $\Spin$-structures are determined by the collections
\begin{equation*}\begin{split}
&\big\{\wt\bI_{n+1;1}^{-1}\wt\io_{n+1;n}''(\wt{g}_{\al\be})\wt\bI_{n+1;1}\!:
U_{\al\be}\!\lra\!\Spin(n\!+\!1)\big\}_{\al,\be\in\cA}
\qquad\hbox{and}\\
&\big\{\wt\io_{n+1;n}''(\wt\bI_{n;1}^{-1}\wt{g}_{\al\be}\wt\bI_{n;1})\!:
U_{\al\be}\!\lra\!\Spin(n\!+\!1)\big\}_{\al,\be\in\cA}
\end{split}\end{equation*}
of transition data.
By the first equation in~\eref{SO2Spin2_e9} with $(n,m)$ replaced by $(n\!+\!1,n)$, 
these two collections are the same.
This establishes the first equality in~\eref{SpinPinStab_e2}.
\end{proof}

\begin{proof}[{\bf{\emph{Proof of SpinPin~\ref{SpinPinObs_prop} property}}}]
For $n\!\in\!\Z^+$, let 
$$\ga_n\lra \G(n)  \qquad\hbox{and}\qquad \wt\ga_n\lra\wt\G(n)$$
be the real tautological $n$-plane bundle over the infinite Grassmannian of 
real $n$-planes and 
the oriented tautological $n$-plane bundle over the infinite Grassmannian of 
oriented $n$-planes.\\

\ref{PinObs_it} By Proposition~\ref{SnakeLG_prp}, the short exact sequence
$$\{\1\}\lra \Z_2\lra\Pin^{\pm}(n)\lra \O(n)\lra \{\1\}$$
of Lie groups induces a commutative diagram
\BE{PinSpin2_e5}\begin{split}
\xymatrix{ \wcH^1\big(\G(n);\Pin^{\pm}(n)\big)\ar[r]\ar[d]^{f^*}& 
\wcH^1\big(\G(n);\O(n)\big)\ar[r]^>>>>>>{\wch\de_1}\ar[d]^{f^*}& 
H^2\big(\G(n);\Z_2\big) \ar[d]^{f^*}\\
\wcH^1\big(Y;\Pin^{\pm}(n)\big)\ar[r]& 
\wcH^1\big(Y;\O(n)\big) \ar[r]^>>>>>>{\wch\de_1}& H^2(Y;\Z_2)}
\end{split}\EE
of pointed sets for every continuous map $f\!:Y\!\lra\!\G(n)$.
By the last part of Section~\ref{SpinPindfn_subs}
and the exactness of the bottom row in this diagram,
a rank~$n$ real vector bundle~$V$ over~$Y$ admits a $\Pin^{\pm}$-structure 
if and only if the image under~$\wch\de_1$ 
of the equivalence class~$[\O(V)]$ of $\O(V)$ in $\wcH^1(Y;\O(n))$
vanishes.\\

By \cite[Theorem~7.1]{MiSt}, $H^2(\G(n);\Z_2)$ is  generated by $w_2(\ga_n)$ 
and~$w_1^2(\ga_n)$. 
By \cite[Theorem~5.6]{MiSt}, for every rank~$n$ real vector bundle $V$ over~$Y$
there exists a continuous map $f\!:Y\!\lra\!\G(n)$ such that $V\!\approx\!f^*\ga_n$.
Along with the commutativity of~\eref{PinSpin2_e5}, these statements imply that there exist 
$a_n^{\pm},b_n^{\pm}\!\in\!\Z_2$ such~that
$$\wch\de_1\big([\O(V)]\big)=a_n^{\pm}w_2(V)+b_n^{\pm}w_1^2(V)$$
for every rank~$n$ real vector bundle~$V$ over every paracompact locally $H^2$-simple space~$Y$.
By the SpinPin~\ref{SpinPinStab_prop} property, $a_n^{\pm},b_n^{\pm}$ do not depend on~$n$;
we thus denote them by $a^{\pm},b^{\pm}$.
By Examples~\ref{Pin1pStr_eg} and~\ref{Pin1mStr_eg}, $b^+\!=\!0$ and $b^-\!=\!1$.
By Example~\ref{Spin2Str_eg} and the already established SpinPin~\ref{Pin2SpinRed_prop} property,  
$a^+\!=\!1$ and $a^-\!=\!1$.\\

\ref{SpinObs_it} By Proposition~\ref{SnakeLG_prp},
the short exact sequence~\eref{SpinSOseq_e} induces a commutative diagram
\BE{PinSpin2_e9}\begin{split}
\xymatrix{ \wcH^1\big(\wt\G(n);\Spin(n)\big)\ar[r]\ar[d]^{f^*}& 
\wcH^1\big(\wt\G(n);\SO(n)\big)\ar[r]^>>>>>>{\wch\de_1}\ar[d]^{f^*}& 
H^2\big(\wt\G(n);\Z_2\big) \ar[d]^{f^*}\\
\wcH^1\big(Y;\Spin(n)\big)\ar[r]& \wcH^1\big(Y;\SO(n)\big) 
\ar[r]^>>>>>>{\wch\de_1}& H^2(Y;\Z_2)}
\end{split}\EE
of pointed sets for every continuous map $f\!:Y\!\lra\!\wt\G(n)$.
By the last part of Section~\ref{SpinPindfn_subs}
and the exactness of the bottom row in this diagram,
a rank~$n$ oriented vector bundle~$V$ over~$Y$ admits a $\Spin$-structure if and only if
the image under~$\wch\de_1$ of the equivalence class~$[\SO(V)]$ of $\SO(V)$ in $\wcH^1(Y;\SO(n))$
vanishes.\\

By \cite[Theorem~12.4]{MiSt}, $H^2(\wt\G(n);\Z_2)$ is  generated by $w_2(\wt\ga_n)$. 
By \cite[Theorem~5.6]{MiSt}, for every rank~$n$ oriented vector bundle $V$ over~$Y$
there exists a continuous map $f\!:Y\!\lra\!\wt\G(n)$ such that $V\!\approx\!f^*\wt\ga_n$.
Along with the commutativity of~\eref{PinSpin2_e9}, these statements imply that there exists 
$a_n\!\in\!\Z_2$ such~that
$$\wch\de_1\big([\SO(V)]\big)=a_nw_2(V)$$
for every rank~$n$ oriented vector bundle~$V$ over every paracompact locally $H^2$-simple space~$Y$.
By the SpinPin~\ref{SpinPinStab_prop} property, $a_n$ does not depend on~$n$;
we thus denote it by~$a$.
By Example~\ref{Spin2Str_eg}, $a\!=\!1$.
\end{proof}

\begin{proof}[{\bf{\emph{Proof of SpinPin~\ref{SpinPinCorr_prop} property}}}] 
Let $V_{\pm}$ be as in~\eref{Vpmdfn_e}, $\fo_V^{\pm}$ be the canonical orientations
on these vector bundles, $n\!=\!\rk\,V$,  $n_{\pm}\!=\!\rk\,V_{\pm}$, and 
\BE{PinSpin4_e3}\begin{split}
\xymatrix{ \Pin^{\pm}(n) \ar[rr]^{\wt\io_n^{\pm}} \ar[d]_{q_n^{\pm}} &&
\Spin(n_{\pm})\ar[d]^{q_{n^{\pm}}} \\
\O(n)  \ar[rr]^{\io_n^{\pm}} && \SO(n_{\pm}) }
\end{split}\EE
be as in Corollary~\ref{Pin2Spin_crl}.
By the SpinPin~\ref{SpinPinObs_prop} property, 
$(V_{\pm},\fo_V^{\pm})$ admits a $\Spin$-structure if and only
if $V$ admits a $\Pin^{\pm}$-structure.\\

We can thus assume that there exists a $\Pin^{\pm}$-structure 
$\fp\!\equiv\!(\Pin^{\pm}(V),q_V)$ on~$V$.
Transition data
$$\big\{\wt{g}_{\al\be}\!:U_{\al\be}\!\lra\!\Pin^{\pm}(n)\big\}_{\al,\be\in\cA}$$
for $\Pin^{\pm}(V)$ determines the transition data
\begin{gather*}
\big\{g_{\al\be}\!\equiv\!q_n^{\pm}(\wt{g}_{\al\be})\!:
U_{\al\be}\!\lra\!\O(n)\big\}_{\al,\be\in\cA}, \quad
\big\{\io_n^{\pm}(g_{\al\be})\!:
U_{\al\be}\!\lra\!\SO(n_{\pm})\big\}_{\al,\be\in\cA},\\
\hbox{and}\qquad
\big\{\wt\io_n^{\pm}(\wt{g}_{\al\be})\!:
U_{\al\be}\!\lra\!\Spin(n_{\pm})\big\}_{\al,\be\in\cA}
\end{gather*}
for $\O(V)$, $\SO(V_{\pm})$, and a principal $\Spin(n_{\pm})$-bundle 
$$\Co_V^{\pm}(\fp)\equiv\big(\Spin(V_{\pm},\fo_{\pm}),q_{V_{\pm}}\big).$$
By the commutativity of~\eref{PinSpin4_e3}, there is a commutative diagram
$$\xymatrix{  \Pin^{\pm}(V) \ar[d]_{q_V}\ar[r] & \Spin(V_{\pm},\fo_{\pm})\ar[d]^{q_{V_{\pm}}}\\
\O(V) \ar[r]& \SO(V_{\pm},\fo_{\pm})}$$
so that the top and bottom arrows are equivariant with respect to
the actions of $\Pin^{\pm}(n)$ and~$\O(n)$, respectively.\\

The equivalence class of $\Spin$-structures on~$V_{\pm}$ represented by
the right vertical arrow in the above diagram is determined by
the equivalence class of $\Pin^{\pm}$-structures on~$V$ represented by
the left vertical arrow.
By the first column in~\eref{Pin2Spin_e},
the resulting maps~\eref{SpinPinCorr_e} are equivariant with respect to 
the actions of $H^1(Y;\Z_2)$ provided by Proposition~\ref{PinSpin_prp1}.
Since these actions are free and transitive, these maps are bijective.
By the above construction and the proof of the SpinPin~\ref{SpinPinStab_prop} property,
the $\Spin$-structures 
$\Co_{\tau_Y\oplus V}^{\pm}(\St_V^{\pm}(\fp))$ and 
$\St_{V_{\pm}}^{\pm}(\Co_V^{\pm}(\fp))$ on~\eref{SpinPinCorr_e0} 
are determined by the collections
\begin{equation*}\begin{split}
&\big\{\wt\io_{n+1}^{\pm}\big(\wt\io_{n+1;n}''(\wt{g}_{\al\be})\!\big)\!:
U_{\al\be}\!\lra\!\Spin(n_{\pm}\!+\!1)\big\}_{\al,\be\in\cA}
\qquad\hbox{and}\\
&\big\{\wt\io_{n_{\pm}+1;n_{\pm}}\big(\wt\io_n^{\pm}(\wt{g}_{\al\be})\!\big)\!:
U_{\al\be}\!\lra\!\Spin(n_{\pm}\!+\!1)\big\}_{\al,\be\in\cA},
\end{split}\end{equation*}
of transition data.
By~\eref{StCorrComm_e} with $(n,m)$ replaced by $(n\!+\!1,n)$, these two collections are the same.
This establishes the last claim of the SpinPin~\ref{SpinPinCorr_prop} 
property.
\end{proof}

\begin{eg}\label{Pin2Spin_eg}
Let $V$ be the infinite Mobius band line bundle of
Examples~\ref{Pin1pMB_eg} and~\ref{Pin1mMB_eg}.
Thus, 
\begin{alignat*}{3}
\SO(V_-,\fo_V^-)&\!=\!
\big([0,1]\!\times\!\SO(2)\big)\!/\!\!\sim\lra\R\P^1, &\quad
\big(0,\bI_{2;2}A\big)&\sim(1,A), &\quad [t,A]&\lra\big[\ne^{\pi\fI t}\big],\\
\SO(V_+,\fo_V^+)&\!=\!
\big([0,1]\!\times\!\SO(4)\big)\!/\!\!\sim\lra\R\P^1, &\quad
\big(0,\bI_{4;4}A\big)&\sim(1,A), &\quad [t,A]&\lra\big[\ne^{\pi\fI t}\big],
\end{alignat*}
with $\bI_{2;2}\!\equiv\!-\bI_2$ and $\bI_{4;4}\!\equiv\!-\bI_4$
as in~\eref{In24dfn_e}.
Under the conventions specified by~\eref{Pin1iden_e} and~\eref{Pin2Spindfn_e}, 
the $\Pin^{\pm}$-structure $\fp_0^{\pm}(V)$ on~$V$
constructed in Examples~\ref{Pin1pMB_eg} and~\ref{Pin1mMB_eg}
corresponds to the $\OSpin$-structure \hbox{$\os_0(V_{\pm},\fo_V^{\pm})$} given~by 
\begin{gather*}
\Spin_0\big(V_{\pm},\fo_V^{\pm}\big)=
\big([0,1]\!\times\!\Spin(3\!\pm\!1)\big)\!/\!\!\sim\,\lra\,\SO\big(V_{\pm},\fo_V^{\pm}\big),\\
\big(0,\wt\bI_{3\pm1;3\pm1}\wt{A}\big)\sim(1,\wt{A}), \qquad 
q_{(3\pm1)V}\big([t,\wt{A}]\big)=\big[t,q_{3\pm1}(\wt{A})\big],
\end{gather*}
on $(V_{\pm},\fo_V^{\pm})$,
with the homomorphism~$q_{3\pm1}$ given by~\eref{Spindfn_e}.
\end{eg}

\subsection{Short exact sequences}
\label{ShExSeq_subs}

We conclude the verification of the properties of Spin- and Pin-structures
in the perspective of Definition~\ref{PinSpin_dfn} by establishing 
the SpinPin~\ref{SpinPinSES_prop} property and 
using it in combination with Examples~\ref{negaautPin_eg} and~\ref{negaaut_eg} 
to obtain~\eref{SpinPinStr_e}.

\begin{proof}[{\bf{\emph{Proof of SpinPin~\ref{SpinPinSES_prop} property}}}] 
A short exact sequence~$\ce$ of vector bundles over a paracompact space~$Y$ 
as in~\eref{SpinPinSES_e0}  
determines a homotopy class of isomorphisms
$$V\approx V'\!\oplus\!V''$$
so that $\io$ is the inclusion as the first component
on the right-hand side above and $\fj$ is the projection to the second component.
Thus, it is sufficient to establish the SpinPin~\ref{SpinPinSES_prop} property
for the direct sum exact sequences as in~\eref{SpinPinDS_e0}.\\

From the definition of the induced orientation before
the statement of the SpinPin~\ref{SpinPinSES_prop} property,
it is immediate~that
\BE{DSCorr1_e3}\ov\fo'\fo''=\fo'\ov\fo''=
\ov{\fo'\fo''}\in\fO\big(V'\!\oplus\!V''\big)
\quad\hbox{and}\quad \fo_1'\big(\fo_2'\fo''\big)
=(\fo_1'\fo_2')\fo''
\in\fO\big(V_1'\!\oplus\!V_2'\oplus\!V\big)\EE
for all $\fo'\!\in\!\fO(V')$, $\fo_1'\!\in\!\fO(V_1')$, $\fo_2'\!\in\!\fO(V_2')$, 
$\fo''\!\in\!\fO(V'')$, and vector bundles $V',V_1',V_2',V$ over~$Y$.
For vector bundles~$V'$ and~$V''$ over~$Y$, we define
$$\io_{V',V''}'\!: V'\lra V'\!\oplus\!V'' \qquad\hbox{and}\qquad
\io_{V',V''}''\!: V''\lra V'\!\oplus\!V''$$
to be the canonical inclusions and  
$$\io_{V',V''}\!: \O(V')\!\times_Y\!\O(V'')\lra \O\big(V'\!\oplus\!V''\big)$$
to be the induced map.\\

Let $V'$ and $V''$ be vector bundles over~$Y$ of ranks~$m$ and $n\!-\!m$,
respectively, and \hbox{$V\!=\!V'\!\oplus\!V''$}.
Suppose 
$$\fo'\in\!\fO(V'), \quad \os'\equiv\big(\Spin(V',\fo'),q_{V'}\big), 
\quad\hbox{and}\quad  \fp''\equiv\big(\Pin^{\pm}(V''),q_{V''}\big)$$
are an orientation on~$V'$, an $\OSpin$-structure on~$V'$, and a $\Pin^{\pm}$-structure on~$V''$,
respectively.
Choose collections
\BE{DSCorr1_e7}\big\{\wt{s}_{\al}'\!:U_{\al}\!\lra\!\Spin(V',\fo')\big\}_{\al\in\cA}
\quad\hbox{and}\quad
\big\{\wt{s}_{\al}''\!:U_{\al}\!\lra\!\Pin^{\pm}(V'')\big\}_{\al\in\cA}\EE
of local sections so that  $\{U_{\al}\}_{\al\in\cA}$ is an open cover of~$Y$.
Let 
\BE{DSCorr1_e8}\big\{\wt{g}_{\al\be}'\!:U_{\al\be}\!\lra\!\Spin(m)\big\}_{\al,\be\in\cA}
\quad\hbox{and}\quad
\big\{\wt{g}_{\al\be}''\!:U_{\al\be}\!\lra\!\Pin^{\pm}(n\!-\!m)\big\}_{\al,\be\in\cA}\EE
be the associated transition data for 
the principal $\Spin(m)$-bundle $\Spin(V',\fo')$ and 
for the principal $\Pin^{\pm}(m)$-bundle $\Pin^{\pm}(V'')$, respectively.
The collections 
$$\big\{g_{\al\be}'\!=\!q_m(\wt{g}_{\al\be}')\!:
U_{\al\be}\!\lra\!\SO(m)\big\}_{\al,\be\in\cA} ~~\hbox{and}~~
\big\{g_{\al\be}''\!=\!q_{n-m}^{\pm}(\wt{g}_{\al\be}'')\!:
U_{\al\be}\!\lra\!\O(n\!-\!m)\big\}_{\al,\be\in\cA}$$
then consist of the transition data for the principal $\SO(m)$-bundle $\SO(V',\fo')$
and for the principal $\O(m\!-\!n)$-bundle $\O(V')$ induced by the collections
$$\big\{s_{\al}'\!\equiv\!q_{V'}\!\circ\!\wt{s}_{\al}'\!:U_{\al}\!\lra\!\SO(V',\fo')
\big\}_{\al\in\cA} \quad\hbox{and}\quad
\big\{s_{\al}''\!\equiv\!q_{V''}\!\circ\!\wt{s}_{\al}''\!:
U_{\al}\!\lra\!\O(V'')\big\}_{\al\in\cA}$$ 
of local trivializations of these bundles.\\

The transition data for the principal $\O(n)$-bundle determined by 
the collection
\BE{DSCorr1_e9}\big\{s_{\al}\!\equiv\!\io_{V',V''}(s_{\al}',s_{\al}'')\!:
U_{\al}\!\lra\!\O(V)\big\}_{\al\in\cA}\EE
of local sections is given by
$$\big\{g_{\al\be}\!\equiv\!\io_{n;m}(g_{\al\be}',g_{\al\be}'')\!:
U_{\al\be}\!\lra\!\O(n)\big\}_{\al,\be\in\cA}\,.$$
The collection 
\BE{DSCorr1_e11} \big\{\wt{g}_{\al\be}\!\equiv\!\wt\io_{n;m}(\wt{g}_{\al\be}',\wt{g}_{\al\be}'')\!:
U_{\al\be}\!\lra\!\Pin^{\pm}(n)\big\}_{\al,\be\in\cA}\EE
then determines a $\Pin^{\pm}$-structure~$\fp$ on~$V$,
\begin{gather*}
q_V\!:\Pin^{\pm}(V)\equiv\bigg(\bigsqcup_{\al\in\cA}\!U_{\al}\!\times\!\Pin^{\pm}(n)
\!\!\bigg)\!\!\Big/\!\!\!\sim \,\lra \O(V),
\quad q_V\big([y,\wt{A}]\big)=s_{\al}(y)\!\cdot\!q_n(\wt{A}),\\
U_{\al}\!\times\!\Pin^{\pm}(n)\!\ni\!\big(y,\wt{g}_{\al\be}(y)\wt{A}\big)\sim
\big(y,\wt{A}\big) \!\in\! U_{\be}\!\times\!\Pin^{\pm}(n)
\quad\forall\,\big(y,\wt{A}\big)\!\in\!U_{\al\be}\!\times\!\Pin^{\pm}(n),
~\al,\be\!\in\!\cA.
\end{gather*}
By the reasoning in Appendix~\ref{CechPB_subs}, 
the equivalence class of $\fp\!\equiv\!\llrr{\os',\fp''}_{\oplus}$ is independent 
of the choices of the collections~\eref{DSCorr1_e7} and depends only
on the equivalence classes of~$\os'$ and~$\fp''$.\\
 
Along with the definition of $\fo\!\equiv\!\fo'\fo''\!\in\!\fO(V)$,
the above construction with~$\fp''$ and $\Pin^{\pm}$ replaced by an $\OSpin$-structure 
\hbox{$\os''\!\equiv\!(\Spin(V'',\fo''),q_{V''})$} on~$V''$ and by $\Spin$, respectively, 
yields an $\OSpin$-structure 
$$\os\!\equiv\!\llrr{\os',\os''}_{\oplus}\equiv\big(\Spin(V,\fo),q_V\big)$$ 
on~$V$;
its equivalence class depends only on the equivalence classes of~$\os'$ and~$\os''$.
By~\eref{wtioSpin_e2} and~\eref{Pindfn_e2},
$$\wt\io_{n;m}\big(\wt{A}',\wt{A}''\wh\bI_{n-m}\big)
=\wt\io_{n;m}\big(\wt{A}'\wh\bI_m,\wt{A}''\big)
=\wt\io_{n;m}\big(\wt{A}',\wt{A}''\big)\wh\bI_n
\quad\forall\,(\wt{A}',\wt{A}'')\!\in\!\Spin(m)\!\times\!\Pin^{\pm}(n\!-\!m).$$
Along with the proof of Proposition~\ref{PinSpin_prp1},
this implies that the two induced maps~$\llrr{\cdot,\cdot}_{\ce}$ in~\eref{SpinPinSESdfn_e0}
are $H^1(Y;\Z_2)$-biequivariant.\\

Suppose the vector bundles $V'$ and $V''$ over~$Y$ are split as direct sums
of oriented line bundles, \hbox{$\os'\!=\!\os_0(V')$}, and $\os''\!=\!\os_0(V'')$.
We can take~$\cA$ to consist of a single element~$\al$
and $s_{\al}'$ and~$s_{\al}''$ to be sections of~$\SO(V',\fo')$ and~$\SO(V'',\fo'')$
respecting the ordered splittings and the orientations of all line bundle components.
The induced section~$s_{\al}$ of~$\SO(V,\fo'\fo'')$ in~\eref{DSCorr1_e9}  
then respects the induced ordered splitting and the orientations of all line bundle 
components of~$(V,\fo'\fo'')$.
This implies~\ref{DSsplit_it}.\\

The property~\eref{DSorient1_e} follows immediately from the construction.
If $\fp''$ is a $\Pin^{\pm}$-structure on~$V''$, $\fo''\!\in\!\fO(V'')$,
and $\os''$ is the $\OSpin$-Spin on~$V''$ obtained from~$\fp''$ and~$\fo''$
as in the proof of the SpinPin~\ref{Pin2SpinRed_prop} property, then
the collection $\{\wt{s}_{\al}''\}_{\al\in\cA}$ in~\eref{DSCorr1_e7} can be chosen
as in the construction for~$\os''$, i.e.~with values~in 
$$\Spin(V'',\fo'')\equiv\Pin^{\pm}(V'')\big|_{\SO(V'',\fo'')}\,.$$
The $\Pin^{\pm}$- and $\OSpin$-structures on~$V$ obtained in this way then satisfy
$$\OSpin(V)=\Pin^{\pm}(V)\big|_{\SO(V,\fo)}\,.$$
This establishes~\eref{DSorient2_e}.\\

Suppose $V_1',V_2',V''$ and $\os_1',\os_2',\fp''$ are as in~\ref{DSassoc_it}
and $m_1,m_2,n\!\in\!\Z^+$ are the ranks of $V_1'$, $V_2'$, and 
\hbox{$V_1'\!\oplus\!V_2'\!\oplus\!V''$}, respectively.
Let $\{s_{1;\al}'\}_{\al\in\cA}$, $\{s_{2;\al}'\}_{\al\in\cA}$, and
$\{s_{\al}''\}_{\al\in\cA}$ be collections of local sections 
associated with $\os_1',\os_2',\fp''$ as in~\eref{DSCorr1_e7} and 
$\{g_{1;\al\be}'\}_{\al,\be\in\cA}$, $\{g_{2;\al\be}'\}_{\al,\be\in\cA}$,
and $\{g_{\al\be}''\}_{\al,\be\in\cA}$ be the induced transition data.
The $\Pin^{\pm}$-structures $\llrr{\os_1',\llrr{\os_2',\fp''}_{\oplus}}_{\oplus}$
and $\llrr{\llrr{\os_1',\os_2'}_{\oplus},\fp''}_{\oplus}$ are 
then given by the transitions data
\begin{equation*}\begin{split}
&\big\{\wt\io_{n;m_1}\big(g_{1;\al\be}',\wt\io_{n-m_1;m_2}
(g_{2;\al\be}',g_{\al\be}'')\!\big)\!:U_{\al\be}\!\lra\!\Pin^{\pm}(n)\big\}_{\al,\be\in\cA}
\quad\hbox{and}\\
&\big\{\wt\io_{n;m_1+m_2}\big(\wt\io_{m_1+m_2;m_1}(g_{1;\al\be}',g_{2;\al\be}'),
g_{\al\be}''\big)\!:U_{\al\be}\!\lra\!\Pin^{\pm}(n)\big\}_{\al,\be\in\cA},
\end{split}\end{equation*}
respectively.
By~\eref{wtioassoci_e}, these two collections are the same.
This implies~\ref{DSassoc_it}.\\

With $\os'\!=\!\os_0(\tau_Y)$, the local sections~$\wt{s}_{\al}'$ in~\eref{DSCorr1_e7}
can be chosen so that $\wt{g}_{\al\be}'\!=\!\wt\bI_1$ for all $\al,\be\!\in\!\cA$.
If $\{\wt{g}_{\al\be}\}_{\al,\be}$ is the transition data for
a $\Pin^{\pm}$-structure $\fp$ determined by a collection $\{s_{\al}\}_{\al\in\cA}$ of local sections,
then the $\Pin^{\pm}$-structure $\llrr{\os_0(\tau_Y),\fp}_{\oplus}$ is determined 
by the transition data
$$\big\{\wt\io_{n+1;1}(\wt\bI_1,\wt{g}_{\al\be})\!:
U_{\al\be}\!\lra\!\Pin^{\pm}(n\!+\!1)\big\}_{\al,\be\in\cA}.$$
This collection is the same as the collection specifying the $\Pin^{\pm}$-structure $\St_V^{\pm}(\fp)$
in the proof of the SpinPin~\ref{SpinPinStab_prop} property.
This establishes~\ref{DSstab_it}.\\

Suppose $V$, $\fo$, and~$\fp$ are as in~\ref{DSPin2Spin_it}, \hbox{$n\!=\!\rk\,V$},
and $\os$ is the $\OSpin$-structure on~$V$ obtained from~$\fp$ and~$\fo$
as in the proof of the SpinPin~\ref{Pin2SpinRed_prop} property. 
We can choose sections~$\wt{s}_{\al}'$ for $(V,\fo)$ and $\wt{s}_{\al}''$
for $(2\!\pm\!1)\la(V)$ in~\eref{DSCorr1_e7} so~that
the associated transition functions $\wt{g}_{\al\be}''$ are constant and equal~$\wt\bI_{2\pm1}$.
Since
$$\wt\io_{n+1;n}\big(\wt{A},\wt\bI_1\big)
=\wt\io_{n+1;n}'\big(\wt{A}\big)=\wt\io_n^{\,-}(\wt{A})
\quad\hbox{and}\quad
\wt\io_{n+3;n}\big(\wt{A},\wt\bI_3\big)
=\wt\io_{n+3;n}'\big(\wt{A}\big)=\wt\io_n^{\,+}(\wt{A})$$
for all $\wt{A}\!\in\!\Spin(n)$,
the transition data for $\llrr{\os,\os_0((2\!\pm\!1)\la(V,\fo))}_{\oplus}$ 
in the above construction 
is the same as for $\Co_V^{\pm}(\fp)$ 
in the proof of the SpinPin~\ref{SpinPinCorr_prop} property.
This establishes~\ref{DSPin2Spin_it}.\\

It remains to establish~\ref{DSEquivSum_it}.
We continue with the setup and notation in the paragraph containing~\eref{DSCorr1_e7}
and in the paragraph immediately after it.
Define 
$$\big\{\eta_{\al\be}\!:U_{\al\be}\!\lra\!\Z_2\big\}_{\al,\be\in\cA}
\qquad\hbox{by}\quad \det g_{\al\be}''\!=\!(-1)^{\eta_{\al\be}}\,.$$
By Corollary~\ref{w1_crl}, the \v{C}ech cohomology class~$[\eta]$ 
represented by $\{\eta_{\al\be}\}$ is~$w_1(V'')$.
Let 
\BE{DSCorr1_e21}\ov\fs'\equiv\big(\Spin(V',\ov\fo'),\ov{q}_V'\big)\EE
be the $\Spin$-structure on $(V',\ov\fo')$ constructed from
the $\Spin$-structure~$\fs$ on  $(V',\fo')$ 
in the proof of the SpinPin~\ref{OSpinRev_prop} property.\\

\noindent
The analogues of the first collection in~\eref{DSCorr1_e8} and 
of the collection~\eref{DSCorr1_e11} for~$\ov\fs'$ are
\BE{DSCorr1_e23}\begin{split}
&\big\{\wh{g}_{\al\be}'\!\equiv\!\wt\bI_{m;1}^{-1}\wt{g}_{\al\be}'\wt\bI_{m;1}
\!:U_{\al\be}\!\lra\!\Spin(m)\big\}_{\al,\be\in\cA}
\qquad\hbox{and}\\
&\big\{\wh{g}_{\al\be}\!\equiv\!\wt\io_{n;m}(\wh{g}_{\al\be}',\wt{g}_{\al\be}'')\!:
U_{\al\be}\!\lra\!\Pin^{\pm}(n)\big\}_{\al,\be\in\cA}\,,
\end{split}\EE
respectively.
By~\eref{wtiocomm_e4},
\BE{DSCorr1_e25}  \wh{g}_{\al\be}=\wt\bI_{n;(m)}^{-1}
\big(\wt{g}_{\al\be}\wh\bI_n^{\eta_{\al\be}}\big)\wt\bI_{n;(m)}
\qquad\forall~\al,\be\!\in\!\cA.\EE
The $\OSpin$-structure~\eref{DSCorr1_e21} is described by
\begin{gather*}
\ov{q}_{V'}\!:
\Spin(V',\ov\fo')\equiv
\bigg(\bigsqcup_{\al\in\cA}\!U_{\al}\!\times\!\Spin(m)\!\!\bigg)\!\!\Big/\!\!\!\sim
\,\lra  \SO(V',\ov\fo'),
~~ \ov{q}_{V'}\big([y,\wt{A}]\big)=s_{\al}'(y)\!\cdot\!\bI_{m;1}q_m(\wt{A}),\\
U_{\al}\!\times\!\Spin(m)\!\ni\!\big(y,\wh{g}_{\al\be}'(y)\wt{A}\big)\sim
\big(y,\wt{A}\big) \!\in\! U_{\be}\!\times\!\Spin(m)
\quad\forall\,\big(y,\wt{A}\big)\!\in\!U_{\al\be}\!\times\!\Spin(m),
~\al,\be\!\in\!\cA.
\end{gather*}
The collections~\eref{DSCorr1_e23} arise from the local sections of $\Spin(V',\ov\fo')$ given~by
$$\wh{s}_{\al}'\!:U_{\al}\lra U_{\al}\!\times\!\Spin(m)\subset\Spin(V',\ov\fo'),
\qquad \wh{s}_{\al}'(y)=[y,\wt\bI_m]\,.$$
The analogue of the collection in~\eref{DSCorr1_e9} determined by these local sections is 
$$\big\{\ov{s}_{\al}\!\equiv\!\io_{V',V''}(s_{\al}'\!\cdot\!\bI_{m;1},s_{\al}'')
\!=\!s_{\al}\!\cdot\!\bI_{n;(m)}\!:
U_{\al}\!\lra\!\O(V)\big\}_{\al\in\cA}\,.$$
The $\Pin^{\pm}$-structure $\wh\fp\!\equiv\!\llrr{\ov\os',\fp''}_{\oplus}$ 
on~$V$ is then given~by 
\begin{gather*}
\wh{q}_V\!:\wh\Pin^{\pm}(V)\equiv\bigg(\bigsqcup_{\al\in\cA}\!U_{\al}\!\times\!\Pin^{\pm}(n)
\!\!\bigg)\!\!\Big/\!\!\!\sim \,\lra \O(V),
~~ \wh{q}_V\big([y,\wt{A}]\big)=s_{\al}(y)\!\cdot\!\bI_{n;(m)}q_n(\wt{A}),\\
U_{\al}\!\times\!\Pin^{\pm}(n)\!\ni\!\big(y,\wh{g}_{\al\be}(y)\wt{A}\big)\sim
\big(y,\wt{A}\big) \!\in\! U_{\be}\!\times\!\Pin^{\pm}(n)
\quad\forall\,\big(y,\wt{A}\big)\!\in\!U_{\al\be}\!\times\!\Pin^{\pm}(n),
~\al,\be\!\in\!\cA.
\end{gather*}

\vspace{.15in}

By the proof of Proposition~\ref{PinSpin_prp1}, the $\Pin^{\pm}$-structure
$$\eta\!\cdot\!\fp\equiv\big(\Pin_{\eta}^{\pm}(V),q_{\eta}\big)$$
on~$V$ is described by the displayed equation below~\eref{DSCorr1_e11} with~$\wt{g}_{\al\be}$
replaced by~$\wt{g}_{\al\be}\wh\bI_n^{\eta_{\al\be}}$.
By~\eref{DSCorr1_e25}, the~map
$$\wt\Psi\!:\Pin_{\eta}^{\pm}(V)\lra \wh\Pin^{\pm}(V),\quad
\wt\Psi\big([y,\wt{A}]\big) = \big[y,\wt\bI_{n;(m)}^{-1}\wt{A}\big]
\quad\forall\,\big(y,\wt{A}\big)\!\in\!U_{\al}\!\times\!\Pin^{\pm}(m),~\al\!\in\!\cA,$$
is thus well-defined. 
It is $\Pin^{\pm}(n)$-equivariant and satisfies $q_{\eta}\!=\!\wh{q}_V\!\circ\!\wt\Psi$.
We conclude that \hbox{$[\eta\!\cdot\!\fp]\!=\![\wh\fp]$} in~$\cP^{\pm}(V)$.
This establishes the first property in~\ref{DSEquivSum_it}.\\

Let $n\!=\!n\!+\!2\pm\!1$ as before.
The $\Spin$-structures $\Co_V^{\pm}\big(\llrr{\os',\fp''}_{\oplus})$ and 
$\llrr{\os',\Co_{V''}^{\pm}(\fp'')}_{\oplus}$ are specified by 
the transition data
\begin{equation*}\begin{split}
&\big\{\wt\io_n^{\pm}\big(\wt\io_{n;m}(g_{\al\be}',g_{\al\be}'')\!\big)\!:
U_{\al\be}\!\lra\!\Spin(n_{\pm})\big\}_{\al,\be\in\cA}
\quad\hbox{and}\\
&\big\{\wt\io_{n_{\pm};m}\big(g_{\al\be}',\wt\io_{n-m}^{\pm}(g_{\al\be}'')\big)\!:
U_{\al\be}\!\lra\!\Spin(n_{\pm})\big\}_{\al,\be\in\cA},
\end{split}\end{equation*}
respectively.
By~\eref{StCorrComm_e}, these two collections are the same.
This establishes the second property in~\ref{DSEquivSum_it}.
\end{proof}

\begin{proof}[{\bf{\emph{Completion of proof of SpinPin~\ref{SpinPinStr_prop} property}}}] 
In light of Proposition~\ref{PinSpin_prp1}, it remains to establish~\eref{SpinPinStr_e}.
It is immediate that this identity holds if the rank of~$V'$ is~0. 
By Example~\ref{negaautPin_eg}, 
it also holds if the rank of~$V''$ is~0.
We thus assume that the ranks of~$V'$ and~$V''$ are positive.\\

We define automorphisms of vector bundles $V',V'',V'\!\oplus\!V''$ by
\begin{alignat*}{3}
\Psi'\!: V'&\lra V', &\qquad \Psi''\!: V''&\lra V'', &\qquad
\wh\Psi\!: V'\!\oplus\!V''&\lra V'\!\oplus\!V'', \\
\Psi'(v')&= -v', &\qquad \Psi''(v'')&=-v'', &\qquad \wh\Psi(v',v'')&=(v',-v'').
\end{alignat*}
By Example~\ref{negaautPin_eg} and the SpinPin~\ref{Pin2SpinRed_prop} property, 
\BE{DSCorr1crl_e2} \Psi'^*\os'=\begin{cases}\os',&\hbox{if}~\rk\,V'\!\in\!2\Z;\\
\ov\os',&\hbox{if}~\rk\,V'\!\not\in\!2\Z;
\end{cases}\EE
these two statements were also established in Examples~\ref{negaaut_eg0} and~\ref{negaaut_eg}
directly.
By the naturality of the $\llrr{\cdot,\cdot}_{\oplus}$-operation of 
the SpinPin~\ref{SpinPinSES_prop} property,
\BE{DSCorr1crl_e3} \Psi^*\llrr{\os',\fp''}_{\oplus}=\llrr{\Psi'^*\os',\fp''}_{\oplus}, 
\qquad
\wh\Psi^*\llrr{\os',\fp''}_{\oplus}=\llrr{\os',\Psi''^*\fp''}_{\oplus}\EE
for all $\os'\!\in\!\OSpin(V')$ and $\fp''\!\in\!\OSpin(V'')$.\\

Suppose $V'$ is orientable and $V''$ admits a $\Pin^{\pm}$-structure~$\fp''$.
By the SpinPin~\ref{SpinPinObs_prop} property, $V'$ then admits an $\OSpin$-structure
if and only if $V'\!\oplus\!V''$ admits a $\Pin^{\pm}$-structure.
Along with the $H^1(Y;\Z_2)$-equivariance of $\llrr{\cdot,\cdot}_{\oplus}$
in the first input
and Proposition~\ref{PinSpin_prp1}, this implies that the~map
$$\OSp(V')\lra \cP^{\pm}(V'\!\oplus\!V''), \qquad \os'\lra \llrr{\os',\fp''}_{\oplus},$$
is a bijection.
It is thus sufficient to assume that $\fp\!=\!\llrr{\os',\fp''}_{\oplus}$ in this case.
By the first statement in~\eref{DSCorr1crl_e3}, \eref{DSCorr1crl_e2}, and 
the first statement in the SpinPin~\ref{SpinPinSES_prop}\ref{DSEquivSum_it} property,
$$\Psi^*\llrr{\os',\fp''}=\big(\!(\rk\,V')w_1(V'')\big)\!\cdot\!\llrr{\os',\fp''}\,.$$
Thus, \eref{SpinPinStr_e} holds if $V'$ is orientable and $V''$ admits 
a $\Pin^{\pm}$-structure.\\

\noindent
Suppose $V'$ is orientable and admits an $\OSpin$-structure~$\os'$.
By the SpinPin~\ref{SpinPinObs_prop} property, $V''$ then admits a $\Pin^{\pm}$-structure
if and only if $V'\!\oplus\!V''$ does.
Along with the $H^1(Y;\Z_2)$-equivariance of $\llrr{\cdot,\cdot}_{\oplus}$
in the second input and  Proposition~\ref{PinSpin_prp1},
this implies that the~map
$$\cP^{\pm}(V'')\lra \cP^{\pm}(V'\!\oplus\!V''), \qquad \fp''\lra \llrr{\os',\fp''}_{\oplus},$$
is an $H^1(Y;\Z_2)$-equivariant bijection.
Combining this with the second statement in~\eref{DSCorr1crl_e3}
and Example~\ref{negaautPin_eg}, we obtain
\begin{equation*}\begin{split}
\wh\Psi^*\llrr{\os',\fp''}_{\oplus}=\llrr{\os',\Psi''^*\fp''}_{\oplus}
&=\bllrr{\os',\big((\rk\,V''\!-\!1)w_1(V'')\big)\!\cdot\!\fp}_{\oplus}\\
&=\big((\rk\,V''\!-\!1)w_1(V'')\big)\!\cdot\!\bllrr{\os',\fp''}_{\oplus}
\end{split}\end{equation*}
for all $\fp''\!\in\!\cP^{\pm}(V'')$.
Thus,
\BE{DSCorr1crl_e5}\wh\Psi^*\fp=\big((\rk\,V''\!-\!1)w_1(V'')\big)\!\cdot\!\fp
\qquad\forall~\fp\!\in\!\cP^{\pm}(V'\!\oplus\!V'')\,.\EE
By symmetry, this implies~\eref{SpinPinStr_e} if 
$V''$ is orientable and admits an $\OSpin$-structure.\\

\noindent
By Proposition~\ref{PinSpin_prp1}, for every $\fp\!\in\!\cP^{\pm}(V'\!\oplus\!V'')$
there exists a unique $\eta(\fp)\!\in\!H^1(Y;\Z_2)$ such~that 
$$\Psi^*\fp= \eta(\fp)\!\cdot\!\fp\,.$$
For $\os\!\in\!\OSp(V'\!\oplus\!V'')$, let $\eta(\os)\!\equiv\!\eta(\fp)$
if $\os\!=\!\fR_{\fo}^{\pm}(\fp)$ with $\fp\!\in\!\cP^{\pm}(V)$;
by the SpinPin~\ref{Pin2SpinRed_prop} property,
this class is the same for the two possibilities for~$\fp$.
Since $\Psi$ is orientation-preserving on a fiber if and only~if $\rk\,V'$ is even,
$$\Psi^*\os=\eta(\os)\!\cdot\begin{cases}\os,&\hbox{if}~\rk\,V'\!\in\!2\Z;\\
\ov\os,&\hbox{if}~\rk\,V'\!\not\in\!2\Z.\end{cases}$$
Along with the $H^1(Y;\Z_2)$-equivariance of $\llrr{\cdot,\cdot}_{\oplus}$
in the first input and
the first statement in the SpinPin~\ref{SpinPinSES_prop}\ref{DSEquivSum_it} property,
this implies~that 
\BE{DSCorr1crl_e7} \bllrr{\Psi^*\os,\fp'''}_{\oplus}
=\big(\eta(\os)\!+\!(\rk\,V')w_1(V''')\big)\!\cdot\!\llrr{\os,\fp'''}_{\oplus}\EE
for every vector bundle $V'''$ over~$Y$ and every $\fp'''\!\in\!\cP^{\pm}(V''')$.\\

Suppose $V'$ admits a $\Pin^{\pm}$-structure~$\fp'$ and
$V'\!\oplus\!V''$ is orientable; in particular, \hbox{$w_1(V')\!=\!w_1(V'')$}.
By the SpinPin~\ref{SpinPinObs_prop} property and the first assumption, 
$V'\!\oplus\!V''$ admits an $\OSpin$-structure
if and only if $V'\!\oplus\!V''\!\oplus\!V'$ admits a $\Pin^{\pm}$-structure.
Along with the $H^1(Y;\Z_2)$-equivariance of $\llrr{\cdot,\cdot}_{\oplus}$
in the first input and  Proposition~\ref{PinSpin_prp1},
this implies that the~map
$$\OSp(V'\!\oplus\!V'')\lra \cP^{\pm}(V'\!\oplus\!V''\!\oplus\!V'), 
\qquad \os\lra \llrr{\os,\fp'}_{\oplus},$$
is an $H^1(Y;\Z_2)$-equivariant bijection.
By the orientability of $V''\!\oplus\!V'$ and the conclusion below~\eref{DSCorr1crl_e5}
with~$V''$ replaced by $V''\!\oplus\!V'$,
\begin{equation*}\begin{split}
\big\{\Psi\!\oplus\!\id_{V'}\big\}^{\,*}\llrr{\os,\fp'}_{\oplus}
&=\big((\rk\,V'\!-\!1)w_1(V')\!+\!(\rk\,V')w_1(V''\!\oplus\!V')\big)
\!\cdot\!\llrr{\os,\fp'}_{\oplus}\\
&=\big((\rk\,V'\!-\!1)w_1(V')\!\big)\!\cdot\!\bllrr{\os,\fp'}_{\oplus}
\end{split}\end{equation*}
for all $\os\!\in\!\OSp(V'\!\oplus\!V'')$.
By the naturality of the $\llrr{\cdot,\cdot}_{\oplus}$-operation and~\eref{DSCorr1crl_e7}, 
$$\big\{\Psi\!\oplus\!\id_{V'}\big\}^{\,*}\llrr{\os,\fp'}_{\oplus}
=\bllrr{\Psi^*\os,\fp'}_{\oplus}
=\big(\eta(\os)\!+\!(\rk\,V')w_1(V')\big)\!\cdot\!
\llrr{\os,\fp'}_{\oplus}\,.$$
Combining the last two equations, we obtain \hbox{$\eta(\os)\!=\!w_1(V')$}.
Thus, \eref{SpinPinStr_e} holds if $V'$ admits a $\Pin^{\pm}$-structure
and $V'\!\oplus\!V''$ is orientable.\\

We now consider the general case and show that 
\BE{DSCorr1crl_e11}\eta(\fp)=(\rk\,V'\!-\!1)w_1(V')\!+\!(\rk\,V')w_1(V'')
\qquad\forall~\fp\!\in\!\cP^{\pm}(V'\!\oplus\!V'').\EE
By the surjectivity of the Hurewicz homomorphism on~$\pi_1$ \cite[Proposition~7.5.2]{Spanier}
and the Universal Coefficient Theorem for Cohomology \cite[Theorem~53.5]{Mu2},
the homomorphism
\BE{kadfn_e}\ka\!:H^1(Y;\Z_2)\lra\Hom\big(\pi_1(Y),H^1(S^1;\Z_2)\big), \qquad
\big\{\ka(\eta)\big\}\big(\al\!:S^1\!\lra\!Y\big)=\al^*\eta,\EE
is injective.
It is thus sufficient to show that 
$$\al^*\eta(\fp)=(\rk\,V'\!-\!1)\al^*w_1(V')\!+\!(\rk\,V')\al^*w_1(V'')\in H^1(S^1;\Z_2)$$
for every $\al\!\in\!\cL(Y)$.
By the naturality of $\eta(\cdot)$ and $w_1(\cdot)$, this is equivalent~to
\BE{DSCorr1crl_e15}\eta(\al^*\fp)=(\rk\,V'\!-\!1)w_1(\al^*V')\!+\!(\rk\,V')w_1(\al^*V'')
\in H^1(S^1;\Z_2).\EE
By the SpinPin~\ref{SpinPinObs_prop} property, 
the bundles $\al^*V'$, $\al^*V''$, and $\al^*(V'\!\oplus\!V'')$
over~$S^1$ admit a $\Pin^{\pm}$-structure;
each of them admits an $\OSpin$-structure if orientable.
At least one of the three vector bundles is orientable.
Thus, \eref{DSCorr1crl_e15} holds because the three cases under which~\eref{SpinPinStr_e}
is shown to hold above cover all cases if $Y\!=\!S^1$.
\end{proof}

\section{Proof of Theorem~\ref{SpinStrEquiv_thm}\ref{PinSpinPrp_it}: trivializations perspectives}
\label{SpinPin_sec2}

We establish the statements of Section~\ref{SpinPinProp_subs} for the notions of 
$\Spin$-structure and $\Pin^{\pm}$-structure arising from Definitions~\ref{PinSpin_dfn2}
and~\ref{PinSpin_dfn3} in Sections~\ref{SpinPindfn2_subs} and~\ref{SpinPindfn3_subs}, respectively.
The proofs of these statements in the perspective of Definition~\ref{PinSpin_dfn3}
rely heavily on the simple topological observations of Section~\ref{TopolPrelim_subs}.

\subsection{Topological preliminaries}
\label{TopolPrelim_subs}

Let $V$ be a rank~$n$ vector bundle over a topological space~$Y$.
For $\fo\!\in\!\fO(V)$,
we denote by $\Triv(V,\fo)$ the set of homotopy classes of trivializations of
the oriented vector bundle~$(V,\fo)$.
The next lemma is used to establish the last equality in~\eref{SpinPinStab_e2}
for vector bundles~$V$ of ranks~1,2.

\begin{lmm}\label{StabAgr_lmm}
Let $V$ be an even-rank vector bundle over a CW complex~$Y$ of dimension at most~2.
For every $\fo\!\in\!\fO(V)$, the~maps
\begin{gather}\notag
\Triv(\tau_Y\!\oplus\!V,\fo_Y\fo)\lra 
\Triv\big(\tau_Y\!\oplus\!V\!\oplus\!\tau_Y,\fo_Y\fo\fo_Y\big),
\quad [\Phi]\lra[\Phi_1],[\Phi_2],\\
\label{StabAgr_e}
\Phi_1(r_1,v,r_2)=\big(\Phi(r_1,v),r_2\big),~~
\Phi_2(r_1,v,r_2)=\big(r_1,\Phi(r_2,v)\big),
\end{gather}
are the same.
\end{lmm} 

\begin{proof} We can assume that $n\!\equiv\!\rk\,V\!\ge\!2$.
Let
$$\io_1,\io_2\!:\SO(n\!+\!1)\lra \SO(n\!+\!2)$$
be the homomorphisms induced by the inclusions
$$\R^{n+1}\lra\R^{n+2}, \qquad (r,v)\lra(r,v,0),(0,v,r)~~\forall\,(r,v)\!\in\!\R\!\times\!\R^n.$$
Since the Lie group homomorphisms~$\io_1$ and~$\io_2$ differ by an automorphism 
of~$\SO(n\!+\!2)$,
\BE{StabAgr_e3} \io_{1*}\!=\!\io_{2*}\!:\pi_1\big(\SO(n\!+\!1)\big)
\lra\pi_1\big(\SO(n\!+\!2)\big)\!=\!\Z_2\,.\EE
Given a trivialization $\Phi$ of $(\tau_Y\!\oplus\!V,\fo_Y\fo)$,
let $\Phi_1,\Phi_2$ be the trivializations of $(\tau_Y\!\oplus\!V\!\oplus\!\tau_Y,\fo_Y\fo\fo_Y)$
as in~\eref{StabAgr_e}
and \hbox{$h\!:Y\!\lra\!\SO(n\!+\!2)$} be the continuous map so that $\Phi_2\!=\!h\Phi_1$.\\

Since $V$ is orientable, we can assume that its restriction to the 1-skeleton~$Y_1$ of~$Y$ is trivial.
The restriction of~$\Phi$ to~$Y_1$ is then given~by
$$\Phi|_{Y_1}\!:(\tau_Y\!\oplus\!V)\big|_{Y_1}\lra Y_1\!\times\!\R^{n+1},\quad
\Phi|_{Y_1}(y,r,v)=\big(y,\phi(y)(r,v)\big)~~\forall\,(y,r,v)\!\in\!Y_1\!\times\!\R\!\times\!\R^n,$$
for some continuous map $\phi\!:Y_1\!\lra\!\SO(n\!+\!1)$.
Thus, 
$$\Phi_1|_{Y_1}(y,w)=\big(y,(\io_1(\phi(y)))w\big), ~
\Phi_2|_{Y_1}(y,w)=\big(y,(\io_2(\phi(y)))w\big) 
\quad\forall\,(y,w)\!\in\!(\tau_Y\!\oplus\!V\!\oplus\!\tau_Y)\big|_{Y_1},$$
with $\io_2\!\circ\!\phi\!=\!h\!\cdot\!(\io_1\!\circ\!\phi)$.
By~\eref{StabAgr_e3},  $h|_{Y_1}$ is homotopic to a constant map.
Since $\pi_2(\SO(n\!+\!2))$ is trivial, it follows that $h$ itself is homotopic to a constant map.
Thus, the trivializations~$\Phi_1$ and~$\Phi_2$ are homotopic.
\end{proof}

For $\fo,\fo'\!\in\!\fO(V)$, we define a bijection
$$\fR_{\fo',\fo}\!:\Triv(V,\fo)\lra \Triv(V,\fo')$$
by sending the homotopy class of a trivialization~$\Phi$ of $(V,\fo)$
to the homotopy class $\fR_{\fo',\fo}([\Phi])$ containing the trivialization~$h\Phi$
for some locally constant map $h\!:Y\!\lra\!\O(n)$.
The next statement is the key input used to establish~\eref{SpinPinStr_e} and 
the first statement of the SpinPin~\ref{SpinPinSES_prop}\ref{DSEquivSum_it} property
in the perspectives of 
Definitions~\ref{PinSpin_dfn2} and~\ref{PinSpin_dfn3}. 

\begin{lmm}\label{Psiact_lmm}
Let $V'$ and $V''$ be vector bundles over a topological space~$Y$, 
$\Psi$ be as in~\eref{SpinPinStr_e0}, $\fo\!\in\!\fO(V'\!\oplus\!V'')$,
and $\hb\!\in\!\Triv(V'\!\oplus\!V'',\fo)$.
\begin{enumerate}[label=(\alph*),leftmargin=*]

\item\label{PsiactFix_it}  If $Y$ is a CW complex of dimension at most~2, $V'$ is orientable,
and $\rk(V'\!\oplus\!V'')\!\ge\!3$, 
then  $\Psi^*\hb\!=\!\fR_{\Psi^*\fo,\fo}(\hb)$.

\item\label{PsiactRev_it}  If $Y$ is paracompact and locally path connected and $V'$ 
is not orientable, then \hbox{$\Psi^*\hb\!\neq\!\fR_{\Psi^*\fo,\fo}(\hb)$}.

\end{enumerate}
\end{lmm}

\begin{proof} The assumptions in both cases imply that $V'\!\oplus\!V''$ is orientable.
We denote by $m$ and~$n$ the ranks of~$V'$ and~$V'\!\oplus\!V''$, respectively.
Let
$$\Phi\!:V'\!\oplus\!V''\lra Y\!\times\!\R^n \qquad\hbox{and}\qquad h\!:Y\lra\O(n)$$
be a trivialization in~$\hb$ and the continuous map 
so that $\Psi^*\Phi\!=\!h\Phi$, respectively.\\

\ref{PsiactFix_it} Since $V'$ and $V''$ are orientable, 
we can assume that their restrictions to the 1-skeleton~$Y_1$ of~$Y$ are trivial.
The restriction of~$\Phi$ to~$Y_1$ is then given~by
$$\Phi|_{Y_1}\!:(V'\!\oplus\!V'')\big|_{Y_1}\lra Y_1,\qquad
\Phi|_{Y_1}(y,v)=\big(y,\phi(y)v\big),$$
for some continuous map $\phi\!:Y_1\!\lra\!\SO(n)$.
Thus, 
$$\Psi^*\Phi|_{Y_1}(y,v)=\big(y,A\big(A^{-1}\phi(y)A\big)v\!\big),$$
where $A\!\in\!\O(n)$ is the diagonal matrix with the first $m$ diagonal entries equal $-1$
and the remainder equal~$1$.
The~map 
$$A^{-1}\phi A\!:Y_1\lra \SO(n)$$
is homotopic to the identity if and only if $\phi$ is.
Since $\pi_1(\SO(n))\!=\!\Z_2$, it follows that the maps $\phi$ and $A^{-1}\phi A$ are homotopic.
Thus, $h|_{Y_1}$ is homotopic to a constant map.
Since $\pi_2(\SO(n))$ is trivial, it follows that $h$ itself is homotopic to a constant map.\\

\ref{PsiactRev_it} Let $\ga_m\!\lra\!\G(m)$ and $\wt\ga_m\!\lra\!\wt\G(m)$
be the tautological vector bundles over the Grassmannians of real $m$-planes
and oriented $m$-planes, respectively.
By \cite[Theorem~5.6]{MiSt}, $V'\!\approx\!f^*\wt\ga_m$ for some continuous map 
$f\!:Y\!\lra\!\G(m)$.
Since the~map
$$\pi_1\big(\G(m)\big)\lra H^1(S^1;\Z_2), \qquad 
\big[\al\!:S^1\!\lra\!\G(m)\big]\lra\al^*w_1(\ga_m),$$
is an isomorphism and $V'$ is not orientable, 
\cite[Lemma~79.1]{Mu} implies that there exists \hbox{$\al\!\in\!\cL(Y)$}
such that $\al^*w_1(V')\!\neq\!0$ (otherwise, $f$ would lift over the projection
\hbox{$\wt\G(m)\!\lra\!\G(m)$} and $V'$ would be a pullback of
the oriented vector bundle~$\wt\ga_m$).
Thus, it is sufficient to establish the claim for $Y\!=\!\R\P^1$.\\

Let $\ga_{\R;1}\!\lra\!\R\P^1$ be the tautological line bundle.
By the assumptions on~$V'$ and~$V''$, 
$$V'\approx \ga_{\R;1}\!\oplus\!(m\!-\!1)\tau_{\R\P^1}, \qquad
V''\approx \ga_{\R;1}\!\oplus\!(n\!-\!m\!-\!1)\tau_{\R\P^1}\,.$$
By~\ref{PsiactFix_it}, the negation of each $\tau_{\R\P^1}$ component sends 
$\hb$ to~$\fR_{\ov\fo,\fo}(\hb)$.
Thus, we can assume that $m\!=\!1$.
Under the standard identification of $\R^2$ with~$\C$, two trivializations of~$2\ga_{\R;1}$
over~$\R\P^1$ are given~by
\begin{alignat*}{2}
\Phi_+\big([\ne^{\fI\th}],a\ne^{\fI\th},b \ne^{\fI\th}\big)
&= \big([\ne^{\fI\th}],\ne^{\fI\th}(a\!+\!\fI b)\big) &\qquad&\forall~a,b\!\in\!\R,\\
\Phi_-\big([\ne^{\fI\th}],a \ne^{\fI\th},b \ne^{\fI\th}\big)
&=\big([\ne^{\fI\th}],\ne^{-\fI\th}(a\!+\!\fI b)\big) &\qquad&\forall~a,b\!\in\!\R.
\end{alignat*}
They extend to trivializations of $V'\!\oplus\!V''$ by the identity on 
the $\tau_{\R\P^1}$ components and thus determine elements 
$\hb_-,\hb_+$ of $\Triv(V'\!\oplus\!V'',\fo)$,
for one of the two orientations~$\fo$ on~$V'\!\oplus\!V''$.\\

Since the conjugation on~$\C$ identifies the pullback of~$\Phi_+$ by the negation
of the second copy of~$\ga_{\R;1}$ with~$\Phi_-$,
\BE{Psiact_e7}\Psi^*\hb_{\pm}=\fR_{\ov\fo,\fo}\big(\hb_{\mp}\big)\,.\EE
On the other hand, the trivializations~$\Phi_+$ and~$\Phi_-$ differ by the multiplication
by the map
$$\R\P^1\lra S^1=\SO(2), \qquad [\ne^{\fI\th}]\lra \ne^{2\fI\th}\,.$$
Since this map generates $\pi_1(\SO(2))$ and the homomorphism 
\hbox{$\pi_1(\SO(2))\!\lra\!\pi_1(\SO(n))$} 
induced by the inclusion $\R^2\!\lra\!\R^n$ 
is surjective, it follows that $\hb_+\!\neq\!\hb_-$.
Combining this with~\eref{Psiact_e7}, we obtain the claim.
\end{proof}

\noindent
For a bordered surface $\Si$, we denote by
$$[\Si]_{\Z_2}\in H_2(\Si,\prt\Si;\Z_2)$$ 
the $\Z_2$-fundamental homology class of~$(\Si,\prt\Si)$.
The first statement and a special case of the second statement of Lemma~\ref{Z2cycles_lmm}
below  are used to establish the SpinPin~\ref{SpinPinObs_prop}
and~\ref{SpinPinStr_prop} properties in the perspective of Definition~\ref{PinSpin_dfn3}. 
The general case of the second statement and the third statement 
are used for similar purposes in Section~\ref{RelSpinPin_sec}.

\begin{lmm}\label{Z2cycles_lmm}
Let $X$ be a topological space and $Y\!\subset\!X$.
\begin{enumerate}[label=(\alph*),leftmargin=*]

\item\label{bound_it}
For every collection $\al_1,\ldots,\al_k\!\in\!\cL(Y)$
such~that $\al_1\!+\!\ldots\!+\!\al_k$ is a singular \hbox{$\Z_2$-boundary} in~$Y$,
there exists a continuous map $F\!:\!\Si\!\lra\!Y$ from a bordered surface
so that $F|_{\prt\Si}$ is the disjoint union of the loops $\al_1,\ldots,\al_k$.
If $Y$ is path-connected, then $\Si$ can be chosen to be connected.

\item\label{cycles_it} For every $b\!\in\!H_2(X,Y;\Z_2)$,
there exists a continuous map $F\!:(\Si,\prt\Si)\!\lra\!(X,Y)$ from a bordered surface so~that
\BE{Z2cycles_e0a}b=F_*\big([\Si]_{\Z_2}\big)\in H_2(X,Y;\Z_2). \EE

\item\label{relbound_it} For every continuous map $u\!:(\Si,\prt\Si)\!\lra\!(X,Y)$ 
from a bordered surface such that $u_*[\Si]_{\Z_2}$ vanishes in $H_2(X,Y;\Z_2)$,
there exists a continuous map $u'\!:\Si'\!\lra\!Y$ from a bordered surface so that
$\prt\Si\!=\prt\Si'$, $u|_{\prt\Si}\!=\!u'|_{\prt\Si'}$, and the continuous~map
$$u\!\cup\!u'\!:\Si\!\cup\!\Si'\lra X$$
obtained by gluing $u$ and $u'$ along the boundaries of their domains
vanishes in $H_2(X;\Z_2)$.

\end{enumerate}
\end{lmm}

\begin{proof} As in \cite[Section~2.1]{pseudo}, 
we denote by $\De^{\!k}\!\subset\!\R^k$ the standard $k$-simplex and
by $\De_j^{\!k}$ with \hbox{$j\!=\!0,1,\ldots,k$} its facets.
A singular $\Z_2$-chain~$b$ in~$X$ is a sum of finitely singular 2-simplicies,
i.e.~continuous maps  $F_i\!:\De^{\!2}\!\lra\!X$.\\

\ref{bound_it} 
Suppose $\al_1,\ldots,\al_k\!\in\!\cL(Y)$ and $b$ is a 
singular $\Z_2$-chain in~$Y$ as above so~that 
\BE{Spin2a_e5}\prt b=\al_1\!+\!\ldots\!+\!\al_k\in S_1(Y;\Z_2).\EE
We can then identify $k$ of the restrictions $F_i|_{\De^2_j}$ with the maps~$\al_i$
and pair up the remaining restrictions as in \cite[Lemma~2.11]{pseudo}.
By \cite[Theorem~77.5]{Mu}, the topological space~$\Si$ obtained by identifying 
the domains of~$F_i$ along the paired up edges is a closed surface with 
some open disks removed.
The closures of these disks may intersect, but only at the vertices of the glued 2-simplicies.
The maps~$F_i$ induce a continuous map  $F\!:\Si\!\lra\!Y$ so~that
its restriction to the boundary of~$\Si$ is the sum in~\eref{Spin2a_e5}.\\

\begin{figure}
\begin{pspicture}(-.2,-1.8)(10,1.5)
\psset{unit=.4cm}
\psline(2,-3)(8,-3)\psline(2,-3)(5,2.2)\psline(5,2.2)(8,-3)
\psline[arrowsize=6pt]{->}(2,-3)(5.5,-3)\psline[arrowsize=6pt]{->}(8,-3)(6.5,-.4)
\pscircle*(2,-3){.2}\pscircle*(8,-3){.2}\pscircle*(5,2.2){.2}
\rput(2.8,-.3){$\al_1$}\rput(5,-3.8){$a$}\rput(7.2,-.3){$c$}
\rput(1.4,-3.6){$x_1$}\rput(5.1,2.9){$x_1$}
\psline(9,-3)(15,-3)\psline(9,-3)(12,2.2)\psline(12,2.2)(15,-3)
\psline[arrowsize=6pt]{->}(15,-3)(11.5,-3)\psline[arrowsize=6pt]{->}(9,-3)(10.5,-.4)
\pscircle*(9,-3){.2}\pscircle*(15,-3){.2}\pscircle*(12,2.2){.2}
\rput(14.2,-.3){$\al_2$}\rput(12,-3.8){$a$}\rput(9.8,-.3){$c$}
\rput(15.6,-3.6){$x_1$}\rput(12.1,2.9){$x_1$}
\psline(19,-2)(22,3.2)\psline(22,3.2)(22,-.27)\psline(19,-2)(22,-.27)
\psline[arrowsize=6pt]{->}(19,-2)(20.5,-1.13)
\psline[arrowsize=6pt]{->}(22,-.27)(22,1.54)
\pscircle*(19,-2){.2}\pscircle*(22,3.2){.2}\pscircle*(22,-.27){.2}
\rput(19.8,.7){$\al_1$}\rput(18.4,-2.6){$x_1$}\rput(22,3.9){$x_1$}
\psline(24,-3)(30,-3)\psline(24,-3)(27,-1.27)\psline(30,-3)(27,-1.27)
\pscircle*(24,-3){.2}\pscircle*(30,-3){.2}\pscircle*(27,-1.27){.2}
\rput(23.9,-3.7){$x_1$}
\psline[arrowsize=6pt]{->}(24,-3)(27.5,-3)\rput(27,-3.8){$a$}
\psline[arrowsize=6pt]{->}(24,-3)(25.5,-2.13)
\psline[arrowsize=6pt]{->}(30,-3)(28.5,-2.13)
\psline(28,3.2)(31,-2)\psline(28,3.2)(28,-.27)\psline(28,-.27)(31,-2)
\psline[arrowsize=6pt]{->}(31,-2)(29.5,-1.13)
\psline[arrowsize=6pt]{->}(28,-.27)(28,1.54)
\rput(28.2,3.9){$x_1$}
\pscircle*(28,3.2){.2}\pscircle*(31,-2){.2}\pscircle*(28,-.27){.2}
\psline[arrowsize=6pt]{->}(31,-2)(29.5,.6)\rput(30.2,.7){$c$}
\psline(23,3.2)(27,3.2)\psline(23,3.2)(25,-.27)\psline(27,3.2)(25,-.27)
\psline[arrowsize=6pt]{->}(23,3.2)(25.3,3.2)\rput(25.2,3.9){$x_1$}
\psline[arrowsize=6pt]{->}(25,-.27)(24,1.47)
\psline[arrowsize=6pt]{->}(25,-.27)(26,1.47)
\pscircle*(23,3.2){.2}\pscircle*(27,3.2){.2}\pscircle*(25,-.27){.2}
\psline(20,-2.2)(23,-2.8)\psline(20,-2.2)(24.3,-1)\psline(23,-2.8)(24.3,-1)
\psline[arrowsize=6pt]{->}(20,-2.2)(22,-2.6)\rput(21.2,-3.2){$x_1$}
\psline[arrowsize=6pt]{->}(20,-2.2)(22.15,-1.6)
\psline[arrowsize=6pt]{->}(23,-2.8)(23.75,-1.76)
\pscircle*(20,-2.2){.2}\pscircle*(23,-2.8){.2}\pscircle*(24.3,-1){.2}
\psline(32,-3)(38,-3)\psline(32,-3)(35,2.2)\psline(35,2.2)(38,-3)
\psline[arrowsize=6pt]{->}(32,-3)(33.5,-.4)\rput(32.8,-.3){$c$}
\pscircle*(32,-3){.2}\pscircle*(38,-3){.2}\pscircle*(35,2.2){.2}
\rput(37.2,-.3){$\al_2$}\rput(35,-3.8){$a$}
\psline[arrowsize=6pt]{->}(38,-3)(34.5,-3)
\rput(38.6,-3.6){$x_1$}\rput(35.1,2.9){$x_1$}
\end{pspicture}
\caption{Modification of a bounding chain for $\al_1$ and $\al_2$ to produce 
a bounding surface}
\label{simplices_fig}
\end{figure}

The vertices of the restrictions of $F_i|_{\De^2_j}$ corresponding
to some of the maps~$\al_i$ may be identified with the vertices of other such restrictions,
instead of or in addition to being identified with just each other;
see the left side of Figure~\ref{simplices_fig}.
If this happens with one or both of the vertices of the restriction  $F_i|_{\De^2_j}$
corresponding to one of the loops~$\al_i$, we subdivide $F_i$ into 3 singular simplices and insert
two additional singular 2-simplicies each of which takes constant value $\al_i(0)\!=\!\al_i(1)$
along one of the edges; see the right side of Figure~\ref{simplices_fig}.
We then pair up the two constant edges to separate the edge corresponding to~$\al_i$
into its own boundary component, adding an extra handle in the process.
After separating off the edges corresponding to all loops~$\al_i$ into their own boundary components
in this way, we obtain a bordered surface and a continuous map as in~\ref{bound_it}.
If $Y$ is path-connected, we can connect the topological components of~$\Si$ by thin cylinders
and map them to appropriate paths in~$Y$.\\

\ref{cycles_it}
If $b$ above is a cycle in $(X,Y)$, then
the restrictions $F_i|_{\De^{\!2}_j}$ that are not mapped into~$Y$
come in pairs as in \cite[Lemma~2.10]{pseudo}.
The topological space~$\Si$ obtained by identifying 
the domains of~$F_i$ along the paired up edges is a closed surface with 
some open disks removed;
the edges forming its boundary~$\prt\Si$ are mapped into~$Y$.
After a modification as in~\ref{bound_it}, 
$\Si$ can be assumed to be a bordered surface.
The maps~$F_i$ induce a continuous map  \hbox{$F\!:(\Si,\prt\Si)\!\lra\!(X,Y)$} so~that 
$$[b]=F_*\big([\Si]_{\Z_2}\big)\in H_2(X,Y;\Z_2).$$ 
This establishes~\ref{cycles_it}.\\

\ref{relbound_it} Choose a triangulation of $\Si$ and a collection $F_i\!:\De^3\!\lra\!X$
of singular 3-simplicies that bound the 2-chain determined by the triangulation of~$\Si$
in the relative
singular group~$\cS_2(X,Y;\Z_2)$.
It can be assumed that no more than two singular 1-simplices obtained by restricting
the maps~$F_i$ to the edges of~$\De^3$ are the same (same 1-simplicies come precisely in pairs).
The restrictions $F_i|_{\De^{\!3}_j}$ that do not come pairs are either identified with 
the 2-simplices of the triangulation of~$\Si$ or are mapped to~$Y$.
The restrictions of such singular 2-simplices $F_i|_{\De^{\!3}_j}$ to the edges come in pairs.
The topological space~$\Si'$ obtained by identifying the domains of the 2-simplices 
$F_i|_{\De^{\!3}_j}$ that are neither paired nor associated with the 2-simplicies of~$\Si$
along the paired up edges is a closed surface with some open disks removed;
the edges forming its boundary~$\prt\Si'$ correspond to the edges forming 
the boundary~$\prt\Si$ of~$\Si$.
After a modification as in~\ref{bound_it}, $\Si'$ can be assumed to be a bordered surface
with $\prt\Si'\!=\!\prt\Si$ and the induced map $u\!:\Si'\!\lra\!Y$ satisfying 
$u|_{\prt\Si}\!=\!u'|_{\prt\Si'}$.
The glue map $u\!\cup\!u'$ bounds the collection of the singular 3-simplices~$F_i$
and thus vanishes in~$H_2(X;\Z_2)$.
\end{proof}

\subsection{The $\Spin$- and $\Pin$-structures of Definition~\ref{PinSpin_dfn3}}
\label{SpinPindfn3_subs}

\noindent
Throughout this section, the terms $\Spin$-structure and $\Pin^{\pm}$-structure
refer to the notions arising from Definition~\ref{PinSpin_dfn3}.
It is sufficient to establish the statements of Section~\ref{SpinPinProp_subs} 
under the assumption that the base topological space~$Y$ is path-connected;
this will be assumed to be the case.

\begin{proof}[{\bf{\emph{Proof of SpinPin~\ref{SpinPinStr_prop} property}}}] 
By the definition of $\Spin$- and $\Pin^{\pm}$-structures in this perspective,
it is sufficient to establish the claims of this property other than~\eref{SpinPinStr_e}  
for $\Spin$-structures \hbox{$\fs\!\equiv\!(\fs_{\al})_{\al}$} 
on oriented vector bundles~$(V,\fo)$ with 
\hbox{$n\!\equiv\!\rk\,V\!\ge\!3$}.\\

\noindent
By~\eref{pi012SO_e}, the oriented vector bundle $\al^*(V,\fo)$
over~$S^1$ has two homotopy classes of trivializations for every
\hbox{$\al\!\in\!\cL(Y)$}.
Let $\fs$ be a $\Spin$-structure on~$(V,\fo)$ and $\eta\!\in\!H^1(Y;\Z_2)$.
We define the $\Spin$-structure 
\hbox{$\eta\!\cdot\!\fs\!\equiv\!(\eta\!\cdot\!\fs_{\al})_{\al}$} 
on~$(V,\fo)$ by
$$\eta\!\cdot\!\fs_{\al}~\begin{cases}
=\fs_{\al},&\hbox{if}~\al^*\eta\!=\!0\!\in\!H^1(S^1;\Z_2);\\
\neq\fs_{\al},&\hbox{if}~\al^*\eta\!=\!0\!\in\!H^1(S^1;\Z_2);
\end{cases}
\qquad\forall~\al\!\in\!\cL(Y).$$
Suppose $F\!:\Si\!\lra\!Y$ is a continuous map from a connected bordered surface and 
\BE{SpinPinStr_e32}\al_1,\ldots,\al_m\!:S^1\lra Y\EE
are the restrictions of~$F$  to the boundary components of~$\Si$.
Since
$$\sum_{i=1}^m\blr{\al_i^*\eta,[S^1]}=
\blr{F^*\eta,\prt\Si}=0\in\Z_2,$$
the number of boundary components~$\al_i$ of~$F$ such that
$(\eta\!\cdot\!\fs)_{\al}\!\neq\!\fs_{\al}$ is even.
Since the trivializations in~$\fs$ extend to trivializations of $F^*(V,\fo)$, 
Corollary~\ref{X2VB_crl2b} then implies that this is also the case 
for the trivializations in~$\eta\!\cdot\!\fs$.
Thus, $\eta\!\cdot\!\fs$ is indeed a $\Spin$-structure on~$(V,\fo)$.
It is immediate that this construction defines a group action of $H^1(Y;\Z_2)$
on the set of such structures.
By the injectivity of the homomorphism~\eref{kadfn_e},
$\eta\cdot\fs_{\al}\!\neq\!\fs_{\al}$ for some
\hbox{$\al\!\in\!\cL(Y)$} if $\eta\!\neq\!0$ and 
so this action of $H^1(Y;\Z_2)$ is free.\\

Suppose $\fs'\!\equiv\!(\fs'_{\al})_{\al}$ is another $\Spin$-structure on~$(V,\fo)$.
Define
$$\eta\!:\cL(Y)\lra\Z_2, \qquad
\eta(\al)=\begin{cases}0,&\hbox{if}~\fs'_{\al}\!=\!\fs_{\al};\\
1,&\hbox{if}~\fs'_{\al}\!\neq\!\fs_{\al}.\end{cases}$$
This determines a linear map from the $\Z_2$-vector space generated by $\cL(Y)$
to~$\Z_2$.
Suppose \hbox{$F\!:\Si\!\lra\!Y$} is a continuous map from a bordered surface 
with boundary components~\eref{SpinPinStr_e32}.
By Definition~\ref{PinSpin_dfn3}, the homotopy classes of trivializations of $F^*(V,\fo)|_{\prt\Si}$ 
determined by~$\fs$ and $\fs'$ extend to a trivialization of~$F^*(V,\fo)$.
Corollary~\ref{X2VB_crl2b} then implies that 
the number of boundary components~$\al_i$ of~$F$ such that
$\fs'_{\al}\!\neq\!\fs_{\al}$ is even and~so
$$\sum_{i=1}^m\eta(\al_i)=0\in\Z_2\,.$$
Along with the surjectivity of the Hurewicz homomorphism for~$\pi_1$
and Lemma~\ref{Z2cycles_lmm}\ref{bound_it}, 
this implies that~$\eta$ descends to a homomorphism
$$\eta\!: H_1(Y;\Z_2)\lra \Z_2\,.$$
By the Universal Coefficient Theorem for Cohomology,
such a homomorphism corresponds to an element of $H^1(Y;\Z_2)$,
which we still denote by~$\eta$.
By the definition of~$\eta$ and the construction above, $\eta\!\cdot\!\fs\!=\!\fs'$.
Thus, the action of $H^1(Y;\Z_2)$ described above is transitive.\\

By the definition of the above action of $H^1(Y;\Z_2)$ on~$\Sp(V,\fo)$, 
\BE{SpinPinStr_e33}\al^*(\eta\!\cdot\!\fs)~
\begin{cases}=\al^*\fs,&\hbox{if}~\al^*\eta\!=\!0\!\in\!H^1(S^1;\Z_2);\\
\neq\al^*\fs,&\hbox{if}~\al^*\eta\!\neq\!0\!\in\!H^1(S^1;\Z_2);\end{cases}
\quad
\begin{aligned}
\forall~&\fs\!\in\!\Sp(V,\fo),\,\eta\!\in\!H^1(Y;\Z_2),\\
&\al\!\in\!\cL(Y).\end{aligned}\EE
Along with Lemma~\ref{Psiact_lmm}, this implies that
\BE{SpinPinStr_e35}
\Psi^*\fs=w_1(V')\cdot\begin{cases}\fs,&\hbox{if}~\rk\,V'\!\in\!2\Z;\\
\fR_{\ov\fo,\fo}(\fs),&\hbox{if}~\rk\,V'\!\not\in\!2\Z;
\end{cases}
\quad
\begin{aligned}
\forall~&\fo\!\in\!\fO(V'\!\oplus\!V''),\\
&\fs\!\in\!\Sp(V'\!\oplus\!V'',\fo),
\end{aligned}\EE
with $\Psi$ as in~\eref{SpinPinStr_e}.\\

\noindent
Let $V'$ and $V''$ be as in~\eref{SpinPinStr_e} and 
$$(W',W'')=\begin{cases}(V',V''\!\oplus\!(2\!\pm\!1)\la(V')\!\otimes\!\la(V''))
&\hbox{if}~\rk\,V'\!\in\!2\Z;\\
(V'\!\oplus\!(2\!\pm\!1)\la(V')\!\otimes\!\la(V''),V'')
&\hbox{if}~\rk\,V'\!\not\in\!2\Z.
\end{cases}$$
The automorphism $\Psi$ in~\eref{SpinPinStr_e0} induces the analogous automorphism~$\wt\Psi$
on~$W'\!\oplus\!W''$.
Since a $\Pin^{\pm}$-structure $\fp$ on $V'\!\oplus\!V''$ in this perspective
is a $\Spin$-structure $\fs^{\pm}$ on $(W'\!\oplus\!W'',\fo_{V'\oplus V''}^{\pm})$,
the first case in~\eref{SpinPinStr_e35} with~$\Psi$ replaced by~$\wt\Psi$ gives
$$\Psi^*\fp\equiv\wt\Psi^*\fs^{\pm}=w_1(W')\!\cdot\!\fs^{\pm}
\equiv \big(w_1(V')\!+\!(\rk\,V')(w_1(V')\!+\!w_1(V'')\!)\big)\!\cdot\!\fp\,.$$
This establishes~\eref{SpinPinStr_e}.
\end{proof}

\begin{proof}[{\bf{\emph{Proof of SpinPin~\ref{SpinPinObs_prop} property}}}] 
By the definition of $\Spin$- and $\Pin^{\pm}$-structures in this perspective,
it is sufficient to establish this claim for $\Spin$-structures on oriented vector bundles~$(V,\fo)$
with \hbox{$n\!\equiv\!\rk\,V\!\ge\!3$}.\\

Suppose $(V,\fo)$ admits a $\Spin$-structure $\fs\!\equiv\!(\fs_{\al})_{\al}$.
For every continuous map \hbox{$F\!:\Si\!\lra\!Y$} from a closed surface,
the bundle $F^*V$ is then trivializable and~so
$$\blr{w_2(V),F_*([\Si]_{\Z_2})}= \blr{w_2(F^*V),[\Si]_{\Z_2}}  =0. $$
Along with Lemma~\ref{Z2cycles_lmm}\ref{cycles_it} with $(X,Y)$ replaced by $(Y,\eset)$, 
this implies that $\lr{w_2(V),b}\!=\!0$ for every $b\!\in\!H_2(Y;\Z_2)$.
By the Universal Coefficient Theorem for Cohomology \cite[Theorems~53.5]{Mu2},
the homomorphism
$$\ka\!:H^2(Y;\Z_2)\lra \Hom_{\Z_2}\big(H_2(Y;\Z_2),\Z_2\big),\quad
\big\{\ka(\vp)\big\}(b)=\lr{\vp,b},$$
is an isomorphism.
Combining the last two statements, we conclude that $w_2(V)\!=\!0$.\\

Suppose $w_2(V)\!=\!0$.
Choose a collection $\cC\!\equiv\!\{\al_i\}_i$ of loops in~$Y$ that form
a basis for $H_1(Y;\Z_2)$ and a trivialization~$\phi_i$ of $\al_i^*(V,\fo)$
for each loop in~$\cC$.
Given \hbox{$\al\!\in\!\cL(Y)$}, let \hbox{$\al_1,\ldots,\al_k\!\in\!\cC$} 
be so~that 
$$\al\!+\!\al_1\!+\!\ldots\!+\!\al_k\in\prt S_2(Y;\Z_2).$$
By Lemma~\ref{Z2cycles_lmm}\ref{bound_it},
there exists a continuous map \hbox{$F\!:\Si\!\lra\!Y$} from a connected bordered surface 
which restricts~to
\BE{aldecomp_e}\al,\al_1,\ldots,\al_k\!:S^1\lra Y\EE
for some parametrization of $\prt\Si$ by $k\!+\!1$ copies of~$S^1$.
By Corollary~\ref{X2VB_crl2}, there exists a trivialization~$\Phi$ of $F^*V\!\lra\!\Si$
so that its restriction to each~$\al_i$ agrees with~$\phi_i$.
We take the homotopy class~$\fs_{\al}$ of the trivializations for $\al^*V\!\lra\!S^1$ to be
the homotopy class of the restriction of~$\Phi$.\\

Suppose $F'\!:\Si'\!\lra\!Y$ is another continuous map satisfying 
the conditions of the previous paragraph and $\Phi'$ is a trivialization 
of $F'^*V\!\lra\!\Si$ so that its restriction to each~$\al_i$ agrees with~$\phi_i$.
Denote by $\wh\Si$ (resp.~$\wt\Si$) the closed (resp.~bordered) surface 
obtained from~$\Si$ and~$\Si'$ by identifying
them along the boundary components corresponding to~$\al,\al_1,\ldots,\al_k$
(resp.~$\al_1,\ldots,\al_k$).
Thus, $\wt\Si$ has two boundary components, each of which corresponds to~$\al$,
and $\wh\Si$ is a connected surface obtained from $\wt\Si$ by identifying these two boundary components.
The~maps~$F$ and~$F'$ induce continuous~maps
$$\wh{F}\!:\wh\Si\lra Y \qquad\hbox{and}\qquad \wt{F}\!:\wt\Si\lra Y,$$
which restrict to $F$ and~$F'$ over $\Si,\Si'\!\!\subset\!\wh\Si,\wt\Si$.
The trivializations~$\Phi$ and~$\Phi'$ induce a trivialization~$\wt\Phi$ 
of~$\wt{F}^*(V,\fo)$ over~$\wt\Si$.
The bundle~$\wh{F}^*(V,\fo)$ over~$\wh\Si$ is obtained from $\wt\Si\!\times\!\R^n$
by identifying the copies of $\al^*V$ via the clutching map \hbox{$\vph\!:S^1\!\lra\!\SO(n)$}
determined by the difference between the trivializations of $\al^*V$ induced by~$\Phi$
and~$\Phi'$. 
Since 
$$w_2\big(\wh{F}^*V\big)=\wh{F}^*w_2(V)=0\in H^2(\wh\Si;\Z_2),$$
Corollary~\ref{X2VB_crl1b} implies that $\vph$ is homotopically trivial.
Thus,
$\Phi$ and $\Phi'$ determine the same homotopy class~$\fs_{\al}$ 
of trivializations of~$\al^*(V,\fo)$.\\

It remains to verify that $\fs\!\equiv\!(\fs_{\al})_{\al}$
satisfies the condition of Definition~\ref{PinSpin_dfn3}.
Suppose \hbox{$F\!:\Si\!\lra\!Y$} is a continuous map from a connected bordered surface and 
$$\al_1',\ldots,\al_m'\!:S^1\lra Y$$
are the restrictions of~$F$  to the boundary components of~$\Si$.
If $m\!=\!0$, the vector bundle $F^*W$ is trivializable by Corollary~\ref{X2VB_crl1}.
Below we consider the case $m\!\in\!\Z^+$.\\

For each $i\!=\!1,\ldots,m$, let 
$$\al_{ij}\!:S^1\!\lra\!Y~~\hbox{with}~j\!=\!1,\ldots,k_i,\qquad
F_i\!:\Si_i\lra Y, \qquad\hbox{and}\quad \Phi_i\!:F_i^*V\lra\Si_i\!\times\!\R^n$$
be the loops in~$\cC$,
a continuous map from a connected bordered surface, and a trivialization of $F_i^*(V,\fo)$, respectively,
as in the construction of the homotopy class~$\fs_{\al_i'}$ of trivializations
of $\al_i'^*(V,\fo)$ below~\eref{aldecomp_e}.
Thus, $F_i|_{\prt\Si_i}$ is the disjoint union of the loops 
$\al_i',\al_{i1},\ldots,\al_{ik_i}$ and
the restriction of~$\Phi_i$ to each~$\al_{ij}$ agrees with 
the initially chosen trivialization~$\phi_{ij}$.
Let $\phi'_i$ be the restriction of~$\Phi_i$ to 
$$\al_i'^*V\subset F_i^*V,F^*V\,.$$
By Corollary~\ref{X2VB_crl2}, the trivializations $\phi_1',\ldots,\phi_{m-1}'$
extend to a trivialization~$\Phi$ of~$F^*(V,\fo)$.
We show below that the restriction of~$\Phi$ to $\al_m'^*V$ is homotopic to~$\phi_m'$,
thus confirming the condition of Definition~\ref{PinSpin_dfn3}.\\

Since the loops in~$\cC$ are linearly independent in $H_1(Y;\Z_2)$ and
$$\sum_{i=1}^m\sum_{j=1}^{k_i}\al_{ij}=\sum_{i=1}^m\al_i'=\prt F
\in S_1(Y;\Z_2),$$
the loops $\al_{ij}$ come in pairs.
Denote by $\wh\Si$ (resp.~$\wt\Si$) the closed (resp.~bordered) surface 
obtained from $\Si_1,\ldots,\Si_m,\Si$ by identifying them along 
the paired up boundary components corresponding to~$\al_{ij}$ and 
along the boundary components corresponding to 
$\al_1',\ldots,\al_m'$ (resp.~$\al_1',\ldots,\al_{m-1}'$).
Thus, $\wt\Si$ has two boundary components, each of which corresponds to~$\al_m'$,
and $\wh\Si$ is a connected surface obtained from $\wt\Si$ by identifying these two boundary components.\\

\noindent
The~maps~$F_1,\ldots,F_m,$ and~$F$ induce continuous~maps
$$\wh{F}\!:\wh\Si\lra Y \qquad\hbox{and}\qquad \wt{F}\!:\wt\Si\lra Y,$$
which restrict to $F_i$ and~$F$ over $\Si_i,\Si\!\!\subset\!\wh\Si,\wt\Si$.
The trivializations~$\Phi_1,\ldots,\Phi_m$ and~$\Phi$ induce a trivialization~$\wt\Phi$ 
of~$\wt{F}^*V$ over~$\wt\Si$.
The bundle~$\wh{F}^*V$ over~$\wh\Si$ is obtained from $\wt\Si\!\times\!\R^n$
by identifying the copies of $\al_m'^*V$ via the clutching map \hbox{$\vph\!:S^1\!\lra\!\SO(n)$}
determined by the difference between the trivializations of $\al_m'^*V$ induced by~$\Phi_m$
and~$\Phi$. 
Since 
$$w_2\big(\wh{F}^*V\big)=\wh{F}^*w_2(V)=0\in H^2(\wh\Si;\Z_2),$$
Corollary~\ref{X2VB_crl1b} implies that $\vph$ is homotopically trivial.
Thus,
$\Phi_m$ and $\Phi$ determine the same homotopy class of trivializations of~$\al_m'^*V$.
\end{proof}

\begin{proof}[{\bf{\emph{Proof of SpinPin~\ref{OSpinRev_prop},\ref{Pin2SpinRed_prop} properties}}}]
For the purposes of establishing the first of these properties,
we can assume that \hbox{$n\!\equiv\!\rk\,V\!\ge\!3$}.
For $\fo\!\in\!\fO(V)$ and $\fs\!\in\!\Sp(V,\ov\fo)$,
we define the $\Spin$-structure $\ov\fs\!\in\!\Sp(V,\fo)$ by
\BE{OSpinRev_e33}\ov\fs_{\al}
=\big\{\ov\phi\!\equiv\!\bI_{n;1}\phi\!:\phi\!\in\!\fs_{\al}\big\}
\quad\forall\,\al\!\in\!\cL(Y),\EE
with $\bI_{n;1}\!\in\!\O(n)$ as below~\eref{Inmdfn_e}.
If \hbox{$F\!:\Si\!\lra\!Y$} is a continuous map from a bordered surface 
with boundary components~\eref{SpinPinStr_e32} and $\Phi$ is a trivialization
of $F^*(V,\fo)$ so that its restriction to the boundary component $(\prt\Si)_i$ of~$\Si$ 
corresponding to~$\al_i$ lies in~$\fs_{\al_i}$ for every~$i$, then
the restriction of the trivialization $\bI_{n;1}\Phi$ of $F^*(V,\ov\fo)$ to~$(\prt\Si)_i$
lies in~$\ov\fs_{\al_i}$.
Thus, $\ov\fs$ is indeed a $\Spin$-structure on~$(V,\ov\fo)$.
By~\eref{SpinPinStr_e33}, the resulting bijection~\eref{OSpinRev_e}
is $H^1(Y;\Z_2)$-equivariant.
It also satisfies the last condition of the SpinPin~\ref{OSpinRev_prop} property.\\

\noindent
We next describe a bijection between $\Sp(V,\fo)$ and 
$$\cP^{\pm}(V)\equiv \Sp(V\!\oplus\!(2\!\pm\!1)\la(V),\fo_V^{\pm}\big)
=\Sp(V\!\oplus\!(2\!\pm\!1)\tau_Y,\fo(2\!\pm\!1)\fo_Y\big);$$
the last equality holds because the orientation $\fo$ on~$V$ determines 
a canonical homotopy class of trivializations of $\al^*\la(V,\fo)$
for every loop~$\al$ in~$Y$.
By the SpinPin~\ref{SpinPinObs_prop} property,
$(V,\fo)$ admits a $\Spin$-structure if and only if $V$ admits a $\Pin^{\pm}$-structure.
We can thus assume that $(V,\fo)$ admits a $\Spin$-structure.\\

\noindent
If $\rk\,V\!\ge\!3$, we associate 
a $\Spin$-structure \hbox{$\fs\!\equiv\!(\fs_{\al})_{\al}$}
on~$(V,\fo)$ with the $\Spin$-structure \hbox{$\fs^{\pm}\!\equiv\!(\fs_{\al}^{\pm})_{\al}$}
on~$(V_{\pm},\fo_V^{\pm})$ given~by 
\BE{Pin2SpinRed_e33}
\fs_{\al}^{\pm}=\big\{\phi_{\pm}\!\equiv\!\phi\!\oplus\!(2\!\pm\!1)(\det\phi)\!:
\phi\!\in\!\fs_{\al}\big\} \quad\forall\,\al\!\in\!\cL(Y).\EE
If $F$ and $\Phi$ are as below~\eref{OSpinRev_e33}, then
the restriction of the trivialization 
$$\Phi^{\pm}\!\equiv\!\Phi\!\oplus\!(2\!\pm\!1)(\det\Phi)\!:
F^*(V_{\pm},\fo_V^{\pm})\lra \Si\!\times\!\R^{n+2\pm1}$$ 
to the boundary component of~$\Si$ corresponding to~$\al_i$ lies in~$\fs_{\al_i}^{\pm}$.
Thus, $\fs^{\pm}$ is indeed a $\Spin$-structure on~$(V_{\pm},\fo_V^{\pm})$.\\

\noindent
By~\eref{SpinPinStr_e33}, the map 
\BE{SpinPinStr_e36}\Sp(V,\fo)\lra\cP^{\pm}(V), \qquad \fs\lra\fs^{\pm},\EE
is $H^1(Y;\Z_2)$-equivariant.
Along with the SpinPin~\ref{SpinPinStr_prop} property, 
this implies that this map is a bijection.
We take the map~$\fR_{\fo}^{\pm}$ in~\eref{Pin2SpinRed_e} to be its inverse.
If~$\ov\phi$ is a trivialization of $\al^*(V,\ov\fo)$ as in~\eref{OSpinRev_e33},
then
$$\big(\ov\phi\big)_{\!\pm}=\bI_{n+2\pm1;3\pm1}\phi_{\pm}\!:
\al^*\big(V\!\oplus\!(2\!\pm\!1)\la(V)\big) \lra S^1\!\times\!\R^{n+2\pm1}.$$
Since this trivialization is homotopic to~$\phi_{\pm}$,
the last claim in the SpinPin~\ref{Pin2SpinRed_prop} property holds.\\

\noindent
If $\rk\,V\!=\!2$, we first identify the trivializations of $\al^*(\tau_Y\!\oplus\!V)$
with the trivializations of $\al^*(V\!\oplus\!\tau_Y)$ in the obvious way and
then extend them to trivializations of $\al^*(V\!\oplus\!3\tau_Y)$ by the identity
on the last two $\tau_Y$-summands.
If $\rk\,V\!=\!1$, we first identify the trivializations of 
$$\al^*(2\tau_Y\!\oplus\!V)\equiv
\al^*\big(\tau_Y\!\oplus\!(\tau_Y\!\oplus\!V)\big)$$ 
with the trivializations of $\al^*((\tau_Y\!\oplus\!V)\!\oplus\!\tau_Y)$
and of $\al^*(V\!\oplus\!2\tau_Y)$ 
in the obvious orientation-preserving ways and then extend the latter to the trivializations
of $\al^*(V\!\oplus\!3\tau_Y)$  by the identity on the last $\tau_Y$-summand.
These identifications of sections extend to pullbacks by continuous maps~$F$ 
as below~\eref{OSpinRev_e33}.
Thus, the collection $\fs^{\pm}\!\equiv\!(\fs^{\pm}_{\al})$ induced by~$\fs$
via these identifications is indeed a $\Spin$-structure on~$(V_{\pm},\fo_V^{\pm})$.
We note that 
\BE{Pin2SpinRed_e37} \fR_{\fo}^-\!=\!\fR_{\St_V(\fo)}^-\!:
\cP^-(V)\!=\!\cP^-(\tau_Y\!\oplus\!V)\lra\Sp(V,\fo)\!=\!\Sp\big(\St(V,\fo)\!\big)\EE
under the identifications~\eref{SpinPinCorr_e0} 
if $\rk\,V\!=\!1$ and $\fo\!\in\!\fO(V)$.
In all cases, the resulting bijections~\eref{Pin2SpinRed_e} are $H^1(Y;\Z_2)$-equivariant
and satisfy
the last requirement of the SpinPin~\ref{Pin2SpinRed_prop} property.
\end{proof}

\begin{proof}[{\bf{\emph{Proof of SpinPin~\ref{SpinPinStab_prop},\ref{SpinPinCorr_prop} properties}}}]
The two sides of~\eref{SpinPinCorr_e} are the same by definition;
we take this map to be the identity.
The requirement after~\eref{SpinPinCorr_e} is forced by the definition
of the second map of the SpinPin~\ref{SpinPinStab_prop} property in~\eref{SpinPinStab_e35} below.\\

Let $n\!=\!\rk\,V$.
By the SpinPin~\ref{SpinPinObs_prop} property, 
an oriented vector bundle $(V,\fo)$ admits a $\Spin$-structure if and only if 
 $(\tau_Y\!\oplus\!V,\St(\fo))$ does.
In order to establish the SpinPin~\ref{SpinPinStab_prop} property,
we can thus assume that there exists a $\Spin$-structure on~$(V,\fo)$.\\

If $n\!\ge\!3$, we associate 
a $\Spin$-structure \hbox{$\fs\!\equiv\!(\fs_{\al})_{\al}$}
on~$(V,\fo)$ with the $\Spin$-structure \hbox{$\St_V\fs\!\equiv\!(\St_V\fs_{\al})_{\al}$}
on $(\tau_Y\!\oplus\!V,\St(\fo))$ given~by 
\BE{SpinPinStab_e34}\St_V\fs_{\al}=\big\{\St_V\phi\!\equiv\id_{\tau_{S^1}}\!\oplus\!\phi\!:
\phi\!\in\!\fs_{\al}\big\} \quad\forall\,\al\!\in\!\cL(Y).\EE
If \hbox{$F\!:\Si\!\lra\!Y$} is a continuous map from a bordered surface 
with boundary components~\eref{SpinPinStr_e32} and $\Phi$ is a trivialization
of $F^*(V,\fo)$ so that its restriction to the boundary component~$(\prt\Si)_i$ of~$\Si$ 
corresponding to~$\al_i$ lies in~$\fs_{\al_i}$ for every~$i$, then
the restriction of the trivialization 
$$\St_V\Phi\!\equiv\id_{\tau_{\Si}}\!\oplus\!\Phi\!:
F^*\big(\tau_Y\!\oplus\!V,\St(\fo)\big)\lra \Si\!\times\!\R^{n+1}$$ 
to~$(\prt\Si)_i$ lies in~$\St_V\fs_{\al_i}$.
Thus, $\St_V\fs$ is indeed a $\Spin$-structure on~$(\tau_Y\!\oplus\!V,\St(\fo))$.
It is immediate that it satisfies the first condition after~\eref{SpinPinStab_e}.
With~$\ov\phi$ as in~\eref{OSpinRev_e33}, 
$$\St_V\ov\phi\!=\ov{\St_V\phi}\!: \al^*\big(\tau_Y\!\oplus\!V)\lra S^1\!\times\!\R^{n+1}\,.$$
Thus, the resulting first map in~\eref{SpinPinStab_e} satisfies the first equality 
in~\eref{SpinPinStab_e2}.\\

We take the second map in~\eref{SpinPinStab_e} to~be
\BE{SpinPinStab_e35}\begin{split}
\St_V^{\pm}\!\equiv\!\St_{V_{\pm}}\big|_{\Sp(V_{\pm},\fo_V^{\pm})}\!:\,
&\cP^{\pm}(V)\!\equiv\!\Sp(V_{\pm},\fo_V^{\pm})\\
&\lra 
\Sp\big(\St(V_{\pm},\fo_V^{\pm})\big)\!=\!
\Sp\big((\tau_Y\!\oplus\!V)_{\pm},\fo_{\tau_Y\oplus V}^{\pm}\big)
\!\equiv\!\cP^{\pm}\big(\tau_Y\!\oplus\!V\big);
\end{split}\EE
the identifications~\eref{laStiden_e} imply the equality above.
By~\eref{SpinPinStr_e33},
both maps~\eref{SpinPinStab_e} are then $H^1(Y;\Z_2)$-equivariant.
In light of the SpinPin~\ref{SpinPinStr_prop} property,
this implies that they are bijections.
With~$\phi_{\pm}$ as in~\eref{Pin2SpinRed_e33}, 
$$\big(\St_V\phi\big)_{\pm}\!=\St_V^{\pm}\phi_{\pm}\!: 
\al^*(\tau_Y\!\oplus\!V)_{\pm}\!=\!\al^*\big(\tau_Y\!\oplus\!V_{\pm}\big)
\lra S^1\!\times\!\R^{n+3\pm1}\,.$$
Thus, the second equality in~\eref{SpinPinStab_e2} also holds.\\

\noindent
For $n\!=\!1,2$, the two sides of the first map in~\eref{SpinPinStab_e} are 
the same by definition. 
We take this map to be the identity~then and define the second map in~\eref{SpinPinStab_e}
by~\eref{SpinPinStab_e35}.
The second equality in~\eref{SpinPinStab_e2} is  satisfied
by~\eref{Pin2SpinRed_e37} in the $\Pin^-$-case if $n\!=\!1$
and by Lemma~\ref{StabAgr_lmm} in the three remaining cases.
The first map in~\eref{SpinPinStab_e} satisfies the first condition after~\eref{SpinPinStab_e}
and the first equality in~\eref{SpinPinStab_e2} in both cases by definition.
The two maps are $H^1(Y;\Z_2)$-equivariant in all six cases,
by definition in three of the cases and by~\eref{SpinPinStr_e33} in the remaining
three cases.
\end{proof}

\begin{proof}[{\bf{\emph{Proof of SpinPin~\ref{SpinPinSES_prop} property}}}]
For every $\al\!\in\!\cL(Y)$, a short exact sequence~$\ce$ of vector bundles 
over~$Y$ as in~\eref{SpinPinSES_e0}  determines a homotopy class of isomorphisms
$$\al^*V\approx \al^*V'\!\oplus\!\al^*V''$$
so that $\al^*\io$ is the inclusion as the first component
on the right-hand side above and $\al^*\fj$ is the projection to the second component.
Thus, it is sufficient to establish the SpinPin~\ref{SpinPinSES_prop} property
for the direct sum exact sequences as in~\eref{SpinPinDS_e0}.
Furthermore, an orientation~$\fo'$ on~$V'$ determines a homotopy class of trivializations 
of~$\al^*\la(V')$ and thus of isomorphisms
\BE{SpinPinSES_e33}\begin{split} 
\al^*\big(V'\!\oplus\!V''\big)_{\pm}
&\equiv \al^*\big(V'\!\oplus\!V''\!\oplus\!(2\!\pm\!1)\la(V'\!\oplus\!V'')\big)\\
&\approx  \al^*\big(V'\!\oplus\!V''\!\oplus\!(2\!\pm\!1)\la(V'')\big)
\equiv \al^*\big(V'\!\oplus\!V''_{\pm}\big)
\end{split}\EE
for every $\al\!\in\!\cL(Y)$.\\

Let $V'$ and $V''$ be vector bundles over~$Y$ of rank~$n'$ and~$n''$, respectively.
Suppose $n',n''\!\ge\!3$, $\fo'\!\in\!\fO(V')$, and $\fo''\!\in\!\fO(V'')$.
For a $\Spin$-structure \hbox{$\fs'\!\equiv\!(\fs_{\al}')_{\al}$} on~$(V',\fo')$ 
and a $\Spin$-structure \hbox{$\fs''\!\equiv\!(\fs_{\al}'')_{\al}$} on~$(V'',\fo'')$,
we define a $\Spin$-structure 
$\llrr{\fs',\fs''}_{\oplus}\!\equiv\!({\fs_{\al}'\!\oplus\!\fs_{\al}''})_{\al}$
on $(V'\!\oplus\!V'',\fo'\fo'')$ by
\BE{dfn3sum_e}\fs_{\al}'\!\oplus\!\fs_{\al}''=
\big\{\phi'\!\oplus\!\phi''\!:\phi'\!\in\!\fs_{\al}',\,\phi''\!\in\!\fs_{\al}''\big\} 
\quad\forall\,\al\!\in\!\cL(Y).\EE
If \hbox{$F\!:\Si\!\lra\!Y$} is a continuous map from a bordered surface 
with boundary components~\eref{SpinPinStr_e32} and $\Phi',\Phi''$ are trivializations
of $F^*(V',\fo')$ and $F^*(V'',\fo'')$ so that their restrictions to 
the boundary component~$(\prt\Si)_i$ of~$\Si$  corresponding to~$\al_i$ 
lie in~$\fs_{\al_i}'$ and~$\fs_{\al_i}''$, respectively, for every~$i$, 
then the restriction of the trivialization 
$$\Phi'\!\oplus\!\Phi''\!:
F^*(V'\!\oplus\!V'',\fo'\fo'')\lra \Si\!\times\!\R^{n'+n''}$$ 
to~$(\prt\Si)_i$ lies in~$\fs_{\al_i}'\!\oplus\!\fs_{\al_i}''$.
Thus, $\llrr{\fs',\fs''}_{\oplus}$ is indeed a $\Spin$-structure 
on~$(V'\!\oplus\!V'',\fo'\fo'')$.\\

The above construction determines the first map~$\llrr{\cdot,\cdot}_{\ce}$ 
in~\eref{SpinPinSESdfn_e0} if \hbox{$n',n''\!\ge\!3$}.
Along with~\eref{SpinPinSES_e33}, it also determines the second map 
$\llrr{\cdot,\cdot}_{\ce}$
if $n'\!\ge\!3$ and \hbox{$\rk\,V''_{\pm}\!\ge\!3$} (i.e.~not $V_-''$ if $n''\!=\!1$).
By~\eref{SpinPinStr_e33}, both maps are $H^1(Y;\Z_2)$-biequivariant.
These two maps satisfy \ref{DSsplit_it},  \eref{DSorient1_e}, 
and~\ref{DSassoc_it} for vector bundles of ranks at least~3
by definition.
With~$\phi_{\pm}$ as in~\eref{Pin2SpinRed_e33},
$$\big(\phi'\!\oplus\!\phi''\big)_{\pm}\!=\!\phi'\!\oplus\!\phi''_{\pm}\!:
\al^*(V'\!\oplus\!V'')_{\pm}\!=\!
\big(\al^*V'\big)\!\oplus\!\big(\al^*V_{\pm}''\big)\lra S^1\!\times\!\R^{n'+n''+2\pm1}\,.$$
Thus, \eref{DSorient2_e} is also satisfied in these cases.\\

For a vector bundle $V$ over $Y$ and $a\!\in\!\Z^{\ge0}$, let
\BE{SpinPinSES_e34a}\St_V^a\!: \OSp(V)\lra \OSp\big(a\tau_Y\!\oplus\!V\big)
\quad\hbox{and}\quad
\St_V^a\!:\cP^{\pm}(V)\lra \cP^{\pm}\big(a\tau_Y\!\oplus\!V\big)\EE
denote $a$ iterations of the maps~\eref{SpinPinStab_e}.
By the proof of the SpinPin~\ref{SpinPinStab_prop} property and the construction above,
\BE{SpinPinSES_e35}
\St_V^a(\fp) =\llrr{\os_0(a\tau_Y),\fp}_{\oplus} \EE
 for all $a\!\ge\!3$ and $\fp\!\in\!\cP^{\pm}(V)$ if $\rk\,V_{\pm}\!\ge\!3$.\\

For vector bundles $V',V''$ over~$Y$ and $a,b\!\in\!\Z^{\ge0}$, let
$$\Psi_{V',V''}\!:a\tau_Y\!\oplus\!b\tau_Y\!\oplus\!V'\!\oplus\!V''\lra
a\tau_Y\!\oplus\!V'\!\oplus\!b\tau_Y\!\oplus\!V''$$
be the obvious bundle isomorphism.
If $b\!\in\!2\Z^{\ge0}$, then
\BE{SpinPinSES_e35a}
\St_{V'\oplus V''}^{a+b}(\fo'\fo'')=\Psi_{V',V''}^{\,*}
\big(\St_{V'}^a(\fo')\St_{V''}^b(\fo'')\!\big)
\quad\forall\,\fo'\!\in\!\fO(V'),\,\fo''\!\in\!\fO(V'').\EE
If $b\!\in\!2\Z^{\ge0}$ and $V',V''$ split as direct sums of oriented line bundles, then
\BE{SpinPinSES_e35b}
\St_{V'\oplus V''}^{a+b}\big(\os_0(V'\!\oplus\!V'')\big)
=\Psi_{V',V''}^{\,*}\os_0\big(\St_{V'}^a(\os_0(V')\!)\!\oplus\!\St_{V''}^b(\os_0(V'')\!)\big).\EE
If $b\!\in\!2\Z^{\ge0}$ and $\rk\,V',\rk\,V''\!\ge\!3$, then
\BE{SpinPinSES_e35c}\begin{split}
\St_{V'\oplus V''}^{a+b}\big(\llrr{\os',\os''}_{\oplus}\big)
&=\Psi_{V',V''}^{\,*}\bllrr{\St_V^a(\os'),\St_V^b(\os'')}_{\oplus}\,,\\
\St_{V'\oplus V''}^{a+b}\big(\llrr{\os',\fp''}_{\oplus}\big)
&=\Psi_{V',V''}^{\,*}\bllrr{\St_V^a(\os'),\St_V^b(\fp'')}_{\oplus}
\end{split}\EE
for all $\os'\!\in\!\OSp(V')$, $\os''\!\in\!\OSp(V'')$, and $\fp''\!\in\!\cP^{\pm}(V'')$;
this follows immediately from the construction of $\llrr{\cdot,\cdot}_{\oplus}$ above.
By the SpinPin~\ref{SpinPinCorr_prop} property, 
the proofs of the  SpinPin~\ref{SpinPinStab_prop} and~\ref{SpinPinCorr_prop} properties,
and the construction above, 
\BE{SpinPinSES_e35d}
\Co_{b\tau_Y\oplus V}^{\pm}\big(\St_V^b\fp\big)
=\St_V^b\big(\Co_V^{\pm}(\fp)\big)
=\Psi_{V,(2\pm1)\la(V,\fo)}^*
\bllrr{\fR_{\fo}^{\pm}(\fp),\St_V^b\big(\os_0((2\!\pm\!1)\la(V,\fo))\big)}_{\oplus}\EE
for all $b\!\in\!2\Z^{\ge0}$, $\fo\!\in\!\fO(V)$, and $\fp\!\in\!\cP^{\pm}(V)$ if $\rk\,V\!\ge\!3$.\\

For vector bundles $V',V''$ over~$Y$ not both of ranks at least~3,
we define the two maps in~\eref{SpinPinSESdfn_e0}
by~\eref{SpinPinSES_e35c} with $b\!\in\!2\Z^{\ge0}$.
By \eref{SpinPinSES_e35c} for vector bundles~$V',V''$ of ranks at least~3
and the SpinPin~\ref{SpinPinStab_prop} property,
the resulting maps are well-defined and $H^1(Y;\Z_2)$-biequivariant.
Furthermore, the definition of the second map agrees with the definition above
if the rank of~$V''_{\pm}$ is at least~3.
By~\ref{DSsplit_it} for~$V',V''$ of ranks at least~3, 
\eref{SpinPinSES_e35b}, and the first property after~\eref{SpinPinStab_e}, 
the first map satisfies~\ref{DSsplit_it} for all~$V',V''$.
By~\eref{SpinPinSES_e35a}, this map satisfies~\eref{DSorient1_e}.
By~\eref{DSorient2_e} for~$V',V''$ of ranks at least~3
and the second equality in~\eref{SpinPinStab_e2}, 
the two maps satisfy~\eref{DSorient2_e} for all~$V',V''$.
Since the map~\eref{SpinPinCorr_e} is identity under the identification~\eref{SpinPinSES_e33}, 
they satisfy the second property in~\ref{DSEquivSum_it}.
By~\ref{DSassoc_it} for $V_1',V_2',V''$ of ranks at least~3,
they satisfy~\ref{DSassoc_it} for all vector bundles $V_1',V_2',V''$.
By the SpinPin~\ref{SpinPinStab_prop} property and~\eref{SpinPinSES_e35},
the two maps satisfy~\ref{DSstab_it} for all~$V$.
By the SpinPin~\ref{SpinPinCorr_prop} property and~\eref{SpinPinSES_e35d},
they also satisfy~\ref{DSPin2Spin_it} for all~$V$.\\ 

By the SpinPin~\ref{OSpinRev_prop} and~\ref{SpinPinStab_prop} properties and~\eref{SpinPinSES_e35c},  
it is sufficient to establish the first property in~\ref{DSEquivSum_it} under the assumption that 
the ranks of~$V'$ and~$V''$ are at least~3 and the former is odd
(after possibly applying $\St_V^{a+b}$ to both~sides for some $a,b\!\in\!\Z^{\ge0}$
with $b$ even).
It is immediate from the definition of the map~\eref{SPinPullback_e} and 
the construction of the map~\eref{OSpinRev_e} in the perspective of Definition~\ref{PinSpin_dfn3}
that~\eref{DSCorr1crl_e2} holds in this perspective.
Combining this with the first equation in~\eref{DSCorr1crl_e3} and~\eref{SpinPinStr_e},
we obtain
$$\llrr{\ov\os',\fp''}_{\oplus}=\llrr{\Psi'^*\os',\fp''}_{\oplus}
=\Psi^*\llrr{\os',\fp''}_{\oplus}
=w_1(V'')\!\cdot\!\llrr{\os',\fp''}_{\oplus}\,.$$
This establishes the SpinPin~\ref{SpinPinSES_prop} property
for all vector bundles~$V'$ and~$V''$ over~$Y$.
\end{proof}

\subsection{The $\Spin$- and $\Pin$-structures of Definition~\ref{PinSpin_dfn2}}
\label{SpinPindfn2_subs}

Throughout this section, the terms $\Spin$-structure and $\Pin^{\pm}$-structure
refer to the notions arising from Definition~\ref{PinSpin_dfn2}.
As these notions are restricted to CW complexes, 
the properties of Section~\ref{SpinPinProp_subs} are fairly easy to establish
in this setting.\\

For $n\!\in\!\Z^{\ge0}$,
a topological space $Z$ is called \sf{$n$-connected}\gena{nconn@$n$-connected} if
$\pi_k(Z)$ is trivial for every $k\!\in\!\Z^{\ge0}$ with \hbox{$k\!\le\!n$}. 
For topological spaces~$X$ and~$Z$, we denote by $[X,Z]$\nota{XtoZ2@$[X,Z]$} 
the set of
homotopy classes of continuous maps from~$X$ to~$Z$. 
If in addition $Y\!\subset\!X$ and $z_0\!\in\!Z$, let 
$[(X,Y),(Z,z_0)]$\nota{XtoZ1@$[(X,Y),(Z,z_0)]$}  be the set of
homotopy classes of continuous maps from~$X$ to~$Z$ sending~$Y$ to~$z_0$.
If $Z$ is a group, then so are $[X,Z]$ and $[(X,Y),(Z_0,\id)]$.
If $Y$ is a CW complex and $k\!\in\!\Z^{\ge0}$, we denote by~$Y_k$ 
the $k$-skeleton of~$Y$.
Lemma~\ref{CWhomotop_lmm}  is used in the proof of the SpinPin~\ref{SpinPinStr_prop} property below
and in Section~\ref{TopolPrelim2_subs}.

\begin{lmm}\label{CWhomotop_lmm}
Suppose $n\!\in\!\Z^+$, $Z$ is an $(n\!-\!1)$-connected topological space with $|\pi_n(Z)|\!=\!2$
and $\pi_{n+1}(Z)$ trivial, $z_0\!\in\!Z$, and
$$\eta_Z\in H^n\big(Z,\{z_0\};\Z_2\big)\!=\!H^n\big(Z;\Z_2\big)$$
is the generator.
Let $X$ be a CW complex and $Y\!\subset\!X$ be a subcomplex.
The~map
\BE{CWhomotop_e0a}
\big[(X_{n+1},Y_n),(Z,z_0)\big]\lra H^n\big(X_{n+1},Y_n;\Z_2\big),
\qquad [h]\lra h^*\eta_Z,\EE
is then bijective, while the map
 \BE{CWhomotop_e0b}
\big[(X_{n+2},Y_n),(Z,z_0)\big]\lra H^n\big(X_{n+2},Y_n;\Z_2\big),
\qquad [h]\lra h^*\eta_Z,\EE
is surjective.
\end{lmm}

\begin{proof}
By  the Hurewicz Theorem \cite[Proposition~7.5.2]{Spanier} and 
the Universal Coefficient Theorem for Cohomology \cite[Theorem~53.5]{Mu2},
the homomorphism
$$\pi_n(X)\lra H_n(X;\Z)\lra \Z_2, \qquad [f\!:S^n\!\lra\!X]\lra \blr{\eta_Z,f_*[S^n]_{\Z_2}},$$
is an isomorphism.
Along with \cite[Theorem~8.1.15]{Spanier}, this implies~that the~map
\BE{CWhomotop_e3}
\big[(X_n,Y_n),(Z,z_0)\big]\lra H^n\big(X_n,Y_n;\Z_2\big),
\qquad [h]\lra h^*\eta_Z,\EE
is a bijection.\\

Suppose $\de\!=\!1,2$ and $\eta\!\in\!H^n(X_{n+\de},Y_n;\Z_2)$. 
By the bijectivity of~\eref{CWhomotop_e3}, 
there exists a continuous map~$f$ from $(X_n,Y_n)$ to $(Z,z_0)$ so~that
\BE{CWhomotop_e4}\eta\big|_{(X_n,Y_n)}=f^*\eta_Z\in H^n\big(X_n,Y_n;\Z_2\big)\,.\EE
By \cite[Theorem~8.1.17]{Spanier}, $f$ extends to a continuous map~$h$ from $X_{n+1}$ 
to~$Z$.
Since $\pi_{n+1}(Z)$ is trivial, $h$ then extends to a continuous map
from~$X_{n+2}$ to~$Z$.
Thus, there exists a continuous map~$h$ from $(X_n,Y_n)$ to $(Z,z_0)$ so~that
\BE{CWhomotop_e5}
\eta\big|_{(X_n,Y_n)}=(h^*\eta_Z)\big|_{(X_n,Y_n)}\in H^n\big(X_n,Y_n;\Z_2\big)\,.\EE
By the cohomology exact sequence for the triples $Y_n\!\subset\!X_n\!\subset\!X_{n+1}$
and \hbox{$Y_n\!\subset\!X_n\!\subset\!X_{n+2}$}, the restriction homomorphisms
$$H^n\big(X_{n+1},Y_n;\Z_2\big)\lra H^n\big(X_n,Y_n;\Z_2\big)
\quad\hbox{and}\quad 
H^n\big(X_{n+2},Y_n;\Z_2\big)\lra H^n\big(X_n,Y_n;\Z_2\big) $$
are injective.
Along with~\eref{CWhomotop_e5}, this establishes the surjectivity of~\eref{CWhomotop_e0a}
and~\eref{CWhomotop_e0b}.\\

Since $\pi_{n+1}(Z)$ is trivial, any two extensions of a continuous map $f\!:X_n\!\lra\!Z$
to a continuous map \hbox{$h\!:X_{n+1}\!\lra\!Z$} are homotopic with~$h|_{X_n}$ fixed.
By the bijectivity of~\eref{CWhomotop_e3}, any two maps~$f$ from~$(X_n,Y_n)$ to~$(Z,z_0)$
satisfying~\eref{CWhomotop_e4} are homotopic.
Combining the last two statements, we conclude that the map~\eref{CWhomotop_e0a} is injective.
\end{proof}

\begin{crl}\label{SOnH1_crl}
Suppose $n\!\ge\!3$ and $Y$ is a CW complex.
Let $\eta_n\!\in\!H^1(\SO(n);\Z_2)$ be the generator.
The homomorphism 
\BE{Spin1a_e5} \big[Y_2,\SO(n)\big]\lra H^1(Y_2;\Z_2), \qquad [h]\lra h^*\eta_n,\EE
is an isomorphism.
The homomorphism 
$$\big[Y_3,\SO(n)\big]\lra H^1(Y_3;\Z_2), \qquad [h]\lra h^*\eta_n,$$
is surjective.
\end{crl}

\begin{proof} This follows immediately from~\eref{pi012SO_e} and
Lemma~\ref{CWhomotop_lmm} with $n\!=\!1$, $Z\!=\!\SO(n)$, 
and~$(X,Y)$ replaced by~$(Y,\eset)$.
\end{proof}

\begin{proof}[{\bf{\emph{Proof of SpinPin~\ref{SpinPinObs_prop} property}}}]
By the definition of $\Spin$- and $\Pin^{\pm}$-structures in this perspective,
it is sufficient to establish this claim for $\Spin$-structures on oriented vector bundles~$(V,\fo)$
with $\rk\,V\!\ge\!3$.
If the vector bundle $V|_{Y_2}$ is trivializable, then $w_2(V)|_{Y_2}\!=\!0$.
Since the restriction homomorphism
\BE{SpinPinObs2_e2}H^2(Y;\Z_2)\lra H^2(Y_2;\Z_2)\EE
is injective by the cohomology exact sequence for the pair $Y_2\!\subset\!Y$, 
it follows that $w_2(V)\!=\!0$.\\

Suppose $w_2(V)\!=\!0$ and $\fo\!\in\!\fO(V)$.
By the SpinPin~\ref{SpinPinObs_prop} property for perspective of Definition~\ref{PinSpin_dfn},
$(V,\fo)$ then admits a $\Spin$-structure
$$q_V\!:\Spin(V,\fo)\lra\SO(V,\fo)$$
in the sense of Definition~\ref{PinSpin_dfn}.
By~\eref{pi012Spin_e},  
the principal $\Spin(n)$-bundle $\Spin(V,\fo)|_{Y_2}$ is trivializable.
For any section~$\wt{s}$ of $\Spin(V,\fo)|_{Y_2}$, the section $q_V\!\circ\!\wt{s}$ of 
$\SO(V,\fo)|_{Y_2}$
determines a trivialization of $(V,\fo)|_{Y_2}$ and thus 
a $\Spin$-structure on~$(V,\fo)$ in the sense of Definition~\ref{PinSpin_dfn2}.
\end{proof}

\begin{proof}[{\bf{\emph{Proof of SpinPin~\ref{SpinPinStr_prop} property}}}]
By the definition of $\Spin$- and $\Pin^{\pm}$-structures in this perspective,
it is sufficient to establish the claims of this property other than~\eref{SpinPinStr_e}  
for $\Spin$-structures on oriented vector bundles~$(V,\fo)$ with 
\hbox{$n\!\equiv\!\rk\,V\!\ge\!3$}.
By the cohomology exact sequence for the pair $Y_2\!\subset\!Y$, 
the restriction homomorphism
\BE{Spin1a_e3} H^1(Y;\Z_2)\lra H^1(Y_2;\Z_2)\EE
is an isomorphism.\\

If $\fo\!\in\!\fO(V)$ and $(V,\fo)$ admits a trivialization over~$Y_2$,
then the natural action of $[Y_2,\SO(n)]$ on the set of homotopy equivalence classes
of such trivializations given~by
\BE{Spin1a_e4}[Y_2,\SO(n)]\!\times\!\Sp(V,\fo)\lra\Sp(V,\fo), \qquad
\big([h],[\Phi]\big)\lra\big[h\Phi],\EE
is free and transitive.
Combining this action with the isomorphisms~\eref{Spin1a_e3} and~\eref{Spin1a_e5},
we obtain the claims of the SpinPin~\ref{SpinPinStr_prop} property other
 than~\eref{SpinPinStr_e}.\\  

By the naturality of the above action of $H^1(Y;\Z_2)$ on~$\Sp(V,\fo)$ 
with respect to continuous maps, it satisfies~\eref{SpinPinStr_e33} 
and thus~\eref{SpinPinStr_e35}.
The last paragraph in the proof of the SpinPin~\ref{SpinPinStr_prop} property
in the perspective of Definition~\ref{PinSpin_dfn3} now applies verbatim and
implies that
\eref{SpinPinStr_e} holds in the perspective of Definition~\ref{PinSpin_dfn3} as~well.
\end{proof}

\begin{proof}[{\bf{\emph{Proof of SpinPin~\ref{OSpinRev_prop},\ref{Pin2SpinRed_prop} properties}}}]
For the purposes of establishing the first of these properties,
we can assume that \hbox{$n\!\equiv\!\rk\,V\!\ge\!3$}.
For each $\fo\!\in\!\fO(V)$ and a homotopy class $\fs\!\in\!\Sp(V,\fo)$
of trivializations $\Phi$ of~$(V,\fo)$ over the 2-skeleton~$Y_2$ of~$Y$,
we take $\ov\fs\!\in\!\Sp(V,\ov\fo)$ to be the homotopy class
of trivializations of~$(V,\ov\fo)$ over the 2-skeleton~$Y_2$ of~$Y$
given~by 
\BE{OSpinRev_e23}
\ov\fs=\big\{\ov\Phi\!\equiv\!\bI_{n;1}\Phi\!:\Phi\!\in\!\fs\big\},\EE
with $\bI_{n;1}$ as below~\eref{Inmdfn_e}.
By~\eref{SpinPinStr_e33}, the resulting bijection~\eref{OSpinRev_e}
is $H^1(Y;\Z_2)$-equivariant.
It also satisfies the last condition of the SpinPin~\ref{OSpinRev_prop} property.\\

\noindent
We next describe a bijection between $\Sp(V,\fo)$ and 
$$\cP^{\pm}(V)\equiv \Sp(V\!\oplus\!(2\!\pm\!1)\la(V),\fo_V^{\pm}\big)
=\Sp(V\!\oplus\!(2\!\pm\!1)\tau_Y,\fo(2\!\pm\!1)\fo_Y\big);$$
the last equality holds because the orientation $\fo$ on~$V$ determines 
a canonical homotopy class of trivializations of~$\la(V)$.
By the SpinPin~\ref{SpinPinObs_prop} property,
$(V,\fo)$ admits a $\Spin$-structure if and only if $V$ admits a $\Pin^{\pm}$-structure.
We can thus assume that $(V,\fo)$ admits a $\Spin$-structure.\\

\noindent
If $\rk\,V\!\ge\!3$, we identify a homotopy class~$\fs$ of trivializations of 
$V|_{Y_2}$ with the homotopy class of trivializations of $(V_{\pm},\fo_V^{\pm})|_{Y_2}$
given~by
\BE{Pin2SpinRed_e23}
\fs^{\pm}=\big\{\Phi_{\pm}\!\equiv\!\Phi\!\oplus\!(2\!\pm\!1)(\det\Phi)\!:
\Phi\!\in\!\fs\big\}.\EE
By~\eref{SpinPinStr_e33}, the resulting map~\eref{SpinPinStr_e36} is again $H^1(Y;\Z_2)$-equivariant
and thus a bijection.
We take the map~$\fR_{\fo}^{\pm}$ in~\eref{Pin2SpinRed_e} to be its inverse.
If~$\ov\Phi$ is a trivialization of $(V,\ov\fo)|_{Y_2}$ as in~\eref{OSpinRev_e23},
then
$$\big(\ov\Phi\big)_{\!\pm}=\bI_{n+2\pm1;3\pm1}\Phi_{\pm}\!:
\big(V\!\oplus\!(2\!\pm\!1)\la(V)\big)\big|_{Y_2} \lra Y_2\!\times\!\R^{n+2\pm1}.$$
Since this trivialization is homotopic to~$\Phi_{\pm}$,
the last claim in the SpinPin~\ref{Pin2SpinRed_prop} property holds.\\

\noindent
If $\rk\,V\!=\!2$, we first identify the trivializations of $(\tau_Y\!\oplus\!V)|_{Y_2}$
with the trivializations of $(V\!\oplus\!\tau_Y)|_{Y_2}$ in the obvious way and
then extend them to trivializations of $(V\!\oplus\!3\tau_Y)|_{Y_2}$ by the identity
on the last two $\tau_Y$-summands.
If $\rk\,V\!=\!1$, we first identify the trivializations of 
$$(2\tau_Y\!\oplus\!V)\big|_{Y_2}\equiv
\big(\tau_Y\!\oplus\!(\tau_Y\!\oplus\!V)\big)\big|_{Y_2}$$ 
with the trivializations of $((\tau_Y\!\oplus\!V)\!\oplus\!\tau_Y)|_{Y_2}$
 and of $(V\!\oplus\!2\tau_Y)|_{Y_2}$ 
in the obvious orientation-preserving ways and then extend the latter to the trivializations
of $(V\!\oplus\!3\tau_Y)|_{Y_2}$  by the identity on the last $\tau_Y$-summand.
In all cases, the resulting bijections~\eref{Pin2SpinRed_e} are $H^1(Y;\Z_2)$-equivariant
and satisfy~\eref{Pin2SpinRed_e37} and
the last requirement of the SpinPin~\ref{Pin2SpinRed_prop} property.
\end{proof}

\begin{proof}[{\bf{\emph{Proof of SpinPin~\ref{SpinPinStab_prop},\ref{SpinPinCorr_prop} properties}}}]
Let $n\!=\!\rk\,V$.
By the SpinPin~\ref{SpinPinObs_prop} property,
an oriented vector bundle $(V,\fo)$ admits a $\Spin$-structure if and only if
$(\tau_Y\!\oplus\!V,\St_V(\fo))$ does.
In order to establish the SpinPin~\ref{SpinPinStab_prop} property,
we can thus assume that $(V,\fo)$ admits a $\Spin$-structure.\\

If $n\!\ge\!3$, we identify a homotopy class~$\fs$ of trivializations of 
$V|_{Y_2}$ with the homotopy class of trivializations of 
$(\tau_Y\!\oplus\!V,\St_V(\fo))|_{Y_2}$ given~by
\BE{SpinPinStab_e24}\St_V\fs=\big\{\St_V\Phi\!\equiv\id_{\tau_{Y_2}}\!\oplus\!\Phi\!:\Phi\!\in\!\fs\big\}.\EE
It is immediate that $\St_V\fs$ satisfies the first condition after~\eref{SpinPinStab_e}.
With~$\ov\Phi$ as in~\eref{OSpinRev_e23}, 
$$\St_V\ov\Phi\!=\ov{\St_V\Phi}\!: \big(\tau_Y\!\oplus\!V)\big|_{Y_2}\lra Y_2\!\times\!\R^{n+1}\,.$$
Thus, the resulting first map in~\eref{SpinPinStab_e} satisfies the first equality 
in~\eref{SpinPinStab_e2}.
We define the second map in~\eref{SpinPinStab_e} by~\eref{SpinPinStab_e35}.
By~\eref{SpinPinStr_e33},
both maps~\eref{SpinPinStab_e}
are then $H^1(Y;\Z_2)$-equivariant bijections.
With~$\Phi_{\pm}$ as in~\eref{Pin2SpinRed_e23}, 
$$\big(\St_V\Phi\big)_{\pm}\!=\St_V^{\pm}\Phi_{\pm}\!: 
(\tau_Y\!\oplus\!V)_{\pm}\big|_{Y_2}\!=\!\big(\tau_Y\!\oplus\!V_{\pm}\big)|_{Y_2}
\lra Y_2\!\times\!\R^{n+3\pm1}\,.$$
Thus, the second equality in~\eref{SpinPinStab_e2} also holds.\\

For $n\!=\!1,2$, the last paragraph of
the proof of the SpinPin~\ref{SpinPinStab_prop} property
in the perspective of Definition~\ref{PinSpin_dfn3} now applies verbatim.
The two sides of~\eref{SpinPinCorr_e} are again the same by definition;
we take this map to be the identity.
The requirement after~\eref{SpinPinCorr_e} is then again forced by the definition
of the second map of the SpinPin~\ref{SpinPinObs_prop} in~\eref{SpinPinStab_e35}.
\end{proof}

\begin{proof}[{\bf{\emph{Proof of SpinPin~\ref{SpinPinSES_prop} property}}}]
As in the perspective of Definition~\ref{PinSpin_dfn}, it is sufficient 
to establish this property for the direct sum exact sequences as in~\eref{SpinPinDS_e0}.
Furthermore, an orientation~$\fo'$ on~$V'$ determines a homotopy class of trivializations 
of~$\la(V')$ and thus of isomorphisms
\BE{SpinPinSES_e23} \big(V'\!\oplus\!V''\big)_{\pm}
\equiv V'\!\oplus\!V''\!\oplus\!(2\!\pm\!1)\la(V'\!\oplus\!V'')
\approx  V'\!\oplus\!V''\!\oplus\!(2\!\pm\!1)\la(V'')
\equiv V'\!\oplus\!V''_{\pm}\,.\EE

\vspace{.15in}

\noindent
Let $V'$ and $V''$ be vector bundles over~$Y$ of rank~$n'$ and~$n''$, respectively,
$\fo'\!\in\!\fO(V')$, and $\fo''\!\in\!\fO(V'')$.
If $\fs'$ and $\fs''$ are homotopy classes of trivializations of 
$(V',\fo')|_{Y_2}$ and~$(V'',\fo'')_{Y_2}$,
we take $\llrr{\fs',\fs''}$ to be the homotopy class of trivializations of
$(V'\!\oplus\!V'',\fo'\fo'')|_{Y_2}$ given~by
\BE{dfn2sum_e}\llrr{\fs',\fs''}=
\big\{\Phi'\!\oplus\!\Phi''\!:\Phi'\!\in\!\fs',\,\Phi''\!\in\!\fs''\big\}.\EE
This determines the first map in~\eref{SpinPinSESdfn_e0} if \hbox{$n',n''\!\ge\!3$}.
Along with~\eref{SpinPinSES_e23}, this also determines the second map of this property
if $n'\!\ge\!3$ and \hbox{$\rk\,V''_{\pm}\!\ge\!3$} (i.e.~not $V_-''$ if $n''\!=\!1$).
By~\eref{SpinPinStr_e33}, both maps are then $H^1(Y;\Z_2)$-biequivariant.
These two maps satisfy \ref{DSsplit_it},  \eref{DSorient1_e}, 
and~\ref{DSassoc_it} for vector bundles of ranks at least~3
by definition.
With~$\Phi_{\pm}$ as in~\eref{Pin2SpinRed_e23},
$$\big(\Phi'\!\oplus\!\Phi''\big)_{\pm}\!=\!\Phi'\!\oplus\!\Phi''_{\pm}\!:
\big(V'|_{Y_2}\big)\!\oplus\!(V_{\pm}''|_{Y_2}\big)\lra Y_2\!\times\!\R^{n'+n''+2\pm1}\,.$$
Thus, \eref{DSorient2_e} is also satisfied in these cases.\\

For vector bundles $V',V''$ over~$Y$ not both of ranks at least~3,
we define the two maps in~\eref{SpinPinSESdfn_e0}
as in the proof of this property in the perspective
of Definition~\ref{PinSpin_dfn3}, i.e.~by~\eref{SpinPinSES_e35c},
but with all relevant bijections taken between the sets of $\Spin$- and $\Pin^{\pm}$-structures
in the perspective of Definition~\ref{PinSpin_dfn2}.
The entire part of that proof starting with~\eref{SpinPinSES_e34a}
applies verbatim in this case,
even though the domains and targets of the relevant bijections
now refer to different sets. 
\end{proof}

\begin{eg}\label{ga1Rpin_eg2}
Let $\ga_{\R;1}\!\lra\!\R\P^1$, $\fo_{T\R\P^1}$, $\prt_{\th}$, $f$, and~$g$
be as in Example~\ref{ga1Rpin_eg0}, $\Phi_0$ be as in~\eref{ga1Rpin_e}, 
and 
$\os_1(2\ga_{\R;1},\fo_{\ga_{\R;1}}^-)$ be the $\OSpin$-structure 
on $(2\ga_{\R;1},\fo_{\ga_{\R;1}}^-)$ as in Example~\ref{ga1Rpin_eg}.
The composition of a right inverse for the vector bundle homomorphism~$g$ 
with the trivialization of~$T(\R\P^1)$ determined by the vector field~$\prt_{\th}$ is given~by 
$$h\!:\tau_{\R\P^1}\lra 2\ga_{\R;1}, \quad
h\big([x,y],b\big)=
\bigg(\![x,y],b\frac{(-y)(x,y)}{x^2\!+\!y^2},b\frac{x(x,y)}{x^2\!+\!y^2}\bigg)\,.$$
Since
$$\Phi_0\!\circ\!(f\!\oplus\!h)\!:2\tau_{\R\P^1}\!=\!\R\P^1\!\times\!\C\lra \R\P^1\!\times\!\C,
\quad \big([\ne^{\fI\th}],c\big)=\big([\ne^{\fI\th}],\ne^{2\fI\th}c\big),$$
the automorphism $\id\!\oplus\!(\Phi_0\!\circ\!(f\!\oplus\!h)\!)$ of $3\tau_{\R\P^1}$ 
is homotopically non-trivial.
Thus,
\BE{ga1Rpin2_e}
\bllrr{\os_0\big(\tau_{\R\P^1},\fo_{\R\P^1}\big),
\os_0\big(T(\R\P^1),\fo_{T\R\P^1}\big)}_{\eref{RPnEES_e3}}
=\os_1\big(2\ga_{\R;1},\fo_{\ga_{\R;1}}^-\big).\EE
\end{eg}

\vspace{.1in}

For $n\!\in\!\Z^+$, we denote~by 
$$\cN_{\R\P^n}\R\P^1\approx (n\!-\!1)\ga_{\R;1}^* \approx (n\!-\!1)\ga_{\R;1}$$ 
the normal bundle for the standard embedding of~$\R\P^1$ into~$\R\P^n$.
Let
\BE{RP2pin2_e3}0\lra T(\R\P^1)\lra T(\R\P^n)|_{\R\P^1}\lra 
\cN_{\R\P^n}\R\P^1\lra 0\EE
be the associated short exact sequence of vector bundles over~$\R\P^1$.

\begin{eg}\label{RP2pin_eg2}
Let $\os_0(\R\P^1,\fo_{T\R\P^1})$, $\os_0(2\ga_{\R;1},\fo_{\ga_{\R;1}}^-)$,
$\os_1(2\ga_{\R;1},\fo_{\ga_{\R;1}}^-)$, and $\os_1(4\ga_{\R;2},\fo_{\ga_{\R;2}}^+)$
be the $\OSpin$-structures on~$\R\P^1$, $2\ga_{\R;1}$, and $4\ga_{\R;2}$
provided by Examples~\ref{ga1Rpin_eg0}-\ref{RP2pin_eg}.
Let~$\fp_1^-(\R\P^2)$ be the $\Pin^-$-structure on~$\R\P^2$ provided by Example~\ref{RP2pin_eg}
and~$\os_1^-(\R\P^2)$ be the $\OSpin$-structure on $(T\R\P^2)_-$ corresponding to~$\fp_1^-(\R\P^2)$
under the bijection~\eref{SpinPinCorr_e}.
We denote by~$\fp_0^-(\cN_{\R\P^2}\R\P^1)$ the $\Pin^-$-structure 
on $\cN_{\R\P^2}\R\P^1\!\approx\!\ga_{\R;1}$ 
corresponding to $\os_0(2\ga_{\R;1},\fo_{\ga_{\R;1}}^-)$
under the bijection~\eref{SpinPinCorr_e}.
Euler's exact sequences~\eref{RPnEES_e3} for~$\R\P^1$ and~$\R\P^2$ and an identification of 
$\cN_{\R\P^2}\R\P^1$ with~$\ga_{\R;1}^*$ induce the commutative square of 
exact rows and columns of vector bundles over~$\R\P^1$ in Figure~\ref{RP2pin_fig}.
By the definitions of~$\os_1^-(\R\P^2)$ and $\os_1(4\ga_{\R;2},\fo_{\ga_{\R;1}}^+)$
and the SpinPin~\ref{SpinPinCorr_prop} property,
\begin{equation*}\begin{split}
\bllrr{\os_0\big(\tau_{\R\P^1},
\fo_{\R\P^1}\big),\os_1^-(\R\P^2)\big|_{\R\P^1}}_{\eref{RPnEES_e3}}
&=\os_1\big(4\ga_{\R;1},\fo_{\ga_{\R;1}}^+\big)\big|_{\R\P^1}\\
&
=\bllrr{\os_1(2\ga_{\R;1},\fo_{\ga_{\R;1}}^-),\os_0(2\ga_{\R;1},\fo_{\ga_{\R;1}}^-)}_{\oplus}.
\end{split}\end{equation*}
By~\eref{ga1Rpin2_e} and~\eref{OspinDSassoc_e}, 
\begin{equation*}\begin{split}
&\bllrr{\os_1(2\ga_{\R;1},\fo_{\ga_{\R;1}}^-),\os_0(2\ga_{\R;1},\fo_{\ga_{\R;1}}^-)}_{\oplus}\\
&\hspace{1in}
=\bllrr{\!\bllrr{\os_0\big(\tau_{\R\P^1},\fo_{\R\P^1}\big),
\os_0\big(T(\R\P^1),\fo_{T\R\P^1}\big)}_{\eref{RPnEES_e3}},
\os_0(2\ga_{\R;1},\fo_{\ga_{\R;1}}^-)}_{\oplus}\\
&\hspace{1in}
=\bllrr{\os_0\big(\tau_{\R\P^1},\fo_{\R\P^1}\big),\bllrr{\os_0\big(\R\P^1,\fo_{T\R\P^1}\big),
\os_0(2\ga_{\R;1},\fo_{\ga_{\R;1}}^-)}_{\eref{RP2pin2_e3}}}_{\oplus}\,.
\end{split}\end{equation*}
Along with the second statement in the SpinPin~\ref{SpinPinSES_prop}\ref{DSEquivSum_it} property,
the last two equations give
\BE{RP2pin2_e7}
\fp_1^-(\R\P^2)\big|_{\R\P^1}=
\bllrr{\os_0\big(\R\P^1,\fo_{T\R\P^1}\big),
\fp_0^-\big(\cN_{\R\P^2}\R\P^1\big)}_{\eref{RP2pin2_e3}}.\EE
\end{eg}

\begin{figure}
$$\xymatrix{& 0\ar[d] & 0\ar[d] & 0\ar[d]  \\
0\ar[r] & \tau_{\R\P^1}\ar[r]\ar[d]_{\id} & 2\ga_{\R;1}^*\ar[r]\ar[d] & T(\R\P^1)\ar[r]\ar[d] & 0 \\
0\ar[r] & \tau_{\R\P^1}\ar[r]\ar[d] & 3\ga_{\R;1}^*\!\oplus\!\ga_{\R;1}^* \ar[r]\ar[d] 
& T(\R\P^2)\big|_{\R\P^1}\!\oplus\!\ga_{\R;1}^* \ar[r]\ar[d] & 0 \\
0\ar[r] & 0\ar[r] & \ga_{\R;1}^*\!\oplus\!\ga_{\R;1}^*\ar[r]\ar[d] 
&  \cN_{\R\P^2}\R\P^1\!\oplus\!\cN_{\R\P^2}\R\P^1\ar[r]\ar[d] & 0 \\
&& 0 & 0}$$ 
\caption{Commutative square of exact rows and columns of
vector bundles over~$\R\P^1$ used in Example~\ref{RP2pin_eg2};
the vertical map from $3\ga_{\R;1}^*\!\oplus\!\ga_{\R;1}^*$ is the projection
to the last two components.}
\label{RP2pin_fig}
\end{figure}

\begin{eg}\label{RP3spin_eg2}
Let $\os_0(\R\P^1,\fo_{T\R\P^1})$, $\os_0(2\ga_{\R;1},\fo_{\ga_{\R;1}}^-)$,
$\os_1(2\ga_{\R;1},\fo_{\ga_{\R;1}}^-)$, and $\os_1(4\ga_{\R;3},\fo_{\ga_{\R;3}}^+)$
be the $\OSpin$-structures on~$\R\P^1$, $2\ga_{\R;1}$, and $4\ga_{\R;3}$
provided by Examples~\ref{ga1Rpin_eg0}-\ref{RP3spin_eg}
and $\os_0(\cN_{\R\P^3}\R\P^1)$ the $\OSpin$-structure 
on $\cN_{\R\P^3}\R\P^1\!\approx\!2\ga_{\R;1}$ 
corresponding to $\os_0(2\ga_{\R;1},\fo_{\ga_{\R;1}}^-)$.
Euler's exact sequences~\eref{RPnEES_e3} for~$\R\P^1$ and~$\R\P^3$ and 
a standard identification of 
$\cN_{\R\P^3}\R\P^1$ with $2\ga_{\R;1}^*$ induce the commutative square of 
exact rows and columns of vector bundles over~$\R\P^1$ in Figure~\ref{RP2pin_fig}
with the middle and bottom entries in the right column replaced~by 
$T(\R\P^3)\big|_{\R\P^1}$ and $\cN_{\R\P^3}\R\P^1$, respectively.
The same computation as in Example~\ref{RP2pin_eg2} yields
\begin{equation*}\begin{split}
&\bllrr{\os_0\big(\tau_{\R\P^1},
\fo_{\R\P^1}\big),\os_1(\R\P^3)\big|_{\R\P^1}}_{\eref{RPnEES_e3}}\\
&\hspace{1in}
=\bllrr{\os_0\big(\tau_{\R\P^1},\fo_{\R\P^1}\big),\bllrr{\os_0\big(\R\P^1,\fo_{T\R\P^1}\big),
\os_0(2\ga_{\R;1},\fo_{\ga_{\R;1}}^-)}_{\eref{RP2pin2_e3}}}_{\oplus}\,.
\end{split}\end{equation*}
Along with the SpinPin~\ref{SpinPinSES_prop}\ref{DSstab_it} and~\ref{SpinPinStab_prop} properties,
this implies~that
\BE{RP3pin2_e7}
\os_1(\R\P^3)\big|_{\R\P^1}=
\bllrr{\os_0\big(\R\P^1,\fo_{T\R\P^1}\big),
\os_0\big(\cN_{\R\P^3}\R\P^1\big)}_{\eref{RP2pin2_e3}}.\EE
\end{eg}

\section{Equivalence of Definitions~\ref{PinSpin_dfn}, \ref{PinSpin_dfn2}, and~\ref{PinSpin_dfn3}}
\label{thmpf_sec}

\noindent
It remains to establish the second and third statements of 
Theorem~\ref{SpinStrEquiv_thm}.
This is done in Sections~\ref{PinSpin1vs3_subs} and~\ref{PinSpin1vs2_subs}, respectively, 
implementing the approach indicated after the statement of this theorem.
We first confirm a basic observation that underpins this approach.

\begin{proof}[{\bf{\emph{Proof of Lemma~\ref{SpinStrEquiv_lmm}}}}]
Let $n\!=\!\rk\,V$.
As in Section~\ref{SO2Spin_subs}, we denote by $\wt\bI_n\!\in\!\Spin(n)$ 
the identity element and by $\wh\bI_n\!\in\!\Spin(n)$ the other preimage 
of the identity $\bI_n\!\in\!\SO(n)$ under the projection~\eref{Spindfn_e}.
For $r\!=\!0,1$, let 
$$\al_{n;r}\!:[0,1]\lra\SO(n) \qquad\hbox{and}\qquad \wt\al_{n;r}\!:[0,1]\lra\Spin(n)$$
be a loop representing $r$ times the generator of $\pi_1(\SO(n))$ with 
$\al_{n;r}(0),\al_{n;r}(1)\!=\!\bI_n$
and its lift to a path with $\wt\al_{n;r}(0)\!=\!\wt\bI_n$ and 
$\wt\al_{n;r}(1)\!=\!\wh\bI_n^r$.\\

\noindent
Since every orientable vector bundle over $S^1$ is trivializable,
we can assume that \hbox{$(V,\fo)\!=\!n(\tau_{S^1},\fo_{S^1})$}.
By the SpinPin~\ref{SpinPinObs_prop} and~\ref{SpinPinStr_prop} properties 
in Section~\ref{SpinPinProp_subs},
there are two equivalence classes of $\Spin$-structures in 
the sense of Definition~\ref{PinSpin_dfn} on this vector bundle.
They are represented~by
\begin{gather*}
q_{V;r}\!:\Spin_r(V,\fo)\!\equiv\!
\big([0,1]\!\times\!\Spin(n)\big)\!\big/\!\sim\lra S^1\!\times\!\SO(n),\\
\quad (0,\wh\bI_n^r\wt{A})\!\sim\!\big(1,\wt{A}\big)~\forall\,\wt{A}\!\in\!\Spin(n),\quad
q_{V;r}\big([t,\wt{A}]\big)=\big(\ne^{2\pi\fI t},q_n(\wt{A})\big)
~\forall\,(t,\wt{A})\!\in\![0,1]\!\times\!\Spin(n),
\end{gather*}
with $r\!=\!0,1$.
Sections of these $\Spin(n)$-principal bundles and 
the induced sections of $S^1\!\times\!\SO(n)$ are given~by
\begin{alignat}{2}
\label{SpinStrEquiv_e7}
\wt{s}_{V;r}\!:S^1&\lra \Spin_r(V,\fo), &\qquad  
\wt{s}_{V;r}\big(\ne^{2\pi\fI t}\big)&=\big[t,\wt\al_{n;r}(t)\big],\\
\notag
s_{V;r}\!:S^1&\lra S^1\!\times\!\SO(n), &\qquad  
s_{V;r}\big(\ne^{2\pi\fI t}\big)&=\big(\ne^{2\pi\fI t},\al_{n;r}(t)\big).
\end{alignat}
Since the sections $s_{V;0}$ and $s_{V;1}$ are not homotopic,
the homotopy classes of the associated trivializations of $S^1\!\times\!\R^n$ are
different.
\end{proof}

\subsection{Proof of Theorem~\ref{SpinStrEquiv_thm}\ref{SpinStrEquiv1vs3_it}}
\label{PinSpin1vs3_subs}

Let $(V,\fo)$ be a rank $n\!\ge\!3$ oriented vector bundle over a 
paracompact locally contractible space~$Y$.
By the SpinPin~\ref{SpinPinObs_prop} property,
$(V,\fo)$ admits a $\Spin$-structure in the sense of Definition~\ref{PinSpin_dfn}\ref{SpinStrDfn_it}
if and only if it admits a $\Spin$-structure in the sense of 
Definition~\ref{PinSpin_dfn3}.
We can thus assume that $(V,\fo)$ admits a $\Spin$-structure 
\BE{thmpf_e3}q_V\!:\Spin(V,\fo)\lra\SO(V,\fo)\EE
in the sense of Definition~\ref{PinSpin_dfn}\ref{SpinStrDfn_it}.\\

For every $\al\!\in\!\cL(Y)$,
\BE{thmpf_e11}\al^*q_V\!:\al^*\Spin(V,\fo)\lra \al^*\SO(V,\fo)\!=\!\SO\big(\al^*(V,\fo)\!\big)\EE
is a $\Spin$-structure over~$S^1$ in the sense of Definition~\ref{PinSpin_dfn}\ref{SpinStrDfn_it}.
Since $\Spin(n)$ is path-connected, the principal bundle 
$\al^*\Spin(V,\fo)$ admits a section~$\wt{s}_{\al}$.
Since in addition $\Spin(n)$ is simply connected, any two such sections are homotopic.
We take the homotopy class~$\fs_{\al}$ of trivializations of~$\al^*(V,\fo)$ associated to
the $\Spin$-structure~\eref{thmpf_e3} to be the homotopy class represented by 
the trivialization induced by the section $(\al^*q_V)\!\circ\!\wt{s}_{\al}$
of~$\SO(\al^*(V,\fo))$.\\

For every continuous map $F\!:\Si\!\lra\!Y$ from a compact bordered surface,
$$F^*q_V\!:F^*\Spin(V,\fo)\lra F^*\SO(V,\fo)\!=\!\SO(F^*(V,\fo)\!\big)$$
is a $\Spin$-structure over~$\Si$ in the sense of Definition~\ref{PinSpin_dfn}\ref{SpinStrDfn_it}.
Since $\Spin(n)$ is path-connected and simply connected, 
the principal bundle $F^*\Spin(V,\fo)$ admits a section~$\wt{s}_F$.
The restrictions of the trivialization of~$F^*(V,\fo)$ induced by the section 
$(F^*q_V)\!\circ\!\wt{s}_F$ of~$\SO(F^*(V,\fo)\!)$ to the boundary components 
$$\al_1\!\equiv\!F|_{C_1},\ldots,\al_m\!\equiv\!F|_{C_m}\!:S^1\lra Y$$
of~$F$ are induced by the sections $\wt{s}_{\al_i}\!\equiv\!\wt{s}_F|_{C_i}$
of~$F^*\Spin(V,\fo)|_{C_i}$.
Thus, the homotopy classes of trivializations constructed in the previous paragraph
satisfy the condition of Definition~\ref{PinSpin_dfn3}.\\

We conclude that a $\Spin$-structure~\eref{thmpf_e3} on~$(V,\fo)$ 
in the sense of Definition~\ref{PinSpin_dfn}\ref{SpinStrDfn_it} determines
a $\Spin$-structure 
$$(\fs_{\al})_{\al\in\cL(Y)}\equiv\Th_3\big(\Spin(V,\fo),q_V\big)$$ 
on~$(V,\fo)$ in the sense of Definition~\ref{PinSpin_dfn3}.
It is immediate that the latter depends only on the 
equivalence class of the $\Spin$-structure~\eref{thmpf_e3}.
By the construction, the resulting~map 
\BE{Th3dfn_e}\Th_3\!: \OSp_1(V)\lra\OSp_3(V)   \EE 
from the set of equivalence classes of $\OSpin$-structures on~$V$ 
in the sense of Definition~\ref{PinSpin_dfn}\ref{SpinStrDfn_it}
to the set of $\OSpin$-structures on~$V$ in the sense of Definition~\ref{PinSpin_dfn3}
is natural with respect to the pullbacks induced by continuous maps.
We show below~that $\Th_3$ is compatible with the
SpinPin~\ref{SpinPinStr_prop}\ref{SpinStr_it} property, i.e.
\BE{thmpf_e15} \Th_3\big(\eta\!\cdot\!(\Spin(V,\fo),q_V)\big)=  
\eta\!\cdot\!\big(\Th_3(\Spin(V,\fo),q_V)\big)
\qquad\forall~\eta\!\in\!H^1(Y;\Z_2)\,.\EE
Along with the SpinPin~\ref{SpinPinStr_prop}\ref{SpinStr_it} property, 
this implies that~\eref{Th3dfn_e} is a bijection.\\

Suppose $(V,\fo)\!\equiv\!n(\tau_Y,\fo_Y)$ is the trivial rank~$n$ oriented vector bundle
over~$Y$ and $q_V$ is the $\Spin$-structure on~$(V,\fo)$ given by~\eref{osVfodfn_e}.
For each $\al\!\in\!\cL(Y)$, we can then take 
$$\wt{s}_{\al}\!: S^1\lra \al^*\Spin(V,\fo)\!=\!S^1\!\times\!\Spin(n), \qquad
\wt{s}_{\al}(x)=(x,\wt\bI_n).$$
The induced homotopy class~$\fs_{\al}$ of trivializations of $(V,\fo)$ is then
the canonical homotopy class of trivializations of $S^1\!\times\!\R^n$.
For any oriented vector bundle~$(V,\fo)$ split into oriented line bundles,
the map~\eref{Th3dfn_e} thus takes
the distinguished element $\fs_0(V,\fo)$ in the perspective of Definition~\ref{PinSpin_dfn} to
the the distinguished element $\fs_0(V,\fo)$ in the perspective of Definition~\ref{PinSpin_dfn3}.

\begin{proof}[{\bf{\emph{Compatibility with SpinPin~\ref{SpinPinStr_prop}\ref{SpinStr_it} property}}}]
Let $\eta\!\in\!H^1(Y;\Z_2)$ and $\al\!\in\!\cL(Y)$.
By the naturality of the $H^1(Y;\Z_2)$-action of Propositions~\ref{PinSpin_prp1}\ref{Spin1_it},
\BE{thmpf_e17}\al^*\big(\eta\!\cdot\!(\Spin(V,\fo),q_V)\big)= 
\big\{\al^*\eta\big\}\!\cdot\!\big(\al^*\Spin(V,\fo),\al^*q_V\big)\,.\EE
By the SpinPin~\ref{SpinPinStr_prop} property with $Y\!=\!S^1$,
the $\Spin$-structure~\eref{thmpf_e17} is equivalent to 
the $\Spin$-structure~\eref{thmpf_e11} if and only~if $\al^*\eta\!=\!0$.
Along with Lemma~\ref{SpinStrEquiv_lmm}, this implies that the homotopy class 
of trivializations of $\al^*V$ determined by the $\Spin$-structure $\eta\!\cdot\!(\Spin(V,\fo),q_V)$
is the same as~$\fs_{\al}$ if and only~if $\al^*\eta\!=\!0$.
Combining this with~\eref{SpinPinStr_e33}, we obtain~\eref{thmpf_e15}.
\end{proof}

\begin{proof}[{\bf{\emph{Compatibility with SpinPin~\ref{OSpinRev_prop} property}}}]
Suppose $Y\!=\!S^1$ and $(V,\fo)$ is the trivial rank~$n$ oriented vector bundle over~$Y$.
For $r\!=\!0,1$, let $\al_r$, $\wt\al_r$, $q_{V;r}$, $\Spin_r(V,\fo)$, $\wt{s}_{V;r}$,
and~$s_{V;r}$ be as in the proof of Lemma~\ref{SpinStrEquiv_lmm}.
In particular, the homotopy class~$\fs_{\id}$ of trivializations of~$(V,\fo)$ 
associated to the $\OSpin$-structure $(\fo,q_{V;r})$ is represented by 
the trivialization~$\phi_{n;r}$ induced by~$s_{V;r}$.\\

By the paragraph containing~\eref{conjSpin_e},
\eref{OSpinRev_e} takes the $\Spin$-structure $q_{V;r}$ to 
the $\Spin$-structure 
on the oriented vector bundle~$(V,\ov\fo)$ given~by
\begin{gather*}
\ov{q}_{V;r}\!:\Spin_r(V,\ov\fo)\!\equiv\!
S^1\!\times\!\Spin(n)\lra \SO(V,\ov\fo)\!\equiv\!S^1\!\times\!\big(\O(n)\!-\!\SO(n)\!\big),\\
\ov{q}_{V;r}\big(\ne^{2\pi\fI t},\wt{A}\big)=\big(\ne^{2\pi\fI t},\al_{n;r}(t)\bI_{n;1}q_n(\wt{A})\!\big)
\quad\forall\,(t,\wt{A})\!\in\![0,1]\!\times\!\Spin(n).
\end{gather*}
A section of this $\Spin(n)$-principal bundle and 
the induced section of $S^1\!\times\!\SO(n)$ are given~by
$$\ne^{2\pi\fI t}\lra\big(\ne^{2\pi\fI t},\wt\bI_n\big),\qquad
\ne^{2\pi\fI t}\lra\big(\ne^{2\pi\fI t},\al_{n;r}(t)\bI_{n;1}\big)
\quad\forall\,t\!\in\![0,1].$$
The trivialization of $(V,\ov\fo)$ induced by the latter section is~$\bI_{n;1}\phi_{n;r}$.
Along with~\eref{OSpinRev_e33}, this implies that 
the homotopy class~$\ov\fs_{\id}'$ of trivializations of~$(V,\ov\fo)$ 
associated to the $\OSpin$-structure $(\ov\fo,\ov{q}_{V;r})$ is 
the homotopy class~$\ov\fs_{\id}$ 
of trivializations of~$(V,\ov\fo)$ determined by~$\fs_{\id}$.\\

By the naturality of~\eref{OSpinRev_e}, the conclusion of the previous paragraph implies~that
\BE{Th3OSpinRev_e}\Th_3\big(\ov\os\big)=\ov{\Th_3(\os)} \qquad\forall\,\os\!\in\!\OSp_1(V)\EE
for every vector bundle $V$ of rank at least~3 over a paracompact locally contractible space~$Y$.
\end{proof}

\begin{proof}[{\bf{\emph{Compatibility with SpinPin~\ref{SpinPinStab_prop} property}}}]
We continue with the setup at the beginning of the proof of the compatibility 
with the SpinPin~\ref{OSpinRev_prop} property.
By the proof of the SpinPin~\ref{SpinPinStab_prop} property in Section~\ref{SpinPinprp_subs}
and~\eref{wtioSpin_e2},
the first map in~\eref{SpinPinStab_e} takes the $\Spin$-structure~$q_{V;r}$ to
the $\Spin$-structure
$$q_{\St(V);r}\!: \Spin_r\big(\St(V,\fo)\!)\lra S^1\!\times\!\SO(n\!+\!1).$$
A section $\wt{s}_{\St(V);r}$ of this principal $\Spin(n\!+\!1)$-bundle over~$S^1$ 
is obtained by replacing $\wt\al_{n;r}$ with $\io_{n+1;n}''(\wt\al_{n;r})$ in~\eref{SpinStrEquiv_e7}.
The trivialization of $\St(V,\fo)$ induced by the trivialization 
$q_{\St(V);r}\!\circ\!\wt\al_{n;r}$ of $\SO(\St(V,\fo)\!)$ is
the trivialization $\St_V\phi_{n;r}$ in~\eref{SpinPinStab_e34}.
Thus, the homotopy class $(\St_V\fs)'_{\id}$ of trivializations of $\St_V(V,\fo)$ 
associated to the $\OSpin$-structure $\St_V(\fo,q_{V;r})$ is 
the homotopy class~$\St_V\fs_{\id}$ determined by~$\fs_{\id}$.
By the naturality of the first map in~\eref{SpinPinStab_e}, this implies~that
\BE{Th3SpinPinStab_e}
\Th_3\big(\St_V(\os)\!\big)=\St_V\big(\Th_3(\os)\!\big) \qquad\forall\,\os\!\in\!\OSp_1(V)\EE
for every vector bundle $V$ of rank at least~3 over a paracompact locally contractible space~$Y$.
\end{proof}

We define the map~$\Th_3$ in~\eref{Th3dfn_e} for rank~2 vector bundles $V$ over~$Y$ 
and then for rank~1 vector bundles $V$ over~$Y$ by~\eref{Th3SpinPinStab_e}.
By the SpinPin~\ref{SpinPinStab_prop} property and by the already established 
properties of~\eref{Th3dfn_e} for vector bundles of ranks at least~3,
the resulting maps~$\Th_3$ for vector bundles of ranks~1 and~2 are 
natural $H^1(Y;\Z_2)$-equivariant bijections which associate the distinguished 
elements~$\fs_0(V,\fo)$ in the perspectives of Definitions~\ref{PinSpin_dfn} and~\ref{PinSpin_dfn3}
with each other for all oriented vector bundles~$(V,\fo)$
split into oriented line bundles.
These maps are compatible with the SpinPin~\ref{SpinPinStab_prop} property,
i.e.~satisfy~\eref{Th3SpinPinStab_e}, by definition.
Along with the first equality in~\eref{SpinPinStab_e2} and~\eref{Th3OSpinRev_e} 
for vector bundles of rank at least~3, this implies that 
the maps~$\Th_3$ satisfy~\eref{Th3OSpinRev_e} for all vector bundles~$V$.

\begin{proof}[{\bf{\emph{Compatibility with SpinPin~\ref{SpinPinSES_prop} property}}}]
By the naturality of the first map in~\eref{SpinPinSESdfn_e0},
it is sufficient to show~that 
\BE{Th3SpinPinSES_e}\Th_3\big(\llrr{\os',\os''}_{\ce}\big)=
\bllrr{\Th_3(\os'),\Th_3(\os'')}_{\ce}
\qquad\forall\,\os'\!\in\!\OSp_1(V'),\,\os''\!\in\!\OSp_1(V''),\EE
for all short exact sequences~$\ce$ of trivial oriented vector bundles over $Y\!=\!S^1$
as in~\eref{SpinPinSES_e0}.
By the SpinPin~\ref{SpinPinStr_prop}\ref{SpinStr_it} property and 
the $H^1(Y;\Z_2)$-equivariance of~$\Th_3$ and of the first map in~\eref{SpinPinSESdfn_e0}
in the first input,
it is sufficient to check the $Y\!=\!S^1$ case of~\eref{Th3SpinPinSES_e} 
with $\os'$ being the distinguished $\OSpin$-structure $\os_0(m(\tau_{S^1},\fo_{S^1})\!)$.
By~\ref{DSsplit_it}, \ref{DSassoc_it}, and~\ref{DSstab_it} in the SpinPin~\ref{SpinPinSES_prop} property,
$$\bllrr{\os_0(m(\tau_Y,\fo_Y)\!),\os''}_{\oplus}=\St_V^m(\os'')\,.$$
Along with~\eref{Th3SpinPinStab_e}, this implies that~\eref{Th3SpinPinSES_e} holds 
with $\os'\!=\!\os_0(m(\tau_{S^1},\fo_{S^1})\!)$ and 
thus for all short exact sequences~\eref{SpinPinSES_e0} of vector bundles
over a paracompact locally contractible space~$Y$.
\end{proof}

With $\cP_1^{\pm}(V)$ and $\cP_3^{\pm}(V)$ denoting the sets~$\cP^{\pm}(V)$
of the $\Pin^{\pm}$-structures on a vector bundle~$V$ over~$Y$
in the perspectives of Definitions~\ref{PinSpin_dfn} and~\ref{PinSpin_dfn3}, 
respectively, define
\BE{Th3dfn2_e}\Th_3\!: \cP_1^{\pm}(V)\lra\cP_3^{\pm}(V) 
\qquad\hbox{by}\qquad  
\Co_V^{\pm}\big(\Th_3(\fp)\!\big)=\Th_3\big(\Co_V^{\pm}(\fp)\!\big)\,.\EE 
By the SpinPin~\ref{SpinPinCorr_prop} property and the already established 
properties of~\eref{Th3dfn_e},
\eref{Th3dfn2_e} is a natural $H^1(Y;\Z_2)$-equivariant bijection which is 
compatible with the SpinPin~\ref{SpinPinStab_prop} property, 
i.e.~\eref{Th3SpinPinStab_e} holds for $\Pin^{\pm}$-structures~$\fp$ in place of 
the $\OSpin$-structures~$\os$.
The map~$\Th_3$ is compatible with the SpinPin~\ref{SpinPinCorr_prop} property by definition.
By the second statement in the SpinPin~\ref{SpinPinSES_prop}\ref{DSEquivSum_it} property 
and~\eref{Th3SpinPinSES_e}, 
\eref{Th3dfn2_e} is compatible with the second map in~\eref{SpinPinSESdfn_e0},
i.e.~\eref{Th3SpinPinSES_e} holds for $\Pin^{\pm}$-structures~$\fp''$ in place of 
the $\OSpin$-structures~$\os''$.
By the SpinPin~\ref{SpinPinSES_prop}\ref{DSPin2Spin_it} property,
\eref{Th3dfn2_e} is compatible with the  SpinPin~\ref{Pin2SpinRed_prop} property,
i.e.
$$\Th_3\big(\fR_{\fo}^{\pm}(\fp)\!\big)=\fR_{\fo}^{\pm}\big(\Th_3(\fp)\!\big) 
\qquad\forall\,\fp\!\in\!\cP^{\pm}_1(V),\,\fo\!\in\!\fO(V),$$
for every vector bundle $V$ over a paracompact locally contractible space~$Y$.
This concludes the proof of Theorem~\ref{SpinStrEquiv_thm}\ref{SpinStrEquiv1vs3_it}.

\subsection{Proof of Theorem~\ref{SpinStrEquiv_thm}\ref{SpinStrEquiv1vs2_it}}
\label{PinSpin1vs2_subs}

Let $(V,\fo)$ be a rank $n\!\ge\!3$ oriented vector bundle over a CW complex~$Y$.
By the SpinPin~\ref{SpinPinObs_prop} property,
$(V,\fo)$ admits a $\Spin$-structure in the sense of Definition~\ref{PinSpin_dfn}\ref{SpinStrDfn_it}
if and only if it admits a $\Spin$-structure in the sense of 
Definition~\ref{PinSpin_dfn2}.
We can thus assume that $(V,\fo)$ admits a $\Spin$-structure 
in the sense of Definition~\ref{PinSpin_dfn}\ref{SpinStrDfn_it}
as in~\eref{thmpf_e3}.\\

Since $\Spin(n)$ is path-connected and simply connected,  
the principal bundle $\Spin(V,\fo)|_{Y_2}$ admits a section~$\wt{s}_2$.
Since in addition $\pi_2(\Spin(n))$ is trivial, any two such sections are homotopic.
We take the $\Spin$-structure 
\BE{Th2dfn_e0}\fs\equiv \Th_2(\Spin(V,\fo),q_V)\EE 
on~$(V,\fo)$ 
in the sense of Definition~\ref{PinSpin_dfn2} to be
the homotopy class of trivializations of~$(V,\fo)|_{Y_2}$ represented by 
the trivialization~$\Phi$ induced by the section $q_V\!\circ\!\wt{s}_2$
of~$\SO((V,\fo)|_{Y_2})$.\\

It is immediate that $\fs$ depends only on the 
equivalence class of the $\Spin$-structure~\eref{thmpf_e3}.
By the construction, the resulting~map 
\BE{Th2dfn_e}\Th_2\!: \OSp_1(V)\lra\OSp_2(V)   \EE 
from the set of equivalence classes of $\OSpin$-structures on~$V$ 
in the sense of Definition~\ref{PinSpin_dfn}\ref{SpinStrDfn_it}
to the set of $\OSpin$-structures on~$V$ in the sense of Definition~\ref{PinSpin_dfn2}
is natural with respect to the pullbacks induced by continuous maps.
We show below~that $\Th_2$ is compatible with the
SpinPin~\ref{SpinPinStr_prop}\ref{SpinStr_it} property, i.e.
\BE{thmpf_e25} \Th_2\big(\eta\!\cdot\!(\Spin(V,\fo),q_V)\big)=  
\eta\!\cdot\!\big(\Th_2(\Spin(V,\fo),q_V)\big)
\qquad\forall~\eta\!\in\!H^1(Y;\Z_2)\,.\EE
Along with the SpinPin~\ref{SpinPinStr_prop}\ref{SpinStr_it} property, 
this implies that~\eref{Th2dfn_e} is a bijection.\\

Suppose $(V,\fo)\!\equiv\!n(\tau_Y,\fo_Y)$ is the trivial rank~$n$ oriented vector bundle
over~$Y$ and $q_V$ is the $\Spin$-structure on~$(V,\fo)$ given by~\eref{osVfodfn_e}.
We can then take 
$$\wt{s}_2\!: Y_2\lra \Spin(V,\fo)\big|_{Y_2}\!=\!Y_2\!\times\!\Spin(n), \qquad
\wt{s}_2(x)=(x,\wt\bI_n).$$
The induced trivialization of $(V,\fo)$ is then
the canonical trivialization of $Y_2\!\times\!\R^n$.
For any oriented vector bundle~$(V,\fo)$ split into oriented line bundles,
the map~\eref{Th2dfn_e} thus takes
the distinguished element $\fs_0(V,\fo)$ in the perspective of Definition~\ref{PinSpin_dfn} to
the the distinguished element $\fs_0(V,\fo)$ in the perspective of Definition~\ref{PinSpin_dfn2}.

\begin{proof}[{\bf{\emph{Compatibility with SpinPin~\ref{SpinPinStr_prop}\ref{SpinStr_it} property}}}]
Suppose $\eta\!\in\!H^1(Y;\Z_2)$, 
$$\big(\Spin'(V,\fo),q_V'\big)\equiv \eta\!\cdot\!\big(\Spin(V,\fo),q_V\big)$$
is as in~\eref{PinSpin1_e4}, and
$h\!:Y_2\!\lra\!\SO(n)$ is a continuous map such that $\eta|_{Y_2}\!=\!h^*\eta_n$;
see the proof of the SpinPin~\ref{SpinPinStr_prop}\ref{SpinStr_it} property in
the perspective of Definition~\ref{PinSpin_dfn2} in Section~\ref{SpinPindfn2_subs}.
Let $\wt{s}_2'$ be a section of $\Spin'(V,\fo)|_{Y_2}$.
Since $\pi_2(\Spin(n)\!)$ is trivial, \eref{thmpf_e25} holds if and only if 
the trivializations of $\SO((V,\fo)|_{Y_1})$ induced by the restrictions of
$q_V'\!\circ\!\wt{s}_2'$ and $h\!\cdot\!(q_V\!\circ\!\wt{s}_2)$ are homotopic.\\

Let $\al\!\in\!\cL(Y)$.
By~\eref{thmpf_e15} and~\eref{SpinPinStr_e33}, the trivializations of $\SO(\al^*(V,\fo)\!)$ induced~by
$$\al^*\big\{q_V\!\circ\!\wt{s}_2\big\}=\big(\al^*q_V\big)\!\circ\!\big(\al^*\wt{s}_2\big)
\qquad\hbox{and}\qquad
\al^*\big\{q_V'\!\circ\!\wt{s}_2'\big\}=\big(\al^*q_V'\big)\!\circ\!\big(\al^*\wt{s}_2'\big)$$
are homotopic if and only if 
\BE{thmpf_e29}\big\{h\!\circ\!\al\big\}^*\eta_n=\al^*\eta=0\in H^1(S^1;\Z_2).\EE
The trivializations of $\SO(\al^*(V,\fo)\!)$ induced~by
$\al^*\big\{q_V\!\circ\!\wt{s}_2\big\}$ and
$$\al^*\big\{h\!\cdot\!(q_V\!\circ\!\wt{s}_2)\big\}
=(h\!\circ\!\al)\!\cdot\!\al^*\big\{q_V\!\circ\!\wt{s}_2\big\}$$
are homotopic if and only if the loop $h\!\circ\!\al\!:S^1\!\lra\!\SO(n)$ is homotopically trivial.
The latter is also equivalent to~\eref{thmpf_e29}.
Thus, the trivializations of $\SO((V,\fo)|_{Y_1})$ induced by the restrictions of
$q_V'\!\circ\!\wt{s}_2'$ and $h\!\cdot\!(q_V\!\circ\!\wt{s}_2)$ are homotopic.
This establishes~\eref{thmpf_e25}
for every vector bundle $V$ of rank at least~3 over~$Y$.
\end{proof}

\begin{proof}[{\bf{\emph{Compatibility with SpinPin~\ref{OSpinRev_prop} property}}}]
Suppose $Y$, $(V,\fo)$, $\wt{s}_2$, and $\Phi$ are as above and
$$\ov{q}_V\!: \Spin(V,\ov\fo)\lra \SO(V,\ov\fo)$$
is as below~\eref{conjSpin_e}.
Let $\wt{s}_2'$ be a section of $\Spin(V,\ov\fo)|_{Y_2}$.
Since $\pi_2(\Spin(n)\!)$ is trivial,
\BE{Th2OSpinRev_e}\Th_2\big(\Spin(V,\ov\fo),\ov{q}_V\big)=\ov{\Th_2(\Spin(V,\fo),q_V)}\EE
if and only if the trivialization of $\SO((V,\ov\fo)|_{Y_1})$  induced by 
the restriction of $\ov{q}_V\!\circ\!\wt{s}_2'$ is homotopic to~$\bI_{n;1}\Phi|_{Y_1}$;
see~\eref{OSpinRev_e23}.
The latter is the case if and only if the trivialization of $\SO(\al^*(V,\ov\fo)\!)$ induced~by
$$\al^*\big\{\ov{q}_V\!\circ\!\wt{s}_2'\big\}=\big(\al^*\ov{q}_V\big)\!\circ\!\big(\al^*\wt{s}_2'\big)$$
is homotopic to $\al^*(\bI_{n;1}\Phi)\!=\!\bI_{n;1}(\al^*\Phi)$ for every $\al\!\in\!\cL(Y)$.
The last statement is implied by~\eref{Th3OSpinRev_e}.
This establishes~\eref{Th2OSpinRev_e} for every vector bundle $V$ of rank at least~3 over~$Y$.
\end{proof}

\begin{proof}[{\bf{\emph{Compatibility with SpinPin~\ref{SpinPinStab_prop} property}}}]
Suppose $Y$, $(V,\fo)$, $\wt{s}_2$, and $\Phi$ are as above and
$$q_{\St(V)}\!: \Spin\big(\St(V,\fo)\!\big)\lra S^1\!\times\!\SO\big(\St(V,\fo)\!\big)$$
is as in the proof of the SpinPin~\ref{SpinPinStab_prop} property in Section~\ref{SpinPinprp_subs}.
Let $\wt{s}_2'$ be a section of $\Spin(\St(V,\fo)\!)|_{Y_2}$.
Since $\pi_2(\Spin(n)\!)$ is trivial,
\BE{Th2SpinPinStab_e}
\Th_2\big(\St_V(\Spin(V,\fo),q_V)\!\big)=\St_V\big(\Th_2(\Spin(V,\fo),q_V)\!\big)\EE
if and only if the trivialization of $\SO((V,\ov\fo)|_{Y_1})$ induced by 
the restriction of $q_{\St(V)}\!\circ\!\wt{s}_2'$ is homotopic to 
$\id_{\tau_{Y_1}}\!\oplus\!\Phi|_{Y_1}$;
see~\eref{SpinPinStab_e24}.
The latter is the case if and only if the trivialization of $\SO(\al^*\St(V,\fo)\!)$ induced~by
$$\al^*\big\{q_{\St(V)}\!\circ\!\wt{s}_2'\big\}=
\big(\al^*q_{\St(V)}\big)\!\circ\!\big(\al^*\wt{s}_2'\big)$$
is homotopic to 
$\al^*(\id_{\tau_{Y_1}}\!\oplus\!\Phi|_{Y_1})\!=\!\id_{\tau_{S_1}}\!\oplus\!(\al^*\Phi)$ 
for every $\al\!\in\!\cL(Y)$.
The last statement is implied by~\eref{Th3SpinPinStab_e}.
This establishes~\eref{Th2SpinPinStab_e} for every vector bundle $V$ of rank at least~3 over~$Y$.
\end{proof}

We define the map~$\Th_2$ in~\eref{Th2dfn_e} for rank~2 vector bundles $V$ over~$Y$ 
and then for rank~1 vector bundles $V$ over~$Y$ by~\eref{Th2SpinPinStab_e}.
By the same reasoning as in Section~\ref{PinSpin1vs3_subs}, 
the resulting maps~$\Th_2$ are natural $H^1(Y;\Z_2)$-equivariant bijections 
which associate the distinguished elements~$\fs_0(V,\fo)$ in the perspectives of 
Definitions~\ref{PinSpin_dfn} and~\ref{PinSpin_dfn2}
with each other for all oriented vector bundles~$(V,\fo)$
split into oriented line bundles and are 
compatible with the SpinPin~\ref{OSpinRev_prop}  and~\ref{SpinPinStab_prop} properties.

\begin{proof}[{\bf{\emph{Compatibility with SpinPin~\ref{SpinPinSES_prop} property}}}]
By Definition~\ref{PinSpin_dfn2},  the SpinPin~\ref{SpinPinStr_prop}\ref{SpinStr_it} property, and 
the $H^1(Y;\Z_2)$-equivariance of~$\Th_2$ and of the first map in~\eref{SpinPinSESdfn_e0}
in the first input, it is sufficient to show~that 
\BE{Th2SpinPinSES_e}\Th_2\big(\llrr{\os',\os''}_{\ce}\big)=
\bllrr{\Th_2(\os'),\Th_2(\os'')}_{\ce}
\qquad\forall\,\os'\!\in\!\OSp_1(V'),\,\os''\!\in\!\OSp_1(V''),\EE
for all short exact sequences~$\ce$ of trivial oriented vector bundles as in~\eref{SpinPinSES_e0}
with $\os'$ being the distinguished $\OSpin$-structure $\os_0(m(\tau_Y,\fo_Y)\!)$.
By~\ref{DSsplit_it}, \ref{DSassoc_it}, and~\ref{DSstab_it} in the SpinPin~\ref{SpinPinSES_prop} property,
$$\bllrr{\os_0(m(\tau_Y,\fo_Y)\!),\os''}_{\oplus}=\St_V^m(\os'')\,.$$
Along with~\eref{Th2SpinPinStab_e}, this implies that~\eref{Th2SpinPinSES_e} holds 
with $\os'\!=\!\os_0(m(\tau_{S^1},\fo_{S^1})\!)$ and 
thus for all short exact sequences~\eref{SpinPinSES_e0} of vector bundles over~$Y$.
\end{proof}

With $\cP_1^{\pm}(V)$ and $\cP_2^{\pm}(V)$ denoting the sets~$\cP^{\pm}(V)$
of the $\Pin^{\pm}$-structures on a vector bundle~$V$ over~$Y$
in the perspectives of Definitions~\ref{PinSpin_dfn} and~\ref{PinSpin_dfn2}, 
respectively, define
\BE{Th2dfn2_e}\Th_2\!: \cP_1^{\pm}(V)\lra\cP_2^{\pm}(V) 
\qquad\hbox{by}\qquad  
\Co_V^{\pm}\big(\Th_2(\fp)\!\big)=\Th_2\big(\Co_V^{\pm}(\fp)\!\big)\,.\EE
By the same reasoning as in Section~\ref{PinSpin1vs3_subs},
\eref{Th2dfn2_e} is a natural $H^1(Y;\Z_2)$-equivariant bijection which is 
compatible with the 
SpinPin~\ref{Pin2SpinRed_prop}-\ref{SpinPinSES_prop} properties.
This concludes the proof of Theorem~\ref{SpinStrEquiv_thm}\ref{SpinStrEquiv1vs2_it}.

\begin{eg}\label{SpinDfn1to3_eg}
With $V\!\lra\!\R\P^1$ denoting the infinite Mobius band
line bundle of Examples~\ref{Pin1pMB_eg}, \ref{Pin1mMB_eg},
and~\ref{Pin2Spin_eg}, 
\begin{alignat*}{3}
\pi_{V_-}\!:V_-&=\!\big([0,1]\!\times\!\C\big)\!\big/\!\!\sim\lra \R\P^1\!=\!S^1/\Z_2, &\quad
(0,-c)&\sim\!(1,c),&~~\pi_V\big([t,c]\big)&=\big[\ne^{\pi\fI t}\big],\\
\pi_{V_+}\!:V_+&=\!\big([0,1]\!\times\!\C^2\big)\!\big/\!\!\sim\lra \R\P^1\!=\!S^1/\Z_2, &\quad
(0,-\bfc)&\sim\!(1,\bfc),&~~\pi_V\big([t,\bfc]\big)&=\big[\ne^{\pi\fI t}\big].
\end{alignat*}
The orientations $\fo_V^-$ of $V_-$ and $\fo_V^+$ of $V_+$ described below~\eref{Vpmdfn_e2}
are the complex orientations of these bundles.
Let
$$\ga_{\pm}\!:[0,1]\lra\SO(3\!\pm\!1) \qquad\hbox{and}\qquad 
\wt\ga_{\pm}\!:[0,1]\lra\Spin(3\!\pm\!1) $$
be the paths as in the $n\!=\!2$ case of the proof of Lemma~\ref{SO2Spin_lmm} and
their lifts with \hbox{$\wt\ga_{\pm}(0)\!=\!\wt\bI_{3\pm1}$}.
By the conventions in~\eref{I24dfn_e}, $\wt\ga_-(1)\!=\!\wt\bI_{2;2}$ and
$\wt\ga_+(1)\!=\!\bI_{4;4}$.
Thus,
$$\wt{s}_{3\pm1}\!: \R\P^1\!\equiv\!S^1/\Z_2\lra\Spin_0(V_{\pm},\fo_V^{\pm}),\qquad 
\wt{s}_{3\pm1}\big([\ne^{\pi\fI t}]\big)=\big[t,\wt\ga_{\pm}(t)^{-1}\big],$$
is a well-defined continuous section of the $\Spin(3\!\pm\!1)$-bundle 
$\Spin_0(V_{\pm},\fo_V^{\pm})$ of Example~\ref{Pin2Spin_eg}.
The induced trivializations~$\Phi_2$ of~$(V_-,\fo_V^-)$ 
and~$\Phi_4$ of~$(V_+,\fo_V^+)$ are given~by 
\begin{alignat*}{2}
\Phi_2^{-1}\!:\R\P^1\!\times\!\C&\lra V_-, &\qquad
\Phi_2^{-1}\big([\ne^{\pi\fI t}],c\big)&=\big[t,\ne^{-\pi\fI t}c\big],\\
\Phi_4^{-1}\!:\R\P^1\!\times\!\C^2&\lra V_+, &\qquad
\Phi_4^{-1}\big([\ne^{\pi\fI t}],\bfc\big)&=\big[t,\ne^{-\pi\fI t}\bfc\big].
\end{alignat*}
Thus, the homotopy class of~$\Phi_4$ is the canonical $\OSpin$-structure 
$\os_0(V_+,\fo_V^+)$ on~$4V$ of Example~\ref{CanSpin_eg}.
The trivialization~$\Phi_2$ is the composition of 
the trivialization~\eref{ga1Rpin_e} of $2\ga_{\R;1}$
with the direct sum of two copies of the identification~\eref{Pin1mMBiden_e}.
\end{eg}

\begin{rmk}\label{SpinDfn1to3_rmk}
The map $\ga_-$ of Example~\ref{SpinDfn1to3_eg} is denoted by $\tn{rot}$ in the second case 
of \cite[Section~1.2.2]{Wel6b}.
By Example~\ref{SpinDfn1to3_eg}, the OSpin-structure $\os_0(2\ga_{\R;1},\fo_{\ga_{\R;1}}^-)$
of Example~\ref{ga1Rpin_eg} is not the base OSpin-structure on~$2\ga_{\R;1}$
in the second case of \cite[Section~1.2.2]{Wel6b} and in the corresponding case in
\cite[Section~2.2]{Wel6}.
Thus, the base OSpin-structure on~$2\ga_{\R;1}$ in \cite{Wel6,Wel6b} is
the OSpin-structure $\os_1(2\ga_{\R;1},\fo_{\ga_{\R;1}}^-)$ of Example~\ref{ga1Rpin_eg}.
Along with~\eref{RP3pin2_e7}, this implies that 
the base OSpin-structure on~$T(\R\P^3)|_{\R\P^1}$ in \cite{Wel6,Wel6b} is
the OSpin-structure $\os_0(\R\P^3)|_{\R\P^1}$ of Example~\ref{RP3spin_eg}.
While neither choice of OSpin-structure on~$2\ga_{\R;1}$ or
of $\Pin^-$-structure on~$\ga_{\R;1}$ appears more canonical,
some of the advantages $\os_0(2\ga_{\R;1},\fo_{\ga_{\R;1}}^-)$
over $\os_1(2\ga_{\R;1},\fo_{\ga_{\R;1}}^-)$ 
and of the associated $\Pin^-$-structure~$\fp_0^-(\ga_{\R;1})$ on~$\ga_{\R;1}$
over~$\fp_1^-(\ga_{\R;1})$ are described at the end of Section~\ref{OrientPrp_subs1},
after the CROrient~\ref{CRONormal_prop} property.
\end{rmk}

\section{Relative $\Spin$ and $\Pin$-structures}
\label{RelSpinPin_sec}

We recall the definitions of relative Spin- and Pin-structures from the CW perspective
of \cite{FOOO},
introduce a completely intrinsic perspective on these structures,
and state  the second theorem of this chapter, 
that these two perspectives are essentially equivalent, in Section~\ref{RelDfnThm_subs}.
Section~\ref{RelSpinPinProp_subs} describes key 
properties of relative Spin- and Pin-structures;
many of the stated properties resemble the properties of Spin- and Pin-structures collected 
in Section~\ref{RelSpinPinProp_subs}.
We establish these properties for the two perspectives in Sections~\ref{RelSpinPindfn2_subs}
and~\ref{RelSpinPindfn3_subs}, respectively, and
show the two perspectives to be equivalent, when restricted to CW complexes, 
in Section~\ref{RelPinSpin2vs3_subs}.
The proofs of the properties of Section~\ref{RelSpinPinProp_subs}
in the perspective of Definition~\ref{RelPinSpin_dfn2}
rely heavily on the simple topological observations of Section~\ref{TopolPrelim2_subs}.

\subsection{Definitions and main theorem}
\label{RelDfnThm_subs}

We begin with the CW perspective of \cite[Definition~8.1.2]{FOOO} on relative $\Spin$-structures.
Its extension to $\Pin$-structures appearing in \cite[Section~1.2]{Sol} is described
after we introduce a completely intrinsic notion of relative $\Spin$-structures.
We call $(X,Y)$ a \sf{CW pair} if $X$ is a CW complex and $Y$ is a CW subcomplex of~$X$.

\begin{dfn}\label{RelPinSpin_dfn2}
Let $(X,Y)$ be a CW pair and $V$ be a real vector bundle over~$Y$.
\begin{enumerate}[label=(\alph*),leftmargin=*]

\item\label{RelPinStrDfn2_it}
A \sf{relative $\Pin^{\pm}$-structure}~$\fp$\gena{relative Pin-structure!Solomon} 
on~$V$ is a tuple $(E,\fo_E,\fp_{E,V})$
consisting of an oriented vector bundle $(E,\fo_E)$ over the 3-skeleton~$X_3$ of~$X$
and a $\Pin^{\pm}$-structure~$\fp_{E,V}$ on $E|_{Y_2}\!\oplus\!V|_{Y_2}$
in the perspective of Definition~\ref{PinSpin_dfn2}.

\item\label{RelSpinStrDfn2_it} If $\fo\!\in\!\fO(V)$, 
a \sf{relative $\Spin$-structure}~$\fs$\gena{relative Spin-structure!FOOO} 
on~$(V,\fo)$ is a tuple $(E,\fo_E,\fs_{E,V})$
consisting of an oriented vector bundle $(E,\fo_E)$ over~$X_3$ 
and a Spin-structure~$\fs_{E,V}$ on $(E,\fo_E)|_{Y_2}\!\oplus\!(V,\fo)|_{Y_2}$
in the perspective of Definition~\ref{PinSpin_dfn2}.

\end{enumerate}
\end{dfn}

For $m\!\in\!\Z^{\ge0}$, we denote by $\St^m$ and $\St_V^m$ the compositions of
$m$-copies of the maps in~\eref{Stdfn_e} and~\eref{SpinPinStab_e}.
Relative $\Spin$-structures $(E,\fo_E,\fs_{E,V})$ and $(E',\fo_{E'},\fs_{E',V})$
on $(V,\fo)$ as in Definition~\ref{RelPinSpin_dfn2} are called 
\sf{equivalent}\gena{relative Spin-structure!equivalent}
if there exist $m,m'\!\in\!\Z^{\ge0}$ and an isomorphism 
$$\Psi\!:\St^m(E,\fo_E)\!\equiv\!\big(m\tau_{X_3}\!\oplus\!E,\St_E^m(\fo_E)\!\big)
\lra \St^{m'}\!(E',\fo_{E'})\!\equiv\!\big(m'\tau_{X_3}\!\oplus\!E',\St_E^{m'}\!(\fo_{E'})\!\big)$$ 
of oriented vector bundles over~$X_3$ so that the induced isomorphism
$$\Psi|_{m\tau_{Y_2}\oplus E|_{Y_2}}\!\oplus\!\id_{V|_{Y_2}}\!:
\St^m(E,\fo_E)\big|_{Y_2}\!\oplus\!(V,\fo)\big|_{Y_2}
\lra \St^{m'}\!(E',\fo_{E'})\big|_{Y_2}\!\oplus\!(V,\fo)\big|_{Y_2}$$
identifies the Spin-structures $\St_{E|_{Y_2}\oplus V|_{Y_2}}^m\!(\fs_{E,V})$ and
$\St_{E'|_{Y_2}\oplus V|_{Y_2}}^{m'}\!(\fs_{E',V})$.
The equivalence classes of relative $\Spin$-structures on~$(V,\fo)$
for different triangulations  of the pair $Y\!\subset\!X$
can be canonically identified; see \cite[Proposition~8.1.6]{FOOO}.
Notions of \sf{equivalence}\gena{relative Pin-structure!equivalent}
for relative $\Pin^{\pm}$-structures of Definition~\ref{RelPinSpin_dfn2}\ref{RelPinStrDfn2_it}
are defined similarly.
When there is no ambiguity, we will not distinguish between the Spin- and Pin-structures
of Definition~\ref{RelPinSpin_dfn2} and their equivalence classes.
For a relative $\Pin$ or $\Spin$-structure~$\fs$ on~$V$ as in Definition~\ref{RelPinSpin_dfn2},  
let
$$w_2(\fs)\nota{w2@$w_2(\fs),w_2(\os),w_2(\fp)$!FOOO}
=w_2(E)\in H^2(X_3;\Z_2)\!=\!H^2(X;\Z_2);$$
we use the same notation for equivalence classes of such structures.

\begin{rmk}\label{JakePin_rmk}
In \cite[Section~1.2]{Sol}, a relative $\Pin^{\pm}$-structure on
a real vector bundle~$V$ over $Y\!\subset\!X$ is defined (with some typos)
to be a tuple $(E,\fo_E,\fp_{E,V})$
consisting of an oriented vector bundle $(E,\fo_E)$ over~$X_3$ of~$X$
and a $\Pin^{\pm}$-structure~$\fp_{E,V}$ on $(E,\fo_E)|_{Y_2}\!\oplus\!(V,\fo)|_{Y_2}$.
For the primarily qualitative purposes of~\cite{Sol}, 
it is not material whether the latter is taken in 
the perspective of Definition~\ref{PinSpin_dfn}, \ref{PinSpin_dfn2}, or~\ref{PinSpin_dfn3},
but the approach in~\cite{Sol} leads to some difficulties in studying and applying
the disk invariants defined there.
\end{rmk}

Let $(V,\fo)$ be an oriented vector bundle over a compact nonempty one-dimensional manifold~$Y$, 
i.e.~a union of circles, with $\rk_{\R}V\!\ge\!3$.
We call two homotopy classes of trivializations of~$(V,\fo)$ \sf{equivalent} 
if the number of connected components of~$Y$ on which they differ is even. 
We call such an equivalence class a \sf{mod~2 homotopy class of trivializations} of~$(V,\fo)$. 
By~\eref{pi012SO_e}, there are two such classes.
If $Y$ is empty, we define a mod~2 homotopy class of trivializations to be an element of~$\Z_2$.
If $Y$ and $Y'$ are compact one-dimensional manifolds and 
$(V,\fo)$ is an oriented vector bundle over $Y\!\sqcup\!Y'$,
a mod~2 homotopy class~$\fs$ of trivializations of $(V,\fo)|_Y$ and
a mod~2 homotopy class~$\fs'$ of trivializations of $(V,\fo)|_{Y'}$
determine a mod~2 homotopy class~$\fs\fs'$ of trivializations of $(V,\fo)$.
If neither~$Y$ nor~$Y'$ is empty, $\fs\fs'$ is obtained by simply combining 
the trivializations in~$\fs$ and~$\fs'$.
If $Y'$ is empty, we take  $\fs\fs'\!=\!\fs$ if $\fs'\!=\!0$ and 
$\fs\fs'\!\neq\!\fs$ if $\fs'\!\neq\!0$.
We call $\fs\fs'$ the \sf{disjoint union} of the trivializations~$\fs$ and~$\fs'$.\\

Suppose $(V,\fo)$ is an oriented vector bundle over a bordered surface~$\Si$
with \hbox{$\rk_{\R}V\!\ge\!3$}. 
Let \hbox{$\Si^{\circ}\!\subset\!\Si$} be the union of the components~$\Si'$ of~$\Si$ with 
$\prt\Si'\!\neq\!\eset$,
$\Si^{\bu}\!\subset\!\Si$ be the union of the remaining components of~$\Si$, and
$$\fs_{\Si^{\bu}}(V,\fo)=\blr{w_2(V),[\Si^{\bu}]_{\Z_2}}
=\blr{w_2(V),[\Si]_{\Z_2}}\in\Z_2\,.$$
If $\Si^{\circ}\!\neq\!\eset$, let $\fs_{\Si^{\circ}}(V,\fo)$ be 
the mod~2 homotopy class of trivializations of 
$$(V,\fo)\big|_{\prt\Si}=(V,\fo)\big|_{\prt\Si^{\circ}}$$
containing the restriction of a trivialization~$\Phi$ of $(V,\fo)|_{\Si^{\circ}}$;
by Corollary~\ref{X2VB_crl2b}, $\fs_{\Si^{\circ}}(V,\fo)$ is well-defined.
If $\Si^{\circ}\!=\!\eset$, we take $\fs_{\Si^{\circ}}(V,\fo)\!=\!0$.
In either case, 
\BE{fsSiVfodfn_e}\fs_{\Si}(V,\fo)\equiv\fs_{\Si^{\bu}}(V,\fo)\fs_{\Si^{\circ}}(V,\fo)\EE
is a mod~2 homotopy class of trivializations of $(V,\fo)|_{\prt\Si}$.
In particular, \hbox{$\fs_{\eset}(V,\fo)\!=\!0\!\in\!\Z_2$}.\\

Let $Y$ be a subspace of a topological space~$X$. 
We denote by~$\cL_X(Y)$\nota{lyx@$\cL_X(Y)$} 
the collection of all continuous maps 
\hbox{$u\!:(\Si,\prt\Si)\!\lra\!(X,Y)$} from bordered surfaces.
For such a map, let \hbox{$\prt u\!=\!u|_{\prt\Si}$}. 
We call two such maps
$$u\!:(\Si,\prt\Si)\lra(X,Y) \qquad\hbox{and}\qquad
u'\!:(\Si',\prt\Si')\lra(X,Y)$$ 
\sf{relatively equivalent} if $\prt u\!=\!\prt u'$ and 
the continuous~map
$$u\!\cup\!u'\!:\Si\!\cup\!\Si'\lra X$$
obtained by gluing $u$ and $u'$ along the boundaries of their domains
represents the zero class in~$H_2(X;\Z_2)$.

\begin{dfn}\label{RelPinSpin_dfn3}
Let $X$ be a topological space and $(V,\fo)$ be an oriented vector bundle over \hbox{$Y\!\subset\!X$}
with $\rk_{\R}V\!\ge\!3$.
A \sf{relative $\Spin$-structure}~$\fs$ \gena{relative Spin-structure!intrinsic}
 on~$(V,\fo)$ is a collection $(\fs_u)_{u\in\cL_X(Y)}$
of mod~2 homotopy classes~$\fs_u$ of trivializations of $\{\prt u\}^*(V,\fo)$
such that for every continuous map \hbox{$F\!:\Si\!\lra\!Y$} from a bordered surface
relatively equivalent to the disjoint union of continuous maps
$$u_i\!:(\Si_i,\prt\Si_i)\lra(X,Y), \qquad i=1,\ldots,k,$$
the mod~2 homotopy class $\fs_{\Si}(F^*(V,\fo)\!)$ of trivializations of 
$\{\prt F\}^*(V,\fo)$ equals to 
the disjoint union of the mod~2 homotopy classes~$\fs_{u_i}$ of trivializations 
of~$\{\prt u_i\}^*(V,\fo)$.
\end{dfn}

By the $\Si\!=\!\eset$ case of the condition above,
the existence of a relative $\Spin$-structure $\fs$ in the sense of Definition~\ref{RelPinSpin_dfn3}
on an oriented vector bundle $(V,\fo)$ with \hbox{$\rk_{\R}V\!\ge\!3$} implies that
$\fs_{u_i}\!=\!0$ for every continuous map $u_i\!:\Si_i\!\lra\!X$ from a closed surface 
such that $u_{i*}[\Si_i]$ represents the zero element of~$H_2(X;\Z_2)$.
Thus, the~number 
\BE{w2fsdfn_e}\blr{w_2(\fs),u_*[\Si]_{\Z_2}}\nota{w2@$w_2(\fs),w_2(\os),w_2(\fp)$!intrinsic}
\equiv \fs_u\in\Z_2\EE
depends only on the class in $H_2(X;\Z_2)$ determined by a 
continuous map $u\!:\Si\!\lra\!X$ from a closed surface.
Along with Lemma~\ref{Z2cycles_lmm}\ref{cycles_it} with $Y\!=\!\eset$ 
and the Universal Coefficient Theorem
for Cohomology, this implies that $\fs$ determines an element 
\hbox{$w_2(\fs)\!\in\!H^2(X;\Z_2)$}.\\

In the perspectives of Definition~\ref{RelPinSpin_dfn3}, 
a \sf{relative Spin-structure}~$\fs$ on an oriented vector bundle~$(V,\fo)$ over $Y\!\subset\!X$
with $\rk_{\R}V\!\in\!\{1,2\}$ is a relative $\Spin$-structure on the vector bundle
$2\tau_Y\!\oplus\!V$ with the induced orientation in the first case and
on $\tau_Y\!\oplus\!V$ in the second.
A \sf{relative $\Pin^{\pm}$-structure}~$\fp$\gena{relative Pin-structure!intrinsic} 
on a real vector bundle~$V$ over $Y\!\subset\!X$
is a relative $\Spin$-structure on the canonically oriented vector bundle~\eref{Vpmdfn_e}.\\

In either of the two perspectives 
a \sf{relative $\OSpin$-structure}\gena{relative OSpin-structure}
on~$V$ is a pair \hbox{$\os\!\equiv\!(\fo,\fs)$} 
consisting of an orientation~$\fo$ on~$V$ and a relative $\Spin$-structure~$\fs$ on~$(V,\fo)$.
We denote by~$\OSp_X(V)$\nota{ospvx@$\OSp_X(V)$} and
$\cP_X^{\pm}(V)$\nota{pvx@$\cP_X^{\pm}(V)$} 
the sets of relative  $\OSpin$-structures and 
$\Pin^{\pm}$-structures, respectively, on~$V$
(up to equivalence in the perspective of Definition~\ref{RelPinSpin_dfn2}).
For $\fo\!\in\!\fO(V)$,  we denote by $\Sp_X(V,\fo)$\nota{spvox@$\Sp_X(V,\fo)$} 
the set of relative $\Spin$-structures on~$(V,\fo)$ in either perspective.
We identity $\Sp_X(V,\fo)$ with a subset of~$\OSp_X(V)$ in the obvious way.\\

In either perspective, there are natural maps 
\BE{vsSpinPin_e0}\io_X\!:\OSp(V)\lra \OSp_X(V) \qquad\hbox{and}\qquad
\io_X\!:\cP^{\pm}(V)\lra \cP_X^{\pm}(V)\,.\EE
In the perspective of Definitions~\ref{PinSpin_dfn2} and~\ref{RelPinSpin_dfn2}, 
these maps are obtained by taking the bundle~$E$ over~$X$ to be of rank~0.
If the rank of~$V$ is at least~3 and $(\fs_{\al})_{\al\in\cL(Y)}$ is a $\Spin$-structure
on $(V,\fo)$ in the perspective of Definition~\ref{PinSpin_dfn3},
we define the $\Spin$-structure
\BE{Spin2Rel_e3}
(\fs_u)_{u\in\cL_X(Y)}\equiv \io_X\big((\fs_{\al})_{\al\in\cL(Y)}\big)\EE
in the perspective of Definition~\ref{RelPinSpin_dfn3} as follows.
Let \hbox{$u\!:(\Si,\prt\Si)\!\lra\!(X,Y)$} be a continuous map from a bordered surface.
If $\prt\Si\!=\!\eset$, we take $\fs_u\!=\!0$.
If $\prt\Si\!\neq\!\eset$, we take~$\fs_u$ to be 
the disjoint union of the (mod~2) homotopy classes $\fs_{u|_{\prt_r\Si}}$ of
trivializations of $\{u|_{\prt_r\Si}\}^*(V,\fo)$ under any identification of
each topological component~$\prt_r\Si$ of~$\prt\Si$ with~$S^1$.
By the condition at the end of Definition~\ref{PinSpin_dfn3}, 
$\fs_u$~does not depend on the choice of such identifications.
The collection~\eref{Spin2Rel_e3} then satisfies 
the condition at the end of Definition~\ref{RelPinSpin_dfn3} and is thus 
is a relative $\Spin$-structure on~$(V,\fo)$.
This construction induces the first map in~\eref{vsSpinPin_e0} for vector bundles~$V$ of 
ranks~1 and~2 and the second map in~\eref{vsSpinPin_e0} for vector bundles. 

\begin{thm}\label{RelSpinStrEquiv_thm}
Let $X$ be a topological space and $(V,\fo)$ be an oriented vector bundle over \hbox{$Y\!\subset\!X$}.
\begin{enumerate}[label=(\arabic*),leftmargin=*]

\item\label{RelPinSpinPrp_it} The relative OSpin- and Pin-structures on~$V$
in the perspective of Definition~\ref{RelPinSpin_dfn2} 
if $(X,Y)$ is a CW pair 
and in the perspective of Definition~\ref{RelPinSpin_dfn3} 
satisfy all properties of Section~\ref{RelSpinPinProp_subs}.

\item\label{RelSpinStrEquiv2vs3_it} 
If $(X,Y)$ is a CW pair,
there are canonical identifications of the sets $\OSp_X(V)$ in the perspectives of 
Definitions~\ref{RelPinSpin_dfn2} and~\ref{RelPinSpin_dfn3} and
of the sets $\cP_X^{\pm}(V)$ in the two perspectives for every vector bundle~$V$ over~$Y$.
These identifications are intertwined with the identifications of
Theorem~\ref{SpinStrEquiv_thm} via the maps~\eref{vsSpinPin_e0} 
and respect all structures and correspondences of Section~\ref{RelSpinPinProp_subs}.

\end{enumerate}

\end{thm}

\subsection{Properties of relative $\Spin$ and $\Pin$ structures}
\label{RelSpinPinProp_subs} 

Let $X$ be a topological space and $Y\!\subset\!X$. 
An isomorphism $\Psi\!:V'\!\lra\!V$ of vector bundles
over~$Y$ induces bijections 
\BE{RelSPinPullback_e}\Psi_X^*\!: \OSp_X(V)\lra\OSp_X(V') \qquad\hbox{and}\qquad 
\Psi^*\!: \cP_X^{\pm}(V)\lra\cP_X^{\pm}(V')\EE
between the relative $\OSpin$-structures on~$V$ and~$V'$ and 
the $\Pin^{\pm}$-structures on~$V$ and~$V'$
in the perspective of Definition~\ref{RelPinSpin_dfn2} if $(X,Y)$ is a CW pair  
and in the perspective of Definition~\ref{RelPinSpin_dfn3}.
The first map sends a relative $\OSpin$-structure $(E,\fo_E,\os_{E,V})$
in the perspective of Definition~\ref{RelPinSpin_dfn2} to 
$$\big(E,\fo_E,\big\{\id_{E|_{Y_2}}\!\oplus\!\Psi|_{V|_{Y_2}}\big\}^*\os_{E,V}\big)
\in \OSp_X(V'),$$
with $\{\ldots\}^*$ on $\OSp(E|_{Y_2}\!\oplus\!V|_{Y_2})$ as in~\eref{SPinPullback_e}.
In the perspective of Definition~\ref{RelPinSpin_dfn3},
the relative $\OSpin$-structures on~$V'$
are obtained from the relative $\OSpin$-structures on~$V$ by 
pre-composing the relevant trivializations with~$\Psi$.\\

Let $X,X'$ be topological spaces, $Y\!\subset\!X$, $Y'\!\subset\!X'$, and 
$V$ be a vector bundle over~$Y$.
A continuous map \hbox{$f\!:(X',Y')\!\lra\!(X,Y)$} induces maps
\BE{RelfSPinPullback_e}f^*\!: \OSp_X(V)\lra\OSp_{X'}(f^*V) \qquad\hbox{and}\qquad 
f^*\!: \cP_X^{\pm}(V)\lra\cP_{X'}^{\pm}(f^*V)\EE
in the perspective of Definition~\ref{RelPinSpin_dfn2} if $f$ is a CW map between
CW pairs and in the perspective of Definition~\ref{RelPinSpin_dfn3}.
The first map sends a relative $\OSpin$-structure $(E,\fo_E,\os_{E,V})$
in the perspective of Definition~\ref{RelPinSpin_dfn2} to 
$$\big(f^*E,f^*\fo_E,\big\{\id_{E|_{Y_2}}\!\oplus\!\Psi|_{V|_{Y_2}}\big\}^*\os_{E,V}\big)
\in \OSp_X(V'),$$
with $f^*$ on $\OSp(E|_{Y_2}\!\oplus\!V|_{Y_2})$ as in~\eref{fSPinPullback_e}.
In the perspective of Definition~\ref{RelPinSpin_dfn3},
the relative $\OSpin$-structures on~$f^*V$ 
are obtained from the relative $\OSpin$-structures on~$V$ by 
pre-composing the relevant trivializations with~\eref{fVpullPsidfn_e}.\\

We denote by
\BE{YdeXYZ2_e}\de_{X,Y}\!:H^1(Y;\Z_2)\lra H^2(X,Y;\Z_2)\EE
the coboundary homomorphism in the cohomology exact sequence of the pair \hbox{$Y\!\subset\!X$}.
The RelSpinPin properties below apply in either of the perspectives of 
Definitions~\ref{RelPinSpin_dfn2} and~\ref{RelPinSpin_dfn3} on relative
$\Spin$-, $\Pin^{\pm}$-,  and $\OSpin$-structures,
provided $(X,Y)$ is a CW pair in the first perspective.
The naturality properties of the group actions and correspondences below
refer to the commutativity with the pullbacks~\eref{RelSPinPullback_e}
induced by isomorphisms of vector bundles over~$Y$ and 
the pullbacks~\eref{RelfSPinPullback_e} induced by the admissible continuous maps.

\begin{RelSpinPin}[Obstruction to Existence]\label{RelSpinPinObs_prop}
Let $V$ be a vector bundle over $Y\!\subset\!X$ and $\mu\!\in\!H^2(X;\Z_2)$.
\begin{enumerate}[label=(\alph*),leftmargin=*]

\item\label{RelPinObs_it} The vector bundle $V$ admits 
a relative $\Pin^-$-structure (resp.~$\Pin^+$-structure) $\fp$ with $w_2(\fp)\!=\!\mu$
if and only if $\mu|_Y\!=\!w_2(V)\!-\!w_1^2(V)$ (resp.~$\mu|_Y\!=\!w_2(V)$).

\item\label{RelSpinObs_it} If\, $V$ is orientable, $V$ admits 
a relative $\Spin$-structure~$\fs$  with $w_2(\fs)\!=\!\mu$
if and only if \hbox{$\mu|_Y\!=\!w_2(V)$}.

\end{enumerate}
\end{RelSpinPin}

\begin{RelSpinPin}[Affine Structure]\label{RelSpinPinStr_prop}
Let $V$ be a vector bundle over $Y\!\subset\!X$.
\begin{enumerate}[label=(\alph*),leftmargin=*]

\item\label{RelPinStr_it}  If\, $V$ admits a relative $\Pin^{\pm}$-structure,
the group $H^2(X,Y;\Z_2)$ acts freely and transitively on the set~$\cP_X^{\pm}(V)$.

\item\label{RelSpinStr_it}  If $\fo\!\in\!\fO(V)$ and $(V,\fo)$
admits a relative $\Spin$-structure,
the group $H^2(X,Y;\Z_2)$ acts freely and transitively on the set~$\Sp_X(V,\fo)$.

\end{enumerate}
Both actions are natural with respect to pullbacks by continuous maps and
\BE{RelSpinPinStr_e0} w_2(\eta\!\cdot\!\fp)=w_2(\fp)\!+\!\eta|_X, \qquad
w_2(\eta\!\cdot\!\fs)=w_2(\fs)\!+\!\eta|_X\EE
for all $\eta\!\in\!H^2(X,Y;\Z_2)$, $\fp\!\in\!\cP_X^{\pm}(V)$, and $\fs\!\in\!\Sp_X(V,\fo)$.
If $V'$ and $V''$ are vector bundles over~$Y$, with at most one of them possibly of rank~0, 
then the action of the automorphism~$\Psi$ in~\eref{SpinPinStr_e0} 
on $\cP_X^{\pm}(V'\!\oplus\!V'')$ is given~by 
\BE{RelSpinPinStr_e}  \Psi^*\fp= 
\de_{X,Y}\!\big((\rk\,V'\!-\!1)w_1(V')\!+\!(\rk\,V')w_1(V'')\big)\!\cdot\!\fp
\qquad\forall\,\fp\!\in\!\cP_X^{\pm}(V'\!\oplus\!V'').\EE

\end{RelSpinPin}

\begin{RelSpinPin}[Orientation Reversal]\label{RelOSpinRev_prop}
Let $V$ be a real vector bundle over $Y\!\subset\!X$. 
There is a natural $H^2(X,Y;\Z_2)$-equivariant 
involution 
\BE{RelOSpinRev_e}\OSp_X(V)\lra \OSp_X(V), \qquad \os\lra\ov\os,\EE
so that $w_2(\ov\os)\!=\!w_2(\os)$ and 
$\ov\os\!\in\Sp_X(V,\ov\fo)$ for all $\os\!\in\!\Sp_X(V,\fo)$ 
and $\fo\!\in\!\fO(V)$.
\end{RelSpinPin}

\begin{RelSpinPin}[Reduction]\label{RelPin2SpinRed_prop}
Let $V$ be a real vector bundle over $X\!\subset\!Y$.
For every $\fo\!\in\!\fO(V)$, there are natural 
$H^2(X,Y;\Z_2)$-equivariant bijections
\BE{RelPin2SpinRed_e} \fR_{\fo}^{\pm}\!:\cP_X^{\pm}(V)\lra \Sp_X(V,\fo)\EE
so that $w_2(\fR_{\fo}^{\pm}(\cdot))\!=\!w_2(\cdot)$ and 
$\fR_{\ov\fo}^{\pm}(\cdot)\!=\!\ov{\fR_{\fo}^{\pm}(\cdot)}$.
\end{RelSpinPin}

\begin{RelSpinPin}[Stability]\label{RelSpinPinStab_prop}
Let $V$ be a  vector bundle over $Y\!\subset\!X$. 
There are natural $H^2(X,Y;\Z_2)$-equivariant bijections
\BE{RelSpinPinStab_e}
\St_V\!: \OSp_X(V)\lra \OSp_X\big(\tau_Y\!\oplus\!V\big), \qquad
\St_V^{\pm}\!:\cP_X^{\pm}(V)\lra \cP_X^{\pm}\big(\tau_Y\!\oplus\!V\big)\EE
so that $w_2(\St_V(\cdot))\!=\!w_2(\cdot)$,
\BE{RelSpinPinStab_e2}
\St_V(\ov\os)=\ov{\St_V(\os)}~~\forall\,\os\!\in\!\OSp_X(V), \quad
\St_V\!\circ\!\fR_{\fo}^{\pm}\!=\!\fR_{\St_V(\fo)}^{\pm}\!\circ\!\St_V^{\pm}
~~\forall\,\fo\!\in\!\fO(V).\EE 
\end{RelSpinPin}

\begin{RelSpinPin}[Correspondences]\label{RelSpinPinCorr_prop}
Let $V$ be a real vector bundle over $Y\!\subset\!X$. 
There are natural $H^2(X,Y;\Z_2)$-equivariant bijections
\BE{RelSpinPinCorr_e}\Co_V^{\pm}\!:\cP_X^{\pm}(V)\lra\Sp_X\!\big(V_{\pm},\fo_V^{\pm}\big)\EE
so that $w_2\!\circ\!\Co_V^{\pm}\!=\!w_2$ and
$\Co_{\tau_Y\oplus V}^{\pm}\!\circ\!\St_V^{\pm}\!=\!\St_{V_{\pm}}^{\pm}\!\circ\!\Co_V^{\pm}$. 
\end{RelSpinPin}

\begin{RelSpinPin}[Exact Triples]\label{RelSpinPinSES_prop}
Every short exact sequence~$\ce$ of vector bundles over $Y\!\subset\!X$ as in~\eref{SpinPinSES_e0}
determines natural  $H^2(X,Y;\Z_2)$-biequivariant maps
\BE{RelSpinPinSESdfn_e0}\begin{split}
\llrr{\cdot,\cdot}_{\ce}\!:\OSp_X(V')\!\times\!\OSp_X(V'')&\lra\OSp_X(V), \\
\llrr{\cdot,\cdot}_{\ce}\!:\OSp_X(V')\!\times\!\cP_X^{\pm}(V'')&\lra\cP_X^{\pm}(V)
\end{split}\EE
so that the following properties hold. 
\begin{enumerate}[label=(ses\arabic*$_X$),leftmargin=*]

\stepcounter{enumi}

\item\label{RelDSorient_it} If $\ce$ is as in~\eref{SpinPinSES_e0},
$\fo'\!\in\!\fO(V')$, and $\fo''\!\in\!\fO(V'')$, then
\begin{alignat}{2}
\label{RelDSorient1_e}
\llrr{\os',\os''}_{\ce}&\in\OSp_X(V,\fo'_{\ce}\fo'') &\quad
&\forall~\os'\!\in\!\Sp_X(V',\fo'),\,\os''\!\in\!\Sp_X(V'',\fo''),\\
\label{RelDSorient2_e}
\fR_{\fo'_{\ce}\fo''}^{\pm}\big(\llrr{\os',\fp''}_{\ce}\big)
&=\bllrr{\os',\fR_{\fo''}^{\pm}(\fp'')}_{\ce} &\quad
&\forall~\os'\!\in\!\Sp_X(V',\fo'),\,\fp''\!\in\!\cP_X^{\pm}(V'').
\end{alignat}

\item\label{RelDSEquivSum_it} 
If $\ce$ is as in~\eref{SpinPinSES_e0},  $\fo'\!\in\!\fO(V')$, $\os'\!\in\!\Sp_X(V',\fo')$, and
 $\fp''\!\in\!\cP_X^{\pm}(V'')$, then
\begin{gather}\notag
\llrr{\ov\os',\fp''}_{\ce}=\de_{X,Y}\big(w_1(V'')\!\big)\!\cdot\!\llrr{\os',\fp''}_{\ce}\,,
\quad
\Co_V^{\pm}\big(\llrr{\os',\fp''}_{\ce}\big)=
\bllrr{\os',\Co_{V''}^{\pm}(\fp'')}_{\ce_{\fo'}^{\pm}}\,,\\
\label{Relw2ses_e} w_2\big(\llrr{\os',\fp''}_{\ce}\big)=w_2(\os')\!+\!w_2(\fp'').
\end{gather}

\item\label{RelDSassoc_it} 
If $V_1',V_2',V''$ are vector bundles over~$Y$, then
$$\bllrr{\os_1',\llrr{\os_2',\fp''}_{\oplus}\!}_{\oplus}
=\bllrr{\!\llrr{\os_1',\os_2'}_{\oplus},\fp''}_{\oplus}$$
for all $\os_1'\!\in\!\OSp_X(V_1')$, $\os_2'\!\in\!\OSp_X(V_2')$,  and
$\fp''\!\in\!\cP_X^{\pm}(V'')$.

\item\label{RelDSstab_it} 
If $V$ is a vector bundle over~$Y$ and $\fp\!\in\!\cP_X^{\pm}(V)$,
$$\St_V^{\pm}(\fp) =\llrr{\io_X(\os_0(\tau_Y,\fo_Y)\!),\fp}_{\oplus} .$$

\item\label{RelDSPin2Spin_it} 
If $V$ is a vector bundle over~$Y$, $\fo\!\in\!\fO(V)$, 
and $\fp\!\in\!\cP_X^{\pm}(V)$, then
$$\Co_V^{\pm}(\fp)=
\bllrr{\fR_{\fo}^{\pm}(\fp),\io_X\big(\os_0((2\!\pm\!1)\la(V,\fo))\!\big)\!}_{\oplus}\,.$$

\end{enumerate}

\end{RelSpinPin}

\vspace{.1in}

By RelSpinPin~\ref{RelPin2SpinRed_prop}, \eref{RelDSorient2_e},
and the first and third statements in RelSpinPin~\ref{RelSpinPinSES_prop}\ref{RelDSEquivSum_it},
\BE{RelOspinDSEquivSum_e}  
\llrr{\ov\os',\os''}_{\ce}=\llrr{\os',\ov\os''}_{\ce}=\ov{\llrr{\os',\os''}_{\ce}},
\quad w_2\big(\llrr{\os',\os''}_{\ce}\big)=w_2(\os')\!+\!w_2(\os'')\EE
for every short exact sequence~$\ce$ of vector bundles as in~\eref{SpinPinSES_e0},
$\os'\!\in\!\OSp_X(V')$, and $\os''\!\in\!\OSp_X(V'')$.
By RelSpinPin~\ref{RelPin2SpinRed_prop} and~\ref{RelSpinPinSES_prop}\ref{RelDSassoc_it} 
and~\eref{RelDSorient2_e},
\BE{RelOspinDSassoc_e}  
\bllrr{\os_1',\llrr{\os_2',\os''}_{\oplus}\!}_{\oplus}
=\bllrr{\!\llrr{\os_1',\os_2'}_{\oplus},\os''}_{\oplus}\EE
for all vector bundles $V_1',V_2',V''$ over~$Y$, 
$\os_1'\!\in\!\OSp_X(V_1')$, $\os_2'\!\in\!\OSp_X(V_2')$,  and $\os''\!\in\!\OSp_X(V'')$.
By RelSpinPin~\ref{RelPin2SpinRed_prop}, \ref{RelSpinPinStab_prop}, 
and~\ref{RelSpinPinSES_prop}\ref{RelDSstab_it}
and \eref{RelDSorient2_e}, 
\BE{RelOspinDSstab_e}\St_V(\os) =\llrr{\io_X(\os_0(\tau_Y,\fo_Y)\!),\os}_{\oplus}\EE
for every vector bundle~$V$ over~$Y$ and $\os\!\in\!\OSp_X(V)$.\\

Combining the second map in~\eref{RelSpinPinSESdfn_e0} with the canonical isomorphism
of $V'\!\oplus\!V''$ with $V''\!\oplus\!V'$, we obtain
a natural $H^2(X,Y;\Z_2)$-biequivariant map
\BE{RelSpinPinSESdfn_e0b}
\llrr{\cdot,\cdot}_{\ce}\!:\cP_X^{\pm}(V')\!\times\!\OSp_X(V'')\lra\cP_X^{\pm}(V)\,.\EE
By the RelSpinPin~\ref{RelSpinPinSES_prop} property, this map satisfies the obvious
analogues of \eref{RelDSorient2_e}, 
the first and third statements in RelSpinPin~\ref{RelDSEquivSum_it}, 
and RelSpinPin~\ref{RelDSassoc_it} and~\ref{RelDSPin2Spin_it}.

\begin{RelSpinPin}[Compatibility with Spin/Pin-Structures]\label{vsSpinPin_prop}
For every real vector bundle~$V$ over $Y\!\subset\!X$, 
the maps~\eref{vsSpinPin_e0} are equivariant with respect to the homomorphism~\eref{YdeXYZ2_e},
respect all structures and correspondences of 
SpinPin~\ref{OSpinRev_prop}-\ref{SpinPinSES_prop}
and RelSpinPin~\ref{RelOSpinRev_prop}-\ref{RelSpinPinSES_prop},
\BE{vsSpinPin_e}
\io_X\big(\cP^{\pm}(V)\big)=\big\{\fp\!\in\!\cP_X^{\pm}(V)\!: w_2(\fp)\!=\!0\big\}.\EE
\end{RelSpinPin}

\vspace{.1in}

By RelSpinPin~\ref{RelPin2SpinRed_prop} and~\ref{vsSpinPin_prop},
\BE{vsSpinPin2_e}
\io_X\big(\OSp(V)\big)=\big\{\os\!\in\!\OSp_X(V)\!: w_2(\os)\!=\!0\big\}.\EE

\subsection{Proof of Theorem~\ref{RelSpinStrEquiv_thm}\ref{RelPinSpinPrp_it}:  
Definition~\ref{RelPinSpin_dfn3} perspective}
\label{RelSpinPindfn3_subs}

We now establish the statements of Section~\ref{RelSpinPinProp_subs} for the notions of
relative  $\Spin$-structure and $\Pin^{\pm}$-structure arising from Definition~\ref{RelPinSpin_dfn3}.
Similarly to Section~\ref{SpinPindfn3_subs}, 
the proofs of these statements
rely heavily on the topological observations of Section~\ref{TopolPrelim_subs}.
Throughout this section, the terms relative $\Spin$-structure and $\Pin^{\pm}$-structure
refer to the notions arising from Definition~\ref{RelPinSpin_dfn3}.

\begin{proof}[{\bf{\emph{Proof of RelSpinPin~\ref{RelSpinPinStr_prop} property}}}]
We construct natural actions of $H^2(X,Y;\Z_2)$ on $\OSp_X(V)$ and $\cP_X^{\pm}(V)$
so~that 
\BE{RelSpinPinStr_e1}
\io_X\big(\eta\!\cdot\!\os\big)=\de_{X,Y}(\eta)\!\cdot\!\io_X(\os) 
\qquad\hbox{and}\qquad
\io_X\big(\eta\!\cdot\!\fp\big)=\de_{X,Y}(\eta)\!\cdot\!\io_X(\fp) \EE
for all $\eta\!\in\!H^1(Y;\Z_2)$, $\os\!\in\!\OSp(V)$, and $\fp\!\in\!\cP^{\pm}(V)$.
By the definition of relative $\Spin$- and $\Pin^{\pm}$-structures in this perspective,
it is sufficient to describe the first action 
and to establish \eref{RelSpinPinStr_e1} and 
the claims of the RelSpinPin~\ref{RelSpinPinStr_prop} property, 
other than~\eref{RelSpinPinStr_e},
for $\Spin$-structures on oriented vector bundles~$(V,\fo)$
with \hbox{$n\!\equiv\!\rk\,V\!\ge\!3$}.\\

The oriented vector bundle $\{\prt u\}^*(V,\fo)$ has two mod~2 homotopy classes of trivializations
for every \hbox{$u\!\in\!\cL_X(Y)$}.
Let $\fs$ be a relative $\Spin$-structure on~$(V,\fo)$ and $\eta\!\in\!H^2(X,Y;\Z_2)$.
We define the relative $\Spin$-structure \hbox{$\eta\!\cdot\!\fs\!\equiv\!(\eta\!\cdot\!\fs_u)_u$}
on~$(V,\fo)$ by 
\BE{RelSpinPinStr_e33}
\eta\!\cdot\!\fs_u=\begin{cases}=\fs_u,&\hbox{if}~\lr{u^*\eta,[\Si]_{\Z_2}}\!=\!0;\\
\neq\fs_u,&\hbox{if}~\lr{u^*\eta,[\Si]_{\Z_2}}\!\neq\!0;
\end{cases}
\qquad\forall\,\big(u\!:(\Si,\prt\Si)\!\lra\!(X,Y)\!\big)\in\!\cL_X(Y).\EE
Suppose $F\!:\Si\!\lra\!Y$ is a continuous map from a bordered surface relatively equivalent
to the disjoint union of continuous maps
\BE{RelSpinPinStr_e5}u_i\!:(\Si_i,\prt\Si_i)\lra(X,Y), \qquad i=1,\ldots,m,\EE
from bordered surfaces.
In particular,
\begin{gather*}
{u_1}_*([\Si_1]_{\Z_2})+\!\ldots\!+{u_m}_*([\Si_m]_{\Z_2})=0\in H_2(X,Y;\Z_2), \\
\sum_{i=1}^m\blr{u_i^*\eta,[\Si_i]_{\Z_2}}=
\blr{\eta,\sum_{i=1}^m{u_i}_*([\Si_i]_{\Z_2})}=0\in\Z_2.
\end{gather*}
Thus, the number of the maps~$u_i$ such that $\eta\!\cdot\fs_{u_i}\!\neq\!\fs_{u_i}$ is even.
Since the disjoint union of the mod~2 homotopy classes~$\fs_{u_i}$ of trivializations 
of~$\{\prt u_i\}^*(V,\fo)$ equals to
the mod~2 homotopy class $\fs_{\Si}(F^*(V,\fo)\!)$ of trivializations of 
$\{\prt F\}^*(V,\fo)$, the same is the case for the disjoint union of 
the mod~2 homotopy classes~$\eta\!\cdot\fs_{u_i}$.
Thus, $\eta\!\cdot\!\fs$ is indeed a relative $\Spin$-structure on~$(V,\fo)$.\\

It is immediate that the above construction defines a group action of $H^2(X,Y;\Z_2)$
on the set of such structures.
In light of~\eref{SpinPinStr_e33}, this group satisfies 
the first property in~\eref{RelSpinPinStr_e1}.
If \hbox{$u\!:\Si\!\lra\!X$} is a continuous map from a closed surface,
$\fs_u\!=\!\eta\!\cdot\!\fs_u$ if and only if 
$$\blr{\eta|_X,u_*[\Si]_{\Z_2}}=0\in\Z_2.$$
Along with~\eref{w2fsdfn_e}, Lemma~\ref{Z2cycles_lmm}\ref{cycles_it} with $Y\!=\!\eset$, and 
the Universal Coefficient Theorem for Cohomology, 
this implies the second statement in~\eref{RelSpinPinStr_e0}.\\

By Lemma~\ref{Z2cycles_lmm}\ref{cycles_it}, every $b\!\in\!H_2(X,Y;\Z_2)$ 
can be represented by a continuous map
\hbox{$(\Si,\prt\Si)\!\lra\!(X,Y)$} from a bordered surface.
By the Universal Coefficient Theorem for Cohomology, the homomorphism
$$\ka\!:H^2(X,Y;\Z_2)\lra\Hom\big(H_2(X,Y),\Z_2\big), \quad
\big\{\ka(\eta)\big\}(b)=\lr{\eta,b},$$
is bijective.
Thus, \hbox{$\eta\!\cdot\!\fs_u\!\neq\!\fs_u$} for some $u\!\in\!\cL_X(Y)$ if 
$\eta\!\neq\!0$ and so the action of $H^2(X,Y;\Z_2)$ above is~free.\\

Suppose $\fs'\!\equiv\!(\fs_u')_u$ is another relative $\Spin$-structure on~$(V,\fo)$.
Define
$$\eta\!:\cL_X(Y)\lra\Z_2, \qquad 
\eta(u)=\begin{cases}0,&\hbox{if}~\fs_u'\!=\!\fs_u;\\
1,&\hbox{if}~\fs_u'\!\neq\!\fs_u.\end{cases}$$
This determines a linear map from the $\Z_2$-vector space generated by $\cL_X(Y)$ 
to~$\Z_2$. 
Suppose \hbox{$u_1,\ldots,u_m\!\in\!\cL_X(Y)$} and $F$ are as in and above~\eref{RelSpinPinStr_e5}. 
By Definition~\ref{RelPinSpin_dfn3}, 
$$\fs_{u_1}\!\ldots\!\fs_{u_m}=\fs_{\Si}\!\big(F^*(V,\fo)\!\big)=\fs_{u_1}'\!\ldots\!\fs_{u_m}'\,.$$
Thus, the number of maps~$u_i$ such that $\fs_{u_i}\!\neq\!\fs_{u_i}$ is even and~so
$$\sum_{i=1}^m\eta(u_i)=0\in\Z_2\,.$$
Along with~\ref{cycles_it} and~\ref{relbound_it} in Lemma~\ref{Z2cycles_lmm},
this implies that~$\eta$ descends to a homomorphism
$$\eta\!: H_2(X,Y;\Z_2)\lra \Z_2\,.$$
By the Universal Coefficient Theorem for Cohomology,
such a homomorphism corresponds to an element of $H^2(X,Y;\Z_2)$,
which we still denote by~$\eta$.
By the definition of~$\eta$ and the construction above, $\eta\!\cdot\!\fs\!=\!\fs'$.
Thus, the action of $H^2(X,Y;\Z_2)$ described above is transitive.\\

It remains to establish~\eref{RelSpinPinStr_e}.
Let 
$$\eta=(\rk\,V'\!-\!1)w_1(V')\!+\!(\rk\,V')w_1(V'')\in H^1(Y;\Z_2) \qquad\hbox{and}\qquad
 \fp\!\in\!\cP_X^{\pm}(V'\!\oplus\!V'')\,.$$
Let $u\!:(\Si,\prt\Si)\!\lra\!(X,Y)$ be a continuous map from a bordered surface and
$\prt_1\Si,\ldots,\prt_k\Si$ be the components of~$\prt\Si$.
By~\eref{SpinPinStr_e} and the SpinPin~\ref{SpinPinStr_prop}\ref{PinStr_it} property, 
\hbox{$(\Psi^*\fp)_u\!\equiv\!\Psi^*\fp_u$} equals~$\fp_u$
if and only~if
$$\sum_{r=1}^k\!\blr{\{u|_{\prt_r\Si}\}^*\eta,[\prt_r\Si]_{\Z_2}}
\equiv \blr{u^*(\de_{X,Y}\eta),[\Si]_{\Z_2}}\in\Z_2$$
vanishes. This implies~\eref{RelSpinPinStr_e}
\end{proof}

Lemma~\ref{X2VBdisconn_lmm} below is a reformulation of 
Corollaries~\ref{X2VB_crl1b}-\ref{X2VB_crl2b} 
for (possibly) disconnected surfaces. 
It is used in the proof of the RelSpinPin~\ref{RelSpinPinObs_prop} property.
If $Z_-$ and~$Z_+$ are compact one-dimensional manifolds, 
$(\wt{V},\wt\fo)$ is an oriented vector bundle over~$Z_-\!\sqcup\!Z_+$,
\BE{ZZvphdfn_e}\vph\!:\big(\wt{V},\wt\fo\big)|_{Z_-}\lra \big(\wt{V},\wt\fo\big)|_{Z_+}\EE
is an isomorphism covering an identification of~$Z_-$ with~$Z_+$,
and $\fs_-$ and~$\fs_+$ are mod~2 homotopy classes of trivializations of
the two sides in~\eref{ZZvphdfn_e}, 
we define
$$\lr{\fs_-,\fs_+}_{\vph}=\begin{cases}0\!\in\!\Z_2,&\hbox{if}~\fs_-\!=\!\vph^*\fs_+;\\
1\!\in\!\Z_2,&\hbox{if}~\fs_-\!\neq\!\vph^*\fs_+.
\end{cases}$$
If $Z_-\!=\!Z_+$, let 
$$\lr{\fs_-,\fs_+}=\lr{\fs_-,\fs_+}_{\id}\,.$$

\begin{lmm}\label{X2VBdisconn_lmm}
Suppose $\wt\Si$ is a bordered surface, 
$Z_-$ and~$Z_+$ are disjoint unions of components of~$\prt\wt\Si$, 
$(\wt{V},\wt\fo)$ is an oriented vector bundle over~$\wt\Si$ with
\hbox{$\rk_{\R}\wt{V}\!\ge\!3$}, and
\eref{ZZvphdfn_e} is an isomorphism covering an identification of~$Z_-$ with~$Z_+$.
Let $(\wh{V},\wh\fo)$ be the oriented vector bundle over the closed surface~$\wh\Si$
obtained from $(\wt{V},\wt\fo)$ and~$\wt\Si$  by gluing via~$\vph$.
If $\fs_-$, $\fs_+$, and~$\fs_0$ are mod~2 homotopy classes of trivializations of $(\wt{V},\wt\fo)$ over
$Z_-$, $Z_+$, and $\prt\Si\!-\!Z_-\!\cup\!Z_+$, respectively, then
$$\blr{\fs_0,\fs_{\wh\Si}(\wh{V},\wh\fo)}=
\blr{\fs_-,\fs_+}_{\vph}\!+\!\blr{\fs_0\fs_-\fs_+,\fs_{\wt\Si}(\wt{V},\wt\fo)}
\in\Z_2\,.$$
\end{lmm}

\begin{proof}[{\bf{\emph{Proof of RelSpinPin~\ref{RelSpinPinObs_prop} property}}}]
By the definition of relative $\Spin$- and $\Pin^{\pm}$-structures in this perspective,
it is sufficient to establish this claim for relative $\Spin$-structures
on oriented vector bundles~$(V,\fo)$ with $\rk\,V\!\ge\!3$.\\

Suppose $(V,\fo)$ admits a relative $\Spin$-structure $\fs\!\equiv\!(\fs_u)_u$
with $w_2(\fs)\!=\!\mu$.
Let \hbox{$u\!:\Si\!\lra\!Y$} be a continuous map \hbox{$u\!:\Si\!\lra\!Y$} from 
a closed surface.
By~\eref{w2fsdfn_e} and  the condition of Definition~\ref{RelPinSpin_dfn3} applied with 
$u_1,F\!=\!u$,
$$\blr{\mu|_Y,u_*[\Si]_{\Z_2}}=\fs_u=\fs_{\Si}\big(u^*(V,\fo)\!\big)
\equiv\blr{w_2(u^*V),[\Si]_{\Z_2}}=\blr{w_2(V),u_*[\Si]_{\Z_2}}.$$
Along with Lemma~\ref{Z2cycles_lmm}\ref{cycles_it} with $(X,Y)\!=\!(Y,\eset)$ and 
the Universal Coefficient Theorem for Cohomology, 
this implies that $\mu|_Y\!=\!w_2(V)$.\\

Suppose $w_2(V)\!=\!\mu|_Y$.
Choose a collection 
$$\cC\equiv\big\{u_i\!:(\Si_i,\prt\Si_i)\!\lra\!(X,Y)\big\}$$ 
of maps that form a basis for $H_2(X,Y;\Z_2)$ and
a mod~2 homotopy class~$\fs_i$ of trivializations of $\{\prt u_i\}^*(V,\fo)$ 
for each map in~$\cC$.
Given a continuous map $u\!:(\Si,\prt\Si)\!\lra\!(X,Y)$ from a bordered surface, 
let $u_1,\ldots,u_k\!\in\!\cC$ be so~that 
$$u_*\big([\Si]_{\Z_2}\big)\!+\!{u_1}_*\big([\Si_1]_{\Z_2}\big)
\!+\!\ldots\!+\!{u_k}_*\big([\Si_k]_{\Z_2}\big)=\!0 \in H_2(X,Y;\Z_2).$$ 
By Lemma~\ref{Z2cycles_lmm}\ref{relbound_it},  
there exists a continuous map $F\!:\Si_0\!\lra\!Y$ 
from a bordered surface relatively equivalent to the disjoint union of the maps $u,u_1,\ldots,u_k$.
In particular,
$$\prt F=\prt u\sqcup \prt u_1\sqcup\ldots\sqcup \prt u_k\,.$$
We take the mod~2 homotopy class~$\fs_u$ of trivializations of $\{\prt u\}^*(V,\fo)$ to be
such~that 
\BE{RelSpinPinObs_e2}  \fs_u\fs_1\ldots\fs_k=\fs_{\Si_0}\!\big(F^*(V,\fo)\!\big)\,.\EE

\vspace{.10in}

Suppose $F'\!:\Si_0'\!\lra\!Y$ is another continuous map satisfying 
the conditions of the previous paragraph and $\fs_u'$ is 
the mod~2 homotopy class of trivializations of $\{\prt u\}^*(V,\fo)$ such~that 
\BE{RelSpinPinObs_e2b}   \fs_u'\fs_1\ldots\fs_k=\fs_{\Si_0'}\!\big(F'^*(V,\fo)\!\big)\,.\EE
Denote by $\wh\Si$ (resp.~$\wt\Si$) the closed (resp.~bordered) surface 
obtained from~$\Si_0$ and~$\Si_0'$ by identifying them along 
the boundary components corresponding to $\prt\Si,\prt\Si_1,\ldots,\prt\Si_k$
(resp.~$\prt\Si_1,\ldots,\prt\Si_k$).
Thus, the boundary components of $\wt\Si$ split into those forming 
$\prt\Si\!\subset\!\Si_0$ and 
those forming  \hbox{$\prt\Si\!\subset\!\Si_0'$};
$\wh\Si$ is obtained from $\wt\Si$ by identifying these two copies of~$\prt\Si$.
The~maps~$F$ and~$F'$ induce continuous~maps
$$\wh{F}\!:\wh\Si\lra Y \qquad\hbox{and}\qquad \wt{F}\!:\wt\Si\lra Y,$$
which restrict to $F$ and~$F'$ over $\Si_0,\Si_0'\!\!\subset\!\wh\Si,\wt\Si$.
By~\eref{RelSpinPinObs_e2}, \eref{RelSpinPinObs_e2b}, and Lemma~\ref{X2VBdisconn_lmm}, 
\BE{RelSpinPinObs_e4} \fs_u\fs_u'=\fs_{\wt\Si}\!\big(\wt{F}^*(V,\fo)\!\big)\,.\EE
Since each of the maps $F$ and $F'$ is relatively equivalent 
to \hbox{$u\!\sqcup\!u_1\!\ldots\!\sqcup\!u_k$},
$\wh{F}_*[\wh\Si]_{\Z_2}$ vanishes in $H_2(X;\Z_2)$.
Thus,
\BE{RelSpinPinObs_e9}
\fs_{\wh\Si}\big(\wh{F}^*(V,\fo)\!\big)\equiv 
\blr{w_2(\wh{F}^*V),[\wh\Si]_{\Z_2}}
=\blr{\wh{F}^*(\mu|_Y),[\wh\Si]_{\Z_2}}
=\blr{\mu|_Y,\wh{F}_*[\wh\Si]_{\Z_2}}=0\,.\EE
The bundle $\wh{F}^*V$ over $\wh\Si$ is obtained by identifying the two copies of 
$\{\prt u\}^*V$ in~$\wt{F}^*V$.
Along with~\eref{RelSpinPinObs_e4}, \eref{RelSpinPinObs_e9}, and Lemma~\ref{X2VBdisconn_lmm}, 
this implies that $\fs_u\!=\!\fs_u'$.\\

We next verify that $\fs\!\equiv\!(\fs_u)_u$ satisfies 
the condition of Definition~\ref{RelPinSpin_dfn3}.
Suppose $F\!:\Si\!\lra\!Y$ is a continuous map from a bordered surface relatively equivalent
to the disjoint union of continuous~maps
$$u_i'\!:(\Si_i',\prt\Si_i')\lra(X,Y), \qquad i=1,\ldots,m,$$
from bordered surfaces.
Thus,
\BE{RelSpinPinObs_e14a} \prt F=\prt u_1'\sqcup\ldots\sqcup\prt u_m'\,.\EE
We show below that 
\BE{RelSpinPinObs_e15}\fs_{u_1'}\!\ldots\!\fs_{u_m'}=\fs_{\Si}\big(F^*(V,\fo)\!\big),\EE
thus confirming the condition of Definition~\ref{RelPinSpin_dfn3}.\\

For each $i\!=\!1,\ldots,m$, let 
$$u_{ij}\!:\big(\Si_{ij},\prt\Si_{ij}\big)\!\lra\!(X,Y)~~\hbox{with}~j\!=\!1,\ldots,k_i
\qquad\hbox{and}\qquad F_i\!:\Si_{i0}\lra Y$$
be the maps in~$\cC$ and a continuous map from a bordered surface, respectively,
as in the construction of the mod~2 homotopy class~$\fs_{u_i'}$ of trivializations
of $\{\prt u_i'\}^*(V,\fo)$ above.
Let $\fs_{ij}$ be the initially chosen mod~2 homotopy class of trivializations 
of $\prt u_{ij}^*(V,\fo)$.
Thus,
\BE{RelSpinPinObs_e14b} 
\prt F_i=\prt u_i'\sqcup \prt u_{i1}\sqcup\ldots\sqcup \prt u_{ik_i}, \qquad
\fs_{u_i'}\fs_{i1}\!\ldots\!\fs_{ik_i}=\fs_{\Si_{i0}}\!\big(F_i^*(V,\fo)\!\big).\EE

\vspace{.1in}

Since the maps in~$\cC$ are linearly independent in~$H_2(X,Y;\Z_2)$ and 
$$\sum_{i=1}^m\sum_{j=1}^{k_i}\!{u_{ij}}_*\big([\Si_{ij}]_{\Z_2}\big)
=\sum_{i=1}^m\!{u_i'}_*\big([\Si'_i]_{\Z_2}\big)=F_*\big([\Si]_{\Z_2}\big)=0
\in H_2(X,Y;\Z_2),$$
the maps $u_{ij}$ come in pairs.
Along with the first identity in~\eref{RelSpinPinObs_e14b}, this implies 
that the boundary components of $\Si_{10},\ldots,\Si_{m0}$ corresponding to 
the boundaries~$\prt u_{ij}$ of these maps come in pairs as well.
By~\eref{RelSpinPinObs_e14a}, the remaining boundary components of $\Si_{10},\ldots,\Si_{m0}$ 
correspond to the boundary components of~$\Si$.
Denote by $\Si'$ the bordered surface obtained from $\Si_{10},\ldots,\Si_{m0}$ 
by identifying them along  the paired up boundary components corresponding to~$\prt u_{ij}$.
Let $\wh\Si$ be the closed surface obtained from~$\Si$ and~$\Si'$ by identifying them along
the boundary components corresponding to $\prt u_1',\ldots,\prt u_m'$.\\

The~maps~$F,F_1,\ldots,F_m$ induce continuous~maps
\BE{RelSpinPinObs_e18a}F'\!:\Si'\lra Y \qquad\hbox{and}\qquad \wh{F}\!:\wh\Si\lra Y,\EE
which restrict to $F_i$ over $\Si_{0i}\!\subset\!\Si'$;
the second map also restricts to~$F$ over $\Si\!\subset\!\wh\Si$.
By the second identity in~\eref{RelSpinPinObs_e14b} and Lemma~\ref{X2VBdisconn_lmm},
\BE{RelSpinPinObs_e18}\fs_{u_1'}\!\ldots\!\fs_{u_m'}=\fs_{\Si'}\big(F'^*(V,\fo)\!\big), \quad
\fs_{\Si}\big(F^*(V,\fo)\!\big)\big(\fs_{u_1'}\!\ldots\!\fs_{u_m'}\big)=
\fs_{\Si\sqcup\Si'}\!\big(\{F\!\sqcup\!F'\}^*(V,\fo)\!\big)\,.\EE
Since the disjoint union of~$F$ and the maps~$u_{ij}$ (which come in pairs)
is relatively equivalent to $F_1\!\sqcup\!\ldots\!\sqcup\!F_m$,
$\wh{F}_*[\wh\Si]_{\Z_2}$ vanishes in $H_2(X;\Z_2)$.
Thus,
\BE{RelSpinPinObs_e19}
\fs_{\wh\Si}\big(\wh{F}^*(V,\fo)\!\big)\equiv 
\blr{w_2(\wh{F}^*V),[\wh\Si]_{\Z_2}}=\blr{\wh{F}^*(\mu|_Y),[\wh\Si]_{\Z_2}}
=\blr{\mu|_Y,\wh{F}_*[\wh\Si]_{\Z_2}}=0\,.\EE
The bundle $\wh{F}^*V$ over $\wh\Si$ is obtained by identifying the copies of 
$\{\prt u_1\}^*V,\ldots,\{\prt u_m\}^*V$  in $\{F\!\sqcup\!F'\}^*V$.
Along with the second identity in~\eref{RelSpinPinObs_e18}, 
\eref{RelSpinPinObs_e19}, and Lemma~\ref{X2VBdisconn_lmm}, 
this implies~\eref{RelSpinPinObs_e15}.\\

Thus, $\fs$ is a relative $\Spin$-structure on the oriented vector bundle $(V,\fo)$
over $Y\!\subset\!X$ and~so \hbox{$w_2(\fs)|_Y\!=\!w_2(V)$} by the first paragraph of
the proof.
Along with the cohomology exact sequence for $(X,Y)$, this implies~that 
$$\mu=w_2(\fs)\!+\!\eta|_X\in H^2(X;\Z_2)$$
for some $\eta\!\in\!H^2(X,Y;\Z_2)$.
By the RelSpinPin~\ref{RelSpinPinStr_prop} property, the relative $\Spin$-structure 
$\eta\!\cdot\!\fs$ on~$(V,\fo)$ satisfies $w_2(\eta\!\cdot\!\fs)\!=\!\mu$.
\end{proof}

\begin{proof}[{\bf{\emph{Proof of 
RelSpinPin~\ref{RelOSpinRev_prop},\ref{RelPin2SpinRed_prop} properties}}}]
For the purposes of establishing the first of these properties,
we can assume that \hbox{$n\!\equiv\!\rk\,V\!\ge\!3$}.
For $\fo\!\in\!\fO(V)$ and $\fs\!\in\!\Sp_X(V,\fo)$,
we define the relative $\Spin$-structure \hbox{$\ov\fs\!\equiv\!(\fs_u)_u$}
on $(V,\ov\fo)$ by
\BE{RelOSpinRev_e33}\ov\fs_u=\begin{cases}
\{\ov\phi\!\equiv\!\bI_{n;1}\phi\!:\phi\!\in\!\fs_u\},
&\hbox{if}~\prt\Si\!\neq\!\eset;\\
\fs_u,&\hbox{if}~\prt\Si\!=\!\eset;
\end{cases}
\quad\forall\,\big(u\!:(\Si,\prt\Si)\!\lra\!(X,Y)\!\big)\!\in\!\cL_X(Y).\EE
Since this orientation-reversal operation commutes with the disjoint union operation on
the mod~2 homotopy classes of trivializations over one-dimensional manifolds and
$$\ov{\fs_{\Si}(V',\fo')}=\fs_{\Si}(V',\ov\fo')$$
for any oriented vector bundle $(V',\fo')$ over a bordered surface~$\Si$,
the collection $(\ov\fs_u)_u$ satisfies the condition of Definition~\ref{RelPinSpin_dfn3}
(because the collection $(\fs_u)_u$ does) and is thus 
indeed a relative $\Spin$-structure on~$(V,\ov\fo)$.
By~\eref{RelSpinPinStr_e33}, the resulting involution~\eref{RelOSpinRev_e}
is $H^2(X,Y;\Z_2)$-equivariant.
It also satisfies the two conditions after~\eref{RelOSpinRev_e}.
By~\eref{OSpinRev_e33} and~\eref{RelOSpinRev_e33}, 
\BE{RelOSpinRev_e1}
\io_X\big(\ov\os\big)=\ov{\io_X(\os)} 
\qquad\forall\,\os\!\in\!\OSp(V). \EE

\vspace{.1in}

We next describe a bijection between $\Sp_X(V,\fo)$ and 
$$\cP_X^{\pm}(V)\equiv \Sp_X(V_{\pm},\fo_V^{\pm}\big)
=\Sp_X(V\!\oplus\!(2\!\pm\!1)\tau_Y,\fo(2\!\pm\!1)\fo_Y\big);$$
the last equality holds because the orientation $\fo$ on~$V$ determines 
a canonical homotopy class of trivializations of $u^*\la(V)$
for every continuous map $u\!:\!\Si\!\lra\!Y$ from a bordered surface.
By the RelSpinPin~\ref{RelSpinPinObs_prop} property,
$(V,\fo)$ admits a relative $\Spin$-structure if and only if 
$V$ admits a relative $\Pin^{\pm}$-structure.
We can thus assume that $(V,\fo)$ admits a relative $\Spin$-structure.\\

If $\rk\,V\!\ge\!3$, we associate 
a relative $\Spin$-structure \hbox{$\fs\!\equiv\!(\fs_u)_u$}
on~$(V,\fo)$ with the relative $\Spin$-structure \hbox{$\fs^{\pm}\!\equiv\!(\fs_u^{\pm})_u$}
on~$(V_{\pm},\fo_V^{\pm})$ given~by 
\BE{RelPin2SpinRed_e33}
\fs_u^{\pm}=\begin{cases}
\{\phi_{\pm}\!\equiv\!\phi\!\oplus\!(2\!\pm\!1)(\det\phi)\!:
\phi\!\in\!\fs_u\big\}
&\hbox{if}~\prt\Si\!\neq\!\eset;\\
\fs_u,&\hbox{if}~\prt\Si\!=\!\eset;
\end{cases}
\quad\forall\,\big(u\!:(\Si,\prt\Si)\!\lra\!(X,Y)\!\big)\!\in\!\cL_X(Y).\EE
Since this operation of adding determinant factors commutes with the disjoint union operation on
the mod~2 homotopy classes of trivializations over one-dimensional manifolds and
$$\big(\fs_{\Si}(V',\fo')\!\big)^{\pm}=\fs_{\Si}\big(V'_{\pm},\fo_{V'}^{\pm}\big)$$
for any oriented vector bundle $(V',\fo')$ over a bordered surface~$\Si$,
the collection $(\fs_u^{\pm})_u$ satisfies the condition of Definition~\ref{RelPinSpin_dfn3}
(because the collection $(\fs_u)_u$ does) and is thus 
indeed a relative $\Spin$-structure on~$(V_{\pm},\fo_V^{\pm})$.\\

By~\eref{RelSpinPinStr_e33}, the map 
\BE{RelPin2SpinRed_e36}\Sp_X(V,\fo)\lra\cP_X^{\pm}(V), \qquad \fs\lra\fs^{\pm},\EE
is $H^2(X,Y;\Z_2)$-equivariant.
Along with the RelSpinPin~\ref{RelSpinPinStr_prop} property, 
this implies that this map is a bijection.
We take the map~$\fR_{\fo}^{\pm}$ in~\eref{RelPin2SpinRed_e} to be its inverse.
By~\eref{w2fsdfn_e} and the second case in~\eref{RelPin2SpinRed_e33}, 
this inverse satisfies the first condition after~\eref{RelPin2SpinRed_e}.
By~\eref{Pin2SpinRed_e33} and~\eref{RelPin2SpinRed_e33}, 
\BE{RelPin2SpinRed_e1}
\io_X\big(\fR_{\fo}^{\pm}(\fp)\!\big)=
\fR_{\fo}^{\pm}\big(\io_X(\fp)\!\big)  \qquad\forall\,\fp\!\in\!\cP^{\pm}(V).\EE
If~$\ov\phi$ is a trivialization of $\{\prt u\}^*(V,\ov\fo)$ as in~\eref{RelOSpinRev_e33}
and $\prt\Si\!\neq\!\eset$, then
$$\big(\ov\phi\big)_{\!\pm}=\bI_{n+2\pm1;3\pm1}\phi_{\pm}\!:
\{\prt u\}^*\big(V_{\pm}\big) \lra (\prt\Si)\!\times\!\R^{n+2\pm1}.$$
Since this trivialization is homotopic to~$\phi_{\pm}$,
the second condition after~\eref{RelPin2SpinRed_e} is also satisfied.\\

We extend~\eref{RelPin2SpinRed_e36} to vector bundles~$V$ with $\rk\,V\!=\!1,2$
as at the end of the proof of 
the SpinPin~\ref{Pin2SpinRed_prop} property in Section~\ref{SpinPindfn3_subs}.
Similarly to~\eref{Pin2SpinRed_e37},
\BE{RelPin2SpinRed_e37} \fR_{\fo}^-\!=\!\fR_{\St_V(\fo)}^-\!:
\cP_X^-(V)\!=\!\cP_X^-\big(\tau_Y\!\oplus\!V\big)\lra
\Sp_X(V,\fo)\!=\!\Sp_X\big(\St(V,\fo)\!\big)\EE
under the identifications~\eref{SpinPinCorr_e0} if $\rk\,V\!=\!1$ and $\fo\!\in\!\fO(V)$.
The resulting extensions are $H^2(X,Y;\Z_2)$-equivariant
and satisfy~\eref{RelPin2SpinRed_e1} and 
the last requirement of the RelSpinPin~\ref{RelPin2SpinRed_prop} property.
\end{proof}

\begin{proof}[{\bf{\emph{Proof of RelSpinPin~\ref{RelSpinPinStab_prop},\ref{RelSpinPinCorr_prop} 
properties}}}]
The two sides of~\eref{RelSpinPinCorr_e} are the same by definition;
we take this map to be the identity.
The first condition after~\eref{RelSpinPinCorr_e} is then satisfied.
Since the map~\eref{SpinPinCorr_e} in the perspective of Definition~\ref{PinSpin_dfn3}
is also the identity, 
\BE{RelSpinPinCorr_e1}
\io_X\big(\Co_V^{\pm}(\fp)\big)=\Co_V^{\pm}\big(\io_X(\fp)\!\big) 
\qquad\forall\,\fp\!\in\!\cP^{\pm}(V). \EE
By the RelSpinPin~\ref{RelSpinPinObs_prop} property,
an oriented vector bundle $(V,\fo)$ admits a relative $\Spin$-structure if and only if
$(\tau_Y\!\oplus\!V,\St_V(\fo))$ does.
In order to establish the RelSpinPin~\ref{RelSpinPinStab_prop} property,
we can thus assume that $(V,\fo)$ admits a relative $\Spin$-structure.\\

If $\rk\,V\!\ge\!3$, we associate a relative $\Spin$-structure \hbox{$\fs\!\equiv\!(\fs_u)_u$}
on~$(V,\fo)$ with the relative $\Spin$-structure \hbox{$\St_V\fs\!\equiv\!(\St_V\fs_u)_u$}
on $(\tau_Y\!\oplus\!V,\St(\fo))$ given~by 
\BE{RelSpinPinStab_e34}
\St_V\fs_u=\begin{cases}
\{\St_V\phi\!\equiv\id_{\tau_{\prt\Si}}\!\oplus\!\phi\!:\phi\!\in\!\fs_u\},
&\hbox{if}~\prt\Si\!\neq\!\eset;\\
\fs_u,&\hbox{if}~\prt\Si\!=\!\eset;
\end{cases}
\quad\forall\,\big(u\!:(\Si,\prt\Si)\!\lra\!(X,Y)\!\big)\!\in\!\cL_X(Y).\EE
Since this stabilization operation commutes with the disjoint union operation on
the mod~2 homotopy classes of trivializations over one-dimensional manifolds and
$$\St_{V'}\big(\fs_{\Si}(V',\fo')\big)=\fs_{\Si}\big(\St_{V'}(V',\fo')\!\big)$$
for any oriented vector bundle $(V',\fo')$ over a bordered surface~$\Si$,
the collection $(\St_V\fs_u)_u$ satisfies the condition of Definition~\ref{RelPinSpin_dfn3}
(because the collection $(\fs_u)_u$ does) and is thus 
indeed a relative $\Spin$-structure on~$(\tau_Y\!\oplus\!V,\St(\fo))$.
By~\eref{SpinPinStab_e34} and~\eref{RelSpinPinStab_e34},
\BE{RelSpinPinStab_e1}
\io_X\big(\St_V(\os)\big)=\St_V\big(\io_X(\os)\!\big) 
\qquad\forall\,\os\!\in\!\OSp(V). \EE

\vspace{.1in} 

The above construction of the first map in~\eref{RelSpinPinStab_e} with $\rk\,V\!\ge\!3$ 
induces the first map in~\eref{RelSpinPinStab_e} with $\rk\,V\!=\!1,2$ and 
the second map in~\eref{RelSpinPinStab_e} as in the proof of 
the SpinPin~\ref{SpinPinStab_prop},\ref{SpinPinCorr_prop} properties in
Section~\ref{SpinPindfn3_subs}.
The resulting maps thus satisfy~\eref{RelSpinPinStab_e1} for $\OSpin$-
and $\Pin^{\pm}$-structures on all vector bundles~$V$ over~$Y$.
By~\eref{w2fsdfn_e} and the second case in~\eref{RelSpinPinStab_e34},
they also satisfy the first condition after~\eref{RelSpinPinStab_e}. 
The verification that the maps in~\eref{RelSpinPinStab_e} and~\eref{RelSpinPinCorr_e} 
satisfy the remaining stated requirements is identical to the argument in Section~\ref{SpinPindfn3_subs}.
\end{proof}

\begin{proof}[{\bf{\emph{Proof of RelSpinPin~\ref{RelSpinPinSES_prop} property}}}]
For every continuous map $u\!:\Si\!\lra\!Y$ from a bordered surface, 
a short exact sequence~$\ce$ of vector bundles 
over~$Y$ as in~\eref{SpinPinSES_e0}  determines a homotopy class of isomorphisms
$$u^*V\approx u^*V'\!\oplus\!u^*V''$$
so that $u^*\io$ is the inclusion as the first component
on the right-hand side above and $u^*\fj$ is the projection to the second component.
Thus, it is sufficient to establish the RelSpinPin~\ref{RelSpinPinSES_prop} property
for the direct sum exact sequences as in~\eref{SpinPinDS_e0}.
Furthermore, an orientation~$\fo'$ on~$V'$ determines a homotopy class of trivializations 
of~$u^*\la(V')$ and thus of isomorphisms
\BE{RelSpinPinSES_e33}\begin{split} 
u^*\big(V'\!\oplus\!V''\big)_{\pm}
&\equiv u^*\big(V'\!\oplus\!V''\!\oplus\!(2\!\pm\!1)\la(V'\!\oplus\!V'')\big)\\
&\approx  u^*\big(V'\!\oplus\!V''\!\oplus\!(2\!\pm\!1)\la(V'')\big)
\equiv u^*\big(V'\!\oplus\!V''_{\pm}\big)
\end{split}\EE
for every continuous map $u\!:\Si\!\lra\!Y$ from a bordered surface.\\

Let $V'$ and $V''$ be vector bundles over~$Y$ of rank~$n'$ and~$n''$, respectively.
Suppose $n',n''\!\ge\!3$, $\fo'\!\in\!\fO(V')$, and $\fo''\!\in\!\fO(V'')$.
For a relative $\Spin$-structure \hbox{$\fs'\!\equiv\!(\fs_u')_u$} on~$(V',\fo')$ 
and a relative $\Spin$-structure \hbox{$\fs''\!\equiv\!(\fs_u'')_u$} on~$(V'',\fo'')$,
we define a $\Spin$-structure 
$\llrr{\fs',\fs''}_{\oplus}\!\equiv\!({\fs_u'\!\oplus\!\fs_u''})_u$
on $(V'\!\oplus\!V'',\fo'\fo'')$ by
\BE{Reldfn3sum_e}
\fs_u'\!\oplus\!\fs_u''=\begin{cases}
\{\phi'\!\oplus\!\phi''\!:\phi'\!\in\!\fs_u',\,\phi''\!\in\!\fs_u''\},
&\hbox{if}~\prt\Si\!\neq\!\eset;\\
\fs_u'\!+\!\fs_u''\!\in\!\Z_2,&\hbox{if}~\prt\Si\!=\!\eset;
\end{cases}
\quad\forall\,\big(u\!:(\Si,\prt\Si)\!\lra\!(X,Y)\!\big)\!\in\!\cL_X(Y).\EE
Since this direct sum operation commutes with the disjoint union operation on
the mod~2 homotopy classes of trivializations over one-dimensional manifolds and
$$\fs_{\Si}(V',\fo')\!\oplus\!\fs_{\Si}(V'',\fo'')=
\fs_{\Si}\big((V',\fo')\!\oplus\!(V'',\fo'')\big)$$
for any oriented vector bundles $(V',\fo')$ and $(V'',\fo'')$
over a bordered surface~$\Si$,
the collection $({\fs_u'\!\oplus\!\fs_u''})_u$ satisfies the condition 
of Definition~\ref{RelPinSpin_dfn3}
(because the collections $(\fs_u')_u$ and~$(\fs_u'')_u$ do) and is thus 
indeed a relative $\Spin$-structure on~$(V'\!\oplus\!V'',\fo'_{\ce}\fo'')$.
By~\eref{w2fsdfn_e} and the second case in~\eref{Reldfn3sum_e}, 
the resulting first map in~\eref{RelSpinPinSESdfn_e0} satisfies
the second statement in~\eref{RelOspinDSEquivSum_e}.
By~\eref{dfn3sum_e} and~\eref{Reldfn3sum_e},
\BE{RelSpinPinSES_e1}
\io_X\big(\os'\!\oplus\!\os''\big)=\io_X(\os')\!\oplus\!\io_X(\os'')
\qquad\forall\,\os'\!\in\!\OSp(V'),\,\os''\!\in\!\OSp(V''). \EE

\vspace{.1in}

The above construction determines the first map~$\llrr{\cdot,\cdot}_{\ce}$ 
in~\eref{RelSpinPinSESdfn_e0} if \hbox{$n',n''\!\ge\!3$}.
Along with~\eref{RelSpinPinSES_e33}, it also determines the second map 
$\llrr{\cdot,\cdot}_{\ce}$
if $n'\!\ge\!3$ and \hbox{$\rk\,V''_{\pm}\!\ge\!3$} (i.e.~not $V_-''$ if $n''\!=\!1$).
By the second statement in~\eref{RelOspinDSEquivSum_e} and the first 
statement after~\eref{RelSpinPinCorr_e}, this map satisfies~\eref{Relw2ses_e}.
By~\eref{RelSpinPinSES_e1} and~\eref{RelSpinPinCorr_e1}, it also satisfies~\eref{RelSpinPinSES_e1} 
with $\OSp(V'')$ replaced by~$\cP^{\pm}(V'')$.
By~\eref{RelSpinPinStr_e33}, both maps are $H^2(X,Y;\Z_2)$-biequivariant.
These two maps satisfy \eref{RelDSorient1_e} 
and~\ref{RelDSassoc_it} for vector bundles of ranks at least~3
by definition.
With~$\phi_{\pm}$ as in~\eref{RelPin2SpinRed_e33},
$$\big(\phi'\!\oplus\!\phi''\big)_{\pm}\!=\!\phi'\!\oplus\!\phi''_{\pm}\!:
\{\prt u\}^*(V'\!\oplus\!V'')_{\pm}\!=\!
\big(\{\prt u\}^*V'\big)\!\oplus\!\big(\{\prt u\}^*V_{\pm}''\big)
\lra (\prt\Si)\!\times\!\R^{n'+n''+2\pm1}\,.$$
Thus, \eref{RelDSorient2_e} is also satisfied in these cases.\\

For vector bundles $V',V''$ over~$Y$ not both of ranks at least~3,
we define the two maps in~\eref{RelSpinPinSESdfn_e0}
via~\eref{SpinPinSES_e35c} with $b\!\in\!2\Z^{\ge0}$ and 
with $\OSpin$- and $\Pin^{\pm}$-structures replaced by 
relative $\OSpin$- and $\Pin^{\pm}$-structures.
The resulting maps thus satisfy~\eref{RelSpinPinSES_e1} for 
$\OSpin$- and $\Pin^{\pm}$-structures on all vector bundles~$V',V''$ over~$Y$.
By~\eref{Relw2ses_e} for vector bundles of ranks at least~3 and  
the first statements after~\eref{RelPin2SpinRed_e} and~\eref{RelSpinPinStab_e},
they also satisfy~\eref{Relw2ses_e} for all vector bundles~$V',V''$.
The verification that the maps in~\eref{RelSpinPinSESdfn_e0}
satisfy the remaining stated requirements is identical to the argument
in the proof of the SpinPin~\ref{SpinPinSES_prop} property in Section~\ref{SpinPindfn3_subs}.
\end{proof}

\begin{proof}[{\bf{\emph{Proof of RelSpinPin~\ref{vsSpinPin_prop} property}}}]
By~\eref{RelSpinPinStr_e1}, 
the maps~\eref{vsSpinPin_e0} are equivariant with respect to the homomorphism~\eref{YdeXYZ2_e}.
By~\eref{RelOSpinRev_e1}, \eref{RelPin2SpinRed_e1}, \eref{RelSpinPinStab_e1}, 
\eref{RelSpinPinCorr_e1}, and~\eref{RelSpinPinSES_e1}, 
these maps respect all structures and correspondences of 
the SpinPin~\ref{OSpinRev_prop}-\ref{SpinPinSES_prop}
and RelSpinPin~\ref{RelOSpinRev_prop}-\ref{RelSpinPinSES_prop} properties.
By~\eref{w2fsdfn_e} and the definition of the map~\eref{YdeXYZ2_e},
\BE{vsSpinPin_e2} 
\io_X\big(\cP^{\pm}(V)\big)\subset\big\{\fp\!\in\!\cP_X^{\pm}(V)\!: w_2(\fp)\!=\!0\big\}.\EE
Suppose $\fp\!\in\!\cP_X^{\pm}(V)$ and $w_2(\fp)\!=\!0$.
By the RelSpinPin~\ref{RelSpinPinObs_prop} and SpinPin~\ref{SpinPinObs_prop} properties,
this implies that $V$ admits a $\Pin^{\pm}$-structure~$\fp'$.
By the RelSpinPin~\ref{RelSpinPinStr_prop} property and~\eref{vsSpinPin_e2}, 
\BE{vsSpinPin_e3} \eta\!\cdot\!\big(\io_X(\fp')\!\big)=\fp \qquad\hbox{and}\qquad
\eta|_X=0\EE
for some $\eta\!\in\!H^2(X,Y;\Z_2)$.
By the second condition in~\eref{vsSpinPin_e3} and 
the cohomology exact sequence for the pair $X\!\subset\!Y$,
$\eta\!=\!\de_{X,Y}(\eta')$ for some $\eta'\!\in\!H^1(Y;\Z_2)$.
Since the second map in~\eref{vsSpinPin_e0} is equivariant with respect 
to the homomorphism~\eref{YdeXYZ2_e}, the first condition~\eref{vsSpinPin_e3} 
implies~that 
$$\io_X\big(\eta'\!\cdot\!\fp'\big)=\fp\,.$$
Along with~\eref{vsSpinPin_e2}, this establishes~\eref{vsSpinPin_e}.
\end{proof}

\subsection{Topological preliminaries}
\label{TopolPrelim2_subs}

Our proof of Theorem~\ref{RelSpinStrEquiv_thm}\ref{RelPinSpinPrp_it}
for the relative Spin- and Pin- structures of Definition~\ref{RelPinSpin_dfn2} 
in Section~\ref{RelSpinPindfn2_subs}
relies on the standard topological statements recalled in this section.\\

We denote by \hbox{$\wt\ga_m\!\lra\!\wt\G(m)$} 
the tautological vector bundle over the Grassmannian of oriented $m$-planes in~$\R^{\i}$
and by $P_m\!\in\!\wt\G(m)$ the canonically oriented subspace of~$\R^{\i}$ spanned
by the first $m$ coordinate vectors.
For $m',m''\!\in\!\Z^{\ge0}$, choose a continuous~map 
\begin{gather*}
F_{m',m''}\!: \big(\wt\G(m')\!\times\!\wt\G(m''),(P_{m'},P_{m''})\!\big)\lra 
\big(\wt\G(m'\!+\!m''),P_{m'+m''}\big)\\ 
\hbox{s.t.}\qquad
\wt\ga_{m'}\!\times\!\wt\ga_{m''}=F_{m',m''}^*\wt\ga_{m'+m''}\,.
\end{gather*}
If $X$ is a topological space and $Y\!\subset\!X$, then
\BE{GnH2_e2} \big\{F_{m',m''}\!\circ\!(h',h'')\big\}^*w_2\big(\wt\ga_{m'+m''}\big)=
h'^*w_2\big(\wt\ga_{m'}\big)\!+\!h''^*w_2\big(\wt\ga_{m''}\big)
\in H^2(X,Y;\Z_2)\EE
for all continuous maps~$h'$ and~$h''$
from $(X,Y)$ to $(\wt\G(m'),P_{m'})$ and $(\wt\G(m''),P_{m''})$, respectively.
We denote by
$[(X,Y),(\wt\G(m),P_m)]$ the set of homotopy classes of continuous maps 
from~$X$ to~$\wt\G(m)$ which send~$Y$ to~$P_m$.

\begin{lmm}\label{GnH2_lmm}
Suppose $m\!\ge\!3$ and $(X,Y)$ is a CW pair.
The~map 
\BE{GnH2_e} \big[(X_3,Y_2),(\wt\G(m),P_m)\big]\lra H^2(X_3,Y_2;\Z_2), \quad 
[h]\lra h^*w_2(\wt\ga_m),\EE
is a bijection.
\end{lmm}

\begin{proof} The space $\wt\G(m)$ is connected and simply connected.
By the homotopy exact sequence for the fiber bundle
\BE{GnH2_e1}\SO(m) \lra \wt\F(m) \stackrel{q}{\lra} \wt\G(m),\EE
where $\wt\F(m)$ is the space of oriented $m$-frames on~$\R^{\i}$,
$$\pi_2\big(\wt\G(m)\!\big)\approx \pi_1\big(\SO(m)\!\big)\approx\Z_2
\qquad\hbox{and}\qquad 
\pi_3\big(\wt\G(m)\!\big)\approx \pi_2\big(\SO(m)\!\big)=0\,.$$
Thus, the claim follows immediately from Lemma~\ref{CWhomotop_lmm} 
with $n\!=\!2$, $Z\!=\!\wt\G(m)$, and $z_0\!=\!P_m$.
\end{proof}

Let $X$ be a topological space and $Y\!\subset\!X$. 
For each $m\!\in\!\Z^+$, we denote by $\OVB_m(X,Y)$ the set of isomorphism classes
of oriented rank~$m$ vector bundles 
over~$X$ with a homotopy class of trivialization over $Y\!\subset\!X$.
For every continuous map~$h$ from $(X,Y)$ to $(\wt\G(m),P_m)$,
the restriction to~$Y$ of
the oriented vector bundle $h^*\wt\ga_m$ over~$X$ is equipped with 
a trivialization given by the canonical identification
\BE{psihdfn_e}\psi_h\!: h^*\wt\ga_m\big|_Y=Y\!\times\!\wt\ga_m|_{P_m}=Y\!\times\!\R^m\,.\EE
We denote the homotopy class of~$\psi_h$ by~$\fs_h$.
If $h$ is the constant map to~$P_m$, we denote~$\fs_h$ by~$\fs_{X;m}$.
By the proof of \cite[Theorem~5.6]{MiSt}, the~map
\BE{RelSpinPinStr2_e1} \big[(X,Y),(\wt\G(m),P_m)\big]\lra \OVB_m(X,Y), \qquad 
h\lra \big[h^*\wt\ga_m,\fs_h\big],\EE
is a bijection.
We note~that 
\BE{RelSpinPinStr2_e2} 
\big[\{F_{m',m''}(h',h'')\}^*\wt\ga_{m'+m''},\fs_{F_{m',m''}(h',h'')}\big]
=\big[h'^*\wt\ga_{m'}\!\oplus\!h''^*\wt\ga_{m''},\fs_{h'}\!\oplus\!\fs_{h''}\big]\EE
for all $m',m''\!\in\!\Z^+$ and continuous maps~$h'$ and~$h''$
from $(X,Y)$ to $(\wt\G(m'),P_{m'})$ and $(\wt\G(m''),P_{m''})$, respectively.

\begin{crl}\label{GnH2_crl2}
Suppose $m',m''\!\in\!\Z^{\ge0}$, $(X,Y)$ is a CW pair,  and
\begin{gather}\notag
\big[(E_1',\fs_1')\big],\big[(E_2',\fs_2')\big]\in \OVB_{m'}(X,Y), \quad
\big[(E'',\fs'')\big]\in \OVB_{m''}(X,Y) \qquad\hbox{s.t.}\\
\label{GnH2crl_e0}
\big[(E_1',\fs_1')\!\oplus\!(E'',\fs'')\big]\!=\!
\big[(E_2',\fs_2')\!\oplus\!(E'',\fs'')\big]\in \OVB_{m'+m''}(X,Y).
\end{gather}
If $m'\!\ge\!3$, then
\BE{GnH2crl_e}\big[(E_1',\fs_1')|_{(X_3,Y_2)}\big]\!=\!\big[(E_2',\fs_2')|_{(X_3,Y_2)}\big]
\in \OVB_{m'}(X_3,Y_2).\EE
\end{crl}

\begin{proof}
Let $g$ be a continuous map from $(X_3,Y_2)$ to $(\wt\G(m''),P_{m''})$ representing
the preimage of $[(E'',\fs'')|_{(X_3,Y_2)}]$ under the bijection~\eref{RelSpinPinStr2_e1} 
with $m\!=\!m''$.
By~\eref{GnH2_e2} and~\eref{RelSpinPinStr2_e2}, the diagram
$$\xymatrix{ H^2(X_3,Y_2;\Z_2)\ar[d]^{\cdot+g^*w_2(\wt\ga_{m''})}
& \big[(X_3,Y_2),(\wt\G(m'),P_{m'})\big] 
\ar[l]_>>>>>>>>>{\eref{GnH2_e}}^>>>>>>>>>{\approx}
\ar[r]^>>>>>>>>>{\eref{RelSpinPinStr2_e1}}_>>>>>>>>>{\approx}
\ar[d]_{F_{m',m''}\circ(\cdot,g)}
&\OVB_{m'}(X_3,Y_2)\ar[d]_{\cdot\,\oplus[(E'',\fs'')|_{(X_3,Y_2)}]}\\
H^2(X_3,Y_2;\Z_2)& 
\big[(X_3,Y_2),(\wt\G(m'\!+\!m''),P_{m'+m''})\big] 
\ar[l]_>>>>>>{\eref{GnH2_e}}^>>>>>>{\approx} 
\ar[r]^>>>>>>{\eref{RelSpinPinStr2_e1}}_>>>>>>{\approx}
& \OVB_{m'+m''}(X_3,Y_2)}$$
then commutes.
Since the left vertical arrow in this diagram is a bijection,
so is the right vertical arrow.
This implies the claim.
\end{proof}

\begin{prp}\label{GnH2_prp2}
Let $m',m''\!\in\!\Z^{\ge0}$ with $m'\!\ge\!3$ and $(X,Y)$ be a CW pair.
Suppose 
$$\big[(E_1',\fs_1')\big],\big[(E_2',\fs_2')\big]\in\OVB_{m'}(X,Y),$$
$E''$ and $V$ are  oriented vector bundles over~$X$ and~$Y$, respectively,
and $\fs_{E'',V}$ is a homotopy class of trivializations of~$E''|_Y\!\oplus\!V$.
If there exists an isomorphism
\BE{GnH2prp_e}\Psi\!: E_1'\!\oplus\!E''\lra E_2'\!\oplus\!E''
\quad\hbox{s.t.}\quad 
\fs_1'\!\oplus\!\fs_{E'',V}= \big\{\Psi|_{(E_1'\oplus E)''|_Y}\!\oplus\!\id_V\big\}^* 
\big(\fs_2'\!\oplus\!\fs_{E'',V}\big),\EE
then~\eref{GnH2crl_e} holds.
\end{prp}

\begin{proof}
Since the diagram
$$\xymatrix{ H^2(X_3,Y_2;\Z_2)\ar[d]^{\cdot|_{(X_3,Y_1)}}
& \big[(X_3,Y_2),(\wt\G(m'),P_{m'})\big] 
\ar[l]_>>>>>>{\eref{GnH2_e}}^>>>>>>{\approx}
\ar[r]^>>>>>>{\eref{RelSpinPinStr2_e1}}_>>>>>>{\approx}
\ar[d]_{\cdot|_{(X_3,Y_1)}}
&\OVB_{m'}(X_3,Y_2)\ar[d]_{\cdot|_{(X_3,Y_1)}}\\
H^2(X_3,Y_1;\Z_2)& 
\big[(X_3,Y_1),(\wt\G(m'),P_{m'})\big] 
\ar[l]_>>>>>>{\eref{GnH2_e}}^>>>>>>{\approx} 
\ar[r]^>>>>>>{\eref{RelSpinPinStr2_e1}}_>>>>>>{\approx}
& \OVB_{m'}(X_3,Y_1)}$$
commutes and the left vertical arrow is injective, so is the right vertical arrow.
Thus, it is enough to establish~\eref{GnH2crl_e} with~$Y_2$ replaced by~$Y_1$.\\

We can assume that the sum of the ranks of~$E''$ and~$V$ is at least~3.
Let~$\fs''$ and $\fs_V$ be homotopy classes of trivializations of
the oriented vector bundles~$E''|_{Y_1}$ and~$V|_{Y_1}$, respectively, such that 
$\fs''\!\oplus\!\fs_V\!=\!\fs_{E'',V}$.
By~\eref{GnH2prp_e},  \eref{GnH2crl_e0} with $Y\!=\!Y_1$ holds.
By Corollary~\ref{GnH2_crl2}, \eref{GnH2crl_e} with~$Y_2$ replaced by~$Y_1$
thus holds.
\end{proof}

Suppose $(X,Y)$ is a CW pair. 
A homotopy class of trivializations of an oriented rank~$m$ vector bundle~$(E,\fo_E)$ over~$X_3$ determines
a homotopy class of trivializations of~$(E,\fo_E)|_{Y_2}$.
Trivializations of~$(E,\fo_E)|_{Y_2}$ obtained from trivializations of~$(E,\fo_E)$ differ by
the action of $h\!\equiv\!f|_{Y_2}$ as in~\eref{Spin1a_e4} for
some continuous function~$f$ from~$X_3$ to~$\SO(m)$. 
Along with the proof of the SpinPin~\ref{SpinPinStr_prop} property in 
Section~\ref{SpinPindfn2_subs}, 
this implies~that a trivializable vector bundle~$(E,\fo_E)$ over~$X_3$ determines
a collection 
$$H^1(X_3;\Z_2)\fs_0\big(\!(E,\fo_E)|_{Y_2}\big)\subset\Sp\big(\!(E,\fo_E)|_{Y_2}\big)$$
of homotopy classes of trivializations of $(E,\fo_E)|_{Y_2}$ which form an orbit of
$$\big\{\eta|_{Y_2}\!:\eta\!\in\!H^1(X_3;\Z_2)\big\}\subset H^1(Y_2;\Z_2)$$
under the action of $H^1(Y_2;\Z_2)$ on $\Sp((E,\fo_E)|_{Y_2})$ provided
by the SpinPin~\ref{SpinPinStr_prop}\ref{SpinStr_it} property in Section~\ref{SpinPinProp_subs}.
For any $\eta\!\in\!H^1(Y_2;\Z_2)$, let 
$$\eta\!\cdot\!H^1(X_3;\Z_2)\fs_0\big(\!(E,\fo_E)|_{Y_2}\big)
\subset\Sp\big(\!(E,\fo_E)|_{Y_2}\big)$$
be the image of this orbit under the action of~$\eta$.

\begin{lmm}\label{GnH2_lmm2}
Suppose $m\!\ge\!3$,  $(X,Y)$ is a CW pair, $\eta\!\in\!H^1(Y_2;\Z_2)$, and
$h$ is a continuous map from $(X_3,Y_2)$ to $(\wt\G(m),P_m)$.
If
\BE{GnH22_e0a}\de_{X_3,Y_2}(\eta)=h^*w_2(\wt\ga_m)\in H^2\big(X_3,Y_2;\Z_2\big),\EE
then the oriented vector bundle $h^*\wt\ga_m$ over~$X_3$ is trivializable and 
\BE{GnH22_e0b}\fs_h\in\eta\!\cdot\!H^1(X_3;\Z_2)\fs_0\big(h^*\wt\ga_m|_{Y_2}\big)\,.\EE
\end{lmm}

\begin{proof}
By~\eref{GnH22_e0a} and the cohomology exact sequence for $Y_2\!\subset\!X_3$,
$w_2(h^*\wt\ga_m)\!=\!0$ in $H^2(X_3;\Z_2)$.
By the bijectivity of~\eref{RelSpinPinStr2_e1} and~\eref{GnH2_e} with $Y_2\!=\!\eset$,
the vector bundle $h^*\wt\ga_m$ is thus trivializable.
By the bijectivity of~\eref{RelSpinPinStr2_e1},
it is sufficient to establish~\eref{GnH22_e0b}
for any map~$h$ satisfying~\eref{GnH22_e0a}.\\

Let $\eta_m\!\in\!H^1(\SO(m);\Z_2)$ be the generator as before.
We identify the fiber of the fiber bundle~$\wt\F(m)$ in~\eref{GnH2_e1} over $P_m\!\in\!\wt\G(m)$
with~$\SO(m)$.
By \cite[Theorem~7.2.8]{Spanier}, the homomorphism
$$q_*\!:\pi_k\big(\wt\F(m),\SO(m)\!\big)\lra\pi_k\big(\wt\G(m)\!\big)$$
is an isomorphism for every $k\!\in\!\Z^+$.
Since $\wt\F(m)$ is contractible, 
the boundary homomorphism
$$\prt\!:\pi_2\big(\wt\F(m),\SO(m)\!\big)\lra \pi_1\big(\SO(m)\!\big)$$
in the homotopy exact sequence for~\eref{GnH2_e1} is an isomorphism.
Let 
$$\wt\al_2\!:\big(\bD^2,S^1,1\big)\lra \big(\wt\F(m),\SO(m),\bI_m\big)$$
be a continuous map generating 
$$\pi_2\big(\wt\F(m),\SO(m)\!\big)\approx\pi_1\big(\SO(m),\bI_m\!\big)\approx\Z_2.$$
Thus, $q\!\circ\!\wt\al_2$ induces a continuous map~$\al_2$ from $(\P^1,\i)$ to $(\wt\G(m),P_m)$
generating $\pi_2(\wt\G(m)\!)$ and 
$$\al_1\!\equiv\!\wt\al_2|_{S^1}\!: (S^1,1)\lra\big(\SO(m),\bI_m\big)$$
is a loop generating $\pi_1(\SO(m)\!)\!\approx\!\Z_2$.
Thus,
\BE{GnH22_e4} \blr{\eta_m,\al_{1*}[S^1]_{\Z_2}}=1\in\Z_2 \qquad\hbox{and}\qquad
\blr{w_2(\wt\ga_m),\wt\al_{2*}[S^2]_{\Z_2}}=1\in\Z_2.\EE

\vspace{.1in}

Let 
$$\wt\eta\!:C_1(Y_2;\Z_2)\lra \Z_2 \qquad\hbox{and}\qquad
\wt\eta_2\!:C_2(X_3,Y_2;\Z_2)\lra \Z_2$$
be a cocycle representing $\eta\!\in\!H^1(Y_2;\Z_2)$ and 
the cocycle sending each 2-cell $e\!\subset\!X_2$ to
$$\wt\eta_2(e)\equiv  \wt\eta\big((\prt e)\!\cap\!C_1(Y_2;\Z_2)\!\big),$$
respectively.
Thus, $\wt\eta_2$ represents $\de_{X_3,Y_2}(\eta)\!\in\!H^2(X_3,Y_2;\Z_2)$.
We define a continuous map
$$\wt{h}\!:\big(X_3,Y_2\big)\lra \big(\wt\F(m),\SO(m)\!\big)$$
as follows.
Let $f$ be a continuous map sending~$X_0$, 
$X_1\!-\!Y_1$, and each 1-cell $e\!\subset\!Y_1$ with $\wt\eta(e)\!=\!0$ to~$\bI_m$ 
and each 1-cell $e\!\subset\!Y_1$ with $\wt\eta(e)\!=\!1$ to~$\al_1$.
By the first statement in~\eref{GnH22_e4},
\BE{GnH22_e6} \big\{f|_{Y_1}\big\}^*\eta_m=\eta|_{Y_1}\in H^1(Y_1;\Z_2).\EE
Along with \cite[Theorem~8.1.17]{Spanier}, this implies that $f|_{Y_1}$ extends 
to a continuous map from $Y_2\!\cup\!X_1$ to $\SO(m)$.
For the same reason, $f$ extends continuously over
the 2-cells $e\!\subset\!X_2\!-\!Y_2$ with $\wt\eta_2(e)\!=\!0$. 
We extend the resulting map to a continuous map~$\wt{h}$ from~$X_2$ to~$\wt\F(m)$
by sending  each 2-cell $e\!\subset\!X_2\!-\!Y_2$ with $\wt\eta_2(e)\!=\!1$
to~$\wt\al_2$.
Since the space~$\wt\F(m)$ is contractible, the map~$\wt{h}$ extends to~$X_3$.\\

By the cohomology exact sequence for the pair $(Y_2,Y_1)$, the restriction homomorphism
$$H^1(Y_2;\Z_2)\lra H^1(Y_1;\Z_2)$$
is injective.
Along with~\eref{GnH22_e6}, this implies that
\BE{GnH22_e6b} \big\{\wt{h}|_{Y_2}\big\}^*\eta_m=\eta\in H^1(Y_2;\Z_2).\EE
The continuous~map
$$h\!\equiv\!q\!\circ\!\wt{h}\!:(X_2,Y_2)\lra\big(\wt\G(m),P_m\big)$$
sends every 2-cell $e_2\!\subset\!X_2$ with $\wt\eta_2(e)\!=\!0$ to~$P_m$
and every 2-cell $e_2\!\subset\!X_2$ with $\wt\eta_2(e)\!=\!1$ to~$\al_2$.
Along with the second statement in~\eref{GnH22_e4}, this implies~that 
$h$ satisfies~\eref{GnH22_e0a}.\\

The oriented vector bundle $q^*\wt\ga_m$ over $\wt\F(m)$ has a canonical trivialization~$\Psi_m$,
since each point of~$\wt\F(m)$ is an oriented frame for the fiber over~it.
The restriction of~$\Psi_m$ to $\SO(m)\!\subset\!\wt\F(m)$ 
differs from the canonical identification
$$q^*\wt\ga_m\big|_{\SO(m)}=\SO(m)\!\times\!\wt\ga_m\big|_{P_m}
=\SO(m)\!\times\!\R^m$$
by the action as in~\eref{Spin1a_e4} of the identity map on~$\SO(m)$.
Thus, the canonical identification~$\psi_h$ in~\eref{psihdfn_e} of the restriction~of 
$$h^*\wt\ga_m=\wt{h}^*\big(q^*\wt\ga_m)$$
to $Y_2\!\subset\!X_3$ differs from the restriction of the trivialization 
of $h^*\wt\ga_m$ induced by~$\Psi_m$ by the multiplication by $\wt{h}|_{Y_2}\!=\!f$.
Along with~\eref{GnH22_e6b} and the proof of the SpinPin~\ref{SpinPinStr_prop} property 
in Section~\ref{SpinPindfn2_subs}, this establishes~\eref{GnH22_e0b}.
\end{proof}

\begin{crl}\label{GnH2_crl2b}
Let $m\!\ge\!3$,  $(X,Y)$, $\eta$, and $h$ be as in Lemma~\ref{GnH2_lmm2}.
There exists a trivialization~$\Psi_h$ of~$h^*\wt\ga_m$ such~that
$$\fs_h\!\oplus\!\fs_{E'',V}= \big\{\Psi_h|_{(h^*\wt\ga_m)|_{Y_2}}\!\oplus\!
\id_{E''|_{Y_2}\oplus V}\big\}^{\!*}\big(\fs_{X;m}\!\oplus\!
(\eta\!\cdot\!\fs_{E'',V})\!\big)$$
for all oriented vector bundles $E''$ and $V$ over~$X$ and~$Y_2$, respectively,
and homotopy classes $\fs_{E'',V}$ of trivializations of~$E''|_{Y_2}\!\oplus\!V$. 
\end{crl}

\begin{proof} Let $\fs_0(h^*\wt\ga_m|_{Y_2})$ be the homotopy class of trivializations
of $h^*\wt\ga_m|_{Y_2}$ determined by a trivialization~$\Psi$ of~$h^*\wt\ga_m$.
By Lemma~\ref{GnH2_lmm2}, there exists $\eta'\!\in\!H^1(X_3;\Z_2)$ such that 
\BE{GnH2crl2b_e3} \eta\!\cdot\!s_h=(\eta'|_{Y_2})\!\cdot\!\fs_0\big(h^*\wt\ga_m|_{Y_2}\big)
\in \Sp\big(h^*\wt\ga_m|_{Y_2}\big).\EE
By the second statement of Corollary~\ref{SOnH1_crl}, 
$\eta'\!=\!f^*\eta_m$ for some continuous map~$f$ from~$X_3$ to~$\SO(m)$.
Let
$$\Psi_h\!=\!f\Psi\!: h^*\wt\ga_m\lra m\tau_X\,.$$
By~\eref{GnH2crl2b_e3} and the construction of the $H^1(Y;\Z_2)$-action
in the proof of the SpinPin~\ref{SpinPinStr_prop} property in Section~\ref{SpinPindfn2_subs},
$$\big\{\Psi_h|_{(h^*\wt\ga_m)|_{Y_2}}\big\}^{\!*}\fs_{X;m}=\eta\!\cdot\!s_h.$$
Along with the SpinPin~\ref{SpinPinSES_prop} property 
in Section~\ref{SpinPinProp_subs}, this implies~that 
$$\big\{\Psi_h|_{(h^*\wt\ga_m)|_{Y_2}}\!\oplus\!
\id_{E''|_{Y_2}\oplus V}\big\}^{\!*}
\big(\fs_{X;m}\!\oplus\!(\eta\!\cdot\!\fs_{E'',V})\!\big)
=\big(\eta\!\cdot\!s_h\big)\!\oplus\!\big(\eta\!\cdot\!\fs_{E'',V}\big)
=s_h\!\oplus\!\fs_{E'',V}\,.$$
This establishes the claim.
\end{proof}

\subsection{Proof of Theorem~\ref{RelSpinStrEquiv_thm}\ref{RelPinSpinPrp_it}:  
Definition~\ref{RelPinSpin_dfn2} perspective}
\label{RelSpinPindfn2_subs}

We now establish the statements of Section~\ref{RelSpinPinProp_subs} for the notions of
relative $\Spin$-structure and $\Pin^{\pm}$-structure arising from Definition~\ref{RelPinSpin_dfn2}.
Throughout this section, the terms relative $\Spin$-structure and $\Pin^{\pm}$-structure
refer to these notions.

\begin{proof}[{\bf{\emph{Proof of RelSpinPin~\ref{RelSpinPinObs_prop} property}}}]
Suppose $\fo\!\in\!\fO(V)$ and $(V,\fo)$ admits a relative $\Spin$-structure 
$\fs\!\equiv\!(E,\fo_E,\fs_{E,V})$ with $w_2(E)\!=\!\mu$.
Since the vector bundle $(E,\fo_E)|_{Y_2}\!\oplus\!(V,\fo)|_{Y_2}$ then admits a $\Spin$-structure,
$\mu|_{Y_2}\!=\!w_2(E)|_{Y_2}$ by the SpinPin~\ref{SpinPinObs_prop}\ref{SpinObs_it} property
in Section~\ref{SpinPinProp_subs}.
Since the restriction homomorphism~\eref{SpinPinObs2_e2} is injective, 
it follows that $\mu|_Y\!=\!w_2(V)$.\\

Suppose $\mu|_Y\!=\!w_2(V)$.
By Lemma~\ref{GnH2_lmm} with $Y\!=\!\eset$, $\mu|_{X_3}\!=\!w_2(E)$
for some oriented vector bundle $(E,\fo_E)\!\equiv\!h^*\wt\ga_3$ over~$X_3$.
By the SpinPin~\ref{SpinPinObs_prop}\ref{SpinObs_it} property, 
the vector bundle $(E,\fo_E)|_{Y_2}\!\oplus\!(V,\fo)|_{Y_2}$ then admits a $\Spin$-structure.\\

The above establishes the RelSpinPin~\ref{RelSpinPinObs_prop}\ref{RelSpinObs_it} property.
The RelSpinPin~\ref{RelSpinPinObs_prop}\ref{RelPinObs_it} property
is obtained similarly from the SpinPin~\ref{SpinPinObs_prop}\ref{PinObs_it} property
and Lemma~\ref{GnH2_lmm}.
\end{proof}

\begin{proof}[{\bf{\emph{Proof of RelSpinPin~\ref{RelSpinPinStr_prop} property}}}]
By the cohomology exact sequences for the triples \hbox{$Y_2\!\subset\!Y\!\subset\!X$}
and $Y_2\!\subset\!X_3\!\subset\!X$, the restriction homomorphism
\BE{RelSpinPinStr_e1a} H^2\big(X,Y;\Z_2\big)\lra H^2\big(X_3,Y_2;\Z_2\big)\EE
is an isomorphism.
It is thus sufficient to establish the claims of this property with 
$H^2(X,Y;\Z_2)$ replaced by $H^2(X_3,Y_2;\Z_2)$.\\

Suppose $\fs\!\equiv\!(E,\fo_E,\fs_{E,V})$ is a relative $\Spin$-structure on $(V,\fo)$ 
and $\eta\!\in\!H^2(X_3,Y_2;\Z_2)$.
Let $m\!\ge\!3$. By Lemma~\ref{GnH2_lmm}, 
\BE{hetamdfn_e}\eta=h_{\eta;m}^*w_2(\wt\ga_m)\in H^2\big(X_3,Y_2;\Z_2\big)\EE
for some continuous map~$h_{\eta;m}$ from $(X_3,Y_2)$ to $(\wt\G_m,P_m)$.
We define
\BE{Reletafsdfn2_e}\eta\!\cdot\!\fs=\big[h_{\eta;m}^*\wt\ga_m\!\oplus\!(E,\fo_E),
\lr{\fs_{h_{\eta;m}},\fs_{E,V}}_{\oplus}\big],\EE
with $\fs_{h_{\eta;m}}$ as below~\eref{psihdfn_e} and
$\lr{\cdot,\cdot}$ as in~\eref{SpinPinSESdfn_e0}.
By Lemma~\ref{GnH2_lmm}, $[\eta\!\cdot\!\fs]\!\in\!\Sp_X(V,\fo)$ 
does not depend on the choice of~$h_{\eta;m}$.
By~\eref{RelSpinPinStr2_e2}, $[\eta\!\cdot\!\fs]$ does not depend on the choice of~$m$
or the representative~$\fs$ of \hbox{$[\fs]\!\in\!\Sp_X(V,\fo)$}. 
By~\eref{GnH2_e2} and~\eref{RelSpinPinStr2_e2}, the resulting map 
\BE{RelSpinPinStr2_e5}H^2(X_3,Y_2;\Z_2)\!\times\!\Sp_X(V,\fo)\lra \Sp_X(V,\fo)\EE
is an action of $H^2(X_3,Y_2;\Z_2)$ on~$\Sp_X(V,\fo)$.
By~\eref{Reletafsdfn2_e} and~\eref{hetamdfn_e}, 
this action satisfies the second statement in~\eref{RelSpinPinStr_e0}.\\

Suppose $[\eta\!\cdot\!\fs]\!=\![\fs]$.
Thus, there exist $m'\!\in\Z^{\ge0}$ and an isomorphism
$$\Psi\!:\St^{m'}\!\big(h_{\eta;m}^*\wt\ga_m\!\oplus\!(E,\fo_E)\!\big)
\lra \St^{m+m'}\!(E,\fo_E)$$ 
of oriented vector bundles over~$X_3$ so that 
$$\St_{(h_{\eta;m}^*\wt\ga_m\oplus E)|_{Y_2}\oplus V|_{Y_2}}^{m'}\!\!
\big(\blr{\fs_{h_{\eta;m}},\fs_{E,V}}_{\oplus}\big)=
\Big\{\!\Psi\big|_{\St^{m'}\!(h_{\eta;m}^*\wt\ga_m\oplus E)|_{Y_2}}\!\!\oplus\!\id_{V|_{Y_2}}
\!\Big\}^{\!*}
\St_{E|_{Y_2}\oplus V|_{Y_2}}^{m+m'}\!\big(\fs_{E,V}\big).$$
By Proposition~\ref{GnH2_prp2}, this implies that 
$$\big[h_{\eta;m}^*\wt\ga_m,\fs_{h_{\eta;m}}\big]\!=\!
\big[m\tau_{X_3},\fs_{X;m}\big]\in \OVB_m(X_3,Y_2).$$
Since the maps~\eref{GnH2_e} and~\eref{RelSpinPinStr2_e1} are bijections, 
it follows that $\eta\!=\!0$.
Thus, the action of $H^2(X_3,Y_2;\Z_2)$ on~$\Sp_X(V,\fo)$ defined above is free.\\

Suppose $\fs'\!\equiv\!(E',\fo_{E'},\fs_{E',V})$
is another relative $\Spin$-structure on~$(V,\fo)$.
We can assume that the ranks of~$E$ and~$E'$ are the same.
By the RelSpinPin~\ref{RelSpinPinObs_prop} property, 
$$\big(w_2(E)\!-\!w_2(E')\big)\big|_{Y_2}=
w_2(V)\big|_{Y_2}\!-\!w_2(V)\big|_{Y_2}=0\in H^2(Y_2;\Z_2).$$ 
Along with the cohomology exact sequence for the pair $Y_2\!\subset\!X_3$ and 
Lemma~\ref{GnH2_lmm}, this implies~that 
\BE{RelSpinPinStr2_e9}w_2(E)\!-\!w_2(E')=h_2^*w_2\big(\wt\ga_3\big)\in H^2(X_3;\Z_2)\EE
for some continuous map $h_2$ from $(X_3,Y_2)$ to $(\wt\G(3),P_3)$.
Let 
$$\eta_2\!\equiv\!h_2^*w_2\big(\wt\ga_3\big)\in H^2(X_3,Y_2;\Z_2)\,.$$
By~\eref{RelSpinPinStr2_e9}, Lemma~\ref{GnH2_lmm}, and the bijectivity of~\eref{RelSpinPinStr2_e1},
there exists an isomorphism 
$$\Psi_2\!: h_2^*\wt\ga_3\!\oplus\!(E,\fo_E)\lra 3(\tau_X,\fo_X)\!\oplus\!(E',\fo_{E'})$$
of oriented vector bundles over~$X_3$.\\

By the SpinPin~\ref{SpinPinStr_prop}\ref{SpinStr_it} property 
in Section~\ref{SpinPinProp_subs},  there exists $\eta_1\!\in\!H^1(Y_2;\Z_2)$ such~that 
\BE{RelSpinPinStr2_e11}\eta_1\!\cdot\!\big(\fs_{h_2}\!\oplus\!\fs_{E,V}\big)=
\big\{\Psi_2|_{(h_2^*\wt\ga_3\oplus(E,\fo_E))|_{Y_2}}\!\oplus\!\id_{V|_{Y_2}}\big\}^{\!*}
\big(\fs_{X_3;3}\!\oplus\!\fs_{E',V}\big)\,.\EE
Let $h_1$ be a continuous map from $(X_3,Y_2)$ to $(\wt\G(3),P_3)$ so that 
\eref{GnH22_e0a} with $(\eta,h)$ replaced by~$(\eta_1,h_1)$ holds and~$\Psi_1$ 
be a trivialization provided by Corollary~\ref{GnH2_crl2b}.
By~\eref{RelSpinPinStr2_e11}, 
$$\Psi\!\equiv\!\Psi_1\!\oplus\!\Psi_2\!: h_1^*\wt\ga_3\!\oplus\!
\big(h_2^*\wt\ga_3\!\oplus\!(E,\fo_E)\!\big)\lra 
3(\tau_X,\fo_X)\!\oplus\!\big(3(\tau_X,\fo_X)\!\oplus\!(E',\fo_{E'})\!\big)$$
is an isomorphisms of oriented vector bundles over $X_3$ such that
$$\fs_{h_1}\!\oplus\!\fs_{h_2}\!\oplus\!\fs_{E,V}
=\big\{\Psi|_{(h_1^*\wt\ga_3\oplus h_2^*\wt\ga_3\oplus(E,\fo_E))|_{Y_2}}
\!\oplus\!\id_{V|_{Y_2}}\big\}^{\!*}
\big(\fs_{X_3;6}\!\oplus\!\fs_{E,V}\big)\,.$$
By the definition of the $H^2(X,Y;\Z_2)$-action above, this implies that 
$$\big(\de_{X_3,Y_2}(\eta_1)\!+\!\eta\big)\!\cdot\!\fs=\fs'\,.$$
Thus, this action is transitive.
By Corollary~\ref{GnH2_crl2b} with~$E''$ being the rank~0 bundle, 
it also satisfies the first identity in~\eref{RelSpinPinStr_e1}.\\

We define the action of $H^2(X,Y;\Z_2)$ on $\cP_X^{\pm}(X)$ by~\eref{Reletafsdfn2_e},
replacing the Spin structures by $\Pin^{\pm}$-structures.
By the same reasoning as in the relative Spin case, this indeed defines a group action.
By the SpinPin~\ref{SpinPinCorr_prop} property in Section~\ref{SpinPinProp_subs}
and the second statement in the SpinPin~\ref{SpinPinSES_prop}\ref{DSEquivSum_it} property,
the $H^2(X,Y;\Z_2)$-action on $\Sp_X(V_{\pm},\fo_V^{\pm})$ being free and transitive implies
that $H^2(X,Y;\Z_2)$-action on $\cP_X^{\pm}(V)$ is also free and transitive.
The SpinPin~\ref{SpinPinCorr_prop} property and the first identity in~\eref{RelSpinPinStr_e1}
imply the second identity in~\eref{RelSpinPinStr_e1}.
By~\eref{SpinPinStr_e} and Corollary~\ref{GnH2_crl2b}, 
\eref{RelSpinPinStr_e} holds as well.
\end{proof}

\begin{proof}[{\bf{\emph{Proof of 
RelSpinPin~\ref{RelOSpinRev_prop},\ref{RelPin2SpinRed_prop} properties}}}]
Let $\fo\!\in\!\fO(V)$.
For an oriented vector bundle $(E,\fo_E)$, we denote by $\fo_E|_{Y_2}\fo|_{Y_2}$
the orientation on $E|_{Y_2}\!\oplus\!V|_{Y_2}$ induced by~$\fo_E$ and~$\fo$.
With the notation as in~\eref{OSpinRev_e} and~\eref{Pin2SpinRed_e},
we define the maps~\eref{RelOSpinRev_e} and~\eref{RelPin2SpinRed_e} by
\begin{alignat*}{2}
\OSp_X(V)&\lra \OSp_X(V), &\quad
\os\!\equiv\!\big(E,\fo_E,\os_{E,V}\big)&\lra 
\ov\os\!\equiv\!\big(E,\fo_E,\ov{\os_{E,V}}\big),\\
\fR_{\fo}^{\pm}\!:\cP_X^{\pm}(V)&\lra\Sp_X(V,\fo), &\quad
\fp\!\equiv\!\big(E,\fo_E,\fp_{E,V}\big)&\lra 
\big(E,\fo_E,\fR_{\fo_E|_{Y_2}\fo|_{Y_2}}^{\,\pm}(\fp_{E,V})\big).
\end{alignat*}
By~\eref{SpinPinStab_e2}, these two maps are well-defined on the equivalence classes
of relative $\OSpin$- and $\Pin^{\pm}$-structures of Definition~\ref{RelPinSpin_dfn2}.
It is immediate that the two maps preserve~$w_2$ and satisfy~\eref{RelOSpinRev_e1} 
and~\eref{RelPin2SpinRed_e1}.
By~\eref{OspinDSEquivSum_e} and~\eref{DSorient2_e}, 
the two maps are $H^2(X,Y;\Z_2)$-equivariant.
The last requirement in the RelSpinPin~\ref{RelPin2SpinRed_prop} property holds
by the statement after~\eref{Pin2SpinRed_e}.
\end{proof}

\begin{proof}[{\bf{\emph{Proof of RelSpinPin~\ref{RelSpinPinStab_prop},\ref{RelSpinPinCorr_prop} 
properties}}}]
For an isomorphism $\Psi$ of vector bundles $V$ and $W$ over~$Y_2$, let
$$\Psi_{\pm}\!:V_{\pm}\lra W_{\pm} $$
be the induced isomorphism.
For vector bundles $V',V''$ over $Y_2$, let 
\BE{RelSpinPinStab_e21a}\Psi_{V',V''}\!: V'\!\oplus\!V''\lra V''\!\oplus\!V'\EE
be the factor-interchanging isomorphism.
If $V'''$ is another vector bundle over~$Y_2$, then 
\BE{RelSpinPinStab_e21}  \Psi_{V'\oplus V'', V'''}= 
  \big\{\Psi_{V',V'''}\!\oplus\!\id_{V''}\big\}
	\!\circ\!\big\{\id_{V'}\!\oplus\!\Psi_{V'',V'''}\big\}\,.\EE
An orientation $\fo'$ on $V'$ determines
a homotopy class of isomorphisms 
$$\Psi_{\fo',V''}\!: (V',\fo')\!\oplus\!\big(V''_{\pm},\fo_{V''}^{\pm}\big)
\approx \big((V'\!\oplus\!V'')_{\pm},\fo_{V'\oplus V''}^{\pm}\big)$$
of vector bundles over~$Y_2$.
If in addition $\fo''$ is an orientation of~$V''$, then 
\begin{gather}\label{RelSpinPinCorr_e21}  
\Psi_{\fo'\fo'',V'''}\sim \Psi_{\fo',V''\oplus V'''}\!\circ\!
\big\{\id_{V'}\!\oplus\!\Psi_{\fo'',V'''}\big\}\,,\\
\label{RelSpinPinCorr_e21b}
\begin{split}
&\Psi_{\fo'',V'\oplus V'''}\!\circ\!
\big\{\id_{V''}\!\oplus\!\Psi_{\fo',V'''}\big\}\!\circ\!
\big\{\Psi_{V',V''}\!\oplus\!\id_{V'''_{\pm}}\big\}\\
&\hspace{1in}\sim\big\{\Psi_{V',V''}\!\oplus\!\id_V\big\}_{\pm}\!\circ\!
\Psi_{\fo',V''\oplus V'''}\!\circ\!
\big\{\id_{V'}\!\oplus\!\Psi_{\fo'',V'''}\big\}\,.
\end{split}\end{gather}

\vspace{.15in}

With the notation as in~\eref{SpinPinStab_e}, we define the two maps in~\eref{RelSpinPinStab_e} 
by 
\begin{equation*}\begin{split}
\St_V\big(E,\fo_E,\os_{E,V}\big)&=
\big(E,\fo_E,
\big\{\Psi_{E|_{Y_2},\tau_{Y_2}}\!\oplus\!\id_{V|_{Y_2}}\big\}^{\!*}
\!\big(\St_{E|_{Y_2}\oplus V|_{Y_2}}(\os_{E,V})\!\big)\big),\\
\St_V^{\pm}\big(E,\fo_E,\fp_{E,V}\big)&=
\big(E,\fo_E,
\big\{\Psi_{E|_{Y_2},\tau_{Y_2}}\!\oplus\!\id_{V|_{Y_2}}\big\}^{\!*}
\!\big(\St_{E|_{Y_2}\oplus V|_{Y_2}}(\fp_{E,V})\!\big)\big).
\end{split}\end{equation*}
By \eref{RelSpinPinStab_e21}, \eref{OspinDSEquivSum_e}, and~\eref{RelOspinDSassoc_e},
the first map is well-defined on the equivalence classes
of relative $\OSpin$- and $\Pin^{\pm}$-structures of Definition~\ref{RelPinSpin_dfn2}
and equivariant with respect to the $H^2(X,Y;\Z_2)$-actions.
By \eref{RelSpinPinStab_e21} and~\ref{DSstab_it} and~\ref{DSassoc_it}
in the SpinPin~\ref{SpinPinSES_prop} property 
in Section~\ref{SpinPinProp_subs},  so is the second~map.
It is immediate that the two maps preserve~$w_2$ and satisfy~\eref{RelSpinPinStab_e1} 
for $\OSpin$- and $\Pin^{\pm}$-structures on all vector bundles~$V$ over~$Y$.
By~\eref{SpinPinStab_e2}, they also satisfy~\eref{RelSpinPinStab_e2}.\\

With the notation as in~\eref{SpinPinCorr_e},  we define the map in~\eref{RelSpinPinCorr_e} 
by 
$$\Co_V^{\pm}\big(E,\fo_E,\fp_{E,V}\big)
=\big(E,\fo_E, \Psi_{\fo_E|_{Y_2},V|_{Y_2}}^{\,*}\!
\big(\Co_{E|_{Y_2}\oplus V|_{Y_2}}^{\pm}(\fp_{E,V})\!\big)\big).$$
By~\eref{RelSpinPinCorr_e21} and 
the second statement in the SpinPin~\ref{SpinPinSES_prop}\ref{DSEquivSum_it} property, 
this map is well-defined on the equivalence classes
of relative $\Pin^{\pm}$- and $\OSpin$-structures of Definition~\ref{RelPinSpin_dfn2}
and is equivariant with respect to the $H^2(X,Y;\Z_2)$-actions.
It is immediate that this map preserves~$w_2$ and satisfies~\eref{RelSpinPinCorr_e1}.
By~\eref{RelSpinPinCorr_e21b} and the equation after~\eref{SpinPinCorr_e},
it also satisfies the last condition of the RelSpinPin~\ref{RelSpinPinCorr_prop} 
property. 
\end{proof}

\begin{proof}[{\bf{\emph{Proof of RelSpinPin~\ref{RelSpinPinSES_prop} property}}}]
Let $V',V''$ be vector bundles over~$Y$.
With the notation as in~\eref{SpinPinSESdfn_e0} and~\eref{RelSpinPinStab_e21a}, 
we define the first map in~\eref{RelSpinPinSESdfn_e0} by
\begin{equation*}\begin{split}
&\bllrr{(E',\fo_{E'},\os_{E',V'}),(E'',\fo_{E''},\os_{E'',V''})}_{\oplus}\\
&\hspace{.8in}
=\big(E'\!\oplus\!E'',\fo_{E'}\fo_{E''},
\big\{\id_{E'|_{Y_2}}\!\oplus\!\Psi_{E''|_{Y_2},V'|_{Y_2}}\!\oplus\!\id_{V''|_{Y_2}}\big\}^{\!*}
\big(\llrr{\os_{E',V'},\os_{E'',V''}}_{\oplus}\big)\big)
\end{split}\end{equation*}
and the second map by the same expression with the $\OSpin$-structure $\os_{E'',V''}$
replaced by a $\Pin^{\pm}$-structure~$\fp_{E'',V''}$.
By~\eref{RelSpinPinStab_e21}, \eref{RelOspinDSassoc_e}, and
the SpinPin~\ref{SpinPinSES_prop}\ref{DSassoc_it} property in Section~\ref{SpinPinProp_subs},
these two maps are well-defined on the equivalence classes
of relative $\OSpin$- and $\Pin^{\pm}$-structures of Definition~\ref{RelPinSpin_dfn2}
and are bi-equivariant with respect to the $H^2(X,Y;\Z_2)$-actions.
It is immediate that the two maps satisfy~\eref{Relw2ses_e} and~\eref{RelSpinPinSES_e1} 
for $\OSpin$- and $\Pin^{\pm}$-structures on all vector bundles~$V'$ and~$V''$
over~$Y$.
By the first equation in the SpinPin~\ref{SpinPinSES_prop}\ref{DSEquivSum_it} property
and Corollary~\ref{GnH2_crl2b},  
they also satisfy  the first equation in 
the RelSpinPin~\ref{RelSpinPinSES_prop}\ref{RelDSEquivSum_it} property.
The remaining requirements on the two maps in the RelSpinPin~\ref{RelSpinPinSES_prop}
property follow from the corresponding requirements in
the SpinPin~\ref{RelSpinPinSES_prop} property.
\end{proof}

\begin{proof}[{\bf{\emph{Proof of RelSpinPin~\ref{vsSpinPin_prop} property}}}]
The reasoning in the proof of this property 
for the perspective of Definition~\ref{RelPinSpin_dfn3} in Section~\ref{RelSpinPindfn3_subs}
applies verbatim for the perspective of  Definition~\ref{RelPinSpin_dfn2}.
\end{proof}

\subsection{Equivalence of Definitions~\ref{RelPinSpin_dfn2} and~\ref{RelPinSpin_dfn3}}
\label{RelPinSpin2vs3_subs}

Let $(X,Y)$ be a CW pair and $V$ be a vector bundle over~$Y$.
Since the restriction homomorphism~\eref{RelSpinPinStr_e1a} is an isomorphism, 
the RelSpinPin~\ref{RelSpinPinObs_prop} property implies that
$V$ admits a relative $\Pin^{\pm}$-structure in the perspective of 
either Definition~\ref{RelPinSpin_dfn2} or~\ref{RelPinSpin_dfn3} if 
and only if the vector bundle $V|_{Y_2}$ over $Y_2\!\subset\!X_3$
admits a relative $\Pin^{\pm}$-structure in the perspective of 
either Definition~\ref{RelPinSpin_dfn2} or~\ref{RelPinSpin_dfn3}.
If $\fo\!\in\!\fO(V)$, then $(V,\fo)$ admits a relative $\Spin$-structure in the perspective of 
either Definition~\ref{RelPinSpin_dfn2} or~\ref{RelPinSpin_dfn3} if 
and only if the vector bundle $(V,\fo)|_{Y_2}$ over $Y_2\!\subset\!X_3$
admits a relative $\Spin$-structure in the perspective of 
either Definition~\ref{RelPinSpin_dfn2} or~\ref{RelPinSpin_dfn3}.
Furthermore, the natural restriction~maps
$$\cP_X^{\pm}(V)\lra\cP_{X_3}^{\pm}\!\big(V|_{Y_2}\big) \qquad\hbox{and}\qquad
\Sp_{\!X_3\!}(V,\fo)\lra\Sp_{X_3}^{\pm}\!\big((V,\fo)|_{Y_2}\big)$$
in either of the two perspectives
are equivariant with respect to the isomorphism~\eref{RelSpinPinStr_e1a}.
Along with the RelSpinPin~\ref{RelSpinPinStr_prop} property,
this implies that the two maps are bijections.
These bijections are intertwined with the identifications of
Theorem~\ref{SpinStrEquiv_thm} via the maps~\eref{vsSpinPin_e0} 
and respect all structures and correspondences of Section~\ref{RelSpinPinProp_subs}.\\

By the previous paragraph, it is thus sufficient to establish 
Theorem~\ref{RelSpinStrEquiv_thm}\ref{RelSpinStrEquiv2vs3_it} under the assumption 
that the dimensions of the CW complexes~$X$ and~$Y$ are at most~3 and~2, respectively. 
We can also restrict the consideration to vector bundles~$V$ over~$Y$ that admit
a relative $\Pin^{\pm}$-structure in the sense of Definition~\ref{RelPinSpin_dfn2} 
and  vector bundles~$(V,\fo)$ over~$Y$ that admit
a relative $\Spin$-structure in the sense of Definition~\ref{RelPinSpin_dfn2},
as appropriate.\\ 

Suppose $(V,\fo)$ is a rank $n\!\ge\!3$ oriented vector bundle over $Y\!\subset\!X$
and $\fs\!\equiv\!(E,\fo_E,\fs_{E,V})$ is a relative $\Spin$-structure on~$(V,\fo)$
in the sense of Definition~\ref{RelPinSpin_dfn2}.
We can assume that the rank of~$E$ is at least~3 as well.
Let $u\!:(\Si,\prt\Si)\!\lra\!(X,Y)$ be an element of~$\cL_X(Y)$ and
$\fs_u(E,\fo_E)$ be the mod~2 homotopy class of trivializations of $\{\prt u\}^*(E,\fo_E)$
as in~\eref{fsSiVfodfn_e}.
If $\prt\Si\neq\!\eset$,
we denote by $\wh\fs_u$ the mod~2 homotopy class of trivializations of $\{\prt u\}^*(V,\fo)$
so that the mod~2 homotopy class $\fs_u(E,\fo_E)\!\oplus\!\wh\fs_u$ of trivializations~of
\BE{RelPinSpin2vs3_e4}\{\prt u\}^*\big((E,\fo_E)|_Y\!\oplus\!(V,\fo)\!\big)=
\{\prt u\}^*(E,\fo_E)\!\oplus\!\{\prt u\}^*(V,\fo)\lra \prt\Si\EE
is the mod~2 homotopy class determined by $\{\prt u\}^*\fs_{E,V}$.
If $\prt\Si\!=\!\eset$, let
\BE{w2RelCorr_e0}
\wh\fs_u=\blr{w_2(E),u_*[\Si]_{\Z_2}}\in\Z_2\,.\EE

\vspace{.15in}

Suppose $F$ and $u_1,\ldots,u_m$ are as in~\eref{RelSpinPinStr_e5} and just above
and \hbox{$\wh{F}\!:\wh\Si\!\lra\!X$} is the continuous map obtained by gluing these maps
along the boundaries of their domains as in~\eref{RelSpinPinObs_e18a}.
Since $F$ is relatively equivalent to $u_1\!\sqcup\!\ldots\!\sqcup\!u_m$,
$$\fs_{\wh\Si}\big(\wh{F}^*(E,\fo_E)\!\big)\equiv 
\blr{w_2(\wh{F}^*E),[\wh\Si]_{\Z_2}}=\blr{w_2(E),\wh{F}_*[\wh\Si]_{\Z_2}}
=\blr{w_2(E),0}=0\,.$$
Along with Lemma~\ref{X2VBdisconn_lmm}, this implies that
$$\fs_{u_1}\!\big(E,\fo_E\!\big)\!\ldots\!\fs_{u_m}\!\big(E,\fo_E\!\big)
=\fs_{\Si}\big(F^*(E,\fo_E)\!\big)\,.$$
Since $\fs_{E,V}$ is a homotopy class of trivializations of $(E,\fo_E)|_Y\!\oplus\!(V,\fo)$
over~$Y$,
\begin{equation*}\begin{split}
\big(\fs_{u_1}\!(E,\fo_E)\!\ldots\!\fs_{u_m}\!(E,\fo_E)\big)\!\oplus\!
\big(\wh\fs_{u_1}\!\ldots\!\wh\fs_{u_m}\big)
&\equiv \big(\fs_{u_1}\!(E,\fo_E)\!\oplus\!\wh\fs_{u_1}\big)\!\ldots\!
\big(\fs_{u_m}\!(E,\fo_E)\!\oplus\!\wh\fs_{u_m}\big)\\
&= \fs_{\Si}\big(F^*\big((E,\fo_E)|_{Y_2}\!\oplus\!(V,\fo)\!\big)\!\big)\\
&=\fs_{\Si}\big(F^*(E,\fo_E)\!\big)\fs_{\Si}\big(F^*(V,\fo)\!\big).
\end{split}\end{equation*}
By the last two statements, 
$$\wh\fs_{u_1}\!\ldots\!\wh\fs_{u_m}=\fs_{\Si}\big(F^*(V,\fo)\!\big).$$
We conclude that a relative $\Spin$-structure~$\fs$ on~$(V,\fo)$ 
in the sense of Definition~\ref{RelPinSpin_dfn2}\ref{RelSpinStrDfn2_it} determines
a relative $\Spin$-structure 
$$(\wh\fs_u)_{u\in\cL_X(Y)}\equiv\Th_X(\fs)$$ 
on~$(V,\fo)$ in the sense of Definition~\ref{RelPinSpin_dfn3}.\\

Suppose $\eta\!\in\!H^2(X,Y;\Z_2)$.
Let $m$ and $h_{\eta;m}$ be as in~\eref{hetamdfn_e} and $u\!\in\!\cL_X(Y)$ be as above. 
By the bijectivity of~\eref{GnH2_e} and~\eref{RelSpinPinStr2_e1} with $(X,Y)\!=\!(\Si,\prt\Si)$, 
$\fs_u(h_{\eta;m}^*\wt\ga_m)$ is the mod~2 homotopy class of trivializations
of $\{\prt u\}^*(h_{\eta;m}^*\wt\ga_m)$ determined by $\{\prt u\}^*\fs_{h_{\eta;m}}$
if and only~if
\BE{RelPinSpin2vs3_e2}\blr{u^*\eta,[\Si]_{\Z_2}}\!\equiv\!
\blr{w_2(u^*(h_{\eta;m}^*\wt\ga_m)\!\big),[\Si]_{\Z_2}}=0\,.\EE
Thus,  the mod~2 homotopy class $\fs_u(h_{\eta;m}^*\wt\ga_m\!\oplus\!(E,\fo_E)\!)\!\oplus\!\wh\fs_u$
of trivializations~of
$$\{\prt u\}^*\big(\!\big(h_{\eta;m}^*\wt\ga_m\!\oplus\!(E,\fo_E)\!\big)|_Y\!\oplus\!(V,\fo)\!\big)=
\{\prt u\}^*(h_{\eta;m}^*)\!\oplus\!
\big(\{\prt u\}^*(E,\fo_E)\!\oplus\!\{\prt u\}^*(V,\fo)\big)$$
is the mod~2 homotopy class determined by $\{\prt u\}^*\llrr{\fs_{h_{\eta;m}},\fs_{E,V}}_{\oplus}$
if and only if~\eref{RelPinSpin2vs3_e2} holds.
Along with~\eref{Reletafsdfn2_e} and~\eref{RelSpinPinStr_e33}, this implies~that
$$\big(\wh{\eta\!\cdot\!\fs}\big)_u=\eta\cdot\wh\fs_u\,.$$
We conclude the construction of~$\Th_X$ above induces 
a natural $H^2(X,Y;\Z_2)$-equivariant map 
\BE{RelThdfn_e}\Th_X\!: \OSp_{X;2}(V)\lra\OSp_{X;3}(V)   \EE 
from the set of equivalence classes of relative $\OSpin$-structures on
the vector bundle~$V$ over \hbox{$Y\!\subset\!X$} 
in the sense of Definition~\ref{RelPinSpin_dfn2}\ref{RelSpinStrDfn2_it}
to the set of relative $\OSpin$-structures on~$V$ in the sense of Definition~\ref{RelPinSpin_dfn3}.
Along with the RelSpinPin~\ref{RelSpinPinStr_prop}\ref{RelSpinStr_it} property, 
this implies that~\eref{RelThdfn_e} is a bijection.
By~\eref{w2fsdfn_e} and~\eref{w2RelCorr_e0},
\BE{w2RelCorr_e} w_2\big(\Th_X(\os)\!\big)=w_2(\os)
\qquad\forall\,\os\!\in\!\OSp_{X;2}(V).\EE

\vspace{.15in}

Suppose~\eref{thmpf_e3} is a $\Spin$-structure on the oriented vector bundle
$(V,\fo)$ over~$Y$ in the sense of Definition~\ref{PinSpin_dfn}\ref{SpinStrDfn_it}.
Let~$\wt{s}_2$ be a section of $\Spin(V,\fo)$ as above~\eref{Th2dfn_e0} 
and~$\fs$ be the induced $\Spin$-structure on $(V,\fo)$ 
in the sense of Definition~\ref{PinSpin_dfn2}.
For each $\al\!\in\!\cL(Y)$, $\wt{s}_{\al}\!\equiv\!\al^*(q_V\!\circ\!\wt{s}_2)$
is then a section of the $\Spin(n)$-bundle $\al^*\Spin(V,\fo)$ in~\eref{thmpf_e11}
and determines a homotopy class~$\fs_{\al}$ of trivializations of~$\al^*(V,\fo)$.
With the notation as in~\eref{Th3dfn_e}, \eref{Th2dfn_e}, and~\eref{RelThdfn_e}, 
this implies that 
\begin{equation*}\begin{split}
\io_X\big(\Th_3\big(q_V\!:\Spin(V,\fo)\!\lra\!\SO(V,\fo)\!\big)\!\big)
&\equiv\io_X\big((\fs_{\al})_{\al\in\cL(Y)}\big)
=\Th_X\big(\io_X(\fs)\big)\\
&=\Th_X\big(\io_X\big(
\Th_2(q_V\!:\Spin(V,\fo)\!\lra\!\SO(V,\fo)\!)\!\big)\!\big).
\end{split}\end{equation*}
Since $\Th_2$ is a bijection, it follows that 
\BE{RelThdfn_e2} \Th_X\big(\io_X(\os)\!\big)=\io_X\big(\Th_3\big(\Th_2^{-1}(\os)\!\big)\!\big)
\qquad\forall\,\os\!\in\!\OSp_2(Y)\EE
for every vector bundle~$V$ of rank at least~3 over~$Y$.

\begin{proof}[{\bf{\emph{Compatibility with RelSpinPin~\ref{RelOSpinRev_prop} property}}}]
Suppose $(V,\fo)$ is an oriented vector bundle over~$Y$ of rank at least~3, 
$$\fs\equiv\big(E,\fo_E,\fs_{E,V}\big) \qquad\hbox{and}\qquad
\ov\fs\equiv\big(E,\fo_E,\ov{\fs_{E,V}}\big)$$
is a relative $\Spin$-structure on~$(V,\fo)$ in the sense of Definition~\ref{RelPinSpin_dfn2}
and the corresponding relative $\Spin$-structure on~$(V,\ov\fo)$, respectively.
Let $u\!\in\cL_X(Y)$, $\fs_u(E,\fo_E)$, and $\wh\fs_u$ be as before
and $\wh{\ov\fs}_u$ be the mod~2 homotopy class of trivializations of $\{\prt u\}^*(V,\ov\fo)$
induced by~$\fs_u(E,\fo_E)$ and the homotopy class~$\ov{\fs_{E,V}}$ of
trivializations of 
\BE{RelThXOSpinRev_e3}
\{\prt u\}^*\big((E,\fo_E)|_Y\!\oplus\!(V,\ov\fo)\!\big)=
\{\prt u\}^*(E,\fo_E)\!\oplus\!\{\prt u\}^*(V,\ov\fo)\lra \prt\Si\,.\EE
Thus, $\fs_u(E,\fo_E)\!\oplus\!\wh\fs_u$ and $\fs_u(E,\fo_E)\!\oplus\!\wh{\ov\fs}_u$ are
the mod~2 homotopy class of trivializations of~\eref{RelPinSpin2vs3_e4}  
determined by $\{\prt u\}^*\fs_{E,V}$ and 
the mod~2 homotopy class of trivializations of~\eref{RelThXOSpinRev_e3}  
determined by $\{\prt u\}^*\ov{\fs_{E,V}}$, respectively. 
Along with~\eref{OspinDSEquivSum_e}, this implies that
$\wh{\ov\fs}_u\!=\!\ov{\wh\fs}_u$ is the mod~2 homotopy class 
of trivializations of $\{\prt u\}^*(V,\ov\fo)$ corresponding to~$\wh\fs_u$.
Thus,
\BE{RelThXOSpinRev_e}\Th_X\big(\ov\os\big)=\ov{\Th_X(\os)}
\qquad\forall\,\os\!\in\!\OSp_{X;2}(V)\EE
for every vector bundle~$V$ of rank at least~3 over~$Y$.
\end{proof}

\begin{proof}[{\bf{\emph{Compatibility with RelSpinPin~\ref{RelSpinPinStab_prop} property}}}]
Suppose $(V,\fo)$ and $\fs$ are as above, $\Psi_{V',V''}$ is as in~\eref{RelSpinPinStab_e21a}, 
and
$$\fs_{E,V}'\equiv \big\{\Psi_{E|_Y,\tau_Y}\!\oplus\!\id_V\big\}^{\!*}
\!\big(\St_{E|_Y\oplus V}\!(\fs_{E,V})\!\big)$$
is the homotopy class of trivializations of $(E,\fo_E)|_Y\!\oplus\!\St(V,\fo)$ induced by~$\fs$.
Let $u\!\in\cL_X(Y)$, $\fs_u(E,\fo_E)$, and $\wh{\fs}_u$ be as before
and $\wh\fs_u'$ be the mod~2 homotopy class of trivializations of  $\{\prt u\}^*\St(V,\fo)$
induced by~$\fs_u(E,\fo_E)$ and the homotopy class~$\fs'_{E,V}$ of trivializations~of 
\BE{RelSpinPinStab_e3}
\{\prt u\}^*\big((E,\fo_E)|_Y\!\oplus\!\St(V,\fo)\!\big)=
\{\prt u\}^*(E,\fo_E)\!\oplus\!\{\prt u\}^{\!*}\!\big(\St(V,\fo)\!\big)\lra \prt\Si\,.\EE
Thus, $\fs_u(E,\fo_E)\!\oplus\!\wh\fs_u$ and $\fs_u(E,\fo_E)\!\oplus\!\wh{\fs}_u'$ are
the mod~2 homotopy class of trivializations of~\eref{RelPinSpin2vs3_e4}  
determined by $\{\prt u\}^*\fs_{E,V}$ and 
the mod~2 homotopy class of trivializations of~\eref{RelSpinPinStab_e3}  
determined by $\{\prt u\}^*\fs'_{E,V}$, respectively. 
Along with~\eref{OspinDSstab_e} and~\eref{OspinDSassoc_e}, this implies that
$\wh{\fs}'_u\!=\!\St_V(\wh\fs_u)$ is the mod~2 homotopy class of trivializations 
of $\{\prt u\}^*\St(V,\fo)$ corresponding to~$\wh\fs_u$.
Thus,
\BE{RelSpinPinStab5_e}\Th_X\big(\St_V(\os)\!\big)=\St_V\big(\Th_X(\os)\!\big)
\qquad\forall\,\os\!\in\!\OSp_{X;2}(V)\EE
for every vector bundle~$V$ of rank at least~3 over~$Y$.
\end{proof}

We define the map~$\Th_X$ in~\eref{RelThdfn_e} for rank~2 vector bundles $V$ over $Y\!\subset\!X$ 
and then for rank~1 vector bundles $V$ over~$Y$ by~\eref{RelSpinPinStab5_e}.
By the RelSpinPin~\ref{RelSpinPinStab_prop} property, the already established 
properties of~\eref{RelThdfn_e} for vector bundles of ranks at least~3,
and Theorem~\ref{SpinStrEquiv_thm}, 
the resulting maps~$\Th_X$ for vector bundles of ranks~1 and~2 are natural 
$H^2(X,Y;\Z_2)$-equivariant bijections which intertwine the identifications of
Theorem~\ref{SpinStrEquiv_thm} via the maps~\eref{vsSpinPin_e0} and satisfy~\eref{w2RelCorr_e}.
These maps are compatible with the RelSpinPin~\ref{RelSpinPinStab_prop} property,
i.e.~satisfy~\eref{RelSpinPinStab5_e}, by definition.
Along with the first equality in~\eref{RelSpinPinStab_e2} and~\eref{RelThXOSpinRev_e} 
for vector bundles of rank at least~3, this implies that 
the maps~$\Th_X$ satisfy~\eref{RelThXOSpinRev_e} for all vector bundles~$V$ over~$Y$.

\begin{proof}[{\bf{\emph{Compatibility with RelSpinPin~\ref{RelSpinPinSES_prop} property}}}]
Suppose $(V',\fo')$ and $(V'',\fo'')$ are oriented vector bundles over~$Y$ of ranks at least~3,
$$\fs'\equiv\big(E',\fo_{E'},\fs_{E',V'}\big) \qquad\hbox{and}\qquad
\fs''\equiv\big(E'',\fo_{E''},\fs_{E'',V''}\big)$$
are relative $\Spin$-structures on $(V',\fo')$ and $(V'',\fo'')$, respectively,
$\Psi_{V',V''}$ is as in~\eref{RelSpinPinStab_e21a},  and
$$\fs_{E,V}\equiv
\big\{\id_{E'|_Y}\!\oplus\!\Psi_{E''|_Y,V'}\!\oplus\!\id_{V''}\big\}^{\!*}
\big(\fs_{E',V'}\!\oplus\!\fs_{E'',V''}\big)$$
is the homotopy class of trivializations of 
$$(E,\fo_E)|_Y\!\oplus\!(V,\fo) \equiv
\big(\!(E',\fo_{E'})\!\oplus\!(E'',\fo_{E''})\!\big)\big|_Y\!\oplus\!
\big(\!(V',\fo')\!\oplus\!(V'',\fo'')\!\big)$$
induced by~$\fs'$ and~$\fs''$.
For $u\!\in\cL_X(Y)$, let $\fs_u(E',\fo_{E'})$ and $\fs_u(E'',\fo_{E''})$
be the mod~2 homotopy class of trivializations~of
${\prt u}^*(E',\fo_{E'})$ and ${\prt u}^*(E'',\fo_{E''})$, respectively,
as in~\eref{fsSiVfodfn_e}.\\

Let $\wh{\fs}'_u$, $\wh{\fs}''_u$, and  $\wh{\fs}_u$
be the mod~2 homotopy classes of trivializations 
of $\{\prt u\}^*(V',\fo')$ induced by $\fs_u(E',\fo_{E'})$ and~$\fs_{E',V'}$,
of $\{\prt u\}^*(V'',\fo'')$ induced by $\fs_u(E'',\fo_{E''})$ and~$\fs_{E'',V''}$,
and of $\{\prt u\}^*(V,\fo)$ induced by 
$$\fs_u\big(E,\fo_E\big)= \fs_u\big(E',\fo_{E'}\big)\!\oplus\!\fs_u\big(E'',\fo_{E''}\big)$$ 
and the homotopy class $\fs_{E',V'}\!\oplus\!\fs_{E'',V''}$ of trivializations~of 
\BE{RelSpinPinSES_e3}
\{\prt u\}^*\big(\!(E,\fo_E)|_Y\!\oplus\!(V,\fo)\!\big)=
\{\prt u\}^*(E',\fo_{E'})\!\oplus\!\{\prt u\}^*(E'',\fo_{E''})\!\oplus\!
\{\prt u\}^{\!*}(V',\fo')\!\oplus\!\{\prt u\}^{\!*}(V'',\fo''),\EE
respectively.
Thus, 
$$\fs_u(E',\fo_{E'})\!\oplus\!\wh\fs_u',\quad 
\fs_u(E'',\fo_{E''})\!\oplus\!\wh\fs_u'',\quad\hbox{and}\quad
\big(\fs_u(E',\fo_{E''})\!\oplus\!\fs_u(E'',\fo_{E''})\!\big)\!\oplus\!\wh{\fs}_u$$ 
are the mod~2 homotopy classes of trivializations
of $\{\prt uW\}^*((E',\fo_{E'})|_Y\!\oplus\!(V',\fo')\!)$ 
 determined by $\{\prt u\}^*\fs_{E',V'}$,
of $\{\prt u\}^*((E'',\fo_{E''})|_Y\!\oplus\!(V'',\fo'')\!)$ 
 determined by $\{\prt u\}^*\fs_{E'',V''}$, and
of~\eref{RelSpinPinSES_e3} determined by 
$\{\prt u\}^*(\fs_{E',V'}\!\oplus\!\fs_{E'',V''})$, respectively.
Along with~\eref{OspinDSassoc_e}, this implies that
$\wh{\fs}_u\!=\!\wh{\fs}_u'\!\oplus\!\wh{\fs}_u''$ 
is the mod~2 homotopy class of trivializations 
of $\{\prt u\}^*(V,\fo)$ determined by $\wh{\fs}_u'$ and~$\wh{\fs}_u''$.
Thus,
\BE{RelSpinPinSES_e5}\Th_X\big(\llrr{\os',\os''}_{\oplus}\big)=
\bllrr{\Th_X(\os'),\Th_X(\os'')}_{\oplus}
\quad\forall\,\os'\!\in\!\OSp_{X;2}(V'),\,\os''\!\in\!\OSp_{X;2}(V'')\EE
for all vector bundles~$V'$ and~$V''$ of ranks at least~3 over~$Y$.
This identity for arbitrary rank bundles then follows from 
the RelSpinPin~\ref{RelSpinPinStab_prop} property, \eref{RelOspinDSstab_e}, 
and~\eref{RelOspinDSassoc_e}.
\end{proof}

With $\cP_{X;2}^{\pm}(V)$ and $\cP_{X;3}^{\pm}(V)$ denoting the sets~$\cP_X^{\pm}(V)$
of the relative $\Pin^{\pm}$-structures on a vector bundle~$V$ over~$Y$
in the perspectives of Definitions~\ref{RelPinSpin_dfn2} and~\ref{RelPinSpin_dfn3}, 
respectively, define
\BE{RelThXdfn2_e}\Th_X\!: \cP_{X;2}^{\pm}(V)\lra\cP_{X;3}^{\pm}(V) 
\qquad\hbox{by}\qquad  
\Co_V^{\pm}\big(\Th_X(\fp)\!\big)=\Th_X\big(\Co_V^{\pm}(\fp)\!\big)\,.\EE 
By the RelSpinPin~\ref{RelSpinPinCorr_prop} property, the already established 
properties of~\eref{RelThdfn_e}, and Theorem~\ref{SpinStrEquiv_thm}, 
\eref{RelThXdfn2_e} is a natural $H^2(X,Y;\Z_2)$-equivariant bijection which
intertwines the identifications of
Theorem~\ref{SpinStrEquiv_thm} via the maps~\eref{vsSpinPin_e0}.
It preserves~$w_2$ and is 
compatible with the RelSpinPin~\ref{RelSpinPinStab_prop} property, 
i.e.~the map~$\Th_X$ satisfies~\eref{w2RelCorr_e} and~\eref{RelSpinPinStab5_e} 
for relative $\Pin^{\pm}$-structures~$\fp$ in place of 
the relative $\OSpin$-structures~$\os$.
This map is compatible with the RelSpinPin~\ref{RelSpinPinCorr_prop} property,
i.e.~satisfies the second equality in~\eref{RelThXdfn2_e}, by definition.
By the second statement in the RelSpinPin~\ref{RelSpinPinSES_prop}\ref{RelDSEquivSum_it} property 
and~\eref{RelSpinPinSES_e5}, 
\eref{RelThXdfn2_e} is compatible with the second map in~\eref{RelSpinPinSESdfn_e0},
i.e.~\eref{RelSpinPinSES_e5} holds for relative $\Pin^{\pm}$-structures~$\fp''$ in place of 
the relative $\OSpin$-structures~$\os''$.
By the RelSpinPin~\ref{RelSpinPinSES_prop}\ref{RelDSPin2Spin_it} property,
\eref{RelThXdfn2_e} is compatible with the RelSpinPin~\ref{RelPin2SpinRed_prop} property,
i.e.
$$\Th_X\big(\fR_{\fo}^{\pm}(\fp)\!\big)=\fR_{\fo}^{\pm}\big(\Th_X(\fp)\!\big) 
\qquad\forall\,\fp\!\in\!\cP^{\pm}_{X;2}(V),\,\fo\!\in\!\fO(V)$$
for every vector bundle $V$ over~$Y$.
This concludes the proof of Theorem~\ref{RelSpinStrEquiv_thm}\ref{RelSpinStrEquiv2vs3_it}.

\chapter{Orientations for Real CR-Operators}
\label{CROrient_ch}

The so-called \sf{real Cauchy-Riemann operators} 
(or \sf{CR-operators} for short)
are central to Gromov-Witten theory and the Fukaya category literature.
As demonstrated in \cite[Section~8.1]{FOOO}, 
relative $\Spin$-structures induce orientations on the determinants of such operators.
As demonstrated in \cite[Section~3]{Sol}, 
relative $\Pin$-structures induce orientations on certain twisted determinants 
of real CR-operators.
Section~\ref{CRintro_sec} collects properties of the induced orientations in a ready-to-use format.
The construction of these orientations and the verification of 
their claimed properties are detailed in Sections~\ref{BaseCR_sec}-\ref{CROrientPf_sec}.
We formulate all statements in terms of \sf{real bundle pairs}~$(V,\vph)$ over 
symmetric surfaces~$(\Si,\si)$ with separating fixed loci~$\Si^{\si}$.
Analogous statements for vector bundle pairs~$(E,F)$ over bordered surfaces~$(\Si^b,\prt\Si^b)$
follow immediately by doubling $(E,F)$ to a real bundle pair 
as in \cite[Section~3]{XCapsSetup}.\\

Our construction of orientations in the relative $\Spin$ case differs from~\cite{FOOO} 
and is instead inspired by the approach of \cite[Section~6]{Melissa};
this approach is now the standard perspective on orienting real CR-operators in this case.
Our construction of orientations in the relative $\Pin$ case differs from that in~\cite{Sol}
and is instead motivated by the approach of \cite[Section~3]{Ge1}.
This approach combines the $\Spin$ and rank~1 cases and eliminates 
the need for an arbitrary choice of 
a distinguished $\Pin$-structure on a vector bundle over~$S^1$ in each rank.
However, the way we combine  the $\Spin$  and rank~1 cases differs 
from~\cite{Ge1}.\\

Section~\ref{BaseCR_sec} constructs orientations of determinants of
real CR-operators on rank~1 real bundle pairs over~$S^2$
with the standard orientation-reversing involution~$\tau$
and on even-degree real bundle pairs~$(V,\vph)$ over~$(S^2,\tau)$.
As in~\cite{Ge1}, the orientations~$\fo(V,\vph;\fo_{x_1})$ constructed in the former case
depend on the choice of orientation~$\fo_{x_1}$ of~$V^{\vph}_{x_1}$ of the real line bundle~$V^{\vph}$
at a point~$x_1$ in the $\tau$-fixed locus~$S^1$ and 
the choice of half-surface~$\bD^2_+$ of~$(S^2,\tau)$.
As in~\cite{FOOO,Melissa}, the orientations~$\fo_{\os}(V,\vph)$ 
constructed in the latter case depend on
the choice of relative $\OSpin$-structure~$\os$ on the real vector bundle~$V^{\vph}$
over \hbox{$S^1\!\subset\!S^2$} and the choice of half-surface~$\bD^2_+$ of~$(S^2,\tau)$.
We recall the orientation constructions of~\cite{Ge1} and~\cite{Melissa}
in these two cases and establish key properties of the resulting orientations.
Section~\ref{BaseCR_sec} contains the most technical arguments of this chapter.\\

Section~\ref{InterCR_sec} constructs orientations of determinants of
real CR-operators on rank~1 real bundle pairs and on even-degree real bundle pairs~$(V,\vph)$
over arbitrary smooth symmetric surfaces~$(\Si,\si)$ with separating fixed loci~$\Si^{\si}$.
These orientations are obtained from the orientations constructed in 
the two corresponding cases in Section~\ref{BaseCR_sec}.
Similarly to Section~\ref{BaseCR_sec}, 
the orientations~$\fo(V,\vph;\fo_{\x})$ in the former case
depend on the choices of orientation~$\fo_{x_r}$ of~$V^{\vph}_{x_r}$ 
at a point~$x_r$ in each component $S^1_r$ of the $\si$-fixed locus~$\Si^{\si}$,
half-surface~$\Si^b$ of~$(\Si,\si)$,  and ordering on~$\pi_0(\Si^{\si})$.
The orientations~$\fo_{\os}(V,\vph)$ in the latter case depend on
the choices of relative $\OSpin$-structure~$\os$ on~$V^{\vph}$,
half-surface~$\Si^b$ of~$(\Si,\si)$,  and ordering on~$\pi_0(\Si^{\si})$.
We again recall the orientation constructions of~\cite{Ge1} and~\cite{Melissa}
in these two cases and establish key properties of the resulting orientations.
These properties are deduced from the analogous properties in Section~\ref{BaseCR_sec}
primarily through algebraic considerations involving commutative squares of 
short exact sequences.
Section~\ref{InterCR_sec} concludes the verification of the properties of 
orientations of the determinants of real CR-operators induced by 
relative $\OSpin$-structures stated in Section~\ref{CRintro_sec}.\\

Section~\ref{CROrientPf_sec} finally constructs orientations on the determinants 
and twisted determinants of real CR-operators induced by relative $\Pin$-structures.
It uses the RelSpinPin~\ref{RelSpinPinCorr_prop} property on page~\pageref{RelSpinPinCorr_prop} 
to combine the orientations constructed in Section~\ref{InterCR_sec}.
Section~\ref{CROrientPf_sec} also deduces 
the properties of the orientations induced by relative $\Pin$-structures stated 
in Section~\ref{CRintro_sec} from the analogous properties established
in the two settings of Section~\ref{InterCR_sec}.

\section{Introduction}
\label{CRintro_sec}

We assemble the key notions concerning real CR-operators and their determinants 
and state the main (and only) theorem of this chapter in Section~\ref{NRS_subs}.
The properties of orientations of real CR-operators described in 
Sections~\ref{OrientPrp_subs1} and~\ref{OrientPrp_subs2}
include the compatibility of the orientations
with short exact sequences of real bundle pairs and with flat degenerations of the domains.
Some consequences of these properties are deduced in Section~\ref{OrientEg_subs}.

\subsection{Definitions and main theorem}
\label{NRS_subs}

An \sf{exact triple} (short exact sequence)\gena{Fredholm operator!exact triple}
$$0\lra D'\lra D\lra D''\lra 0$$
of Fredholm operators is a commutative diagram
\BE{cTdiag_e}\begin{split}
\xymatrix{ 0 \ar[r]&  X'\ar[d]^{D'}\ar[r]& X\ar[d]^D\ar[r]&  X''\ar[d]^{D''}\ar[r]& 0\\
0 \ar[r]&  Y'\ar[r]& Y\ar[r]&  Y''\ar[r]& 0}
\end{split}\EE
so that the rows are exact sequences of bounded linear homomorphisms between 
Banach vector spaces (over~$\R$) and the columns are Fredholm operators.
If $X''$ is a finite-dimensional vector space and $Y''$ is the zero vector space,
we abbreviate the exact triple~\eref{cTdiag_e} of Fredholm operators as
$$0\lra D'\lra D\lra X''\lra 0.$$

\vspace{.15in}

If $X,Y$ are Banach spaces and $D\!:X\!\lra\!Y$ is a Fredholm operator, let
$$\la(D)\equiv\la(\ker D) \otimes \la(\cok\,D)^*\nota{laD@$\la(D)$}$$
denote the \sf{determinant}\gena{Fredholm operator!determinant} 
of~$D$.
An \sf{orientation}\gena{Fredholm operator!orientation} 
on a Fredholm operator~$D$ 
is an orientation on 
the one-dimensional real vector space~$\la(D)$ or equivalently 
a homotopy class of isomorphisms of~$\la(D)$ with~$\R$.
If $\cH$ is a topological space, a continuous family 
$$\cD\equiv \big\{D_{\t}\!:\t\!\in\!\cH\big\}$$
of Fredholm operators~$D_{\t}$ over~$\cH$  
determines a line bundle~$\la(\cD)$ over~$\cH$, called 
\sf{the determinant line bundle of~$\cD$}\gena{determinant line bundle}; 
see~\cite{detLB}. 
An exact triple~$\ft$ of Fredholm operators as in~\eref{cTdiag_e}
determines a canonical isomorphism
\BE{sum} \Psi_{\ft}\!: \la(D')\otimes\la(D'')\stackrel{\approx}{\lra} \la(D).\EE
For a continuous family of exact triples of Fredholm operators, 
the isomorphisms~\eref{sum} give rise to a canonical isomorphism
between the determinant line bundles; see the Exact Triples property in \cite[Section~2]{detLB}.\\

\noindent
For $k\!\in\!\Z^{\ge0}$, let $[k]\!=\!\{1,\ldots,k\}$.
An \sf{involution}\gena{involution} on a topological space~$X$ is a homeomorphism
$\phi\!:X\!\lra\!X$ such that $\phi\!\circ\!\phi\!=\!\id_X$.
A \sf{symmetric surface}\gena{symmetric surface} 
$(\Si,\si)$ is a closed  oriented (possibly nodal) surface~$\Si$  
with an orientation-reversing involution~$\si$.
If $\Si$ is smooth, the fixed locus~$\Si^{\si}$ of~$\si$ is a disjoint union of circles.
In general, $\Si^{\si}$ consists of isolated points
(called \sf{$E$~nodes}\gena{node!type $E$} in \cite[Section~3.2]{Melissa})
and circles identified at pairs of points (called \sf{$H$~nodes}\gena{node!type $H$} 
in~\cite{Melissa}).
The remaining nodes of~$\Si$, which we call \sf{$C$~nodes}, come in pairs 
$\nod_{ij}^{\pm}\!\not\in\!\Si^{\si}$ interchanged by the involution~$\si$.
A \sf{complex structure} on~$(\Si,\si)$ is a complex structure~$\fj$ on~$\Si$ 
compatible with the orientation of~$\Si$ such~that $\si^*\fj\!=\!-\fj$.
A \sf{symmetric surface with $k$ real marked points and $l$ conjugate pairs of marked 
points} is a~tuple 
\BE{cCsymdfn_e} \cC\equiv\big(\Si,\si,(x_i)_{i\in[k]},(z_i^+,z_i^-)_{i\in[l]}\big),\EE
where $(\Si,\si)$ is a symmetric surface and 
$x_1,\ldots,x_k$, $z_1^+,\ldots,z_l^+$, and $z_1^-,\ldots,z_l^-$ are distinct smooth 
points of~$\Si$ such that $\si(x_i)\!=\!x_i$ for all $i\!\in\![k]$ and 
$\si(z_i^+)\!=\!z_i^-$ for all $i\!\in\![l]$.\\

\noindent
There are two basic kinds of symmetric surfaces: 
connected ones and pairs of connected surfaces interchanged by the involution,
called \sf{doublets} in~\cite{RealGWsII}.
Every symmetric surface~$(\Si,\si)$ is a union of such surfaces,
which we then call the \sf{real elemental components}\gena{elemental component!real} 
and 
the \sf{conjugate elemental components}\gena{elemental component!conjugate} 
of~$(\Si,\si)$, respectively.
A bordered surface $\Si^b$ \sf{doubles} to a symmetric surface $(\Si,\si)$ so~that 
$$\Si^b\subset\Si, \qquad \Si=\Si^b\!\cup\!\si(\Si^b),
 \qquad\hbox{and}\qquad \Si^b\!\cap\!\si(\Si^b)=\Si^{\si}=\prt\Si^b\,.$$
In such a case, we call $\Si^b$ a \sf{half-surface}\gena{half-surface} 
of~$(\Si,\si)$ and orient 
$\Si^{\si}$ as the boundary of~$\Si^b$ as in \cite[p146]{Warner}. 
A smooth symmetric surface $(\Si,\si)$ admits a half-surface if and only if
$\Si^{\si}$ separates~$\Si$ into two halves interchanged by~$\si$,
i.e.~the fixed locus of~$(\Si,\si)$ is \sf{separating}\gena{separating fixed locus}.
A half-surface of a doublet~$(\Si,\si)$ is either of the two connected components of~$(\Si,\si)$.
If the fixed locus of a connected elemental symmetric surface~$(\Si,\si)$ is separating,
$(\Si,\si)$ again contains two half-surfaces.
A choice of a half-surface for an arbitrary symmetric surface~$(\Si,\si)$ corresponds
to a choice of a half-surface for each elemental component of~$(\Si,\si)$.\\

\noindent
Let $(X,\phi)$ be a topological space with an involution.
A \sf{conjugation} on a complex vector bundle $V\!\lra\!X$ 
\sf{lifting} the involution~$\phi$ is a vector bundle homomorphism 
$\vph\!:V\!\lra\!V$ covering~$\phi$ (or equivalently 
a vector bundle homomorphism  $\vph\!:V\!\lra\!\phi^*V$ covering~$\id_X$)
such that the restriction of~$\vph$ to each fiber is anti-complex linear
and $\vph\!\circ\!\vph\!=\!\id_V$.
A \sf{real bundle pair}\gena{real bundle pair} $(V,\vph)\!\lra\!(X,\phi)$   
consists of a complex vector bundle $V\!\lra\!X$ and 
a conjugation~$\vph$ on $V$ lifting~$\phi$.
For example, 
$$\big(X\!\times\!\C^n,\phi\!\times\!\fc\big)\lra(X,\phi),$$
where $\fc\!:\C^n\!\lra\!\C^n$ is the standard conjugation on~$\C^n$,
is a real bundle pair.
In general, \hbox{$V^{\vph}\!\lra\!X^{\phi}$} is a real vector bundle with
\hbox{$\rk_{\R}V^{\vph}\!=\!\rk_{\C}V$}.\\

\noindent
Let $(V,\vph)$ be a real bundle pair over a symmetric surface~$(\Si,\si)$.
A \sf{real Cauchy-Riemann} (or \sf{CR-}) \sf{operator}\gena{real CR-operator} 
on~$(V,\vph)$  
is a linear map of the~form
\BE{CRdfn_e}\begin{split}
D=\bp\!+\!A\!: \Ga(\Si;V)^{\vph}
\equiv&\big\{\xi\!\in\!\Ga(\Si;V)\!:\,\xi\!\circ\!\si\!=\!\vph\!\circ\!\xi\big\}\\
&\lra
\Ga_{\fj}^{0,1}(\Si;V)^{\vph}\equiv
\big\{\ze\!\in\!\Ga(\Si;(T^*\Si,\fj)^{0,1}\!\otimes_{\C}\!V)\!:\,
\ze\!\circ\!\nd\si=\vph\!\circ\!\ze\big\},
\end{split}\EE
where $\bp$ is the holomorphic $\bp$-operator for some holomorphic structure in~$V$ 
which lifts a complex structure~$\fj$ on~$(\Si,\si)$ and is reversed  by~$\vph$ 
and  
$$A\in\Ga\big(\Si;\Hom_{\R}(V,(T^*\Si,\fj)^{0,1}\!\otimes_{\C}\!V) \big)^{\vph}$$ 
is a zeroth-order deformation term. 
The completion of a real CR-operator~$D$ on~$(V,\vph)$
with respect to appropriate norms on its domain and target 
(e.g.~Sobolev $L^p_1$ and~$L^p$-norms, respectively, with $p\!>\!2$), 
which we also denote by~$D$, is Fredholm;
see the proof of \cite[Proposition~3.6]{XCapsSetup}.
If $\Si$ is a smooth connected surface of genus~$g$, then the index of~$D$ is given~by
$$\tn{ind}\,D=(1\!-\!g)\,\rk\,V+\blr{c_1(V),[\Si]}\,.$$
The norms can be chosen so that the evaluation~map
$$\Ga(\Si;V)^{\vph}\lra V_x, \qquad \xi\lra\xi(x),$$
is continuous for each $x\!\in\!\Si$.\\

The space of completions of all real CR-operators on~$(V,\vph)$ is contractible 
with respect to the operator~norm.
This implies that there is a canonical homotopy class of isomorphisms 
between the determinant lines of any two real CR-operators on a real bundle pair~$(V,\vph)$;
we thus denote any such line by~$\la(D_{(V,\vph)})$.
If in addition $\cC$ is a marked symmetric surface as in~\eref{cCsymdfn_e}, let
$$\la_{\cC}^{\R}(V,\vph)=
\la\bigg(\bigoplus_{i=1}^k\!V_{x_i}^{\vph}\bigg)
=\bigotimes_{i=1}^k\!\la(V_{x_i}^{\vph}), \qquad
\wt\la_{\cC}(D)=\la_{\cC}^{\R}(V,\vph)^*\!\otimes\!\la(D);$$
the summands and the factors in the definition of $\la_{\cC}^{\R}(V,\vph)$ 
are {\it not} ordered.
We call $(V,\vph)$ 
\sf{$\cC$-balanced}\gena{Cbalanced@$\cC$-balanced} 
if the parity of the number~$k_r$ of 
the real marked points~$x_i$ carried by every topological component $\Si^{\si}_r$
of~$\Si^{\si}$ is {\it different} from $\lr{w_1(V^{\vph}),[\Si^{\si}_r]_{\Z_2}}$.\\

We call a symmetric surface $(\Si,\si)$ with separating fixed locus
and with choices of a half-surface~$\Si^b$ 
and of an ordering  of the topological components of~$\Si^{\si}$
a \sf{decorated} symmetric surface.

\begin{thm}\label{CROrient_thm}
Suppose $(\Si,\si)$ is a smooth decorated symmetric surface 
and $(V,\vph)$ is a  real bundle pair over~$(\Si,\si)$.
\begin{enumerate}[label=(\alph*),leftmargin=*]

\item\label{CROrientSpin_it} A relative $\OSpin$-structure $\os$ on the vector bundle $V^{\vph}$ over 
$\Si^{\si}\!\subset\!\Si$ determines an orientation~$\fo_{\os}(V,\vph)$ 
of~$\la(D_{(V,\vph)})$.

\item\label{CROrientPin_it} If $\cC$ is a marked symmetric surface with 
the underlying symmetric surface~$(\Si,\si)$ and $(V,\vph)$ is $\cC$-balanced, 
a relative $\Pin^{\pm}$-structure $\fp$ on~$V^{\vph}$ determines 
an orientation $\fo_{\cC;\fp}(V,\vph)$ of $\wt\la_{\cC}(D_{(V,\vph)})$.

\end{enumerate}
These orientations satisfy all properties of 
Sections~\ref{OrientPrp_subs1} and~\ref{OrientPrp_subs2}.
\end{thm}

We denote by $\la_{\C}(V)$ the top exterior power of a complex vector bundle~$V$.
A conjugation~$\vph$ on~$V$ induces a conjugation $\la_{\C}(\vph)$ on~$\la_{\C}(V)$.
For a real bundle pair $(V,\vph)$, we define
\BE{RBPstab_e}\la(V,\vph)=\big(\la_{\C}(V),\la_{\C}(\vph)\!\big),\quad
(V,\vph)_{\pm}\equiv\big(V_{\pm},\vph_{\pm}\!\big)
=(V,\vph)\oplus(2\!\pm\!1)\la(V,\vph).\EE

\vspace{.15in}

Let $\cC$ be a smooth decorated marked symmetric surface and 
$(V,\vph)$ be a $\cC$-balanced real bundle pair over~$\cC$.
By Proposition~\ref{cOpincg_prp}, 
the line $\wt\la_{\cC}(D_{(V,\vph)})$
has a natural orientation $\fo_{\cC}(V,\vph)$ if the rank of~$(V,\vph)$ is~1.
If $(V,\vph)$ is of any rank, a relative  $\Pin^{\pm}$-structure $\fp$ on~$V^{\vph}$
is a relative $\OSpin$-structure 
$$\os_{\pm}\equiv\Co_V^{\pm}(\fp)$$
on~$V_{\pm}^{\vph}$; see the RelSpinPin~\ref{RelSpinPinCorr_prop} property 
on page~\pageref{RelSpinPinCorr_prop}.
We define the orientation $\fo_{\cC;\fp}(D_{(V,\vph)})$ 
of Theorem~\ref{CROrient_thm}\ref{CROrientPin_it} 
at the beginning of Section~\ref{twisteddetLB_subs} via the isomorphism
\BE{PinOrientDfn_e}\la\big(D_{(V,\vph)_{\pm}}\big)\approx \la\big(D_{(V,\vph)}\big)\!\otimes\!
\la\big(D_{\la(V,\vph)}\big)^{\otimes(2\pm1)}\EE
as in~\eref{sum} from the orientation $\fo_{\os_{\pm}}((V,\vph)_{\pm})$
as in Theorem~\ref{CROrient_thm}\ref{CROrientSpin_it} and the orientation 
$\fo_{\cC}(\la(V,\vph))$;
in the plus case, the tensor product of the last two factors of $\la(D_{\la(V,\vph)})$
above has a canonical orientation.
Most properties of $\fo_{\cC;\fp}(V,\vph)$ then follow directly from
the corresponding properties of $\fo_{\os_{\pm}}((V,\vph)_{\pm})$
and~$\fo_{\cC}(V,\vph)$.

\subsection{Properties of orientations: smooth surfaces}
\label{OrientPrp_subs1}

Let $(\Si,\si)$ be a smooth symmetric surface.
If the genera of the connected components of~$\Si$ are $g_1,\ldots,g_N$,
we define 
$$g(\Si) \nota{g0si@$g(\Si)$}
\equiv 1+\!\sum_{q=1}^N\!\big(g_q\!-\!1\big)$$
to be the \sf{genus} of~$\Si$.
If the fixed locus $\Si^{\si}\!\subset\!\Si$ is separating, 
then $g(\Si)\!+\!1\!$ has the same parity as the number $|\pi_0(\Si^{\si})|$ 
of connected components of~$\Si^{\si}$.
For a rank~$n$ real bundle pair $(V,\vph)$ over~$(\Si,\si)$,
we denote by 
$$\deg V\equiv\blr{c_1(V),[\Si]_{\Z}}\in\Z \quad\hbox{and}\quad
W_1(V,\vph)\nota{w1Vvph@$W_1(V,\vph)$}
\equiv\big|\big\{S_r^1\!\in\!\pi_0(\Si^{\si})\!:
w_1(V^{\vph})|_{S_r^1}\!\neq\!0\big\}\big|\in\Z^{\ge0}$$
the \sf{degree} of~$V$ and the number of connected components of~$\Si^{\si}$ 
over which $V^{\vph}$ is not orientable, respectively.
By \cite[Proposition~4.1]{BHH}, these two numbers have the same parity.
Thus, the elements
\begin{equation*}\begin{split}
\vp_{\os}(\Si_*)\nota{wpos@$\vp_{\os}(\Si_*)$}
 &\equiv 
n\,\frac{g(\Si_*)\!+\!|\pi_0(\Si_*^{\si})|\!-\!1}{2}+\frac{\deg V|_{\Si_*}}{2}
+\blr{w_2(\os),[\Si_*]_{\Z_2}}\in\Z_2 \qquad\hbox{and}\\
\vp_{\fp}(\Si_*) \nota{wpp@$\vp_{\fp}(\Si_*)$} 
&\equiv n\,\frac{g(\Si_*)\!+\!|\pi_0(\Si_*^{\si})|\!-\!1}{2}
+\frac{\deg V|_{\Si_*}\!\pm\!W_1((V,\vph)|_{\Si_*})}{2}
+\blr{w_2(\fp),[\Si_*]_{\Z_2}}\in \Z_2
\end{split}\end{equation*}
are well-defined for every relative $\OSpin$-structure~$\os$ 
on the real vector bundle~$V^{\vph}$ over $\Si^{\si}\!\subset\!\Si$,
relative $\Pin^{\pm}$-structure~$\fp$ on~$V^{\vph}$, and 
an elemental component $\Si_*$ of~$(\Si,\si)$.\\

A \sf{real CR-operator} on 
a complex vector bundle~$V$ over a closed oriented (possibly nodal) surface~$\Si$ 
is a linear map of the~form
$$D=\bp\!+\!A\!: \Ga(\Si;V)\lra
\Ga_{\fj}^{0,1}(\Si;V)\!\equiv\!\Ga\big(\Si;(T^*\Si,\fj)^{0,1}\!\otimes_{\C}\!V\big),$$
where $\bp$ is the holomorphic $\bp$-operator for some holomorphic structure in~$V$ 
which lifts a complex structure~$\fj$ on~$(\Si,\si)$ and  
$$A\in\Ga\big(\Si;\Hom_{\R}(V,(T^*\Si,\fj)^{0,1}\!\otimes_{\C}\!V) \big)$$ 
is a zeroth-order deformation term. 
The completion of a real CR-operator~$D$ on~$V$
with respect to appropriate norms on its domain and target, 
which we also denote by~$D$, is Fredholm;
see \cite[Theorem~C.1.10]{MS}.
The space of completions of all real CR-operators on~$V$ is contractible 
with respect to the operator~norm.
The real CR-operator $\bp$ on~$V$ in the notation above is $\C$-linear and its 
kernel and cokernel have canonical orientations;
we call the induced orientation of~$\la(D)$ for any real CR-operator~$D$ on~$V$
the \sf{complex orientation} of~$\la(D)$ and~$D$.\\

If $D$ is a real CR-operator on a real bundle pair $(V,\vph)$ over 
a decorated symmetric surface~$(\Si,\si)$ with $\Si^{\si}\!=\!\eset$, 
i.e.~$(\Si,\si)$ is a union of doublets, and $\Si^b$ is a half-surface of $(\Si,\si)$,
then the homomorphisms
\BE{CROdoubl_e0}\Ga(\Si;V)^{\vph}\lra \Ga\big(\Si^b;V|_{\Si^b}\big),  ~~\xi\lra\xi\big|_{\Si^b},\quad
\Ga_{\fj}^{0,1}(\Si;V)^{\vph}\lra \Ga_{\fj|_{\Si^b}}^{0,1}\big(\Si^b;V|_{\Si^b}\big),  
~~\ze\lra\ze\big|_{\Si^b},\EE
are isomorphisms.
They identify~$D$ with its ``restriction" to a real CR-operator $D^+$
on $V|_{\Si^b}$ and thus induce an isomorphism
\BE{CROdoubl_e} \la(D)\lra \la\big(D^+\big)\EE
between the determinants of real CR-operators.

\begingroup
\def\theCROrient{1$\os$}

\begin{CROrient}[Dependence on Decorated Structure: $\OSpin$]\label{CROos_prop}
Suppose $(\Si,\si)$ is a smooth decorated symmetric surface,
$(V,\vph)$ is a rank~$n$ real bundle pair over~$(\Si,\si)$, and
$\os$ is a relative $\OSpin$-structure on the real vector bundle~$V^{\vph}$
over $\Si^{\si}\!\subset\!\Si$.
\begin{enumerate}[label=(\arabic*),leftmargin=*]

\item\label{osflip_it} 
The orientation~$\fo_{\os}(V,\vph)$ of~$D_{(V,\vph)}$ determined by~$\os$
does not depend on the choice of half-surface $\Si_*^b$ of 
an elemental component~$\Si_*$ of~$(\Si,\si)$
if and only~if $\vp_{\os}(\Si_*)\!=\!0$.

\item\label{osinter_it} The interchange in the ordering of two components of~$\Si^{\si}$
preserves the orientation $\fo_{\os}(V,\vph)$ if and only if $n\!\in\!2\Z$.

\item\label{dltCROos_it} If $\Si^{\si}\!=\!\eset$ and
$\Si^b\!\subset\!\Si$ is the chosen half-surface, then the isomorphism~\eref{CROdoubl_e} is 
orientation-preserving with respect to the orientation $\fo_{\os}(V,\vph)$ on
the left-hand side and the canonical complex orientation on the right-hand side if and
only if \hbox{$\lr{w_2(\os),[\Si^b]_{\Z_2}}\!=\!0$}.

\end{enumerate}

\end{CROrient}
\endgroup

\vspace{.1in}

Suppose in addition that $\cC$ is a decorated marked symmetric surface as in~\eref{cCsymdfn_e}
and $\Si^b\!\subset\!\Si$ is the chosen half-surface.
For a component~$S^1_r$ of~$\Si^{\si}$, let
$$j_1^r(\cC),\ldots,j_{k_r}^r(\cC)\in[k]
\qquad\hbox{with}\quad  j_1^r(\cC)<j_2^r(\cC),\ldots,j_{k_r}^r(\cC)$$ 
be such that the real marked points  
\BE{posorder_e} x_{j_1^r(\cC)},\ldots,x_{j_{k_r}^r(\cC)}\in\Si_r^{\si}\EE
are ordered by their position on~$S_r^1$ with respect to the orientation of~$S^1_r$
induced by~$\Si^b$.
If $\cC'$ is a marked symmetric surface obtained from~$\cC$ by interchanging 
some real marked points, the twisted determinants $\wt\la_{\cC}(D)$ and $\wt\la_{\cC'}(D)$
are the same.
We~set 
$$\binom{-1}{1},\binom{-1}{2}\equiv 0.$$

\begingroup
\def\theCROrient{1$\fp$}

\begin{CROrient}[Dependence on Decorated Structure: $\Pin^{\pm}$]\label{CROp_prop}
Suppose $\cC$ is a smooth decorated marked symmetric surface,
$(V,\vph)$ is a $\cC$-balanced rank~$n$ real bundle pair over~$\cC$, and
$\fp$ is a relative $\Pin^{\pm}$-structure  on the real vector bundle~$V^{\vph}$
over $\Si^{\si}\!\subset\!\Si$.
\begin{enumerate}[label=(\arabic*),leftmargin=*]

\item\label{pflip_it} The orientation $\fo_{\cC;\fp}(V,\vph)$ 
on~$\wt\la_{\cC}(D_{(V,\vph)})$ determined by~$\fp$ 
does not depend on the choice of half-surface $\Si_*^b$ of 
an elemental component~$\Si_*$ of $(\Si,\si)$ if and only~if 
$$\sum_{S_r^1\in\pi_0(\Si_*^{\si})}\!\!\!
\bigg(\!\!\binom{k_r\!-\!1}{1}\!+\!n\binom{k_r\!-\!1}{2}\!\!\!\bigg)+2\Z =
\vp_{\fp}(\Si_*)\in \Z_2\,.$$

\item\label{pinter_it} The interchange in the ordering of two consecutive components~$S_r^1$ and~$S_{r+1}^1$ 
of~$\Si^{\si}$ preserves the orientation $\fo_{\cC;\fp}(V,\vph)$ 
if and only if $(n\!+\!1)(k_rk_{r+1}\!+\!1)\!\in\!2\Z$.

\item\label{pinter2_it} The interchange of two real marked points $x_{j_i^r(\cC)}$ and $x_{j_{i'}^r(\cC)}$
on the same connected component $S_r^1$ of~$\Si^{\si}$
with \hbox{$2\!\le\!i,i'\!\le\!k_r$} preserves $\fo_{\cC;\fp}(V,\vph)$.
The interchange of the real points $x_{j_1^r(\cC)}$ and $x_{j_i^r(\cC)}$
with $2\!\le\!i\!\le\!k_r$ preserves $\fo_{\cC;\fp}(V,\vph)$ if and only~if
$(n\!+\!1)(k_r\!-\!1)(i\!-\!1)\!\in\!2\Z$.

\end{enumerate}

\end{CROrient}
\endgroup

\addtocounter{CROrient}{-1}

\begin{CROrient}[Dependence on $\Spin$/$\Pin$-Structure]\label{CROSpinPinStr_prop}
Suppose $(\Si,\si)$ is a smooth decorated symmetric surface,
$(V,\vph)$ is a  real bundle pair over~$(\Si,\si)$, and $\eta\!\in\!H^2(\Si,\Si^{\si};\Z_2)$.
\begin{enumerate}[label=(\alph*),leftmargin=*]

\item\label{CROSpinStr_it} Let $\os\!\in\!\OSp_{\Si}(V^{\vph})$.
The orientations $\fo_{\os}(V,\vph)$ and $\fo_{\eta\cdot\os}(V,\vph)$ 
are the same if and only~if $\lr{\eta,[\Si^b]_{\Z_2}}\!=\!0$.
The orientations $\fo_{\os}(V,\vph)$ and $\fo_{\ov\os}(V,\vph)$
 are the same if and only if $|\pi_0(\Si^{\si})|$ is even.

\item\label{CROPinStr_it} Let $\cC$ be a marked symmetric surface with 
the underlying symmetric surface~$(\Si,\si)$ so that $(V,\vph)$ is $\cC$-balanced.
Let $\fp\!\in\!\cP_{\Si}^{\pm}(V^{\vph})$.
The orientations $\fo_{\cC;\fp}(V,\vph)$ and $\fo_{\cC;\eta\cdot\fp}(V,\vph)$ 
on~$\wt\la_{\cC}(D_{(V,\vph)})$  are the same if and only~if $\lr{\eta,[\Si^b]_{\Z_2}}\!=\!0$.

\end{enumerate}
\end{CROrient}

\vspace{.15in}

Let $(V,\vph)$ be a real bundle pair over a decorated marked symmetric surface~$\cC$
as in~\eref{cCsymdfn_e} and
$S^1_1,\ldots,S^1_N$ be the connected components of~$\Si^{\si}$ in
the chosen order.
For each component $S^1_r$ of~$\Si^{\si}$, define  
\BE{VvphcCrdfn_e}V^{\vph}_{\cC;r}=
V^{\vph}_{x_{j_1^r(\cC)}}\!\oplus\!\ldots\!\oplus\!V^{\vph}_{x_{j_{k_r}^r(\cC)}}\,.\EE
If $\fo$ is an orientation of~$V^{\vph}$,
we define the orientation $\la_{\cC}^{\R}(\fo)$ on $\la_{\cC}^{\R}(V,\vph)$ to be
the orientation induced by
the restriction of~$\fo$ to each~$V^{\vph}_{x_i}$ via the identification
$$\la_{\cC}^{\R}(V,\vph)=
\la\big(V^{\vph}_{x_{j_1^1(\cC)}}\big)\!\otimes\!\ldots\!\otimes\!
\la\big(V^{\vph}_{x_{j_{k_1}^1(\cC)}}\big)\!\otimes\!\ldots\!\otimes\!
\la\big(V^{\vph}_{x_{j_1^N(\cC)}}\big)\!\otimes\!\ldots\!\otimes\!
\la\big(V^{\vph}_{x_{j_{k_N}^N(\cC)}}\big)\,.$$
In other words, for each component $S^1_r$ of~$\Si^{\si}$ we first orient 
the direct sum~\eref{VvphcCrdfn_e} from the orientation~$\fo$ 
based on the listed order of the factors, 
i.e.~the order in which the marked points are positioned on~$S^1_r$.
We then orient the direct sum 
\BE{VRcCdfn_e}V^{\vph}_{\cC}\equiv V^{\vph}_{\cC;1}\!\oplus\!\ldots\!\oplus\!V^{\vph}_{\cC;N}\EE
over all components of~$\Si^{\si}$ based on the ordering of the components.

\begin{CROrient}[Reduction]\label{CROPin2SpinRed_prop}
Suppose $\cC$ is a decorated marked symmetric surface,
$(V,\vph)$ is a $\cC$-balanced rank~$n$ real bundle pair over~$\cC$, and
$\fo$ is an orientation on~$V^{\vph}$.
For every $\fp\!\in\!\cP_{\Si}^{\pm}(V^{\vph})$, 
the orientation $\fo_{\cC;\fp}(V,\vph)$ on~$\wt\la_{\cC}(D_{(V,\vph)})$ 
corresponds to the homotopy class of isomorphisms of 
$\la(D_{(V,\vph)})$ and $\la_{\cC}^{\R}(V,\vph)$ determined by the orientations
$\fo_{\fR_{\fo}^{\pm}(\fp)}(V,\vph)$  and  $\la_{\cC}^{\R}(\fo)$ if and only~if 
\BE{CROPin2SpinRed_e}\fp\in\cP_{\Si}^-(V^{\vph}),~n\binom{|\pi_0(\Si^{\si})|}{2}\in2\Z
\quad\hbox{or}\quad
\fp\in\cP_{\Si}^+(V^{\vph}),~(n\!+\!1)\!\binom{|\pi_0(\Si^{\si})|}{2}\in2\Z.\EE
\end{CROrient}

\vspace{.1in}

Suppose $(\Si_1,\si_1)$ and $(\Si_2,\si_2)$ are symmetric surfaces,
$(V_1,\vph_1)$ is a real bundle pair on $(\Si_1,\si_1)$, and 
$(V_2,\vph_2)$ is a real bundle pair on $(\Si_2,\si_2)$ of the same rank as~$(V_1,\vph_1)$.
We denote by
\BE{V1V2dfn_e}
(V_1,\vph_1)\!\sqcup\!(V_2,\vph_2)\equiv \big(V_1\!\sqcup\!V_2,\vph_1\!\sqcup\!\vph_2\big)\EE
the induced real bundle pair over the symmetric surface
\BE{Si1Si2dfn_e}(\Si_1,\si_1)\!\sqcup\!(\Si_2,\si_2)\equiv 
\big(\Si_1\!\sqcup\!\Si_2,\si_1\!\sqcup\!\si_2\big).\EE
For relative $\OSpin$-structures (resp.~$\Pin^{\pm}$-structures) 
$\os_1$ (resp.~$\fp_1$) on $V_1^{\vph_1}$ and $\os_2$ (resp.~$\fp_2$) on~$V_2^{\vph_2}$,
we denote by $\os_1\!\sqcup\!\os_2$ (resp.~$\fp_1\!\sqcup\!\fp_2$)
the relative $\OSpin$-structure (resp.~$\Pin^{\pm}$-structure) on 
the real vector bundle
$$(V_1\!\sqcup\!V_2)^{\vph_1\sqcup\vph_2}=V_1^{\vph_1}\!\sqcup\!V_2^{\vph_2}$$
over
\BE{Si1Si2dfn_e2}
 \big(\Si_1\!\sqcup\!\Si_2\big)^{\si_1\sqcup\si_2}
= \Si_1^{\si_1}\!\sqcup\!\Si_2^{\si_2} \subset\Si_1\!\sqcup\!\Si_2.\EE
If $(\Si_1,\si_1)$ and $(\Si_2,\si_2)$ are decorated symmetric surfaces,
we take~\eref{Si1Si2dfn_e} to be the decorated symmetric surface obtained by~defining
$$\big(\Si_1\!\sqcup\!\Si_2\big)^b=\Si_1^b\!\sqcup\!\Si_2^b$$
and ordering the topological components of~\eref{Si1Si2dfn_e2} so that 
the topological components of~$\Si_1^{\si_1}$ in their order are  
followed by the topological components of~$\Si_2^{\si_2}$ in their order.\\

If in addition $\cC_1$ and $\cC_2$ are marked symmetric surfaces with 
the underlying symmetric surfaces~$(\Si_1,\si_1)$ and~$(\Si_2,\si_2)$,
respectively, we denote by $\cC_1\!\sqcup\!\cC_2$ the  marked symmetric surface
with the underlying symmetric surface~\eref{Si1Si2dfn_e} and the same marked
points as~$\cC_1$ and~$\cC_2$ which are ordered consistently with their orderings
in~$\cC_1$ and~$\cC_2$ (but the marked points of~$\cC_1$ do not necessarily precede
the marked points of~$\cC_2$).
Let
\BE{V1V2dfn_e4} 
\Psi_{\cC_1,\cC_2}^{\R}\!:\la_{\cC_1\sqcup\cC_2}^{\R}
\big(\!(V_1,\vph_1)\!\sqcup\!(V_2,\vph_2)\!\big)
\approx
\la_{\cC_1}^{\R}(V_1,\vph_1)\!\otimes\!\la_{\cC_2}^{\R}(V_2,\vph_2)\EE
be the isomorphism induced by the identification 
$$(V_1\!\sqcup\!V_2)^{\vph_1\sqcup\vph_2}_{\cC_1\sqcup\cC_2}
=(V_1)^{\vph_1}_{\cC_1}\!\oplus\!(V_2)^{\vph_2}_{\cC_2}\,.$$
If $(V_1,\vph_1)$ is $\cC_1$-balanced and $(V_2,\vph_2)$ is $\cC_2$-balanced,
then~\eref{V1V2dfn_e} is $(\cC_1\!\sqcup\!\cC_2)$-balanced.\\

For every real CR-operator~$D$ on the real bundle pair~\eref{V1V2dfn_e} over~\eref{Si1Si2dfn_e},
there is a canonical splitting 
\BE{D1D2split_e} D=D\big|_{V_1}\!\oplus\!D\big|_{V_2}\EE
of $D$ into real CR-operators on $(V_1,\vph_1)$ and $(V_2,\vph_2)$.
Along with~\eref{sum}, this splitting determines a homotopy class of isomorphisms
\BE{CRODisjUn_e0} \la\big(D_{(V_1,\vph_1)\sqcup(V_2,\vph_2)}\big)\approx 
\la\big(D_{(V_1,\vph_1)}\big)\!\otimes\!\la\big(D_{(V_2,\vph_2)}\big)\,.\EE
Orientations~$\wt\fo_1$ of $\wt\la_{\cC_1}\!(D_{(V_1,\vph_1)})$ and
$\wt\fo_2$ of $\wt\la_{\cC_2}\!(D_{(V_2,\vph_2)})$ correspond to 
homotopy classes of isomorphisms 
$$\la\big(D_{(V_1,\vph_1)}\big)\approx \la_{\cC_1}^{\R}(V_1,\vph_1)
\qquad\hbox{and}\qquad
\la\big(D_{(V_2,\vph_2)}\big)\approx \la_{\cC_2}^{\R}(V_2,\vph_2),$$
respectively. Along with~\eref{CRODisjUn_e0} and~\eref{V1V2dfn_e4}, 
these homotopy classes of isomorphisms determine a homotopy class of isomorphisms
$$\la\big(D_{(V_1,\vph_1)\sqcup(V_2,\vph_2)}\big)\approx
\la_{\cC_1\sqcup\cC_2}^{\R}\big(\!(V_1,\vph_1)\!\sqcup\!(V_2,\vph_2)\!\big)\,.$$
The latter corresponds to an orientation of 
$\wt\la_{\cC_1\sqcup\cC_2}\!(D_{(V_1,\vph_1)\sqcup(V_2,\vph_2)})$,
which we denote by $\wt\fo_1\!\sqcup\!\wt\fo_2$. 

\begin{CROrient}[Disjoint unions]\label{CRODisjUn_prop}
Suppose $(\Si_1,\si_1)$ and $(\Si_2,\si_2)$ are smooth decorated symmetric surfaces
of genera~$g_1$ and~$g_2$, respectively,
$(V_1,\vph_1)$ is a rank~$n$ real bundle pair on $(\Si_1,\si_1)$, and 
$(V_2,\vph_2)$ is a rank~$n$ real bundle pair on $(\Si_2,\si_2)$.
\begin{enumerate}[label=(\alph*),leftmargin=*]

\item\label{osDisjUn_it} Let $\os_1\!\in\!\OSp_{\Si_1}(V_1^{\vph_1})$ and 
$\os_2\!\in\!\OSp_{\Si_2}(V_2^{\vph_2})$.
The homotopy class of isomorphisms~\eref{CRODisjUn_e0}
induced by the splitting~\eref{D1D2split_e}  respects the orientations 
$$\fo_{\os_1\sqcup\os_2}\!\big(\!(V_1,\vph_1)\!\sqcup\!(V_2,\vph_2)\!\big),\quad
\fo_{\os_1}\!(V_1,\vph_1), \quad\hbox{and}\quad \fo_{\os_2}\!(V_2,\vph_2)$$ 
on the three factors. 

\item\label{pDisjUn_it} Let $\cC_1$ and $\cC_2$ be marked symmetric surfaces with 
the underlying symmetric surfaces~$(\Si_1,\si_1)$ and~$(\Si_2,\si_2)$,
respectively,  so that $(V_1,\vph_1)$ is $\cC_1$-balanced and 
$(V_2,\vph_2)$ is $\cC_2$-balanced.
If \hbox{$\fp_1\!\in\!\cP_{\Si_1}^{\pm}(V^{\vph_1})$} and $\fp_2\!\in\!\cP_{\Si_2}^{\pm}(V^{\vph_2})$,
then
$$\fo_{\cC_1\sqcup\cC_2;\fp_1\sqcup\fp_2}\!\big(\!(V_1,\vph_1)\!\sqcup\!(V_2,\vph_2)\!\big)
=\fo_{\cC_1;\fp_1}\!(V_1,\vph_1)\!\sqcup\!\fo_{\cC_2;\fp_2}\!(V_2,\vph_2)$$
if and only~if 
\BE{CRODisjUn_e0c}\begin{split}
&\fp\in\cP_{\Si}^-(V^{\vph}),~
\big(1\!-\!g_1\!+\!\deg V_1\big)\big(\!(1\!-\!g_2)n\!+\!\deg V_2\big)\in2\Z
\qquad\hbox{or}\\
&\fp\in\cP_{\Si}^+(V^{\vph}),~
(n\!+\!1)\big(1\!-\!g_1\!+\!\deg V_1\big)\big(1\!-\!g_2\big)\in2\Z.
\end{split}\EE

\end{enumerate}
\end{CROrient}

\vspace{.15in}

A short exact sequence 
\BE{RBPsesdfn_e} 0\lra (V',\vph')\lra (V,\vph)\lra (V'',\vph'')\lra0\EE
of real bundle pairs over $(X,\phi)$ induces a short exact sequence
\BE{RBPsesdfn_e2} 0\lra V'^{\vph'}\lra V^{\vph}\lra V''^{\vph''}\lra0\EE
of real vector bundles over~$X^{\phi}$.
If the first sequence is denoted by~$\ce$,  we denote the second one by~$\ce_{\R}$.
By the SpinPinRel~\ref{RelSpinPinSES_prop} property on page~\pageref{RelSpinPinSES_prop}, 
a relative $\OSpin$-structure $\os'\!\in\!\OSp_X(V'^{\vph'})$
on the vector bundle~$V'^{\vph'}$ over $X^{\phi}\!\subset\!X$ and
a relative $\OSpin$-structure $\os''\!\in\!\OSp_X(V''^{\vph''})$
(resp.~$\Pin^{\pm}$-structure $\fp''\!\in\!\cP_X^{\pm}(V''^{\vph''})$)
determine a relative $\OSpin$-structure
$$ \llrr{\os',\os''}_{\ce_{\R}}\in\OSp_X(V^{\vph})$$
(resp.~$\Pin^{\pm}$-structure $\llrr{\os',\fp''}_{\ce_{\R}}\!\in\!\cP_X^{\pm}(V^{\vph})$).\\

Suppose in addition that $(X,\phi)\!=\!(\Si,\si)$ is a smooth symmetric surface  
and $\cC$ is a marked symmetric surface with the underlying symmetric surface~$(\Si,\si)$
as in~\eref{cCsymdfn_e}.
For each $i\!\in\![k]$, the short exact sequence~\eref{RBPsesdfn_e2} determines
an isomorphism
$$\la\big(V^{\vph}_{x_i}\big)\approx 
\la\big(V'^{\vph'}_{x_i}\big)\!\otimes\!\la\big(V''^{\vph''}_{x_i}\big).$$
Putting these isomorphisms together, we obtain an identification
\BE{RBPsesdfn_e4} \Psi_{\cC}^{\R}\!:\la_{\cC}^{\R}(V,\vph)\approx
\la_{\cC}^{\R}(V',\vph')\!\otimes\!\la_{\cC}^{\R}(V'',\vph'')\,.\EE
If the real vector bundle~$V'^{\vph'}$ is orientable, then
the real bundle pair $(V,\vph)$ is $\cC$-balanced if and only if
the real bundle pair $(V'',\vph'')$ is $\cC$-balanced.
If $\fo'$ is an orientation on~$V'^{\vph'}$ and $\os'\!\in\!\Spin_{\Si}(V'^{\vph'},\fo')$,
let 
$$\la_{\cC}^{\R}(\os')\equiv \la_{\cC}^{\R}(\fo')$$
be the orientation on $\la_{\cC}^{\R}(V',\vph')$ defined as above
the CROrient~\ref{CROPin2SpinRed_prop} property.\\

Along with~\eref{sum}, the short exact sequence~$\ce$ of real bundle pairs over~$(\Si,\si)$ 
as in~\eref{RBPsesdfn_e} determines a homotopy class of isomorphisms
\BE{CROses_e0} \la\big(D_{(V,\vph)}\big)\approx 
\la\big(D_{(V',\vph')}\big)\!\otimes\!\la\big(D_{(V'',\vph'')}\big)\,.\EE
Orientations~$\fo_D'$ of $\la(D_{(V',\vph')})$, 
$\la_{\cC}^{\R}(\os')$ of $\la_{\cC}^{\R}(V',\vph')$,
and $\wt\fo''$ of $\wt\la_{\cC}(D_{(V'',\vph'')})$ determine 
homotopy classes of isomorphisms 
$$\la\big(D_{(V',\vph')}\big)\approx \la_{\cC}^{\R}(V',\vph')
\qquad\hbox{and}\qquad
\la\big(D_{(V'',\vph'')}\big)\approx \la_{\cC}^{\R}(V'',\vph''),$$
respectively. 
Along with~\eref{CROses_e0} and~\eref{RBPsesdfn_e4}, 
these homotopy classes of isomorphisms determine a homotopy class of isomorphisms
$$\la\big(D_{(V,\vph)}\big)\approx \la_{\cC}^{\R}(V,\vph)\,.$$
The latter corresponds to an orientation of $\wt\la_{\cC}(D_{(V,\vph)})$,
which we denote by $(\fo'_D\la_{\cC}^{\R}(\os')\!)_{\ce}\wt\fo''$.

\begin{CROrient}[Exact Triples]\label{CROSpinPinSES_prop}
Suppose $(\Si,\si)$ is a smooth decorated symmetric surface,
$\ce$ is a short exact sequence of real bundle pairs over~$(\Si,\si)$ 
as in~\eref{RBPsesdfn_e}, and $\os'\!\in\!\OSp_{\Si}(V'^{\vph'})$.
\begin{enumerate}[label=(\alph*),leftmargin=*]

\item\label{CROsesSpin_it} Let $\os''\!\in\!\OSp_{\Si}(V''^{\vph''})$.
The homotopy class of isomorphisms~\eref{CROses_e0} induced by~$\ce$ respects the orientations 
$$\fo_{\llrr{\os',\os''}_{\ce_{\R}}}\!(V,\vph),\quad \fo_{\os'}(V',\vph'), 
\quad\hbox{and}\quad \fo_{\os''}(V'',\vph'')$$ 
if and only~if 
\BE{CROses_e0b}\big(\rk\,V'\big)\big(\rk\,V''\big)\binom{|\pi_0(\Si^{\si})|}{2}\in2\Z.\EE

\item\label{CROsesPin_it} Let $\cC$ be a marked symmetric surface with 
the underlying symmetric surface~$(\Si,\si)$ so that $(V'',\vph'')$ is $\cC$-balanced.
If $\fp''\!\in\!\cP_{\Si}^{\pm}(V''^{\vph''})$, then 
$$\fo_{\cC;\llrr{\os',\fp''}_{\ce_{\R}}}\!(V,\vph)=
\big(\fo_{\cC;\os'}(V',\vph')\la_{\cC}^{\R}(\os')\!\big)_{\ce}\fo_{\cC;\fp''}(V'',\vph'')$$
if and only~if 
\BE{CROsesPin_e0}\big(\rk\,V'\big)\!\big(\rk\,V''\!+\!1\big)\!\binom{|\pi_0(\Si^{\si})|}{2}
+\big(\rk\,V'\big)\!\big(\rk\,V''\big)\binom{k}{2}\in2\Z.\EE

\end{enumerate}
\end{CROrient}

\vspace{.15in}

There are two equivalence classes of orientation-reversing involutions on
the Riemann sphere \hbox{$S^2\!=\!\P^1$}.
They are represented by the involutions 
\BE{tauetadfn_e}\tau\!:\P^1\lra\P^1, \quad z\lra 1/\bar{z}, \quad\hbox{and}\quad
\eta\!:\P^1\lra\P^1, \quad z\lra -1/\bar{z}.\EE
The fixed-point locus of~$\eta$ is empty.
The fixed-point locus of~$\tau$ is the unit circle $S^1\!\subset\!\C$;
it separates~$\P^1$ into the unit disk~$\bD^2_+\!\subset\!\C$ and its complement~$\bD^2_-$.
A choice of half-surface~$\bD^2_+$ of~$(S^2,\tau)$ is equivalent to a choice
of one of the two natural embeddings of~$\C$ as a subspace of~$S^2$.\\

Let $(V,\vph)$ be a rank~1 degree~$a$ real bundle pair over~$(S^2,\tau)$.
If $a\!\ge\!-1$, $x_1,\ldots,x_{a+1}\!\in\!S^1$ are distinct points,
and $D_{(V,\vph)}$ is a real CR-operator on~$(V,\vph)$, then
$D_{(V,\vph)}$ is surjective and the evaluation homomorphism
\BE{cOevdfn_e0} \ev\!: \ker D_{(V,\vph)}\lra  
V_{x_1}^{\vph}\!\oplus\!\ldots\!\oplus\!V_{x_{a+1}}^{\vph}, \qquad
\xi\lra\big(\xi(x_1),\ldots,\xi(x_{a+1})\big),\EE
is an isomorphism; see Section~\ref{OrientEg_subs}.
If $a\!\in\!2\Z$, the real line bundle~$V^{\vph}$ over~$S^1$ is orientable.
For an orientation~$\fo$ on~$V^{\vph}$, we then let
\BE{osLBdfn_e}\fo_0(V,\vph;\fo)\equiv \fo_{\io_{S^2}(\os_0(V^{\vph},\fo))}(V,\vph)\EE
denote the orientation of~$D_{(V,\vph)}$ determined by the image
$$\io_{S^2}\big(\os_0(V^{\vph},\fo)\!\big)\in\OSp_{S^2}\big(V^{\vph}\big)$$
of  the canonical $\OSpin$-structure $\os_0(V^{\vph},\fo)$ on $(V^{\vph},\fo)$
under the first map in~\eref{vsSpinPin_e0} with $X\!=\!S^2$.\\

If $a\!\not\in\!2\Z$, the real line bundle~$V^{\vph}$ over~$S^1$ is not orientable.
The two $\Pin^{\pm}$-structures on~$V^{\vph}$, $\fp_0^{\pm}(V^{\vph})$ and $\fp_1^{\pm}(V^{\vph})$, 
are then distinguished by Examples~\ref{Pin1pMB_eg} and~\ref{Pin1mMB_eg}
from the classical perspective and by Example~\ref{SpinDfn1to3_eg}
from the trivializations perspectives.
The distinction between $\fp_0^-(V^{\vph})$ and $\fp_1^-(V^{\vph})$ depends
on the choice of half-surface~$\bD^2_+$ of~$(S^2,\tau)$,
while the distinction between $\fp_0^+(V^{\vph})$ and $\fp_1^+(V^{\vph})$ is independent 
of such a choice.
If $\cC$ is a marked symmetric surface as in~\eref{cCsymdfn_e} so that 
$(\Si,\si)$ is $S^2$ with the involution~$\tau$, $k\!\in\!2\Z$, and $r\!=\!0,1$,
let
$$\fo_{\cC;r}^{\pm}(V,\vph)\equiv \fo_{\cC;\io_{S^2}(\fp_r^{\pm}(V^{\vph}))}(V,\vph)$$
denote the orientation on $\wt\la_{\cC}(D_{(V,\vph)})$ determined by the image
\BE{PinNormDfn_e}\io_{S^2}\big(\fp_r^{\pm}(V^{\vph})\!\big)\in\cP_{S^2}^{\pm}\big(V^{\vph}\big)\EE
of $\fp_r^{\pm}(V^{\vph})$
under the second map in~\eref{vsSpinPin_e0} with $X\!=\!S^2$.

\begin{CROrient}[Normalizations]\label{CRONormal_prop}
Let $(V,\vph)$ be a rank~1 real bundle pair over~$(S^2,\tau)$.
\begin{enumerate}[label=(\alph*),leftmargin=*]

\item\label{CROnormSpin_it} If $\deg V\!=\!0$ and $\fo\!\in\!\fO(V^{\vph})$, 
the isomorphism~\eref{cOevdfn_e0} with $a\!=\!0$ respects the orientations $\fo_0(V,\vph;\fo)$ and~$\fo$
for every $x_1\!\in\!S^1$.

\item\label{CROnormPin_it} If $\deg V\!=\!1$,
$\cC$ is a marked symmetric surface with the underlying symmetric surface~$(S^2,\tau)$ 
and $k\!=\!2$ real marked points, and
$\bD^2_+$ is the distinguished half-surface of~$(S^2,\tau)$, then
$\fo_{\cC;0}^{\pm}(V,\vph)$ is the orientation of $\wt\la_{\cC}(D_{(V,\vph)})$
corresponding to the homotopy class of the isomorphism~\eref{cOevdfn_e0} with $a\!=\!1$.
\end{enumerate}

\end{CROrient}

\vspace{.15in}

In the setting of the CROrient~\ref{CRONormal_prop}\ref{CROnormSpin_it} property,  
the $\OSpin$-structure $\os_0(V^{\vph},\fo)$ on the oriented real line bundle~$(V^{\vph},\fo)$
over~$S^1$ is the homotopy class of trivializations of 
$2\tau_{S^1}\!\oplus\!V^{\vph}$ induced by 
the unique homotopy classes of trivializations of~$\tau_{S^1}$ and~$V^{\vph}$ as oriented line bundles.
It is thus not surprising that a systematic orienting procedure for CR-operators 
on even-degree real bundle pairs over~$(S^2,\tau)$ with relative $\OSpin$-structures,
such as the one described at the beginning of Section~\ref{tauspinorient_subs},   
satisfies the  CROrient~\ref{CRONormal_prop}\ref{CROnormSpin_it} property.
This is immediate from the first case of~\eref{SpinP1rk1_e} and
is noted in Proposition~\ref{tauspinorient_prp}\ref{SpinvsCan_it}.\\

In the setting of the CROrient~\ref{CRONormal_prop}\ref{CROnormPin_it} property, 
the $\Pin^-$- and $\Pin^+$-structures $\fp_0^-(V^{\vph})$ and $\fp_0^+(V^{\vph})$
on the unorientable real line bundle~$V^{\vph}$ over~$S^1$
correspond to the homotopy classes $\os_0(2\ga_{\R;1},\fo_{\ga_{\R;1}}^-)$
and $\os_0(4\ga_{\R;1},\fo_{\ga_{\R;1}}^+)$ of trivializations of the canonically oriented vector
bundles
$$\tau_{\R\P^1}\!\oplus\!2\ga_{\R;1}\lra\R\P^1  \qquad\hbox{and}\qquad 
4\ga_{\R;1}\lra\R\P^1,$$
respectively.
The latter is the canonical homotopy class of trivializations provided 
by Example~\ref{CanSpin_eg}\ref{CanSpin_it2}.
The former is the homotopy class of trivializations determined by 
the unique homotopy class of trivializations of~$\tau_{\R\P^1}$ as an oriented line bundle
and the homotopy class of trivializations of~$2\ga_{\R;1}$ represented by the trivialization~$\Phi_0$
in~\eref{ga1Rpin_e}.
According to \cite[Proposition~3.5]{RealGWsII}, the orientation on the~line
$$\la(D_{2(V,\vph)})=\la(2D_{(V,\vph)})=\la(D_{(V,\vph)})^{\otimes2}$$
induced by the OSpin-structure $\os_0(2\ga_{\R;1},\fo_{\ga_{\R;1}}^-)$ 
is the canonical orientation of this line as a square of another line; 
see also Proposition~\ref{tauspinorient_prp}\ref{PinCan_it}.
Along with the construction of twisted orientations of Theorem~\ref{CROrient_thm}\ref{CROrientPin_it}
at the beginning of Section~\ref{twisteddetLB_subs},
this immediately implies both cases of
the  CROrient~\ref{CRONormal_prop}\ref{CROnormPin_it} property.\\

By the CROrient~\ref{CRONormal_prop} property, the orientations~$\fo_0(V,\vph;\fo)$
and $\fo_{\cC;0}^{\pm}(V,\vph)$ match intrinsic evaluation orientations 
for rank~1 real bundle pairs~$(V,\vph)$ over $(S^1,\tau)$ of degrees~0 and~1.
Corollary~\ref{NormEv_crl} extends this property to 
rank~1 real bundle pairs over $(S^1,\tau)$ of arbitrary degrees.

\subsection{Properties of orientations: degenerations}
\label{OrientPrp_subs2}

The compatibility of the orientations on the determinants of real CR-operators 
induced by the relative $\Spin$- and $\Pin^{\pm}$-structures with degenerations
of smooth domains is essential for
any study of the structures in open and real Gromov-Witten theories 
that depend on node-splitting properties in the spirit of \cite[2.2.6]{KM}.
We next describe the relevant degeneration settings and characterize 
the behavior of the induced orientations in these settings by
the CROrient~\ref{CROdegenC_prop} and~\ref{CROdegenH3_prop} properties.\\

Let $\cC\!\equiv\!(\Si,(z_i)_{i\in[l]})$ be a nodal marked surface,
i.e.~$\Si$ is a closed oriented (possibly nodal) surface and
$z_1,\ldots,z_l$ are distinct smooth points of~$\Si$.
A \sf{flat family of deformations} of~$\cC$ is a~tuple
\BE{CcCdfn_e}\big(\pi\!:\cU\!\lra\!\De,(s_i\!:\De\!\lra\!\cU)_{i\in[l]}\big),\EE
where $\cU$ is a smooth manifold, 
$\De\!\subset\!\C^N$ is a ball around~$0$, 
and $\pi,s_1,\ldots,s_l$ are smooth maps such~that
\begin{enumerate}[label=$\bullet$,leftmargin=*]

\item $\Si_{\t}\!\equiv\!\pi^{-1}(\t)$ 
is a closed oriented (possibly nodal) surface for each $\t\!\in\!\De$
and $\pi$ is a submersion outside of the nodes of the fibers of~$\pi$,

\item for every $\t^*\!\equiv\!(t_j^*)_{j\in[N]}\!\in\!\De$ and 
every node $z^*\!\in\!\Si_{\t^*}$, 
there exist $i\!\in\!\{1,\ldots,N\}$ with $t_i^*\!=\!0$,
neighborhoods $\De_{\t^*}$ of~$\t^*$ in~$\De$ and $\cU_{z^*}$ of~$z^*$ in~$\cU$, 
and a diffeomorphism
$$\Psi\!:\cU_{z^*}\lra \big\{\big((t_j)_{j\in[N]},x,y\big)\!\in\!\De_{\t^*}\!\times\!\C^2\!:
xy\!=\!t_i\big\}$$
such that the composition of~$\Psi$ with the projection to~$\De_{\t^*}$ equals~$\pi|_{\cU_{z^*}}$,

\item $\pi\!\circ\!s_i\!=\!\id_{\De}$ and $s_1(\t),\ldots\!s_l(\t)\!\in\!\Si_{\t}$
are distinct smooth points for every $\t\!\in\!\De$,

\item $(\Si_0,(s_i(0))_{i\in[l]})\!=\!\cC$.

\end{enumerate} 
We call an almost complex structure~$\fj_{\cU}$ on~$\cU$ \sf{admissible} if
the projection~$\pi$ is $(\fj_{\C^N},\fj_{\cU})$-holomorphic and 
for every node $z^*\!\in\!\Si_{\t^*}$ the diffeomorphism~$\Psi$ above can be chosen
to be holomorphic.\\

Let $\cC$ be a nodal marked symmetric surface as in~\eref{cCsymdfn_e}.
A \sf{flat family of deformations} of~$\cC$ is a tuple 
\BE{cUsymmdfn_e} \fF\!\equiv\!\big(\pi\!:\cU\!\lra\!\De,\wt\fc\!:\cU\!\lra\!\cU,
(s_i^{\R}\!:\De\!\lra\!\cU)_{i\in[k]},
(s_i^+\!:\De\!\lra\!\cU,s_i^-\!:\De\!\lra\!\cU)_{i\in[l]}\big)\EE
such that
$(\pi,(s_i^{\R})_{i\in[k]},(s_i^+,s_i^-)_{i\in[l]})$
is a flat family of deformations of~$\cC$,
$\wt\fc$ is an involution on~$\cU$ lifting 
the standard involution~$\fc$ on~$\De$ and restricting to~$\si$ over
\hbox{$\Si\!\equiv\!\Si_0$}, and
$$s_i^{\R}\!\circ\!\fc=\wt\fc\!\circ\!s_i^{\R}~~
\forall\,i\!\in\![k], \qquad
s_i^+\!\circ\!\fc=\wt\fc\!\circ\!s_i^-~~\forall\,i\!\in\![l].$$
We call an almost complex structure~$\fj_{\cU}$ on~$(\cU,\wt\fc)$ \sf{admissible}
if $\fj_{\cU}$ is admissible for~$\cU$ and \hbox{$\wt\fc^*\fj_{\cU}\!=\!-\fj_{\cU}$}.\\

If $\fF$ is as in~\eref{cUsymmdfn_e}, let $\si_{\t}\!=\!\wt\fc|_{\Si_{\t}}$ for each
parameter~$\t$ in $\De_{\R}\!\equiv\!\De\!\cap\!\R^N$.
The tuple
$$\cC_{\t}\equiv\big(\Si_{\t},\si_{\t},(s_i^{\R}(\t))_{i\in[k]},
(s_i^+(\t),s_i^-(\t))_{i\in[l]}\big)$$
is then a nodal marked symmetric surface.
If $\fj_{\cU}$ is an admissible almost complex structure on~$(\cU,\wt\fc)$,
then $\fj_{\t}\!\equiv\!\fj_{\cU}|_{\Si_{\t}}$ is a complex structure on~$\Si_{\t}$.
We denote by $\De_{\R}^*\!\subset\!\De_{\R}$ the subspace of elements $\t\!\in\!\De_{\R}$
so that the surface~$\Si_{\t}$ is smooth.
If $(\Si_0,\si_0)$ has only conjugate pairs of nodes, $\De_{\R}\!-\!\De_{\R}^*$ is
a subspace of~$\De_{\R}$ of codimension~2;
if $(\Si_0,\si_0)$ has real nodes, $\De_{\R}\!-\!\De_{\R}^*$ is
a subspace of codimension~1.\\

\noindent
Let $(\pi\!:\cU\!\lra\!\De,\wt\fc\!:\cU\!\lra\!\cU)$ be a flat family of deformations
of a (possibly) nodal symmetric surface and 
$(V,\vph)$ be a real bundle pair over~$(\cU,\wt\fc)$.
For each $\t\!\in\!\De_{\R}$,
$$(V_{\t},\vph_{\t}\big)\equiv \big(V,\vph)|_{\Si_{\t}}$$
is then a real bundle pair over $(\Si_{\t},\si_{\t})$.
Suppose in addition that $\fj_{\cU}$ is  an admissible almost complex structure on~$(\cU,\wt\fc)$,
$\na$ is a $\vph$-compatible (complex-linear) connection in~$V$, and 
\BE{AcUdfn_e} A\in\Ga\big(\cU;\Hom_{\R}(V,(T^*\cU,\fj_{\cU})^{0,1}\otimes_{\C}\!V)\big)^{\vph}.\EE
The restrictions of~$\na$ and~$A$ to each fiber~$(\Si_{\t},\si_{\t})$
of~$\pi$ with $\t\!\in\!\De_{\R}$ then determine a real CR-operator
\BE{DVvphdfn_e}D_{(V,\vph);\t}\!: \Ga\big(\Si_{\t};V_{\t}\big)^{\vph_{\t}}\lra
\Ga_{\fj_{\t}}^{0,1}\big(\Si_{\t};V_{\t}\big)^{\vph_{\t}}\EE
on $(V_{\t},\vph_{\t})$.
Let
\BE{detLB_e}\pi_{(V,\vph)}\!:\la\big(\cD_{(V,\vph)}\!\big)\!\equiv\!
\bigsqcup_{\t\in\De_{\R}}\!\!
\{\t\}\!\times\!\la\big(D_{(V,\vph);\t}\!\big)\lra\De_{\R}\,.\EE
The set $\la(\cD_{(V,\vph)})$ carries natural topologies so that the projection
$\pi_{(V,\vph)}$ is a real line bundle; see \cite[Appendix~A]{RealGWsI}.
These topologies in particular satisfy the following properties:
\begin{enumerate}[label=(D\arabic*),leftmargin=*]

\item\label{dNatI_it} a homotopy class of continuous isomorphisms 
$\Psi\!:(V_1,\vph_1)\!\lra\!(V_2,\vph_2)$ of real bundle pairs over $(\cU,\wt\fc)|_{\De_{\R}}$
determines a homotopy class of isomorphisms
$$\la(\cD_{\Psi})\!:\la\big(\cD_{(V_1,\vph_1)}\!\big)\lra \la\big(\cD_{(V_2,\vph_2)}\!\big)$$
of line bundles over $\De_{\R}$;

\item\label{dses_it} the isomorphisms \eref{sum} determine a homotopy class of isomorphisms
$$\la\!\big(\cD_{(V_1,\vph_1)\oplus(V_2,\vph_2)}\!\big)\approx
\la\big(\cD_{(V_1,\vph_1)}\!\big)\otimes \la\big(\cD_{(V_2,\vph_2)}\!\big)$$
of line bundles over $\De_{\R}$ for all real bundle pairs 
$(V_1,\vph_1)$ and $(V_2,\vph_2)$ over~$(\cU,\wt\fc)$.

\end{enumerate}

\vspace{.1in}

For a flat family of deformations of a (possibly) nodal marked symmetric surface 
as in~\eref{cUsymmdfn_e} and $\na$ and $A$ as above, we similarly define 
$$\wt\la_{\fF}\big(\cD_{(V,\vph)}\!\big)
=\bigg(\bigotimes_{i=1}^k\!s_i^{\R*}\la(V^{\vph})\bigg)^{\!\!*}
\!\!\otimes_{\R}\la\big(\cD_{(V,\vph)}\!\big)
\equiv\bigsqcup_{\t\in\De_{\R}}\!\!
\{\t\}\!\times\!\wt\la_{\cC_{\t}}\!\big(D_{(V,\vph)}\!\big).$$
This set inherits a topology from $\la\big(D_{(V,\vph)}\!\big)$, 
which makes it into a real line bundle over~$\De_{\R}$.\\

We describe the behavior of the orientations induced by relative $\Spin$- and 
$\Pin$-structures under degenerations to two types of nodal surfaces.
First, let 
\BE{cCdegC_e} \cC_0\equiv\big(\Si_0,\si_0,(x_i)_{i\in[k]},(z_i^+,z_i^-)_{i\in[l]}\big)\EE
be a marked symmetric surface which contains precisely one conjugate pair $(\nod^+,\nod^-)$
of nodes and no other nodes.
Let 
$$ \wt\cC_0\equiv\big(\wt\Si_0,
\wt\si_0,(x_i)_{i\in[k]},(z_i^+,z_i^-)_{i\in[l+2]}\big)$$
be the normalization of~$\cC_0$ so that~$\cC_0$ is obtained from $\wt\cC_0$ by identifying
$z_{l+1}^{\pm}$ with $z_{l+2}^{\pm}$ into the node~$\nod^{\pm}$; see Figure~\ref{Cdeg_fig}.
If $(\Si_0,\si_0)$ is a decorated symmetric surface, we take~$\nod^+$ to be the node 
in the chosen half-surface~$\Si^b$.
In this case,  the decorated structure on~$(\Si_0,\si_0)$ induces a decorated structure
on $(\wt\Si_0,\wt\si_0)$ so that~$z_{l+1}^+$ and~$z_{l+2}^+$ lie on 
on the chosen half-surface~$\wt\Si_0^b$.\\

\begin{figure}
\begin{pspicture}(.3,-3)(10,2.5)
\psset{unit=.4cm}
\psellipse(8,0)(3,2)
\psarc(8,2.5){.5}{180}{360}\psarc(8,-2.5){.5}{0}{180}
\psarc(10.5,2.5){2}{0}{180}\psarc(10.5,-2.5){2}{180}{360}
\psarc(11,2.5){3.5}{0}{180}\psarc(11,-2.5){3.5}{180}{360}
\psline(12.5,2.5)(12.5,-2.5)\psline(14.5,2.5)(14.5,-2.5)
\psarc[linewidth=.04](8,.7){1}{210}{330} 
\psarc[linewidth=.04](8,-.7){1}{30}{150} 
\pscircle*(8,2){.2}\pscircle*(8,-2){.2}
\rput(8,1.3){\sm{$\nod^+$}}\rput(8,-1.3){\sm{$\nod^-$}}
\psline{<->}(4,4)(4,-4)\rput(3.2,0){$\si_0$}
\rput(9.5,-7){$\cC_0$}
\psellipse(28,0)(3,2)
\psarc[linewidth=.04](28,.7){1}{210}{330} 
\psarc[linewidth=.04](28,-.7){1}{30}{150} 
\pscircle*(28,2){.2}\pscircle*(28,-2){.2}
\rput(28,2.9){\sm{$z_{l+1}^+$}}\rput(28,-2.9){\sm{$z_{l+1}^-$}}
\pscircle(34.5,0){3}
\pscircle*(34.5,3){.2}\pscircle*(34.5,-3){.2}
\rput(34.5,3.9){\sm{$z_{l+2}^+$}}\rput(34.5,-3.9){\sm{$z_{l+2}^-$}}
\psline{<->}(24,4)(24,-4)\rput(23.2,0){$\wt\si_0$}
\rput(32,-7){$\wt\cC_0$}
\end{pspicture}
\caption{A marked symmetric surface~$\cC_0$ with a conjugate pair of nodes~$\nod^{\pm}$
and its normalization~$\wt\cC_0$}
\label{Cdeg_fig}
\end{figure}

A real bundle pair $(V_0,\vph_0)$ over~$(\Si_0,\si_0)$ lifts to a real bundle pair  
$(\wt{V}_0,\wt\vph_0)$ over~$(\wt\Si_0,\wt\si_0)$.
A real CR-operator~$D_0$ on $(V_0,\vph_0)$ lifts to a real CR-operator~$\wt{D}_0$ 
on $(\wt{V}_0,\wt\vph_0)$ so that there is a natural exact triple
\BE{Cdegses_e}0\lra D_0\lra \wt{D}_0\lra V_0\big|_{\nod^+}\lra0,  \quad
\wt\xi\lra\wt\xi\big(z_{l+2}^+\!\big)\!-\!\wt\xi\big(z_{l+1}^+\!\big),\EE
of Fredholm operators. It induces an isomorphism
\BE{Cdegses_e2} \la(D_0)\!\otimes\!\la\big(V_0|_{\nod^+}\big)\approx\la(\wt{D}_0)\,.\EE
The real line $\la(V_0|_{\nod^+})$ is oriented by the complex orientation of~$V_0|_{\nod^+}$.
An orientation of~$\wt{D}_0$ thus induces an orientation of~$D_0$.
This orientation does not change if the points $z_{l+1}^+$ and $z_{l+2}^+$
are interchanged because the real dimension of $V_0|_{\nod^+}$ is even.
If $\os_0$ 
is a relative $\OSpin$-structure on the real vector bundle~$V_0^{\vph_0}$ 
over $\Si_0^{\si_0}\!\subset\!\Si_0$ and 
$\wt\os_0$ is its lift to a relative $\OSpin$-structure on 
the real vector bundle~$\wt{V}_0^{\wt\vph_0}$ over $\wt\Si_0^{\wt\si_0}\!\subset\!\wt\Si_0$, 
we denote by $\fo_{\os_0}(V_0,\vph_0)$ the orientation of~$D_0$ induced by
the orientation $\fo_{\wt\os_0}(\wt{V}_0,\wt\vph_0)$ of~$\wt{D}_0$  via~\eref{Cdegses_e2}.
We call $\fo_{\os_0}(V_0,\vph_0)$ the \sf{intrinsic orientation induced by~$\os_0$}.
If $\fp_0$ is a relative $\Pin^{\pm}$-structure on~$V_0^{\vph_0}$ and
$(V_0,\vph_0)$ is $\cC_0$-balanced, we similarly define the \sf{intrinsic orientation}
$\fo_{\cC_0;\fp_0}(V_0,\vph_0)$ on $\wt\la_{\cC_0}(D_0)$  \sf{induced by~$\fp_0$}.\\

Suppose in addition~\eref{cUsymmdfn_e} is a flat family of deformations of~$\cC_0$
and $(V,\vph)$ is a real bundle pair over $(\cU,\wt\fc)$ extending~$(V_0,\vph_0)$. 
The decorated structure on~$(\Si_0,\si_0)$ determines a subspace $\cU^b\!\subset\!\cU$ so that 
$\cU^b$ is a manifold with boundary, \hbox{$\Si_{\t}^b\!\equiv\!\cU^b\!\cap\!\Si_{\t}$} 
is a bordered, possibly nodal, surface for every \hbox{$\ft\!\in\!\De$},
$$\cU=\cU^b\!\cup\!\wt\fc(\cU^b),
 \qquad\hbox{and}\qquad \cU^b\!\cap\!\wt\fc(\cU^b)=\prt\cU^b\,.$$
The ordering of the components of~$\Si_0^{\si_0}$ likewise induces an ordering
of the components of~$\Si_{\t}^{\si_{\t}}$ for   every $\t\!\in\!\De_{\R}$.
Thus, a decorated structure on $(\Si_0,\si_0)$ determines 
a decorated structure on $(\Si_{\t},\si_{\t})$ for every $\t\!\in\!\De_{\R}$
in a continuous manner.\\

A relative $\OSpin$-structure~$\os_0$ on the real vector bundle~$V_0^{\vph_0}$ over
$\Si_0^{\si_0}\!\subset\!\Si_0$ extends to a relative $\OSpin$-structure~$\os$
on the real vector bundle~$V^{\vph}$ over~$\cU^{\wt\fc}\!\subset\!\cU$.
The latter in turn restricts to a relative $\OSpin$-structure~$\os_{\t}$ 
on the real vector bundle~$V_{\t}^{\vph_{\t}}$ over $\Si_{\t}^{\si_{\t}}\!\subset\!\Si_{\t}$ 
and thus determines an orientation $\fo_{\os_{\t}}(V_{\t},\vph_{\t})$ of $D_{(V,\vph);\t}$ 
for every~$\t\!\in\!\De_{\R}^*$ which varies continuously with~$\t$,
i.e.~an orientation of the restriction of the real line bundle~\eref{detLB_e} to~$\De_{\R}^*$.
Since the codimension of \hbox{$\De_{\R}\!-\!\De_{\R}^*$} in~$\De_{\R}$ is~2 in this case,
the last orientation extends continuously over the entire base~$\De_{\R}$.
We denote the restriction of this orientation to $D_{(V_0,\vph_0)}$
by $\fo_{\os_0}'(V_0,\vph_0)$ and call it the \sf{limiting orientation induced by~$\os$}.
If $\fp_0$ is a relative $\Pin^{\pm}$-structure on~$V_0^{\vph_0}$ and
$(V_0,\vph_0)$ is $\cC_0$-balanced, we similarly define 
the \sf{limiting orientation} $\fo_{\cC_0;\fp_0}'(V_0,\vph_0)$
on $\wt\la_{\cC_0}(D_{(V_0,\vph_0)})$  \sf{induced by~$\fp_0$}.
These limiting orientations depend only on~$(V_0,\vph_0)$,
$\os_0$, and~$\fp_0$, and not on~$(\cU,\wt\fc)$ or~$(V,\vph)$.

\setcounter{temp}{\value{CROrient}}
\stepcounter{temp}

\begingroup
\def\theCROrient{\thetemp$C$}
\begin{CROrient}[Degenerations: $C$ nodes]\label{CROdegenC_prop}
Suppose $\cC_0$ is a decorated marked symmetric surface which contains precisely 
one conjugate pair $(\nod^+,\nod^-)$ of nodes and no other nodes and
$(V_0,\vph_0)$ is a real bundle pair over~$\cC_0$.
\begin{enumerate}[label=(\alph*),leftmargin=*]

\item\label{CdegenSpin_it} 
The intrinsic and limiting orientations of $D_{(V_0,\vph_0)}$ 
induced by a relative $\OSpin$-structure on the vector bundle~$V_0^{\vph_0}$ over
$\Si_0^{\si_0}\!\subset\!\Si_0$ are the same.

\item\label{CdegenPin_it} If $(V_0,\vph_0)$ is $\cC_0$-balanced,
the intrinsic and limiting orientations on $\wt\la_{\cC_0}(D_{(V_0,\vph_0)})$ 
induced by a relative $\Pin^{\pm}$-structure on~$V_0^{\vph_0}$
are the same.
 
\end{enumerate}
\end{CROrient}
\endgroup

\vspace{.1in}

Suppose~\eref{cCdegC_e} is a marked symmetric surface which contains precisely 
one $H$~node~$\nod$ (i.e.~a non-isolated point of~$\Si^{\si}$) and no other nodes 
so that~$\nod$ splits a topological component~$\Si_{\bu}$ of~$\Si$ into
two irreducible surfaces, $(\Si_{\bu1},\si_{\bu1})$ and~$(\Si_{\bu2},\si_{\bu2})$.
Such a node, which we call an \sf{H3~node}, is one of the three possible types
of $H$~nodes a symmetric surface~$(\Si,\si)$ may have;
see \cite[Section~3.2]{RealGWsI}.
Let
\BE{wtcCH3deg_e} \wt\cC_0\equiv\big(\wt\Si_0,
\wt\si_0,(\nod_1,\nod_2,(x_i)_{i\in[k]}),(z_i^+,z_i^-)_{i\in[l]}\big)\EE
be the normalization of~$\cC_0$ so that 
$\Si_{\bu1},\Si_{\bu2}\!\subset\!\wt\Si_0$ and
$\cC_0$ is obtained from $\wt\cC_0$ by identifying
the real marked points $\nod_1\!\in\!\Si_{\bu1}$ and $\nod_2\!\in\!\Si_{\bu2}$
into the node~$\nod$; see Figure~\ref{H3deg_fig}.
We view these two points as the first two real marked points of~$\wt\cC_0$ 
when specifying orientations in the $\Pin^{\pm}$-case below.
Let $S^1_{\bu1}\!\subset\!\Si_{\bu1}$ and $S^1_{\bu2}\!\subset\!\Si_{\bu2}$ be 
the components of $\wt\Si_0^{\wt\si_0}$ containing~$\nod_1$ and~$\nod_2$, respectively.
Let $k_{\bu1}$ and $k_{\bu2}$ be the numbers of real marked points 
(not including~$\nod_1$ and~$\nod_2$) carried by these components.
We denote by $\wt\cC_1'$ and~$\wt\cC_2'$ the marked curves obtained from~$\wt\cC_0$ by 
dropping~$\nod_1$ and~$\nod_2$, respectively.
If $(\Si_0,\si_0)$ is a decorated surface.
we define a decorated structure on the normalization $(\wt\Si_0,\wt\si_0)$ by 
\begin{enumerate}[label=$\bu$,leftmargin=*]

\item taking the distinguished half-surface \hbox{$\wt\Si_0^b\!\subset\!\wt\Si_0$} to
be the preimage of the distinguished half-surface \hbox{$\Si_0^b\!\subset\!\Si_0$} and

\item replacing the position~$r_{\bu}$ of the topological component of~$\Si_0^{\si_0}$ containing
the node~$\nod$ (i.e.~the wedge $S_{\bu1}^1\!\cup_{\nod}\!S_{\bu2}^1$)
in the ordering of~$\pi_0(\Si_0^{\si_0})$ by 
$S^1_{\bu1}$ followed by~$S^1_{\bu2}$.

\end{enumerate}

\vspace{.1in}

\begin{figure}
\begin{pspicture}(.3,-2)(10,2)
\psset{unit=.4cm}
\psellipse(8,0)(3,2)
\psarc[linewidth=.04](8,.7){1}{210}{330} 
\psarc[linewidth=.04](8,-.7){1}{30}{150} 
\pscircle(14,0){3}
\pscircle*(11,0){.2}\rput(11.8,0.1){\sm{$\nod$}}
\psline{<->}(4,4)(4,-4)\rput(3.2,0){$\si_0$}
\rput(9.5,-4){$\cC_0$}
\psellipse(28,0)(3,2)
\psarc[linewidth=.04](28,.7){1}{210}{330} 
\psarc[linewidth=.04](28,-.7){1}{30}{150} 
\pscircle*(31,0){.2}\rput(30.1,0.1){\sm{$\nod_1$}}
\pscircle(34.5,0){3}
\pscircle*(31.5,0){.2}\rput(32.6,0.1){\sm{$\nod_2$}}
\psline{<->}(24,4)(24,-4)\rput(23.2,0){$\wt\si_0$}
\rput(32,-4){$\wt\cC_0$}
\end{pspicture}
\caption{A marked symmetric surface~$\cC_0$ with a real $H3$ node~$\nod$
and its normalization~$\wt\cC_0$}
\label{H3deg_fig}
\end{figure}

A real bundle pair $(V_0,\vph_0)$ over~$(\Si_0,\si_0)$ lifts to a real bundle pair  
$(\wt{V}_0,\wt\vph_0)$ over~$(\wt\Si_0,\wt\si_0)$.
A real CR-operator~$D_0$ on $(V_0,\vph_0)$ lifts to a real CR-operator~$\wt{D}_0$ 
on $(\wt{V}_0,\wt\vph_0)$ so that there is a natural exact triple
\BE{Rdegses_e}0\lra D_0\lra \wt{D}_0\lra V_0^{\vph_0}\big|_{\nod}\lra0,  \qquad
\wt\xi\lra\wt\xi(\nod_2)\!-\!\wt\xi(\nod_1),\EE
of Fredholm operators. It induces an isomorphism
\BE{Rdegses_e2} \la(D_0)\!\otimes\!\la\big(V_0^{\vph_0}|_{\nod}\big)\approx\la(\wt{D}_0)\,.\EE
Orientations on $V_0^{\vph_0}$ and~$\wt{D}_0$ thus induce an orientation on~$D_0$.
If $\os_0\!\in\!\Spin_{\Si}(V^{\vph},\fo)$ 
is a relative $\OSpin$-structure on the real vector bundle~$V_0^{\vph_0}$ 
over $\Si_0^{\si_0}\!\subset\!\Si_0$ and 
$\wt\os_0$ is its lift to a relative $\OSpin$-structure on 
the real vector bundle~$\wt{V}_0^{\wt\vph_0}$ over $\wt\Si_0^{\wt\si_0}\!\subset\!\wt\Si_0$, 
we denote by $\fo_{\os_0}(V_0,\vph_0)$ the orientation of~$D_0$ induced by
the orientations~$\fo_{\nod}$ of~$V_0^{\vph_0}|_{\nod}$
and $\fo_{\wt\os_0}(\wt{V}_0,\wt\vph_0)$ of~$\wt{D}_0$  via~\eref{Rdegses_e2}.
We call $\fo_{\os_0}(V_0,\vph_0)$ the \sf{intrinsic orientation induced by~$\os_0$}.
If $(V_0,\vph_0)$ is $\cC_0$-balanced, then $(\wt{V}_0,\wt\vph_0)$ is either 
$\wt\cC_1'$-balanced or $\wt\cC_2'$-balanced (but not both).
If the $i$-th case holds, the orientation $\fo_{\wt\cC_i';\wt\fp_0}(\wt{V}_0,\wt\vph_0)$ 
induced by 
the lift~$\wt\fp_0$ of a relative $\Pin^{\pm}$-structure~$\fp$ on~$V_0^{\vph_0}$
to~$\wt{V}_0^{\wt\vph_0}$ determines a homotopy class of isomorphisms
\BE{Rdegses_e4}
\la(\wt{D}_0)\approx \la_{\wt\cC_i'}^{\R}(\wt{V}_0,\wt\vph_0)
\approx \la_{\cC_0}^{\R}(V_0,\vph_0)\!\otimes\!\la\big(V_0^{\vph_0}|_{\nod}\big).\EE
Combining~\eref{Rdegses_e2} and~\eref{Rdegses_e4}, 
we obtain an orientation $\fo_{\cC_0;\fp_0}(V_0,\vph_0)$ on $\wt\la_{\cC_0}(D_0)$.
We call $\fo_{\cC_0;\fp_0}(V_0,\vph_0)$ the \sf{intrinsic orientation induced by~$\fp_0$}.\\  

If~\eref{cUsymmdfn_e} is a flat family of deformations of~$\cC_0$,
the topological components of the fixed locus of $(\Si_{\t},\si_{\t})$ for each $\t\!\in\!\De_{\R}$
correspond to the components of $\Si_0^{\si_0}$ and thus inherit an ordering
from the chosen ordering of~$\pi_0(\Si_0^{\si_0})$.
Let $\cU_{\bu}\!\subset\!\cU$ be the topological component containing~$\Si_{\bu}$.
The restriction of~$\pi$ to $(\cU\!-\!\cU_{\bu})|_{\De_{\R}}$ is a topologically trivial
fiber bundle of symmetric surfaces and 
a half-surface of $\Si_0\!-\!\Si_{\bu}$ determines a half-surface of 
every fiber of this restriction.
For each fiber~$\Si_{\t;\bu}$ of the restriction of~$\pi$ to~$\cU_{\bu}|_{\De_{\R}}$,
we take $\Si_{\t;\bu}^b\!\subset\!\Si_{\t;\bu}$ to be 
the half-surface so that 
for every point~$z$ in $\Si_{\t;\bu}^b\!-\!\Si_{\t}^{\si_{\t}}$ there exists a path
$$\al\!:[0,1]\lra\cU|_{\De_{\R}}\!-\!\cU^{\wt\fc} \qquad\hbox{s.t.}\quad
\al(0)\in\Si_{\bu}^b\!-\!\Si_{\bu2},\quad \al(1)=z.$$
Thus, a decorated structure on $(\Si_0,\si_0)$ determines 
a decorated structure on $(\Si_{\t},\si_{\t})$ for every $\t\!\in\!\De_{\R}^*$
in a continuous manner.
The space $\De_{\R}^*$ has two topological components.
We denote by $\De_{\R}^+\!\subset\!\De_{\R}^*$ the topological component so that 
for every $\t\!\in\!\De_{\R}^+$ and $z\!\in\!\Si_{\t;\bu}^b\!-\!\Si_{\t}^{\si_{\t}}$
there exists a path
$$\al\!:[0,1]\lra\cU|_{\De_{\R}}\!-\!\cU^{\wt\fc} \qquad\hbox{s.t.}\quad
\al(0)\in\Si_{\bu}^b\!-\!\Si_{\bu1},\quad \al(1)=z.$$
We denote the other component of~$\De_{\R}^*$ by~$\De_{\R}^-$.\\

Suppose in addition that $(V,\vph)$ is a real bundle pair over $(\cU,\wt\fc)$.
As above the CROrient~\ref{CROdegenC_prop} property, 
a decorated structure on~$(\Si_0,\si_0)$ and
a relative $\OSpin$-structure~$\os_0$ on the real vector bundle~$V_0^{\vph_0}$ over
$\Si_0^{\si_0}\!\subset\!\Si_0$ determine an orientation of the restriction of 
the real line bundle~\eref{detLB_e} to~$\De_{\R}^*$.
The restrictions of this orientation to $\la(\cD_{(V,\vph)})|_{\De_{\R}^+}$ 
and $\la(\cD_{(V,\vph)})|_{\De_{\R}^-}$ extend
continuously to orientations $\fo_{\os_0}^+(V_0,\vph_0)$ and $\fo_{\os_0}^-(V_0,\vph_0)$, 
respectively, of $\la(D_{(V_0,\vph_0)})$;
we call them the \sf{limiting orientations induced by~$\os_0$}.
If $\fp_0$ is a relative $\Pin^{\pm}$-structure on~$V_0^{\vph_0}$ and
$(V_0,\vph_0)$ is $\cC_0$-balanced, we similarly define 
the \sf{limiting orientations}
$\fo_{\cC_0;\fp_0}^+(V_0,\vph_0)$ and $\fo_{\cC_0;\fp_0}^-(V_0,\vph_0)$
of $\wt\la_{\cC_0}(D_{(V_0,\vph_0)})$  \sf{induced by~$\fp_0$}.
All four limiting orientations depend only on~$(V_0,\vph_0)$,
$\os_0$, and~$\fp_0$, and not on~$(\cU,\wt\fc)$ or~$(V,\vph)$.\\

With the notation as above, let $r(\cC_0)\!=\!1$ if either $k_{\bu2}\!=\!0$ or
\BE{rcC0dfn_e}\inf\!\big\{i\!\in\![k]\!:x_i\!\in\!S^1_{\bu1}\big\}<
\min\!\big\{i\!\in\![k]\!:x_i\!\in\!S^1_{\bu2}\big\};\EE
otherwise, we take $r(\cC_0)\!=\!2$.
For $r\!=\!1,2$, $S^1_{\bu r}\!\subset\!\Si_{\bu r}$ is oriented as the boundary of~$\Si_{\bu r}^b$.
Let $j_r'(\cC_0)\!\in\!\Z^{\ge0}$ be the number of real marked points that lie 
on the oriented arc of~$S_{\bu r}^1$ between the nodal point~$\nod_r$ 
and  the real marked point~$x_i$ with the minimal index $i\!\in\![k]$ on~$S^1_{\bu r}$;
if $k_{\bu r}\!=\!0$, we take $j_r'(\cC_0)\!=\!0$.
Define
$$\de_{\R}(\cC_0)
=j_{r(\cC_0)}'(\cC_0)\!+\!1\!+\!\big(r(\cC_0)\!-\!1\big)k_{\bu1}\,.$$
If in addition $(V_0,\vph_0)$ is a $\cC_0$-balanced rank~$n$ real bundle pair, let
\begin{gather*}
r(V_0,\vph_0)=\begin{cases}1,&\hbox{if}~
k_{\bu1}\!+\!2\Z\!\neq\!\lr{w_1(V_0^{\vph_0}),[S^1_{\bu1}]_{\Z_2}},\,
k_{\bu2}\!+\!2\Z\!=\!\lr{w_1(V_0^{\vph_0}),[S^1_{\bu2}]_{\Z_2}};\\
2,&\hbox{if}~
k_{\bu1}\!+\!2\Z\!=\!\lr{w_1(V_0^{\vph_0}),[S^1_{\bu1}]_{\Z_2}},\,
k_{\bu2}\!+\!2\Z\!\neq\!\lr{w_1(V_0^{\vph_0}),[S^1_{\bu2}]_{\Z_2}};
\end{cases}\\
\de_{\R}(V_0,\vph_0)=\big(k_{\bu r(V_0,\vph_0)}\!-\!1\big)
\big(j_{r(V_0,\vph_0)}'(\cC_0)\!+\!r(V_0,\vph_0)\!\big)\,.
\end{gather*}

\begingroup
\def\theCROrient{\thetemp$H3$}
\begin{CROrient}[Degenerations: $H3$ nodes]\label{CROdegenH3_prop}
Suppose $\cC_0$ is a decorated marked symmetric surface which contains precisely 
one $H3$ node and no other nodes and
$(V_0,\vph_0)$  is a rank~$n$ real bundle pair over~$\cC_0$.
\begin{enumerate}[label=(\alph*),leftmargin=*]

\item\label{HdegenSpin_it} The limiting orientation $\fo_{\os_0}^+(V_0,\vph_0)$
of $D_{(V_0,\vph_0)}$ 
induced by a relative $\OSpin$-structure~$\os_0$ on the real vector bundle~$V_0^{\vph_0}$ over
$\Si_0^{\si_0}\!\subset\!\Si_0$ is the same as the intrinsic orientation
$\fo_{\os_0}(V_0,\vph_0)$ if and only if $n(|\pi_0(\Si_0^{\si_0})|\!-\!r_{\bu})$ is even.

\item\label{HdegenPin_it} If $(V_0,\vph_0)$ is $\cC_0$-balanced,
the limiting orientation $\fo_{\cC_0;\fp_0}^+(V_0,\vph_0)$ on $\wt\la_{\cC_0}(D_{(V_0,\vph_0)})$ 
induced by a relative $\Pin^{\pm}$-structure~$\fp_0$ on~$V_0^{\vph_0}$
 is the same as the intrinsic orientation $\fo_{\cC_0;\fp_0}(V_0,\vph_0)$
if and only~if 
\begin{equation*}\begin{split}
&(n\!+\!1)\bigg(\!\!\big(k_{r_{\bu}}\!-\!1\!\big)\de_{\R}(\cC_0)\!+\!\de_{\R}(V_0,\vph_0)
+\!\!\!\!\!\!\!\sum_{\begin{subarray}{c}S_r^1\in\pi_0(\Si_0^{\si_0})\\ 
r>r_{\bu}\end{subarray}}\!\!\!\!\!\!\!\!\!\blr{w_1(V_0^{\vph_0}),[S_r^1]_{\Z_2}}
\!\!\bigg)\\
&\hspace{1.5in}=n\big(r(V_0,\vph_0)\!-\!1\big)\!+\!2\Z
+\begin{cases}
nk,&\hbox{if}~\fp_0\!\in\!\cP^-_{\Si_0}(V_0^{\vph_0});\\
(n\!+\!1)k,&\hbox{if}~\fp_0\!\in\!\cP^+_{\Si_0}(V_0^{\vph_0}).
\end{cases}
\end{split}\end{equation*}
\end{enumerate}
\end{CROrient}
\endgroup

\addtocounter{CROrient}{-1}

\begin{rmk}\label{CROdegenH3_rmk}
By CROrient~\ref{CROos_prop}\ref{osflip_it} and~\ref{CROdegenH3_prop}\ref{HdegenSpin_it},
the limiting orientation $\fo_{\os_0}^-(V_0,\vph_0)$ is the same as the intrinsic orientation
$\fo_{\os_0}(V_0,\vph_0)$ if and only~if 
$$ n\big(|\pi_0(\Si_0^{\si_0})|\!-\!r_{\bu}\big)+2\Z =\vp_{\os_0}(\Si_{\bu2})\in\Z_2.$$ 
By CROrient~\ref{CROp_prop}\ref{pflip_it} and~\ref{CROdegenH3_prop}\ref{HdegenPin_it}, 
the limiting orientation $\fo_{\cC_0;\fp_0}^-(V_0,\vph_0)$
is the same as the intrinsic orientation $\fo_{\cC_0;\fp_0}(V_0,\vph_0)$
if and only~if
\begin{equation*}\begin{split}
&(n\!+\!1)\bigg(\!\!\big(k_{r_{\bu}}\!-\!1\!\big)\de_{\R}(\cC_0)\!+\!\de_{\R}(V_0,\vph_0)
+\!\!\!\!\!\!\!\sum_{\begin{subarray}{c}S_r^1\in\pi_0(\Si_0^{\si_0})\\ 
r>r_{\bu}\end{subarray}}\!\!\!\!\!\!\!\!\!\blr{w_1(V_0^{\vph_0}),[S_r^1]_{\Z_2}}
\!\!\bigg)\\
&\hspace{1in}+
\!\!\!\!\!\sum_{\begin{subarray}{c}S_r^1\in\pi_0(\Si_{\bu2}^{\si_{\bu2}})\\ 
S_r^1\neq S_{\bu2}^1\end{subarray}}\!\!\!
\bigg(\!\!\!\binom{k_r\!-\!1}{1}\!+\!n\binom{k_r\!-\!1}{2}\!\!\!\bigg)
+r(V_0,\vph_0)\!\bigg(\!\!k_{\bu2}\!+\!n\binom{k_{\bu2}}{2}\!\!\!\bigg)\\
&\hspace{1in}
+\big(r(V_0,\vph_0)\!-\!1\big)n
\bigg(\!\!\!\binom{k_{\bu2}\!-\!1}{1}\!+\!
\binom{k_{\bu2}\!-\!1}{2}\!\!\!\bigg)
+\big(r(\cC_0)\!-\!1\big)(n\!+\!1)(k_{\bu1}\!-\!1)(k_{\bu2}\!-\!1) \\
&\hspace{.5in}=\vp_{\fp_0}(\Si_{\bu2})+n\big(r(V_0,\vph_0)\!-\!1\big)+\begin{cases}
nk,&\hbox{if}~\fp_0\!\in\!\cP^-_{\Si_0}(V_0^{\vph_0});\\
(n\!+\!1)k,&\hbox{if}~\fp_0\!\in\!\cP^+_{\Si_0}(V_0^{\vph_0}).
\end{cases}
\end{split}\end{equation*}
The last term on the third line above accounts for the effect of the change in
$\Si_{\bu2}^b$ on $j_{r(\cC_0)}'(\cC_0)$ if $r(\cC_0)\!=\!2$.  
The preceding term includes the  effect of this change on
$j_{r(V_0,\vph_0)}'(\cC_0)$ if $r(V_0,\vph_0)\!=\!2$.
\end{rmk}

\subsection{Some implications}
\label{OrientEg_subs}

With the setup as in the CROrient~\ref{CROdegenC_prop} property,
suppose $(\Si_{01},\si_{01})$ and $(\Si_{02},\si_{02})$ are symmetric surfaces so~that 
\BE{CROdegenC_e3} 
\Si_{02}^{\si_{02}}=\eset, \qquad
\big(\Si_0,\si_0\big)=
\big(\Si_{01},\si_{01}\big)\!\sqcup\!\big(\Si_{02},\si_{02}\big);\EE
see Figure~\ref{Csplit_fig}.
We denote by $(V_1,\vph_1)$ and $(V_2,\vph_2)$ the restrictions of
the real bundle pair~$(V_0,\vph_0)$ to $(\Si_{01},\si_{01})$ and $(\Si_{02},\si_{02})$,
respectively, and by~$\os_1$ (resp~$\fp_1$)
the restriction of the relative $\OSpin$-structure~$\os_0$ (resp.~$\Pin^{\pm}$-structure~$\fp_0$)
to relative $\OSpin$-structure (resp.~$\Pin^{\pm}$-structure)
on the real vector bundle~$V_0^{\vph_0}$ over $\Si_0^{\si_0}\!\subset\!\Si_0$
to the real vector bundle $V_1^{\vph_1}$  over $\Si_{01}^{\si_{01}}\!\subset\!\Si_{01}$.
Let $D_0$ and $\wt{D}_0$ be as in~\eref{Cdegses_e} and
$\wt{D}_{01}$, $\wt{D}_{02}$, and $\wt{D}_{02}^b$  be  the restrictions of~$\wt{D}_0$
to real CR-operators on $(V_1,\vph_1)$, $(V_2,\vph_2)$, and~$V_2$, respectively.\\

\begin{figure}
\begin{pspicture}(-.3,-2.5)(10,2.3)
\psset{unit=.4cm}
\psellipse(8,0)(3,2)
\psarc[linewidth=.04](8,.7){1}{210}{330} 
\psarc[linewidth=.04](8,-.7){1}{30}{150} 
\pscircle(8,3.5){1.5}\pscircle(8,-3.5){1.5}
\pscircle*(8,2){.2}\pscircle*(8,-2){.2}
\rput(8,1.3){\sm{$\nod^+$}}\rput(8,-1.3){\sm{$\nod^-$}}
\psline{<->}(4,4)(4,-4)\rput(3.2,0){$\si_0$}
\rput(8,-6){$\cC_0$}
\psellipse(28,0)(3,2)
\psarc[linewidth=.04](28,.7){1}{210}{330} 
\psarc[linewidth=.04](28,-.7){1}{30}{150} 
\pscircle*(28,2){.2}\pscircle*(28,-2){.2}
\rput(28,2.9){\sm{$z_{l+1}^+$}}\rput(28,-2.9){\sm{$z_{l+1}^-$}}
\pscircle(33,2){1.5}\pscircle(33,-2){1.5}
\pscircle*(33,.5){.2}\pscircle*(33,-.5){.2}
\rput(33,1.5){\sm{$z_{l+2}^+$}}\rput(33,-1.4){\sm{$z_{l+2}^-$}}
\psline{<->}(24,4)(24,-4)\rput(23.2,0){$\wt\si_0$}
\rput(28.2,-6){$\wt\cC_{01}$}\rput(33.2,-6){$\wt\cC_{02}$}
\end{pspicture}
\caption{A marked symmetric surface~$\cC_0$ with a conjugate pair of nodes~$\nod^{\pm}$
as in Corollary~\ref{CROdegenC_crl} and its normalization~$\wt\cC_0$}
\label{Csplit_fig}
\end{figure}

The orientation $\fo_{\os_1}\!(V_1,\vph_1)$ of $\wt{D}_{01}$ and
the complex orientation of $\wt{D}_{02}^b$ determine 
an orientation $\fo_{\os_0}^{\C}(\wt{V}_0,\wt\vph_0)$ of~$\wt{D}_0$ 
via the isomorphisms
\BE{CintrvsCsplit_e}\la\big(\wt{D}_0\big)\approx \la\big(\wt{D}_{01}\big)\!\otimes\!\la\big(\wt{D}_{02}\big)
\approx \la\big(\wt{D}_{01}\big)\!\otimes\!\la\big(\wt{D}_{02}^b\big)\EE
as in~\eref{CRODisjUn_e0} and~\eref{CROdoubl_e}
and thus an orientation 
$\fo_{\os_0}^{\C}(V_0,\vph_0)$ of~$D_0$ via the isomorphism~\eref{Cdegses_e2}.
We call $\fo_{\os_0}^{\C}(V_0,\vph_0)$ the \sf{$\C$-split orientation induced by~$\os_0$}.
Suppose in addition that $(V_0,\vph_0)$ is $\cC_0$-balanced and
\BE{CROdegenC_e3b}\wt\cC_0 \equiv \wt\cC_{01}\!\sqcup\!\wt\cC_{02}\EE
is the decomposition induced by the decomposition in~\eref{CROdegenC_e3}.
The orientation $\fo_{\wt\cC_{01};\fp_1}\!(V_1,\vph_1)$ of $\wt\la_{\wt\cC_{01}}\!(\wt{D}_{01})$ 
and the complex orientation of $\wt{D}_{02}^b$ then determine 
an orientation $\fo_{\wt\cC_0;\fp_0}^{\C}\!(\wt{V}_0,\wt\vph_0)$ on
$\wt\la_{\wt\cC_0}(\wt{D}_0)$ via the isomorphisms~\eref{CintrvsCsplit_e}
and thus an orientation 
$\fo_{\cC_0;\fp_0}^{\C}(V_0,\vph_0)$ on $\wt\la_{\cC_0}\!(D_0)$ via 
the isomorphism~\eref{Cdegses_e2}.
We call $\fo_{\cC_0;\fp_0}^{\C}(V_0,\vph_0)$ the \sf{$\C$-split orientation induced by~$\fp_0$}.\\

The CROrient~\ref{CRODisjUn_prop}, \ref{CROPin2SpinRed_prop},  and
\ref{CROos_prop}\ref{dltCROos_it} properties compare
 the intrinsic orientation $\fo_{\os_0}(V_0,\vph_0)$ 
(resp. $\fo_{\cC_0;\fp_0}(V_0,\vph_0)$)
with the $\C$-split orientation $\fo_{\os_0}^{\C}(V_0,\vph_0)$
(resp.~$\fo_{\cC_0;\fp_0}^{\C}(V_0,\vph_0)$). 
Combining these properties with the CROrient~\ref{CROdegenC_prop} property, 
we obtain the following conclusion.

\begin{crl}\label{CROdegenC_crl}
Suppose $\cC_0$ and $(V_0,\vph_0)$ are as in  
the CROrient~\ref{CROdegenC_prop} property and \eref{CROdegenC_e3} holds.
\begin{enumerate}[label=(\alph*),leftmargin=*]

\item\label{CdegenSpin2_it} The limiting orientation $\fo_{\os_0}'(V_0,\vph_0)$
of $D_{(V_0,\vph_0)}$ 
induced by a relative $\OSpin$-structure~$\os_0$ on the real vector bundle~$V_0^{\vph_0}$ over
$\Si_0^{\si_0}\!\subset\!\Si_0$ is the same as the $\C$-split orientation
$\fo_{\os_0}^{\C}(V_0,\vph_0)$ 
if and only if \hbox{$\lr{w_2(\os_0),[\Si_{02}^b]_{\Z_2}}\!=\!0$}.

\item\label{CdegenPin2_it} If $(V_0,\vph_0)$ is $\cC_0$-balanced and~\eref{CROdegenC_e3b} holds, then
the limiting orientation $\fo_{\cC_0;\fp_0}'(V_0,\vph_0)$ on $\wt\la_{\cC_0}(D_{(V_0,\vph_0)})$ 
induced by a relative $\Pin^{\pm}$-structure~$\fp_0$ on~$V_0^{\vph_0}$
is the same as the $\C$-split orientation $\fo_{\cC_0;\fp_0}^{\C}(V_0,\vph_0)$
if and only~if \hbox{$\lr{w_2(\fp_0),[\Si_{02}^b]_{\Z_2}}\!=\!0$}.
\end{enumerate}
\end{crl}

\vspace{.1in}

With the setup as in the CROrient~\ref{CROdegenH3_prop} property,
let $(\Si_{01},\si_{01})$ and $(\Si_{02},\si_{02})$ be symmetric surfaces so~that 
\BE{CROdegen_e3} \Si_{\bu1}\subset\Si_{01}, \qquad
\Si_{\bu2}\subset\Si_{02}, \qquad
\big(\Si_0,\si_0\big)=
\big(\Si_{01},\si_{01}\big)\!\sqcup\!\big(\Si_{02},\si_{02}\big)\EE
and $(V_1,\vph_1)$ and $(V_2,\vph_2)$ be as below~\eref{CROdegenC_e3}.
We denote by $\os_1$ and~$\os_2$ (resp~$\fp_1$ and~$\fp_2$) 
the restrictions of the relative $\OSpin$-structure~$\os_0$ (resp.~$\Pin^{\pm}$-structure~$\fp_0$)
on the real vector bundle~$V_0^{\vph_0}$ over $\Si_0^{\si_0}\!\subset\!\Si_0$
to relative $\OSpin$-structures (resp.~$\Pin^{\pm}$-structures)
on the real vector bundles $V_1^{\vph_1}$ and $V_2^{\vph_2}$
over $\Si_{01}^{\si_{01}}\!\subset\!\Si_{01}$ and $\Si_{02}^{\si_{02}}\!\subset\!\Si_{02}$,
respectively.
Let $D_0$ and $\wt{D}_0$ be as in~\eref{Rdegses_e} and $\wt{D}_{01}$ and $\wt{D}_{02}$ be 
the restrictions of~$\wt{D}_0$ to real CR-operators on $(V_1,\vph_1)$ and $(V_2,\vph_2)$, 
respectively.\\

The orientations $\fo_{\os_1}(V_1,\vph_1)$ and $\fo_{\os_2}(V_2,\vph_2)$ 
determine an orientation $\fo_{\os}^{\sp}(\wt{V}_0,\wt\vph_0)$ of~$\wt{D}_0$ 
via the isomorphism
\BE{H3intrvssplit_e} 
\la\big(\wt{D}_0\big)\approx \la\big(\wt{D}_{01}\big)\!\otimes\!\la\big(\wt{D}_{02}\big) \EE
as in~\eref{CRODisjUn_e0} and thus an orientation 
$\fo_{\os_0}^{\sp}(V_0,\vph_0)$ of~$D_0$ via the isomorphism~\eref{Rdegses_e2}.
We call $\fo_{\os_0}^{\sp}(V_0,\vph_0)$ the \sf{split orientation induced by~$\os_0$}.
Suppose in addition that $(V_0,\vph_0)$ is $\cC_0$-balanced, $i\!=\!r(V,\vph)$, and
$$\wt\cC_i' \equiv \wt\cC_{i1}'\!\sqcup\!\wt\cC_{i2}'$$
is the decomposition induced by the decomposition in~\eref{CROdegen_e3}.
The orientations $\fo_{\wt\cC_{i1}';\fp_1}(V_1,\vph_1)$ of 
$\wt\la_{\wt\cC_{i1}'}(\wt{D}_{01})$ and
$\fo_{\wt\cC_{i2}';\fp_2}(V_2,\vph_2)$ of  $\wt\la_{\wt\cC_{i2}'}(\wt{D}_{02})$
determine an orientation $\fo_{\wt\cC_i';\fp_0}^{\sp}(\wt{V}_0,\wt\vph_0)$ on
$\wt\la_{\wt\cC_i'}(\wt{D}_0)$ via the isomorphisms~\eref{H3intrvssplit_e} and~\eref{V1V2dfn_e4} 
and thus an orientation 
$\fo_{\cC_0;\fp_0}^{\sp}(V_0,\vph_0)$ on $\wt\la_{\cC_0}(D_0)$ via 
the isomorphisms~\eref{Rdegses_e2} and~\eref{Rdegses_e4}.
We call $\fo_{\cC_0;\fp_0}^{\sp}(V_0,\vph_0)$ the \sf{split orientation induced by~$\fp_0$}.\\

The CROrient~\ref{CRODisjUn_prop} property compares 
the intrinsic orientation $\fo_{\os_0}(V_0,\vph_0)$ (resp.~$\fo_{\cC_0;\fp_0}(V_0,\vph_0)$)
with the split orientation $\fo_{\os_0}^{\sp}(V_0,\vph_0)$
(resp.~$\fo_{\cC_0;\fp_0}^{\sp}(V_0,\vph_0)$) if 
the decomposition in~\eref{CROdegen_e3} respects the decorated structures, 
i.e.~the elements of~$\pi_0(\Si_0^{\si_0})$ are ordered so~that 
all elements of~$\pi_0(\Si_{01}^{\si_{01}})$ precede all elements
of~$\pi_0(\Si_{02}^{\si_{02}})$. 
Combining the CROrient~\ref{CRODisjUn_prop} and~\ref{CROdegenH3_prop} properties, 
we obtain the following conclusion.

\begin{crl}\label{CROdegenH3_crl}
Suppose $\cC_0$ and $(V_0,\vph_0)$ are as in  
the CROrient~\ref{CROdegenH3_prop} property and
the decomposition in~\eref{CROdegen_e3} respects the decorated structures.
\begin{enumerate}[label=(\alph*),leftmargin=*]

\item\label{HdegenSpin2_it} The limiting orientation $\fo_{\os_0}^+(V_0,\vph_0)$
of $D_{(V_0,\vph_0)}$ 
induced by a relative $\OSpin$-structure~$\os_0$ on the real vector bundle~$V_0^{\vph_0}$ over
$\Si_0^{\si_0}\!\subset\!\Si_0$ is the same as the split orientation
$\fo_{\os_0}^{\sp}(V_0,\vph_0)$ 
if and only if $n(|\pi_0(\Si_0^{\si_0})|\!-\!r_{\bu})$ is even.

\item\label{HdegenPin2_it} If $(V_0,\vph_0)$ is $\cC_0$-balanced,
the limiting orientation $\fo_{\cC_0;\fp_0}^+(V_0,\vph_0)$ on $\wt\la_{\cC_0}(D_{(V_0,\vph_0)})$ 
induced by a relative $\Pin^{\pm}$-structure~$\fp_0$ on~$V_0^{\vph_0}$
 is the same as the split orientation $\fo_{\cC_0;\fp_0}^{\sp}(V_0,\vph_0)$
if and only~if 
\begin{equation*}\begin{split}
&(n\!+\!1)\bigg(\!\!\big(k_{r_{\bu}}\!-\!1\!\big)\de_{\R}(\cC_0)\!+\!\de_{\R}(V_0,\vph_0)
+\!\!\!\!\!\!\!\sum_{\begin{subarray}{c}S_r^1\in\pi_0(\Si_0^{\si_0})\\ 
r>r_{\bu}\end{subarray}}\!\!\!\!\!\!\!\!\!\blr{w_1(V_0^{\vph_0}),[S_r^1]_{\Z_2}}
\!\!\bigg)+n\big(r(V_0,\vph_0)\!-\!1\big)\\
&\hspace{.7in}
=\begin{cases}
nk\!+\!(1\!-\!g_1\!+\!\deg V_1)((1\!-\!g_2)n\!+\!\deg V_2)\!+\!2\Z,
&\hbox{if}~\fp_0\!\in\!\cP^-_{\Si_0}(V_0^{\vph_0});\\
(n\!+\!1)(k\!+\!(1\!-\!g_1\!+\!\deg V_1)(1\!-\!g_2)\!)\!+\!2\Z,
&\hbox{if}~\fp_0\!\in\!\cP^+_{\Si_0}(V_0^{\vph_0}).
\end{cases}
\end{split}\end{equation*}
\end{enumerate}
\end{crl}

\vspace{.1in}

The CROrient~\ref{CRODisjUn_prop} property and Remark~\ref{CROdegenH3_rmk} similarly
yield a comparison between the limiting orientation $\fo_{\os_0}^-(V_0,\vph_0)$
(resp.~$\fo_{\cC_0;\fp_0}^-(V_0,\vph_0)$) with the split orientation 
$\fo_{\os_0}^{\sp}(V_0,\vph_0)$ (resp.~$\fo_{\cC_0;\fp_0}^{\sp}(V_0,\vph_0)$).
We make a special note of the comparisons of $\fo_{\cC_0;\fp_0}^{\pm}(V_0,\vph_0)$ 
with the split orientation in the case of $H3$ degenerations of~$S^2$.
In this case, we denote the numbers $k_{\bu1}$ and $k_{\bu2}$ 
of the CROrient~\ref{CROdegenH3_prop} property and Remark~\ref{CROdegenH3_rmk}
by~$k_1$ and~$k_2$, respectively.

\begin{crl}\label{CROdegenP1H3_crl}
Let $\cC_0$ be a decorated marked symmetric surface as in~\eref{cCdegC_e}
so that $(\Si_0,\si_0)$ consists of two copies of $(S^2,\tau)$ joined
at a single $H3$ node.
Suppose 
$(V_0,\vph_0)$ is a $\cC_0$-balanced rank~$n$ real bundle pair over~$(\Si_0,\si_0)$
and $\fp_0$ is a relative $\Pin^{\pm}$-structure on the real vector bundle~$V_0^{\vph_0}$
over $\Si_0^{\si_0}\!\subset\!\Si_0$.
\begin{enumerate}[label=(\alph*),leftmargin=*]

\item\label{P1HdegenPin_it1} 
The limiting orientation $\fo_{\cC_0;\fp_0}^+(V_0,\vph_0)$ on $\wt\la_{\cC_0}(D_{(V_0,\vph_0)})$ 
is the same as the split orientation $\fo_{\cC_0;\fp_0}^{\sp_0}(V_0,\vph_0)$
if and only~if 
\begin{equation*}\begin{split}
&(n\!+\!1)\!\big(\!(k\!-\!1)\de_{\R}(\cC_0)\!+\!\de_{\R}(V_0,\vph_0)
\!\big)+n\big(r(V_0,\vph_0)\!-\!1\big)\\
&\hspace{1.2in}
=\begin{cases}
nk\!+\!(1\!+\!\deg V_1)(n\!+\!\deg V_2)\!+\!2\Z,
&\hbox{if}~\fp_0\!\in\!\cP^-_{\Si_0}(V_0^{\vph_0});\\
(n\!+\!1)(k\!+\!1\!+\!\deg V_1)\!+\!2\Z,
&\hbox{if}~\fp_0\!\in\!\cP^+_{\Si_0}(V_0^{\vph_0}).
\end{cases}
\end{split}\end{equation*}

\item\label{P1Hdegen_it2} 
The limiting orientation $\fo_{\cC_0;\fp_0}^-(V_0,\vph_0)$ on $\wt\la_{\cC_0}(D_{(V_0,\vph_0)})$ 
is the same as the split orientation $\fo_{\cC_0;\fp_0}^{\sp}(V_0,\vph_0)$
if and only~if 
\begin{equation*}\begin{split}
&(n\!+\!1)\big(\!(k\!-\!1)\de_{\R}(\cC_0)\!+\!\de_{\R}(V_0,\vph_0)
\!\big)
\!+\!r(V_0,\vph_0)\!\bigg(\!\!k_2\!+\!n\binom{k_2}{2}\!\!\!\bigg)\\
&\hspace{.5in}
+\!\big(r(V_0,\vph_0)\!-\!1\big)n
\bigg(\!\!1\!+\!\binom{k_2\!-\!1}{1}\!+\!
\binom{k_2\!-\!1}{2}\!\!\!\bigg)
\!+\!\big(r(\cC_0)\!-\!1\big)(n\!+\!1)(k_1\!-\!1)(k_2\!-\!1) \\
&\hspace{.3in}
=\blr{w_2(\fp_0),[\Si_*]_{\Z_2}}\!+\!\binom{\deg V_2}{2}
\!+\!\begin{cases}
nk\!+\!(1\!+\!\deg V_1)(n\!+\!\deg V_2),
&\hbox{if}~\fp_0\!\in\!\cP^-_{\Si_0}(V_0^{\vph_0});\\
(n\!+\!1)(k\!+\!1\!+\!\deg V_1)\!+\!\deg V_2,
&\hbox{if}~\fp_0\!\in\!\cP^+_{\Si_0}(V_0^{\vph_0}).
\end{cases}
\end{split}\end{equation*}
\end{enumerate}
\end{crl}

\subsection{Orientations and evaluation isomorphisms}
\label{OrientEvalIsom_subs}

Let $\cC$ be a marked symmetric surface as in~\eref{cCsymdfn_e} and 
$(V,\vph)$ be a real bundle pair over~$(\Si,\si)$. 
Define
$$\ev_{\cC}\!: \Ga(\Si;V)^{\vph}\lra 
V_{\cC}\!\equiv\!\bigoplus_{i=1}^k V_{x_i}^{\vph} \oplus
\bigoplus_{i=1}^lV_{z_i^+},  ~~
\ev_{\cC}(\xi)=\big(\!\big(\xi(x_i)\!\big)_{\!i\in[k]},\big(\xi(z_i^+)\!\big)_{\!i\in[l]}\big).$$
We identify an orientation of $\wt\la_{\cC}(D_{(V,\vph)})$ with a homotopy class of isomorphisms
$$\la(D_{(V,\vph)})\lra \la(V_{\cC})$$
via the complex orientations of the fibers~$V_{z_i^+}$ of~$V$.
The conjugation $\vph$ on~$V$ induces a conjugation on the complex vector bundle
$$V(\!-\!\cC)\equiv V\otimes \bigotimes_{i=1}^k\cO_{\Si}(-x_i)
\otimes \bigotimes_{i=1}^l\cO_{\Si}(-z_i^+\!-\!z_i^-),$$
which we denote in the same way.
By \cite[Lemma~2.4.1]{Sh}, a real CR-operator~$D$ on $(V,\vph)$ induces 
a real CR-operator $D_{-\cC}$ on $(V(\!-\!\cC),\vph)$ so that the restriction sequence
\BE{Vrestses_e} 0\lra D_{-\cC} \lra D\xlra{\ev_{\cC}} V_{\cC} \lra0\EE
of Fredholm operators is exact.
We note the following.

\begin{lmm}\label{rk1eval_lmm}
Suppose $\cC$ is a marked symmetric surface as in~\eref{cCsymdfn_e} with 
the underlying symmetric surface~$(S^2,\tau)$ and
$D$ is a real CR-operator on rank~1 real bundle pair $(V,\vph)$ over~$(S^2,\tau)$.
\begin{enumerate}[label=(\arabic*),leftmargin=*]

\item\label{rk1eval_it1} The operator $D$ is surjective if $\deg V\!\ge\!-1$ and 
injective if $\deg V\!\le\!-1$.

\item\label{rk1eval_it2}  The homomorphism 
\BE{cOevdfn_e} \ev_{\cC}\!: \ker D\lra  V_{\cC}
\!\equiv\!\bigoplus_{i=1}^k V_{x_i}^{\vph} \oplus
\bigoplus_{i=1}^lV_{z_i^+},  ~~
\ev_{\cC}(\xi)=\big(\!\big(\xi(x_i)\!\big)_{\!i\in[k]},\big(\xi(z_i^+)\!\big)_{\!i\in[l]}\big),\EE
is surjective if $k\!+\!2l\!\le\!\deg V\!+\!1$ and is an isomorphism if 
$k\!+\!2l\!=\!\deg V\!+\!1$.

\end{enumerate}
\end{lmm}

\begin{proof} By the proof of \cite[Proposition~3.6]{XCapsSetup},
the dimensions of the kernel and cokernel of the real CR-operator~$D$ on $(V,\vph)$ 
are halves of the dimensions of the kernel and cokernel of the associated real CR-operator
on~$V$ (i.e.~without restricting to the subspaces of real sections).
By \cite[Theorem~$1'$]{HLS97}, the latter is surjective if $\deg V\!\ge\!-1$ and 
injective if $\deg V\!\le\!-1$.
This establishes~\ref{rk1eval_it1}.
Along with the exactness of the restriction sequence~\eref{Vrestses_e},
\ref{rk1eval_it1} in turn implies~\ref{rk1eval_it2}.
\end{proof}

If $V^{\vph}$ is orientable and $k\!\not\in\!2\Z$, let 
$$\fo_{\cC;0}^{\pm}(V,\vph)\equiv \fo_{\cC;\io_{S^2}(\fp_0^{\pm}(V^{\vph}))}(V,\vph)$$
denote the orientation on $\wt\la_{\cC}(D_{(V,\vph)})$ determined by the image
$$\io_{S^2}\big(\fp_0^{\pm}(V^{\vph})\!\big)\in\cP_{S^2}^{\pm}\big(V^{\vph}\big)$$
of  the canonical $\Pin^{\pm}$-structure $\fp_0^{\pm}(V^{\vph})$ on~$V^{\vph}$
as in Remark~\ref{PinpCan_lmm} under the second map in~\eref{vsSpinPin_e0} with $X\!=\!S^2$.
If $V^{\vph}$ is not orientable, 
$\bD^2_+$ is a half-surface of~$(S^2,\tau)$, and $k\!\in\!2\Z$,
let $\fo_{\cC;0}^{\pm}(V,\vph)$ be the orientation of $\wt\la_{\cC}(D_{(V,\vph)})$ 
as above the CROrient~\ref{CRONormal_prop} property.\\

A half-surface $\bD^2_+$ of~$(S^2,\tau)$ determines an orientation of the fixed locus~$S^1$ 
of~$(S^2,\tau)$
(the counterclockwise rotation if $\bD^2_+$ is identified with the unit disk 
around $0\!\in\!\C$).
We say that distinct points $x_1,\ldots,x_k\in\!S^1$ are \sf{ordered by position
with respect to~$\bD^2_+$} if every $x_i$ with $i\!\in\![k\!-\!1]$ 
directly precedes~$x_{i+1}$ as $S^1$ is traversed in the direction determined 
by its orientation as the boundary of~$\bD^2_+$.

\begin{crl}\label{NormEv_crl}
Suppose $\cC$ is a marked symmetric surface as in~\eref{cCsymdfn_e} and 
$(V,\vph)$ is a rank~1 real bundle pair over~$\cC$ so~that 
\BE{NormEvCrl_e}(\Si,\si)=(S^2,\tau) \qquad\hbox{and}\qquad
k\!+\!2l\!=\!\deg V\!+\!1.\EE
If $\bD^2_+$ is a half-surface of~$(S^2,\tau)$, $z_i^+\!\in\!\bD^2_+$ for every $i\!\in\![l]$, 
and the real marked points 
$x_1,\ldots,x_k\in\!S^1$ are ordered by position with respect to~$\bD^2_+$,
then the isomorphism~\eref{cOevdfn_e} lies in the homotopy class determined
by the orientation $\fo_{\cC;0}^{\pm}(V,\vph)$ of  $\wt\la_{\cC}(D_{(V,\vph)})$.
\end{crl}

\begin{proof} Let $a\!=\!\deg V$.
By~\eref{NormEvCrl_e}, $a\!\ge\!-1$.
The $(k,l)\!=\!(1,0),(2,0)$ cases of this claim are 
the two statements of the CROrient~\ref{CRONormal_prop} property.
We first show below that the claim in the $l\!=\!0$ cases implies the claim in all cases.
We then give two proofs of the $l\!=\!0$ cases.
The first proof makes use of Proposition~\ref{cOpincg_prp},
which in turn is used in the proof of Theorem~\ref{CROrient_thm};
the second proceeds directly from
the properties of the orientations of Theorem~\ref{CROrient_thm}
stated in Sections~\ref{OrientPrp_subs1} and~\ref{OrientPrp_subs2}.\\

(1) Let $\cC_0\!\equiv\!(\Si_0,\si_0)$ be the connected symmetric surface
which contains precisely one conjugate pair $(\nod^+,\nod^-)$ of nodes and
consists of three irreducible components, each of which is isomorphic to~$S^2$.
Let $\P^1_0\!\subset\!\Si_0$ be the component preserved by~$\si_0$ and
$\bD^2_+,\bD^2_-\!\subset\!\P^1_0$ be the two half-surfaces.
 We label the two conjugate nodes so that $\nod^+\!\in\!\bD^2_+$ and denote by 
$\P^1_+$ the conjugate component of~$\Si_0$ also containing~$\nod^+$;
see the first diagram in Figure~\ref{P1CH3split_fig}.
Let $(\cU,\wt\fc)$ be a flat family of deformations of~$\cC_0$ as in~\eref{cUsymmdfn_e}
over the unit ball $\De\!\subset\!\C^2$ around the origin so that 
\BE{3compcUdfn_e}\big(\Si_{\t_1},\si_{\t_1}\big)=\big(S^2,\tau\big)\EE
for some $\t_1\!\in\!\De_{\R}$.
The decorated structure on $(S^2,\tau)$ then induces a decorated structure 
on the fiber $(\Si_{\t},\si_{\t})$ of $(\cU,\wt\fc)$ over every $\t\!\in\!\De_{\R}$.\\

Let $(\wt{V},\wt\vph)$ be a rank~1 real bundle pair over $(\cU,\wt\fc)$ so~that 
\BE{pinchcond_e3a}\deg\wt{V}\big|_{\P^1_0}=k\!-\!1, \quad  \deg\wt{V}\big|_{\P^1_+}=l,
\quad \big(\wt{V}_{\t_1},\wt\vph_{\t_1}\big)=(V,\vph).\EE
Such a real bundle pair can be obtained by extending an appropriate real bundle pair from
the central fiber and deforming the extension so that it satisfies 
the last condition in~\eref{pinchcond_e3a}.
The $\Pin^{\pm}$-structure $\fp_0^{\pm}(V^{\vph})$ on~$V^{\vph}$ extends to 
a  $\Pin^{\pm}$-structure~$\wt\fp^{\pm}$ on~$\wt{V}^{\wt\vph}$ which restricts to 
$\fp_0^{\pm}(\wt{V}_{\t}^{\wt\vph_{\t}})$ for every $\t\!\in\!\De_{\R}$. 
We denote the restriction of $(\wt{V}_0,\wt\vph_0)$ to~$\P^1_0$ by $(\wt{V}_{00},\wt\vph_{00})$
and the restriction of~$\wt\fp^{\pm}$ to~$\wt{V}_0^{\wt\vph_0}$ by~$\wt\fp_0^{\pm}$.\\

\begin{figure}
\begin{pspicture}(-.3,-2)(10,2.3)
\psset{unit=.4cm}
\pscircle(8,0){2}\pscircle(8,3.5){1.5}\pscircle(8,-3.5){1.5}
\pscircle*(8,2){.2}\pscircle*(8,-2){.2}
\rput(8.1,2.7){\sm{$\nod^+$}}\rput(8,-2.7){\sm{$\nod^-$}}
\psline{<->}(4,4)(4,-4)\rput(3.2,0){$\si_0$}
\rput(10.2,4.5){$\P^1_+$}\rput(10.2,-4.5){$\P^1_-$}
\rput(5.8,1.8){$\bD^2_+$}\rput(5.8,-1.8){$\bD^2_-$}
\rput(10.8,0){$\P^1_0$}
\psarc[linewidth=.03](8,3.46){4}{240}{300} 
\psarc[linewidth=.03,linestyle=dashed](8,-3.46){4}{60}{120} 
\pscircle(27,0){3}\pscircle(33,0){3}
\pscircle*(30,0){.2}\rput(30,3){\sm{$\nod$}}
\psline[linewidth=.03]{->}(30,2.3)(30,.7)
\rput(27.2,3.8){$\P^1_1$}\rput(33.2,3.8){$\P^1_2$}
\rput(24.4,3){$\bD^2_{1+}$}\rput(24.5,-2.6){$\bD^2_{1-}$}
\rput(35.5,2.9){$\bD^2_{2+}$}\rput(35.7,-2.9){$\bD^2_{2-}$}
\psarc[linewidth=.03](27,5.2){6}{240}{300} 
\psarc[linewidth=.03,linestyle=dashed](27,-5.2){6}{60}{120} 
\psarc[linewidth=.03](33,5.2){6}{240}{300} 
\psarc[linewidth=.03,linestyle=dashed](33,-5.2){6}{60}{120} 
\rput(27,-.2){\sm{$S^1_1$}}\rput(33,-.2){\sm{$S^1_2$}}
\psline{<->}(23,4)(23,-4)\rput(22.2,0){$\si_0$}
\end{pspicture}
\caption{The three- and two-component symmetric surfaces~$\cC_0$ of
the proof of Corollary~\ref{NormEv_crl}}
\label{P1CH3split_fig}
\end{figure}

Let \hbox{$\cD\!\equiv\!\{D_{\t}\}$} be a family of real CR-operators on
$(\wt{V}_{\t},\wt\vph_{\t})$ as in~\eref{DVvphdfn_e} so~that $D_{\t_1}\!=\!D$.
By Lemma~\ref{rk1eval_lmm}\ref{rk1eval_it1}, each operator~$D_{\t}$ is then surjective; their kernels form 
a vector bundle
\BE{picKdfn_e}\pi_{\cK}\!:\cK\!(\cD)\!\equiv\! 
\bigsqcup_{\t\in\De_{\R}}\!\!
\{\t\}\!\times\!\big(\!\ker D_{\t}\big) \lra\De_{\R}.\EE
We denote by $D_{00}$ and $D_0^+$ the restrictions of $D_0$ to 
real CR-operators on $(\wt{V}_{00},\wt\vph_{00})$ and~$\wt{V}|_{\P^1_+}$, respectively. 
The exact triple 
\BE{Cndses_e}\begin{split}
&\qquad 0\lra D_0 \lra D_{00}\!\oplus\!D_0^+
\lra \wt{V}_{\nod^+}\lra 0, \\
&(\xi_-,\xi_0,\xi_+)\lra (\xi_0,\xi_+), \quad
 (\xi_0,\xi_+)\lra\xi_+(\nod^+)\!-\!\xi_0(\nod^+), 
\end{split}\EE
of Fredholm operators then determines an isomorphism
\BE{Cndses_e2}
\la(D_0)\!\otimes\!\la\big(\wt{V}_{\nod^+}\big)\approx 
\la\big(D_{00}\big)\!\otimes\!\la\big(D_0^+\big).\EE

\vspace{.15in}

Let $s_1^+,\ldots,s_l^+$ be disjoint sections of $\cU\!-\!\cU^{\wt\fc}$ and
$s_1^{\R},\ldots,s_k^{\R}$ be disjoint sections of $\cU^{\wt\fc}$ over~$\De_{\R}$ 
so~that
$$s_i^{\R}(\t_1)=x_i~~\forall\,i\!\in\![k], \qquad
s_i^+(0)\in\P^1_+\!-\!\big\{\nod^{\pm}\big\},~~s_i(\t_1)=z_i^+~~\forall\,i\!\in\![l].$$
By Lemma~\ref{rk1eval_lmm}\ref{rk1eval_it2}, the bundle homomorphism
\BE{evcUdfn_e}\begin{split}
\ev_{\cU}\!:\cK\!(\cD)\lra \wt{V}_{\cU}^{\R}\!\oplus\!\wt{V}_{\cU}^{\C}
\equiv\bigoplus_{i=1}^k s_i^{\R*}\wt{V}^{\wt\vph} \oplus
\bigoplus_{i=1}^ls_i^{+*}\wt{V},\\
\ev_{\cU}(\t,\xi)=\big(\t,\big(\xi(s_i^{\R}(\t))\!\big)_{\!i\in[k]},
\big(\xi(s_i^+(\t))\!\big)_{\!i\in[l]}\big),
\end{split}\EE
is then an isomorphism restricting to the isomorphism~$\ev_{\cC}$ in~\eref{cOevdfn_e} over~$\t_1$.
Thus, $\ev_{\cC}$ lies in the homotopy class determined
by the orientation $\fo_{\cC;0}^{\pm}(V,\vph)$ of  $\wt\la_{\cC}(D_{(V,\vph)})$ if
and only~if the restriction~$\ev_{\cC_0}$ of~$\ev_{\cU}$ over~0 lies in the homotopy class determined
by the limiting orientation $\fo_{\cC_0;\wt\fp_0^{\pm}}'(\wt{V}_0,\wt\vph_0)$ of 
$\wt\la_{\cC_0}(D_0)$.
By Corollary~\ref{CROdegenC_crl}, the last orientation is the same as 
the $\C$-split orientation $\fo_{\cC_0;0}^{\C\pm}(\wt{V}_0,\wt\vph_0)$
determined by the $\Pin^{\pm}$-structure 
$\fp_0^{\pm}(\wt{V}_{00}^{\wt\vph_{00}})\!=\!\fp_0^{\pm}(V^{\vph})$.\\

By the first two conditions in~\eref{pinchcond_e3a} and Lemma~\ref{rk1eval_lmm}\ref{rk1eval_it2}, 
the evaluation homomorphisms
\begin{alignat*}{2}
\ev_{\cC_0}^0\!:\ker D_{00}&\lra  \wt{V}_{\cU}^{\R}\big|_0,  &\qquad
\ev_{\cC_0}^0(\xi)&=\big(\big(\xi(s_i^{\R}(0))\!\big)_{\!i\in[k]}\big),\\
\ev_{\cC_0}^+\!:\ker D_0^+&\lra  \wt{V}_{\nod^+}\!\oplus\!\wt{V}_{\cU}^{\C}\big|_0, 
&\qquad
\ev_{\cC_0}^+(\xi)&=\big(\xi(\nod^+),\big(\xi(s_i^+(0))\!\big)_{\!i\in[l]}\big),
\end{alignat*}
are isomorphisms.
The last isomorphism determines the complex orientation of~$D_0^+$.
The $\C$-split orientation $\fo_{\cC_0;0}^{\C\pm}(\wt{V}_0,\wt\vph_0)$
is thus determined by the homotopy class of the isomorphism~$\ev_{\cC_0}^+$,
the orientation $\fo_{\cC_{00};0}^{\pm}(V_{00},\vph_{00})$ of~$D_{00}$,
and the complex orientation of~$\la(\wt{V}_{\nod^+})$ via the isomorphism~\eref{Cndses_e2}.
Along with the conclusion of the previous paragraph, this implies that
the isomorphism~$\ev_{\cC}$ in~\eref{cOevdfn_e} lies  in the homotopy class determined
by the orientation~$\fo_{\cC;0}^{\pm}(V,\vph)$ if and only if 
the isomorphism~$\ev_{\cC_0}^0$ lies  in the homotopy class determined
by the orientation~$\fo_{\cC_{00};0}^{\pm}(V_{00},\vph_{00})$.
Thus, the claim of the corollary holds for $(k,l)$ if and only if it holds for $(k,0)$.
In light of the CROrient~\ref{CRONormal_prop} property, this establishes
the claim if $k\!=\!1$ or $k\!=\!2$.\\

(2) By Proposition~\ref{cOpincg_prp}, the claim of the corollary holds for $(k,l)$ 
if and only if it holds for $(k\!+\!2l,0)$.
Along with the conclusion of the previous paragraph, this implies that the claim holds
if $k\!>\!0$ or $l\!>\!0$.
The $(k,l)\!=\!(0,0)$ case then follows from the $(k,l)\!=\!(0,1)$ case 
via the if and only if statement of the previous paragraph.
Below we obtain the same conclusions by using Corollary~\ref{CROdegenP1H3_crl} instead of 
Proposition~\ref{cOpincg_prp}.\\

We consider the behavior of the relevant orientations under flat degenerations of~$(S^2,\tau)$ 
to a connected symmetric surface $\cC_0\!\equiv\!(\Si_0,\si_0)$
with one node~$\nod$.
This node is then of type~$H3$, i.e.~it is a non-isolated point of~$\Si_0^{\si_0}$
separating~$\Si_0$ into two copies of~$(\P^1,\tau)$, which we denote by~$\P^1_1$ and~$\P^1_2$.
For $r\!=\!1,2$, let $\bD^2_{r\pm}\!\subset\!\P^1_r$ be the two distinguished half-surfaces and
$$S^1_r\equiv \bD^2_{r+}\!\cap\!\bD^2_{r-}\subset\P^1_r$$
be the fixed locus of the involution~$\si_0|_{\P^1_r}$;
see the second diagram in Figure~\ref{P1CH3split_fig}.
Let $\cU$ be a flat family of deformations of~$\cC_0$ as in~\eref{cUsymmdfn_e}
over the unit ball $\De\!\subset\!\C$ around the origin satisfying~\eref{3compcUdfn_e}.
We choose a decorated structure on $(\Si_0,\si_0)$ so that $\t_1\!\in\!\De_{\R}^*$
lies in the topological component $\De_{\R}^+\!\subset\!\De_{\R}^*$ 
defined as above the CROrient~\ref{CROdegenH3_prop} property.
This decorated structure then induces a decorated structure 
on the fiber $(\Si_{\t},\si_{\t})$ of $(\cU,\wt\fc)$ over every $\t\!\in\!\De_{\R}$.\\

Suppose $k\!\ge\!1$. 
Let $(\wt{V},\wt\vph)$ be a rank~1 real bundle pair over $(\cU,\wt\fc)$ so~that 
\BE{pinchcond_e3c}\deg\wt{V}\big|_{\P^1_1}=k\!+\!2l\!-\!2, \quad  
\deg\wt{V}\big|_{\P^1_2}=1,
\quad \big(\wt{V}_{\t_1},\wt\vph_{\t_1}\big)=(V,\vph)\,.\EE
We denote the restrictions of $(\wt{V}_0,\wt\vph_0)$  to
$\P^1_1$ and $\P^1_2$ by $(V_1,\vph_1)$ and $(V_2,\vph_2)$, respectively.
The $\Pin^{\pm}$-structure $\fp_0^{\pm}(V^{\vph})$ on~$V^{\vph}$ extends to 
a  $\Pin^{\pm}$-structure $\wt\fp^{\pm}$ on~$\wt{V}^{\wt\vph}$ which restricts~to
$\fp_0^{\pm}(\wt{V}_{\t}^{\wt\vph_{\t}})$ for every $\t\!\in\!\De_{\R}^+$
and to $\fp_2^{\pm}\!\equiv\!\fp_0^{\pm}(V_2^{\vph_2})$ on~$V_2^{\vph_2}$.
Then,
\BE{pinchcond_e5a}
\fp_1^-\!\equiv\!\wt\fp^-\big|_{V_1^{\vph_1}}=
\begin{cases}
\fp_0^-(V_1^{\vph_1}),&\hbox{if}~\deg V\!\not\in\!2\Z;\\
\fp_1^-(V_1^{\vph_1}),&\hbox{if}~\deg V\!\in\!2\Z;
\end{cases}
\qquad
\fp_1^+\!\equiv\!\wt\fp^+\big|_{V_1^{\vph_1}}=\fp_0^+(V_1^{\vph_1})\,.\EE
We denote the restriction of~$\wt\fp^{\pm}$ to~$\wt{V}_0^{\wt\vph_0}$ by~$\wt\fp_0^{\pm}$.\\

Let \hbox{$\cD\!\equiv\!\{D_{\t}\}$} be a family of real CR-operators on
$(\wt{V}_{\t},\wt\vph_{\t})$ as in~\eref{DVvphdfn_e} so~that $D_{\t_1}\!=\!D$.
By  Lemma~\ref{rk1eval_lmm}\ref{rk1eval_it1}, each operator~$D_{\t}$ is again surjective; 
their kernels form  a vector bundle $\cK(\cD)$ as before.
We denote by $D_1$ and $D_2$ the restrictions of $D_0$ to 
real CR-operators on $(V_1,\vph_1)$ and~$(V_2,\vph_2)$, respectively. 
The exact~triple 
\BE{Rndses_e}
0\lra D_0 \lra D_1\!\oplus\!D_2 \lra V_{\nod}^{\vph}\lra 0, \quad
(\xi_1,\xi_2)\lra\xi_2(\nod)\!-\!\xi_1(\nod),\EE
of Fredholm operators then determines an isomorphism
\BE{Rndses_e2} \la(D_0)\!\otimes\!\la\big(V_{\nod}^{\vph}\big)
\approx\la(D_1)\!\otimes\!\la(D_2).\EE

\vspace{.15in}

Let $s_1^+,\ldots,s_l^+$ be disjoint sections of $\cU\!-\!\cU^{\wt\fc}$ and
$s_1^{\R},\ldots,s_k^{\R}$ be disjoint sections of $\cU^{\wt\fc}$ over~$\De_{\R}$ 
so~that
$$s_i^{\R}(0)\in S^1_1\!-\!\{\nod\}~~\forall\,i\!\in\![k\!-\!1],\quad 
s_k^{\R}(0)\in S^1_2\!-\!\{\nod\}, \quad
s_i^{\R}(\t_1)=x_i~~\forall\,i\!\in\![k],\quad 
s_i(\t_1)=z_i^+~~\forall\,i\!\in\![l],$$
and the points $\nod,x_1,\ldots,x_{k-1}$ of $S^1_1$ are ordered by position 
with respect to~$\bD^2_{1+}$.
By Lemma~\ref{rk1eval_lmm}\ref{rk1eval_it2},
the bundle homomorphism $\ev_{\cU}$ as in~\eref{evcUdfn_e} is again an isomorphism
restricting to the isomorphism~$\ev_{\cC}$ in~\eref{cOevdfn_e} over~$\t_1$.
Thus, $\ev_{\cC}$ lies in the homotopy class determined
by the orientation $\fo_{\cC;0}^{\pm}(V,\vph)$ of  $\wt\la_{\cC}(D_{(V,\vph)})$ if
and only~if the restriction~$\ev_{\cC_0}$ of~$\ev_{\cU}$ over~0 lies in the homotopy class determined
by the limiting orientation $\fo_{\cC_0;\wt\fp_0^{\pm}}^+(\wt{V}_0,\wt\vph_0)$ of $\wt\la_{\cC_0}(D_0)$.
By Corollary~\ref{CROdegenP1H3_crl}\ref{P1HdegenPin_it1}, 
\BE{pinchcond_e5b}
\fo_{\cC_0;\wt\fp_0^-}^+(\wt{V}_0,\wt\vph_0)
\begin{cases}=\fo_{\cC_0;\wt\fp_0^-}^{\sp}(\wt{V}_0,\wt\vph_0),
&\hbox{if}~\deg V\!\not\in\!2\Z;\\
\neq \fo_{\cC_0;\wt\fp_0^-}^{\sp}(\wt{V}_0,\wt\vph_0),
&\hbox{if}~\deg V\!\in\!2\Z; \end{cases}\quad
\fo_{\cC_0;\wt\fp_0^+}^+(\wt{V}_0,\wt\vph_0)
=\fo_{\cC_0;\wt\fp_0^+}(\wt{V}_0,\wt\vph_0),\EE
where $\fo_{\cC_0;\wt\fp_0^{\pm}}(\wt{V}_0,\wt\vph_0)$ is the split orientation
determined by the $\Pin^{\pm}$-structures $\fp_1^{\pm}$ on~$V_1^{\vph_1}$ 
and~$\fp_2^{\pm}$ on~$V_2^{\vph_2}$.\\

Let $\wt\cC_{11}'$ denote the irreducible component $\P^1_1$  of $\Si_0$
with its involution and the marked points ($k\!-\!1$ real points and $l$ conjugate pairs
of points).
Let $\wt\cC_{12}'$ denote the irreducible component $\P^1_2$  of $\Si_0$
with its involution, the marked point $x_k\!\in\!S^1_2$,
and the nodal point~$\nod$.
By the first two conditions in~\eref{pinchcond_e3c} and Lemma~\ref{rk1eval_lmm}\ref{rk1eval_it2}, 
the evaluation homomorphisms
\begin{alignat*}{2}
\ev_{\wt\cC_{11}'}\!:\ker D_1&\lra 
\bigoplus_{i=1}^{k-1}\!\wt{V}_{s_i^{\R}(0)}^{\wt\vph}
\!\oplus\!\wt{V}_{\cU}^{\C}, &\quad
\ev_{\wt\cC_{11}'}(\xi)&=\big(\!\big(\xi(s_i^{\R}(0))\!\big)_{\!i\in[k-1]},
\big(\xi(s_i^+(0))\!\big)_{\!i\in[l]}\big),\\
\ev_{\wt\cC_{12}'}\!:\ker D_2&\lra 
\wt{V}_{s_k^{\R}(0)}^{\wt\vph}\!\oplus\!\wt{V}_{\nod}^{\wt\vph}, &\quad
\ev_{\wt\cC_{12}'}(\xi)&=\big(\xi(s_k^{\R}(0)),\xi(\nod^+)\!\big),
\end{alignat*}
are isomorphisms.
By the CROrient~\ref{CRONormal_prop}\ref{CROnormPin_it} property, 
the second isomorphism above lies in the homotopy class determined
by the orientation 
$$\fo_{\wt\cC_{12}';\io_{\P_2^1}(\fp_2^{\pm})}(V_2,\vph_2)=\fo_{\wt\cC_{12}';0}^{\pm}(V_2,\vph_2)$$
on $\wt\la_{\wt\cC_{12}'}(D_2)$.\\

The split orientation $\fo_{\cC_0;\wt\fp_0^{\pm}}^{\sp}(\wt{V}_0,\wt\vph_0)$ is thus determined 
by the homotopy class of the isomorphism~$\ev_{\wt\cC_{12}'}$ and
by the orientation
$\fo_{\wt\cC_{11}';\io_{\P_1^1}(\fp_1^{\pm})}(V_1,\vph_1)$ on 
$\wt\la_{\wt\cC_{11}'}(D_1)$ via the isomorphism~\eref{Rndses_e2}.
By~\eref{pinchcond_e5a} and the CROrient~\ref{CROSpinPinStr_prop}\ref{CROPinStr_it} property, 
$$\fo_{\wt\cC_{11}';\io_{\P_1^1}(\fp_1^-)}(V_1,\vph_1)
\begin{cases}
=\fo_{\wt\cC_{11}';0}^-(V_1,\vph_1),&\hbox{if}~\deg V\!\not\in\!2\Z;\\
\neq\fo_{\wt\cC_{11}';0}^-(V_1,\vph_1),&\hbox{if}~\deg V\!\in\!2\Z;
\end{cases}
\quad
\fo_{\wt\cC_{11}';\io_{\P_1^1}(\fp_1^+)}(V_1,\vph_1)
=\fo_{\wt\cC_{11}';0}^+(V_1,\vph_1)\,.$$
Along with~\eref{pinchcond_e5b} and the conclusion of the previous paragraph, 
this implies that the isomorphism~$\ev_{\cC}$ in~\eref{cOevdfn_e} lies 
in the homotopy class determined by the orientation~$\fo_{\cC;0}^{\pm}(V,\vph)$ 
if and only if 
the isomorphism~$\ev_{\wt\cC_{11}'}$ lies  in the homotopy class determined
by the orientation~$\fo_{\wt\cC_{11}';0}^{\pm}(V_1,\vph_1)$.\\

Thus, the claim of the corollary holds for $(k,l)$ if and only if it holds for
$(k\!-\!1,l)$. 
Along with the already obtained conclusion for the $k\!=\!1,2$ cases of the claim,
this establishes the full statement of the claim.
\end{proof}

\begin{eg}\label{CRORientTP1_eg}
Let $D^2_+\!\subset\!S^2$ be the disk cut out by the fixed locus $S^1$ 
of the involution~$\tau$ on~$S^2$ which contains the origin~0 in $\C\!\subset\!S^2$
and $\cC$ be the marked symmetric surface consisting of~$S^2$ with one
real marked point.
We denote by $\prt_r$ the usual unit outward radial vector field on~$\C^*$ and
by $\fo_{S^1}$ the usual counterclockwise orientation on~$TS^1$.
Let
$\bp$ be the real CR-operator on the rank~1 real bundle $(T\P^1,\nd\tau)$ over~$(\P^1,\tau)$
determined by the holomorphic $\bp$-operator on~$T\P^1$ and
$\fo_0(T\P^1,\nd\tau;\fo_{S^1})$ be the orientation of 
$$\la(\bp)=\la(\ker\bp)$$
induced by $\fo_{S^1}$ as above the CROrient~\ref{CRONormal_prop} property
in Section~\ref{OrientPrp_subs1}.
By the CROrient~\ref{CROPin2SpinRed_prop} property,
the orientation $\fo_{\cC;0}^{\pm}(T\P^1,\nd\tau)$ on $\wt\la_{\cC}(\bp)$
corresponds to the homotopy class of isomorphisms of 
$\la(\bp)$ and $T_{x_1}S^1$ determined by the orientations $\fo_0(T\P^1,\nd\tau;\fo_{S^1})$  
and~$\fo_{S^1}$, 
respectively.
Along with Corollary~\ref{NormEv_crl}, this implies that the isomorphism
$$\ker \bp\lra T_0\P^1\!\oplus\! T_{x_1}S^1, \qquad
\xi\lra\big(\xi(0),\xi(x_1)\!\big),$$
respects the orientation $\fo_0(T\P^1,\nd\tau;\fo_{S^1})$  on $\ker\bp$, 
the complex orientation on~$T_0\P^1$, and
the counterclockwise orientation~$\fo_{S^1}$ on~$T_{x_1}S^1$.
Thus, the homomorphism
$$(\ker\bp)\!\oplus\!\R\lra T_0\P^1\!\oplus\!T_z\P^1, \qquad
(\xi,t)\lra \big(\xi(0),\xi(z)\!-\!t\prt_r\big),$$
is an isomorphism which respects the orientation $\fo_0(T\P^1,\nd\tau;\fo_{S^1})$ on $\ker\bp$,
the standard orientation on~$\R$, 
and the complex orientations on $T_0\P^1$ and $T_z\P^1$ for every $z\!\in\!\C^*$.
\end{eg}

\section{Base cases}
\label{BaseCR_sec}

In Section~\ref{tauorient_subs}, we follow the perspective in \cite[Section~3]{Ge1}
to describe orientations $\fo_{k,l}(V,\vph;\fo_{x_1})$ on the determinants of
real CR-operators on rank~1 real bundle pairs $(V,\vph)$ over~$(S^2,\tau)$,
with $k,l\!\in\!\Z^{\ge0}$ and \hbox{$k\!\not\cong\!\deg V$} mod~2.
We establish a number of properties of these orientations in Section~\ref{tauorient_subs}
and study their behavior 
under degenerations of $(S^2,\tau)$ to nodal symmetric surfaces 
in Sections~\ref{tauorient_subs2a} and~\ref{tauorient_subs2b}. 
The main results of  Sections~\ref{tauorient_subs2a} and~\ref{tauorient_subs2b},
Propositions~\ref{rk1Cdegen_prp} and~\ref{cOlim_prp},
are the analogues of the CROrient~\ref{CROdegenH3_prop} and~\ref{CROdegenC_prop} properties 
for the orientations $\fo_{k,l}(V,\vph;\fo_{x_1})$ of Section~\ref{tauorient_subs}. 
They imply in particular that the orientations $\fo_{k,l}(V,\vph;\fo_{x_1})$ do not depend
on the admissible choices of~$(k,l)$;  see Proposition~\ref{cOpincg_prp}\ref{cOpincg_it3}.
Outside of Sections~\ref{tauorient_subs}-\ref{tauorient_subs2b}, we thus denote these
orientations simply by~$\fo(V,\vph;\fo_{x_1})$.
In particular, the statements of Propositions~\ref{cOpincg_prp}\ref{cOpincg_it1},
\ref{rk1Cdegen_prp}, and~\ref{cOlim_prp} apply to the orientations~$\fo(V,\vph;\fo_{x_1})$.\\

In the first part of Section~\ref{tauspinorient_subs}, 
we recall the construction of \cite[Section~6.5]{Melissa} for orienting the determinants 
of real CR-operators on even-degree real bundle pairs~$(V,\vph)$ over~$(S^2,\tau)$
from relative $\OSpin$-structures~$\os$ on the real vector bundles~$V^{\vph}$ 
over the $\tau$-fixed locus $S^1\!\subset\!S^2$.
In the remainder of this section and in Section~\ref{tauspinorient_subs2}, 
we confirm that the resulting orientations~$\fo_{\os}(V,\vph)$ satisfy the applicable properties 
of Section~\ref{OrientPrp_subs1}.
While this can be done through geometric arguments similar to those
in Sections~\ref{tauorient_subs}-\ref{tauorient_subs2b},
we instead deduce some of these properties from the analogous properties
of the orientations of  Section~\ref{tauorient_subs} via 
the CROrient~\ref{CROSpinPinSES_prop}\ref{CROsesSpin_it} property 
for even-degree real bundle pairs~$(V,\vph)$ 
over~$(S^2,\tau)$, which we establish directly in Section~\ref{tauspinorient_subs2}.
We use the perspective of Definition~\ref{RelPinSpin_dfn3} on relative $\OSpin$-structures
throughout Sections~\ref{tauspinorient_subs} and~\ref{tauspinorient_subs2}.

\subsection{Line bundles over $(S^2,\tau)$: construction and properties}
\label{tauorient_subs}

Suppose $\cC$ is a marked symmetric surface as in~\eref{cCsymdfn_e} so that 
$(\Si,\si)$ is $S^2$ with the involution~$\tau$ and
$D$ is a real CR-operator on a rank~1 real bundle pair~$(V,\vph)$
over~$(S^2,\tau)$.
By Lemma~\ref{rk1eval_lmm},  the operator~$D$ is surjective and
the evaluation homomorphism in~\eref{cOevdfn_e} is an isomorphism~if 
\BE{klcond_e}  k\!+\!2l=\deg V\!+\!1\,.\EE
If $k\!=\!0$, the target $V_{\cC}$ of this isomorphism is oriented by the complex orientation of~$V$.
If $k\!>\!0$, an orientation~$\fo_{x_1}$ on $V_{x_1}$ determines an orientation
on $V_{\cC}$ by transporting~$\fo_{x_1}$
to each $x_i$ with $2\!\le\!i\!\le\!k$ along the positive direction of~$S^1$
with respect to the disk \hbox{$\bD^2_+\!\subset\!S^2$}.
Along with the complex orientations of the summands~$V_{z_i^+}$,
$\fo_{x_1}$ thus determines an orientation on~$V_{\cC}$.
If~\eref{klcond_e} holds, the resulting orientation on~$V_{\cC}$ determines 
an orientation~$\fo_{\cC}(\fo_{x_1})$ on the determinant
$$\la(D)\approx\la\big(\ker D\big)$$
of the real CR-operator~$D$ via the isomorphism~\eref{cOevdfn_e}.\\

We denote by $\fo_{k,l}(\fo_{x_1})$ the orientation $\fo_{\cC}(\fo_{x_1})$  
for a marked symmetric surface~$\cC$ as in~\eref{cCsymdfn_e} so that its real
points $x_1,\ldots,x_k$ are ordered by position with respect to~$\bD^2_+$
and $z_i^+\!\in\!\bD^2_+$ for every $i\!\in\![l]$.
Since the space of such~$\cC$ is path-connected,
the orientation~$\fo_{k,l}(\fo_{x_1})$ does not depend on the choice of~$\cC$ in~it.
We call $\fo_{k,l}(\fo_{x_1})$ the \sf{$(k,l)$-evaluation orientation} of~$D$.
It is straightforward to see that this orientation satisfies the following properties.

\begin{lmm}\label{cOev_lmm}
Suppose $k,l\!\in\!\Z^{\ge0}$ and
$(V,\vph)$ is a rank~1 real bundle pair  over~$(S^2,\tau)$
so that~\eref{klcond_e} holds.  
Let $\fo_{x_1}$ be an orientation of~$V^{\vph}$ at a point $x_1$ in $S^1\!\subset\!S^2$.
\begin{enumerate}[label=(\arabic*),leftmargin=*]

\item\label{cOevorient_it} 
If $\ov\fo_{x_1}$ is the orientation of~$V^{\vph}_{x_1}$ opposite to~$\fo_{x_1}$, then
$$\fo_{k,l}(\ov{\fo_{x_1}})=\begin{cases} 
\fo_{k,l}(\fo_{x_1}),&\hbox{if}~\deg V\!\not\in\!2\Z;\\
\ov{\fo_{k,l}(\fo_{x_1})},&\hbox{if}~\deg V\!\in\!2\Z.\end{cases}$$

\item\label{cOevorient_it2}  
The orientation $\fo_{k,l}(\fo_{x_1})$ does not depend on the choice 
of the half-surface $\bD^2_+\!\subset\!S^2$ if and only if
\hbox{$\deg V\!\cong\!0,3$} mod~4.

\end{enumerate}
\end{lmm}

\vspace{.1in}

The orientation $\fo_{k,l}(\fo_{x_1})$ induces an orientation $\fo_{k,l}(V,\vph;\fo_{x_1})$ 
of a real CR-operator~$D$ on any rank~1 real bundle pair 
$(V,\vph)$ with 
\BE{pinchcond_e} k\!+\!2l \cong \deg V\!+\!1  \qquad\hbox{mod}~~2\EE
as follows.
Let $\cC_0$, $(\nod^+,\nod^-)$, $\P^1_0,\P^1_{\pm}$, $\bD^2_{\pm}$, 
$(\cU,\wt\fc)$, and~$\t_1$ be as above~\eref{3compcUdfn_e},
$$a_0=k\!+\!2l\!-\!1, \qquad\hbox{and}\qquad a'=(\deg V\!-\!a_0)/2.$$
By~\eref{pinchcond_e}, $a'\!\in\!\Z$.
Let $(\wt{V},\wt\vph)$ be a rank~1 real bundle pair over $(\cU,\wt\fc)$
so that 
\BE{pinchcond_e3}\deg\wt{V}\big|_{\P^1_0}=a_0,
\quad  \deg\wt{V}\big|_{\P^1_+}=a',
\quad \big(\wt{V}_{\t_1},\wt\vph_{\t_1}\big)=(V,\vph).\EE
For an orientation~$\fo_{x_1}$ of $V^{\vph}$ at 
a point~$x_1$ in the fixed locus~$S^1\!\subset\!\Si_{\t_1}$ and
a point~$x_1'$ in the fixed locus~$S^1\!\subset\!\Si_0$, 
let $\fo_{x_1'}$ denote the orientation of $\wt{V}^{\wt\vph}$ at~$x_1'$ 
obtained by transferring~$\fo_{x_1}$ along a path in~$\cU^{\wt\fc}$ from~$x_1$ to~$x_1'$.
If the degree of~$V$ is even (i.e.~the real line bundle~$V^{\vph}$ over $S^1$ is orientable), 
then~$\fo_{x_1'}$ does not depend on the choice of this path.\\ 

Let \hbox{$\cD\!\equiv\!\{D_{\t}\}$} be a family of real CR-operators on
$(\wt{V}_{\t},\wt\vph_{\t})$ as in~\eref{DVvphdfn_e} so~that $D_{\t_1}\!=\!D$.
As above~\eref{Cndses_e}, we denote the restriction of $(\wt{V}_0,\wt\vph_0)$ 
to~$\P^1_0$ by $(\wt{V}_{00},\wt\vph_{00})$ and the restrictions of $D_0$ to 
real CR-operators on $(\wt{V}_{00},\wt\vph_{00})$ and~$\wt{V}|_{\P^1_+}$
by $D_{00}$ and $D_0^+$, respectively. 
The orientation $\fo_{k,l}(\fo_{x_1'})$ of~$D_{00}$
and the complex orientations of~$D_0^+$ and $\wt{V}_{\nod^+}$
determine an orientation of~$D_0$ via the isomorphism~\eref{Cndses_e2}  and 
thus an orientation of the line bundle~$\la(\cD)$ over~$\De_{\R}$.
The latter restricts to an orientation $\fo_{k,l}(V,\vph;\fo_{x_1})$ of~$\la(D)$,
which we call the \sf{$(k,l)$-intrinsic orientation of~$D$}.
By Lemma~\ref{cOev_lmm}\ref{cOevorient_it}, $\fo_{k,l}(V,\vph;\fo_{x_1})$ does not depend
on the choice of the path from~$x_1$ to~$x_1'$ above even if the degree of~$V$ is~odd.
In light of Proposition~\ref{cOpincg_prp}\ref{cOpincg_it3} below, 
which is established in Section~\ref{tauorient_subs2b}, 
we denote the orientations $\fo_{k,l}(V,\vph;\fo_{x_1})$ by $\fo(V,\vph;\fo_{x_1})$
after Section~\ref{tauorient_subs2b} and call $\fo(V,\vph;\fo_{x_1})$
the \sf{intrinsic orientation of~$D$}.

\begin{prp}\label{cOpincg_prp}
Suppose $(V,\vph)$ is a rank~1 real bundle pair  over~$(S^2,\tau)$ and
$D$ is a real CR-operator on~$(V,\vph)$.
Let $\fo_{x_1}$ be an orientation of~$V^{\vph}$ at a point $x_1$ in $S^1\!\subset\!S^2$.
\begin{enumerate}[label=(\arabic*),leftmargin=*]

\item\label{cOpincg_it1} For all $k,l\!\in\!\Z^{\ge0}$ satisfying~\eref{pinchcond_e},
the orientation $\fo_{k,l}(V,\vph;\fo_{x_1})$ of~$D$
satisfies the two properties of Lemma~\ref{cOev_lmm}.

\item\label{cOpincg_it2} 
For all $k,l\!\in\!\Z^{\ge0}$ satisfying~\eref{klcond_e},
the $(k,l)$-evaluation and $(k,l)$-intrinsic orientations  $\fo_{k,l}(\fo_{x_1})$
and $\fo_{k,l}(V,\vph;\fo_{x_1})$ of~$D$ are the~same. 

\item\label{cOpincg_it3} The orientation $\fo_{k,l}(V,\vph;\fo_{x_1})$ 
of~$D$ does not depend on
the choice of $k,l\!\in\!\Z^{\ge0}$ satisfying~\eref{pinchcond_e}.

\end{enumerate}
\end{prp}

\begin{proof}[{\bf{\emph{Proof of Proposition~\ref{cOpincg_prp}\ref{cOpincg_it1}}}}]
The property of Lemma~\ref{cOev_lmm}\ref{cOevorient_it} for the orientation $\fo_{k,l}(V,\vph;\fo_{x_1})$
follows from this property for the orientation~$\fo_{k,l}(\fo_{x_1})$ because
the choice of~$\fo_{x_1}$ affects the orientation only on the first term on the right-hand side
of~\eref{Cndses_e2} and the parities of the degrees of~$V$ and $\wt{V}\big|_{\P^1_0}$
are the same.
The change in the choice of the half-surface of~$(S^2,\tau)$ acts by 
the complex conjugation on the complex orientations of $\wt{V}_{\nod^+}$ and~$D_0^+$.
The complex dimension of~$\wt{V}_{\nod^+}$ is~1;
the complex index of $D_0^+$ is $a'\!+\!1$.
Combining this with Lemma~\ref{cOev_lmm}\ref{cOevorient_it2},
we find that the orientation $\fo_{k,l}(V,\vph;\fo_{x_1})$ does not depend on the choice 
of the half-surface $\bD^2_+\!\subset\!S^2$ if and only if 
the number
\BE{cOpincg1_e3}-1+(a'\!+\!1)+\frac{a_0(a_0\!+\!1)}2=
\frac{(\deg V)(\deg V\!+\!1)}{2}-2\big(a_0a'\!+\!a'^2)\EE
is even.
Thus, the orientation $\fo_{k,l}(V,\vph;\fo_{x_1})$ satisfies
the property of Lemma~\ref{cOev_lmm}\ref{cOevorient_it2}.
\end{proof}

Proposition~\ref{cOpincg_prp}\ref{cOpincg_it2} is a special case of  
Lemma~\ref{cOpincg_lmm}\ref{cOpincgl_it1} below.
We establish Proposition~\ref{cOpincg_prp}\ref{cOpincg_it3} in Section~\ref{tauorient_subs2b}
based on the following observation about orientations of the determinants of real CR-operators
on degree~2 rank~1 real bundle pairs over~$(S^2,\tau)$.

\begin{lmm}\label{cOpincg_lmm2}
With the notation as in Lemma~\ref{cOev_lmm}, $\fo_{3,0}(\fo_{x_1})\!=\!\fo_{1,1}(\fo_{x_1})$.
\end{lmm}

\begin{proof}
By \cite[Proposition~2.1]{RBP}, every rank~1 degree~$a$ real bundle pair $(V,\vph)$ 
over~$(\P^1,\tau)$ is isomorphic to the holomorphic line bundle $\cO_{\P^1}(a)$ 
with the natural lift~$\wt\tau_a$ of~$\tau$.
Thus, we can assume~that $(V,\vph)$ is  $(\cO_{\P^1}(2),\wt\tau_2)$, 
$D$ is the holomorphic $\bp$-operator, and 
$$x_1=1, \qquad x_2=\fI, \qquad x_3=-1, \qquad z_1^+=0.$$
The evaluation homomorphisms 
\BE{cOpincg_e3}\ker D\lra V^{\vph}_{x_1}\!\oplus\!V^{\vph}_{x_2}\!\oplus\!V^{\vph}_{x_3}
\qquad\hbox{and}\qquad
\ker D\lra V^{\vph}_{x_1}\!\oplus\!V_{z_1^+}\EE
are the isomorphisms determining the orientations 
$\fo_{3,0}(\fo_{x_1})$ and $\fo_{1,1}(\fo_{x_1})$, respectively.
An orientation~$\fo_{x_1}$ of~$V^{\vph}_{x_1}$ determines an orientation~$\fo_1$ of
the target of the first isomorphism;
combined with the complex orientation of~$V_{z_1^+}$, it also determines 
an orientation~$\fo_2$ of the target of the second isomorphism.
By definition,
the first (resp.~second) isomorphism above is orientation-preserving with respect to 
the orientation~$\fo_{3,0}(\fo_{x_1})$ (resp.~$\fo_{1,1}(\fo_{x_1})$) on its domain
and the orientation~$\fo_1$ (resp.~$\fo_2$) on its target.\\

\noindent
Choose a trivialization of $V|_{\bD^2_+}$ so that $V^{\vph}|_{S^1}$ is identified
with the real line subbundle
$$\La\equiv\big\{(z,v)\!\in\!S^1\!\times\!\C\!:v/z\!\in\!\R\big\}$$
of $S^1\!\times\!\C$; see \cite[Corollary~C.3.9]{MS}.
We can assume that the orientation~$\fo_{x_1}$ of~$V^{\vph}_{x_1}$ agrees with 
the standard orientation of $\R\!\subset\!\C$ under this identification.
The~sets
\BE{cOpincg_e7}
(1,0,0),(0,\fI,0),(0,0,-1)\in V^{\vph}_{x_1}\!\oplus\!V^{\vph}_{x_2}\!\oplus\!V^{\vph}_{x_3}
\quad\hbox{and}\quad
(1,0),(0,1),(0,\fI)\in V^{\vph}_{x_1}\!\oplus\!V_{z_1^+}\EE
are then oriented bases with respect to the orientations~$\fo_1$ and~$\fo_2$.\\

\noindent
The functions $z,1\!+\!z^2,\fI\!-\!\fI z^2$ form a basis for $\ker D$;
see Step~2 in the proof of \cite[Theorem~C.4.1]{MS}.
With respect to this basis and the bases~\eref{cOpincg_e7},
the isomorphisms~\eref{cOpincg_e3} are given by the matrices
$$\left(\!\!\begin{array}{ccc}1 & 2 & 0 \\ 1 & 0 & 2 \\ 1 & -2 & 0 \end{array}\!\!\right)
\qquad\hbox{and}\qquad
\left(\!\!\begin{array}{ccc} 1 & 2 & 0 \\ 0 & 1 & 0 \\ 0 & 0 & 1\end{array}\!\!\right),$$
respectively.
Since these matrices have determinants of the same sign, 
the orientation~$\fo_{3,0}(\fo_{x_1})$ on $\ker D$ corresponding to~$\fo_1$ and
the orientation~$\fo_{1,1}(\fo_{x_1})$ corresponding to~$\fo_2$ are the same.
\end{proof}

\begin{lmm}\label{cOpincg_lmm}
Suppose $(V,\vph)$ is a rank~1  real bundle pair  over~$(S^2,\tau)$,
$D$ is a real CR-operator on~$(V,\vph)$, and $k,l\!\in\!\Z^{\ge0}$ satisfy~\eref{pinchcond_e}.
Let $\fo_{x_1}$ be an orientation of~$V^{\vph}$ at a point $x_1$ in $S^1\!\subset\!S^2$.
\begin{enumerate}[label=(\arabic*),leftmargin=*]

\item\label{cOpincgl_it1} If $k\!+\!2l\!\le\!\deg V\!+\!1$, then the orientations 
$\fo_{k,l}(V,\vph;\fo_{x_1})$ and $\fo_{k,(\deg V+1-k)/2}(\fo_{x_1})$ 
are the same.

\item\label{cOorientcomp_it} The  $(k,l)$-intrinsic orientation
 $\fo_{k,l}(V,\vph;\fo_{x_1})$  does not depend
on the choice of $l\!\in\!\Z^{\ge0}$.

\end{enumerate}
\end{lmm}

\begin{proof}[{\bf{\emph{Proof of Lemma~\ref{cOpincg_lmm}\ref{cOpincgl_it1}}}}]
We continue with the notation above Proposition~\ref{cOpincg_prp}, taking~$\cC$
as in the definition of the $(k,l)$-evaluation orientation above Lemma~\ref{cOev_lmm}.
In this case, $a'\!\ge\!0$ and each operator~$D_{\t}$ is surjective.
Their kernels form a vector bundle~$\cK(\cD)$ as in~\eref{picKdfn_e}.
Let $s_1^+,\ldots,s_{l+a'}^+$ be disjoint sections of $\cU\!-\!\cU^{\wt\fc}$ and
$s_1^{\R},\ldots,s_k^{\R}$ be disjoint sections of $\cU^{\wt\fc}$ over~$\De_{\R}$ 
so~that
\begin{alignat*}{4}
s_i^{\R}(\t_1)&=x_i&~~&\forall\,i\!\in\![k],&\qquad
s_i^+(\t_1)&=z_i^+&~~&\forall\,i\!\in\![l\!+\!a'],\\
s_i^+(0)&\in\P^1_0\!-\!\big\{\nod^{\pm}\big\}&~~&\forall\,i\!\in\![l],  &\qquad 
s_{l+i}^+(0)&\in\P^1_+-\!\big\{\nod^+\big\}&~~&\forall\,i\!\in\![a'].
\end{alignat*}
By the assumption that $a'\!\ge\!0$, the evaluation homomorphisms
\BE{cOpincg_e2}\begin{aligned}
\ev_{\cC_0}^0\!:\ker D_{00}&\lra 
\bigoplus_{i=1}^k\wt{V}^{\wt\vph}_{s_i^{\R}(0)}\!\oplus\!\bigoplus_{i=1}^l\!\wt{V}_{s_i^+(0)}, 
&
\ev_{\cC_0}^0(\xi)&=
\big(\big(\xi(s_i^{\R}(0))\!\big)_{\!i\in[k]},\big(\xi(s_i^+(0))\!\big)_{\!i\in[l]}\big),\\
\ev_{\cC_0}^+\!:\ker D_0^+&\lra 
\wt{V}_{\nod^+}\!\oplus\!\bigoplus_{i=l+1}^{a'}\!\!\!\wt{V}_{s_{l+i}^+(0)}, &
\ev_{\cC_0}^+(\xi)&=\big(\xi(\nod^+),\big(\xi(s_{l+i}^+(0))\!\big)_{\!i\in[a']}\big),
\end{aligned}\EE
are isomorphisms.
Along with the sentence containing~\eref{cOevdfn_e}, this implies that 
\eref{evcUdfn_e} with~$l$ replaced~by 
$$l\!+\!a'\equiv \frac{\deg V\!+\!1\!-\!k}{2}$$
is an isomorphism of vector bundles over~$\De_{\R}$.\\

Since each operator $D_{\t}$ is surjective, the last isomorphism induces
an isomorphism
\BE{cOpincglmm_e3}
\la(\cD)\lra \bigotimes_{i=1}^k s_i^{\R*}\wt{V}^{\wt\vph} \otimes
\bigotimes_{i=1}^{l+a'}\!s_i^{+*}\!\la(\wt{V})\EE
of line bundles over~$\De_{\R}$.
The first product on the right-hand side above is oriented by transferring 
the orientation~$\fo_{x_1}$ from the point~$x_1$ in $S^1\!\subset\!\Si_{\t_1}$ 
along the positive direction of $S^1\!\subset\!\bD^2_+$ to the points~$x_i$
with $i\!=\!2,\ldots,k$ and then along the section~$s_i^{\R}$ to an orientation
$\fo_{x_i;\t}$ of $\wt{V}^{\wt\vph}_{s_i^{\R}(\t)}$ for each $\t\!\in\!\De_{\R}$.
The second product is oriented by the complex orientation of~$\wt{V}$.
By definition, the first isomorphism in~\eref{cOpincg_e2} preserves
the $(k,l)$-evaluation $\fo_{k,l}(\fo_{x_i;0})$.
The second isomorphism in~\eref{cOpincg_e2} preserves the complex orientation 
on its domain.
Thus, the restriction of~\eref{cOpincglmm_e3} to the fiber over~$\t_1$
is orientation-preserving with respect to 
the $(k,l)$-intrinsic orientation $\fo_{k,l}(V,\vph;\fo_{x_1})$ on its domain.
By definition, this restriction is also orientation-preserving with respect 
to the $(k,l\!+\!a')$-evaluation orientation~$\fo_{k,l+a'}(\fo_{x_1})$.
This establishes the first claim of the lemma.
\end{proof}

Lemma~\ref{cOpincg_lmm}\ref{cOorientcomp_it} is proved in Section~\ref{tauorient_subs2a}.

\subsection{Line bundles over $C$ degenerations of $(S^2,\tau)$}
\label{tauorient_subs2a}

We now describe the behavior of the orientations of Proposition~\ref{cOpincg_prp}
under flat degenerations of~$(S^2,\tau)$ as in the CROrient~\ref{CROdegenC_prop} property 
on page~\pageref{CROdegenC_prop}.
Suppose $k,l\!\in\!\Z^{\ge0}$,
$\cC_0\!\equiv\!(\Si_0,\si_0)$ is a connected decorated symmetric surface
as above~\eref{3compcUdfn_e} and in the left diagram in
Figure~\ref{P1CH3split_fig}, 
$(V_0,\vph_0)$ is a rank~1 real bundle pair over $(\Si_0,\si_0)$ such~that
\BE{Cklcond_e} k\!+\!2l\cong\deg V_0|_{\P^1_0}\!+\!1 \cong
\deg V_0\!+\!1 \mod2,\EE
and $\fo_{x_1}$ is an orientation of $V_0^{\vph_0}$ at a point~$x_1$ of the fixed locus
$S^1\!=\!\Si_0^{\si_0}$ of~$\si_0$.
With $\nod^{\pm}$ and $\P^1_0,\P^1_{\pm}$ as above~\eref{3compcUdfn_e}, 
we denote the restriction of $(V_0,\vph_0)$  to~$\P^1_0$ by $(V_{00},\vph_{00})$.\\

Let $D_0$ be a real CR-operator on $(V_0,\vph_0)$.
We denote its restrictions to real CR-operators on $(V_{00},\vph_{00})$ and~$V_0|_{\P^1_+}$
by $D_{00}$ and $D_0^+$, respectively. 
The orientation 
$$\fo_{0;0}(\fo_{x_1})\equiv\fo_{k,l}(V_{00},\vph_{00};\fo_{x_1})$$ 
of $D_{00}$
and the complex orientations of $D_0^+$ and $V_0|_{\nod^+}$
determine an orientation 
\BE{rk1Csplit_e}
\fo_0(\fo_{x_1})\equiv \fo_{k,l}\big(V_0,\vph_0;\fo_{x_1}\big)\EE
of $D_0$ via the isomorphism~\eref{Cndses_e2}.
In an analogy with the $\C$-split orientation of Corollary~\ref{CROdegenC_crl},
we call \eref{rk1Csplit_e} the \sf{$(k,l)$-split orientation} of~$D_0$
induced by~$\fo_{x_1}$.\\

Suppose in addition~\eref{cUsymmdfn_e} is a flat family of deformations of~$\cC_0$,
$(V,\vph)$ is a real bundle pair over $(\cU,\wt\fc)$ extending~$(V_0,\vph_0)$,
$s_1^{\R}$ is a section of $\cU^{\wt\fc}$ over~$\De_{\R}$ with $s_1^{\R}(0)\!=\!x_1$,
and  \hbox{$\cD\!\equiv\!\{D_{\t}\}$} is a family of real CR-operators on
$(V_{\t},\vph_{\t})$ as in~\eref{DVvphdfn_e} extending~$D_0$.
The decorated structure on~$\cC_0$ induces a decorated structure on
the fiber $(\Si_{\t},\si_{\t})$ of~$\pi$ for every $\t\!\in\!\De_{\R}$ as 
above the CROrient~\ref{CROdegenC_prop} property.
The orientation~$\fo_{x_1}$ of $V_0^{\vph_0}|_{x_1}$ induces an orientation
of the line bundle $s_1^{\R*}V^{\vph}$ over $\De_{\R}$, which in turn restricts to
an orientation $\fo_{x_1;\t}$ of $V_{\t}^{\vph_{\t}}$ at $s_1^{\R}(\t)$ for each $\t\!\in\!\De_{\R}$.
The orientations
$$\fo_{\t}(\fo_{x_1})\equiv
\fo_{k,l}\big(V_{\t},\vph_{\t};\fo_{x_1;\t}\!\big), \qquad \t\!\in\!\De_{\R}^*,$$
of $D_{\t}$ depend continuously on~$\t$ and extend to an orientation
\BE{rk1Clim_e}
\fo_0'(\fo_{x_1})\equiv\fo_{k,l}'\big(V_0,\vph_0;\fo_{x_1}\big)\EE
of $D_0$.
In an analogy with the limiting orientation of the CROrient~\ref{CROdegenC_prop} property,
we call \eref{rk1Clim_e} the \sf{limiting orientation} of~$D_0$
induced by~$\fo_{x_1}$.
It depends only on~$(V_0,\vph_0)$ and $\fo_{x_1}$,
and not on~$(\cU,\wt\fc)$, $(V,\vph)$, or~$s_1^{\R}$.

\begin{prp}\label{rk1Cdegen_prp}
Suppose $\cC_0$ is a connected decorated symmetric surface with precisely 
one pair of conjugate nodes and each component isomorphic to~$S^2$
$(V_0,\vph_0)$  is a rank~1 real bundle pair over~$(\Si_0,\si_0)$.
Let $k,l\!\in\!\Z^{\ge0}$ be so that \eref{Cklcond_e} holds
and $D_0$ be a real CR-operator on~$(V_0,\vph_0)$.
The $(k,l)$-split and limiting orientations, \eref{rk1Csplit_e} and~\eref{rk1Clim_e},
of~$D_0$ are the~same.
\end{prp}

We deduce this proposition from Lemma~\ref{cOevComp_lmm} below.
It shows that certain orientations $\fo_{k,l}'(V,\vph;\fo_{x_1})$
of a real CR-operator over a smooth symmetric surface
constructed via a flat family of deformations of a symmetric surface with two pairs 
of conjugate nodes, instead of one as above Proposition~\ref{cOpincg_prp}, 
are in fact the same as $\fo_{k,l}(V,\vph;\fo_{x_1})$.\\

Suppose $(V,\vph)$, $k$, $l$, $x_1$, $\fo_{x_1}$, and $D$ are as above Lemma~\ref{cOev_lmm}.
Let $\cC_0\!\equiv\!(\Si_0,\si_0)$ be a connected symmetric surface 
which contains precisely two conjugate pairs of nodes and
consists of one invariant component and two pairs of conjugate components attached to~it, 
all of which are isomorphic to~$S^2$;
see the left diagram in Figure~\ref{P1Cdegen_fig2a}.
Let $\P^1_0\!\subset\!\Si_0$ be the component preserved by~$\si_0$,
$\bD^2_+,\bD^2_-\!\subset\!\P^1_0$ be the two half-surfaces,
and $\nod_1^+,\nod_2^+\!\in\!\bD^2_+$ be the two nodes on~$\bD^2_+$.
We denote by $\P^1_{1+}$ and $\P^1_{2+}$ the other irreducible 
components containing~$\nod_1^+$ and~$\nod_2^+$, respectively.
Let $(\cU,\wt\fc)$ be a flat family of deformations of~$\cC_0$ as in~\eref{cUsymmdfn_e}
over the unit ball $\De\!\subset\!\C^4$ around the origin satisfying~\eref{3compcUdfn_e}
for some $\t_1\!\in\!\De_{\R}$.\\

Let $a_0\!=\!k\!+\!2l\!-\!1$ be as before and $a_1',a_2'\!\in\!\Z$ be such that
\BE{degsplitcond_e}a_0+2(a_1'\!+a_2')=\deg V.\EE
Let $(\wt{V},\wt\vph)$ be a rank~1 real bundle pair over $(\cU,\wt\fc)$ so that 
\BE{pinchcond5_e3}\deg\wt{V}\big|_{\P^1_0}=a_0,
\quad  \deg\wt{V}\big|_{\P^1_{1+}}=a_1',\quad  
\deg\wt{V}\big|_{\P^1_{2+}}=a_2',
\quad \big(\wt{V}_{\t_1},\wt\vph_{\t_1}\big)=(V,\vph)\EE
and \hbox{$\cD\!\equiv\!\{D_{\t}\}$} be a family of real CR-operators on
$(\wt{V}_{\t},\wt\vph_{\t})$ as in~\eref{DVvphdfn_e} so~that $D_{\t_1}\!=\!D$.
The orientation~$\fo_{x_1}$ of $V^{\vph}$ induces an orientation~$\fo_{x_1'}$
of $\wt{V}^{\wt\vph}$
at a point~$x_1'$ in the fixed locus~$S^1\!\subset\!\Si_0$ as before.\\

We denote the restriction of $(\wt{V}_0,\wt\vph_0)$  to~$\P^1_0$ by $(\wt{V}_{00},\wt\vph_{00})$ and 
the restrictions of $D_0$ to real CR-operators on $(\wt{V}_{00},\wt\vph_{00})$, $\wt{V}_0|_{\P^1_{1+}}$,
and $\wt{V}_0|_{\P^1_{2+}}$ by $D_{00}$, $D_{01}^+$, and $D_{02}^+$, respectively. 
The exact~triple 
\begin{gather}
\label{ses5comp_e}
0\lra D_0 \lra 
D_{00}\!\oplus\!D_{01}^+\!\oplus\!D_{02}^+
\lra \wt{V}_{\nod_1^+}\!\oplus\!\wt{V}_{\nod_2^+}\lra 0, \\
\notag
\big(\xi_{2-},\xi_{1-},\xi_0,\xi_{1+},\xi_{2+}\big)\lra \big(\xi_0,\xi_{1+},\xi_{2+}\big), \\
\notag
 \big(\xi_0,\xi_{1+},\xi_{2+}\big)\lra
\big(\xi_{1+}(\nod_1^+)\!-\!\xi_0(\nod_1^+),\xi_{2+}(\nod_2^+)\!-\!\xi_0(\nod_2^+)\big),
\end{gather}
of Fredholm operators then determines an isomorphism
\BE{ses5comp_e2}
\la(D_0)\!\otimes\!\la\big(\wt{V}_{\nod_1^+}\big)\!\otimes\!\la\big(\wt{V}_{\nod_2^+}\big)
\approx  \la\big(D_{00}\big)\!\otimes\!\la\big(D_{01}^+\big)
\!\otimes\!\la\big(D_{02}^+\big).\EE
The $(k,l)$-evaluation orientation $\fo_{k,l}(\fo_{x_1'})$ of $D_{00}$
and the complex orientations of
$D_{01}^+$, $D_{02}^+$, $\wt{V}_{\nod_1^+}$, and $\wt{V}_{\nod_2^+}$
determine an orientation of $D_0$ via the isomorphism~\eref{ses5comp_e2}  and 
thus an orientation $\wt\fo_{k,l}'(\wt{V},\wt\vph;\fo_{x_1})$
of the line bundle~$\la(\cD)$ over~$\De_{\R}$.
The latter restricts to an orientation $\fo_{k,l}'(V,\vph;\fo_{x_1})$ of~$\la(D)$.

\begin{lmm}\label{cOevComp_lmm}
Suppose $(V,\vph)$ is a rank~1 real bundle pair  over~$(S^2,\tau)$,
$D$ is a real CR-operator on~$(V,\vph)$,
$k,l\!\in\!\Z^{\ge0}$ satisfy~\eref{pinchcond_e},
and $a_1',a_2'\!\in\!\Z$ satisfy~\eref{degsplitcond_e}.
Let $\fo_{x_1}$ be an orientation of~$V^{\vph}$ at a point $x_1$ in $S^1\!\subset\!S^2$.
The orientation $\fo_{k,l}'(V,\vph;\fo_{x_1})$ of~$D$ defined above
is the same as the $(k,l)$-intrinsic orientation of Proposition~\ref{cOpincg_prp}.
\end{lmm}

\begin{figure}
\begin{pspicture}(-.3,-3.8)(10,2.5)
\psset{unit=.4cm}
\pscircle(8,-2){2}
\pscircle(5.53,.47){1.5}\pscircle(5.53,-4.47){1.5}
\pscircle(10.47,.47){1.5}\pscircle(10.47,-4.47){1.5}
\pscircle*(6.59,-.59){.2}\pscircle*(6.59,-3.41){.2}
\pscircle*(9.41,-.59){.2}\pscircle*(9.41,-3.41){.2}
\rput(5.7,0){\sm{$\nod_1^+$}}\rput(5.7,-4){\sm{$\nod_1^-$}}
\rput(10.3,0){\sm{$\nod_2^+$}}\rput(10.3,-4){\sm{$\nod_2^-$}}
\psline{<->}(3,2)(3,-6)\rput(2.2,-2){$\si_0$}
\rput(6,2.6){$\P^1_{1+}$}\rput(6,-6.7){$\P^1_{1-}$}
\rput(11,2.6){$\P^1_{2+}$}\rput(11,-6.7){$\P^1_{2-}$}
\rput(8.1,.6){\sm{$\bD^2_+$}}\rput(8.1,-4.7){\sm{$\bD^2_-$}}
\rput(10.8,-2){$\P^1_0$}
\psarc[linewidth=.03](8,1.46){4}{240}{300} 
\psarc[linewidth=.03,linestyle=dashed](8,-5.46){4}{60}{120} 
\pscircle(30,0){2}\pscircle(32.47,2.47){1.5}\pscircle(32.47,-2.47){1.5}
\pscircle*(31.41,1.41){.2}\pscircle*(31.41,-1.41){.2}
\rput(32.3,2){\sm{$\nod_{\t}^+$}}\rput(32.3,-2){\sm{$\nod_{\t}^-$}}
\rput(33,0){$\P^1_{\t;0}$}\rput(33,4.6){$\P^1_{\t+}$}\rput(33,-4.7){$\P^1_{\t-}$}
\psarc[linewidth=.03](30,3.46){4}{240}{300} 
\psarc[linewidth=.03,linestyle=dashed](30,-3.46){4}{60}{120} 
\psline{<->}(35,4)(35,-4)\rput(35.8,0){$\si_{\t}$}
\psline{->}(18,-8)(18,6)\psline{->}(18,-8)(32,-8)
\rput(30.5,-8.8){\sm{$\C^2$}}\rput{90}(17.4,4.4){\sm{$\C^2$}}
\psline[linewidth=.03]{->}(28.5,-2.5)(26,-7)
\psline[linewidth=.03]{->}(12.5,-5.5)(17,-7.5)
\pscircle*(18,-8){.2}\pscircle*(22,-8){.2}
\pscircle*(22,-4){.2}
\rput(22.2,-7.4){$\t_0$}\rput(21.6,-4.4){$\t_1$}
\pscircle(23,2){2.5}\psarc[linewidth=.03](23,6.33){5}{240}{300} 
\psarc[linewidth=.03,linestyle=dashed](23,-2.33){5}{60}{120} 
\psline[linewidth=.03]{->}(23.2,-1)(22.3,-3.7)
\end{pspicture}
\caption{Deformations of the five-component symmetric surface $(\Si_0,\si_0)$
over a unit ball in $\C^4$ in the proof of Proposition~\ref{rk1Cdegen_prp}}
\label{P1Cdegen_fig2a}
\end{figure}

\begin{proof}[{\bf{\emph{Proof of Proposition~\ref{rk1Cdegen_prp}}}}]
In order to distinguish the nodal symmetric surfaces $\cC_0$,
their invariant components~$\P^1_0$,
and the families of deformations $(\cU,\wt\fc)$ over~$\De$ appearing above
the statements of Proposition~\ref{rk1Cdegen_prp} and Lemma~\ref{cOevComp_lmm},
we denote the objects appearing above Lemma~\ref{cOevComp_lmm} by 
$\cC_0'$, $\P'^1_0$, $(\cU',\wt\fc')$, and~$\De'$, respectively, in this proof.
We can choose $(\cU',\wt\fc')$ so that $\cC_{\t_0}'$ for some $\t_0\!\in\!\De_{\R}'$ 
is the nodal symmetric surface~$\cC_0$  as in the statement of the proposition
obtained from~$\cC_0'$ by smoothing the node~$\nod_1^+$; see Figure~\ref{P1Cdegen_fig2a}.
Let $\t_1\!\in\!\De_{\R}'^*$ and 
$s_1'^{\R}$ be a section of $\cU'^{\wt\fc'}$ over~$\De_{\R}'$ with $s_1'^{\R}(\t_0)\!=\!x_1$.\\

With the notation as above the statement of the proposition, let
$$a_0=k\!+\!2l\!-\!1, \qquad a_1'=\big(\!\deg V_0|_{\P^1_0}\!-\!a_0\big)/2, \qquad 
a_2'=\deg V_0|_{\P^1_+}\,.$$
For the purposes of applying Lemma~\ref{cOevComp_lmm},
we choose $(\wt{V},\wt\vph)$ and \hbox{$\cD\!\equiv\!\{D_{\t}\}$} so that 
$$(\wt{V}_{\t_0},\wt\vph_{\t_0})=(V_0,\vph_0) \qquad\hbox{and}\qquad 
\qquad D_{\t_0}=D_0.$$
For each $\t\!\in\!\De_{\R}'$, let $\fo_{x_1;\t}$ be the orientation 
of $\wt{V}^{\wt\vph}$ at $s_1'^{\R}(\t)$ induced by~$\fo_{x_1}$
via the line bundle $s_1'^{\R*}\wt{V}^{\wt\vph}$.
By definition and the evenness of the (real) dimension of~$\wt{V}_{\nod_1^+}$, 
the restriction of the orientation $\wt\fo_{k,l}'(\wt{V},\wt\vph;\fo_{x_1;\t_1})$
of the line bundle~$\la(\cD)$ over~$\De_{\R}'$ to $\la(D_{\t_0})$
is the $(k,l)$-split orientation~$\fo_0(\fo_{x_1})$ in~\eref{rk1Csplit_e}.
By Lemma~\ref{cOevComp_lmm}, the restriction of $\wt\fo_{k,l}'(\wt{V},\wt\vph;\fo_{x_1;\t_1})$
to $\la(D_{\t})$ for each $\t\!\in\!\De_{\R}'^*$ is the $(k,l)$-intrinsic orientation
$\fo_{k,l}(\wt{V}_{\t}',\wt\vph_{\t}';\fo_{x_1;\t})$.
Thus, the restriction of $\wt\fo_{k,l}'(\wt{V},\wt\vph;\fo_{x_1;\t_1})$ to $\la(D_{\t_0})$
is the limiting orientation~$\fo_0'(\fo_{x_1})$ in~\eref{rk1Clim_e}.
We conclude that the orientations~\eref{rk1Csplit_e} and~\eref{rk1Clim_e} are the same.
\end{proof}

\begin{proof}[{\bf{\emph{Proof of Lemma~\ref{cOpincg_lmm}\ref{cOorientcomp_it}}}}]
Suppose $l,\wt{l}\!\in\!\Z^{\ge0}$ with $\wt{l}\!\ge\!l$.
Let $\cC_0$, $\P_0^1$, $\P^1_+$, $(\cU,\wt\fc)$, $\De_{\R}$, and $\t_1$ 
be as above Proposition~\ref{cOpincg_prp} and 
$s_1^{\R}$ be a section of $\cU^{\wt\fc}$ over~$\De_{\R}$ with $s_1^{\R}(\t_1)\!=\!x_1$.
We can choose $(\wt{V},\wt\vph)$ and \hbox{$\cD\!\equiv\!\{D_{\t}\}$}
as above Proposition~\ref{cOpincg_prp} so~that
$$\deg\wt{V}\big|_{\P^1_0}=k\!+\!2\wt{l}\!-\!1,\quad 
\deg\wt{V}\big|_{\P^1_+}=\big(\!\deg V\!-\!\deg\wt{V}|_{\P^1_0}\big)/2, \quad
\big(\wt{V}_{\t_1},\wt\vph_{\t_1}\big)=(V,\vph), 
\quad D_{\t_1}=D\,.$$
For each $\t\!\in\!\De_{\R}$, let $\fo_{x_1;\t}$ be the orientation 
of $\wt{V}^{\wt\vph}$ at $s_1^{\R}(\t)$ induced by~$\fo_{x_1}$
via the line bundle $s_1^{\R*}\wt{V}^{\wt\vph}$.\\

Since $\wt{l}\!\ge\!l$, 
the definitions of the $(k,l)$- and $(k,\wt{l})$-split orientations of~$D_0$ 
and Lemma~\ref{cOpincg_lmm}\ref{cOpincgl_it1} imply~that 
$$\fo_{k,l}\big(\wt{V}_0,\wt\vph_0;\fo_{x_1;0}\big)
=\fo_{k,\wt{l}}\big(\wt{V}_0,\wt\vph_0;\fo_{x_1;0}\big).$$ 
Along with Proposition~\ref{rk1Cdegen_prp}, this in turn implies that the corresponding
limiting orientations of~$D_0$ are equal, i.e.
$$\fo_{k,l}'\big(\wt{V}_0,\wt\vph_0;\fo_{x_1;0}\big)
=\fo_{k,\wt{l}}'\big(\wt{V}_0,\wt\vph_0;\fo_{x_1;0}\big).$$
Thus, the $(k,l)$- and $(k,\wt{l})$-intrinsic orientations 
$\fo_{k,l}(\wt{V}_{\t},\wt\vph_{\t};\fo_{x_1;\t})$ and 
$\fo_{k,\wt{l}}(\wt{V}_{\t},\wt\vph_{\t};\fo_{x_1;\t})$
of~$D_{\t}$ with $\t\!\in\!\De_{\R}^*$ are the same.
Taking $\t\!=\!\t_1$, we obtain the claim.
\end{proof}

\begin{proof}[{\bf{\emph{Proof of Lemma~\ref{cOevComp_lmm}}}}] 
Let $\cC_0'\!\equiv\!(\Si_0',\si_0')$ be a connected symmetric surface 
which contains precisely three conjugate pairs of nodes and
consists of one invariant component~$\P^1_0$ and three pairs of conjugate components, 
all of which are isomorphic to~$S^2$, so that one pair
of the conjugate components is attached to~$\P^1_0$ and 
the other two pairs are attached to the first pair;
see the bottom diagram in Figure~\ref{cOevCompPf_fig}.
We denote by $\bD^2_+,\bD^2_-\!\subset\!\P^1_0$ the two half-surfaces,
by $\nod_0^+\!\in\!\bD^2_+$ the node on~$\bD^2_+$, 
by $\P^1_{0+}\!\subset\!\Si_0'$ the other irreducible component containing~$\nod_0^+$,
by $\nod_1^+,\nod_2^+\!\in\!\P^1_{0+}$ the other nodes on~$\P^1_{0+}$,
and by  $\P^1_{1+}$ and $\P^1_{2+}$ the other irreducible 
components containing~$\nod_1^+$ and~$\nod_2^+$, respectively.\\
 
Let  $(\cU',\wt\fc')$ be a flat family of deformations of~$\cC_0'$ as in~\eref{cUsymmdfn_e}
over the unit ball $\De'\!\subset\!\C^6$ around the origin satisfying~\eref{3compcUdfn_e}
for some $\t_1\!\in\!\De_{\R}'$ so that 
$$\cC_{\t_4}'\equiv \big(\Si_{\t_4}',\si_{\t_4}'\big)
\qquad\hbox{and}\qquad
\cC_{\t_2}'\equiv \big(\Si_{\t_2}',\si_{\t_2}'\big)$$
for some $\t_4,\t_2\!\in\!\De_{\R}'$ are 
a five-component decorated symmetric surface as on the left side of Figure~\ref{P1Cdegen_fig2a}
obtained from~$\cC_0'$ by smoothing the conjugate pair of nodes on~$\P^1_0$  and
a three-component decorated symmetric surface as on the left side of Figure~\ref{P1CH3split_fig}
obtained from~$\cC_0'$ by smoothing the other two pairs of nodes,
respectively; see Figure~\ref{cOevCompPf_fig}.
We denote~by 
$$\De_{\R;4}'^*,\De_{\R;2}'^*\subset\De_{\R}'$$
the subspaces parametrizing all five-component symmetric surfaces 
 obtained from~$\cC_0'$
by smoothing the conjugate pair of nodes on~$\P^1_0$
and all three-component  symmetric surfaces  
obtained from~$\cC_0'$ by smoothing the other two pairs of nodes,
respectively.
Let $s_1'^{\R}$ be a section of $\cU'^{\wt\fc'}$ over~$\De_{\R}'$ with 
\hbox{$s_1'^{\R}(\t_1)\!=\!x_1$}.\\

Let $(\wt{V}',\wt\vph')$ be a rank~1 real bundle pair over $(\cU',\wt\fc')$
so~that 
$$\deg\wt{V}'\big|_{\P^1_0}=a_0, ~~ \deg\wt{V}'\big|_{\P^1_{0+}}=0, ~~
~~  \deg\wt{V}'\big|_{\P^1_{1+}}=a_1',~~
\deg\wt{V}'\big|_{\P^1_{2+}}=a_2',
~~ \big(\wt{V}'_{\t_1},\wt\vph'_{\t_1}\big)=(V,\vph).$$
For each $\t\!\in\!\De_{\R}'$, let $\fo_{x_1;\t}$ be the orientation 
of $\wt{V}'^{\wt\vph'}$ at $s_1'^{\R}(\t)$ induced by~$\fo_{x_1}$
via the line bundle $s_1'^{\R*}\wt{V}'^{\wt\vph'}$.\\

\begin{figure}
\begin{pspicture}(-1,-6.4)(10,4.5)
\psset{unit=.4cm}
\pscircle(8,5){2}
\pscircle(5.53,7.47){1.5}\pscircle(5.53,2.53){1.5}
\pscircle(10.47,7.47){1.5}\pscircle(10.47,2.53){1.5}
\pscircle*(6.59,6.41){.2}\pscircle*(6.59,3.59){.2}
\pscircle*(9.41,6.41){.2}\pscircle*(9.41,3.59){.2}
\rput(5.5,7){\sm{$\nod_{\t;1}^+$}}\rput(5.7,3){\sm{$\nod_{\t;1}^-$}}
\rput(10.4,7){\sm{$\nod_{\t;2}^+$}}\rput(10.3,3){\sm{$\nod_{\t;2}^-$}}
\psline{<->}(3,9)(3,1)\rput(2.2,5){$\si_{\t}'$}
\rput(6,9.6){$\P^1_{\t;1+}$}\rput(6,.3){$\P^1_{\t;1-}$}
\rput(11,9.6){$\P^1_{\t;2+}$}\rput(11,.3){$\P^1_{\t;2-}$}
\rput(11,5){$\P^1_{\t;0}$}
\psarc[linewidth=.03](8,8.46){4}{240}{300} 
\psarc[linewidth=.03,linestyle=dashed](8,1.54){4}{60}{120} 
\psline[linewidth=.03]{->}(12.5,5)(17,0)
\pscircle(30,0){2}\pscircle(30,4){2}\pscircle(30,-4){2}
\pscircle*(30,2){.2}\pscircle*(30,-2){.2}
\rput(30,2.8){\sm{$\nod_{\t}^+$}}\rput(30,-2.8){\sm{$\nod_{\t}^-$}}
\rput(33,0){$\P^1_{\t;0}$}\rput(33,3.9){$\P^1_{\t;+}$}\rput(33,-3.8){$\P^1_{\t;-}$}
\psarc[linewidth=.03](30,3.46){4}{240}{300} 
\psarc[linewidth=.03,linestyle=dashed](30,-3.46){4}{60}{120} 
\psline{<->}(35,4)(35,-4)\rput(35.8,0){$\si_{\t}'$}
\psline{->}(18,-8)(18,6)\psline{->}(18,-8)(32,-8)
\rput(30,-9){\sm{$\De_{\R;2}'^*\!\subset\!\C^4$}}
\rput{90}(17,4){\sm{$\De_{\R;4}'^*\!\subset\!\C^2$}}
\psline[linewidth=.03]{->}(27.5,-1.5)(26,-7)
\pscircle*(18,-8){.2}\pscircle*(22,-8){.2}
\pscircle*(18,-4){.2}\pscircle*(22,-4){.2}
\rput(18.7,-4){$\t_4$}\rput(22,-7.2){$\t_2$}\rput(21.6,-4.4){$\t_1$}
\pscircle(23,2){2.5}\psarc[linewidth=.03](23,6.33){5}{240}{300} 
\psarc[linewidth=.03,linestyle=dashed](23,-2.33){5}{60}{120} 
\psline[linewidth=.03]{->}(22.5,-1)(22.1,-3.5)
\pscircle(8,-8){1.5}\pscircle(8,-5){1.5}\pscircle(8,-11){1.5}
\pscircle(5.88,-2.88){1.5}\pscircle(10.12,-2.88){1.5}
\pscircle(5.88,-13.12){1.5}\pscircle(10.12,-13.12){1.5}
\pscircle*(8,-6.5){.2}\pscircle*(8,-9.5){.2}
\pscircle*(6.94,-3.94){.2}\pscircle*(6.94,-12.06){.2}
\pscircle*(9.06,-3.94){.2}\pscircle*(9.06,-12.06){.2}
\rput(8,-5.7){\sm{$\nod_0^+$}}\rput(8,-10.3){\sm{$\nod_0^-$}}
\rput(6.3,-3.2){\sm{$\nod_1^+$}}\rput(6.1,-12.7){\sm{$\nod_1^-$}}
\rput(10.1,-3.3){\sm{$\nod_2^+$}}\rput(10.1,-12.6){\sm{$\nod_2^-$}}
\rput(3.5,-2.5){$\P^1_{1+}$}\rput(3.5,-13.5){$\P^1_{1-}$}
\rput(12.6,-2.5){$\P^1_{2+}$}\rput(12.6,-13.5){$\P^1_{2-}$}
\rput(6,-7){\sm{$\bD^2_+$}}\rput(6,-9.1){\sm{$\bD^2_-$}}
\rput(10.3,-8){$\P^1_0$}\rput(10.5,-5.3){$\P^1_{0+}$}\rput(10.5,-10.6){$\P^1_{0-}$}
\psarc[linewidth=.03](8,-5.41){3}{240}{300} 
\psarc[linewidth=.03,linestyle=dashed](8,-10.59){3}{60}{120} 
\psline{<->}(2,-4)(2,-12)\rput(1.2,-8){$\si_0'$}
\rput(8,-15.5){$\cC_0'$}
\psline[linewidth=.03]{->}(12.5,-8)(17,-8)
\end{pspicture}
\caption{Deformations of~$\cC_0'$ in the family $(\cU',\wt\fc')$ over 
$\De'\!\subset\!\C^6$ of the proof of Lemma~\ref{cOevComp_lmm}}
\label{cOevCompPf_fig}
\end{figure}

Let \hbox{$\cD'\!\equiv\!\{D_{\t}'\}$} be a family of real CR-operators on 
$(\wt{V}_{\t}',\wt\vph_{\t}')$ 
as in~\eref{DVvphdfn_e} so that $D_{\t_1}'$ is the operator~$D$ in the statement of 
the lemma. 
We denote by $(\wt{V}'_{00},\wt\vph'_{00})$ the restriction of $(\wt{V}',\wt\vph')$ 
to~$\P^1_0$ and 
 by~$D_{00}'$ the restriction of $D_0'$ to a real CR-operator on $(\wt{V}_{00}',\wt\vph_{00}')$.
For $r\!=\!0,1,2$, let $D_r'^+$ be the restriction of~$D_0'$
to a real CR-operator on $\wt{V}_0|_{\P^1_{r+}}$. 
The exact~triple 
\begin{gather*}
0\lra D'_0 \lra 
D'_{00}\!\oplus\!D_0'^+\!\oplus\!D_1'^+\!\oplus\!D_2'^+\!
\lra \wt{V}'_{\nod_0^+}\!\oplus\!\wt{V}'_{\nod_1^+}\!\oplus\!\wt{V}'_{\nod_2^+}\lra 0, \\
\notag
\big(\xi_{2-},\xi_{1-},\xi_{0-},\xi_0,\xi_{0+},\xi_{1+},\xi_{2+}\big)\lra
\big(\xi_0,\xi_{0+},\xi_{1+},\xi_{2+}\big),\\
\notag
\big(\xi_0,\xi_{0+},\xi_{1+},\xi_{2+}\big)\lra 
\big(\xi_{0+}(\nod_0^+)\!-\!\xi_0(\nod_0^+),\xi_{1+}(\nod_1^+)\!-\!\xi_{0+}(\nod_1^+),
\xi_{2+}(\nod_2^+)\!-\!\xi_{0+}(\nod_2^+)\big),
\end{gather*}
of Fredholm operators then determines an isomorphism
\BE{cOevComp_e4}\la(D_0')\!\otimes\!\la\big(\wt{V}'_{\nod_0^+}\big)
\!\otimes\!\la\big(\wt{V}'_{\nod_1^+}\big)\!\otimes\!
\la\big(\wt{V}'_{\nod_2^+}\big)
\approx  \la\big(D'_{00}\big)\!\otimes\!\la\big(D_0'^+\big)
\!\otimes\!\la\big(D_1'^+\big)\!\otimes\!\la\big(D_2'^+\big).\EE
The $(k,l)$-evaluation orientation $\fo_{k,l}(\fo_{x_1;0})$ of $D'_{00}$
and the complex orientations of $D_0'^+$, $D_1'^+$, $D_2'^+$, $\wt{V}'_{\nod_0^+}$,
$\wt{V}'_{\nod_1^+}$, and $\wt{V}'_{\nod_2^+}$
determine an orientation of $D'_0$ via the isomorphism~\eref{cOevComp_e4}
and  thus an orientation~$\fo_{\cU'}(\fo_{x_1})$ of the line bundle~$\la(\cD')$ 
over~$\De'_{\R}$.\\

The restriction of $(\cU',\wt\fc')$ to $\De_{\R;4}'$ is the product of
two conjugates pairs of~$\P^1$ with a family smoothing 
$$\P^1_{0;0}\equiv \P^1_{0-}\!\cup\!\P^1_0\!\cup\!\P^1_{0+}$$
into a single irreducible component $\P^1_{\t;0}\!\subset\!\Si_{\t}'$
as above Proposition~\ref{cOpincg_prp};
see the top left diagram in Figure~\ref{cOevCompPf_fig}.
The other irreducible components of~$\Si_{\t}'$ with $\t\!\in\!\De_{\R;4}'^*$ 
and its nodes correspond to the conjugate components and conjugate nodes of~$\Si_0'$;
we denote them by~$\P^1_{\t;r\pm}$ and~$\nod_{\t;r}^{\pm}$, respectively.
For each $\t\!\in\!\De_{\R;4}'$, 
we denote the restriction of $(\wt{V}_{\t}',\wt\vph_{\t}')$  to~$\P^1_{\t;0}$ by 
$(\wt{V}'_{\t;0},\wt\vph'_{\t;0})$ 
and 
the restrictions of $D_{\t}'$ to real CR-operators on $(\wt{V}'_{\t;0},\wt\vph'_{\t;0})$, 
$\wt{V}_{\t}'|_{\P^1_{\t;1+}}$,
and $\wt{V}_0|_{\P^1_{\t;2+}}$ by $D_{\t;0}'$, $D_{\t1}'^+$, and $D_{\t2}'^+$, respectively.
The exact triple~\eref{Cndses_e} of Fredholm operators with~$D_0$ replaced 
by~$D'_{0;0}$ and 
the exact triple~\eref{ses5comp_e} of Fredholm operators with~$D_0$ replaced
by~$D'_{\t_4}$ induce isomorphisms
\BE{cOevComp_e7}\begin{split}
\la\big(D'_{0;0}\big)\!\otimes\!\la\big(\wt{V}'_{\nod_0^+}\big)&\approx 
\la\big(D'_{00}\big)\!\otimes\!\la\big(D_0'^+\big),\\
\la(D_{\t_4}')\!\otimes\!
\la\big(\wt{V}'_{\nod_{\t_4;1}^+}\big)\!\otimes\!\la\big(\wt{V}'_{\nod_{\t_4;2}^+}\big) &\approx 
\la\big(D_{\t_4;0}'\big)\!\otimes\!\la\big(D_{\t_41}'^+\big)
\!\otimes\!\la\big(D_{\t_42}'^+\big).
\end{split}\EE
The $(k,l)$-evaluation orientation $\fo_{k,l}(\fo_{x_1;0})$ of $D_{00}'$
and the complex orientations of $D_0'^+$ and $\wt{V}'_{\nod_0^+}$
determine an orientation on $D_{0;0}'$
via the first isomorphism in~\eref{cOevComp_e7} 
and thus an orientation of $D_{\t;0}'$ for each $\t\!\in\!\De_{\R;4}'$.
By definition, the latter with $\t\!=\!\t_4$ is the $(k,l)$-intrinsic orientation
$\fo_{k,l}(\wt{V}'_{\t_4;0},\wt\vph'_{\t_4;0};\fo_{x_1;\t_4})$.
By Lemma~\ref{cOpincg_lmm}\ref{cOpincgl_it1}, this orientation is
the same as the $(k,l)$-evaluation orientation $\fo_{k,l}(\fo_{x_1;\t_4})$.
Since the (real) dimension of~$\wt{V}'_{\nod_0^+}$ is even, 
this implies that the second isomorphism in~\eref{cOevComp_e7} is orientation-preserving 
with respect to the restriction of the orientation $\fo_{\cU'}(\fo_{x_1})$ to $\la(D_{\t_4}')$,
the orientation~$\fo_{k,l}(\fo_{x_1;\t_4})$ of $D_{\t_4;0}'$,
and the complex orientations of the remaining factors.\\

The restriction of $(\cU',\wt\fc')$ to $\De_{\R;2}'$ is the product of 
$(\P^1_0,\tau)$ with a family smoothing the conjugate~pair
$$\P^1_{0;\pm}\equiv \P^1_{0\pm}\cup \P^1_{1\pm}\cup \P^1_{2\pm}\subset\Si_0'$$
into a conjugate pair of irreducible components $\P^1_{\t;\pm}\!\subset\!\Si_{\t}'$;
see the last diagram in Figure~\ref{cOevCompPf_fig}.
The other irreducible component of~$\Si_{\t}'$ with $\t\!\in\!\De_{\R;2}^*$ 
and its nodes correspond to~$\P^1_0$ and the nodes~$\nod_0^{\pm}$;
we denote them by~$\P^1_{\t0}$ and~$\nod_{\t}^{\pm}$, respectively.
For each $\t\!\in\!\De_{\R;2}'$, 
we denote the restriction of $(\wt{V}_{\t}',\wt\vph_{\t}')$  to~$\P^1_{\t0}$ by 
$(\wt{V}'_{\t0},\wt\vph'_{\t0})$ and 
the restrictions of $D_{\t}'$ to real CR-operators on $(\wt{V}'_{\t0},\wt\vph'_{\t0})$ and 
$\wt{V}_{\t}'|_{\P^1_{\t;+}}$ by $D_{\t0}'$ and $D_{\t;+}'^+$, respectively.
The exact~triple 
\begin{equation*}\begin{split}
0\lra D'_{0;+} \lra 
&D_0'^+\!\oplus\!D_1'^+\!\oplus\!D_2'^+
\lra \wt{V}'_{\nod_1^+}\!\oplus\!\wt{V}'_{\nod_2^+}\lra 0, \\
&\big(\xi_{0+},\xi_{1+},\xi_{2+}\big)\lra
\big(\xi_{1+}(\nod_1^+)\!-\!\xi_{0+}(\nod_1^+),\xi_{2+}(\nod_2^+)\!-\!\xi_{0+}(\nod_2^+)\big),
\end{split}\end{equation*}
of Fredholm operators and the exact triple~\eref{Cndses_e} of Fredholm operators 
with~$D_0$ replaced by~$D_{\t_2}'$ induce isomorphisms
\BE{cOevComp_e9}\begin{split}
\la\big(D_{0;+}'\big)
\!\otimes\!\la\big(\wt{V}'_{\nod_1^+}\big)\!\otimes\!\la\big(\wt{V}'_{\nod_2^+}\big)
&\approx 
\la\big(D_0'^+\big)\!\otimes\!\la\big(D_1'^+\big)
\!\otimes\!\la\big(D_2'^+\big),\\
\la(D_{\t_2}')\!\otimes\!
\la\big(\wt{V}'_{\nod_{\t_2}^+}\big)&\approx 
\la\big(D'_{\t_20}\big)\!\otimes\!\la\big(D_{\t_2;+}'\big).
\end{split}\EE
The complex orientations of $D_r'^+$ for $r\!=\!0,1,2$
and $\wt{V}'_{\nod_r^+}$ for $r\!=\!1,2$
determine an orientation of $D_{0;+}'$ via
the first isomorphism in~\eref{cOevComp_e9} 
and thus an orientation of $D_{\t;+}'$ for each $\t\!\in\!\De_{\R;2}'$.
The latter is the complex orientation of~$D_{\t;+}'$.
This implies that the second isomorphism in~\eref{cOevComp_e9} is orientation-preserving 
with respect to the restriction of the orientation $\fo_{\cU'}(\fo_{x_1})$ to $\la(D_{\t_2}')$,
the $(k,l)$-evaluation orientation~$\fo_{k,l}(\fo_{x_1;\t_2})$ of $D'_{\t_20}$,
and the complex orientations of the remaining factors.\\

By the conclusion above regarding the second isomorphism in~\eref{cOevComp_e7},
the restriction of the orientation $\fo_{\cU'}(\fo_{x_1})$ to $\la(D_{\t_1}')$
is the orientation $\fo_{k,l}'(V,\vph;\fo_{x_1})$ constructed above Lemma~\ref{cOevComp_lmm}.
By the conclusion above regarding the second isomorphism in~\eref{cOevComp_e9},
the restriction of the orientation $\fo_{\cU'}(\fo_{x_1})$ to $\la(D_{\t_1}')$
is the $(k,l)$-intrinsic orientation $\fo_{k,l}(V,\vph;\fo_{x_1})$  
constructed above Proposition~\ref{cOpincg_prp}.
Thus, the orientations
$\fo_{k,l}(V,\vph;\fo_{x_1})$ and $\fo_{k,l}'(V,\vph;\fo_{x_1})$ are the same.  
\end{proof}

\subsection{Line bundles over $H3$ degenerations of $(S^2,\tau)$}
\label{tauorient_subs2b}

We next describe the behavior of the orientations of Proposition~\ref{cOpincg_prp}
under flat degenerations of~$(S^2,\tau)$ as in the CROrient~\ref{CROdegenH3_prop} property 
on page~\pageref{CROdegenH3_prop}.
Let $\cC_0\!\equiv\!(\Si_0,\si_0)$ be a connected decorated symmetric surface
as above~\eref{pinchcond_e3c} and in the right diagram in
Figure~\ref{P1CH3split_fig} and 
$\nod$, $\P^1_1,\P^1_2$, $\bD^2_{1\pm}\bD^2_{2\pm}$, and $S^1_1,S^1_2$ 
be as above~\eref{pinchcond_e3c}.
Suppose $k_1,l_1,k_2,l_2\!\in\!\Z^{\ge0}$ and
$(V_0,\vph_0)$ is a rank~1 real bundle pair over $(\Si_0,\si_0)$ such~that
\BE{k1l1k2l2_e} k_1\!+\!k_2\neq0,\quad
 k_1\!+\!2l_1\cong\deg V_0|_{\P^1_1}\!+\!1~~\tn{mod}~2,
\quad k_2\!+\!2l_2\cong\deg V_0|_{\P^1_2}\!+\!1~~\tn{mod}~2\,.\EE
For $r\!=\!1,2$, we denote the restriction of $(V_0,\vph_0)$ to $\P^1_r$ by $(V_r,\vph_r)$.
Let $\fo_{\nod}$ be an orientation of $V_0^{\vph_0}|_{\nod}$.\\

Let $D_0$ be a real CR-operator on $(V_0,\vph_0)$.
We denote its restrictions to real CR-operators on $(V_1,\vph_1)$ and
$(V_2,\vph_2)$ by~$D_1$ and~$D_2$, respectively.
The orientation~$\fo_{\nod}$ determines orientations
\BE{foRnd_e1} 
\fo_{0;1}\big(\fo_{\nod}\big)\equiv\fo_{k_1,l_1}\!\big(V_1,\vph_1;\fo_{\nod}\big)
\quad\hbox{and}\quad
\fo_{0;2}\big(\fo_{\nod}\big)\equiv\fo_{k_2,l_2}\!\big(V_2,\vph_2;\fo_{\nod}\big)\EE
of $D_1$ and~$D_2$ as above Proposition~\ref{cOpincg_prp}.
These two orientations and the orientation~$\fo_{\nod}$ induce an orientation 
\BE{cOlsplor_e} \fo_0(\fo_{\nod})\equiv 
\fo_{(k_1,l_1),(k_2,l_2)}\big(V_0,\vph_0;\fo_{\nod}\big)\EE
of~$D_0$ via the isomorphism~\eref{Rndses_e2}.
We call it the \sf{split orientation} of~$D_0$.\\

Suppose in addition~\eref{cUsymmdfn_e} is a flat family of deformations of~$\cC_0$,
$(V,\vph)$ is a real bundle pair over $(\cU,\wt\fc)$ extending~$(V_0,\vph_0)$,
$s_1^{\R}$ is a section of $\cU^{\wt\fc}$ over~$\De_{\R}$ with 
$$s_1^{\R}(0)\in S^1_1\!-\!\{\nod\},$$
and  \hbox{$\cD\!\equiv\!\{D_{\t}\}$} is a family of real CR-operators on
$(V_{\t},\vph_{\t})$ as in~\eref{DVvphdfn_e} extending~$D_0$.
The decorated structure on~$\cC_0$ induces a decorated structure on
the fiber $(\Si_{\t},\si_{\t})$ of~$\pi$ for every $\t\!\in\!\De_{\R}$ 
and determines open subspaces $\De_{\R}^{\pm}\!\subset\!\De_{\R}$
as above the CROrient~\ref{CROdegenH3_prop} property.
The orientation~$\fo_{\nod}$ determines an orientation $\fo_{\nod;0}$ of 
$V^{\vph}$ at $s_1^{\R}(0)$ by translation along the positive direction of  $S^1_1\!\subset\!\bD^2_{1+}$
and thus an orientation  $\fo_{\nod;\t}$ of $V^{\vph}$ at $s_1^{\R}(\t)$ for 
every $\t\!\in\!\De_{\R}$ via the line bundle $s_1^{\R*}V^{\vph}$.
For each $\t\!\in\!\De_{\R}^*$, let
\BE{foRnd_e} 
\fo_{\t}(\fo_{\nod})\equiv 
\fo_{k_1+k_2-1,l_1+l_2}\big(V_{\t},\vph_{\t};\fo_{\nod;\t}\big)\EE
be the intrinsic orientation of $D_{\t}$ as  above Proposition~\ref{cOpincg_prp}.
We denote~by 
\BE{cOlimor_e} \fo_0^+(\fo_{\nod})\equiv 
\fo_{k_1+k_2-1,l_1+l_2}^+\big(V_0,\vph_0;\fo_{\nod}\big)\EE
the orientation of~$D_0$ obtained as the continuous extension of the orientations~\eref{foRnd_e}
with $\t\!\in\!\De_{\R}^+$.
We call~\eref{cOlimor_e} the \sf{limiting orientation} of~$D_0$.

\begin{prp}\label{cOlim_prp}
Suppose $\cC_0$ is a connected decorated symmetric surface consisting of two copies
of $(S^2,\tau)$ and one node~$\nod$
and $(V_0,\vph_0)$  is a rank~1 real bundle pair over~$(\Si_0,\si_0)$.
Let $k_1,l_1,k_2,l_2\!\in\!\Z^{\ge0}$ be so that \eref{k1l1k2l2_e} holds
and $D_0$ be a real CR-operator on~$(V_0,\vph_0)$.
The split and limiting orientations, \eref{cOlsplor_e} and~\eref{cOlimor_e}, of~$D_0$
are the same if and only~if
\BE{cOlimcond_e}\big(\!\deg(V_0|_{\P^1_1})\!+\!1\big)\deg(V_0|_{\P^1_2})\in2\Z.\EE
\end{prp}

\begin{proof}[{\bf{\emph{Proof of Proposition~\ref{cOpincg_prp}\ref{cOpincg_it3}}}}]
In light of Lemma~\ref{cOpincg_lmm}\ref{cOorientcomp_it}, it remains 
to show that the orientation $\fo_{k,l}(V,\vph;\fo_{x_1})$ of~$D$
does not depend on the choice of $k\!\in\!\Z^{\ge0}$ satisfying~\eref{pinchcond_e}.
We assume that $k\!\ge\!2$ and show below that the orientations $\fo_{k,l}(V,\vph;\fo_{x_1})$
and $\fo_{k-2,l+1}(V,\vph;\fo_{x_1})$ are the~same.\\

We can choose $\cC_0$, $(\cU,\wt\fc)$, $s_1^{\R}$, $(V,\vph)$, and $\cD\!\equiv\!\{D_{\t}\}$ 
as above Proposition~\ref{cOlim_prp}  and just above so that \eref{3compcUdfn_e} 
holds for some $\t_1\!\in\!\De_{\R}^+$, $s_1^{\R}(\t_1)\!=\!x_1$,
$$(V_2,\vph_2)\equiv (V,\vph)|_{\P^1_2}$$ 
is a degree~2 rank~1 real bundle pair,
$(V_{\t_1},\vph_{\t_1})$ is the rank~1 real bundle pair 
in the statement of Proposition~\ref{cOpincg_prp}, and $D_{\t_1}\!=\!D$.
With the notation as before, let $\fo_{\nod}$ be the orientation of~$V_0^{\vph_0}|_{\nod}$
so that $\fo_{\nod;\t_1}\!=\!\fo_{x_1}$.
By Proposition~\ref{cOlim_prp}, the orientations 
\BE{cOpincg_it3_e3a}\fo_0(\fo_{\nod})\equiv \fo_{(k-2,l),(3,0)}\big(V_0,\vph_0;\fo_{\nod}\big)
\quad\hbox{and}\quad
\fo_0^+(\nod)\equiv \fo_{k,l}^+\big(V_0,\vph_0;\fo_{\nod}\big)\EE
of~$D_0$ are the same if and only if the orientations 
\BE{cOpincg_it3_e3b}\fo_0'(\fo_{\nod})\equiv \fo_{(k-2,l),(1,1)}\big(V_0,\vph_0;\fo_{\nod}\big)
\quad\hbox{and}\quad
\fo_0'^+(\fo_{\nod})\equiv \fo_{k-2,l+1}^+\big(V_0,\vph_0;\fo_{\nod}\big)\EE
of~$D_0$ are the same.\\

By Proposition~\ref{cOpincg_prp}\ref{cOpincg_it2} and Lemma~\ref{cOpincg_lmm2}, 
$$\fo_{3,0}\big(V_2,\vph_2;\fo_{\nod}\big)=\fo_{1,1}\big(V_2,\vph_2;\fo_{\nod}\big).$$
Thus, the split orientations $\fo_0(\nod)$ and $\fo_0'(\nod)$ above are the same.
Along with~\eref{cOpincg_it3_e3a} and~\eref{cOpincg_it3_e3b},
this implies that the limiting orientations $\fo_0^+(\nod)$ and $\fo_0'^+(\nod)$ of~$D_0$
are also the same and~so 
$$\fo_{k,l}\big(V_{\t},\vph_{\t};\fo_{\nod;\t}\big)
=\fo_{k-2,l+1}\big(V_{\t},\vph_{\t};\fo_{\nod;\t}\big)
\qquad\forall~\t\!\in\!\De_{\R}^+\,.$$
Taking $\t\!=\!\t_1$, we obtain the claim.
\end{proof}

\begin{lmm}\label{cOlim_lmm}
Proposition~\ref{cOlim_prp} holds if the two congruences in~\eref{k1l1k2l2_e}
are equalities.
\end{lmm}

\begin{proof}[{\bf{\emph{Proof of Proposition~\ref{cOlim_prp}}}}]
Let 
\BE{cOlim_e0}k=k_1\!+\!k_2\!-\!1, \quad l=l_1\!+\!l_2, \quad
a_1=k_1\!+\!2l_1\!-\!1, \quad  a_2=k_2\!+\!2l_2\!-\!1.\EE
By~\eref{k1l1k2l2_e}, the numbers
$$a_1'\equiv \big(\!\deg V_0|_{\P_1^1}\!-\!a_1\big)\big/2 \qquad \hbox{and}\qquad
a_2'\equiv \big(\!\deg V_0|_{\P_2^1}\!-\!a_2\big)\big/2$$
are integers.
Let $\cC_0'\!\equiv\!(\Si_0',\si_0')$
be a connected decorated symmetric surface consisting of two copies, $\cC_{10}$ and~$\cC_{20}$,
of the three-component curve as on the left side of Figure~\ref{P1CH3split_fig}
joined at a real node~$\nod_0$;
see the bottom diagram in Figure~\ref{cOlimPf_fig}.
For $r\!=\!1,2$, we denote by $\P^1_{r0}\!\subset\!\Si_{r0}$ the irreducible component 
in the $r$-th copy preserved by the involution $\si_0'$,
by $\bD^2_{r+},\bD^2_{r-}\!\subset\!\P^1_{r0}$ the two half-surfaces,
by $\nod_r^+\!\in\!\bD^2_{r+}$ the node on~$\bD^2_{r+}$, and
by $\P^1_{r+}$ the conjugate component of~$\Si_0'$ also containing~$\nod_r^+$.\\

Let  $(\cU',\wt\fc')$ be a flat family of deformations of~$\cC_0'$ as in~\eref{cUsymmdfn_e}
over the unit ball $\De'\!\subset\!\C^5$ around the origin so~that
$\cC_{\t_0}'\!\equiv\!(\Si_{\t_0}',\si_{\t_0}')$ for some $\t_0\!\in\!\De_{\R}'$ 
is the two-component decorated curve~$\cC_0$ of Proposition~\ref{cOlim_prp}.
Let $\De_{\R}'^*\!\subset\!\De_{\R}'$ be the subspace parametrizing smooth symmetric surfaces
and \hbox{$\De_{\R}'^+,\De_{\R}'^-\!\subset\!\De_{\R}'^*$} be its two topological components 
distinguished as above~\eref{Rdegses_e}.
We denote by $\De_{\R;1}'^*\!\subset\!\De_{\R}'$ and $\De_{\R;4}'^*\!\subset\!\De_{\R}'$
the subspaces parametrizing symmetric surfaces with one real node 
and two pairs of conjugate nodes, respectively, and by~$\De_{\R;1}'$ and~$\De_{\R;4}'$
their closures.
Let $s_1'^{\R}$ be a section of $\cU'^{\wt\fc'}$ over~$\De_{\R}'$ so~that 
$s_1'^{\R}(\t)$ is not the real node of~$\Si_{\t}'$ for all $\t\!\in\!\De_{\R;1}'$ and 
$s_1'^{\R}(\t_0)\!\in\!S^1_1$.\\

We can choose a rank~1 real bundle pair $(\wt{V}',\wt\vph')$  over $(\cU',\wt\fc')$
so~that 
$$\begin{aligned}
\deg\wt{V}'\big|_{\P^1_{10}}&=a_1,&\quad  \deg\wt{V}'\big|_{\P^1_{20}}&=a_2, \\
\deg\wt{V}'\big|_{\P^1_{1+}}&=a_1',&\quad  \deg\wt{V}'\big|_{\P^1_{2+}}&=a_2',
\end{aligned}\qquad
\big(\wt{V}'_{\t_0},\wt\vph'_{\t_0}\big)=(V_0,\vph_0).$$
The orientation~$\fo_{\nod}$ of $V_0^{\vph_0}$ at the node~$\nod$ of~$\Si_{\t_0}'$
induces an orientation~$\fo_{\nod_{\t}}$ of $\wt{V}'^{\wt\vph'}$ at
the real node~$\nod_{\t}$ of each symmetric surface~$\Si_{\t}'$
with $\t\!\in\!\De_{\R;1}'$.
It also determines an orientation of the line bundle 
$s_1'^{\R*}\wt{V}'^{\wt\vph'}$ over~$\De_{\R}'$ as above Proposition~\ref{cOlim_prp}.
We denote the restriction of this orientation to the fiber over $\t\!\in\!\De_{\R}'$ 
by~$\fo_{\nod;\t}$.\\

\begin{figure}
\begin{pspicture}(-1,-6)(10,3.5)
\psset{unit=.4cm}
\pscircle(8,0){2}
\pscircle(5.53,2.47){1.5}\pscircle(5.53,-2.47){1.5}
\pscircle(10.47,2.47){1.5}\pscircle(10.47,-2.47){1.5}
\pscircle*(6.59,1.41){.2}\pscircle*(6.59,-1.41){.2}
\pscircle*(9.41,1.41){.2}\pscircle*(9.41,-1.41){.2}
\rput(5.5,2){\sm{$\nod_{\t;1}^+$}}\rput(5.7,-2){\sm{$\nod_{\t;1}^-$}}
\rput(10.4,2){\sm{$\nod_{\t;2}^+$}}\rput(10.3,-2){\sm{$\nod_{\t;2}^-$}}
\psline{<->}(3,4)(3,-4)\rput(2.2,0){$\si_{\t}'$}
\rput(6,4.6){$\P^1_{\t;1+}$}\rput(6,-4.7){$\P^1_{\t;1-}$}
\rput(11,4.6){$\P^1_{\t;2+}$}\rput(11,-4.7){$\P^1_{\t;2-}$}
\rput(11,0){$\P^1_{\t;0}$}
\psarc[linewidth=.03](8,3.46){4}{240}{300} 
\psarc[linewidth=.03,linestyle=dashed](8,-3.46){4}{60}{120} 
\pscircle(28,0){2.5}\psarc[linewidth=.03](28,4.33){5}{240}{300} 
\psarc[linewidth=.03,linestyle=dashed](28,-4.33){5}{60}{120} 
\pscircle(33,0){2.5}\psarc[linewidth=.03](33,4.33){5}{240}{300} 
\psarc[linewidth=.03,linestyle=dashed](33,-4.33){5}{60}{120} 
\pscircle*(30.5,0){.2}
\rput(30.7,3){\sm{$\nod_{\t}$}}\psline[linewidth=.03]{->}(30.5,2.3)(30.5,.7)
\psline{<->}(36,4)(36,-4)\rput(36.8,0){$\si_{\t}'$}
\psline{->}(18,-8)(18,6)\psline{->}(18,-8)(35,-8)
\rput(33,-9){\sm{$\De_{\R;1}'^*\!\subset\!\C^4$}}
\rput{90}(17,4){\sm{$\De_{\R;4}'^+\!\subset\!\C$}}\rput(25,-5){$\De_{\R}'^+$}
\psline[linewidth=.03]{->}(30.5,-2)(30.5,-7)
\psline[linewidth=.03]{->}(12.5,0)(17,0)
\pscircle*(18,-8){.2}\pscircle*(22,-8){.2}
\rput(28.3,3.2){$\P^1_{\t;1}$}\rput(33.3,3.2){$\P^1_{\t;2}$}
\pscircle*(22,-4){.2}
\rput(22.1,-7.2){$\t_0$}
\pscircle(22,2){3}
\psarc[linewidth=.03](22,7.2){6}{240}{300} 
\psarc[linewidth=.03,linestyle=dashed](22,-3.2){6}{60}{120} 
\psline[linewidth=.03]{->}(22,-1.5)(22,-3.5)
\pscircle(6.5,-11){1.5}\pscircle(4.38,-8.88){1.5}\pscircle(4.38,-13.12){1.5}
\pscircle*(5.44,-9.94){.2}\pscircle*(5.44,-12.06){.2}
\psarc[linewidth=.03](6.5,-8.41){3}{240}{300} 
\psarc[linewidth=.03,linestyle=dashed](6.5,-13.59){3}{60}{120} 
\rput(4.7,-9.2){\sm{$\nod_1^+$}}\rput(4.4,-12.7){\sm{$\nod_1^-$}}
\rput(7,-9){\sm{$\bD^2_{1+}$}}\rput(7,-13.2){\sm{$\bD^2_{1-}$}}
\rput(2.7,-7.1){$\P^1_{1+}$}\rput(2.8,-14.3){$\P^1_{1-}$}
\rput(4.3,-11){\sm{$\P^1_{10}$}}
\pscircle*(8,-11){.2}
\rput(8.2,-7){\sm{$\nod_0$}}\psline[linewidth=.03]{->}(8,-7.5)(8,-10.5)
\pscircle(9.5,-11){1.5}\pscircle(11.62,-8.88){1.5}\pscircle(11.62,-13.12){1.5}
\pscircle*(10.56,-9.94){.2}\pscircle*(10.56,-12.06){.2}
\psarc[linewidth=.03](9.5,-8.41){3}{240}{300} 
\psarc[linewidth=.03,linestyle=dashed](9.5,-13.59){3}{60}{120}
\rput(11.5,-9.3){\sm{$\nod_2^+$}}\rput(11.6,-12.5){\sm{$\nod_2^-$}}
\rput(9.2,-8.8){\sm{$\bD^2_{2+}$}}\rput(9.2,-13.2){\sm{$\bD^2_{2-}$}}
\rput(13.7,-7.5){$\P^1_{2+}$}\rput(13.7,-14.5){$\P^1_{2-}$}
\rput(11.8,-11){\sm{$\P^1_{20}$}}
\psline{<->}(1,-8)(1,-14)\rput(.2,-11){$\si_0'$}
\rput(8,-15){$\cC_0'$}
\psline[linewidth=.03]{->}(13,-11)(17.3,-8.5)
\end{pspicture}
\caption{Deformations of~$\cC_0'$ in the family $(\cU',\wt\fc')$ over 
$\De'\!\subset\!\C^5$ of the proof of Lemma~\ref{cOevComp_lmm}}
\label{cOlimPf_fig}
\end{figure}

Let \hbox{$\cD'\!\equiv\!\{D_{\t}'\}$} be a family of real CR-operators on 
$(\wt{V}_{\t}',\wt\vph_{\t}')$ 
as in~\eref{DVvphdfn_e} so that $D_{\t_0}'$ is the operator~$D_0$ in the statement of 
the proposition.
For $r\!=\!1,2$, we denote the restriction of $(\wt{V}_0',\wt\vph_0')$ 
to $\P^1_{r0}$ by $(\wt{V}_{0r}',\vph_{0r}')$ and 
the restrictions of~$D_0'$ to real CR-operators on $(\wt{V}_{0r}',\vph_{0r}')$
and $\wt{V}_0'|_{\P^1_{r+}}$ by~$D_{0 r}'$ and~$D_r'^+$, respectively. 
The exact~triple 
\begin{gather*}
0\lra D'_0 \lra 
D'_{01}\!\oplus\!D'_{02}\!\oplus\!D_1'^+\!\oplus\!D_2'^+\!
\lra \wt{V}'^{\wt\vph'}_{\nod_0}\!\oplus\!\wt{V}'_{\nod_1^+}\!\oplus\!\wt{V}'_{\nod_2^+}\lra 0, \\
\notag
\big(\xi_{1-},\xi_{01},\xi_{1+},\xi_{2-},\xi_{02},\xi_{2+}\big)\lra
\big(\xi_{01},\xi_{1+},\xi_{02},\xi_{2+}\big),\\
\notag
\big(\xi_{01},\xi_{1+},\xi_{02},\xi_{2+}\big)\lra 
\big(\xi_{02}(\nod_0)\!-\!\xi_{01}(\nod_0),\xi_{1+}(\nod_1^+)\!-\!\xi_{01}(\nod_1^+),
\xi_{2+}(\nod_2^+)\!-\!\xi_{02}(\nod_2^+)\big),
\end{gather*}
of Fredholm operators then determines an isomorphism
$$\la(D_0')\!\otimes\!\la\big(\wt{V}'^{\wt\vph'}_{\nod_0}\big)
\!\otimes\!\la\big(\wt{V}'_{\nod_1^+}\big)\!\otimes\!
\la\big(\wt{V}'_{\nod_2^+}\big)
\approx  \la\big(D'_{01}\big)\!\otimes\!\la\big(D'_{02}\big)
\!\otimes\!\la\big(D_1'^+\big)\!\otimes\!\la\big(D_2'^+\big).$$
The evaluation orientations $\fo_{k_1,l_1}(\fo_{\nod_0})$ of $D'_{01}$
and  $\fo_{k_2,l_2}(\fo_{\nod_0})$ of $D'_{02}$,
the orientation~$\fo_{\nod_0}$ of $\wt{V}'^{\wt\vph'}_{\nod_0}$, 
and the complex orientations of 
$D_1'^+$, $D_2'^+$, $\wt{V}'_{\nod_1^+}$, and 
$\wt{V}'_{\nod_2^+}$ determine an orientation of $D'_0$
via the above isomorphism  and 
thus an orientation~$\fo_{\cU'}(\fo_{\nod})$ of the line bundle~$\la(\cD')$ over~$\De'_{\R}$.\\

Let $\De_{\R;4}'^+\!\subset\!\De_{\R;4}'^*$ be the intersection of $\De_{\R;4}'^*$
with the closure of~$\De_{\R}'^+$ in~$\De_{\R}'$.
The restriction of $(\cU',\wt\fc')$ to $\De_{\R;4}'$ is the product of
two conjugates pairs of~$\P^1$ with a family of symmetric surfaces as above 
Proposition~\ref{cOlim_prp} smoothing 
$$\P^1_{0;0}\equiv \P^1_{10}\!\cup_{\nod_0}\!\P^1_{20}$$
into 
a single irreducible component $\P^1_{\t;0}\!\subset\!\Si_{\t}'$ for each $\t\!\in\!\De_{\R;4}'^*$;
see the top left diagram in Figure~\ref{cOlimPf_fig}.
The other irreducible components of~$\Si_{\t}'$ with $\t\!\in\!\De_{\R;4}'^*$ 
and its nodes correspond to the conjugate components and conjugate nodes of~$\Si_0'$;
we denote them by~$\P^1_{\t;r\pm}$ and~$\nod_{\t;r}^{\pm}$, respectively.
For each $\t\!\in\!\De_{\R;4}'^*$,  we denote the restriction of $(\wt{V}_{\t}',\wt\vph_{\t}')$ 
to $\P^1_{\t;0}$ by $(\wt{V}_{\t;0}',\vph_{\t;0}')$ and 
the restrictions of~$D_{\t}'$ to real CR-operators on $(\wt{V}_{\t;0}',\vph_{\t;0}')$
$\wt{V}_{\t}'|_{\P^1_{\t;r+}}$ by~$D_{\t;0}'$ and~$D_{\t;r}'^+$, respectively. 
The exact triple~\eref{Rndses_e} of Fredholm operators with~$D_0$ replaced 
by~$D_{0;0}'$ and 
the exact triple~\eref{ses5comp_e} of Fredholm operators with~$D_0$ replaced
by~$D_{\t}'$ with $\t\!\in\!\De_{\R;4}'^*$ induce isomorphisms
\BE{cOlim_e7}\begin{split}
\la\big(D_{0;0}'\big)\!\otimes\!\la\big(\wt{V}_{\nod_0}'^{\wt\vph'}\big)
&\approx\la\big(D_{01}'\big)\!\otimes\!\la\big(D_{02}'\big),\\
\la(D_{\t}')\!\otimes\!
\la\big(\wt{V}_{\nod_{\t;1}^+}\big)\!\otimes\!\la\big(\wt{V}_{\nod_{\t;2}^+}\big)
&\approx 
\la\big(D_{\t;0}'\big)\!\otimes\!\la\big(D_{\t;1}'^+\big)\!\otimes\!\la\big(D_{\t;2}'^+\big),
\end{split}\EE
respectively.
The orientations $\fo_{k_1,l_1}(\fo_{\nod_0})$ of $D'_{01}$,
$\fo_{k_2,l_2}(\fo_{\nod_0})$ of $D'_{02}$, and
$\fo_{\nod_0}$ of $\wt{V}_{\nod_0}'^{\wt\vph'}$
determine an orientation of $D_{0;0}'$
via  the first isomorphism in~\eref{cOlim_e7}  and
thus an orientation~$\fo_{\t;0}'(\fo_{\nod})$ of $D_{\t;0}'$
for each $\t\!\in\!\De_{\R;4}'^*$.\\

Let for $\t\!\in\!\De_{\R;4}'^+$.
By Proposition~\ref{cOpincg_prp}\ref{cOpincg_it2} and Lemma~\ref{cOlim_lmm},
the last orientation is the same as 
the evaluation orientation $\fo_{k,l}(\fo_{\nod;\t})$
if and only if~\eref{cOlimcond_e} holds.
Since the indices of~$D_1'^+$ and~$D_2'^+$ and 
the (real) dimensions of $\wt{V}'_{\nod_1^+}$ and $\wt{V}'_{\nod_2^+}$ are even,
this implies that the second isomorphism in~\eref{cOlim_e7} with $\t\!\in\!\De_{\R;4}'^+$
is orientation-preserving with respect to 
the restriction of the orientation $\fo_{\cU'}(\fo_{\nod})$ to $\la(D_{\t}')$,
the orientation $\fo_{k,l}(\fo_{\nod;\t})$ of $D_{\t;0}'$, 
and the complex orientations of the remaining factors
if and only if~\eref{cOlimcond_e} holds.
Combining this with the definition of the limiting orientation~\eref{cOlimor_e}, 
we conclude that the restriction of the orientation $\fo_{\cU'}(\fo_{\nod})$ to $\la(D_{\t_0}')$
is the limiting orientation~\eref{cOlimor_e} if and only if \eref{cOlimcond_e} holds.\\

The restriction of $(\cU',\wt\fc')$ to $\De_{\R;1}'$ is the product of
families of symmetric surfaces as above Proposition~\ref{cOpincg_prp} smoothing 
the three-component curves
$$\P^1_{0;1}\equiv \Si_{10}\equiv \P^1_{1-}\!\cup\! \P^1_{10}\!\cup\!\P^1_{1+}
\qquad\hbox{and}\qquad
\P^1_{0;2}\equiv \Si_{20} \equiv \P^1_{2-}\!\cup\! \P^1_{20}\!\cup\!\P^1_{2+}$$
into irreducible components $\P^1_{\t;1}\!\subset\!\Si_{\t}'$ and 
$\P^1_{\t;2}\!\subset\!\Si_{\t}'$, respectively;
see the right diagram in Figure~\ref{cOlimPf_fig}.
For $\t\!\in\!\De_{\R;1}'$ and $r\!=\!1,2$, 
we denote the restriction of $(\wt{V}_{\t}',\wt\vph_{\t}')$ 
to $\P^1_{\t;r}$ by $(\wt{V}_{\t;r}',\vph_{\t;r}')$ and 
the restriction of~$D_{\t}'$ to a real CR-operator on $(\wt{V}_{\t;r}',\vph_{\t;r}')$
by~$D_{\t;r}'$. 
The exact triple~\eref{Cndses_e} of Fredholm operators with~$D_0$ replaced 
by~$D_{0;r}'$ for $r\!=\!1,2$ and 
the exact triple~\eref{Rndses_e} of Fredholm operators with~$D_0$ replaced
by~$D_{\t}'$ with $\t\!\in\!\De_{\R;1}'^*$ induce isomorphisms
\BE{cOlim_e9}
\la\big(D_{0;r}'\big)\!\otimes\!\la\big(\wt{V}_{\nod_r^+}'\big)
\approx\la\big(D_{0r}'\big)\!\otimes\!\la\big(D_r'^+\big),\quad
\la(D_{\t}')\!\otimes\!\la\big(\wt{V}^{\wt\vph}_{\nod_{\t}}\big)
\approx 
\la\big(D_{\t;1}'\big)\!\otimes\!\la\big(D_{\t;2}'\big),\EE
respectively.
The evaluation orientation $\fo_{k_r,l_r}(\fo_{\nod_0})$ of $D'_{0r}$
and the complex orientations of $D_r'^+$ and $\wt{V}'_{\nod_r^+}$
determine an orientation on $D_{0;r}'$
via the first isomorphism in~\eref{cOlim_e9}  and
thus an orientation~$\fo_{\t;r}(\fo_{\nod})$ on $D_{\t;r}'$
for each $\t\!\in\!\De_{\R;1}'^*$.
By the definition of the $(k_r,l_r)$-intrinsic orientation above Proposition~\ref{cOpincg_prp},
$$\fo_{\t;r}\big(\fo_{\nod}\big)=
\fo_{k_r,l_r}\big(\wt{V}_{\t;r}',\wt\vph_{\t;r}';\fo_{\nod_{\t}}\big).$$
Since the dimensions of $\wt{V}'_{\nod_1^+}$ and $\wt{V}'_{\nod_2^+}$ are even,
this implies that the second isomorphism in~\eref{cOlim_e9} with $\t\!\in\!\De_{\R;4}'^*$
is orientation-preserving with respect to 
the restriction of the orientation $\fo_{\cU'}(\fo_{\nod})$ to $\la(D_{\t}')$,
the intrinsic orientations $\fo_{k_r,l_r}(\wt{V}_{\t;r},\wt\vph_{\t;r};\fo_{\nod_{\t}})$
of~$D_{\t;r}'$, 
and the orientation~$\fo_{\nod_{\t}}$ of $\wt{V}^{\wt\vph}_{\nod_{\t}}$.
Combining the $\t\!=\!\t_0$ case of this statement
with the definition of the split orientation above Proposition~\ref{cOlim_prp}, we conclude that 
the restriction of the orientation $\fo_{\cU'}(\fo_{\nod})$ to $\la(D_{\t_0}')$
is the split orientation~\eref{cOlsplor_e}.\\

The conclusions of the last two paragraphs establish Proposition~\ref{cOlim_prp}.
\end{proof}

\begin{proof}[{\bf{\emph{Proof of Lemma~\ref{cOlim_lmm}}}}] 
We continue with the notation in~\eref{cOlim_e0} and above Proposition~\ref{cOlim_prp}.
In this~case,  
$$a_1=\deg V_1 \qquad\hbox{and}\qquad a_2=\deg V_2\,.$$
Since \hbox{$a_1,a_2,a_1+\!a_2\!\ge\!-1$}, 
the real CR-operators $D_1$, $D_2$, $D_0$, 
and~$D_{\t}$ with $\t\!\in\!\De_{\R}^*$ are surjective
and~\eref{Rndses_e2} reduces~to 
\BE{cOlim_0b}\la(\ker D_0)\!\otimes\!\la\big(V_0^{\vph_0}|_{\nod}\big)\approx
\la\big(\!\ker D_1\big)\!\otimes\!\la\big(\!\ker D_2\big)\,.\EE
By definition of the split orientation,
this isomorphism respects the orientations
$\fo_0(\fo_{\nod})$ of $\ker D_0$,  $\fo_{\nod}$ of~$V_0^{\vph_0}|_{\nod}$,
$\fo_{0;1}(\fo_{\nod})$ of $\ker D_1$, and
$\fo_{0;2}(\fo_{\nod})$ of $\ker D_2$.\\ 

Let $s_1^+,\ldots,s_l^+$ be disjoint sections of $\cU\!-\!\cU^{\wt\fc}$ and
$s_1^{\R},\ldots,s_{k+1}^{\R}$ be disjoint sections of $\cU^{\wt\fc}$ over~$\De_{\R}$ 
so~that
\begin{alignat*}{4}
s_i^+(0)&\in\bD^2_{1+}&~~&\forall\,i\!\in\![l_1], &\qquad
x_i\equiv s_i^{\R}(0)&\in S^1_1\!-\!\{\nod\}&~~&\forall\,i\!\in\![k_1],\\
s_{l_1+i}^+(0)&\in\bD^2_{2+}&~~&\forall\,i\!\in\![l_2], &\qquad
x_i\equiv s_i^{\R}(0)&\in S^1_2\!-\!\{\nod\}&~~&\forall\,i\!\in\![k\!+\!1]\!-\!\big[\max(k_1,1)\big],
\end{alignat*}
the points $\nod,x_1,\ldots,x_{k_1}$ of $S^1_1$ are ordered by position with 
respect to~$\bD^2_{1+}$,
and the points $\nod,x_{k_1+1},\ldots,x_{k+1}$ of~$S^1_2$ 
are ordered by position with respect to~$\bD^2_{2+}$.
We denote the complex orientations of 
$$V_1^{\C}\equiv\bigoplus_{i=1}^{l_1}V_0\big|_{s_i^+(0)} \qquad\hbox{and}\qquad
V_2^{\C}\equiv\bigoplus_{i=1}^{l_2}V_0\big|_{s_{l_1+i}^+(0)}$$
by $\fo_1^{\C}$ and $\fo_2^{\C}$, respectively.
By Proposition~\ref{cOpincg_prp}\ref{cOpincg_it2}, 
the orientations~$\fo_{0;1}(\fo_{\nod})$, $\fo_{0;2}(\fo_{\nod})$, and~$\fo_0^+(\fo_{\nod})$ 
in~\eref{foRnd_e1} and~\eref{cOlimor_e} agree with 
the orientations obtained via evaluation homomorphisms as in~\eref{cOevdfn_e} at
\begin{enumerate}[label=$\bu$,leftmargin=*]

\item $k_1$ points of $S^1_1$ and $s_1^+(0),\ldots,s_{l_1}^+(0)$,

\item $k_2$ points of $S^1_2$ and $s_{l_1+1}^+(0),\ldots,s_l^+(0)$, and

\item $k$ points of $S^1_{\t}\!\equiv\!\Si_{\t}^{\si_{\t}}$ 
and $s_1^+(\t),\ldots,s_l^+(\t)$ for $\t\!\in\!\De_{\R}^+$,

\end{enumerate}
respectively.\\

Suppose $a_2\!\in\!2\Z$. 
Thus, $k_2\!\not\in\!2\Z$ and $V_0^{\vph_0}$ is orientable along~$S^1_2$.
Let~$\fo_1^{\R}$ and~$\fo_2^{\R}$ denote the orientations~of
\BE{cOlim_e15}
V_1^{\R}\equiv\bigoplus_{i=1}^{k_1}\!V_0^{\vph_0}\big|_{x_i} \qquad\hbox{and}\qquad
V_2^{\R}\equiv\bigoplus_{i=2}^{k_2}\!V_0^{\vph_0}\big|_{x_{k_1+i}}\EE
obtained by translating the orientation $\fo_{\nod;0}$ of $V_0^{\vph_0}|_{x_1}$
(or equivalently~$\fo_{\nod}$ of $V_0^{\vph_0}|_{\nod}$)
along the positive direction of~$S_1^1$ and 
by translating $\fo_{\nod}$ along either direction of~$S_2^1$, 
respectively.
The evaluation homomorphisms
\begin{equation*}\begin{split}
\big(\!\ker D_1,\fo_{0;1}(\fo_{\nod})\!\big)&\lra 
\big(V_1^{\R},\fo_1^{\R}\big)\!\oplus\!\big(V_1^{\C},\fo_1^{\C}\big),\\
\big(\!\ker D_2,\fo_{0;2}(\fo_{\nod})\!\big)&\lra 
\big(V_2^{\R},\fo_2^{\R}\big)\!\oplus\!\big(V_0^{\vph_0}|_{\nod},\fo_{\nod}\!\big)
\!\oplus\!\big(V_2^{\C},\fo_2^{\C}\big), \\
\big(\!\ker D_0,\fo_0^+(\fo_{\nod})\!\big) &\lra
\big(V_1^{\R},\fo_1^{\R}\big)\!\oplus\!\big(V_2^{\R},\fo_2^{\R}\big)
\!\oplus\!\big(V_1^{\C},\fo_1^{\C}\big)\!\oplus\!\big(V_2^{\C},\fo_2^{\C}\big)
\end{split}\end{equation*}
are then orientation-preserving isomorphisms
(since the number of summands in the second sum in~\eref{cOlim_e15} is even,
the orientation of $\fo_2^{\R}$ in fact does not depend on the choice of~$\fo_{\nod}$).  
Since the real dimension of $V_2^{\C}$ is even,
this and the sentence following~\eref{cOlim_0b} imply that 
\hbox{$\fo_0^+(\fo_{\nod})\!=\!\fo_0(\fo_{\nod})$}.\\

Suppose $a_2\!\not\in\!2\Z$.
Thus, $a_1\!\cong\!k$ mod~2, $k_2\!\in\!2\Z$, 
and the orientation $\fo_2(\fo_{\nod})$ does not depend on the choice of~$\fo_{\nod}$
by Lemma~\ref{cOev_lmm}\ref{cOevorient_it}.
If $k_1\!\neq\!0$, let~$\fo_1^{\R}$ and~$\fo_2^{\R}$ denote the orientations~of
$$V_1^{\R}\equiv\bigoplus_{i=2}^{k_1}\!V_0^{\vph_0}\big|_{x_i} \qquad\hbox{and}\qquad
V_2^{\R}\equiv\bigoplus_{i=1}^{k_2}V_0^{\vph_0}\big|_{x_{k_1+i}}$$
obtained by translating $\fo_{\nod}$ along the positive direction of~$S_1^1$ 
and the positive direction of~$S_2^1$, respectively.
The evaluation homomorphisms
\begin{equation*}\begin{split}
\big(\!\ker D_1,\fo_{0;1}(\fo_{\nod})\!\big)&\lra 
\big(V_0^{\vph_0}|_{\nod},\fo_{\nod}\!\big)\!\oplus\!
\big(V_1^{\R},\fo_1^{\R}\big)\!\oplus\!\big(V_1^{\C},\fo_1^{\C}\big),\\
\big(\!\ker D_2,\fo_{0;2}(\fo_{\nod})\!\big)&\lra 
\big(V_2^{\R},\fo_2^{\R}\big)\!\oplus\!\big(V_2^{\C},\fo_2^{\C}\big), \\
\big(\!\ker D_0,\fo_0^+(\fo_{\nod})\!\big) &\lra
\big(V_1^{\R},\fo_1^{\R}\big)\!\oplus\!\big(V_2^{\R},\fo_2^{\R}\big)
\!\oplus\!\big(V_1^{\C},\fo_1^{\C}\big)\!\oplus\!\big(V_2^{\C},\fo_2^{\C}\big)
\end{split}\end{equation*}
are then orientation-preserving isomorphisms.
Combining this with the sentence following~\eref{cOlim_0b}
and taking into account the negative sign in~\eref{Rndses_e}, 
we conclude that the orientations~$\fo_0^+(\fo_{\nod})$ and~$\fo_0(\fo_{\nod})$ 
of~$D_0$ are the opposite if $k\!\in\!2\Z$ (and thus $a_1\!\in\!2\Z$) and 
the same if  $k\!\not\in\!2\Z$ (and thus $a_1\!\not\in\!2\Z$).\\

If $a_2\!\not\in\!2\Z$ and $k_1\!=\!0$, then $a_1\!\not\in\!2\Z$ and
the orientation of~$V_0^{\vph_0}|_{\nod}$ obtained by translating 
$\fo_{x_1}(\fo_{\nod})$ along the positive direction of~$S^1_1$
is the opposite~$\ov{\fo_{\nod}}$ of~$\fo_{\nod}$.
We take~$V_2^{\R}$ as in~\eref{cOlim_e15} and denote by
$\fo_2^{\R}(\fo_{\nod})$ and $\fo_2^{\R}(\ov{\fo_{\nod}})$ the orientations of~$V_2^{\R}$
obtained by translating $\fo_{\nod}$ and $\ov{\fo_{\nod}}$, respectively,
along the positive direction of~$S_2^1$.
The evaluation homomorphisms
\begin{gather*}
\big(\!\ker D_1,\fo_{0;1}(\fo_{\nod})\!\big)\lra\big(V_1^{\C},\fo_1^{\C}\big), ~~
\big(\!\ker D_2,\fo_{0;2}(\fo_{\nod})\!\big)\lra 
\big(V_0^{\vph_0}|_{\nod},\fo_{\nod}\!\big)\!\oplus\!
\big(V_2^{\R},\fo_2^{\R}(\fo_{\nod})\!\big)\!\oplus\!\big(V_2^{\C},\fo_2^{\C}\big),\\
\big(\!\ker D_0,\fo_0^+(\fo_{\nod})\!\big) \lra
\big(V_2^{\R},\fo_2^{\R}(\ov{\fo_{\nod}})\!\big)
\!\oplus\!\big(V_1^{\C},\fo_1^{\C}\big)\!\oplus\!\big(V_2^{\C},\fo_2^{\C}\big)
\end{gather*}
are then orientation-preserving isomorphisms.
Since the real dimension of $V_2^{\C}$ is even,
this and the sentence following~\eref{cOlim_0b} imply that 
\hbox{$\fo_0^+(\fo_{\nod})\!=\!\fo_0(\fo_{\nod})$}.
\end{proof}

\subsection{Even-degree bundles over $(S^2,\tau)$: construction and properties}
\label{tauspinorient_subs}

Let $(V,\vph)$ be a real bundle pair over $(S^2,\tau)$ and 
$\os\!\in\!\OSpin_{S^2}(V^{\vph})$ be a relative $\OSpin$-structure on
the real vector bundle $V^{\vph}$ over $S^1\!\subset\!S^2$.
This  implies that $V^{\vph}$ is orientable and thus
the degree of~$V$ is even; see \cite[Proposition~4.1]{BHH}.
We show below that~$\os$ determines an orientation~$\fo_{\os}(V,\vph)$ of
a real CR-operator~$D$ on~$(V,\vph)$.\\

Let $\bD^2_+\!\subset\!S^2$ be a half-surface of $(S^2,\tau)$ as before
and $\cC_0$, $(\cU,\wt\fc)$, and $\t_1\!\in\!\De_{\R}^*$ 
be as in and above~\eref{3compcUdfn_e}.
We first suppose~that 
$$n\!=\!\rk\, V\ge 3.$$ 
By Definition~\ref{RelPinSpin_dfn3},
$\os$ and the embedding of $\bD^2_+$ into~$S^2$ thus
determine a homotopy class $\os(\bD^2_+)$  of isomorphisms $V^{\vph}\!\approx\!S^1\!\times\!\R^n$.
A rank~$n$ real bundle pair $(\wt{V},\wt\vph)$ over $(\cU,\wt\fc)$ so~that 
$$\big(\wt{V},\wt\vph\big)\big|_{\P^1_0}=
\big(\P^1_0\!\times\!\C^n,\tau\!\times\!\fc\big) \quad\hbox{and}\quad
\big(\wt{V}_{\t_1},\wt\vph_{\t_1}\big)=(V,\vph)$$
determines a homotopy class of trivializations of~$V^{\vph}$ from
the canonical trivialization of $\wt{V}^{\wt\vph}|_{\Si_0^{\si_0}}$.
Choose such $(\wt{V},\wt\vph)$ so that this homotopy class is~$\os(\bD^2_+)$. 
Such a real bundle pair can be obtained by extending a trivialization
of~$V^{\vph}$ in the homotopy class~$\os(\bD^2_+)$ to a trivialization of~$V$ 
over a $\tau$-invariant tubular neighborhood~$U$ of $S^1\!\subset\!S^2$ 
and then collapsing two circles in~$U$ interchanged by~$\tau$ to the nodes
of~$\Si_0$; see the proof of \cite[Proposition~3.1]{Ge1}.\\

\noindent
Let \hbox{$\cD\!=\!\{D_{\t}\}$} be a family of real CR-operators on
$(\wt{V}_{\t},\wt\vph_{\t})$ as in~\eref{DVvphdfn_e} so~that $D_{\t_1}\!=\!D$
and the real CR-operator $D_{00}$ on $(\P^1_0\!\times\!\C^n,\tau\!\times\!\fc)$
induced by restricting~$D_0$ is the standard $\bp$-operator.
The exact triple~\eref{Cndses_e} of Fredholm operators
then induces an isomorphism as in~\eref{Cndses_e2}.
The canonical identification of the kernel of the surjective operator~$D_{00}$ with~$\R^n$,
the complex orientation of the real CR-operator~$D_0^+$ on $\wt{V}|_{\P^1_+}$ induced
by~$D_0$, and the complex orientation of~$\wt{V}_{\nod_+}$
determine an orientation $\fo_{\os}(\wt{V}_0,\wt\vph_0)$ 
on~$\la(D_0)$ via this isomorphism and thus an orientation on the line bundle $\la(\cD)$
over~$\De_{\R}$.
The latter restricts to an orientation~$\fo_{\os}(V,\vph)$ of~$\la(D)$.\\

\noindent
If $n\!=\!2$ and $\fo$ is an orientation of~$V^{\vph}$, 
the group $\pi_1(S^1)\!\approx\!\Z$ acts freely and transitively on the homotopy classes 
of trivializations of~$(V^{\vph},\fo)$.
The $2\Z$-orbits of this action are the collections of these homotopy classes
that induce the same trivializations of the stabilization
$$\St(V^{\vph}) \equiv \tau_{S^1}\!\oplus\!V^{\vph}\lra S^1\,.$$
A relative $\OSpin$-structure~$\os$ and $\bD^2_+$ thus determine a collection $\os(\bD^2_+)$
of homotopy classes of isomorphisms
$V^{\vph}\!\approx\!S^1\!\times\!\R^2$.
Taking any of these homotopy classes and proceeding as in the $n\!\ge\!3$ case,
we obtain an orientation~$\fo_{\os}(V,\vph)$ of~$D$.
By Corollary~\ref{SpinP1rk1_crl} in Section~\ref{tauspinorient_subs2}, 
$\fo_{\os}(V,\vph)$ does not depend on the homotopy class 
in the collection~$\os(\bD^2_+)$\\

If $n\!=\!1$, the orientation~$\fo$ of the real vector bundle
$$\St^2(V^{\vph}) \equiv \tau_{S^1}\!\oplus\!\tau_{S^1}\!\oplus\!V^{\vph}\lra S^1$$
in the pair $\os(\bD^2_+)\!\equiv\!(\fo,\fs)$ determines
a homotopy class~$\fs_{\os}(V^{\vph})$ of trivializations of the real line bundle~$V^{\vph}$ over~$S^1$,
which does not depend on the choice of~$\bD^2_+$.
Let $x_1\!\in\!S^1$ and $\fo_{x_1}$ be  the orientation of~$V_{x_1}^{\vph}$ determined by~$\fo$.
The above construction for the~pair
\BE{wchosdfn_e}\wch\os=\big(\fo,\fs_{\os}(V^{\vph})\!\big)\EE
determines the intrinsic orientation
$$\fo\big(V,\vph;\fo_{x_1}\big)\equiv\fo_{1,0}\big(V,\vph;\fo_{x_1}\big)$$ 
of Section~\ref{tauorient_subs}.
We~take 
\BE{SpinP1rk1_e}\fo_{\os}(V,\vph)=\begin{cases}\fo(V,\vph;\fo_{x_1}),&
\hbox{if}~\os(\bD^2_+)\!=\!\St_{V^{\vph}}^2(\wch\os);\\
\ov{\fo(V,\vph;\fo_{x_1})},&\hbox{if}~
\os(\bD^2_+)\!\neq\!\St_{V^{\vph}}^2(\wch\os).
\end{cases}\EE
In particular, $\fo_{\os}(V,\vph)$ satisfies the property in 
Proposition~\ref{tauspinorient_prp}\ref{SpinvsCan_it} below.
Along with  Proposition~\ref{cOpincg_prp}\ref{cOpincg_it2},
this implies that it also satisfies the CROrient~\ref{CRONormal_prop}\ref{CROnormSpin_it} 
property on page~\pageref{CROnormSpin_it}.\\

Suppose $(V,\vph)$ is a rank~1 odd-degree real bundle pair over $(S^2,\tau)$.
In particular, the real line bundle~$V^{\vph}$ over $S^1\!\subset\!S^2$ is not orientable.
Let $D$ be a real CR-operator on~$(V,\vph)$.
We denote~by 
$$\fo_0^{\pm}\big((3\!\pm\!1)(V,\vph)\big)\equiv 
\fo^{\pm}_{\io_{S^2}(\os_0((3\pm1)V^{\vph},\fo_V^{\pm}))}\big((3\!\pm\!1)(V,\vph)\big)$$
the orientation on 
\BE{CanSpinOrien_e}\la\big((3\!\pm\!1)D\big)\approx \la(D)^{\otimes(3\pm1)} \EE
determined by the image
$$\io_{S^2}\big(\os_0((3\!\pm\!1)V^{\vph},\fo_V^{\pm})\big)
\in\OSp_{S^2}\big((3\!\pm\!1)V^{\vph}\big)$$
of the $\OSpin$-structure $\os_0((3\!\pm\!1)V^{\vph},\fo_V^{\pm})$ of
Examples~\ref{Pin2Spin_eg} and~\ref{SpinDfn1to3_eg}
under the first map in~\eref{vsSpinPin_e0} with \hbox{$X\!=\!S^2$}.
The right-hand side of~\eref{CanSpinOrien_e} also carries a canonical orientation;
it is obtained by taking any orientation on the first factor of~$\la(D)$ 
and the same orientation on the other factor(s) of~$\la(D)$.
Since
$$\os_0(4V^{\vph},\fo_V^+)
=\bllrr{\os_0(2V^{\vph},\fo_V^-),\os_0(2V^{\vph},\fo_V^-)}_{\oplus},$$
the orientation $\fo_0^+(4(V,\vph))$ is the canonical orientation of~$\la(4D)$.
By \cite[Proposition~3.5]{RealGWsII}, the orientation $\fo_0^-(2(V,\vph))$
is also the canonical orientation of~$\la(2D)$;
for the sake of completion, we recall one of the proofs of this fact in~\cite{RealGWsII}
at the end of this section.
We thus obtain the following.

\begin{prp}\label{tauspinorient_prp}
The orientations $\fo_{\os}(V,\vph)$ of the determinants of real CR-operators
on even-degree real bundle pairs $(V,\vph)$  over~$(S^2,\tau)$ constructed 
 above satisfy the CROrient~\ref{CROos_prop}\ref{osflip_it},
\ref{CROSpinPinStr_prop}\ref{CROSpinStr_it}, \ref{CROSpinPinSES_prop}\ref{CROsesSpin_it},
and~\ref{CRONormal_prop}\ref{CROnormSpin_it} properties of Section~\ref{OrientPrp_subs1}.
If $(V,\vph)$ is a rank~1 real bundle pair over~$(S^2,\tau)$ and
\begin{enumerate}[label=(\arabic*),leftmargin=*]

\item\label{SpinvsCan_it} 
$\fo$ is an orientation on~$V^{\vph}$, then the orientation $\fo_0(V,\vph;\fo)$ as in~\eref{osLBdfn_e}
is the intrinsic orientation $\fo(V,\vph;\fo_{x_1})$ of 
Proposition~\ref{cOpincg_prp} for the restriction 
of~$\fo$ to $V^{\vph}_{x_1}$ with $x_1\!\in\!S^1$;

\item\label{PinCan_it} 
the degree of~$V$ is~1, then the orientation $\fo_0^{\pm}((3\!\pm\!1)(V,\vph))$ 
is the canonical orientation on~\eref{CanSpinOrien_e}.

\end{enumerate}
\end{prp}

\begin{rmk}\label{PinCan_rmk}
By the $(S^2,\tau)$ case of the CROrient~\ref{CROdegenC_prop}\ref{CdegenSpin_it} property
established in Section~\ref{SpinOrient_subs} and
Proposition~\ref{tauspinorient_prp}\ref{PinCan_it},
the conclusion of the latter holds for all rank~1 odd-degree real bundles
$(V,\vph)$ over~$(S^2,\tau)$.
\end{rmk}

The orienting construction described above Proposition~\ref{tauspinorient_prp}
depends on the $\OSpin$-structure $\os(\bD^2_+)$
on the real vector bundle~$V^{\vph}$ induced by the relative $\OSpin$-structure~$\os$, 
rather than~$\os$ itself.
For an $\OSpin$-structure~$\os$ on~$V^{\vph}$, we thus denote by $\fo_{\os}(V,\vph)$
the orientation of~$D$ induced by~$\os$.
The next statement 
plays a crucial role in establishing Proposition~\ref{genspinorient_prp}\ref{genspindfn_it}.

\begin{prp}\label{SpinFOOO_prp}
Suppose $(V,\vph)$ is a real bundle pair  over~$(S^2,\tau)$ and
$\os$ and $\os'$ are relative $\OSpin$-structures on~$V^{\vph}$
inducing the same orientation.
If $\os(\bD^2_+)\!\neq\!\os'(\bD^2_+)$, 
\hbox{$\fo_{\os}(V,\vph)\!\neq\!\fo_{\os'}(V,\vph)$}.
\end{prp}

We prove this proposition and the CROrient~\ref{CROSpinPinSES_prop}\ref{CROsesSpin_it} property
for real bundle pairs over $(S^2,\tau)$ in Section~\ref{tauspinorient_subs2},
{\it without} referring to 
the CROrient~\ref{CROos_prop}\ref{osflip_it} or \ref{CROSpinPinStr_prop}\ref{CROSpinStr_it}
properties.
We instead establish these two properties below using Proposition~\ref{SpinFOOO_prp} and
the CROrient~\ref{CROSpinPinSES_prop}\ref{CROsesSpin_it} property.

\begin{proof}[{\bf{\emph{Proof of CROrient~\ref{CROos_prop}\ref{osflip_it} and
\ref{CROSpinPinStr_prop}\ref{CROSpinStr_it} properties for $(S^2,\tau)$}}}]
Let $(V,\vph)$ be a rank~$n$ real bundle pair over $(S^2,\tau)$,
$\os$ be a relative $\OSpin$-structure on the real vector bundle $V^{\vph}$ over $S^1\!\subset\!S^2$,
and $D$ be a real CR-operator on~$(V,\vph)$.
The collections $\os(\bD^2_+)$ and $\os(\bD^2_-)$ of trivializations 
determined by the two half-surfaces
$\bD^2_+$ and $\bD^2_-$ of~$(S^2,\tau)$ induce the same orientation on~$V^{\vph}$.
By the compatibility condition in Definition~\ref{RelPinSpin_dfn3} 
and~\eref{w2fsdfn_e}, these two classes are the same if and only if $w_2(\os)$ vanishes.
If \hbox{$\os(\bD^2_+)\!=\!\os(\bD^2_-)$},
the same reasoning as in the proof of Proposition~\ref{cOpincg_prp}\ref{cOpincg_it1} implies
that the orientations $\fo_{\os}(V,\vph)$ of~$D$ induced by~$\bD^2_+$ and~$\bD^2_-$ 
are the same if  and only if the number
$$-n+\big((\deg V)/2\!+\!n\big)= \frac{\deg V}{2}$$
is even.
Combining the last two statements with Proposition~\ref{SpinFOOO_prp}, 
we conclude the orientations $\fo_{\os}(V,\vph)$ induced by~$\bD^2_+$ and~$\bD^2_-$ 
are the same if  and only~if
$$\frac{\deg V}{2}+\blr{w_2(\os),[S^2]_{\Z_2}} \equiv \vp_{\os}(S^2)$$
vanishes in~$\Z_2$. 
This establishes the CROrient~\ref{CROos_prop}\ref{osflip_it} property for $(S^2,\tau)$.\\

Let $\eta\!\in\!H^2(S^2,S^1;\Z_2)$.
The relative $\OSpin$-structures $\os$ and $\eta\!\cdot\!\os$ determine
the same orientation~$\fo$ on~$V^{\vph}$.
By the RelSpinPin~\ref{RelSpinPinStr_prop} property on page~\pageref{RelSpinPinStr_prop}, 
the collections $\os(\bD^2_+)$ and $\eta\!\cdot\!\os(\bD^2_+)$ of trivializations 
determined by $\os$, $\eta\!\cdot\!\os$, and~$\bD^2_+$
are the same if and only if $\eta|_{\bD^2_+}$ vanishes.
Along with Proposition~\ref{SpinFOOO_prp}, this implies the first statement 
of the CROrient~\ref{CROSpinPinStr_prop}\ref{CROSpinStr_it} property for $(S^2,\tau)$.\\

The collections $\os(\bD^2_+)$ and $\ov\os(\bD^2_+)$ of trivializations determined
by $\os$, $\ov\os$, and~$\bD^2_+$ satisfy
\BE{SpinP1_e5}\ov\os(\bD^2_+) =\ov{\os(\bD^2_+)} .\EE
Suppose first that $n\!=\!1$. 
By~\eref{SpinP1_e5},  
$$\os(\bD^2_+)\!=\!\St_{V^{\vph}}^2(\wch\os) \qquad\Llra\qquad
\ov\os(\bD^2_+)\!=\!\St_{V^{\vph}}^2\big(\wch{\ov\os}\big).$$
Along with~\eref{SpinP1rk1_e} and Proposition~\ref{cOpincg_prp}\ref{cOpincg_it1},
this gives  
\BE{SpinP1_e7} \fo_{\ov\os}(V,\vph)=\ov{\fo_{\os}(V,\vph)}.\EE
This establishes the second statement of the CROrient~\ref{CROSpinPinStr_prop}\ref{CROSpinStr_it} 
property for the rank~1 real bundle pairs over~$(S^2,\tau)$.
Below we deduce the general case from the rank~1 case via 
the CROrient~\ref{CROSpinPinSES_prop}\ref{CROsesSpin_it} property.\\

We denote by $\os_1$ the standard relative $\OSpin$-structure on 
the trivial line bundle $S^1\!\times\!\R$ over $S^1\!\subset\!S^2$ and
by $\bp$ the standard real CR-operator on the real bundle pair 
$(S^2\!\times\!\C,\tau\!\times\!\fs)$.
Let $(L,\phi)$ be a rank~1 real bundle pair over~$(S^2,\tau)$ of even degree,
$\os_L$ be a relative $\OSpin$-structure on the real bundle~$L^{\phi}$ 
over $S^1\!\subset\!S^2$,
and $\bp_L$ be a real CR-operator on~$(L,\phi)$.
For $n\!\in\!\Z^+$, define  
$$(V_n,\vph_n)=(S^2\!\times\!\C^n,\tau\!\times\!\fc), \qquad
(V_{n;L},\vph_{n;L})=\big(V_{n-1},\vph_{n-1}\big)\!\oplus\!(L,\phi).$$
Let $D_n$ be the real CR-operator on $(V_{n;L},\vph_{n;L})$ given by the $n$-fold direct 
sum of the operators~$\bp$ on each factor and \hbox{$D_{n;L}\!=\!D_{n-1}\!\oplus\!\bp_L$}.
Define
\BE{SpinP1_e8}
\os_{n;L}=\St^{n-1}_{L^{\phi}}(\os_L)=\llrr{\os_1,\os_{n-1;L}}_{\oplus}.\EE
By the RelSpinPin~\ref{RelSpinPinStab_prop} and~\ref{RelSpinPinSES_prop} properties 
in Section~\ref{RelSpinPinProp_subs},
$$\ov{\os_{n;L}}=\St^{n-1}_{L^{\phi}}\big(\ov{\os_L}\big)
=\llrr{\os_1,\ov{\os_{n-1;L}}}_{\oplus}.$$
The second equalities in the last two equations hold for $n\!\ge\!2$.\\

Let $n\!\ge\!2$. 
By the second equality in~\eref{SpinP1_e8} and 
the CROrient~\ref{CROSpinPinSES_prop}\ref{CROsesSpin_it} property for $(S^2,\tau)$, 
the natural isomorphism
$$\la(D_{n;L})\approx \la(\bp)\!\otimes\!\la(D_{n-1;L})$$
respects the orientations $\fo_{\os_{n;L}}(V_{n;L},\vph_{n;L})$, 
$\fo_{\os_1}(V_1,\vph_1)$, and $\fo_{\os_{n-1;L}}(V_{n-1;L},\vph_{n-1;L})$.
It also respects the orientations $\fo_{\ov{\os_{n;L}}}(V_{n;L},\vph_{n;L})$, 
$\fo_{\os_1}(V_1,\vph_1)$, and $\fo_{\ov{\os_{n-1;L}}}(V_{n-1;L},\vph_{n-1;L})$.
Along with \eref{SpinP1_e7} with $(V,\vph)\!=\!(V_{1;L},\phi_{1;L})$, this implies~that 
$$\fo_{\ov{\os_{n;L}}}(V_{n;L},\vph_{n;L})= 
\ov{\fo_{\os_{n;L}}(V_{n;L},\vph_{n;L})}\qquad\forall\,n\!\in\!\Z^+.$$
By the RelSpinPin~\ref{RelSpinPinStab_prop} property,
every relative $\OSpin$-structure on the real vector bundle $V_{n;L}^{\vph_{n;L}}$ 
over $S^1\!\subset\!S^2$ equals~$\os_{n;L}$ for some 
relative $\OSpin$-structure $\os_L$ on~$L^{\phi}$.
By \cite[Proposition~4.1]{BHH}, every real bundle pair $(V,\vph)$ over~$(S^2,\tau)$
is isomorphic to $(V_{n;L},\vph_{n;L})$ for some rank~1 real bundle pair $(L,\phi)$ over~$(S^2,\tau)$.
The last three statements imply 
the second statement of the CROrient~\ref{CROSpinPinStr_prop}\ref{CROSpinStr_it} property.
\end{proof}

\begin{proof}[{\bf\emph{Proof of Proposition~\ref{tauspinorient_prp}\ref{PinCan_it} 
for $\fo_0^-(2(V,\vph))$}}]
Three proofs of this statement appear in~\cite{RealGWsII}.
The first one, reproduced below, determines the relevant real holomorphic sections  
and trivializations explicitly.
The second proof uses the comparisons of different orientations
on the moduli spaces of real lines obtained in~\cite{Teh}.
The last one deduces the claim directly from the fixed-edge equivariant
contribution determined in~\cite{Teh}.\\

We first note that there is a canonical isomorphism of the real line bundle
$$L_1^{\R}\!\equiv\!\big\{\big(\ne^{\fI\th},a\ne^{\fI\th/2}\big)\!\in\!S^1\!\times\!\C\!:
a\!\in\!\R\big\}\lra S^1$$
with the real line bundle $\ga_{\R;1}\!\lra\!\R\P^1$.
The composition of the trivialization~\eref{ga1Rpin_e} of $2\ga_{\R;1}$
with the direct sum of two copies of this canonical identification is given~by
\BE{PinCan_e3}\wt\Phi_0\!: 2L_1^{\R}\lra S^1\!\times\!\C, \quad
\wt\Phi_0\big(\!(\ne^{\fI\th},a\ne^{\fI\th/2}),(\ne^{\fI\th},b\ne^{\fI\th/2})\!\big)
=\big(\ne^{\fI\th},(a\!+\!\fI b)\ne^{\fI\th/2}\big).\EE

\vspace{.15in}

Let $S^2_{\bu}\!=\!S^2\!-\!\{1\}$.
The holomorphic map
$$h\!: B\!\equiv\!\big\{t\!\in\!\C\!:\,|t|\!<\!1\big\}\lra\C\subset S^2, \qquad
h(t)=\ne^{\fI t}\,,$$
is injective and intertwines the standard conjugation on~$B$ with $\tau$ on~$S^2$.
We can assume~that $D$ is the standard $\bar\partial$-operator on the holomorphic rank~1 
real bundle pair~$(V,\vph)$ given~by
\begin{gather*}
V =\big(h(B)\!\times\!\C\sqcup S^2_{\bu}\!\times\!\C\big)\big/\!\!\sim, \quad
\big(h(t),tc\big)\sim\big(h(t),c\big)~~\forall~(t,c)\!\in\!\big(B\!-\!\{0\}\!\big)\!\times\!\C,\\
\vph\big([h(t),c]\big)=\big[h(\bar{t}),\bar{c}\big]~~\forall~(t,c)\!\in\!B\!\times\!\C, \quad
\vph\big([z,c]\big)=\big[\tau(z),\bar{c}\big]~~\forall~(z,c)\!\in\!S^2_{\bu}\!\times\!\C.
\end{gather*}
The space $\ker D$ of real holomorphic sections of~$V$ is then generated~by 
the sections~$s_1$ and~$s_2$ described~by
$$s_1(z)=1, \quad s_2(z)=\fI\frac{1\!+\!z}{1\!-\!z} \qquad\forall~z\!\in\!S^2_{\bu}.$$ 
The canonical orientation for $\la(2D)\!=\!\la(\ker 2D)$ is then determined by the basis 
$$s_{11}\equiv (s_1,0), \quad s_{12}\equiv (s_2,0), 
\quad s_{21}\equiv (0,s_1),\quad s_{22}\equiv (0,s_2),$$
for the kernel of the surjective operator~$2D$.\\

The trivialization 
\begin{gather*}
\Psi_2\!:2V|_{S^2_{\bu}-\{0,\i\}}\lra \big(S^2_{\bu}\!-\!\{0,\i\}\big)\!\times\!\C^2, \\
\Psi_2\big([z,c_1],[z,c_2]\big)=\big(z,\fI(z\!-\!z^{-1})c_1\!-\!z^{-1}(1\!-\!z)^2c_2,
z^{-1}(1\!-\!z)^2c_1\!+\!\fI(z\!-\!z^{-1})c_2\big),
\end{gather*}
of $2V|_{S_{\bu}^2-\{0,\i\}}$ extends to a trivialization~$\wt\Psi_2$ of $2V|_{S^2-\{0,\i\}}$.
Since $\Psi_2$ intertwines $2\vph$ with the standard lift $\tau\!\times\!\fc$ of 
$\tau|_{S^2_{\bu}-\{0,\i\}}$ to a conjugation on the trivial bundle
\hbox{$(S^2_{\bu}\!-\!\{0,\i\})\!\times\!\C^2$},
$\wt\Psi_2$ intertwines~$2\vph$ with the conjugation $\tau\!\times\!\fc$ on the trivial bundle
\hbox{$(S^2\!-\!\{0,\i\})\!\times\!\C^2$}.
We note that
\begin{alignat*}{2}
\big\{\wt\Psi_2s_{11}\big\}(z)&=\big(\fI z^{-1}(z^2\!-\!1),z^{-1}(1\!-\!z)^2\big),
&\quad
\big\{\wt\Psi_2s_{12}\big\}(z)&=\big(z^{-1}(1\!+\!z)^2,\fI z^{-1}(1\!-\!z^2)\big),\\
\big\{\wt\Psi_2s_{21}\big\}(z)&=\big(-\!z^{-1}(1\!-\!z)^2,\fI z^{-1}(z^2\!-\!1)\big),
&\quad
\big\{\wt\Psi_2s_{22}\big\}(z)&=\big(-\!\fI z^{-1}(1\!-\!z^2),z^{-1}(1\!+\!z)^2\big).
\end{alignat*}

\vspace{.15in}

We define a trivialization $\Psi$ of~$V$ over $S^2\!-\!\{\i\}\!\supset\!\bD^2_+$ by
\begin{alignat*}{2}
\Psi\big([h(t),c]\big)&=\big(h(t),2\fI\frac{\ne^{\fI t}\!-\!1}{t}c\big) &\qquad
&\forall~(t,c)\!\in\!B\!\times\!\C, \\
\Psi\big([z,c]\big)&=\big(z,2\fI(z\!-\!1)c\big) &\qquad
&\forall~(z,c)\!\in\!\big(S^2_{\bu}\!-\!\{\i\}\big)\!\times\!\C.
\end{alignat*}
This trivialization satisfies $\Psi(V^{\vph})\!=\!L_1^{\R}$.
Under the identification of~$\C$ in the target of~\eref{PinCan_e3} with~$\R^2$ 
in the target of~$\wt\Psi_2$,
the restriction of this trivialization of~$2V|_{S^2-\{0,\i\}}$ to~$2V^{\vph}$ is
the composition of the trivialization~\eref{PinCan_e3} with $\Psi\!\oplus\!\Psi$.
Thus, this restriction lies in the homotopy class~$\os_0(2V^{\vph},\fo_V^-)$ 
of trivializations of~$2V^{\vph}$.
It follows that the orientation~$\fo_0^-(2(V,\vph))$ of~$\la(2D)$ is obtained from
the isomorphism
$$\ker 2D\lra \R\!\oplus\!\R\oplus\big\{\Res_{z=0}(\wt\Psi_2\xi)\!:\,
\xi\!\in\!\ker 2D\big\}, \quad
\xi\lra \big(\{\wt\Psi_2(\xi)\}(1),\Res_{z=0}(\wt\Psi_2(\xi)\!)\!\big).$$
The last space above is a complex subspace of~$\C^2$.\\

Under the above isomorphism, the basis $s_{11},s_{12},s_{21},s_{22}$ is sent~to
$$(0,0;-\fI,1), \quad (4,0;1,\fI), \quad (0,0;-1,-\fI), \quad
(0,4;-\fI,1).$$
Thus, an oriented basis for the target of the above isomorphism is given~by
$$(4,0;0,0), \quad (0,4;0,0), \quad (0,0;-\fI,1), \quad (0,0;1,\fI).$$
The change of basis matrix from the first basis to this one is given~by
$$\left( \begin{array}{cccc} 0& 1& 0& 0\\ 0& 0& 0& 1\\ 1& 0& 0& 1\\
0& 1& -1& 0\end{array}\right).$$
The determinant of this matrix is $+1$.
\end{proof}

\subsection{Even-degree bundles over degenerations of $(S^2,\tau)$ and exact triples}
\label{tauspinorient_subs2}

In this section, we establish Proposition~\ref{SpinFOOO_prp} and 
the CROrient~\ref{CROSpinPinSES_prop}\ref{CROsesSpin_it} property 
on page~\pageref{CROSpinPinSES_prop} for even-degree real bundle pairs over~$(S^2,\tau)$.
The former is immediate from~\eref{SpinP1rk1_e} for rank~1 real bundle pairs;
the latter is straightforward for the short exact sequences of real bundle pairs
of rank at least~2, as noted in Step~1 of the proof of 
the CROrient~\ref{CROSpinPinSES_prop}\ref{CROsesSpin_it} property below.
The next statement is a key ingredient in establishing Proposition~\ref{SpinFOOO_prp}
for real bundle pairs of rank at least~2.

\begin{lmm}\label{SpinFOOO_lmm}
There exist $\OSpin$-structures $\os$ and~$\os'$ on 
the trivial oriented rank~2 vector bundle $S^1\!\times\!\R^2$ over~$S^1$
such~that 
$$\fo_{\os}\big(S^2\!\times\!\C^2,\tau\!\times\!\fc\big)\neq
\fo_{\os'}\big(S^2\!\times\!\C^2,\tau\!\times\!\fc\big)\,.$$
\end{lmm}

\begin{proof}
By \cite[Proposition~8.1.7]{FOOO}, there exist a rank~2 real bundle pair $(V,\vph)$ 
over $S^1\!\times\!(S^2,\tau)$ with $c_1(V)\!=\!0$ and $w_1(V^{\vph})\!=\!0$
and a family 
$\cD\!\equiv\!\{D_s\}_{s\in S^1}$ of real CR-operators on the real bundle pairs
$$(V_s,\vph_s)\equiv (V,\vph)\big|_{\{s\}\times(S^2,\tau)}$$
so that the determinant line bundle~$\la(\cD)$ over $S^1$ is not orientable.
By the construction in~\cite{FOOO} (as well as by \cite[Proposition~4.1]{BHH}),
there exists an automorphism~$\Psi$ of the trivial rank~2 real bundle pair 
$(S^2\!\times\!\C^2,\tau\!\times\!\fc)$ covering the identity on~$S^2$ so~that 
$$(V,\vph)=\big([0,1]\!\times\!(S^2\!\times\!\C^2,\tau\!\times\!\fc)\big)\big/\!\sim,
\qquad \big(0,\Psi(z,v)\big)\sim\big(0,(z,v)\!\big).$$
This automorphism restricts to an automorphism $\Psi_{\R}$ of 
the real vector bundle $S^1\!\times\!\R^2$ over~$S^1$.
Since the vector bundle~$V^{\vph}$ is orientable, $\Psi_{\R}$ is orientation-preserving.\\

Let $\fo$ and $\fs$ be the natural orientation and 
the natural homotopy class of trivializations, respectively, of 
the real vector bundle $S^1\!\times\!\R^2$ over~$S^1$.
Since $\Psi_{\R}$ is orientation-preserving, 
the homotopy class $\fs'\!\equiv\!\Psi_{\R}^*\fs$ preserves the orientation~$\fo$.
Let $\wt\cD\!\equiv\!\{\wt{D}_s\}_{s\in[0,1]}$ be the family of real CR-operators on 
$(S^2\!\times\!\C^2,\tau\!\times\!\fc)$ induced by the family~$\cD$.
Since the orienting construction at the beginning of Section~\ref{tauspinorient_subs}
applied to $\la(\wt{D}_s)$ depends continuously on~$s$, 
it determines a continuously varying orientation~$\fo_s$  of $\la(\wt{D}_s)$ 
from the homotopy class~$\fs$ of trivializations of $S^1\!\times\!\R^2$.
Since the line bundle $\la(\cD)$ is not orientable, the orientations~$\fs_1$ and~$\Psi^*\fs_0$
of the determinant of $\wt{D}_1\!=\!\Psi^*\wt{D}_0$ induced by the
$\OSpin$-structures \hbox{$\os\!\equiv\!(\fo,\fs)$} and $\os'\!\equiv\!(\fo,\fs')$ 
on the trivial rank~2 vector bundle $S^1\!\times\!\R^2$ over~$S^1$ are different.
\end{proof}

We also use Lemma~\ref{3by3_lmm} below.
It is the Exact Squares property for the determinants of Fredholm operators;
its specialization to finite-dimensional vector spaces is a straightforward 
linear algebra observation.

\begin{lmm}[{\cite[(2.27)]{detLB}}]\label{3by3_lmm}
Let $D_{ij}$ with $i,j\!\in\![3]$ be Fredholm operators with orientations 
on their determinants.
If the rows and columns in the diagram in Figure~\ref{3by3_fig} are 
exact triples of Fredholm operators and this diagram commutes,
then the total number of rows and columns in this diagram for which 
the associated isomorphism~\eref{sum} respects the orientations 
is congruent to $\ind(D_{13})\ind(D_{31})$ mod~2.
\end{lmm}

\begin{figure}
$$\xymatrix{& 0\ar[d] & 0\ar[d] & 0\ar[d]  \\
0\ar[r] & D_{11}\ar[r]\ar[d] & D_{12}\ar[r]\ar[d] & D_{13}\ar[r]\ar[d] & 0 \\
0\ar[r] & D_{21}\ar[r]\ar[d] & D_{22}\ar[r]\ar[d] & D_{23}\ar[r]\ar[d] & 0 \\
0\ar[r] & D_{31}\ar[r]\ar[d] & D_{32}\ar[r]\ar[d] & D_{33}\ar[r]\ar[d] & 0 \\
& 0 & 0 & 0}$$ 
\caption{Commutative square of exact rows and columns of Fredholm operators 
or vector spaces for the statement of Lemma~\ref{3by3_lmm}}
\label{3by3_fig}
\end{figure}

\begin{proof}[{\bf{\emph{Proof of CROrient~\ref{CROSpinPinSES_prop}\ref{CROsesSpin_it} property
for $(S^2,\tau)$, Step~1}}}]
Suppose $\ce$ is a short exact sequence of real bundle pairs over~$(S^2,\tau)$ 
as in~\eref{RBPsesdfn_e}, $\os'$ and~$\os''$ are relative $\OSpin$-structures
on the real vector bundles~$V'^{\vph'}$ and~$V''^{\vph''}$ over $S^1\!\subset\!S^2$,
and 
$$\os\equiv\llrr{\os',\os''}_{\ce_{\R}}$$
is the induced relative $\OSpin$-structure on~$V^{\vph}$.
Let $\cC_0$, $(\cU,\wt\fc)$, and $\t_1\!\in\!\De_{\R}^*$ be as in the orienting construction 
at the beginning of Section~\ref{tauspinorient_subs},
$(\wt{V}',\wt\vph')$ and $(\wt{V}',\wt\vph'')$ be real bundle pairs over $(\cU,\wt\fc)$
associated with the collections $\os'(\bD^2_+)$ and $\os''(\bD^2_+)$ of trivializations,
and \hbox{$\cD'\!\equiv\!\{D_{\t}'\}$} and
\hbox{$\cD''\!\equiv\!\{D_{\t}''\}$} be associated families
of real CR-operators.
Let
\BE{P1SpinSES_e3}(\wt{V},\wt\vph)=(\wt{V}',\wt\vph')\!\oplus\!(\wt{V}'',\wt\vph'')
\lra (\cU,\wt\fc)\EE
and $D_{\t}\!\equiv\!D_{\t}'\!\oplus\!D_{\t}''$ 
for each $\t\!\in\!\De_{\R}$.\\

\noindent
The exact triple~\eref{Cndses_e} induces  exact triples of Fredholm operators given
by the rows in the diagram of Figure~\ref{P1SpinSES_fig}.
The splitting~\eref{P1SpinSES_e3} induces exact triples of Fredholm operators given
by the columns in this diagram.
The isomorphisms~\eref{sum} induced by the right column and 
the triple formed by the second summands in the middle column respect 
the complex orientations of all associated determinants.
The isomorphism~\eref{sum} induced by the first summands in the middle column
respects the canonical orientations of the associated determinant lines.
Since the (real) index of $D_0'^+$ is even,
Lemma~\ref{3by3_lmm} implies that the isomorphism~\eref{sum} induced 
by the middle column respects the direct sum orientations.\\

\begin{figure}
$$\xymatrix{& 0\ar[d] & 0\ar[d] & 0\ar[d]  \\
0\ar[r] & D_0'\ar[r]\ar[d] & D_{00}'\!\oplus\!D_0'^+\ar[r]\ar[d] 
& \wt{V}_{\nod^+}'\ar[r]\ar[d] & 0 \\
0\ar[r] & D_0\ar[r]\ar[d] & D_{00}\!\oplus\!D_0^+
\ar[r]\ar[d]  & \wt{V}_{\nod^+}\ar[r]\ar[d] & 0 \\
0\ar[r] & D_0''\ar[r]\ar[d] & D_{00}''\!\oplus\!D_0''^+ \ar[r]\ar[d]
& \wt{V}_{\nod^+}''\ar[r]\ar[d] & 0 \\
& 0 & 0 & 0}$$ 
\caption{Commutative square of exact rows and columns of Fredholm operators 
for the proof of the CROrient~\ref{CROSpinPinSES_prop}\ref{CROsesSpin_it} 
property for $(S^2,\tau)$}
\label{P1SpinSES_fig}
\end{figure}

Suppose that 
\begin{enumerate}[label=(\arabic*),leftmargin=*]

\item\label{P1SpinSES_it1} either $\rk\,V'\!\ge\!2$ or $\rk\,V'\!=\!1$ with 
$\os'(\bD^2_+)\!=\!\St_{V'^{\vph'}}^2(\wch\os')$, and

\item\label{P1SpinSES_it2} either $\rk\,V''\!\ge\!2$ or $\rk\,V''\!=\!1$ with 
$\os''(\bD^2_+)\!=\!\St_{V''^{\vph''}}^2(\wch\os'')$.

\end{enumerate}
The real bundle pair~\eref{P1SpinSES_e3} and the family $\cD\!\equiv\!\{D_{\t}\}$ 
are then associated with the collection $\os(\bD^2_+)$ of trivializations.
The rows in the diagram of Figure~\ref{P1SpinSES_fig} respect 
the orientations $\fo_{\os'}(\wt{V}_0',\wt\vph'_0)$ of~$D_0'$,
$\fo_{\os}(\wt{V}_0,\wt\vph_0)$ of~$D_0$, 
$\fo_{\os''}(\wt{V}_0'',\wt\vph_0'')$ of~$D_0''$,
and the canonical orientations on the remaining terms in this case.
Since the (real) dimension of~$\wt{V}_{\nod^+}'$ is even,
Lemma~\ref{3by3_lmm} implies that the isomorphism~\eref{CROses_e0}
induced by the left column in this diagram respects the orientations 
$\fo_{\os'}(\wt{V}'_0,\wt\vph'_0)$, $\fo_{\os}(\wt{V}_0,\wt\vph_0)$, and 
$\fo_{\os''}(\wt{V}''_0,\wt\vph_0'')$.
By the construction at the beginning of Section~\ref{tauspinorient_subs},
this in turn implies that the isomorphism~\eref{CROses_e0} respects the orientations
$\fo_{\os}(V,\vph)$, $\fo_{\os'}(V',\vph')$, and $\fo_{\os''}(V'',\vph'')$
under the assumptions~\ref{P1SpinSES_it1} and~\ref{P1SpinSES_it2} above. 
\end{proof}

\begin{crl}\label{SpinP1rk1_crl}
Suppose $(V,\vph)$ is a rank~2 real bundle pair over~$(S^2,\tau)$,
$\os$ is a relative $\OSpin$-structure on the real vector bundle~$V^{\vph}$,
and $D$ is a real CR-operator on~$(V,\vph)$.
The orientation~$\fo_{\os}(V,\vph)$ of~$D$ does not depend
on the choice of the homotopy class in the collection $\os(\bD^2_+)$
of such classes determined by~$\os$ and~$\bD^2_+$.
\end{crl}

\begin{proof}
We denote by $\os_1$ the standard relative $\OSpin$-structure on 
the vector bundle $S^1\!\times\!\R$ over $S^1\!\subset\!S^2$ and by 
$\bp$ be the standard real CR-operator on the real bundle pair
$(S^2\!\times\!\C,\tau\!\times\!\fc)$ over~$(S^2,\tau)$.
Let
$$(V',\vph')=(S^2\!\times\!\C,\tau\!\times\!\fc)\!\oplus\!(V,\vph), \qquad
\os'=\llrr{\os_1,\os}_{\oplus}\in\OSp_{S^2}\!\big(\tau_{S^1}\!\oplus\!V^{\vph}\big).$$
Since $\os_1(\bD^2_+)\!=\!\St_{V^{\vph}}^2(\wch\os_1)$,
the conclusion of Step~1 of the proof of the CROrient~\ref{CROSpinPinSES_prop}\ref{CROsesSpin_it} 
property for $(S^2,\tau)$ implies that the natural isomorphism
$$\la(\bp\!\oplus\!D)\approx \la(\bp)\!\otimes\!\la(D)$$
respects the orientations $\fo_{\os'}(V',\vph')$, $\fo_{\os_1}(S^2\!\times\!\C,\tau\!\times\!\fc)$,
and~$\fo_{\os}(V,\vph)$.
Since the first two of these orientations do not depend 
on the choice of the homotopy class in the collection $\os(\bD^2_+)$,
neither does the third.
\end{proof}

\begin{proof}[{\bf{\emph{Proof of Proposition~\ref{SpinFOOO_prp}}}}]
This is immediate from~\eref{SpinP1rk1_e} if $\rk\,V\!=\!1$.
We can thus assume that $\rk\,V\!\ge\!2$.
Let $\os_1$ and $\bp$ be as in the proof of Corollary~\ref{SpinP1rk1_crl},
$(L,\phi)$ be a rank~1 real bundle pair over~$(\P^1,\tau)$ of even degree,
$\os_L$ be a relative $\OSpin$-structure on the real bundle~$L^{\phi}$ 
over $S^1\!\subset\!S^2$ such that \hbox{$\os_L(\bD^2_+)\!=\!\St_{L^{\phi}}^2(\wch\os_L)$},
and $\bp_L$ be a real CR-operator on~$(L,\phi)$.
For $n\!\in\!\Z^+$, let  
$$(V_n,\vph_n)=(S^2\!\times\!\C^n,\tau\!\times\!\fc), \qquad
(V_{n;L},\vph_{n;L})=\big(V_{n-1},\vph_{n-1}\big)\!\oplus\!(L,\phi).$$
Let $D_n$ be the real CR-operator on $(V_{n;L},\vph_{n;L})$ given by the $n$-fold direct 
sum of the operators~$\bp$ on each factor and \hbox{$D_{n;L}\!=\!D_{n-1}\!\oplus\!\bp_L$}.\\

By Lemma~\ref{SpinFOOO_lmm}, there exist relative $\OSpin$-structures~$\os_2$ and~$\os_2'$ 
on the vector bundle $S^1\!\times\!\R^2$  which induce
the same orientation, but
\BE{SpinFOOO_e5}\fo_{\os_2}(V_2,\vph_2)\neq\fo_{\os_2'}(V_2,\vph_2)\,.\EE
For $n\!\ge\!3$, define
\BE{osnLdfn_e}\begin{split}
\os_{n;L}&=\St^{n-3}_{V_{3;L}^{\vph_{3;L}}}\big(\llrr{\os_2,\os_L}_{\oplus}\big)
=\llrr{\os_1,\os_{n-1;L}}_{\oplus},\\
\os_{n;L}'&=\St^{n-3}_{V_{3;L}^{\vph_{3;L}}}\big(\llrr{\os_2',\os_L}_{\oplus}\big)
=\llrr{\os_1,\os_{n-1;L}'}_{\oplus}\,.
\end{split}\EE
The second equalities on the two lines above hold for $n\!\ge\!4$.
We use the $n\!=\!3$ cases of these identities to define 
relative $\OSpin$-structures~$\os_{2;L}$ and~$\os_{2;L}'$ 
on the vector bundle $S^1\!\times\!\R^2$.
The RelSpinPin~\ref{RelSpinPinStr_prop} and~\ref{RelSpinPinSES_prop} properties 
in Section~\ref{RelSpinPinProp_subs} ensure that this is possible.\\

By~\eref{osnLdfn_e} and 
the conclusion of Step~1 of the proof of the CROrient~\ref{CROSpinPinSES_prop}\ref{CROsesSpin_it} 
property for $(S^2,\tau)$,  the natural isomorphisms
$$\la(D_{3;L})\approx \la(D_2)\!\otimes\!\la(\bp_L) \qquad\hbox{and}\qquad 
\la(D_{3;L})\approx \la(\bp)\!\otimes\!\la(D_{2;L})$$
respect the orientations 
$\fo_{\os_{3;L}}(V_{3;L},\vph_{3;L})$, $\fo_{\os_2}(V_2,\vph_2)$, $\fo_{\os_L}(L,\phi)$, 
$\fo_{\os_1}(V_1,\vph_1)$, and $\fo_{\os_{2;L}}(V_{2;L},\vph_{2;L})$.
They also respect the orientations 
$\fo_{\os_{3;L}'}(V_{3;L},\vph_{3;L})$, $\fo_{\os_2'}(V_2,\vph_2)$, $\fo_{\os_L}(L,\phi)$, 
$\fo_{\os_1}(V_1,\vph_1)$, and\linebreak $\fo_{\os_{2;L}'}(V_{2;L},\vph_{2;L})$.
Along with~\eref{SpinFOOO_e5}, this implies that 
\BE{SpinFOOO_e9}
\fo_{\os_{3;L}}(V_{3;L},\vph_{3;L})\neq \fo_{\os_{3;L}'}(V_{3;L},\vph_{3;L}), \qquad
\fo_{\os_{2;L}}(V_{2;L},\vph_{2;L})\neq \fo_{\os_{2;L}'}(V_{2;L},\vph_{2;L}).\EE

\vspace{.1in}

For each $n\!\ge\!4$, the natural isomorphism
$$\la(D_{n;L})\approx \la(\bp)\!\otimes\!\la(D_{n-1;L})$$
respects the orientations $\fo_{\os_{n;L}}(V_{n;L},\vph_{n;L})$, 
$\fo_{\os_1}(V_1,\vph_1)$, and $\fo_{\os_{n-1;L}}(V_{n-1;L},\vph_{n-1;L})$.
It also respects the orientations $\fo_{\os_{n;L}'}(V_{n;L},\vph_{n;L})$, 
$\fo_{\os_1}(V_1,\vph_1)$, and $\fo_{\os_{n-1;L}'}(V_{n-1;L},\vph_{n-1;L})$.
Along with \eref{SpinFOOO_e9}, this implies~that 
\BE{SpinFOOO_e11}
\fo_{\os_{n;L}}(V_{n;L},\vph_{n;L})\neq \fo_{\os_{n;L}'}(V_{n;L},\vph_{n;L})
\qquad\forall\,n\!\ge\!2.\EE

\vspace{.1in}

By \cite[Proposition~4.1]{BHH}, every real bundle pair $(V,\vph)$ over~$(S^2,\tau)$
is isomorphic to $(V_{n;L},\vph_{n;L})$ for some rank~1 real bundle pair $(L,\phi)$ over~$(S^2,\tau)$.
By~\eref{SpinFOOO_e11}, the real vector bundle $V^{\vph}$ over $S^1\!\subset\!S^2$ 
thus admits relative $\OSpin$-structures~$\os$ and~$\os'$ inducing the same orientation
such that $\fo_{\os}(V,\vph)\!\neq\!\fo_{\os'}(V,\vph)$. 
Since $H_1(S^1;\Z_2)\!\approx\!\Z_2$, there are only two possibilities 
for the collection $\os(\bD^2_+)$ of trivializations of~$V^{\vph}$.
Since the orientation $\fo_{\os}(V,\vph)$ of~$D$ is determined by~$\os(\bD^2_+)$,
the last three statements establish the claim.
\end{proof}

\begin{proof}[{\bf{\emph{Proof of CROrient~\ref{CROSpinPinSES_prop}\ref{CROsesSpin_it} property
for $(S^2,\tau)$, Step~2}}}]
We continue with the notation and setup in Step~1 of the proof.
Suppose $\rk\,V'\!=\!1$,  $\os'(\bD^2_+)\!\neq\!\St_{V'^{\vph'}}^2(\wch\os')$, 
and~\ref{P1SpinSES_it2} holds.
Let $\os_{V'}$ be a relative $\OSpin$-structure on the real bundle~$V'^{\vph'}$ 
over $S^1\!\subset\!S^2$ such that \hbox{$\os_{V'}(\bD^2_+)\!=\!\St_{V'^{\vph'}}^2(\wch\os_{V'})$}.
By the RelSpinPin~\ref{RelSpinPinStr_prop} and~\ref{RelSpinPinSES_prop} properties
in Section~\ref{RelSpinPinProp_subs},
$$\os(\bD^2_+)\equiv \llrr{\os',\os''}_{\ce_{\R}}(\bD^2_+)
\neq  \llrr{\os_{V'},\os''}_{\ce_{\R}}(\bD^2_+)\,.$$
Proposition~\ref{SpinFOOO_prp} then implies that
\BE{P1SpinSES_e7} \fo_{\os'}(V',\vph')\neq\fo_{\os_{V'}}(V',\vph')
\qquad\hbox{and}\qquad 
\fo_{\os}(V,\vph)\neq\fo_{\llrr{\os_{V'},\os''}_{\ce_{\R}}}(V,\vph).\EE
By the conclusion of Step~1 of the proof,
the isomorphism~\eref{CROses_e0} respects the two orientations on the right-hand sides
of the inequalities in~\eref{P1SpinSES_e7} and the orientation $\fo_{\os_{V''}}(V'',\vph'')$.
Thus, it also respects  the two orientations on the left-hand sides
and the orientation~$\fo_{\os_{V''}}(V'',\vph'')$.\\

If $\rk\,V''\!=\!1$,  $\os''(\bD^2_+)\!\neq\!\St_{V''^{\vph''}}^2(\wch\os'')$, 
and~\ref{P1SpinSES_it1} holds, the above argument applies with the roles of $\os'$ and~$\os''$
reversed.
Suppose
$$\rk\,V',\rk\,V''\!=\!1, \qquad \os'(\bD^2_+)\!\neq\!\St_{V'^{\vph'}}^2(\wch\os'), \qquad
\os''(\bD^2_+)\!\neq\!\St_{V''^{\vph''}}^2(\wch\os'').$$
By the RelSpinPin~\ref{RelSpinPinStr_prop} and~\ref{RelSpinPinSES_prop} properties,
$$\os(\bD^2_+)\equiv \llrr{\os',\os''}_{\ce_{\R}}(\bD^2_+)
=  \llrr{\os_{V'},\os_{V''}}_{\ce_{\R}}(\bD^2_+).$$
By Proposition~\ref{SpinFOOO_prp},
\BE{P1SpinSES_e9} \fo_{\os'}(V',\vph')\neq\fo_{\os_{V'}}(V',\vph')
\qquad\hbox{and}\qquad 
\fo_{\os''}(V'',\vph'')\neq\fo_{\os_{V''}}(V'',\vph'').\EE
By the conclusion of Step~1 of the proof,
the isomorphism~\eref{CROses_e0} respects the two orientations on the right-hand sides
of the inequalities in~\eref{P1SpinSES_e9} and the orientation $\fo_{\os}(V,\vph)$.
Thus, it also respects  the two orientations on the left-hand sides
and the orientation~$\fo_{\os}(V,\vph)$.
\end{proof}

\section{Intermediate cases}
\label{InterCR_sec}

In Section~\ref{LB_subs}, we use the intrinsic orientations $\fo(V,\vph;\fo_{x_1})$ 
 of real CR-operators on rank~1 real bundle pairs $(V,\vph)$ over~$(S^2,\tau)$
constructed in Section~\ref{tauorient_subs}
to define orientations $\fo(V,\vph;\fo_{\x})$ 
of real CR-operators on rank~1 real bundle pairs $(V,\vph)$
over arbitrary decorated smooth symmetric surfaces~$(\Si,\si)$.
We deduce a number of properties of the orientations $\fo(V,\vph;\fo_{\x})$
in Sections~\ref{LB_subs} and~\ref{LB_subs2} from the properties of 
the orientations $\fo(V,\vph;\fo_{x_1})$ established in 
Sections~\ref{tauorient_subs}-\ref{tauorient_subs2b}.
In Section~\ref{SpinOrient_subs}, we similarly use the orientations $\fo_{\os}(V,\vph)$ 
of real CR-operators on real bundle pairs~$(V,\vph)$ over~$(S^2,\tau)$
induced by relative $\OSpin$-structures~$\os$ on the real vector bundles~$V^{\vph}$ 
over the $\tau$-fixed locus $S^1\!\subset\!S^2$
to  define orientations $\fo_{\os}(V,\vph)$ 
on the determinants of real CR-operators on real bundle pairs~$(V,\vph)$
over arbitrary decorated smooth symmetric surfaces~$(\Si,\si)$
from relative $\OSpin$-structures~$\os$ on the real vector bundles~$V^{\vph}$ 
over the $\si$-fixed locus $\Si^{\si}\!\subset\!\Si$.
We then deduce all properties of these orientations stated in 
Sections~\ref{OrientPrp_subs1} and~\ref{OrientPrp_subs2}
from the already established properties of the orientations $\fo(V,\vph;\fo_{\x})$
of Section~\ref{LB_subs} and the orientations $\fo_{\os}(V,\vph)$  
of Section~\ref{tauspinorient_subs}.
This concludes the proof of the part of Theorem~\ref{CROrient_thm} concerning orientations
induced by relative $\OSpin$-structures.

\subsection{Orientations for line bundle pairs}
\label{LB_subs}

Suppose $(\Si,\si)$ is a smooth decorated symmetric surface,
$S^1_r\!\subset\!\Si^{\si}$ is the $r$-th connected component of~$\Si^{\si}$
with respect to the chosen order on~$\pi_0(\Si^{\si})$, and 
$$\x\!\equiv\!(x_r)_{S^1_r\in\pi_0(\Si^{\si})}$$
is a tuple of points  so that $x_r\!\in\!S^1_r$.
Let $(V,\vph)$ be a rank~1 real bundle pair over $(\Si,\si)$ and
\BE{fobfxfn_e}\fo_{\x}\equiv \big(\fo_{x_r}\big)_{S^1_r\in\pi_0(\Si^{\si})}\EE 
be a tuple  of orientations of $V_{x_r}^{\vph}$.
We then define an orientation $\fo(V,\vph;\fo_{\x})$
of a real CR-operator~$D$ on $(V,\vph)$ as follows.\\

Let $N\!\equiv\!|\pi_0(\Si^{\si})|$ be the number of connected components of
the fixed locus of~$(\Si,\si)$.
For each $S^1_r$, choose a $\si$-invariant tubular neighborhood $U_r\!\subset\!\Si$ of~$S^1_r$
so that the closures~$\ov{U_r}$ of such neighborhoods are disjoint closed tubular neighborhoods
of~$S^1_r$.
Let \hbox{$\cC_0\!\equiv\!(\Si_0,\si_0)$}
be the nodal symmetric surface obtained from $(\Si,\si)$
by collapsing each of the boundary components of each~$\ov{U_r}$ to a single point;
see Figure~\ref{genorient_fig}.
We denote by $\P^1_r\!\subset\!\Si_0$ the irreducible component containing~$S^1_r$ 
and set $\tau_r\!=\!\si_0|_{\P^1_r}$.
Let $\Si_0^{\C}$ be the union of the remaining irreducible components of~$\Si_0$
and $\si_0^{\C}\!=\!\si_0|_{\Si_0^{\C}}$.
The choice of half-surface $\Si^b$ of~$\Si$ determines a distinguished half-surface
$\bD^2_{r+}$ of $(\P^1_r,\tau_r)$ for each~$r$ and 
a distinguished half-surface $\Si_0^+$ of $(\Si_0^{\C},\si_0^{\C})$ so
that all nodal points of the latter are shared with the nodal points of the former.
We denote the complementary half-surfaces by~$\bD^2_{r-}$ and~$\Si_0^-$, 
respectively, and the unique nodal point on~$\bD^2_{r\pm}$ by~$\nod_r^{\pm}$.\\

Let $(\cU,\wt\fc)$ be a family of deformations of~$\cC_0$ over 
$\De\!\subset\!\C^{2N}$,  $s_1^{\R},\ldots,s_N^{\R}$ be sections 
of $\cU^{\wt\fc}$ over~$\De_{\R}$, and $(\wt{V},\wt\vph)$ be a rank~1 real bundle pair
over $(\cU,\wt\fc)$ so~that 
\BE{Sisit1cond_e}\big(\Si_{\t_1},\si_{\t_1}\big)=(\Si,\si), \quad
s_r^{\R}(\t_1)=x_r~~\forall\,r\!\in\![N], 
\quad\hbox{and}\quad
\big(\wt{V}_{\t_1},\wt\vph_{\t_1}\big)=(V,\vph)\EE
for some $\t_1\!\in\!\De_{\R}$.
For each $\t\!\in\!\De_{\R}$, $\fo_{\x}$ determines a tuple
$$\fo_{\x;\t}\equiv\big(\fo_{x_r;\t}\big)_{S^1_r\in\pi_0(\Si^{\si})}$$
of orientations of $\wt{V}^{\wt\vph}$ at $s_r^{\R}(\t)$ via 
the line bundles $s_r^{\R*}\wt{V}^{\wt\vph}$ over~$\De_{\R}$.
If $w_1(V^{\vph})|_{S^1_r}\!=\!0$, then~$\fo_{x_r;\t}$ does not depend on the choice of 
the section~$s_r^{\R}$.\\

\begin{figure}
\begin{pspicture}(-.3,-3.5)(10,3)
\psset{unit=.4cm}
\psellipse(10,0)(5,7)
\psarc[linewidth=.04](13,0){4}{130}{230}\psarc[linewidth=.04](7,0){4}{-50}{50} 
\psarc[linewidth=.03](7,3.46){4}{240}{300} 
\psarc[linewidth=.03,linestyle=dashed](7,-3.46){4}{60}{120}
\psarc[linewidth=.03](13,3.46){4}{240}{300} 
\psarc[linewidth=.03,linestyle=dashed](13,-3.46){4}{60}{120}
\psarc[linewidth=.04](10,6.7){1}{210}{330}\psarc[linewidth=.04](10,5.3){1}{30}{150} 
\psarc[linewidth=.04](10,-5.3){1}{210}{330}\psarc[linewidth=.04](10,-6.7){1}{30}{150}
\psarc[linewidth=.03](7.4,5.46){4.1}{240}{300} 
\psarc[linewidth=.03,linestyle=dashed](7.4,-1.46){4.1}{60}{120}
\psarc[linewidth=.03](12.6,5.46){4.1}{240}{300} 
\psarc[linewidth=.03,linestyle=dashed](12.6,-1.46){4.1}{60}{120}
\psarc[linewidth=.03](7.4,1.46){4.1}{240}{302} 
\psarc[linewidth=.03,linestyle=dashed](7.4,-5.46){4.1}{60}{120}
\psarc[linewidth=.03](12.6,1.46){4.1}{238}{302} 
\psarc[linewidth=.03,linestyle=dashed](12.6,-5.46){4.1}{60}{120}
\rput(6.9,0){\sm{$S^1_1$}}\rput(12.9,0){\sm{$S^1_2$}}
\rput(10,4.3){\sm{$\Si^b$}}
\rput(7.2,2){\sm{$\prt\ov{U_1}$}}\rput(12.8,2){\sm{$\prt\ov{U_2}$}}
\rput(7.2,-2){\sm{$\prt\ov{U_1}$}}\rput(12.8,-2){\sm{$\prt\ov{U_2}$}}
\psline{<->}(4,4)(4,-4)\rput(3.3,0){$\si$}
\rput(10,-8.5){$(\Si,\si)$}
\pscircle(23,0){2}\pscircle(29,0){2}
\psarc[linewidth=.03](23,3.46){4}{240}{300} 
\psarc[linewidth=.03,linestyle=dashed](23,-3.46){4}{60}{120}
\psarc[linewidth=.03](29,3.46){4}{240}{300} 
\psarc[linewidth=.03,linestyle=dashed](29,-3.46){4}{60}{120}
\psarc(23,3){1}{180}{360}\psarc(23,-3){1}{0}{180}
\psarc(29,3){1}{180}{360}\psarc(29,-3){1}{0}{180}
\psarc(26,3){2}{0}{180}\psarc(26,-3){2}{180}{360}
\psarc(26,3){4}{0}{180}\psarc(26,-3){4}{180}{360}
\pscircle*(23,2){.2}\pscircle*(23,-2){.2}
\pscircle*(29,2){.2}\pscircle*(29,-2){.2}
\psarc[linewidth=.04](26,6.7){1}{210}{330}\psarc[linewidth=.04](26,5.3){1}{30}{150} 
\psarc[linewidth=.04](26,-5.3){1}{210}{330}\psarc[linewidth=.04](26,-6.7){1}{30}{150} 
\rput(23,2.7){\sm{$\nod_1^+$}}\rput(23,-2.7){\sm{$\nod_1^-$}}
\rput(29,2.7){\sm{$\nod_2^+$}}\rput(29,-2.7){\sm{$\nod_2^-$}}
\rput(22.9,0){\sm{$S^1_1$}}\rput(28.9,0){\sm{$S^1_2$}}
\rput(25,1.9){\sm{$\bD^2_{1+}$}}\rput(27,2.1){\sm{$\bD^2_{2+}$}}
\rput(25.1,-2){\sm{$\bD^2_{1-}$}}\rput(27.1,-1.8){\sm{$\bD^2_{2-}$}}
\rput(20.3,0){$\P^1_1$}\rput(31.7,0){$\P^1_2$}
\rput(22.2,6){$\Si_0^+$}\rput(22.2,-6){$\Si_0^-$}
\psline{<->}(33,4)(33,-4)\rput(33.8,0){$\si_0$}
\rput(26,-8.5){$\cC_0$}
\end{pspicture}
\caption{A decorated smooth symmetric surface~$(\Si,\si)$ 
and its associated degeneration as below~\eref{fobfxfn_e}}
\label{genorient_fig}
\end{figure}

Let \hbox{$\cD\!\equiv\!\{D_{\t}\}$} be a family of real CR-operators on
$(\wt{V}_{\t},\wt\vph_{\t})$ as in~\eref{DVvphdfn_e} so~that $D_{\t_1}\!=\!D$.
For each $r\!\in\![N]$, we denote by $(\wt{V}_{0r},\wt\vph_{0r})$ 
the restriction of $(\wt{V}_0,\wt\vph_0)$ to~$\P^1_r$ and 
by $D_{0r}$ the real CR-operator on  $(\wt{V}_{0r},\wt\vph_{0r})$ induced by~$D_0$.
Let $D_0^+$ be the real CR-operator on~$\wt{V}_0|_{\Si_0^+}$ induced by~$D_0$.
The exact~triple 
\begin{gather}\label{geCndses_e}
0\lra D_0 \lra \bigoplus_{S^1_r\in\pi_0(\Si^{\si})}\!\!\!\!\!\!\!\!D_{0r}
\oplus\!D_0^+\lra 
\bigoplus_{S^1_r\in\pi_0(\Si^{\si})}\!\!\!\!\!\!\!\!\wt{V}_{\nod_r^+}\lra 0, \\
\notag
\big(\xi_-,(\xi_r)_{S^1_r\in\pi_0(\Si^{\si})} ,\xi_+\big)\lra
\big((\xi_r)_{S^1_r\in\pi_0(\Si^{\si})} ,\xi_+\big), \\
\notag
\big((\xi_r)_{S^1_r\in\pi_0(\Si^{\si})} ,\xi_+\big)\lra
 \big(\xi_+(\nod_r^+)\!-\!\xi_r(\nod_r^+)\big)_{S^1_r\in\pi_0(\Si^{\si})}, 
\end{gather}
of Fredholm operators then determines an isomorphism
\BE{geCndses_e2}
\la(D_0)\!\otimes\!\bigotimes_{S^1_r\in\pi_0(\Si^{\si})}\!\!\!\!\!\!\!\!
\la\big(\wt{V}_{\nod_r^+}\big)\approx 
\bigotimes_{S^1_r\in\pi_0(\Si^{\si})}\!\!\!\!\!\!\!\!\la\big(D_{0r}\big)
\otimes\!\la\big(D_0^+\big).\EE
For each $S^1_r\!\in\!\pi_0(\Si^{\si})$, let 
$\fo(\wt{V}_{0r},\wt\vph_{0r};\fo_{x_r;0})$ be the intrinsic orientation
of $D_{0r}$ as above Proposition~\ref{cOpincg_prp}.
Along with the chosen order on~$\pi_0(\Si^{\si})$ and 
the complex orientations of $D_0^+$ and $\wt{V}_{\nod_r^+}$,
these orientations determine an orientation $\fo(\wt{V}_0,\wt\vph_0;\fo_{\x})$
of $\la(D_0)$ via the isomorphism~\eref{geCndses_e2} and 
thus an orientation of the line bundle~$\la(\cD)$ over~$\De_{\R}$.
The latter restricts to an orientation $\fo(V,\vph;\fo_{\x})$ of~$\la(D)$.
By Proposition~\ref{cOpincg_prp}\ref{cOpincg_it1}, $\fo(V,\vph;\fo_{\x})$ does not depend
on the choice of the section~$s_r^{\R}$ above even if 
$w_1(V^{\vph})|_{S^1_r}\!\neq\!0$.
We call $\fo(V,\vph;\fo_{\x})$ the \sf{intrinsic orientation} of~$D$. 

\begin{prp}\label{gencOpincg_prp}
Suppose $(\Si,\si)$ is a smooth decorated symmetric surface,
$(V,\vph)$ is a rank~1 real bundle pair over~$(\Si,\si)$, and
$D$ is a real CR-operator on~$(V,\vph)$.
Let~$\fo_{\x}$ be a tuple of orientations  of~$V^{\vph}$ at 
points~$x_r$ in $S_r^1\!\subset\!\Si^{\si}$ as in~\eref{fobfxfn_e}.
\begin{enumerate}[label=(\arabic*),leftmargin=*]

\item\label{gencOpincg_it0} The intrinsic orientation $\fo(V,\vph;\fo_{\x})$ of~$D$
does not depend on the choice of the real bundle pair $(\wt{V},\wt\vph)$ over~$(\cU,\wt\fc)$
satisfying the last condition in~\eref{Sisit1cond_e}.

\item\label{gencOpincg_it1a} 
The orientation $\fo(V,\vph;\fo_{\x})$ does not depend on the choice of half-surface $\Si_*^b$ of 
an elemental component~$\Si_*$ of $(\Si,\si)$ if and only~if 
$$\frac{g(\Si_*)\!+\!|\pi_0(\Si_*^{\si})|\!-\!1}{2}
+\frac{\deg V|_{\Si_*}\!+\!W_1((V,\vph)|_{\Si_*})}{2}\in 2\Z.$$

\item\label{gencOpincg_it1b}  
The interchange in the ordering of two consecutive components~$S_r^1$ and~$S_{r+1}^1$ 
of~$\Si^{\si}$ reverses the orientation $\fo(V,\vph;\fo_{\x})$ 
if and only if $w_1(V^{\vph})|_{S^1_r},w_1(V^{\vph})|_{S^1_{r+1}}\!=\!0$.

\item\label{gencOpincg_it2} Reversing the component orientation $\fo_{x_r}$ in~\eref{fobfxfn_e}
preserves the orientation $\fo(V,\vph;\fo_{\x})$ if and only if 
$w_1(V^{\vph})|_{S^1_r}\!\neq\!0$.

\item\label{gencOdisjun_it} Suppose $(\Si,\si)\!=\!(\Si_1,\si_1)\!\sqcup\!(\Si_2,\si_2)$
is a decomposition of~$(\Si,\si)$ into decorated symmetric surfaces as below~\eref{Si1Si2dfn_e2},
$(V_1,\vph_1)$ and $(V_2,\vph_2)$ are the restrictions of $(V,\vph)$
to $\Si_1$ and~$\Si_2$, respectively, and 
\hbox{$\fo_{\x}\!=\!\fo_{\x}^{(1)}\fo_{\x}^{(2)}$} is the associated decomposition of~$\fo_{\x}$.
If $D_1$ and~$D_2$ are the real CR-operators on $(V_1,\vph_1)$ and $(V_2,\vph_2)$ induced by~$D$,
then the isomorphism
$$\la(D)\approx\la\big(D_1\big)\!\otimes\!\la\big(D_2\big)$$
induced by~\eref{D1D2split_e} respects the orientations 
$\fo(V,\vph;\fo_{\x})$, $\fo(V_1,\vph_1;\fo_{\x}^{(1)})$,
and $\fo(V_2,\vph_2;\fo_{\x}^{(2)})$.

\end{enumerate}
\end{prp}

\vspace{.15in}

We establish Proposition~\ref{gencOpincg_prp}\ref{gencOpincg_it0}, 
as well as Propositions~\ref{genrk1Cdegen_prp} and~\ref{genrk1H3degen_prp} stated below,
in Section~\ref{LB_subs2}.

\begin{proof}[{\bf{\emph{Proof of 
Proposition~\ref{gencOpincg_prp}\ref{gencOpincg_it1a}-\ref{gencOdisjun_it}}}}]
We continue with the notation above the statement of the proposition.
For each $S^1_r\!\in\!\pi_0(\Si^{\si})$, let 
$$\ep_r(V,\vph)=\begin{cases}0,&\hbox{if}~w_1(V^{\vph})|_{S^1_r}\!=\!0;\\
1,&\hbox{if}~w_1(V^{\vph})|_{S^1_r}\!\neq\!0.
\end{cases}$$
With the notation as in~\ref{gencOpincg_it1a}, 
let $\Si_{0*}\!\subset\!\Si_0$ be the elemental component corresponding to $\Si_*\!\subset\!\Si$ and
$$\Si_{0*}^+=\Si_{0*}\!\cap\!\Si_0^+\,.$$
We note~that 
\BE{gencOpincg_e2}\begin{split}
g(\Si_{*0}^+)=\frac{g(\Si_*)\!-\!|\pi_0(\Si_*^{\si})|\!+\!1}{2}, &\quad
\deg\wt{V}_0|_{\Si_{*0}^+}= 
\frac{\deg V|_{\Si_*}-\sum\limits_{S^1_r\in\pi_0(\Si_*^{\si})}\!\!\!\!\!\!\deg\wt{V}_{0r}}{2},\\
\ind_{\C}D_0^+|_{\Si_{*0}^+}&=1\!-\!g(\Si_{*0}^+)\!+\!\deg\wt{V}_0|_{\Si_{*0}^+}.
\end{split}\EE
The change in the choice of the half-surface $\Si_*^b$ of $(\Si_*,\si)$ acts by 
the complex conjugation on the complex orientations of $D_0^+|_{\Si_{*0}^+}$
and each $\wt{V}_{\nod_r^+}$ with $S^1_r\!\in\!\pi_0(\Si_*^{\si})$.
By Proposition~\ref{cOpincg_prp}\ref{cOpincg_it1}, this change preserves
the orientation $\fo(\wt{V}_r,\wt\vph_r;\fo_{x_r;0})$ of $D_{0r}$
with $S^1_r\!\in\!\pi_0(\Si_*^{\si})$ if and only~if 
$$\frac{\deg\wt{V}_{0r}\!+\!\ep_r(V,\vph)}{2}\in2\Z\,.$$
Thus, the orientations $\fo(\wt{V}_0,\wt\vph_0;\fo_{\x})$ and
$\fo(V,\vph;\fo_{\x})$ do not depend on the choice of half-surface $\Si_*^b$ 
of $(\Si_*,\si)$ if and only~if
$$\ind_{\C}D_0^+|_{\Si_{*0}^+}-\big|\pi_0(\Si_*^{\si})\big|+
\sum_{S^1_r\in\pi_0(\Si_*^{\si})}\!\!\!\!\!\!\frac{\deg\wt{V}_{0r}\!+\!\ep_r(V,\vph)}{2}
\in2\Z.$$
Combining this with~\eref{gencOpincg_e2}, we obtain~\ref{gencOpincg_it1a}.\\

Since the index of $D_{0r}$ is odd if and only if $w_1(V^{\vph})|_{S^1_r}\!=\!0$, 
\ref{gencOpincg_it1b} follows from the Direct Sum property for 
the determinants of Fredholm operators; see \cite[Section~2]{detLB}.
Proposition~\ref{cOpincg_prp}\ref{cOpincg_it1} immediately implies~\ref{gencOpincg_it2}.\\

It remains to establish~\ref{gencOdisjun_it}.
Let $(\cU_1,\wt\fc_1)$, $(\cU_2,\wt\fc_2)$, 
$(\wt{V}_1,\wt\vph_1)$, $(\wt{V}_2,\wt\vph_2)$,
$\fo_{\x;\t}^{(1)}$, $\fo_{\x;\t}^{(2)}$, $D_{1;\t}$, and $D_{2;\t}$ 
be as in the construction of the orientations 
$\fo(V_1,\vph_1;\fo_{\x}^{(1)})$ and $\fo(V_2,\vph_2;\fo_{\x}^{(2)})$
above Proposition~\ref{gencOpincg_prp}.
The orientation $\fo(V,\vph;\fo_{\x})$ 
on $D$ is obtained via this construction applied with
$$\big(\cU,\wt\fc\big)=
\big(\cU_1,\wt\fc_1\big)\!\!\sqcup\!\!\big(\cU_2,\wt\fc_2\big),~~
\big(\wt{V},\wt\vph\big)=\big(\wt{V}_1,\wt\vph_1\big)
\!\!\sqcup\!\!\big(\wt{V}_2,\wt\vph_2\big),~~
 \fo_{\x;\t}=\fo_{\x;\t}^{(1)}\fo_{\x;\t}^{(2)},~~
D_{\t}=D_{1;\t}\!\sqcup\!D_{2;\t}.$$

\vspace{.1in}

The exact triple~\eref{geCndses_e} of Fredholm operators induces 
the exact triples of Fredholm operators given by the rows in the diagram of 
Figure~\ref{rk1CRODisjUn_fig}.
The decompositions
$$\Si_0=\Si_{1;0}\!\sqcup\!\Si_{2;0}, \quad
\Si_0^+=\Si_{1;0}^+\!\sqcup\!\Si_{2;0}^+, \quad\hbox{and}\quad
\pi_0\big(\Si^{\si}\big)=\pi_0\big(\Si_1^{\si_1}\big)\!\sqcup\!\pi_0\big(\Si_2^{\si_2}\big)$$
induce the exact triples of Fredholm operators given by the columns in this diagram.
By definition, the rows in this diagrams respect the orientations
\BE{LHSorient_e}
\fo\big(\wt{V}_{1;0},\wt\vph_{1;0};\fo_{\x}^{(1)}\big), \quad 
\fo\big(\wt{V}_0,\wt\vph_0;\fo_{\x}\big), \quad\hbox{and}\quad 
\fo\big(\wt{V}_{2;0},\wt\vph_{2;0};\fo_{\x}^{(2)}\big)\EE 
of the operators in the left column, the orientations 
\BE{midorient_e}
\fo\big(\wt{V}_{0r},\wt\vph_{0r};\fo_{x_{r;0}}\big)=
\begin{cases}\fo\big(\wt{V}_{0r},\wt\vph_{0r};\fo_{x_r;0}^{(1)}\big),
&\hbox{if}~S^1_r\!\in\!\pi_0(\Si_1^{\si_1});\\
\fo\big(\wt{V}_{0r},\wt\vph_{0r};\fo_{x_r;0}^{(2)}\big),
&\hbox{if}~S^1_r\!\in\!\pi_0(\Si_2^{\si_2});
\end{cases}\EE
of the operators in the direct sums in the middle column,
and the complex orientations of the remaining terms.\\

\begin{figure}
$$\xymatrix{& 0\ar[d] & 0\ar[d] & 0\ar[d]  \\
0\ar[r] & D_{1;0}\ar[r]\ar[d] & 
\bigoplus\limits_{S^1_r\in\pi_0(\Si_1^{\si_1})}\!\!\!\!\!\!\!\!D_{0r}
\oplus\!D_{1;0}^+ \ar[r]\ar[d]  & 
\bigoplus\limits_{S^1_r\in\pi_0(\Si_1^{\si_1})}\!\!\!\!
\wt{V}_{\nod_r^+}\ar[r]\ar[d] & 0 \\
0\ar[r] & D_0\ar[r]\ar[d] & 
\bigoplus\limits_{S^1_r\in\pi_0(\Si^{\si})}\!\!\!\!\!\!\!\!D_{0r}
\oplus\!D_0^+ \ar[r]\ar[d]  & 
\bigoplus\limits_{S^1_r\in\pi_0(\Si^{\si})}\!\!\!\!\!\!\wt{V}_{\nod_r^+}\ar[r]\ar[d] & 0 \\
0\ar[r] & D_{2;0}\ar[r]\ar[d] & 
\bigoplus\limits_{S^1_r\in\pi_0(\Si_2^{\si_2})}\!\!\!\!\!\!\!\!D_{0r}
\oplus\!D_{2;0}^+ \ar[r]\ar[d]  & 
\bigoplus\limits_{S^1_r\in\pi_0(\Si_2^{\si_2})}\!\!\!\!\wt{V}_{\nod_r^+}\ar[r]\ar[d] & 0 \\
& 0 & 0 & 0}$$ 
\caption{Commutative square of exact rows and columns of Fredholm operators 
for the proof of Proposition~\ref{gencOpincg_prp}\ref{gencOdisjun_it}}
\label{rk1CRODisjUn_fig}
\end{figure}

The right column and the exact triple formed by the last summands in the middle column
respect the complex orientations of the corresponding operators.  
Since the decomposition of~$(\Si,\si)$ into $(\Si_1,\si_1)$ and $(\Si_2,\si_2)$
respects the orderings of the components of the fixed loci,
the exact triple formed by the direct sums in the middle column
respects the ordered direct sum orientations.
Since the (real) index of $D_{1;0}^+$ is even,
Lemma~\ref{3by3_lmm} and the last two statements imply that 
the middle column respects the direct sum orientations.
Since the (real) dimensions of $\wt{V}_{\nod_r^+}^{(1)}$ are even,
Lemma~\ref{3by3_lmm} again, the last sentence in the previous paragraph,
and the conclusions in this paragraph imply that 
the left column in Figure~\ref{rk1CRODisjUn_fig} respects the orientations~\eref{LHSorient_e}.
The claim in~\ref{gencOdisjun_it} now follows from the continuity of 
the isomorphisms
$$\la\big(D_{\t}\big)\approx \la\big(D_{\t;1}\big)\!\otimes\!\la\big(D_{\t;2}\big)$$
with respect to $\t\!\in\!\De_{\R}$.
\end{proof}

We next describe the behavior of the orientations of Proposition~\ref{gencOpincg_prp} 
under flat degenerations of~$(\Si,\si)$ to nodal symmetric surfaces as
in the CROrient~\ref{CROdegenC_prop} and~\ref{CROdegenH3_prop} properties 
of Section~\ref{OrientPrp_subs2}.
We first suppose that $\cC_0$ is a decorated symmetric surface with 
one conjugate pair of nodes~$\nod^{\pm}$ as in~\eref{cCdegC_e}
and in the top left diagram of Figure~\ref{genCdegen_fig} on page~\pageref{genCdegen_fig}
so that $x_r\!\in\!S^1_r$ for each $r\!\in\![N]$.
We also suppose that 
$(V_0,\vph_0)$  is a rank~1 real bundle pair over~$(\Si_0,\si_0)$
and $\fo_{\x}$ is a tuple of orientations of $V_0^{\vph_0}|_{x_r}$
as in~\eref{fobfxfn_e}.
Let $\wt\cC_0$, $(\wt\Si_0,\wt\si_0)$, and $\wt\Si_0^b$ be as below~\eref{cCdegC_e},
$(\wt{V}_0,\wt\vph_0)$ be the lift of $(V_0,\vph_0)$ to a real bundle pair over $(\wt\Si_0,\wt\si_0)$,
and $\wt\fo_{\x}$ be the lift of $\fo_{\x}$ to a tuple of orientations 
of fibers of~$\wt{V}_0^{\wt\vph_0}$.\\

Let $D_0$ be a real CR-operator on $(V_0,\vph_0)$.
We denote its lift to a real CR-operator on $(\wt{V}_0,\wt\vph_0)$ by~$\wt{D}_0$. 
The orientation 
\BE{Cdegenorient_e3}\wt\fo_0\big(\fo_{\x}\big)\equiv\fo\big(\wt{V}_0,\wt\vph_0;\wt\fo_{\x}\big)\EE 
of $\wt{D}_0$ and the complex orientation of $V_0|_{\nod^+}$ determine an orientation 
\BE{rk1Csplit_e4}\fo_0\big(\fo_{\x}\big)\equiv \fo\big(V_0,\vph_0;\fo_{\x}\big)\EE
of $D_0$ via the isomorphism~\eref{Cdegses_e2}.
In an analogy with the intrinsic orientation of the CROrient~\ref{CROdegenC_prop} property,
we call \eref{rk1Csplit_e4} the \sf{intrinsic orientation} of~$D_0$
induced by~$\fo_{\x}$.\\

Suppose in addition that $(\cU,\wt\fc)$ is a flat family of deformations of~$\cC_0$
as in~\eref{cUsymmdfn_e},
$(V,\vph)$ is a real bundle pair over $(\cU,\wt\fc)$ extending~$(V_0,\vph_0)$,
$s_1^{\R},\ldots,s_N^{\R}$ are sections of $\cU^{\wt\fc}$ over~$\De_{\R}$ with 
\hbox{$s_r^{\R}(0)\!=\!x_r$} for all $r\!\in\![N]$,
and  \hbox{$\cD\!\equiv\!\{D_{\t}\}$} is a family of real CR-operators on
$(V_{\t},\vph_{\t})$ as in~\eref{DVvphdfn_e} extending~$D_0$.
The decorated structure on~$\cC_0$ induces a decorated structure on
the fiber $(\Si_{\t},\si_{\t})$ of~$\pi$ for every $\t\!\in\!\De_{\R}$ as 
above the CROrient~\ref{CROdegenC_prop} property.
The orientation~$\fo_{x_r}$ of $V_0^{\vph_0}|_{x_r}$ induces an orientation
of the line bundle $s_r^{\R*}V^{\vph}$ over $\De_{\R}$, which in turn restricts to
an orientation $\fo_{x_r;\t}$ of $V^{\vph}$ at $s_r^{\R}(\t)$ for each $\t\!\in\!\De_{\R}$.
For each $\t\!\in\!\De_{\R}^*$, the~tuple
\BE{fobfxfn_e2}\fo_{\x;\t}\equiv \big(\fo_{x_r;\t}\big)_{S^1_r\in\pi_0(\Si_0^{\si_0})}\EE
of orientations of $V_{\t}^{\vph_{\t}}|_{s_r^{\R}(\t)}$ determines an orientation
$$\fo_{\t}\equiv \fo\big(V_{\t},\vph_{\t};\fo_{\x;\t}\!\big)$$
of $D_{\t}$ as above Proposition~\ref{gencOpincg_prp}.
These orientations vary continuously with~$\t$ and extend to an orientation
\BE{rk1Clim_e4}\fo_0'(\fo_{\x})\equiv \fo'\big(V_0,\vph_0;\fo_{\x}\big)\EE
of~$D_0$ as above the CROrient~\ref{CROdegenC_prop} property on page~\pageref{CROdegenC_prop}.
In an analogy with the limiting orientation of the CROrient~\ref{CROdegenC_prop} property,
we call \eref{rk1Clim_e4} the \sf{limiting orientation} of~$D_0$
induced by~$\fo_{\x}$.

\begin{prp}\label{genrk1Cdegen_prp}
Suppose $\cC_0$ is a decorated marked symmetric surface as in~\eref{cCdegC_e} 
which contains precisely one conjugate pair $(\nod^+,\nod^-)$ of nodes and no other nodes
and carries precisely one real marked point~$x_r$ on each connected component~$S^1_r$ 
of~$\Si_0^{\si_0}$. 
Let $(V_0,\vph_0)$ be a rank~1 real bundle pair over~$(\Si_0,\si_0)$,
$\fo_{\x}$ be a tuple of orientations of $V_0^{\vph_0}|_{x_r}$
as in~\eref{fobfxfn_e}, and $D_0$ be a real CR-operator on $(V_0,\vph_0)$.
The intrinsic and limiting orientations, \eref{rk1Csplit_e4} and~\eref{rk1Clim_e4},
of~$D_0$ are the same.
\end{prp}

We next suppose that $\cC_0$ is a decorated symmetric surface with 
one $H3$~node~$\nod$ as above~\eref{wtcCH3deg_e}
and in the top left diagram of Figure~\ref{genH3degen_fig} on page~\pageref{genH3degen_fig}
so that $x_r\!\in\!S^1_r$ for each $r\!\in\![N]$ different from the index~$r_{\bu}$
of the singular topological component of~$\Si_0^{\si_0}$ and $x_{r_{\bu}}\!=\!\nod$
(in this case, we allow a ``marked point" to be a node).
We also suppose that $(V_0,\vph_0)$ is a rank~1 real bundle pair over~$(\Si_0,\si_0)$
and $\fo_{\x}$ is a tuple of orientations of $V_0^{\vph_0}$ at points $x_r\!\in\!S^1_r$
as in~\eref{fobfxfn_e}.
Let $\wt\cC_0$, $(\wt\Si_0,\wt\si_0)$, $\wt\Si_0^b$, 
$\nod_1\!\in\!S^1_{\bu1}$, and $\nod_2\!\in\!S^1_{\bu2}$ be as in and below~\eref{wtcCH3deg_e},
$(\wt{V}_0,\wt\vph_0)$ be the lift of $(V_0,\vph_0)$ to a real bundle pair over $(\wt\Si_0,\wt\si_0)$,
and 
\BE{wtfoH3gendfn_e}\wt\fo_{\x}\equiv \big(\fo_{\nod_1},\fo_{\nod_2},
\big(\fo_{x_r}\big)_{S^1_r\in\pi_0(\Si_0^{\si_0}),r\neq r_{\bu}}\big)\EE
be the tuple of orientations of the fibers of~$\wt{V}_0^{\wt\vph_0}$
at the points $\nod_1,\nod_2$ and $x_r$ of $\wt\Si_0^{\wt\si_0}$.\\

Let $D_0$ be a real CR-operator on $(V_0,\vph_0)$.
We denote its lift to a real CR-operator on $(\wt{V}_0,\wt\vph_0)$ by~$\wt{D}_0$. 
The orientation 
\BE{rk1H3split_e4a}\wt\fo_0\big(\fo_{\x}\big)\equiv\fo\big(\wt{V}_0,\wt\vph_0;\wt\fo_{\x}\big)\EE 
of $\wt{D}_0$ and the orientation $\fo_{\nod}\!\equiv\!\fo_{x_{r_{\bu}}}$ 
of $V_0^{\vph_0}|_{\nod}$ determine an orientation 
\BE{rk1H3split_e4}\fo_0\big(\fo_{\x}\big)\equiv \fo\big(V_0,\vph_0;\fo_{\x}\big)\EE
of $D_0$ via the isomorphism~\eref{Rdegses_e2}.
We call \eref{rk1H3split_e4} the \sf{intrinsic orientation} of~$D_0$
induced by~$\fo_{\x}$.\\

Suppose in addition that $(\cU,\wt\fc)$ is a flat family of deformations of~$\cC_0$
as in~\eref{cUsymmdfn_e},
$(V,\vph)$ is a real bundle pair over $(\cU,\wt\fc)$ extending~$(V_0,\vph_0)$,
$s_1^{\R},\ldots,s_N^{\R}$ are sections of $\cU^{\wt\fc}$ over~$\De_{\R}$ such~that
\BE{H3sectcond_e}s_r^{\R}(0)=x_r~~\forall\,r\!\in\![N]\!-\!\big\{r_{\bu}\big\} \qquad\hbox{and}\qquad
s_{r_{\bu}}^{\R}(0)\in S^1_{\bu1}\!-\!\big\{\nod\big\},\EE
and  \hbox{$\cD\!\equiv\!\{D_{\t}\}$} is a family of real CR-operators on
$(V_{\t},\vph_{\t})$ as in~\eref{DVvphdfn_e} extending~$D_0$.
The decorated structure on~$\cC_0$ induces a decorated structure on
the fiber $(\Si_{\t},\si_{\t})$ of~$\pi$ for every \hbox{$\t\!\in\!\De_{\R}$} 
and determines open subspaces $\De_{\R}^{\pm}\!\subset\!\De_{\R}$
as above the CROrient~\ref{CROdegenH3_prop} property on page~\pageref{CROdegenH3_prop}.
For each $r\!\neq\!r_{\bu}$ and $\t\!\in\!\De_{\R}$, 
$\fo_{x_r}$ determines an orientation~$\fo_{x_r;\t}$ of $V^{\vph}_{s_r^{\R}(\t)}$
as above Proposition~\ref{genrk1Cdegen_prp}.
The orientation~$\fo_{\nod}$ determines an orientation~$\fo_{x_{r_{\bu}};\t}$ of 
$V^{\vph}$ at $s_{r_{\bu}}^{\R}(\t)$ as above Proposition~\ref{cOlim_prp} 
(with $S_1^1$ there replaced by~$S^1_{\bu1}$).
For each $\t\!\in\!\De_{\R}^*$,
the resulting tuple~\eref{fobfxfn_e2} of orientations of $V_{\t}^{\vph_{\t}}|_{s_r^{\R}(\t)}$ 
determines an orientation
\BE{H3degenorient_e3}
\fo_{\t}\equiv \fo\big(V_{\t},\vph_{\t};\fo_{\x;\t}\!\big)\EE
of $D_{\t}$ as above Proposition~\ref{gencOpincg_prp}.
We denote~by 
\BE{cOlimor_e4} \fo_0^+\big(\fo_{\x}\big)\equiv 
\fo_0^+\big(V_0,\vph_0;\fo_{\x}\big)\EE
the orientation of~$D_0$ obtained as the continuous extension of the orientations~\eref{H3degenorient_e3}
with $\t\!\in\!\De_{\R}^+$.
We call~\eref{cOlimor_e4} the \sf{limiting orientation} of~$D_0$.
Define
\BE{W1Vvphrbudfn_e}W_1(V_0,\vph_0)_{r_{\bu}}=\big|\big\{S_r^1\!\in\!\pi_0(\Si^{\si})\!:
r\!>\!r_{\bu},\,w_1(V_0^{\vph_0})|_{S_r^1}\!\neq\!0\big\}\big|.\EE

\begin{prp}\label{genrk1H3degen_prp}
Suppose $\cC_0$ is a decorated marked symmetric surface as in~\eref{cCdegC_e} 
which contains precisely one~$H3$ node~$\nod$ and no other nodes
and carries precisely one real marked point~$x_r$ on each smooth connected component~$S^1_r$ 
of~$\Si_0^{\si_0}$ and a marked point $x_{r_{\bu}}\!=\!\nod$ on the nodal connected 
component of~$\Si_0^{\si_0}$.
Let $(V_0,\vph_0)$ be a rank~1 real bundle pair over~$(\Si_0,\si_0)$,
$\fo_{\x}$ be a tuple of orientations of $V_0^{\vph_0}|_{x_r}$
as in~\eref{fobfxfn_e}, and $D_0$ be a real CR-operator on $(V_0,\vph_0)$.
The intrinsic and limiting  orientations, \eref{rk1H3split_e4} and~\eref{cOlimor_e4},
of~$D_0$ are the same if and only~if
\BE{genrk1H3degen_e0}\big(\lr{w_1(V_0^{\vph_0}),[S^1_{\bu1}]_{\Z_2}}\!+\!1\big)
\lr{w_1(V_0^{\vph_0}),[S^1_{\bu2}]_{\Z_2}}
+\big|\pi_0(\Si_0^{\si_0})\big|\!-\!r_{\bu}\!-\!W_1(V_0,\vph_0)_{r_{\bu}} 
=0\in \Z_2.\EE 
\end{prp}

\subsection{Proofs of Propositions~\ref{gencOpincg_prp}\ref{gencOpincg_it0}, 
\ref{genrk1Cdegen_prp}, and~\ref{genrk1H3degen_prp}}
\label{LB_subs2}

In all three proofs, we continue with the notation in the statement of the corresponding
proposition and just above.

\begin{proof}[{\bf{\emph{Proof of Proposition~\ref{gencOpincg_prp}\ref{gencOpincg_it0}}}}]
The substance of this claim
is that the orientation $\fo(V,\vph;\fo_{\x})$ does not depend on 
the choice of the restriction of $(\wt{V},\wt\vph)$ to each real component~$(\P^1_r,\tau_r)$
of~$(\Si_0,\si_0)$.
In light of Proposition~\ref{gencOpincg_prp}\ref{gencOpincg_it1b}, it is sufficient to show that
this is the case for $r\!=\!1$.
Thus, suppose that $(\wh{V},\wh\vph)$ is a rank~1 real bundle pair
over $(\cU,\wt\fc)$ so~that 
$$\big(\wh{V}_{\t_1},\wh\vph_{\t_1}\big)\big|_{\Si}=(V,\vph)
\quad\hbox{and}\quad
\big(\wh{V}_{0r},\wh\vph_{0r}\big)\!\equiv\!\big(\wh{V},\wh\vph\big)\big|_{\P^1_r}=
\big(\wt{V},\wt\vph\big)\big|_{\P^1_r}\!\equiv\!\big(\wt{V}_{0r},\wt\vph_{0r}\big)
~~\forall\,r\!\ge\!2.$$

\vspace{.15in}

Let $U_1'\!\subset\!U_1$ be a $\si$-invariant tubular neighborhood of~$S^1_1$
so that $\ov{U_1'}\!\subset\!U_1$ is a closed tubular neighborhood
of~$S^1_1$.
Let
\BE{gencOpincg0_e2}\wh\cC_0 \equiv \big(\wh\Si_0,\wh\si_0\big)
\qquad\hbox{and}\qquad
\cC_0' \equiv \big(\Si_0',\si_0'\big)\EE
be the nodal symmetric  surfaces obtained from $(\Si,\si)$
by collapsing each of the boundary components of 
$\ov{U_1'},\ov{U_2},\ldots,\ov{U_N}$ and 
$\ov{U_1'},\ov{U_1},\ov{U_2},\ldots,\ov{U_N}$, respectively,
to a single point; see the left diagrams in Figure~\ref{OSpinDfn_fig}.
We define $\P^1_r\!\subset\!\Si_0'$ and~$\tau_r'$ as below~\eref{fobfxfn_e}.
The choice of half-surface $\Si^b$ of~$\Si$ again determines a distinguished half-surface
$\bD^2_{r+}$ of $(\P^1_r,\tau_r)$ for $r\!\in\![N]$.
Let $\P^1_{\pm}\!\subset\!\Si_0'$ be the irreducible component sharing a node~$\nod^{\pm}$ 
with~$\bD^2_{1\pm}$.
We label the remaining nodes of~$\Si_0'$ by~$\nod_r^{\pm}$ so that $\nod_1^{\pm}\!\in\!\P^1_{\pm}$
and $\nod_r^{\pm}\!\in\!\bD^2_{r\pm}$ for $r\!\ge\!2$.
We denote by
$\Si_0'^{\C}$ the union of the remaining irreducible components of~$\Si_0'$
and by $\Si_0'^+\!\subset\!\Si_0'^{\C}$ its half-surface distinguished by~$\Si^b$.\\

\begin{figure}
\begin{pspicture}(-1,-6.7)(10,5.3)
\psset{unit=.4cm}
\psellipse(5,5)(2,1)\pscircle(11,5){2}
\psarc[linewidth=.03](5,8.46){4}{240}{300} 
\psarc[linewidth=.03,linestyle=dashed](5,1.54){4}{60}{120}
\psarc[linewidth=.03](11,8.46){4}{240}{300} 
\psarc[linewidth=.03,linestyle=dashed](11,1.54){4}{60}{120}
\psarc(5,7){1}{180}{360}\psarc(5,3){1}{0}{180}
\psarc(11,8){1}{180}{360}\psarc(11,2){1}{0}{180}
\psarc(8,8){2}{0}{180}\psarc(8,2){2}{180}{360}
\psarc(8,8){4}{0}{180}\psarc(8,2){4}{180}{360}
\pscircle*(5,6){.2}\pscircle*(5,4){.2}
\pscircle*(11,7){.2}\pscircle*(11,3){.2}
\psline(4,7)(4,8)\psline(6,7)(6,8)\psline(4,3)(4,2)\psline(6,3)(6,2)
\psarc[linewidth=.04](8,11.7){1}{210}{330}\psarc[linewidth=.04](8,10.3){1}{30}{150} 
\psarc[linewidth=.04](8,-.3){1}{210}{330}\psarc[linewidth=.04](8,-1.7){1}{30}{150} 
\rput(2.8,7.2){\sm{$\nod_{\t;1}^+$}}\rput(2.8,2.8){\sm{$\nod_{\t;1}^-$}}
\psline[linewidth=.03]{->}(2.9,6.5)(4.4,6.05)
\psline[linewidth=.03]{->}(2.8,3.4)(4.4,3.95)
\rput(11,6.2){\sm{$\nod_{\t;2}^+$}}\rput(11,3.8){\sm{$\nod_{\t;2}^-$}}
\rput(2.2,5){$\P^1_{\t;1}$}\rput(13.9,5){$\P^1_{\t;2}$}
\rput(4.2,11){$\Si_{\t}'^+$}\rput(4.1,-1){$\Si_{\t}'^-$}
\psline{<->}(1,9)(1,1)\rput(.2,5){$\si_{\t}'$} 
\pscircle(23,5){2}\pscircle(29,5){2}
\psarc[linewidth=.03](23,8.46){4}{240}{300} 
\psarc[linewidth=.03,linestyle=dashed](23,1.54){4}{60}{120}
\psarc[linewidth=.03](29,8.46){4}{240}{300} 
\psarc[linewidth=.03,linestyle=dashed](29,1.54){4}{60}{120}
\psarc(23,8){1}{180}{360}\psarc(23,2){1}{0}{180}
\psarc(29,8){1}{180}{360}\psarc(29,2){1}{0}{180}
\psarc(26,8){2}{0}{180}\psarc(26,2){2}{180}{360}
\psarc(26,8){4}{0}{180}\psarc(26,2){4}{180}{360}
\pscircle*(23,7){.2}\pscircle*(23,3){.2}
\pscircle*(29,7){.2}\pscircle*(29,3){.2}
\psarc[linewidth=.04](26,11.7){1}{210}{330}\psarc[linewidth=.04](26,10.3){1}{30}{150} 
\psarc[linewidth=.04](26,-.3){1}{210}{330}\psarc[linewidth=.04](26,-1.7){1}{30}{150} 
\rput(23,6.2){\sm{$\nod_{\t;1}^+$}}\rput(23,3.8){\sm{$\nod_{\t;1}^-$}}
\rput(29,6.2){\sm{$\nod_{\t;2}^+$}}\rput(29,3.8){\sm{$\nod_{\t;2}^-$}}
\rput(20,5){$\P^1_{\t;1}$}\rput(32,5){$\P^1_{\t;2}$}
\rput(22.2,11){$\Si_{\t}'^+$}\rput(22.2,-1){$\Si_{\t}'^-$}
\psline{<->}(33,9)(33,1)\rput(33.8,5){$\si_{\t}'$}
\psline{->}(18,-8)(18,6)\psline{->}(18,-8)(35,-8)
\rput(33,-9){\sm{$\De_{\R;1}'^*\!\subset\!\C^2$}}
\rput{90}(17,4){\sm{$\De_{\R;2}'^*\!\subset\!\C^2$}}\rput(25,-5){$\De_{\R}'^*$}
\psline[linewidth=.03]{->}(28,-2.5)(28,-7.2)
\psline[linewidth=.03]{->}(12.5,0)(17,0)
\pscircle*(18,-8){.2}\pscircle*(22,-8){.2}
\pscircle*(18,-4){.2}\rput(18.8,-4){$\t_0'$}
\pscircle*(22,-4){.2}\rput(22.1,-7.2){$\t_0$}
\rput(21.7,-4.5){$\t_1$}\rput(23.5,-3.2){\sm{$(\Si,\si,\fj)$}}
\rput(17.2,-4){\sm{$\wh\cC_0$}}
\psellipse(5,-10)(2,1)\pscircle(11,-10){2}
\pscircle(5,-8.5){.5}\pscircle(5,-11.5){.5}
\psarc[linewidth=.03](5,-6.54){4}{240}{300} 
\psarc[linewidth=.03,linestyle=dashed](5,-13.46){4}{60}{120}
\psarc[linewidth=.03](11,-6.54){4}{240}{300} 
\psarc[linewidth=.03,linestyle=dashed](11,-13.46){4}{60}{120}
\psarc(5,-7){1}{180}{360}\psarc(5,-13){1}{0}{180}
\psarc(11,-7){1}{180}{360}\psarc(11,-13){1}{0}{180}
\psarc(8,-7){2}{0}{180}\psarc(8,-13){2}{180}{360}
\psarc(8,-7){4}{0}{180}\psarc(8,-13){4}{180}{360}
\pscircle*(5,-8){.2}\pscircle*(5,-12){.2}
\pscircle*(5,-9){.2}\pscircle*(5,-11){.2}
\pscircle*(11,-8){.2}\pscircle*(11,-12){.2}
\psarc[linewidth=.04](8,-3.3){1}{210}{330}\psarc[linewidth=.04](8,-4.7){1}{30}{150} 
\psarc[linewidth=.04](8,-15.3){1}{210}{330}\psarc[linewidth=.04](8,-16.7){1}{30}{150} 
\rput(5,-7.3){\sm{$\nod_1^+$}}\rput(5,-12.7){\sm{$\nod_1^-$}}
\rput(11,-7.3){\sm{$\nod_2^+$}}\rput(11,-12.7){\sm{$\nod_2^-$}}
\rput(4.9,-10){\sm{$S^1_1$}}\rput(10.9,-10){\sm{$S^1_2$}}
\rput(7.6,-9.2){\sm{$\bD^2_{1+}$}}\rput(9,-7.9){\sm{$\bD^2_{2+}$}}
\rput(7.6,-10.9){\sm{$\bD^2_{1-}$}}\rput(9.1,-11.8){\sm{$\bD^2_{2-}$}}
\rput(6.1,-8.5){$\P^1_{\!+}$}\rput(6.1,-11.5){$\P^1_{\!-}$}
\rput(2.8,-7.5){\sm{$\nod^+$}}\rput(2.8,-12.5){\sm{$\nod^-$}}
\psline[linewidth=.03]{->}(3,-8.2)(4.4,-8.8)
\psline[linewidth=.03]{->}(3,-12)(4.4,-11.2)
\rput(2.3,-10){$\P^1_1$}\rput(13.7,-10){$\P^1_2$}
\rput(4.2,-4){$\Si_0'^+$}\rput(4.2,-16){$\Si_0'^-$}
\psline{<->}(1,-6)(1,-14)\rput(.2,-10){$\si_0'$}
\rput(14,-13){$\cC_0'$}
\psline[linewidth=.03]{->}(13,-8)(17.3,-8)
\end{pspicture}
\caption{Deformations of~$\cC_0'$ in the family $(\cU',\wt\fc')$ over 
$\De'\!\subset\!\C^{2N+2}$ for the proof of Proposition~\ref{gencOpincg_prp}\ref{gencOpincg_it0};
$\t_1\!\in\!\De_{\R}'^*$ does not lie in the span of $\De_{\R;1}'$ and $\De_{\R;2}'$ unless $N\!=\!1$}
\label{OSpinDfn_fig}
\end{figure}

Let  $(\cU',\wt\fc')$ be a flat family of deformations of~$\cC_0'$ as in~\eref{cUsymmdfn_e}
over the unit ball $\De'\!\subset\!\C^{2N+2}$ around the origin 
satisfying the first condition in~\eref{Sisit1cond_e} for some $\t_1\!\in\!\De_{\R}'$ so that 
$$\cC_{\t_0}'\equiv \big(\Si_{\t_0}',\si_{\t_0}'\big)
\qquad\hbox{and}\qquad
\cC_{\t_0'}'\equiv \big(\Si_{\t_0'}',\si_{\t_0'}'\big)$$
for some $\t_0,\t_0'\!\in\!\De_{\R}'$ are the symmetric surfaces~$\cC_0$ 
as below~\eref{fobfxfn_e} 
and $\wh\cC_0$ as in~\eref{gencOpincg0_e2}
obtained from~$\cC_0'$ by smoothing the conjugate pair of nodes~$\nod^{\pm}$
and  the conjugate pair of nodes~$\nod_1^{\pm}$, respectively.
We denote~by 
$$\De_{\R;1}'^*,\De_{\R;2}'^*\subset\De_{\R}'$$
the subspaces parametrizing all symmetric surfaces obtained from~$\cC_0'$
by smoothing the conjugate pair of nodes~$\nod^{\pm}$
and  the conjugate pair of nodes~$\nod_1^{\pm}$, respectively, and by 
\hbox{$\De_{\R;1}',\De_{\R;2}'\!\subset\!\De_{\R}'$} their closures.\\

Let $s_1'^{\R},\ldots,s_N'^{\R}$ be sections of $\cU'^{\wt\fc'}$ over~$\De_{\R}'$ and
$(\wt{V}',\wt\vph')$ be a rank~1 real bundle pair over $(\cU',\wt\fc')$ so~that
\begin{gather}\notag
s_r'^{\R}(\t_1)=x_r\quad\forall\,r\!\in\![N],\\
\label{gencOevComp_e4} \quad
\big(\wt{V}_{\t_0}',\wt\vph_{\t_0}'\big)=\big(\wt{V}_0,\wt\vph_0\big), \quad
\big(\wt{V}'_{\t_0'},\wt\vph'_{\t_0'}\big)=
\big(\wh{V}_0,\wh\vph_0\big), \quad 
\big(\wt{V}_{\t_1}',\wt\vph_{\t_1}'\big)=(V,\vph).
\end{gather}
For each $\t\!\in\!\De_{\R}'$, $\fo_{\x}$ determines a tuple
$$\fo_{\x;\t}'\equiv\big(\fo_{x_r;\t}'\big)_{S^1_r\in\pi_0(\Si^{\si})}$$
of orientations of $\wt{V}'^{\wt\vph'}$ at $s_r'^{\R}(\t)$ via 
the line bundles $s_r'^{\R*}\wt{V}'^{\wt\vph'}$ over~$\De_{\R}'$.\\

Let \hbox{$\cD'\!=\!\{D_{\t}'\}$} be a family of real CR-operators on 
$(\wt{V}_{\t}',\wt\vph_{\t}')$  as in~\eref{DVvphdfn_e} 
so that $D_{\t_1}'$ is the operator~$D$ in the statement of the proposition.
For each $r\!\in\![N]$, we denote by $(\wt{V}'_{0r},\wt\vph'_{0r})$ 
the restriction of $(\wt{V}'_0,\wt\vph'_0)$ to~$\P^1_r$ and 
by $D_{0r}'$ the real CR-operator on  $(\wt{V}'_{0r},\wt\vph'_{0r})$ induced by~$D_0'$.
Let $D_{01}^+$ and $D_0^+$ be the real CR-operators on~$\wt{V}'_0|_{\P^1_+}$ 
and~$\wt{V}'_0|_{\Si_0'^+}$, respectively, induced by~$D_0'$.
The exact triple 
\begin{gather*}
0\lra D_0' \lra 
\bigoplus_{r\in[N]}\!\!\!D_{0r}' \!\oplus\!
D_{01}'^+\!\oplus\!D_0'^+\lra 
\bigoplus_{r\in[N]}\!\!\!\wt{V}_{\nod_r^+}'\oplus\!\wt{V}_{\nod^+}'\lra 0, \\
\big(\xi_{\Si_0'^-},\xi_{\P_-^1},(\xi_r)_{r\in[N]},\xi_{\P_+^1},\xi_{\Si_0'^+}\big)\lra
\big((\xi_r)_{r\in[N]},\xi_{\P_+^1},\xi_{\Si_0'^+}\big), \\
\big((\xi_r)_{r\in[N]},\xi_{\P_+^1},\xi_{\Si_0'^+}\big)\lra
\big(\xi_{\P_+^1}(\nod^+)\!-\!\xi_1(\nod^+),
\xi_{\Si_0'^+}(\nod_1^+)\!-\!\xi_{\P_+^1}(\nod_1^+),
 \big(\xi_{\Si_0'^+}(\nod_r^+)\!-\!\xi_r(\nod_r^+)\!\big)_{r>1}\big), 
\end{gather*}
of Fredholm operators then determines an isomorphism
$$\la(D_0')\!\otimes\!
\bigotimes_{r\in[N]}\!\!
\la\big(\wt{V}'_{\nod_r^+}\big)\otimes\!\la\big(\wt{V}'_{\nod^+}\big)\approx
\bigotimes_{r\in[N]}\!\!\la\big(D_{0r}'\big)
\!\otimes\! \la\big(D_{01}'^+\big)\!\otimes\!\la\big(D_0'^+\big).$$
The orientations $\fo(\wt{V}'_{0r},\wt\vph'_{0r};\fo_{x_r;0}')$  of $D_{0r}'$
and the complex orientations of $D_{01}'^+$, $D_0'^+$, $\wt{V}'_{\nod_r^+}$, and
$\wt{V}'_{\nod^+}$  determine an orientation of $D'_0$
via the above isomorphism and 
thus an orientation~$\fo_{\cU'}(\fo_{\x})$ of the line bundle~$\la(\cD')$ 
over~$\De'_{\R}$.\\

The restriction of $(\cU',\wt\fc')$ to $\De_{\R;1}'$ is the product of a family smoothing 
$$\P^1_{0;1}\equiv \P^1_-\!\cup\!\P^1_1\!\cup\!\P^1_+$$
into a single irreducible component $\P^1_{\t;1}\!\subset\!\Si_{\t}'$
as above Proposition~\ref{cOpincg_prp} with 
the union of remaining irreducible components of~$\Si_0'$;
see the right diagram in Figure~\ref{OSpinDfn_fig}.
The other components of~$\Si_{\t}'$ with $\t\!\in\!\De_{\R;1}^*$ 
and its nodes correspond to the components $\P^1_r$ of $\Si_0'$ with $r\!\ge\!2$ and~$\Si_0'^{\pm}$
and to the nodes~$\nod_r^{\pm}$ with $r\!\ge\!1$;
we denote them by~$\P^1_{\t;r}$, $\Si_{\t}'^{\pm}$, and~$\nod_{\t;r}^{\pm}$, respectively.
For $\t\!\in\!\De_{\R;1}'$ and $r\!\in\![N]$, we denote by $(\wt{V}'_{\t;r},\wt\vph'_{\t;r})$ 
the restriction of $(\wt{V}'_{\t},\wt\vph'_{\t})$ to~$\P^1_{\t;r}$ and 
by $D_{\t;r}'$ the real CR-operator on  $(\wt{V}'_{\t;r},\wt\vph'_{\t;r})$ induced by~$D_{\t}'$.
Let $D_{\t}^+$ be the real CR-operator on~$\wt{V}'_{\t}|_{\Si_{\t}'^+}$ induced by~$D_{\t}'$.
The exact triple~\eref{Cndses_e} of Fredholm operators with~$D_0$ replaced 
by~$D_{0;1}'$ and 
the exact triple~\eref{geCndses_e} of Fredholm operators with~$D_0$ replaced
by~$D'_{\t_0}$ induce isomorphisms
\BE{gencOevComp_e7}\begin{split}
\la\big(D_{0;1}'\big)\!\otimes\!\la\big(\wt{V}'_{\nod^+}\big)&\approx 
\la\big(D_{01}'\big)\!\otimes\!\la\big(D_{01}'^+\big),\\
\la(D_{\t_0}')\!\otimes\!\bigotimes_{r\in[N]}\!\!
\la\big(\wt{V}'_{\nod_{\t_0;r}^+}\big)&\approx 
\bigotimes_{r\in[N]}\!\!\la\big(D_{\t_0;r}'\big)
\otimes\!\la\big(D_{\t_0}'^+\big).
\end{split}\EE
The orientation $\fo(\wt{V}'_{01},\wt\vph'_{01};\fo_{x_1;0}')$
of $D_{01}'$ and the complex orientations of $D_{01}'^+$ and $\wt{V}'_{\nod^+}$ 
determine an orientation of $D_{0;1'}$ via the first isomorphism in~\eref{gencOevComp_e7} 
and thus an orientation of $D_{\t;1}'$ for each $\t\!\in\!\De_{\R;1}'$;
the induced orientation of~$D_{0;1}'$ is the $\C$-split orientation in the terminology 
of Proposition~\ref{rk1Cdegen_prp}.
By this proposition and the first assumption in~\eref{gencOevComp_e4}, 
the induced orientation of~$D_{\t_0;1}'$ is thus the intrinsic orientation
$$\fo\big(\wt{V}'_{\t_0;1},\wt\vph'_{\t_0;1};\fo_{x_1;\t_0}'\big)
=\fo\big(\wt{V}_{01},\wt\vph_{01};\fo_{x_1;0}\big).$$
Since the (real) index of~$D_{01}'^+$
the (real) dimension of~$\wt{V}'_{\nod^+}$ are even, 
this implies that the second isomorphism in~\eref{gencOevComp_e7} is orientation-preserving 
with respect to the restriction of the orientation~$\fo_{\cU'}(\fo_{\x})$ to $\la(D_{\t_0}')$,
the intrinsic orientations
$$\fo\big(\wt{V}'_{\t_0;r},\wt\vph'_{\t_0;r};\fo_{x_r;\t_0}'\big)
=\fo\big(\wt{V}_{0r},\wt\vph_{0r};\fo_{x_r;0}\big)$$
of $D_{\t_0;r}'\!=\!D_{0r}$, 
and the complex orientations of~$D_{\t_0}'^+$ and~$\wt{V}'_{\nod_{\t_0;r}^+}$.\\

The restriction of $(\cU',\wt\fc')$ to $\De_{\R;2}'$ is the product of a family smoothing 
the conjugate~pair
$$\Si_{0;\pm}'\equiv \P^1_{\pm}\cup\Si_0'^{\pm}\subset\Si_0'$$
into a conjugate pair $\Si_{\t;\pm}'\!\subset\!\Si_{\t}'$ with 
the union of remaining irreducible components of~$\Si_0'$;
see the top left diagram in Figure~\ref{OSpinDfn_fig}.
The other components of~$\Si_{\t}'$ with $\t\!\in\!\De_{\R;2}^*$ 
and its nodes correspond to the components $\P^1_r$ of $\Si_0'$ with $r\!\ge\!1$
and to the nodes~$\nod^{\pm}$ and $\nod_r^{\pm}$ with $r\!\ge\!2$;
we denote them by~$\P^1_{\t r}$, $\nod_{\t;1}^{\pm}$, and~$\nod_{\t;r}^{\pm}$, respectively.
For $\t\!\in\!\De_{\R;2}'$ and $r\!\in\![N]$, we denote by $(\wt{V}'_{\t r},\wt\vph'_{\t r})$ 
the restriction of $(\wt{V}'_{\t},\wt\vph'_{\t})$ to~$\P^1_{\t r}$ and 
by $D_{\t r}'$ the real CR-operator on  $(\wt{V}'_{\t r},\wt\vph'_{\t r})$ induced by~$D_{\t}'$.
Let $D_{\t;+}$ be the real CR-operator on~$\wt{V}'_{\t}|_{\Si_{\t;+}'}$ induced by~$D_{\t}'$.
The exact~triple 
$$0\lra D'_{0;+}\lra D_{01}'^+\!\oplus\!D_0'^+\lra \wt{V}'_{\nod_1^+}\lra 0, 
~~ \big(\xi_{\P^1_+},\xi_{\Si_0'^+}\big)\lra
\xi_{\Si_0'^+}(\nod_1^+)\!-\!\xi_{\P^1_+}(\nod_1^+),$$
of Fredholm operators and the exact triple~\eref{geCndses_e} of 
Fredholm operators with~$D_0$ replaced by~$D'_{\t_0'}$ induce isomorphisms
\BE{gencOevComp_e9}\begin{split}
\la\big(D_{0;+}'\big)\!\otimes\!\la\big(\wt{V}'_{\nod_1^+}\big)
&\approx 
\la\big(D_{01}'^+\big)\!\otimes\!\la\big(D_0'^+\big),\\
\la(D_{\t_0'}')\!\otimes\!\bigotimes_{r\in[N]}\!\!
\la\big(\wt{V}'_{\nod_{\t_0';r}^+}\big)&\approx 
\bigotimes_{r\in[N]}\!\!\la\big(D_{\t_0' r}'\big)
\otimes\!\la\big(D_{\t_0';+}'\big).
\end{split}\EE
The complex orientations of $D_{01}'^+$, $D_0'^+$, and $\wt{V}'_{\nod_1^+}$
determine an orientation on $D_{0;+}'$ via  the first isomorphism in~\eref{gencOevComp_e9}
and thus an orientation of $D_{\t;+}'$ for each $\t\!\in\!\De_{\R;2}'$.
The latter is the complex orientation of this real CR-operator.
Since the dimension of $\wt{V}'_{\nod_1^+}$ is even,
this statement and the second assumption in~\eref{gencOevComp_e4}
imply that the second isomorphism in~\eref{gencOevComp_e9} is orientation-preserving 
with respect to the restriction of the orientation~$\fo_{\cU'}(\fo_{\x})$ to $\la(D_{\t_0'}')$,
the intrinsic orientations 
$$\fo(\wt{V}'_{\t_0'r},\wt\vph'_{\t_0'r};\fo_{x_r;\t_0'}')
=\fo\big(\wh{V}_{0r},\wh\vph_{0r};\fo_{x_r;0}\big)$$ 
of $D_{\t_0' r}'$, and the complex orientations of~$D_{\t_0';+}'$ 
and~$\wt{V}'_{\nod_{\t_0';r}^+}$.\\

By the conclusion above regarding the second isomorphism in~\eref{gencOevComp_e7}
and the last assumption in~\eref{gencOevComp_e4},
the restriction of the orientation~$\fo_{\cU'}(\fo_{\x})$ to $\la(D_{\t_1}')\!=\!\la(D)$
is the orientation $\fo(V,\vph;\fo_{\x})$ of~$D$ constructed above Proposition~\ref{gencOpincg_prp}
from the real bundle pair $(\wt{V},\wt\vph)$ over the nodal curve~$\cC_0$ below~\eref{fobfxfn_e}.
By the conclusion above regarding the second isomorphism in~\eref{gencOevComp_e9}
and the last assumption in~\eref{gencOevComp_e4},
this restriction 
is the orientation $\fo(V,\vph;\fo_{\x})$ constructed above Proposition~\ref{gencOpincg_prp}
from the real bundle pair $(\wh{V},\wh\vph)$ over the nodal curve~$\cC_0$ below~\eref{fobfxfn_e}.
Thus, the rank~1 real bundle pairs $(\wt{V},\wt\vph)$ and $(\wh{V},\wh\vph)$ over~$(\cU,\wt\fc)$ 
determine the same orientation of~$D$.
\end{proof}

\begin{proof}[{\bf{\emph{Proof of Proposition~\ref{genrk1Cdegen_prp}}}}]
We recall that the nodal symmetric surface $(\Si_0,\si_0)$ 
with one conjugate pair of nodes $\nod^{\pm}$ in this case is
obtained from its normalization  $(\wt\Si_0,\wt\si_0)$ by identifying
$z_1^+$ with $z_2^+$ into the node~$\nod^+$ and $z_1^-$ with $z_2^-$ into the node~$\nod^-$.
In particular, the complement of~$\nod^+,\nod^-$ in~$\Si$ is canonically identified
with the complement of $z_1^+,z_1^-,z_2^+,z_2^-$ in~$\wt\Si$.\\

For each $r\!\in\![N]$, let $U_r\!\subset\!\Si_0,\wt\Si_0$ 
be a $\si$-invariant tubular neighborhood of~$S^1_r$
so that the closures~$\ov{U_r}$ of such neighborhoods are disjoint closed tubular neighborhoods
of~$S^1_r$ not containing the points $\nod^+,\nod^-$ and $z_1^+,z_1^-,z_2^+,z_2^-$.
Let \hbox{$\cC_0'\!\equiv\!(\Si_0',\si_0')$}
be the nodal symmetric surface obtained from $(\Si_0,\si_0)$
by collapsing each of the boundary components of each~$\ov{U_r}$ to a single point;
see the bottom left diagram in Figure~\ref{genCdegen_fig}.
We~define 
$$\P^1_r,\Si_0'^{\C}\subset\Si_0', \qquad \bD^2_{r\pm}\subset\P^1_r, \qquad
\Si_0'^{\pm}\subset\Si_0'^{\C}, \qquad\hbox{and}\qquad \nod_r^{\pm}\in \bD^2_{r\pm},\Si_0'^{\pm}$$
as below~\eref{fobfxfn_e}.
In this case, $\nod^{\pm}\!\in\!\Si_0'^{\pm}$.
We denote by~$\wt\Si_0'^+$ and $\wt\Si_0'$ the closed surfaces obtained from~$\Si_0'^+$ and~$\Si_0'$,
respectively, by replacing the conjugate pair of nodes~$\nod^{\pm}$ with two conjugate
pairs~$z_{0;1}^{\pm}$ and~$z_{0;2}^{\pm}$ of marked points.\\

\begin{figure}
\begin{pspicture}(-1,-8)(10,5.3)
\psset{unit=.4cm}
\psarc[linewidth=.03](5,6.73){2}{240}{300} 
\psarc[linewidth=.03,linestyle=dashed](5,3.27){2}{60}{120}
\psarc[linewidth=.03](11,6.473){2}{240}{300} 
\psarc[linewidth=.03,linestyle=dashed](11,3.27){2}{60}{120}
\psarc(8,8){2}{0}{180}\psarc(8,2){2}{180}{360}
\psarc(7.65,8){3.65}{90}{180}\psarc(8.35,8){3.65}{0}{90}
\psarc(7.65,2){3.65}{180}{270}\psarc(8.35,2){3.65}{270}{360}
\psline(4,2)(4,8)\psline(6,2)(6,8)\psline(10,2)(10,8)\psline(12,2)(12,8)
\psarc[linewidth=.04](8,11.7){1}{210}{330}\psarc[linewidth=.04](8,10.3){1}{30}{150} 
\psarc[linewidth=.04](8,-.3){1}{210}{330}\psarc[linewidth=.04](8,-1.7){1}{30}{150} 
\pscircle*(8,11.3){.2}\psarc[linewidth=.04](8,11.475){.175}{180}{360}
\pscircle*(8,-1.3){.2}\psarc[linewidth=.04](8,-1.475){.175}{0}{180}
\psarc[linewidth=.04](7.65,11.475){.175}{0}{90}\psarc[linewidth=.04](8.35,11.475){.175}{90}{180}
\psarc[linewidth=.04](7.65,-1.475){.175}{270}{360}\psarc[linewidth=.04](8.35,-1.475){.175}{180}{270}
\rput(8.3,12.3){\sm{$\nod_{\t}^+$}}\rput(8.1,-2.2){\sm{$\nod_{\t}^-$}}
\rput(4.2,11){$\Si_{\t}'^+$}\rput(4.1,-1){$\Si_{\t}'^-$}
\psline{<->}(2,9)(2,1)\rput(1.2,5){$\si_{\t}'$} 
\psline[linewidth=.03]{->}(13,2)(17,0)
\pscircle(23,5){2}\pscircle(29,5){2}
\psarc[linewidth=.03](23,8.46){4}{240}{300} 
\psarc[linewidth=.03,linestyle=dashed](23,1.54){4}{60}{120}
\psarc[linewidth=.03](29,8.46){4}{240}{300} 
\psarc[linewidth=.03,linestyle=dashed](29,1.54){4}{60}{120}
\psarc(23,8){1}{180}{360}\psarc(23,2){1}{0}{180}
\psarc(29,8){1}{180}{360}\psarc(29,2){1}{0}{180}
\psarc(26,8){2}{0}{180}\psarc(26,2){2}{180}{360}
\psarc(26,8){4}{0}{180}\psarc(26,2){4}{180}{360}
\pscircle*(23,7){.2}\pscircle*(23,3){.2}
\pscircle*(29,7){.2}\pscircle*(29,3){.2}
\psarc[linewidth=.04](26,11.7){1}{210}{330}\psarc[linewidth=.04](26,10.3){1}{30}{150} 
\psarc[linewidth=.04](26,-.3){1}{210}{330}\psarc[linewidth=.04](26,-1.7){1}{30}{150} 
\rput(23,6.2){\sm{$\nod_{\t;1}^+$}}\rput(23,3.8){\sm{$\nod_{\t;1}^-$}}
\rput(29,6.2){\sm{$\nod_{\t;2}^+$}}\rput(29,3.8){\sm{$\nod_{\t;2}^-$}}
\rput(20,5){$\P^1_{\t;1}$}\rput(32,5){$\P^1_{\t;2}$}
\rput(22.2,11){$\Si_{\t}'^+$}\rput(22.2,-1){$\Si_{\t}'^-$}
\psline{<->}(33,9)(33,1)\rput(33.8,5){$\si_{\t}'$}
\psline{->}(18,-8)(18,6)\psline{->}(18,-8)(35,-8)
\rput(33,-9){\sm{$\De_{\R;1}'^*\!\subset\!\C^2$}}
\rput{90}(17,3.5){\sm{$\De_{\R;2}'^*\!\subset\!\C^{2N}$}}\rput(25,-5){$\De_{\R}'^*$}
\psline[linewidth=.03]{->}(28,-2.5)(28,-7.2)
\pscircle*(18,-8){.2}\pscircle*(22,-8){.2}\pscircle*(18,-4){.2}
\rput(18.8,-4){$\t_0'$}\rput(22.1,-7.2){$\t_0$}
\rput(16.2,-4){\sm{Prp~\ref{genrk1Cdegen_prp}}}
\pscircle(5,-12){2}\pscircle(11,-12){2}
\psarc[linewidth=.03](5,-8.54){4}{240}{300} 
\psarc[linewidth=.03,linestyle=dashed](5,-15.46){4}{60}{120}
\psarc[linewidth=.03](11,-8.54){4}{240}{300} 
\psarc[linewidth=.03,linestyle=dashed](11,-15.46){4}{60}{120}
\psarc(5,-9){1}{180}{360}\psarc(5,-15){1}{0}{180}
\psarc(11,-9){1}{180}{360}\psarc(11,-15){1}{0}{180}
\psarc(8,-9){2}{0}{180}\psarc(8,-15){2}{180}{360}
\psarc(7.65,-9){3.65}{90}{180}\psarc(8.35,-9){3.65}{0}{90}
\psarc(7.65,-15){3.65}{180}{270}\psarc(8.35,-15){3.65}{270}{360}
\pscircle*(5,-10){.2}\pscircle*(5,-14){.2}\pscircle*(11,-10){.2}\pscircle*(11,-14){.2}
\psarc[linewidth=.04](8,-5.3){1}{210}{330}\psarc[linewidth=.04](8,-6.7){1}{30}{150} 
\psarc[linewidth=.04](8,-17.3){1}{210}{330}\psarc[linewidth=.04](8,-18.7){1}{30}{150} 
\pscircle*(8,-5.7){.2}\psarc[linewidth=.04](8,-5.525){.175}{180}{360}
\pscircle*(8,-18.3){.2}\psarc[linewidth=.04](8,-18.475){.175}{0}{180}
\psarc[linewidth=.04](7.65,-5.525){.175}{0}{90}\psarc[linewidth=.04](8.35,-5.525){.175}{90}{180}
\psarc[linewidth=.04](7.65,-18.475){.175}{270}{360}\psarc[linewidth=.04](8.35,-18.475){.175}{180}{270}
\rput(8.3,-4.7){\sm{$\nod^+$}}\rput(8.1,-19.1){\sm{$\nod^-$}}
\rput(5,-9.3){\sm{$\nod_1^+$}}\rput(5,-14.7){\sm{$\nod_1^-$}}
\rput(11,-9.3){\sm{$\nod_2^+$}}\rput(11,-14.7){\sm{$\nod_2^-$}}
\rput(4.9,-12){\sm{$S^1_1$}}\rput(10.9,-12){\sm{$S^1_2$}}
\rput(7.2,-10.1){\sm{$\bD^2_{1+}$}}\rput(9,-9.9){\sm{$\bD^2_{2+}$}}
\rput(7.3,-13.8){\sm{$\bD^2_{1-}$}}\rput(9.1,-13.8){\sm{$\bD^2_{2-}$}}
\rput(2.3,-12){$\P^1_1$}\rput(13.7,-12){$\P^1_2$}
\rput(4.2,-6){$\Si_0'^+$}\rput(4.2,-18){$\Si_0'^-$}
\psline{<->}(1,-8)(1,-16)\rput(.2,-12){$\si_0'$}
\rput(14,-15){$\cC_0'$}
\psline[linewidth=.03]{->}(13,-10)(17.3,-8.3)
\end{pspicture}
\caption{Deformations of~$\cC_0'$ in the family $(\cU',\wt\fc')$ over 
$\De'\!\subset\!\C^{2N+2}$ for the proof of Proposition~\ref{genrk1Cdegen_prp}}
\label{genCdegen_fig}
\end{figure}

Let  $(\cU',\wt\fc')$ be a flat family of deformations of~$\cC_0'$ as in~\eref{cUsymmdfn_e}
over the unit ball $\De'\!\subset\!\C^{2N+2}$ around the origin so~that 
$$\cC_{\t_0}'\equiv \big(\Si_{\t_0}',\si_{\t_0}'\big)
\qquad\hbox{and}\qquad
\cC_{\t_0'}'\equiv \big(\Si_{\t_0'}',\si_{\t_0'}'\big)$$
for some $\t_0,\t_0'\!\in\!\De_{\R}'$ are the symmetric surfaces~$\cC_0$ 
below~\eref{fobfxfn_e} and in the statement of the proposition
obtained from~$\cC_0'$ by smoothing the conjugate pair of nodes~$\nod^{\pm}$
and the conjugate pairs of nodes~$\nod_r^{\pm}$, respectively.
We denote~by 
$$\De_{\R;1}'^*,\De_{\R;2}'^*\subset\De_{\R}'$$
the subspaces parametrizing all symmetric surfaces 
obtained from~$\cC_0'$ by smoothing the conjugate pair of nodes~$\nod^{\pm}$
and the conjugate pairs of nodes~$\nod_r^{\pm}$, respectively,  and by 
\hbox{$\De_{\R;1}',\De_{\R;2}'\!\subset\!\De_{\R}'$} their closures.\\

Let $s_1'^{\R},\ldots,s_N'^{\R}$ be sections of $\cU'^{\wt\fc'}$ over~$\De_{\R}'$ and
$(V',\vph')$ be a rank~1 real bundle pair over $(\cU',\wt\fc')$ so~that
$$s_r'^{\R}(\t_0')=x_r~~\forall\,r\!\in\![N] \qquad\hbox{and}\qquad
\big(V'_{\t_0'},\vph'_{\t_0'}\big)=(V_0,\vph_0).$$
For each $\t\!\in\!\De_{\R}'$, $\fo_{\x}$ determines a tuple
$$\fo_{\x;\t}'\equiv\big(\fo_{x_r;\t}'\big)_{S^1_r\in\pi_0(\Si^{\si})}$$
of orientations of $V'^{\vph'}$ at $s_r'^{\R}(\t)$ via 
the line bundles $s_r'^{\R*}V'^{\vph'}$ over~$\De_{\R}'$.
We denote by $(\wt{V}_0',\wt\vph_0')$ the lift of $(V_0',\vph_0')$ 
to~$(\wt\Si_0',\wt\si_0')$.\\

Let \hbox{$\cD'\!=\!\{D_{\t}'\}$} be a family of real CR-operators on 
$(V'_{\t},\vph'_{\t})$ as in~\eref{DVvphdfn_e} 
so that $D_{\t_0'}'$ is the operator~$D_0$ in the statement of the proposition.
For each $r\!\in\![N]$, we denote by $(V'_{0r},\vph'_{0r})$ 
the restriction of $(V'_0,\vph'_0)$ to~$\P^1_r$ and 
by $D_{0r}'$ the real CR-operator on  $(V'_{0r},\vph'_{0r})$ induced by~$D_0'$.
Let $D_0'^+$ be the real CR-operator on~$V'_0|_{\Si_0'^+}$ induced by~$D_0'$
and $\wt{D}_0'$ and $\wt{D}_0'^+$ be the lifts of $D_0'$ and $D_0'^+$, respectively,
to real CR-operators on $(\wt{V}_0',\wt\vph_0')$ and~$\wt{V}_0'|_{\wt\Si_0'^+}$.
The exact triple 
\begin{gather*}
0\lra D_0' \lra 
\bigoplus_{r\in[N]}\!\!\!D_{0r}' \!\oplus\!\wt{D}_0'^+\lra 
\bigoplus_{r\in[N]}\!\!\!V_{\nod_r^+}' \oplus\!V_{\nod^+}'\!\lra 0, \\
\big(\xi_-,(\xi_r)_{r\in[N]} ,\xi_+\big)\lra
\big((\xi_r)_{r\in[N]},\xi_+\big), \\
\big((\xi_r)_{r\in[N]},\xi_+\big)\lra
\big(\big(\xi_+(\nod_r^+)\!-\!\xi_r(\nod_r^+)\!\big)_{r\in[N]},
\xi_+(z_2^+)\!-\!\xi_+(z_1^+)\big), 
\end{gather*}
of Fredholm operators then determines an isomorphism
$$\la(D_0')\!\otimes\!\bigotimes_{r\in[N]}\!\!\la\big(V'_{\nod_r^+}\big)
\otimes\!\la\big(V'_{\nod^+}\big)\approx
\bigotimes_{r\in[N]}\!\!\la\big(D_{0r}'\big)\!\otimes\!\la\big(\wt{D}_0'^+\big).$$
The orientations $\fo(V'_{0r},\vph'_{0r};\fo_{x_r;0}')$ of $D_{0r}'$
and the complex orientations of $V'_{\nod_r^+}$, $V'_{\nod^+}$,
and $\wt{D}_0'^+$ determine an orientation of $D'_0$
via the above isomorphism  and 
thus an orientation~$\fo_{\cU'}(\fo_{\x})$ of the line bundle~$\la(\cD')$ 
over~$\De'_{\R}$.\\

The restriction of $(\cU',\wt\fc')$ to $\De_{\R;1}'$ is the product of a family smoothing 
the conjugate~pair~$\Si_0'^{\pm}$
into a conjugate pair  $\Si_{\t}'^{\pm}\!\subset\!\Si_{\t}'$ with 
the union of remaining irreducible components $(\P^1_r,\tau_r)$ of~$\Si_0'$.
The other components of~$\Si_{\t}'$ with $\t\!\in\!\De_{\R;1}^*$ 
and its nodes correspond to the components $\P^1_r$ of $\Si_0'$
and to the nodes~$\nod_r^{\pm}$;
we denote them by~$\P^1_{\t r}$ and~$\nod_{\t;r}^{\pm}$, respectively.
For $\t\!\in\!\De_{\R;1}'$ and $r\!\in\![N]$, we denote by $(V'_{\t r},\vph'_{\t r})$ 
the restriction of $(V'_{\t},\vph'_{\t})$ to~$\P^1_{\t r}$ and 
by $D_{\t r}'$ the real CR-operator on  $(V'_{\t r},\vph'_{\t r})$ induced by~$D_{\t}'$.
Let $D_{\t}'^+$ be the real CR-operator on~$V'_0|_{\Si_{\t}'^+}$ induced by~$D_0'$.
The exact~triple 
$$0\lra D_0'^+\lra \wt{D}_0'^+\lra V'_{\nod^+}\lra 0, 
\qquad \xi\lra\xi\big(z_{2;0}^+\big)\!-\!\xi\big(z_{1;0}^+\big),$$
of Fredholm operators and the exact triple~\eref{geCndses_e} of 
Fredholm operators with~$D_0$ replaced by~$D'_{\t_0}$ induce isomorphisms
\BE{gencOevComp2_e9}\begin{split}
\la\big(D_0'^+\big)\!\otimes\!\la\big(V'_{\nod^+}\big)
&\approx  \la\big(\wt{D}_0'^+\big),\\
\la(D_{\t_0}')\!\otimes\!\bigotimes_{r\in[N]}\!\!
\la\big(\wt{V}'_{\nod_{\t_0;r}^+}\big)&\approx 
\bigotimes_{r\in[N]}\!\!\la\big(D_{\t_0 r}'\big)
\otimes\!\la\big(D_{\t_0}'^+\big).
\end{split}\EE
The complex orientations of $\wt{D}_0'^+$ and $V'_{\nod^+}$
determine an orientation of $D_0'^+$ via the first isomorphism in~\eref{gencOevComp2_e9} 
and thus an orientation of $D_{\t}'^+$ for each $\t\!\in\!\De_{\R;1}'$.
The latter is the complex orientation of~$D_{\t}'^+$.
This implies that 
the second isomorphism in~\eref{gencOevComp2_e9} is orientation-preserving 
with respect to the restriction of the orientation~$\fo_{\cU'}(\fo_{\x})$ to $\la(\wt{D}_{\t_0'}')$,
the intrinsic orientations $\fo(\wt{V}'_{\t_0 r},\wt\vph'_{\t_0 r};\fo_{x_r;\t_0}')$ 
of $D_{\t_0 r}'$,  and the complex orientations of $D_{\t_0}'^+$ and $\wt{V}'_{\nod_{\t_0;r}^+}$.
Thus, the restriction of~$\fo_{\cU'}(\fo_{\x})$ to $\la(D_{\t}')$ with $\t\!\in\!\De_{\R}'^*$
is the intrinsic orientation $\fo(V'_{\t},\vph'_{\t'};\fo_{\x;\t}')$
defined above Proposition~\ref{gencOpincg_prp}.
This in turn implies that the restriction of~$\fo_{\cU'}(\fo_{\x})$ to $\la(D_{\t_0'}')$
is the limiting orientation~$\fo_0'(\fo_{\x})$ of \hbox{$D_{\t_0'}'\!=\!D_0$}
in~\eref{rk1Clim_e4}.\\

The restriction of $(\cU',\wt\fc')$ to $\De_{\R;2}'$ is a family deforming~$(\Si_0',\si_0')$
to the symmetric surface in the statement of the proposition and in~\eref{cCdegC_e}.
For each $\t\!\in\!\De_{\R;2}'$, we denote by $(\wt\Si_{\t}',\wt\si_{\t}')$ 
the symmetric surface obtained from~$(\Si_{\t}',\si_{\t}')$ by replacing 
the conjugate pair of nodes~$\nod_{\t}^{\pm}$ corresponding to~$\nod^{\pm}$ 
with two conjugate pairs~$z_{\t;1}^{\pm}$ and~$z_{\t;2}^{\pm}$ of marked points,
by $(\wt{V}_{\t}',\wt\vph_{\t}')$ the lift of $(V_{\t}',\vph_{\t}')$ to $(\wt\Si_{\t}',\wt\si_{\t}')$,
and by~$\wt{D}_{\t}'$ the real CR-operator on $(V_{\t}',\vph_{\t}')$
induced by~$D_{\t}'$.
The exact triple~\eref{geCndses_e} of Fredholm operators with~$D_0$ replaced
by~$\wt{D}'_0$ induces an isomorphism
\BE{gencOevComp2_e7}
\la(\wt{D}_0')\!\otimes\!\bigotimes_{r\in[N]}\!\!
\la\big(V'_{\nod_r^+}\big)\approx 
\bigotimes_{r\in[N]}\!\!\la\big(D_{0r}'\big)
\otimes\!\la\big(\wt{D}_0'^+\big).\EE
The orientations $\fo(V_{0r}',\vph_{0r}';\fo_{x_r;0}')$ of $D_{0;r}'$
and the complex orientations of $\wt{D}_0'^+$ and $V'_{\nod_r^+}$ 
determine an orientation of $\wt{D}_0'$ via the isomorphism~\eref{gencOevComp2_e7} 
and thus an orientation of $\wt{D}_{\t}'$ for each $\t\!\in\!\De_{\R;2}'$.
By definition, this orientation for $\t\!=\!\t_0'$ is the orientation 
$\wt\fo_0(\fo_{\x})$ of $\wt{D}_{\t_0'}'\!=\!\wt{D}_0$ in~\eref{Cdegenorient_e3}. 
Since the dimension of~$V'_{\nod^+}$ is even,  this implies that
the restriction of~$\fo_{\cU'}(\fo_{\x})$ to $\la(D_{\t_0'}')$
is the intrinsic orientation~$\fo_0(\fo_{\x})$ of $D_{\t_0'}'\!=\!D_0$ in~\eref{rk1Csplit_e4}.
Combining this with the conclusion of the previous paragraph, we conclude that
$\fo_0(\fo_{\x})\!=\!\fo_0'(\fo_{\x})$.
\end{proof}

\begin{proof}[{\bf{\emph{Proof of Proposition~\ref{genrk1H3degen_prp}}}}]
In this case, the nodal symmetric surface $(\Si_0,\si_0)$ with one H3~node is
obtained from its normalization  $(\wt\Si_0,\wt\si_0)$ by identifying
marked points $\nod_1\!\in\!S^1_{\bu1}$ and $\nod_2\!\in\!S^1_{\bu2}$
into the node~$\nod$.
In light of Proposition~\ref{gencOpincg_prp}\ref{gencOpincg_it1b},
we can assume that $r_{\bu}\!=\!N$.\\

For $r\!=\!1,2$, let $U_{\bu r}\!\subset\!\Si_{\bu r}$ 
be a $\si$-invariant tubular neighborhood of~$S^1_{\bu r}$
so that its closure~$\ov{U_{\bu r}}$ is a closed tubular neighborhood of~$S^1_{\bu r}$.
Let $\cC_0'\!\equiv\!(\Si_0',\si_0')$ be the nodal symmetric surface obtained from $(\Si_0,\si_0)$
by collapsing each of the boundary components of~$\ov{U_{\bu1}}$ and~$\ov{U_{\bu2}}$ 
to a single point; see the bottom left diagram in Figure~\ref{genH3degen_fig}.
We denote by $\P^1_{\bu r}\!\subset\!\Si_0'$ the irreducible component containing~$S^1_{\bu r}$,
by~$\Si_0'^c\!\subset\!\Si_0'$ the union of the remaining irreducible components,
and by $\tau_{\bu r}'$ and $\si_0'^c$ the restrictions of~$\si_0'$ to~$\P^1_{\bu r}$ and~$\Si_0'^c$,
respectively.
The decorated structure on~$(\Si_0,\si_0)$ induces decorated structures on
the symmetric surfaces $(\P^1_{\bu r},\tau_{\bu r}')$ and~$(\Si_0'^c,\si_0'^c)$.
We denote by $\nod_{\bu r}^+\!\in\!\P_{\bu r}^1$ the node contained 
in the interior of the distinguished half-surface~$\bD^2_{\bu r}$ of~$\P_{\bu r}^1$
and by $\nod_{\bu r}^-\!\in\!\P_{\bu r}^1$ its conjugate.\\

\begin{figure}
\begin{pspicture}(-1,-7.5)(10,4.5)
\psset{unit=.4cm}
\psellipse(5.5,3)(4.5,7)\psellipse(12,3)(2,6)
\psarc[linewidth=.04](8,3){4}{130}{230}\psarc[linewidth=.04](2,3){4}{-50}{50} 
\psarc[linewidth=.03](2.5,5.6){3}{240}{300} 
\psarc[linewidth=.03,linestyle=dashed](2.5,.4){3}{60}{120}
\psarc[linewidth=.03](8,6.46){4}{240}{300} 
\psarc[linewidth=.03,linestyle=dashed](8,-.46){4}{60}{120}
\psarc[linewidth=.03](12,6.46){4}{240}{300} 
\psarc[linewidth=.03,linestyle=dashed](12,-.46){4}{60}{120}
\pscircle*(10,3){.2}\rput(12,-1){\sm{$\nod_{\t}$}}
\psline[linewidth=.03]{->}(12,-.5)(10.3,2.7)
\psarc[linewidth=.04](5.5,9.7){1}{210}{330}\psarc[linewidth=.04](5.5,8.3){1}{30}{150} 
\psarc[linewidth=.04](5.5,-2.3){1}{210}{330}\psarc[linewidth=.04](5.5,-3.7){1}{30}{150} 
\rput(2.5,4){\sm{$S^1_1$}}
\rput(8.1,4.1){\sm{$S^1_{\!\bu1}$}}\rput(12,4.1){\sm{$S^1_{\!\bu2}$}}
\psline{<->}(0,7)(0,-1)\rput(-.8,3){$\si_{\t}'$}
\psline[linewidth=.03]{->}(14.2,1)(17.5,0)
\psarc(27,3){2}{90}{270}\psarc(31,3){2}{-90}{90}\pscircle(31,7){2}\pscircle(31,-1){2}
\psline(27,5)(31,5)\psline(27,1)(31,1)
\psarc[linewidth=.03](21,4.73){2}{240}{300} 
\psarc[linewidth=.03,linestyle=dashed](21,1.27){2}{60}{120}
\psarc[linewidth=.03](27,6.46){4}{240}{270} 
\psarc[linewidth=.03,linestyle=dashed](27,-.46){4}{90}{120}
\psarc[linewidth=.03](31,6.46){4}{270}{300} 
\psarc[linewidth=.03,linestyle=dashed](31,-.46){4}{60}{90}
\psline[linewidth=.03](27,2.46)(31,2.46)
\psline[linewidth=.03,linestyle=dashed](27,3.54)(31,3.54)
\psline(20,6)(20,0)\psline(22,6)(22,0)
\psarc(27,6){1}{180}{360}\psarc(27,0){1}{0}{180}
\psarc(24,6){2}{0}{180}\psarc(24,0){2}{180}{360}
\psarc(24,6){4}{0}{180}\psarc(24,0){4}{180}{360}
\pscircle*(27,5){.2}\pscircle*(27,1){.2}\pscircle*(31,5){.2}\pscircle*(31,1){.2}
\psarc[linewidth=.04](24,9.7){1}{210}{330}\psarc[linewidth=.04](24,8.3){1}{30}{150} 
\psarc[linewidth=.04](24,-2.3){1}{210}{330}\psarc[linewidth=.04](24,-3.7){1}{30}{150} 
\rput(27.2,4.3){\sm{$\nod_{\t;\bu1}^+$}}\rput(27,1.8){\sm{$\nod_{\t;\bu1}^-$}}
\rput(31.2,4.3){\sm{$\nod_{\t;\bu2}^+$}}\rput(31,1.8){\sm{$\nod_{\t;\bu2}^-$}}
\rput(21.1,4){\sm{$S^1_1$}}\rput(28.9,3){\sm{$S^1_2$}}
\rput(34,3){$\P^1_{\t;\bu}$}
\psline{<->}(35,7)(35,-1)\rput(35.8,3){$\si_{\t}'$}
\psline{->}(18,-8)(18,6)\psline{->}(18,-8)(35,-8)
\rput(33,-9){\sm{$\De_{\R;1}'^+\!\subset\!\C$}}
\rput{90}(17,3.5){\sm{$\De_{\R;2}'^*\!\subset\!\C^4$}}\rput(25,-6){$\De_{\R}'^*$}
\psline[linewidth=.03]{->}(28,-2.5)(28,-7.5)
\pscircle*(18,-8){.2}\pscircle*(22,-8){.2}\pscircle*(18,-4){.2}
\rput(18.8,-4){$\t_0'$}\rput(22.1,-7.2){$\t_0$}
\rput(16.2,-4){\sm{Prp~\ref{genrk1H3degen_prp}}}
\pscircle(8,-12){2}\pscircle(12,-12){2}\pscircle(12,-8){2}\pscircle(12,-16){2}
\psarc[linewidth=.03](2,-10.27){2}{240}{300} 
\psarc[linewidth=.03,linestyle=dashed](2,-13.73){2}{60}{120}
\psarc[linewidth=.03](8,-8.54){4}{240}{300} 
\psarc[linewidth=.03,linestyle=dashed](8,-15.46){4}{60}{120}
\psarc[linewidth=.03](12,-8.54){4}{240}{300} 
\psarc[linewidth=.03,linestyle=dashed](12,-15.46){4}{60}{120}
\psline(1,-9)(1,-15)\psline(3,-9)(3,-15)
\psarc(8,-9){1}{180}{360}\psarc(8,-15){1}{0}{180}
\psarc(5,-9){2}{0}{180}\psarc(5,-15){2}{180}{360}
\psarc(5,-9){4}{0}{180}\psarc(5,-15){4}{180}{360}
\pscircle*(8,-10){.2}\pscircle*(8,-14){.2}\pscircle*(12,-10){.2}\pscircle*(12,-14){.2}
\pscircle*(10,-12){.2}\rput(10,-14){\sm{$\nod$}}
\psline[linewidth=.03]{->}(10,-13.7)(10,-12.7)
\psarc[linewidth=.04](5,-5.3){1}{210}{330}\psarc[linewidth=.04](5,-6.7){1}{30}{150} 
\psarc[linewidth=.04](5,-17.3){1}{210}{330}\psarc[linewidth=.04](5,-18.7){1}{30}{150} 
\rput(8.2,-10.7){\sm{$\nod_{\bu1}^+$}}\rput(8,-13.2){\sm{$\nod_{\bu1}^-$}}
\rput(12.2,-10.7){\sm{$\nod_{\bu2}^+$}}\rput(12,-13.2){\sm{$\nod_{\bu2}^-$}}
\rput(2.1,-11){\sm{$S^1_1$}}
\rput(7.9,-12){\sm{$S^1_{\!\bu1}$}}\rput(11.9,-12){\sm{$S^1_{\!\bu2}$}}
\rput(6,-9.9){\sm{$\bD^2_{\bu1}$}}\rput(14.1,-10.1){\sm{$\bD^2_{\bu2}$}}
\rput(5.2,-12){$\P^1_{\bu1}$}\rput(14.9,-12){$\P^1_{\bu2}$}
\psline{<->}(0,-8)(0,-16)\rput(-.8,-12){$\si_0'$}
\rput(15,-14){$\cC_0'$}
\psline[linewidth=.03]{->}(14.2,-9.3)(17.5,-8.2)
\end{pspicture}
\caption{Deformations of~$\cC_0'$ in the family $(\cU',\wt\fc')$ over 
$\De'\!\subset\!\C^5$ for the proof of Proposition~\ref{genrk1H3degen_prp}}
\label{genH3degen_fig}
\end{figure}

Let  $(\cU',\wt\fc')$ be a flat family of deformations of~$\cC_0'$ as in~\eref{cUsymmdfn_e}
over the unit ball $\De'\!\subset\!\C^5$ around the origin so~that 
$$\cC_{\t_0}'\equiv \big(\Si_{\t_0}',\si_{\t_0}'\big)
\qquad\hbox{and}\qquad
\cC_{\t_0'}'\equiv \big(\Si_{\t_0'}',\si_{\t_0'}'\big)$$
for some $\t_0,\t_0'\!\in\!\De_{\R}'$ are 
a symmetric surface obtained from~$\cC_0'$ by smoothing the real node~$\nod$
and the symmetric  surface~$\cC_0$ 
in the statement of the proposition obtained from~$\cC_0'$ by smoothing 
the conjugate pairs of nodes~$\nod_{\bu r}^{\pm}$, respectively.
Let $\De_{\R}'^*\!\subset\!\De_{\R}'$ be the subspace parametrizing smooth symmetric surfaces
and \hbox{$\De_{\R}'^+,\De_{\R}'^-\!\subset\!\De_{\R}'^*$} be its two topological components 
distinguished as above~\eref{Rdegses_e}.
Denote~by 
$$\De_{\R;1}'^*,\De_{\R;2}'^*\subset\De_{\R}'$$
the subspaces parametrizing all symmetric surfaces 
obtained from~$\cC_0'$ by smoothing the real node~$\nod$
and the conjugate pairs of nodes~$\nod_{\bu r}^{\pm}$, respectively,  and by 
\hbox{$\De_{\R;1}',\De_{\R;2}'\!\subset\!\De_{\R}'$} their closures.
For each $\t\!\in\!\De_{\R;2}'$, let $\nod_{\t}$ be the real node of~$\Si_{\t}'$.\\

Let  $s_1'^{\R},\ldots,s_N'^{\R}$ be sections of $\cU'^{\wt\fc'}$ over~$\De_{\R}'$ and
$(V',\vph')$ be a rank~1 real bundle pair over $(\cU',\wt\fc')$ so~that
$$s_r'^{\R}(\t_0')=x_r\quad\forall\,r\!\in\![N\!-\!1], \quad
s_N^{\R}(\t_0')\in S^1_{\bu1}, \quad
s_N^{\R}(\t)\neq\nod_{\t}~~\forall\,\t\!\in\!\De_{\R;2}',\quad
\big(V'_{\t_0'},\vph'_{\t_0'}\big)=\big(V_0,\vph_0\big).$$
For each $\t\!\in\!\De_{\R}'$ and $r\!\in\![N]$, 
the orientation~$\fo_{x_r}$ of $V_0^{\vph_0}|_{_{x_r}}$ 
induces an orientation~$\fo_{x_r;\t}'$ of~$V'^{\vph'}$ at $s_r'^{\R}(\t)$ 
via the section~$s_r'^{\R}$ as above Proposition~\ref{gencOpincg_prp}.
For $\t\!\in\!\De_{\R;2}'$, $\fo_{x_N}$ also induces an orientation~$\fo_{\nod_{\t}}'$
of~$V'^{\vph'}_{\nod_{\t}}$ so that $\fo_{x_N;\t}'$ is obtained from~$\fo_{\nod_{\t}}'$
by translation along the positive direction of $S^1_{\bu1}\!\subset\!\bD^2_{\bu1}$. 
Let 
\begin{gather*}
\fo_{\x;\t}'=\big(\fo_{x_r;\t}'\big)_{r\in[N]}, ~~
\fo_{\x;\t}'^c=\big(\fo_{x_r;\t}'\big)_{r\in[N-1]} \qquad\forall\,\t\!\in\!\De_{\R}', \\
\wt\fo_{\x;\t}'=\big(\fo_{\nod_{\t}}',\fo_{\nod_{\t}}',\big(\fo_{x_r;\t}'\big)_{r\in[N-1]}\big)
\qquad\forall\,\t\!\in\!\De_{\R;2}'.
\end{gather*}
We denote by $(V_{\bu1}',\vph_{\bu1}')$, $(V_{\bu2}',\vph_{\bu2}')$, 
and $(V_0'^c,\vph_0'^c)$ the restrictions of $(V_0',\vph_0')$
to $\P^1_{\bu1}$, $\P^1_{\bu2}$, and $\Si_0'^c$, respectively.\\

Let \hbox{$\cD'\!=\!\{D_{\t}'\}$} be a family of real CR-operators on 
$(V'_{\t},\vph'_{\t})$ as in~\eref{DVvphdfn_e} 
so that $D_{\t_0'}'$ is the operator~$D_0$ in the statement of the proposition.
We denote~by 
\BE{gencOevComp2H3_e6}
\fo_{\t_0}'\big(\fo_{\x;\t_0}'\big)\equiv \fo'\big(V'_{\t_0},\vph'_{\t_0};\fo_{\x;\t_0}'\big)\EE
the orientation of~$D_{\t_0}$ continuously extending the orientations
$\fo(V'_{\t},\vph'_{\t};\fo_{\x;\t})$ with $\t\!\in\!\De_{\R}'^*$
as above the statement of Proposition~\ref{genrk1Cdegen_prp};
this is the analogue of the limiting orientation~\eref{rk1Clim_e4} 
for the rank~1 real bundle $(V'_{\t_0},\vph'_{\t_0})$ over the symmetric surface 
$\cC_{\t_0}'$ with two conjugate pairs of~nodes.\\

Let  $D'_{\bu1}$, $D'_{\bu2}$, and $D_0'^c$ be
the real CR-operators on $(V_{\bu1}',\vph_{\bu1}')$, $(V_{\bu2}',\vph_{\bu2}')$, 
and $(V_0'^c,\vph_0'^c)$, respectively, induced by~$D_0'$.
The exact~triple 
\begin{gather*}
0\lra D_0' \lra D_0'^c\!\oplus\!D'_{\bu1}\!\oplus\!D'_{\bu2}\lra
 V_{\nod_{\bu1}^+}'\!\oplus\!V_{\nod_{\bu2}^+}'\!\oplus\!V_{\nod_0}'^{\vph'}\lra 0, \\
(\xi,\xi_1,\xi_2)\lra
\big(\xi(\nod_{\bu1}^+)\!-\!\xi_1(\nod_{\bu1}^+),\xi(\nod_{\bu2}^+)\!-\!\xi_2(\nod_{\bu2}^+),
\xi_2(\nod)\!-\!\xi_1(\nod)\big), 
\end{gather*}
of Fredholm operators then determines an isomorphism
$$\la(D_0')\!\otimes\!\la\big(V'_{\nod_{\bu1}^+}\big)\!\otimes\!\la\big(V'_{\nod_{\bu2}^+}\big)
\!\otimes\!\la\big(V_{\nod_0}'^{\vph'}\big)\approx
\la\big(D_0'^c\big)\!\otimes\!\la\big(D'_{\bu1}\big)
\!\otimes\!\la\big(D'_{\bu2}\big).$$
The orientations $\fo(V_0'^c,\vph_0'^c;\fo_{\x;0}'^c)$  of $D_0'^c$, 
$\fo(V'_{\bu r},\vph'_{\bu r};\fo_{\nod_0}')$ of $D'_{\bu r}$, 
and $\fo_{\nod_0}'$ of $V_{\nod_0}'^{\vph'}$
and the complex orientations of~$V'_{\nod_{\bu r}^+}$
determine an orientation of~$D'_0$ via the above isomorphism  and 
thus an orientation~$\fo_{\cU'}(\fo_{\x})$ of the line bundle~$\la(\cD')$ 
over~$\De'_{\R}$.\\

Let $\De_{\R;1}'^+\!\subset\!\De_{\R;1}'^*$ be the intersection of $\De_{\R;1}'^*$
with the closure of~$\De_{\R}'^+$ in~$\De_{\R}'$.
The restriction of $(\cU',\wt\fc')$ to $\De_{\R;1}'$ is the product of a family smoothing 
$$\P^1_{0;\bu}\equiv \P^1_{\bu1}\!\cup\!\P^1_{\bu2}$$
into an irreducible component $\P^1_{\t;\bu}\!\subset\!\Si_{\t}'$ with~$\Si_0'^c$.
For $\t\!\in\!\De_{\R;1}'$,
we denote by $\Si_{\t}'^c\!\subset\!\Si_{\t}'$ the union of the irreducible components
other than~$\P^1_{\t;\bu}$ and by $\nod_{\t;\bu1}^{\pm},\nod_{\t;\bu2}^{\pm}$ the nodes
corresponding to $\nod_{\bu1}^{\pm},\nod_{\bu2}^{\pm}$, respectively.
Let $(V'_{\t;\bu},\vph'_{\t;\bu})$ and $(V_{\t}'^c,\vph_{\t}'^c)$
be the restrictions of $(V'_{\t},\vph'_{\t})$ to $\P^1_{\t;\bu}$ and $\Si_{\t}'^c$,
respectively, and $D'_{\t;\bu}$ and $D_{\t}'^c$ be the real CR-operators on 
$(V'_{\t;\bu},\vph'_{\t;\bu})$ and $(V_{\t}'^c,\vph_{\t}'^c)$, respectively, 
induced by~$D'_{\t}$.
The exact triple~\eref{Rndses_e} of Fredholm operators with $D_0$ replaced 
by $D_{0;\bu}'$ and the exact triple 
\begin{gather*}
0\lra D_{\t_0}' \lra  D_{\t_0}'^c\!\oplus\!D_{\t_0;\bu}' \lra
  V'_{\nod_{\t_0;\bu1}^+}\!\oplus\!V'_{\nod_{\t_0;\bu2}^+}\lra 0, \\
(\xi,\xi_{\bu})\lra
\big(\xi(\nod_{\t_0;\bu1}^+)\!-\!\xi_{\bu}(\nod_{\t_0;\bu1}^+),
\xi(\nod_{\t_0;\bu2}^+)\!-\!\xi_{\bu}(\nod_{\t_0;\bu2}^+)\big), 
\end{gather*}
of Fredholm operators induce isomorphisms
\BE{gencOevComp2H3_e9}\begin{split}
\la\big(D_{0;\bu}'\big)\!\otimes\!\la\big(V_{\nod_0}'^{\vph'}\big)
&\approx\la\big(D'_{\bu1}\big)\!\otimes\!\la\big(D'_{\bu2}\big),\\
\la(D_{\t_0}')\!\otimes\!\la\big(V'_{\nod_{{\t_0};\bu1}^+}\big)
\!\otimes\!\la\big(V'_{\nod_{\t_0;\bu2}^+}\big)
&\approx\la\big(D_{\t_0}'^c\big)\!\otimes\!\la\big(D_{\t_0;\bu}'\big).
\end{split}\EE
The orientations $\fo(V'_{\bu r},\vph'_{\bu r};\fo_{\nod_0}')$ 
of $D'_{\bu r}$ and $\fo_{\nod_0}'$ of $V_{\nod_0}'^{\vph'}$ determine 
an orientation of $D_{0;\bu}'$ via the first isomorphism in~\eref{gencOevComp2H3_e9} 
and thus an orientation $\fo_{\bu;\t}'$ on 
$D'_{\t;\bu}$ for every $\t\!\in\!\De_{\R;1}'$;
the induced orientation of~$D_{0;\bu}'$ is the split orientation in the terminology 
of Proposition~\ref{cOlim_prp}.
By this proposition and the first assumption in~\eref{gencOevComp_e4}, 
the orientation $\fo_{\bu;\t}'$ with $\t\!\in\!\De_{\R;1}'^+$
agrees with the intrinsic orientation 
$$\fo_{\bu;\t}\equiv \fo\big(V'_{\t;\bu},\vph'_{\t;\bu};\fo_{x_N;\t}'\big)$$ 
if and only~if
\BE{gencOevComp2H3_e11}\big(\lr{w_1(V^{\vph}),[S^1_{\bu1}]_{\Z_2}}\!+\!1\big)
\blr{w_1(V^{\vph}),[S^1_{\bu2}]_{\Z_2}}=0\in \Z_2.\EE

\vspace{.15in}

We now assume that $\t_0\!\in\!\De_{\R;1}'^+$.
The conclusion of the previous paragraph then implies that 
the second isomorphism in~\eref{gencOevComp2H3_e9} respects 
the restriction~$\fo_{\t_0}$ of the orientation~$\fo_{\cU'}(\fo_{\x})$ to $\la(D_{\t_0}')$,
the orientations $\fo(V_{\t_0}'^c,\vph_{\t_0}'^c;\fo_{\x;\t_0}'^c)$ 
of $D_{\t_0}'^c$ and $\fo_{\bu;\t_0}$ of $D'_{\t_0;\bu}$,
and the complex orientations of $V'_{\nod_{{\t_0};\bu r}^+}$
if and only if~\eref{gencOevComp2H3_e11} holds.
By the assumption that $r_{\bu}\!=\!N$, 
the decomposition $\Si_{\t_0}'^c\!\sqcup\!\P^1_{\t_0;\bu}$ respects 
the orderings of the topological components of the fixed loci, 
i.e.~$S^1_r$ with $r\!\in\![N\!-\!1]$ and $S^1_{\t_0;r_{\bu}}$.
Combining these two statements with two applications of Proposition~\ref{genrk1Cdegen_prp}, 
we conclude that the orientation~$\fo_{\t_0}$ of~$D_{\t_0}'$
is the analogue of the intrinsic orientation
$$\fo_0\big(\fo_{\x;\t_0}'\big)\equiv \fo\big(V'_{\t_0},\vph'_{\t_0};\fo_{\x;\t_0}'\big)$$
in~\eref{rk1Csplit_e4} 
for the rank~1 real bundle $(V'_{\t_0},\vph'_{\t_0})$ over the symmetric surface 
$\cC_{\t_0}'$ with two conjugate pairs of~nodes 
if and only if~\eref{gencOevComp2H3_e11} holds.
Along with two applications of Proposition~\ref{genrk1Cdegen_prp}, this implies that 
the orientation $\fo_{\t_0}$ of~$D_{\t_0}'$ agrees with~\eref{gencOevComp2H3_e6}
if and only if~\eref{gencOevComp2H3_e11} holds.
Thus, the restriction of~$\fo_{\cU'}(\fo_{\x})$ to $\la(D_{\t}')$ with $\t\!\in\!\De_{\R}'^+$
is the intrinsic orientation $\fo(V'_{\t},\vph'_{\t'};\fo_{\x;\t}')$
defined above Proposition~\ref{gencOpincg_prp} if and only if~\eref{gencOevComp2H3_e11} holds.
This in turn implies that the restriction of~$\fo_{\cU'}(\fo_{\x})$ to $\la(D_{\t_0'}')$
is the limiting orientation~$\fo_0^+(\fo_{\x})$ in~\eref{cOlimor_e4}
if and only if~\eref{gencOevComp2H3_e11} holds.\\

The restriction of $(\cU',\wt\fc')$ to $\De_{\R;2}'$ is a family deforming~$(\Si_0',\si_0')$
to the symmetric surface in the statement of the proposition and above~\eref{wtcCH3deg_e}.
For each $\t\!\in\!\De_{\R;2}'$, we denote by $(\wt\Si_{\t}',\wt\si_{\t}')$ 
the symmetric surface obtained from~$(\Si_{\t}',\si_{\t}')$ by replacing 
the real node~$\nod_{\t}$ with real marked points $\nod_{\t;1}\!\in\!S^1_{\bu1}$
and $\nod_{\t;2}\!\in\!S^1_{\bu2}$,
by $(\wt{V}_{\t}',\wt\vph_{\t}')$ the pullback of $(V_{\t}',\vph_{\t}')$ to
 $(\wt\Si_{\t}',\wt\si_{\t}')$,
and by~$\wt{D}_{\t}'$ the real CR-operator on $(\wt{V}_{\t}',\wt\vph_t')$
induced by~$D_{\t}'$.
The exact~triple 
\begin{gather*}
0\lra \wt{D}_0' \lra 
D_0'^c\!\oplus\!D'_{\bu1}\!\oplus\!D'_{\bu2}\lra
 V'_{\nod_{\bu1}^+}\!\oplus\!V'_{\nod_{\bu2}^+}\lra 0, \\
(\xi,\xi_1,\xi_2)\lra
\big(\xi(\nod_{\bu1}^+)\!-\!\xi_1(\nod_{\bu1}^+),\xi(\nod_{\bu2}^+)\!-\!\xi_2(\nod_{\bu2}^+)\big), 
\end{gather*}
of Fredholm operators then determines an isomorphism
$$\la(\wt{D}_0')\!\otimes\!\la\big(V'_{\nod_{\bu1}^+}\big)
\!\otimes\!\la\big(V'_{\nod_{\bu2}^+}\big)
\approx\la\big(D_0'^c\big)\!\otimes\!\la\big(D'_{\bu1}\big)
\!\otimes\!\la\big(D'_{\bu2}\big).$$
The orientations $\fo(V_0'^c,\vph_0'^c;\fo_{\x;0}'^c)$ of $D_0'^c$ and
$\fo(V'_{\bu r},\wt\vph'_{\bu r};\fo_{\nod;0}')$  of $D'_{\bu r}$
and the complex orientations of~$V'_{\nod_{\bu r}^+}$
determine an orientation of~$\wt{D}'_0$ via the above isomorphism and 
thus an orientation $\wt\fo_{\t}$ of~$\wt{D}'_{\t}$ for every $\t\!\in\!\De_{\R;2}'$.\\

By the assumption that $r_{\bu}\!=\!N$, 
the decomposition $\Si_0'^c\!\sqcup\!\P^1_{\bu1}\!\sqcup\!\P^1_{\bu2}$ respects 
the orderings of the topological components of the fixed loci of~$(\wt\Si_0,\wt\si_0)$, 
i.e.~$S^1_r$ with $r\!\in\![N\!-\!1]$, $S^1_{\bu1}$, and~$S^1_{\bu2}$. 
Along with Proposition~\ref{gencOpincg_prp}\ref{gencOdisjun_it} applied twice, 
this implies~that the orientation~$\wt\fo_0$ of $\wt{D}'_0$ above
is the analogue of the intrinsic orientation
$$\fo_0\big(\wt\fo_{\x;0}'\big)\equiv \fo\big(V'_0,\vph'_0;\wt\fo_{\x;0}'\big)$$
in~\eref{rk1Csplit_e4} 
for the rank~1 real bundle $(\wt{V}'_0,\wt\vph'_0)$ over the symmetric surface 
$(\wt\Si_0',\wt\si_0')$ with two conjugate pairs of~nodes.
Combining this with two applications of Proposition~\ref{genrk1Cdegen_prp}, we conclude that 
$\wt\fo_{\t_0'}$ is the orientation 
$$\wt\fo_0\big(\fo_{\x;\t_0'}\big)\equiv
\fo\big(\wt{V}_{\t_0}',\wt\vph_{\t_0}';\wt\fo_{\x;\t_0'}'\big)
=\fo\big(\wt{V}_0,\wt\vph_0;\wt\fo_{\x}\big)$$
of $\wt{D}'_{\t_0}\!=\!\wt{D}_0$ in~\eref{rk1H3split_e4a}.
Thus, the restriction of~$\fo_{\cU'}(\fo_{\x})$ to $\la(D_{\t_0'}')$
is the intrinsic orientation~$\fo_0(\fo_{\x})$ of $D'_{\t_0}\!=\!D_0$
in~\eref{rk1H3split_e4}.
Along with the conclusion regarding the limiting orientation~$\fo_0^+(\fo_{\x})$ above,
this implies that~$\fo_0(\fo_{\x})$ and the limiting orientation~$\fo_0^+(\fo_{\x})$ of~$D_0$ 
are the same if and only if~\eref{gencOevComp2H3_e11} holds
and establishes the $r_{\bu}\!=\!N$ case of the proposition.
\end{proof}

\subsection{Orientations from $\OSpin$-structures}
\label{SpinOrient_subs}

Suppose $(\Si,\si)$ is a smooth decorated symmetric surface, 
$(V,\vph)$ is a real bundle pair over~$(\Si,\si)$,  and 
\hbox{$\os\!\in\!\OSpin_{\Si}(V^{\vph})$} is a relative $\OSpin$-structure on
the real vector bundle $V^{\vph}$ over $\Si^{\si}\!\subset\!\Si$.
We show below that~$\os$ determines an orientation~$\fo_{\os}(V,\vph)$ 
of every real CR-operator~$D$ on~$(V,\vph)$.\\

Let $\cC_0$, $(\P^1_r,\tau_r)$, $S^1_r$, $\bD^2_{r\pm}$, $\Si_0^+\!\subset\!\Si_0^{\C}$,
$\nod_r^{\pm}$, $(\cU,\wt\fc)$, $\t_1\!\in\!\De_{\R}$, $(\wt{V},\wt\vph)$,
$(\wt{V}_{0r},\wt\vph_{0r})$,
\hbox{$\cD\!\equiv\!\{D_{\t}\}$}, $D_{0r}$, and $D_0^+$ be as below~\eref{fobfxfn_e}. 
By the compatibility condition in Definition~\ref{RelPinSpin_dfn3}, a tuple 
\BE{wtos0dfn_e}\wt\os_0\equiv \big((\os_{0;r})_{S^1_r\in\pi_0(\Si^{\si})},\os_{0;\C}\big)
\in \prod_{S^1_r\in\pi_0(\Si^{\si})}\!\!\!\!\!\!\!\!
\OSp_{\P^1_r}\big(\wt{V}^{\wt\vph}|_{S^1_r}\big)
\times\!\OSp_{\Si_0^{\C}}\big(\wt{V}^{\wt\vph}|_{\eset}\big)\EE
of relative $\OSpin$-structures on the restrictions of $\wt{V}_0^{\wt\vph_0}$ to
the $\si_0$-fixed loci contained in $\P^1_r$ and $\Si_0^{\C}$ determines 
a relative $\OSpin$-structure~$\os_0$ on $\wt{V}_0^{\wt\vph_0}$ and thus 
a relative $\OSpin$-structure~$\os_{\t}$ on the restriction of~$\wt{V}^{\wt\vph}$
to $\Si_{\t}^{\si_{\t}}\!\subset\!\Si_{\t}$ for each $\t\!\in\!\De_{\R}$.
Fix a tuple $\wt\os_0$ as in~\eref{wtos0dfn_e} so that the induced relative $\OSpin$-structure
$\os_{\t_1}$ is~$\os$.\\

The exact triple~\eref{geCndses_e} of Fredholm operators again induces 
an isomorphism~\eref{geCndses_e2} of the determinants of the associated real CR-operators.
For each $S^1_r\!\in\!\pi_0(\Si^{\si})$, let 
$\fo_{\os_{0;r}}(\wt{V}_{0r},\wt\vph_{0r})$ be the orientation
of $D_{0r}$ as above Proposition~\ref{tauspinorient_prp}.
We denote by $\fo_{\os_{0;\C}}(\wt{V}|_{\Si_0^+})$
the complex orientation of $D_0^+$ if 
$$\os_{0;\C}(\Si_0^+)\equiv \blr{w_2(\os_{0;\C}),[\Si_0^+]_{\Z_2}}\in\Z_2$$ 
vanishes and the opposite orientation otherwise.
Along with the chosen order on~$\pi_0(\Si^{\si})$ and 
the complex orientations of $\wt{V}_{\nod_r^+}$,
these orientations determine an orientation $\fo_{\os}(\wt{V}_0,\wt\vph_0)$ 
of $D_0$  via the isomorphism~\eref{geCndses_e2}
 and thus an orientation~$\fo_{\os}(V,\vph)$ of $D\!=\!D_{\t_1}$.\\

If the real bundle pair $(V,\vph)$ is of rank~1 and $\fo$ is an orientation on~$V^{\vph}$, 
we let
\BE{osLBdfn_e2}\fo_0(V,\vph;\fo)\equiv \fo_{\io_{\Si}(\os_0(V^{\vph},\fo))}(V,\vph)\EE
denote the orientation of~$D$ determined by the image
$$\io_{\Si}\big(\os_0(V^{\vph},\fo)\!\big)\in\OSp_{\Si}\big(V^{\vph}\big)$$
of  the canonical $\OSpin$-structure $\os_0(V^{\vph},\fo)$ on $(V^{\vph},\fo)$
under the first map in~\eref{vsSpinPin_e0} with $X\!=\!\Si$.

\begin{prp}\label{genspinorient_prp}
Suppose $(V,\vph)$ is a real bundle pair over a smooth decorated symmetric surface  
$(\Si,\si)$, $\os\!\in\!\OSpin_{\Si}(V^{\vph})$ is a relative $\OSpin$-structure on
the real vector bundle $V^{\vph}$ over $\Si^{\si}\!\subset\!\Si$,
and $D$ is a real CR-operator on~$(V,\vph)$.
\begin{enumerate}[label=(\arabic*),leftmargin=*]

\item\label{genspindfn_it} 
The orientation $\fo_{\os}(V,\vph)$ of~$D$ 
constructed above does not depend the choice of an admissible tuple~\eref{wtos0dfn_e}.

\item\label{genrk1spindfn_it} If $(V,\vph)$ is of rank~1
and $\fo$ is an orientation on~$V^{\vph}$,
the orientation~\eref{osLBdfn_e2} of~$D$ is the same 
as the intrinsic orientation $\fo(V,\vph;\fo_{\x})$ of  Proposition~\ref{gencOpincg_prp}
for a tuple~$\fo_{\x}$ of orientations as in~\eref{fobfxfn_e} induced by 
the orientation~$\fo$.

\item\label{tauspindfn_it}  If $(\Si,\si)\!=\!(S^2,\tau)$,
the orientation $\fo_{\os}(V,\vph)$ of~$D$ is the same as the orientation constructed above 
Proposition~\ref{tauspinorient_prp}. 

\end{enumerate}  
The orientations $\fo_{\os}(V,\vph)$ of real CR-operators
on real bundle pairs $(V,\vph)$ over smooth decorated symmetric surfaces  
$(\Si,\si)$ constructed above satisfy all applicable CROrient properties 
of Sections~\ref{OrientPrp_subs1} and~\ref{OrientPrp_subs2}.
\end{prp}

\begin{proof}[{\bf{\emph{Proof of 
Proposition~\ref{genspinorient_prp}\ref{genspindfn_it}-\ref{tauspindfn_it}}}}]
\ref{genspindfn_it} Suppose 
$$\wt\os_0'\equiv \big((\os_{0;r}')_{S^1_r\in\pi_0(\Si^{\si})},\os_{0;\C}'\big)
\in \prod_{S^1_r\in\pi_0(\Si^{\si})}\!\!\!\!\!\!\!\!
\OSp_{\P^1_r}\big(\wt{V}^{\wt\vph}|_{S^1_r}\big)
\times\!\OSp_{\Si_0^{\C}}\big(\wt{V}^{\wt\vph}|_{\eset}\big)$$
is another tuple of relative $\OSpin$-structures inducing the relative $\OSpin$-structure~$\os$
as below~\eref{wtos0dfn_e}.
For each $S^1_r\!\in\!\pi_0(\Si^{\si})$, let $\os_{0;r}(\bD^2_{r+})$ and 
$\os_{0;r}'(\bD^2_{r+})$ be the collections of trivializations of $\wt{V}^{\wt\vph}$
induced by the relative $\OSpin$-structures~$\os_{0;r}$ and~$\os_{0;r}'$,
respectively, and by the inclusion $\bD^2_{r+}\!\subset\!\P^1_r$.
Define 
$$\ep_r=\begin{cases}0\!\in\!\Z_2,&\hbox{if}~
\os_{0;r}(\bD^2_{r+})\!=\!\os_{0;r}'(\bD^2_{r+});\\
1\!\in\!\Z_2,&\hbox{if}~
\os_{0;r}(\bD^2_{r+})\!\neq\!\os_{0;r}'(\bD^2_{r+});
\end{cases} \quad
\ep_{\C}=\os_{0;\C}(\Si_0^+)\!-\!\os_{0;\C}'(\Si_0^+)\!\in\!\Z_2.$$

\vspace{.15in}

Since $\wt\os_0$ and $\wt\os_0'$ induce the same relative $\OSpin$-structure
on the restriction of~$\wt{V}^{\wt\vph}$ to $\Si_{\t}^{\si_{\t}}\!\subset\!\Si_{\t}$
for each $\t\!\in\!\De_{\R}^*$,
the equivalence classes~$\os_0(\Si_0^b)$ and~$\os_0'(\Si_0^b)$ of trivializations of 
the restriction of~$\wt{V}^{\wt\vph}$ to $\Si_0^{\si_0}\!\subset\!\Si_0$
determined by~$\os_0$ and~$\os_0'$, respectively, and the normalization~map
$$u\!:\wt\Si_0^b\!\equiv\!\bigsqcup_{S^1_r\in\pi_0(\Si^{\si})}\!\!\!\!\!\!\!\!\bD^2_+
\sqcup\!\Si_0^+\lra\Si_0$$
are the same. Thus,
$$\sum_{S^1_r\in\pi_0(\Si^{\si})}\!\!\!\!\!\!\!\!\ep_r+\ep_{\C}=0\in\Z_2.$$
Along with Proposition~\ref{SpinFOOO_prp}, 
this implies that the number of factors on the right-hand side of
the isomorphism~\eref{geCndses_e2} for which the orientations determined
by~$\wt\os_0$ and~$\wt\os_0'$ differ is even.
Thus, the orientations of~$D_0$ and $D\!=\!D_{\t_1}$ determined by~$\wt\os_0$ and~$\wt\os_0'$ 
via the isomorphism~\eref{geCndses_e2} are the~same.\\

\ref{genrk1spindfn_it}
For the purposes of constructing the orientation~\eref{osLBdfn_e2}, 
we can choose the tuple in~\eref{wtos0dfn_e} so~that 
$$\os_{0;r}= \io_{\P^1_r}\big(\os_0\big((V^{\vph},\fo)|_{S^1_r}\big)\!\big)
~~\forall\,S^1_r\!\in\!\pi_0(\Si^{\si}) 
\qquad\hbox{and}\qquad \os_{0;\C}(\Si_0^+)=0.$$
The definition of the orientation~\eref{osLBdfn_e2} above Proposition~\ref{genspinorient_prp}
then becomes identical to that 
of the intrinsic orientation of Proposition~\ref{gencOpincg_prp} 
for a tuple~$\fo_{\x}$ of orientations as in~\eref{fobfxfn_e} 
induced by the orientation~$\fo$ of~$V^{\vph}$.\\

\ref{tauspindfn_it} In this case, $|\pi_0(\Si^{\si})|\!=\!1$. 
There are canonical homotopy classes of identifications of
$(\P^1_1,\tau_1,\bD^2_{1+})$ with $(\Si,\si,\Si^b)\!\equiv\!(S^2,\tau,\bD^2_+)$
and of the restriction of~$\wt{V}^{\wt\vph}$ to $\Si_0^{\si_0}$ with the real
vector bundle $V^{\vph}$ over $S^1\!\subset\!\Si$.
We can choose the tuple in~\eref{wtos0dfn_e} so~that $\os_{0;1}$ is identified with $\os$
via these homotopy classes and  
$$\os_{0;\C}(\Si_0^+)\equiv\os_{0;\C}(\P^1_+)=0\in\Z_2.$$
The definition of the orientation $\fo_{\os}(V,\vph)$ above Proposition~\ref{genspinorient_prp}
then becomes identical to that above Proposition~\ref{tauspinorient_prp}. 
\end{proof}

\begin{proof}[{\bf{\emph{Proof of CROrient~\ref{CROos_prop},
\ref{CROSpinPinStr_prop}\ref{CROSpinStr_it}, and 
\ref{CRONormal_prop}\ref{CROnormSpin_it} properties}}}]
We continue with the notation above the statement of Proposition~\ref{genspinorient_prp}.
With the notation as in the statement of the CROrient~\ref{CROos_prop}\ref{osflip_it} property
and in the proof of Proposition~\ref{gencOpincg_prp}\ref{gencOpincg_it1a},
the first two identities in~\eref{gencOpincg_e2} still apply;
the third becomes
\BE{genos_e2}
\ind_{\C}D_0^+=n\big(1\!-\!g(\Si_{*0}^+)\!\big)\!+\!\deg\wt{V}|_{\Si_{*0}^+}.\EE
The change in the choice of the half-surface $\Si_*^b$ of $(\Si_*,\si)$ acts by 
the complex conjugation on the complex orientations of $D_0^+|_{\Si_{*0}^+}$
and each $\wt{V}_{\nod_r^+}$ with $S^1_r\!\in\!\pi_0(\Si_*^{\si})$.
By the CROrient~\ref{CROos_prop}\ref{osflip_it} property for $(\P^1,\tau)$ provided by
Proposition~\ref{tauspinorient_prp}, 
this change preserves the orientation 
$\fo_{\os_{0;r}}(\wt{V}_{0r},\wt\vph_{0r})$ of $D_{0r}$
with $S^1_r\!\in\!\pi_0(\Si_*^{\si})$ if and only if $\vp_{\os_{0;r}}(\P^1_r)$ vanishes
in~$\Z_2$.
Thus, the orientations $\fo_{\os}(\wt{V}_0,\wt\vph_0)$ and $\fo_{\os}(V,\vph)$ 
do not depend on the choice of half-surface $\Si_*^b$  of $(\Si_*,\si)$ if and only~if
\begin{equation*}\begin{split}
&\ind_{\C}D_0^+|_{\Si_{*0}^+}-n\big|\pi_0(\Si_*^{\si})\big|
+\blr{w_2(\os_{0;\C}),[\Si_{0*}^{\C}]_{\Z_2}}
+\sum_{S^1_r\in\pi_0(\Si_*^{\si})}\!\!\!\!\!\!\!\vp_{\os_{0;r}}(\P^1_r)\\
&\hspace{.5in}
=\ind_{\C}D_0^+|_{\Si_{*0}^+}-n\big|\pi_0(\Si_*^{\si})\big|+
\sum_{S^1_r\in\pi_0(\Si_*^{\si})}\!\!\!\!\!\!\frac{\deg\wt{V}_{0r}}{2}
+\blr{w_2(\os),[\Si_*]_{\Z_2}}
\end{split}\end{equation*}
vanishes in~$\Z_2$.
Combining this with~\eref{genos_e2} and the first two identities in~\eref{gencOpincg_e2}, we obtain
the CROrient~\ref{CROos_prop}\ref{osflip_it} property.\\

Since the parity of the index of $D_{0r}$ is the same as the parity of~$n$,
the CROrient~\ref{CROos_prop}\ref{osinter_it} property follows from the Direct Sum property for 
the determinants of Fredholm operators; see \cite[Section~2]{detLB}.
The CROrient~\ref{CROos_prop}\ref{dltCROos_it} property holds by definition.
Propositions~\ref{genspinorient_prp}\ref{tauspindfn_it} and~\ref{tauspinorient_prp} imply
the CROrient~\ref{CRONormal_prop}\ref{CROnormSpin_it} property.
The second statement of the CROrient~\ref{CROSpinPinStr_prop}\ref{CROSpinStr_it} property
follows immediately from its $(S^2,\tau)$ case provided by Proposition~\ref{tauspinorient_prp}.\\

Let $\eta\!\in\!H^2(\Si,\Si^{\si};\Z_2)$.
Since the composition of the retraction isomorphism and the restriction homomorphism
$$ H^2\big(\Si_0,\Si_0^{\si_0};\Z_2\big)\approx
  H^2\big(\cU,\cU^{\wt\fc};\Z_2\big) \lra H^2\big(\Si,\Si^{\si};\Z_2\big)$$
is surjective, some $\eta_0\!\in\!H^2(\Si_0,\Si_0^{\si_0};\Z_2)$ is mapped to~$\eta$ 
under this composition.
The~tuple 
$$\eta_0\!\cdot\!\wt\os_0
\equiv\big((\eta_0|_{(\P^1_r,S^1_r)}\!\cdot\!\os_{0;r})_{S^1_r\in\pi_0(\Si_0^{\si_0})},
\eta_0|_{\Si_0^{\C}}\!\cdot\!\os_{0;\C}\big)$$
then induces the relative $\OSpin$-structure $\eta\!\cdot\!\os$ as below~\eref{wtos0dfn_e}.
By the first statement of the CROrient~\ref{CROSpinPinStr_prop}\ref{CROSpinStr_it} property
for $(S^2,\tau)$
provided by Proposition~\ref{tauspinorient_prp} and the definition of 
the orientation $\fo_{\os_{0;\C}}(\wt{V}|_{\Si_0^+})$
of $D_0^+$,
the orientations $\fo_{\eta\cdot\os}(\wt{V}_0,\wt\vph_0)$ and $\fo_{\eta\cdot\os}(V,\vph)$ 
are the same as $\fo_{\os}(\wt{V}_0,\wt\vph_0)$ and $\fo_{\os}(V,\vph)$,
respectively, if and only~if
$$\sum_{S^1_r\in\pi_0(\Si_*^{\si})}\!\!\!\!\!\!
\blr{\eta_0|_{(\P^1_r,S^1_r)},[\bD^2_{r+}]_{\Z_2}}
+\blr{\eta_0|_{\Si_0^{\C}},[\Si_0^+]_{\Z_2}}=\blr{\eta,[\Si^b]_{\Z_2}}$$
vanishes in~$\Z_2$.
This establishes the first statement of 
the CROrient~\ref{CROSpinPinStr_prop}\ref{CROSpinStr_it} property.
\end{proof}

\begin{proof}[{\bf{\emph{Proof of CROrient~\ref{CRODisjUn_prop}\ref{osDisjUn_it} property}}}]
Suppose $(\Si_1,\si_1)$, $(\Si_2,\si_2)$, $(V_1,\vph_1)$, $(V_2,\vph_2)$, $\os_1$,
$\os_2$, 
\begin{gather*}
D_1\equiv  D_{(V_1,\vph_1)}, \qquad  D_2\equiv  D_{(V_2,\vph_2)}, \\
(V,\vph)\equiv (V_1,\vph_1)\!\sqcup\!(V_2,\vph_2), \quad
\os\!\equiv\!\os_1\!\sqcup\!\os_2, \quad\hbox{and}\quad
D\equiv D_1\!\sqcup\!D_2
\end{gather*}
are as in the statement of this property on page~\pageref{CRODisjUn_prop}.
Let $(\cU_1,\wt\fc_1)$, $(\cU_2,\wt\fc_2)$, 
$(\wt{V}_1,\wt\vph_1)$, $(\wt{V}_2,\wt\vph_2)$, 
$$\wt\os_0^{(1)}\equiv \big((\os_{0;r}^{(1)})_{S^1_r\in\pi_0(\Si_1^{\si_1})},\os_{0;\C}^{(1)}\big),
\qquad
\wt\os_0^{(2)}\equiv \big((\os_{0;r}^{(2)})_{S^1_r\in\pi_0(\Si_2^{\si_2})},\os_{0;\C}^{(2)}\big),$$
$D_{1;\t}$, and $D_{2;\t}$ 
be as in the construction of the orientations 
$\fo_{\os_1}\!(V_1,\vph_1)$ and $\fo_{\os_2}\!(V_2,\vph_2)$ above Proposition~\ref{genspinorient_prp}.
The orientation $\fo_{\os}(V,\vph)$ 
of~$D$ is obtained via this construction applied with
\begin{gather*}
\big(\cU,\wt\fc\big)=
\big(\cU_1,\wt\fc_1\big)\!\!\sqcup\!\!\big(\cU_2,\wt\fc_2\big),
\quad
\big(\wt{V},\wt\vph\big)=\big(\wt{V}_1,\wt\vph_1\big)
\!\!\sqcup\!\!\big(\wt{V}_2,\wt\vph_2\big),
\quad
D_{\t}=D_{1;\t}\!\sqcup\!D_{2;\t}, \\
\wt\os_0=\big((\os_{0;r}^{(1)})_{S^1_r\in\pi_0(\Si_1^{\si_1})},
(\os_{0;r}^{(2)})_{S^1_r\in\pi_0(\Si_2^{\si_2})},
\os_{0;\C}^{(1)}\!\sqcup\!\os_{0;\C}^{(2)}\big).
\end{gather*}
In particular,
$$\os_{0;\C}(\Si_0^+)=
\os_{0;\C}^{(1)}\big(\Si_{1;0}^+\big)\!+\!
\os_{0;\C}^{(2)}\big(\Si_{2;0}^+\big)\in\Z_2\,.$$

\vspace{.15in}

The proof now proceeds via the diagram of Figure~\ref{rk1CRODisjUn_fig},
as in the proof of Proposition~\ref{gencOpincg_prp}\ref{gencOdisjun_it},  
with the orientations~\eref{LHSorient_e} replaced by 
$$\fo_{\os_1}\!\big(\wt{V}_{1;0},\wt\vph_{1;0}\big), 
\quad
\fo_{\os}\big(\wt{V}_0,\wt\vph_0\big),
\quad\hbox{and}\quad
\fo_{\os_2}\!\big(\wt{V}_{2;0},\wt\vph_{2;0}\big),$$
respectively, the orientations~\eref{midorient_e} replaced by 
$$\fo_{\os_{0;r}}(\wt{V}_{0r},\wt\vph_{0r})=
\begin{cases}\fo_{\os_{0;r}^{(1)}}(\wt{V}_{0r},\wt\vph_{0r}),
&\hbox{if}~S^1_r\!\in\!\pi_0(\Si_1^{\si_1});\\
\fo_{\os_{0;r}^{(2)}}(\wt{V}_{0r},\wt\vph_{0r}),
&\hbox{if}~S^1_r\!\in\!\pi_0(\Si_2^{\si_2});
\end{cases}$$
and the complex orientations on the determinants of $D_{1;0}^+$,
$D_0^+$, and $D_{2;0}^+$  replaced~by  
$$\fo_{\os_{0;\C}^{(1)}}\big(\wt{V}_1|_{\Si_{1;0}^+}\big),
\qquad
\fo_{\os_{0;\C}}\big(\wt{V}|_{\Si_0^+}\big),
\qquad\hbox{and}\qquad
\fo_{\os_{0;\C}^{(2)}}\big(\wt{V}_2|_{\Si_{2;0}^+}\big),$$
respectively.
\end{proof}

\begin{proof}[{\bf{\emph{Proof of CROrient~\ref{CROSpinPinSES_prop}\ref{CROsesSpin_it} property}}}]
Suppose $(\Si,\si)$, $\ce$, $\os'$, $\os''$,
$$\os\equiv\llrr{\os',\os''}_{\ce_{\R}}, \quad
 D'\equiv D_{(V',\vph')}, \quad\hbox{and}\quad D''\equiv D_{(V'',\vph'')}$$
are as in the statement of this property on page~\pageref{CROSpinPinSES_prop}.
Let $(\cU,\wt\fc)$,  $(\wt{V}',\wt\vph')$, $(\wt{V}'',\wt\vph'')$, 
$$\wt\os_0'\equiv \big((\os_{0;r}')_{S^1_r\in\pi_0(\Si^{\si})},\os_{0;\C}'\big),
\qquad
\wt\os_0''\equiv \big((\os_{0;r}'')_{S^1_r\in\pi_0(\Si^{\si})},\os_{0;\C}''\big),$$
$D_{\t}'$, and $D_{\t}''$ be associated objects as
in the construction above Proposition~\ref{genspinorient_prp}.
The orientation $\fo_{\os}(V,\vph)$ of
$D'\!\oplus\!D''$ is obtained via this construction applied with
\begin{gather*}
\big(\wt{V},\wt\vph\big)=\big(\wt{V}',\wt\vph'\big)
\!\oplus\!\big(\wt{V}'',\wt\vph''\big),
\qquad
D_{\t}=D_{\t}'\!\oplus\!D_{\t}'', \\
\wt\os_0=\big((\os_{0;r}\!\equiv\!\llrr{\os_{0;r}',\os_{0;r}''}_{\oplus})_{S^1_r\in\pi_0(\Si^{\si})},
\os_{0;\C}\!\equiv\!\llrr{\os_{0;\C}',\os_{0;\C}''}_{\oplus}\big).
\end{gather*}
The exact triple~\eref{geCndses_e} of Fredholm operators induces 
the exact triples of Fredholm operators given by the rows in the diagram of 
Figure~\ref{CROSpinPinSES_fig}.
The short exact sequence
$$0\lra \big(\wt{V}_0',\wt\vph'_0\big)\lra 
\big(\wt{V}_0,\wt\vph_0\big)\lra \big(\wt{V}''_0,\wt\vph''_0\big)\lra0$$
of real bundle pairs 
induces the exact triples of Fredholm operators given by the columns in this diagram.\\

\begin{figure}
$$\xymatrix{& 0\ar[d] & 0\ar[d] & 0\ar[d]  \\
0\ar[r] & D_0'\ar[r]\ar[d] & 
\bigoplus\limits_{S^1_r\in\pi_0(\Si^{\si})}\!\!\!\!\!\!\!\!D_{0r}'
\oplus\!D_0'^+ \ar[r]\ar[d]  & 
\bigoplus\limits_{S^1_r\in\pi_0(\Si^{\si})}\!\!\!\!
\wt{V}'_{\nod_r^+}\ar[r]\ar[d] & 0 \\  
0\ar[r] & D_0\ar[r]\ar[d] & 
\bigoplus\limits_{S^1_r\in\pi_0(\Si^{\si})}\!\!\!\!\!\!\!\!D_{0r}
\oplus\!D_0^+ \ar[r]\ar[d]  & 
\bigoplus\limits_{S^1_r\in\pi_0(\Si^{\si})}\!\!\!\!\!\!\wt{V}_{\nod_r^+}\ar[r]\ar[d] & 0 \\
0\ar[r] & D_0''\ar[r]\ar[d] & 
\bigoplus\limits_{S^1_r\in\pi_0(\Si^{\si})}\!\!\!\!\!\!\!\!D_{0r}''
\oplus\!D_0''^+ \ar[r]\ar[d]  & 
\bigoplus\limits_{S^1_r\in\pi_0(\Si^{\si})}\!\!\!\wt{V}''_{\nod_r^+}\ar[r]\ar[d] & 0 \\
& 0 & 0 & 0}$$ 
\caption{Commutative square of exact rows and columns of Fredholm operators 
for the proof of the CROrient~\ref{CROSpinPinSES_prop}\ref{CROsesSpin_it} property}
\label{CROSpinPinSES_fig}
\end{figure}

By definition, the rows in this diagrams respect the orientations
\BE{sesLHSorient_e}
\fo_{\os_0'}\big(\wt{V}_0',\wt\vph_0'\big), \qquad 
\fo_{\os_0}\big(\wt{V}_0,\wt\vph_0\big), \qquad\hbox{and}\qquad 
\fo_{\os_0''}\big(\wt{V}_0'',\wt\vph_0''\big)\EE 
of the operators in the left column, the orientations 
\BE{sesmidorient_e}
\fo_{\os_{0;r}'}\big(\wt{V}_{0r}',\wt\vph_{0r}'\big), \qquad 
\fo_{\os_{0;r}}\big(\wt{V}_{0r},\wt\vph_{0r}\big), \qquad\hbox{and}\qquad 
\fo_{\os_{0;r}''}\big(\wt{V}_{0r}'',\wt\vph_{0r}''\big)\EE
of the operators in the direct sums in the middle column, the orientations
\BE{sesmidorient_e2}
\fo_{\os_{0;\C}'}\big(\wt{V}'|_{\Si_0^+}\big), \qquad 
\fo_{\os_{0;\C}}\big(\wt{V}|_{\Si_0^+}\big), \qquad\hbox{and}\qquad 
\fo_{\os_{0;\C}''}\big(\wt{V}''|_{\Si_0^+}\big)\EE
on the last summands in middle column,
and the complex orientations of the vector spaces in the right column.
The last column respects the complex orientations;
the exact triple formed by the last summands in the middle column respects 
the orientations~\eref{sesmidorient_e2}, since
$$\os_{0;\C}(\Si_0^+)=\os_{0;\C}'(\Si_0^+)\!+\!\os_{0;\C}''(\Si_0^+).$$
The exact triple formed by the $r$-th summands in the middle column respects 
the orientations~\eref{sesmidorient_e}.
Since 
$$\ind\,D_{0r}'=\rk\,V', \qquad \ind\,D_{0r}''=\rk\,V'', \qquad\hbox{and}\qquad
\ind\,D_0'^+\in2\Z,$$
the last two statements and Lemma~\ref{3by3_lmm} imply that the middle column respects
the direct sum orientations if and only if~\eref{CROses_e0b} holds.
Since the (real) dimensions of $\wt{V}'_{\nod_r^+}$ are even,
another application of Lemma~\ref{3by3_lmm} combined with the statements implies~that 
the left column respects the orientations~\eref{sesLHSorient_e}
if and only if~\eref{CROses_e0b} holds.
The claim now follows from the continuity of 
the isomorphisms
$$\la\big(D_{\t}\big)\approx \la\big(D_{\t}'\big)\!\otimes\!\la\big(D_{\t}''\big)$$
with respect to $\t\!\in\!\De_{\R}$.
\end{proof}

\begin{proof}[{\bf{\emph{Proof of CROrient~\ref{CROdegenC_prop}\ref{CdegenSpin_it} property}}}]
Let $\cC_0$ be as in~\eref{cCdegC_e} with $k,l\!=\!0$. 
We denote by $\os_1$ the standard relative $\OSpin$-structure on 
the trivial line bundle $\Si_0^{\si_0}\!\times\!\R$ over $\Si_0^{\si_0}\!\subset\!\Si_0$ and
by $\bp$ the standard real CR-operator on the real bundle pair 
$(\Si_0\!\times\!\C,\si_0\!\times\!\fc)$ over $(\Si_0,\si_0)$.
Let $(L,\phi)$ be a rank~1 real bundle pair over~$(\Si_0,\si_0)$ so that 
the real line bundle~$L^{\phi}$ over~$\Si_0^{\si_0}$  is orientable,
$\os_L$ be a relative $\OSpin$-structure on 
the real line bundle~$L^{\phi}$ over~$\Si_0^{\si_0}\!\subset\!\Si_0$,
and $\bp_L$ be a real CR-operator on~$(L,\phi)$.
For $n\!\in\!\Z^+$, let  
\begin{gather}\notag
(V_n,\vph_n)=(\Si_0\!\times\!\C^n,\si_0\!\times\!\fc), \qquad
\os_{n;L}=\St_{L^{\phi}}^{n-1}(\os_L)=\llrr{\os_1,\os_{n-1;L}}_{\oplus},\\
\label{SpinPin1Spin_e7}
(V_{n;L},\vph_{n;L})=\big(V_{n-1},\vph_{n-1}\big)\!\oplus\!(L,\phi)
=\big(V_1,\vph_1\big)\!\oplus\!(V_{n-1;L},\vph_{n-1;L}).
\end{gather}
Let $D_n$ be the real CR-operator on $(V_{n;L},\vph_{n;L})$ given by the $n$-fold direct 
sum of the operators~$\bp$ on each factor and \hbox{$D_{n;L}\!=\!D_{n-1}\!\oplus\!\bp_L$}.\\

We first note that the CROrient~\ref{CROdegenC_prop}\ref{CdegenSpin_it} holds
for $\os_1$ and~$\os_{1;L}\!\equiv\!\os_L$, i.e.
\BE{CdegenSpin_e3} 
\fo_{\os_1}\big(V_1,\vph_1\big)=\fo_{\os_1}'\big(V_1,\vph_1\big) 
\qquad\hbox{and}\qquad
\fo_{\os_L}(L,\phi)=\fo_{\os_L}'(L,\phi).\EE
The first equality follows immediately from Propositions~\ref{genspinorient_prp}\ref{genrk1spindfn_it} 
and~\ref{genrk1Cdegen_prp}.
So does the second if $\os_L\!=\!\io_{\Si_0}(\os_0(L^{\phi},\fo))$ for some orientation~$\fo$
on~$L^{\phi}$.
In general, 
$$\os_L=\eta_0\!\cdot\!\io_{\Si_0}\big(\os_0(L^{\phi},\fo)\big)$$ 
for some orientation~$\fo$ on~$V^{\vph}$ and some 
$$\eta_0\in H^2\big(\Si_0,\Si_0^{\si_0};\Z_2\big)\approx
H^2\big(\wt\Si_0,\wt\Si_0^{\wt\si_0};\Z_2\big), H^2\big(\cU,\cU^{\wt\fc};\Z_2\big);$$
see the RelSpinPin~\ref{RelSpinPinStr_prop} property on page~\pageref{RelSpinPinStr_prop}.
By the CROrient~\ref{CROSpinPinStr_prop}\ref{CROSpinStr_it} property,
$$\fo_{\os_L}(L,\phi)=\fo_{\io_{\Si_0}(\os_0(L^{\phi},\fo))}(L,\phi)
\qquad \big(\hbox{resp.}~
\fo_{\os_L}'(L,\phi)=\fo_{\io_{\Si_0}(\os_0(L^{\phi},\fo))}'(L,\phi)\big)$$
if and only if $\lr{\eta,\Si_0^b}\!=\!0$.
Along with the second statement in~\eref{CdegenSpin_e3} with~$\os_L$ replaced by 
$\io_{\Si_0}(\os_0(L^{\phi},\fo))$,
this implies the second statement itself.\\

\begin{figure}
$$\xymatrix{& 0\ar[d] & 0\ar[d] & 0\ar[d]  \\
0\ar[r] & \bp\ar[r]\ar[d] & \wt\bp \ar[r]\ar[d]& V_1|_{\nod^+}\ar[r]\ar[d] & 0 \\
0\ar[r] & D_{n;L}\ar[r]\ar[d] & \wt{D}_{n;L}
\ar[r]\ar[d]  & V_{n;L}|_{\nod^+}\ar[r]\ar[d] & 0 \\
0\ar[r] & D_{n-1;L}\ar[r]\ar[d] & \wt{D}_{n-1;L}\ar[r]\ar[d] 
& V_{n-1;L}|_{\nod^+}\ar[r]\ar[d] & 0 \\
& 0 & 0 & 0}$$ 
\caption{Commutative square of exact rows and columns of Fredholm operators 
for the proof of the property of the CROrient~\ref{CROdegenC_prop}\ref{CdegenSpin_it}
property}
\label{SpinCdegen_fig}
\end{figure}

The exact triple~\eref{Cdegses_e} induces exact triples of Fredholm operators given
by the rows in the diagram of Figure~\ref{SpinCdegen_fig}.
The splitting in~\eref{SpinPin1Spin_e7} induces exact triples of Fredholm operators given
by the columns in this diagram.
The right column respects the complex orientations of the associated vector spaces.
By the CROrient~\ref{CROSpinPinSES_prop}\ref{CROsesSpin_it} property,
the left (resp.~middle) column respects  the orientations 
\begin{alignat*}{3}
&\fo_{\os_1}'\big(V_1,\vph_1\big),&\quad &\fo_{\os_{n;L}}'\big(V_{n;L},\vph_{n;L}\big), 
&\quad&\hbox{and}\quad\fo_{\os_{n-1;L}}'\big(V_{n-1;L},\vph_{n-1;L}\big),\\
\big(\hbox{resp.}\quad&\fo_{\wt\os_1}\big(\wt{V}_1,\wt\vph_1\big),
&\quad&\fo_{\wt\os_{n;L}}\big(\wt{V}_{n;L},\wt\vph_{n;L}\big),
&\quad&\hbox{and}\quad
\fo_{\wt\os_{n-1;L}}\big(\wt{V}_{n-1;L},\wt\vph_{n-1;L}\big)~\big) 
\end{alignat*}
if and only if $(n\!-\!1)\binom{|\pi_0(\Si_0^{\si_0})|}{2}$ is even.\\

By definition, the rows in Figure~\ref{SpinCdegen_fig} respect 
the orientations 
$\fo_{\os}(V,\vph)$, $\fo_{\wt\os}(\wt{V},\wt\vph)$, and
the complex orientation on $V|_{\nod^+}$ with
$$\big(V,\vph,\os\big)=\big(V_1,\vph_1,\os_1\big),
\big(V_{n;L},\vph_{n;L},\os_{n;L}\big),
\big(V_{n-1;L},\vph_{n-1;L},\os_{n-1;L}\big),$$
depending on the row.
Along with the conclusion of the previous paragraph, \eref{CdegenSpin_e3}, 
Lemma~\ref{3by3_lmm}, and the evenness of the real dimension of~$V_1|_{\nod^+}$, 
this implies~that 
$$\fo_{\os_{n;L}}\big(V_{n;L},\vph_{n;L}\big)
=\fo_{\os_{n;L}}'\big(V_{n;L},\vph_{n;L}\big)
\qquad\forall\,n\!\in\!\Z^+\,.$$
By the RelSpinPin~\ref{RelSpinPinStab_prop} property,
every relative $\OSpin$-structure on the real vector bundle $V_{n;L}^{\vph_{n;L}}$ 
over $\Si_0^{\si_0}\!\subset\!\Si_0$ equals~$\os_{n;L}$ for some 
relative $\OSpin$-structure $\os_L$ on~$L^{\phi}$.
By \cite[Theorem~1.1]{RBP}, every real bundle pair $(V_0,\vph_0)$ over~$(\Si_0,\si_0)$
with $V_0^{\vph_0}$ orientable
is isomorphic to $(V_{n;L},\vph_{n;L})$ for some rank~1 real bundle pair 
$(L,\phi)$ over~$(\Si_0,\si_0)$ with $L^{\phi}$ orientable.
The last three statements imply 
the CROrient~\ref{CROdegenC_prop}\ref{CdegenSpin_it} property.
\end{proof}

\begin{proof}[{\bf{\emph{Proof of CROrient~\ref{CROdegenH3_prop}\ref{HdegenSpin_it} property}}}]
Let $\cC_0$ and $r_{\bu}$ be as in the statement of this property.
In light of the CROrient~\ref{CROos_prop}\ref{osinter_it} property, 
we can assume that $r_{\bu}\!=\!|\pi_0(\Si_0^{\si_0})|$. 
We take $\os_1$, $\bp$, $(L,\phi)$, $\os_L$, and $\bp_L$ as in the proof
of the CROrient~\ref{CROdegenC_prop}\ref{CdegenSpin_it} property and 
again define $(V_n,\vph_n)$, $(V_{n;L},\vph_{n;L})$, $\os_{n;L}$,
$D_n$, and~$D_{n;L}$ as in~\eref{SpinPin1Spin_e7} and just above and below.
We denote by $\fo_{1;\nod}$ the orientation of~$V_1^{\vph_1}|_{\nod}$ induced by~$\os_1$ and
by $\fo_{n;L;\nod}$ the orientation of~$V_{n;L}^{\vph_{n;L}}|_{\nod}$ induced 
by~$\os_{n;L}$.
By the same reasoning as below~\eref{CdegenSpin_e3}, 
with Proposition~\ref{genrk1H3degen_prp} used in place of~\ref{genrk1Cdegen_prp},
\BE{H3degenSpin_e3} 
\fo_{\os_1}\big(V_1,\vph_1\big)=\fo_{\os_1}^+\big(V_1,\vph_1\big) 
\qquad\hbox{and}\qquad
\fo_{\os_L}(L,\phi)=\fo_{\os_L}^+(L,\phi).\EE

\vspace{.15in}

\begin{figure}
$$\xymatrix{& 0\ar[d] & 0\ar[d] & 0\ar[d]  \\
0\ar[r] & \bp\ar[r]\ar[d] & \wt\bp \ar[r]\ar[d]& V_1^{\vph_1}|_{\nod}\ar[r]\ar[d] & 0 \\
0\ar[r] & D_{n;L}\ar[r]\ar[d] & \wt{D}_{n;L}
\ar[r]\ar[d]  & V_{n;L}^{\vph_{n;L}}|_{\nod}\ar[r]\ar[d] & 0 \\
0\ar[r] & D_{n-1;L}\ar[r]\ar[d] & \wt{D}_{n-1;L}\ar[r]\ar[d] 
& V_{n-1;L}^{\vph_{n-1;L}}|_{\nod}\ar[r]\ar[d] & 0 \\
& 0 & 0 & 0}$$ 
\caption{Commutative square of exact rows and columns of Fredholm operators 
for the proof of the property of the CROrient~\ref{CROdegenH3_prop}\ref{HdegenSpin_it}
property}
\label{SpinH3degen_fig}
\end{figure}

The exact triple~\eref{Rdegses_e} induces exact triples of Fredholm operators given
by the rows in the diagram of Figure~\ref{SpinH3degen_fig}.
The splitting in~\eref{SpinPin1Spin_e7} induces exact triples of Fredholm operators given
by the columns in this diagram.
The right column respects the orientations $\fo_{1;\nod}$, $\fo_{n;L;\nod}$, and~$\os_{n-1;L;\nod}$ 
of the associated vector spaces.
By the CROrient~\ref{CROSpinPinSES_prop}\ref{CROsesSpin_it} property,
the left (resp.~middle) column respects the orientations 
\begin{alignat*}{3}
&\fo_{\os_1}^+\big(V_1,\vph_1\big), &\quad & 
\fo_{\os_{n;L}}^+\big(V_{n;L},\vph_{n;L}\big), 
&\quad&\hbox{and}\quad\fo_{\os_{n-1;L}}^+\big(V_{n-1;L},\vph_{n-1;L}\big)\\
\big(\hbox{resp.}\quad &\fo_{\wt\os_1}\big(\wt{V}_1,\wt\vph_1\big),
&\quad&\fo_{\wt\os_{n;L}}\big(\wt{V}_{n;L},\wt\vph_{n;L}\big),
&\quad&\hbox{and}\quad \fo_{\wt\os_{n-1;L}}\big(\wt{V}_{n-1;L},\wt\vph_{n-1;L}\big)~\big)
\end{alignat*}
if and only if 
$(n\!-\!1)\binom{|\pi_0(\Si_0^{\si_0})|}{2}$ 
(resp.~$(n\!-\!1)\binom{|\pi_0(\Si_0^{\si_0})|+1}{2}$)
is even.\\

By definition, the rows in Figure~\ref{SpinH3degen_fig} respect 
the orientations 
$\fo_{\os}(V,\vph)$, $\fo_{\wt\os}(\wt{V},\wt\vph)$, and
the orientation~$\fo_{V;\nod}$ on $V^{\vph}|_{\nod}$ with
\begin{equation*}\begin{split}
\big(V,\vph,\os,\fo_{V;\nod}\big)=&\big(V_1,\vph_1,\os_1,\fo_{1;\nod}\big),
\big(V_{n;L},\vph_{n;L},\os_{n;L},\fo_{n;L;\nod}\big),\\
&\big(V_{n-1;L},\vph_{n-1;L},\os_{n-1;L},\fo_{n-1;L;\nod}\big),
\end{split}\end{equation*}
depending on the row.
Since 
\begin{equation*}\begin{split}
&(n\!-\!1)\binom{|\pi_0(\Si_0^{\si_0})|}{2}+(n\!-\!1)\binom{|\pi_0(\Si_0^{\si_0})|+1}{2}
+\big(\!\dim V_1^{\vph_1}|_{\nod}\big)\big(\!\ind\,D_{n-1;L}\big)\\
&\qquad\cong   (n\!-\!1)\big|\pi_0(\Si_0^{\si_0})\big|+(n\!-\!1)\big|\pi_0(\Si_0^{\si_0})\big|
\cong0 \mod2,
\end{split}\end{equation*}
the last two statements, \eref{H3degenSpin_e3}, and
Lemma~\ref{3by3_lmm} imply~that 
$$\fo_{\os_{n;L}}\big(V_{n;L},\vph_{n;L}\big)
=\fo_{\os_{n;L}}^+\big(V_{n;L},\vph_{n;L}\big)
\qquad\forall\,n\!\in\!\Z^+\,.$$
Similarly to the end of the proof of the CROrient~\ref{CROdegenC_prop}\ref{CdegenSpin_it} property,
this in turn implies the CROrient~\ref{CROdegenH3_prop}\ref{HdegenSpin_it} property. 
\end{proof}

\section{Orientations for twisted determinants}
\label{CROrientPf_sec}

Let $\cC$ be a smooth decorated marked symmetric surface as in~\eref{cCsymdfn_e} and
$(V,\vph)$ be a real bundle pair over~$\cC$.
In Section~\ref{LB_subs}, we define an orientation~$\fo_{\cC}(V^{\vph},\fo_{\x})$
of 
\BE{lacCRdfn2_e}\la_{\cC}^{\R}(V,\vph)\equiv
\la\bigg(\bigoplus_{i=1}^k\!V_{x_i}^{\vph}\bigg)
=\bigotimes_{i=1}^k\!\la(V_{x_i}^{\vph}).\EE
This orientation generally depends on a tuple $\fo_{\x}$ of orientations~$\fo_{x_r^*}$
of~fibers of~$V^{\vph}$ at one point~$x_r^*$ of each topological component~$S^1_r$
of the fixed locus~$\Si^{\si}$ of~$\cC$ and on
the decorated structure of~$\cC$.
In Section~\ref{detOrient_subs}, we define an orientation $\fo_{\fp}(V,\vph;\fo_{\x})$
of every real CR-operator~$D$ on the real bundle pair~$(V,\vph)$.
This orientation a prior depends on a relative $\Pin$-structure~$\fp$
on the real vector bundle~$V^{\vph}$ over $\Si^{\si}\!\subset\!\Si$,
on a tuple~$\fo_{\x}$ of orientations as above, and
on the decorated structure of~$\cC$.
We establish properties of the orientations $\fo_{\cC}(V^{\vph},\fo_{\x})$ and
$\fo_{\fp}(V,\vph;\fo_{\x})$, analogous to the CROrient properties
of Sections~\ref{OrientPrp_subs1} and~\ref{OrientPrp_subs2}, in Section~\ref{TargetOrient_subs} and 
in Sections~\ref{detOrient_subs} and~\ref{DegenSES_subs},
respectively.\\

The orientations $\fo_{\cC}(V^{\vph},\fo_{\x})$ and $\fo_{\fp}(V,\vph;\fo_{\x})$
determine an orientation $\fo_{\cC;\fp}(V,\vph;\fo_{\x})$ of the twisted determinant 
\BE{wtlacCDdfn_e} \wt\la_{\cC}(D)\equiv\la_{\cC}^{\R}(V,\vph)^*\!\otimes\!\la(D)\EE
of~$D$.
In Section~\ref{twisteddetLB_subs}, we readily deduce properties of the orientations
$\fo_{\cC;\fp}(V,\vph;\fo_{\x})$ from the properties of the orientations 
$\fo_{\cC}(V^{\vph},\fo_{\x})$ and $\fo_{\fp}(V,\vph;\fo_{\x})$ established
in Sections~\ref{TargetOrient_subs}-\ref{DegenSES_subs}.
By Corollary~\ref{Vvphorient_crl2}\ref{rk1horientcrl_it}, the orientation 
$\fo_{\cC;\fp}(V,\vph;\fo_{\x})$ does not depend on the choice of~$\fo_{\x}$
if $(V,\vph)$ is $\cC$-balanced;
we then denote it by~$\fo_{\cC;\fp}(V,\vph)$. 
In this case, the properties established in Section~\ref{twisteddetLB_subs} reduce
to the CROrient properties of Sections~\ref{OrientPrp_subs1} and~\ref{OrientPrp_subs2}
concerning orientations induced by $\Pin$-structures.
This completes the proof of Theorem~\ref{CROrient_thm}\ref{CROrientPin_it}.\\

Given a symmetric surface $(\Si,\si)$, possibly nodal, we denote by~$N$
the number of the topological components of its fixed locus~$\Si^{\si}$.

\subsection{Orientations of the twisting target}
\label{TargetOrient_subs}

Suppose $\cC$ is a smooth decorated marked symmetric surface as in~\eref{cCsymdfn_e}
and $(V,\vph)$ is a real bundle pair over~$\cC$.
For each topological component~$S^1_r$ of the $\si$-fixed locus $\Si^{\si}$,  let
$$k_r\in\{0\}\!\sqcup\![k] \qquad\hbox{and}\qquad
j_1^r(\cC),\ldots,j_{k_r}^r(\cC)\in[k]$$ 
be as in~\eref{posorder_e}.
We take $x_r^*\!\in\!S^1_r$ to be the real marked point $x_{j_1^r(\cC)}$ if $k_r\!>\!0$
and any point in~$S^1_r$ if $k_r\!=\!0$.
Fix a tuple
\BE{foxcCdfn_e}\fo_{\x}\equiv \big(\fo_{x_r^*}\big)_{S^1_r\in\pi_0(\Si^{\si})}\EE
of orientations $\fo_{x_r^*}$ of $V^{\vph}$ at $x_r^*$.\\

For each real marked point $x_i\!\in\!S^1_r$, the orientation $\fo_{x_r^*}$ induces
an orientation~$\fo_{x_i}$ of $V_{x_i}^{\vph}$ by the transfer along the positive direction of~$S^1_r$
determined by the chosen half-surface~$\Si^b$ of~$\Si$.
Along with the ordered decompositions~\eref{VvphcCrdfn_e} and~\eref{VRcCdfn_e}, 
the orientations~$\fo_{x_i}$ in turn determine an orientation~$\fo_{\cC}(V^{\vph},\fo_{\x})$
of~\eref{lacCRdfn2_e}.
If $\ce$ is a short exact sequence of real bundle pairs over~$(\Si,\si)$ and
\BE{foxcCdfn_e2}
\fo_{\x}'\equiv\big(\fo_{x_r^*}'\big)_{S^1_r\in\pi_0(\Si^{\si})}
\qquad\hbox{and}\qquad
\fo_{\x}''\equiv \big(\fo_{x_r^*}''\big)_{S^1_r\in\pi_0(\Si^{\si})}\EE
are tuples of orientations for fibers of $(V',\vph')$ 
and~$(V'',\vph'')$, respectively, as in~\eref{foxcCdfn_e},  we denote~by 
$$\llrr{\fo_{\x}',\fo_{\x}''}_{\ce_{\R}}\equiv \big(\fo_{x_r^*}\big)_{S^1_r\in\pi_0(\Si^{\si})}$$ 
the tuple of orientations for fibers of $(V,\vph)$ so that the restriction of 
the induced exact sequence~$\ce_{\R}$ in~\eref{RBPsesdfn_e2} to each~$x_r^*$
respects the orientations  $\fo_{x_r^*}'$, $\fo_{x_r^*}$, and $\fo_{x_r^*}''$.
The following observations are straightforward.

\begin{lmm}\label{Vvphorient_lmm}
Let $(V,\vph)$ be a rank~$n$ real bundle pair over 
a smooth decorated marked symmetric surface~$\cC$ and 
$\fo_{\x}$ be a tuple of orientations of~$V^{\vph}$ at points $x_{r}^*$
in $S^1_r\!\subset\!\Si^{\si}$ as in~\eref{foxcCdfn_e}.
\begin{enumerate}[label=(\arabic*),leftmargin=*]

\item\label{Vvphflip_it} The orientation $\fo_{\cC}(V^{\vph},\fo_{\x})$
of~$\la_{\cC}^{\R}(V,\vph)$ does not depend on the choice of half-surface $\Si_*^b$ of 
an elemental component~$\Si_*$ of $(\Si,\si)$ if and only~if 
$$\sum_{S_r^1\in\pi_0(\Si_*^{\si})}\!\!\!
\bigg(\!\!\blr{w_1(V^{\vph}),[S^1_r]_{\Z_2}}\binom{k_r\!-\!1}{1}
\!+\!n\binom{k_r\!-\!1}{2}\!\!\!\bigg) =0\in \Z_2\,.$$

\item\label{Vvphinter_it} The interchange in the ordering of two consecutive components~$S_r^1$ 
and~$S_{r+1}^1$  of~$\Si^{\si}$ preserves the orientation $\fo_{\cC}(V^{\vph},\fo_{\x})$
if and only if $nk_rk_{r+1}\!\in\!2\Z$.

\item\label{Vvphinter2_it} 
The interchange of two real marked points $x_{j_i^r(\cC)}$ and $x_{j_{i'}^r(\cC)}$
on the same connected component $S_r^1$ of~$\Si^{\si}$
with \hbox{$2\!\le\!i,i'\!\le\!k_r$} preserves $\fo_{\cC}(V^{\vph},\fo_{\x})$.
The combination of the interchange of the real points $x_{j_1^r(\cC)}$ and $x_{j_i^r(\cC)}$
with $2\!\le\!i\!\le\!k_r$ and the replacement of the component 
\hbox{$\fo_{x_r^*}\!\equiv\!\fo_{x_{j_1^r(\cC)}}$}
in~\eref{foxcCdfn_e} by $\fo_{x_{j_i^r(\cC)}}$ preserves $\fo_{\cC}(V^{\vph},\fo_{\x})$ if and only~if
\BE{Vvphinter2_e}\big(n(k_r\!-\!1)\!+\!\blr{w_1(V^{\vph}),[S^1_r]_{\Z_2}}\big)(i\!-\!1)=0\in \Z_2.\EE

\item\label{Vvphorient_it} 
The reversal of the component orientation~$\fo_{x_r^*}$ in~\eref{foxcCdfn_e}
preserves $\fo_{\cC}(V^{\vph},\fo_{\x})$ if and only if $k_r\!\in\!2\Z$.

\end{enumerate}
\end{lmm}

\begin{lmm}\label{Vvphorient_lmm1b}
Suppose $\cC$ is a smooth decorated marked symmetric surface, 
$\ce$ is a short exact sequence of real bundle pairs over~$\cC$ 
as in~\eref{RBPsesdfn_e}, and $\fo_{\x}'$ and $\fo_{\x}''$
are tuples of orientations for fibers of~$(V',\vph')$ 
and~$(V'',\vph'')$, respectively, as in~\eref{foxcCdfn_e2}.
The isomorphism~\eref{RBPsesdfn_e4} respects the orientations 
$\fo_{\cC}(V'^{\vph'},\fo_{\x}')$,
$\fo_{\cC}(V^{\vph},\llrr{\fo_{\x}',\fo_{\x}''}_{\ce_{\R}})$, 
and $\fo_{\cC}(V''^{\vph''},\fo_{\x}'')$
if and only~if \hbox{$(\rk\,V')(\rk\,V'')\binom{k}{2}$} is~even.
\end{lmm}

Suppose $\cC_0$ is a decorated marked symmetric surface as in~\eref{cCdegC_e} 
which contains precisely one conjugate pair $(\nod^+,\nod^-)$ of nodes and no other nodes
and $(V_0,\vph_0)$ is a real bundle pair over~$\cC_0$.
The fixed locus $\Si_0^{\si_0}\!\subset\!\Si_0$ in this case is 
a disjoint union of circles. 
Given a tuple~$\fo_{\x}$ of orientations for fibers of~$V_0^{\vph_0}$
as in~\eref{foxcCdfn_e}, we define an orientation~$\fo_{\cC_0}(V_0^{\vph_0},\fo_{\x})$
of 
\BE{lacC0Rdfn2_e}\la_{\cC_0}^{\R}(V_0,\vph_0)\equiv
\la\bigg(\bigoplus_{i=1}^k\!V_0^{\vph_0}|_{x_i}\bigg)
=\bigotimes_{i=1}^k\!\la\big(V_0^{\vph_0}|_{x_i}\big)\EE
as above Lemma~\ref{Vvphorient_lmm}. 
We call it the \sf{intrinsic orientation} of $\la_{\cC_0}^{\R}(V_0,\vph_0)$ induced
by~$\fo_{\x}$.
We define the \sf{limiting orientation}~$\fo_{\cC_0}'(V_0^{\vph_0},\fo_{\x})$ 
of~\eref{lacC0Rdfn2_e} induced by~$\fo_{\x}$ similarly to the construction 
of the limiting orientation above Proposition~\ref{genrk1Cdegen_prp} with $s_r^{\R}(0)\!=\!x_r^*$.
In this case, it is immediate that the orientations~$\fo_{\cC_0}(V_0^{\vph_0},\fo_{\x})$
and~$\fo_{\cC_0}'(V_0^{\vph_0},\fo_{\x})$ are the same.\\

Suppose $\cC_0$ is a decorated marked symmetric surface as in~\eref{cCdegC_e} 
which contains precisely one~$H3$ node~$\nod$ and no other nodes and
$(V_0,\vph_0)$ is a real bundle pair over~$\cC_0$.
Let $S^1_{\bu1}\!\subset\!\Si_{\bu1}$, $S^1_{\bu2}\!\subset\!\Si_{\bu2}$,
$r_{\bu},k_{\bu1},k_{\bu2}$, $\nod_1,\nod_2$, $\wt\cC_0,\wt\cC_1',\wt\cC_2'$, 
and $(\wt{V}_0,\wt\vph_0)$ be as in the paragraph containing~\eref{wtcCH3deg_e}
and just below and 
$r(\cC_0)$, $j_1'(\cC_0),j_2'(\cC_0)$, and $\de_{\R}(\cC_0)$
be as in the sentence containing~\eref{rcC0dfn_e}
and just below.
The fixed locus $\Si_0^{\si_0}\!\subset\!\Si_0$ in this case is 
a disjoint union of circles $S^1_r$ with $r\!\neq\!r_{\bu}$ and 
of the wedge $S^1_{\bu1}\!\cup_{\nod}\!S^1_{\bu2}$ of two circles.
We take $x_r^*\!\in\!S^1_r$ for $r\!\neq\!r_{\bu}$ as above~\eref{foxcCdfn_e}
and $x_{r_{\bu}}^*\!=\!\nod$.
For $r\!=\!1,2$, let 
$$x_{\bu r}=\begin{cases}x_i,&\hbox{if}~
i\!\equiv\!\inf\big\{j\!\in\![k]\!:x_j\!\in\!S^1_{\bu r}\big\}\!<\!\i;\\
\nod_r,&\hbox{otherwise}.\end{cases}$$
Given a tuple of orientations of $V_0^{\vph_0}|_{x_r^*}$ as in~\eref{foxcCdfn_e},
denote by~$\fo_{\nod_1}$ and~$\fo_{\nod_2}$ the induced orientations of~$\wt{V}_0^{\wt\vph_0}$
at~$\nod_1$ and~$\nod_2$, respectively.
Let $\fo_{\bu1}$ be the orientation of $V_0^{\vph_0}|_{x_{\bu1}}$ 
obtained by transferring~$\fo_{\nod}$ along the positive direction of~$S^1_{\bu1}$
determined by the chosen half-surface~$\Si^b$ of~$\Si$.
We denote by $\fo_{\bu2}$ (resp.~$\fo_{\bu2}'$) the orientation of $V_0^{\vph_0}|_{x_{\bu2}}$ 
obtained by transferring~$\fo_{\nod}$ along the positive direction of~$S^1_{\bu2}$
(resp.~transferring~$\fo_{\nod}$ first around~$S^1_{\bu1}$ back
to an orientation~$\fo_{\nod}'$ of  $V_0^{\vph_0}|_{\nod}$
and then transferring~$\fo_{\nod}'$ along the positive direction of~$S^1_{\bu2}$).\\

The tuples
$$\wt\fo_{\x;1}\equiv\big(\fo_{\bu1},\fo_{\nod_2},
(\fo_{x_r^*})_{S^1_r\in\pi_0(\Si_0^{\si_0}),r\neq r_{\bu}}\big)
\quad\hbox{and}\quad
\wt\fo_{\x;2}\equiv\big(\fo_{\nod_1},\fo_{\bu2},
(\fo_{x_r^*})_{S^1_r\in\pi_0(\Si_0^{\si_0}),r\neq r_{\bu}}\big)$$
of orientations of fibers of $\wt{V}_0^{\wt\vph_0}$
determine orientations 
\BE{H3orientlift_e}\wt\fo_{\cC_0;1}\big(V_0^{\vph_0},\fo_{\x}\big)
\equiv \fo_{\wt\cC_1'}\big(\wt{V}_0^{\wt\vph_0},\wt\fo_{\x;1}\big)
\quad\hbox{and}\quad
\wt\fo_{\cC_0;2}\big(V_0^{\vph_0},\fo_{\x}\big)
\equiv \fo_{\wt\cC_2'}\big(\wt{V}_0^{\wt\vph_0},\wt\fo_{\x;2}\big)\EE
of the lines
\BE{H3orientisom_e}\la_{\wt\cC_1'}^{\R}(\wt{V}_0,\wt\vph_0)\approx
\la_{\cC_0}^{\R}(V_0,\vph_0)\!\otimes\!\la\big(V_0^{\vph_0}|_{\nod}\big)
\quad\hbox{and}\quad
\la_{\wt\cC_2'}^{\R}(\wt{V}_0,\wt\vph_0)\approx
\la_{\cC_0}^{\R}(V_0,\vph_0)\!\otimes\!\la\big(V_0^{\vph_0}|_{\nod}\big),\EE
respectively.
We denote the orientations of $\la_{\cC_0}^{\R}(V_0,\vph_0)$ induced by 
the orientations~\eref{H3orientlift_e} and $\fo_{\nod}$ via the isomorphisms~\eref{H3orientisom_e}
by $\fo_{\cC_0;1}(V_0^{\vph_0},\fo_{\x})$ and 
$\fo_{\cC_0;2}(V_0^{\vph_0},\fo_{\x})$, respectively.\\

Suppose in addition that $(\cU,\wt\fc)$ is a flat family of deformations of $\cC_0$
as in~\eref{cUsymmdfn_e} over a disk $\De\!\subset\!\C$,
$\De_{\R}^+\!\subset\!\De_{\R}^*$ is as above~\eref{Rdegses_e}, 
and $(V,\vph)$ is a real bundle pair over $(\cU,\wt\fc)$ extending~$(V_0,\vph_0)$.
Let $s_1^{\R},\ldots,s_N^{\R}$ be
sections of $\cU^{\wt\fc}$ over~$\De_{\R}$ so~that 
$$s_r^{\R}(0)=x_r^*~~\forall\,r\!\neq\!r_{\bu}, \qquad
s_{r_{\bu}}^{\R}(0)=x_{\bu r(\cC_0)}.$$
For each $\t\!\in\!\De_{\R}$ and $r\!\neq\!r_{\bu}$, 
the orientation $\fo_{x_r^*}$ induces an orientation $\fo_{x_r;\t}$ of $V^{\vph}$
at $s_r^{\R}(\t)$ via the vector bundle $s_r^{\R*}V^{\vph}$.
If $r(\cC_0)\!=\!1$ (resp.~$r(\cC_0)\!=\!2$),
the orientation~$\fo_{\bu1}$ (resp.~$\fo_{\bu2}'$)  similarly induces an orientation 
$\fo_{x_{r_{\bu}}^*;\t}$ of $V^{\vph}$ at~$s_{r_{\bu}}^{\R}(\t)$. 
Let 
$$\fo_{\x;\t}=\big(\fo_{x_r^*}\big)_{S^1_r\in\pi_0(\Si^{\si})}, 
\qquad \fo_{\t}=\fo_{\cC_{\t}}\big(V_{\t}^{\vph_{\t}},\fo_{\x;\t}\big).$$
We denote~by 
\BE{cCOlimor_e4} \fo_{\cC_0}^+\big(\fo_{\x}\big)\equiv 
\fo_{\cC_0}^+\big(V_0^{\vph_0},\fo_{\x}\big)\EE
the orientation of~$\la_{\cC_0}^{\R}(V_0,\vph_0)$ 
obtained as the continuous extension of the orientations~$\fo_{\t}$
with $\t\!\in\!\De_{\R}^+$.\\

We define
$$\ep_{\R}\big(\cC_0;V_0^{\vph_0}\big)
=\blr{w_1(V_0^{\vph_0}),[S^1_{\bu1}]_{\Z_2}\!+\![S^1_{\bu2}]_{\Z_2}}
\big(j_{r(\cC_0)}'(\cC_0)\!+\!(r(\cC_0)\!-\!1)k_{\bu1}\big).$$
For $r^*\!=\!1,2$, let 
\begin{gather*}
\de_{r^*}(\cC_0)=\big(k_{\bu r^*}\!-\!1\big)
\big(j_{r^*}'(\cC_0)\!+\!r^*\!\big),\\
\ep_{r^*}\big(\cC_0;V_0^{\vph_0}\big)=
\blr{w_1(V_0^{\vph_0}),[S^1_{\bu r^*}]_{\Z_2}}j_{r^*}'(\cC_0)
+\blr{w_1(V_0^{\vph_0}),\big[S^1_{\bu1}]_{\Z_2}}k_{\bu2}.
\end{gather*}

\begin{lmm}\label{VvphorientH3deg_lmm1}
Suppose $\cC_0$ is a decorated marked symmetric surface which contains precisely 
one $H3$ node and no other nodes,
$(V_0,\vph_0)$  is a rank~$n$ real bundle pair over~$\cC_0$,
and $\fo_{\x}$ is  a tuple of orientations of~$V_0^{\vph_0}$ at points 
$x_r^*\!\in\!S^1_r$ for $r\!\neq\!r_{\bu}$ and $x_{r_{\bu}}^*\!=\!\nod$.
Let $r^*\!=\!1,2$.
The intrinsic orientation~$\fo_{\cC_0;r^*}(V_0^{\vph_0},\fo_{\x})$  
and the limiting orientation~\eref{cCOlimor_e4}
of~$\la_{\cC_0}^{\R}(V_0,\vph_0)$ are the same if and only~if
\begin{equation*}\begin{split}
&n\bigg(\!\!\big(k_{r_{\bu}}\!-\!1\!\big)\de_{\R}(\cC_0)\!+\!\de_{r^*}(\cC_0)
+\!\!\!\!\!\!\!\sum_{\begin{subarray}{c}S_r^1\in\pi_0(\Si_0^{\si_0})\\ 
r>r_{\bu}\end{subarray}}\!\!\!\!\!\!\!k_r~\bigg)
\!+\!\ep_{\R}\big(\cC_0;V_0^{\vph_0}\big)\!+\!\ep_{r^*}\big(\cC_0;V_0^{\vph_0}\big)\\
&\hspace{4in}
=n\big(r^*\!-\!1\big)\!+\!2\Z\,.
\end{split}\end{equation*}
\end{lmm}

\begin{proof} In light of Lemma~\ref{Vvphorient_lmm}\ref{Vvphinter_it},
we can assume that $r_{\bu}\!=\!|\pi_0(\Si_0^{\si_0})|$.
By the definitions of intrinsic and limiting orientation, we can also assume
that $k_r\!=\!0$ for all $r\!\neq\!r_{\bu}$.
Let 
$$\ep_1=\blr{w_1(V_0^{\vph_0}),[S^1_{\bu1}]_{\Z_2}}, \qquad
\ep_2=\blr{w_1(V_0^{\vph_0}),[S^1_{\bu2}]_{\Z_2}}\,.$$
We denote by $\fo_{\nod}''$ (resp.~$\fo_{\nod}^{\dag}$)
the orientation  of $V_0^{\vph_0}|_{\nod}$ 
obtained by transferring the orientation~$\fo_{\nod}$ (resp.~$\fo_{\nod}'$) 
around~$S^1_{\bu2}$ back to~$\nod$.\\

Let $r\!=\!1,2$. 
The decorated structure of~$\cC_0$ determines an orientation on~$S^1_{\bu r}$.
If $k_{\bu r}\!\neq\!0$, this decorated structure thus determines an oriented
arc $S^1_{\bu r;1}\!\subset\!S^1_{\bu r}$ from the node~$\nod$ to 
the marked point $x_{\bu r}\!\in\!S^1_{\bu r}$
and an oriented arc $S^1_{\bu r;2}$ from this marked point 
back to~$\nod$.
The first (resp.~second) arc carries $j_1'(\cC_0)$ (resp.~$k_{\bu r}\!-\!1\!-\!j_1'(\cC_0)$)
marked points~$x_i$ of~$\cC_0$ other than~$x_{\bu r}$.
Let $V_{r;1}^{\R}$ (resp.~$V_{r;2}^{\R}$) be the ordered direct sum of 
the $j_r'(\cC_0)$ (resp.~$k_{\bu r}\!-\!j_r'(\cC_0)$) fibers of $V_0^{\vph_0}$ 
at the marked points $x_i\!\in\!S^1_{\bu r;1}$ other than $x_{\bu r}$
(resp.~$x_i\!\in\!S^1_{\bu r;2}$ including $x_{\bu r}$)
in the order they appear on the arc.
If $k_{\bu r}\!=\!0$, we take  $V_{r;1}^{\R}$ and $V_{r;2}^{\R}$ to 
be the zero vector spaces.
The orientation $\fo_{\nod}$ (resp.~$\fo_{\nod}'$, $\fo_{\nod}''$, $\fo_{\nod}^{\dag}$) 
of $V_0^{\vph_0}|_{\nod}$ determines an orientation~$\fo_{x_i}$ 
(resp.~$\fo_{x_i}'$, $\fo_{x_i}''$, $\fo_{x_i}^{\dag}$)
of $V_0^{\vph_0}$ at each marked point $x_i\!\in\!S^1_{\bu r}$
by the transfer in the positive direction of~$S^1_{\bu r}$.
If in addition $s\!=\!1,2$, let $\fo_{r;s}$
(resp.~$\fo_{r;s}'$, $\fo_{r;s}''$, $\fo_{r;s}^{\dag}$)
be the resulting orientation of~$V_{r;s}^{\R}$.\\

By definition and the assumptions above,
$$\big(V_0^{\vph_0}\big)_{\cC_0}\equiv \big(V_0^{\vph_0}\big)_{\cC_0;r_{\bu}}
=V_{1;1}^{\R}\!\oplus\!V_{1;2}^{\R}\!\oplus\!V_{2;1}^{\R}\!\oplus\!V_{2;2}^{\R}\,.$$
The two intrinsic orientations of the lemma are described~by
\begin{alignat}{3}\label{VvphorientH3deg1_e3a}
\fo_{\cC_0;1}\big(V_0^{\vph_0},\fo_{\x}\big)&=
\fo_{1;2}\!\oplus\!\fo_{1;1}'\!\oplus\!\fo_{2;1}\!\oplus\!\fo_{2;2}
&\qquad&\Llra&\qquad nk_{\bu2}&\in\!2\Z,\\
\label{VvphorientH3deg1_e3b}
\fo_{\cC_0;2}\big(V_0^{\vph_0},\fo_{\x}\big)&=
\fo_{1;1}\!\oplus\!\fo_{1;2}\!\oplus\!\fo_{2;2}\!\oplus\!\fo_{2;1}''
&\qquad&\Llra&\qquad nk&\in\!2\Z.
\end{alignat}
The limiting orientation of the lemma is described by 
\BE{VvphorientH3deg1_e5}
\fo_{\cC_0}^+\big(\fo_{\x}\big)=\begin{cases}
\fo_{1;2}\!\oplus\!\fo_{2;1}'\!\oplus\!\fo_{2;2}'\!\oplus\!\fo_{1;1}^{\dag},
&\hbox{if}~r(\cC_0)\!=\!1;\\
\fo_{2;2}'\!\oplus\!\fo_{1;1}^{\dag}\!\oplus\!\fo_{1;2}^{\dag}\!\oplus\!\fo_{2;1}'',
&\hbox{if}~r(\cC_0)\!=\!2.
\end{cases}\EE
Combining \eref{VvphorientH3deg1_e3a}-\eref{VvphorientH3deg1_e5} with
\begin{alignat*}{2}
\fo_{1;1}^{\dag}&=\begin{cases}\fo_{1;1}'&\hbox{iff}~j_1'(\cC_0)\ep_2\!=\!0;\\
\fo_{1;1}&\hbox{iff}~j_1'(\cC_0)(\ep_1\!+\!\ep_2)\!=\!0;
\end{cases}
&\qquad
\fo_{2;1}''&=\begin{cases}\fo_{2;1}&\hbox{iff}~j_2'(\cC_0)\ep_2\!=\!0;\\
\fo_{2;1}'&\hbox{iff}~j_2'(\cC_0)(\ep_1\!+\!\ep_2)\!=\!0;
\end{cases}\\
\fo_{1;2}^{\dag}&=\fo_{1;2}\quad\hbox{iff}~
\big(k_{\bu1}\!-\!j_1'(\cC_0)\!\big)(\ep_1\!+\!\ep_2)\!=\!0,
&\qquad
\fo_{2;2}'&=\fo_{2;2}\quad\hbox{iff}~
\big(k_{\bu2}\!-\!j_2'(\cC_0)\!\big)\ep_1\!=\!0,
\end{alignat*}
we obtain the claim.
\end{proof}

\subsection{Orientations of real CR-operators}
\label{detOrient_subs}

Suppose $(\Si,\si)$ is a smooth decorated symmetric surface and
$(V,\vph)$ is a real bundle over~$(\Si,\si)$. 
A tuple~$\fo_{\x}$ of orientations as in~\eref{fobfxfn_e} determines 
an orientation $\fo(\la(V,\vph);\fo_{\x})$
of every real CR-operator~$D_{\la}$  on the rank~1 real bundle pair $\la(V,\vph)$
as above Proposition~\ref{cOpincg_prp}.
A relative $\Pin^{\pm}$-structure~$\fp$ on the real vector bundle $V^{\vph}$ 
over $\Si^{\si}\!\subset\!\Si$ is a relative $\OSpin$-structure on
the real vector bundle
$$V_{\pm}^{\vph_{\pm}}= V^{\vph}\!\oplus\!(2\!\pm\!1)\la(V^{\vph})$$
compatible with the canonical orientation of this vector bundle.
Via the construction at the beginning of Section~\ref{SpinOrient_subs}, 
$\fp$ thus determines an orientation $\fo_{\fp}(V_{\pm},\vph_{\pm})$ 
of every real CR-operator~$D_{\pm}$ on
the real bundle pair~$(V_{\pm},\vph_{\pm})$. 
Along with the canonical homotopy class of isomorphisms~\eref{PinOrientDfn_e}, 
the orientations $\fo(\la(V,\vph);\fo_{\x})$ and $\fo_{\fp}(V_{\pm},\vph_{\pm})$ determine 
an orientation $\fo_{\fp}(V,\vph;\fo_{\x})$ 
of every real CR-operator $D$ on the real bundle pair~$(V,\vph)$. \\

Suppose $(V,\vph)$ is a rank~1 odd-degree real bundle pair over $(S^2,\tau)$.
In particular, the real line bundle~$V^{\vph}$ over $S^1\!\subset\!S^2$ is not orientable.
Let $D$ be a real CR-operator on $(V,\vph)$.
For an orientation $\fo_{x_1}$ of~$V^{\vph}$ at a point~$x_1$ in $S^1\!\subset\!S^2$,
we denote~by 
\BE{Pin0Orient_e}\fo_0^{\pm}\big(V,\vph;\fo_{x_1}\big)\equiv 
\fo_{\io_{S^2}(\fp_0^{\pm}(V^{\vph}))}\big(V,\vph;\fo_{x_1}\big)\EE
the orientation of $D$ determined by the image~\eref{PinNormDfn_e}
of the $\Pin^{\pm}$-structure $\fp_0^{\pm}(V^{\vph})$ of
Examples~\ref{Pin1pMB_eg}, \ref{Pin1mMB_eg}, and~\ref{SpinDfn1to3_eg}
under the second map in~\eref{vsSpinPin_e0} with \hbox{$X\!=\!S^2$}.

\begin{lmm}\label{Vvphorient_lmm2}
Suppose $(\Si,\si)$ is a smooth decorated symmetric surface, 
$(V,\vph)$ is a rank~$n$ real bundle pair over~$(\Si,\si)$,
$\fp$ is a relative $\Pin^{\pm}$-structure on the real vector bundle $V^{\vph}$ 
over $\Si^{\si}\!\subset\!\Si$,
and $D$ is a real CR-operator on~$(V,\vph)$.
Let $\fo_{\x}$ be a tuple of orientations of~$V^{\vph}$ at points $x_r$
in $S^1_r\!\subset\!\Si^{\si}$ as in~\eref{fobfxfn_e}.
\begin{enumerate}[label=(\arabic*),leftmargin=*]

\item\label{rk1flip_it} The orientation $\fo_{\fp}(V,\vph;\fo_{\x})$
of~$D$ 
does not depend on the choice of half-surface $\Si_*^b$ of 
an elemental component~$\Si_*$ of $(\Si,\si)$ if and only~if 
$\vp_{\fp}(\Si_*)\!=\!0$.

\item\label{kr1inter_it} The interchange in the ordering of two consecutive components~$S_r^1$ 
and~$S_{r+1}^1$  of~$\Si^{\si}$ preserves the orientation $\fo_{\fp}(V,\vph;\fo_{\x})$
if and only~if 
$$n\!+\!2\Z\neq
\big(\blr{w_1(V^{\vph}),[S^1_r]_{\Z_2}}\!+\!1\big)
\big(\blr{w_1(V^{\vph}),[S^1_{r+1}]_{\Z_2}}\!+\!1\big)\in \Z_2.$$

\item\label{rk1horient_it} 
The reversal of the component orientation~$\fo_{x_r}$ in~\eref{fobfxfn_e}
preserves $\fo_{\fp}(V,\vph;\fo_{\x})$ if and only if 
$w_1(V^{\vph})|_{S^1_r}\!\neq\!0$.

\item\label{rk1eta_it} If in addition $\eta\!\in\!H^2(\Si,\Si^{\si};\Z_2)$,
then the orientations $\fo_{\fp}(V,\vph;\fo_{\x})$ and $\fo_{\eta\cdot\fp}(V,\vph;\fo_{\x})$ 
are the same if and only~if $\lr{\eta,[\Si^b]_{\Z_2}}\!=\!0$.

\end{enumerate}
\end{lmm}

\begin{proof}
The first statement follows from the CROrient~\ref{CROos_prop}\ref{osflip_it} property 
applied  to $(V_{\pm},\vph_{\pm})$ and
Proposition~\ref{gencOpincg_prp}\ref{gencOpincg_it1a} applied to~$\la(V,\vph)$.
The second claim follows from the CROrient~\ref{CROos_prop}\ref{osinter_it} property 
applied  to $(V_{\pm},\vph_{\pm})$ and
Proposition~\ref{gencOpincg_prp}\ref{gencOpincg_it1b} applied to~$\la(V,\vph)$.
The third claim follows from Proposition~\ref{gencOpincg_prp}\ref{gencOpincg_it2}
applied to $\la(V,\vph)$.
The last statement follows from the RelSpinPin~\ref{RelSpinPinCorr_prop} property
on page~\pageref{RelSpinPinCorr_prop}
and the first part of the CROrient~\ref{CROSpinPinStr_prop}\ref{CROSpinStr_it} property
applied to~$(V_{\pm},\vph_{\pm})$.
\end{proof}

\begin{lmm}\label{Vvphorient_lmm3}
Suppose $(\Si,\si)$, $(V,\vph)$, $\fp$, and $D$ are as in Lemma~\ref{Vvphorient_lmm2}.
\begin{enumerate}[label=(\arabic*),leftmargin=*]

\item\label{Pin2SpinRed_it} If $\fo$ is an orientation on~$V^{\vph}$ and
$\fo_{\x}$ is a tuple of orientations of~$V^{\vph}$ at points $x_r$
in $S^1_r\!\subset\!\Si^{\si}$ as in~\eref{fobfxfn_e} obtained by restricting~$\fo$,
then the orientations $\fo_{\fp}(V,\vph;\fo_{\x})$ and $\fo_{\fR_{\fo}^{\pm}(\fp)}(V,\vph)$ 
of $D$ are the same if and only if \eref{CROPin2SpinRed_e} holds.

\item\label{CROnorm_it} If $(V,\vph)$ is a rank~1 degree~1 real bundle pair over~$(S^2,\tau)$
and $\fo_{x_1}$ is an orientation of~$V^{\vph}$ at a point~$x_1$ in $S^1\!\subset\!S^2$, 
then the orientation~\eref{Pin0Orient_e} is the intrinsic orientation
$\fo(V,\vph;\fo_{x_1})$ of Proposition~\ref{cOpincg_prp}.

\item\label{CRODisjUn_it} If $\fo_{\x}$ is a tuple of orientations of~$V^{\vph}$  
at points $x_r$ in $S^1_r\!\subset\!\Si^{\si}$ as in~\eref{fobfxfn_e} and
\eref{Si1Si2dfn_e} is a decomposition of~$(\Si,\si)$ into 
decorated symmetric surfaces
of genera~$g_1$ and~$g_2$, respectively, then the isomorphism~\eref{CRODisjUn_e0} 
respects the orientations 
\BE{CRODisjUn_e2}\fo_{\fp}\equiv \fo_{\fp}\big(V,\vph;\fo_{\x}\big), ~~ 
\fo_{\fp;1}\equiv\fo_{\fp|_{\Si_1}}\!\big(V|_{\Si_1},\vph|_{\Si_1};\fo_{\x}|_{\Si_1}\big), 
~~
\fo_{\fp;2}\equiv \fo_{\fp|_{\Si_2}}\!\big(V|_{\Si_2},\vph|_{\Si_2};\fo_{\x}|_{\Si_2}\big),\EE
if and only~if \eref{CRODisjUn_e0c} with $V_1\!\equiv\!V|_{\Si_1}$ and $V_2\!\equiv\!V|_{\Si_2}$
holds.

\end{enumerate}
\end{lmm}

\begin{proof} \ref{Pin2SpinRed_it} 
By the RelSpinPin~\ref{RelSpinPinSES_prop}\ref{RelDSPin2Spin_it} 
and CROrient~\ref{CROSpinPinSES_prop}\ref{CROsesSpin_it} properties
(the latter applied three times if \hbox{$\fp\!\in\!\cP_{\Si}^+(V^{\vph})$}), 
the canonical homotopy class of isomorphisms~\eref{PinOrientDfn_e}
respects the orientations $\fo_{\fp}(V_{\pm},\vph_{\pm})$ of $D_{(V,\vph)_{\pm}}$,
$\fo_{\fR_{\fo}^{\pm}(\fp)}(V,\vph)$ of $D_{(V,\vph)}$,
and
\BE{Pin2SpinRed_e3}\fo_0\big(\la(V,\vph);\la(\fo)\big)
\equiv\fo_{\io_{\Si}(\os_0(\la(V,\vph),\la(\fo)))}\big(\la(V,\vph)\big)\EE
of $D_{\la(V,\vph)}$, where $\io_{\Si}$ is the first map in~\eref{vsSpinPin_e0} 
with \hbox{$X\!=\!\Si$},  
if and only if \eref{CROPin2SpinRed_e} holds.
By Proposition~\ref{genspinorient_prp}\ref{genrk1spindfn_it}, 
the orientation~\eref{Pin2SpinRed_e3} equals $\fo(\la(V,\vph),\fo_{\x})$.
These two statements imply the claim.\\

\ref{CROnorm_it} In this case, the real bundle pair $(V,\vph)_{\pm}$
is $3\!\pm\!1$ copies of~$(V,\vph)$.
By Proposition~\ref{tauspinorient_prp}\ref{PinCan_it}, the isomorphism~\eref{CanSpinOrien_e}
respects the orientation
$$\fo_0^{\pm}\big((3\!\pm\!1)(V,\vph)\big)=
\fo_{\io_{S^2}(\fp_0^{\pm}(V^{\vph}))}\big((3\!\pm\!1)(V,\vph)\big)$$
of $(3\!\pm\!1)D$ induced by the relative $\OSpin$-structure $\io_{S^2}(\fp_0^{\pm}(V^{\vph}))$
on the real vector bundle $V_{\pm}^{\vph_{\pm}}$ over $S^1\!\subset\!S^2$ and
the intrinsic orientation $\fo(V,\vph;\fo_{x_1})$ on the factors~$\la(D)$
on the right-hand side of~\eref{CanSpinOrien_e}.
By definition, this means that the orientation~\eref{Pin0Orient_e} on the first factor 
of~$\la(D)$ equals~$\fo(V,\vph;\fo_{x_1})$.\\

\ref{CRODisjUn_it} For $r\!=\!1,2$, we denote by $V_{\pm;r}$,
$\vph_r$, $\vph_{\pm;r}$, $\fp_r$, $D_r$, $D_{\la;r}$, $D_{\pm;r}$, and $\fo_{\x}^{(r)}$ 
the restrictions of~$V_{\pm}$, $\vph$, $\vph_{\pm}$, $\fp$, $D$,
a real CR-operator $D_{\la}$ on $\la(V,\vph)$, the real CR-operator
\hbox{$D_{\pm}\!\equiv\!D\!\oplus\!D_{\la}$} on $(V_{\pm},\vph_{\pm})$, 
and the tuple~$\fo_{\x}$ in~\eref{fobfxfn_e}, respectively, to~$\Si_r$.
The splitting in~\eref{RBPstab_e} induces the exact triples of Fredholm 
operators given by the rows in the commutative diagram of Figure~\ref{CRODisjUn_fig}.
The splitting~\eref{D1D2split_e} induces the exact triples of Fredholm 
operators given by the left and middle columns in this diagram and 
by the individual summands in the right column.
The isomorphism~\eref{CRODisjUn_e0} is induced by the left column in this diagram.\\

\begin{figure}
$$\xymatrix{& 0\ar[d] & 0\ar[d] & 0\ar[d]  \\
0\ar[r] & D_1\ar[r]\ar[d] & D_{\pm;1}\ar[r]\ar[d] 
& (2\!\pm\!1)D_{\la;1}\ar[r]\ar[d] & 0 \\
0\ar[r] & D\ar[r]\ar[d] & D_{\pm} \ar[r]\ar[d]  & (2\!\pm\!1)D_{\la}\ar[r]\ar[d] & 0 \\
0\ar[r] & D_2\ar[r]\ar[d] & D_{\pm;2} \ar[r]\ar[d]& (2\!\pm\!1)D_{\la;2}\ar[r]\ar[d] & 0 \\
& 0 & 0 & 0}$$ 
\caption{Commutative square of exact rows and columns of Fredholm operators 
for the proof of Lemma~\ref{Vvphorient_lmm3}\ref{CRODisjUn_it}}
\label{CRODisjUn_fig}
\end{figure}

By definition, the middle (resp.~top/bottom) row respects the orientations 
$\fo_{\fp}$, $\fo_{\fp}(V_{\pm},\vph_{\pm})$, and
$\fo(\la(V,\vph);\fo_{\x})$
(resp.~$\fo_{\fp;r}$, $\fo_{\fp_r}(V_{\pm;r},\vph_{\pm;r})$, and
$\fo(\la(V_r,\vph_r);\fo_{\x}^{(r)})$).
By the CROrient~\ref{CRODisjUn_prop}\ref{osDisjUn_it} property
(resp.~Proposition~\ref{gencOpincg_prp}\ref{gencOdisjun_it}),
the middle column (resp.~the exact triple formed by the individual summands in the right column),
respects the orientations 
$\fo_{\fp_1}(V_{\pm;1},\vph_{\pm;1})$, $\fo_{\fp}(V_{\pm},\vph_{\pm})$,
and $\fo_{\fp_2}(V_{\pm;2},\vph_{\pm;2})$
(resp.\,$\fo(\la(V_1,\vph_1);\fo_{\x}^{(1)})$,
$\fo(\la(V,\vph);\fo_{\x})$, and $\fo(\la(V_2,\vph_2);\fo_{\x}^{(2)})$).
Since
\BE{CRODisjUn_e5}\big(\ind\,D_{\la;1}\big)\big(\ind\,D_2\big)=
\big(1\!-\!g_1\!+\!\deg V_1\big)\big(\!(1\!-\!g_2)n\!+\!\deg V_2\big),\EE
Lemma~\ref{3by3_lmm} and the last two sentences imply that the left column in the diagram
respects the orientations~\eref{CRODisjUn_e2} 
in the $\Pin^-$-case if and only if the number~\eref{CRODisjUn_e5} is even.\\

Since the exact triple formed by the individual summands in the right column in 
Figure~\ref{CRODisjUn_fig} respects the orientations $\fo(\la(V_1,\vph_1);\fo_{\x}^{(1)})$,
$\fo(\la(V,\vph);\fo_{\x})$, and $\fo(\la(V_2,\vph_2);\fo_{\x}^{(2)})$,
Lemma~\ref{3by3_lmm} implies that the right column in the $\Pin^+$-case respects
the direct sum orientations if and only if the number 
$$\big(\ind\,D_{\la;1}\big)\big(\ind\,D_{\la;2}\big)=
\big(1\!-\!g_1\!+\!\deg V_1\big)\big(1\!-\!g_2\!+\!\deg V_2\big)$$
is even.
Since
$$\big(\ind\,3D_{\la;1}\big)\big(\ind\,D_2\big)\cong
\big(\ind\,D_{\la;1}\big)\big(\ind\,D_2\big) \mod2,$$
Lemma~\ref{3by3_lmm} again, \eref{CRODisjUn_e5}, and the two sentences above~\eref{CRODisjUn_e5}
imply that the left column 
respects the orientations~\eref{CRODisjUn_e2}
in the $\Pin^+$-case if and only if the number
$$\big(\ind\,D_{\la;1}\big)\big(\ind\,D_{\la;2}\big)
+\big(\ind\,3D_{\la;1}\big)\big(\ind\,D_2\big)\cong
(n\!+\!1)\big(1\!-\!g_1\!+\!\deg V_1\big)(1\!-\!g_2)\mod2$$
is even.
\end{proof}

\subsection{Degenerations and exact triples}
\label{DegenSES_subs}

We now describe the behavior of the orientations $\fo_{\fp}(V,\vph;\fo_{\x})$ 
under flat degenerations of~$(\Si,\si)$ to nodal symmetric surfaces as
in the CROrient~\ref{CROdegenC_prop} and~\ref{CROdegenH3_prop} properties 
of Section~\ref{OrientPrp_subs2}.
We first suppose that $\cC_0$ is a decorated symmetric surface with 
one conjugate pair of nodes~$\nod^{\pm}$ as in~\eref{cCdegC_e}
and in the top left diagram of Figure~\ref{genCdegen_fig} on page~\pageref{genCdegen_fig}
so that $x_r\!\in\!S^1_r$ for each $r\!\in\![N]$.
We also suppose that $(V_0,\vph_0)$ is a real bundle pair over~$(\Si_0,\si_0)$.
Let $\wt\cC_0$, $(\wt\Si_0,\wt\si_0)$, and $\wt\Si_0^b$ be as below~\eref{cCdegC_e},
$(\wt{V}_0,\wt\vph_0)$ be the lift of $(V_0,\vph_0)$ to a real bundle pair over $(\wt\Si_0,\wt\si_0)$.
A relative $\Pin^{\pm}$-structure~$\fp_0$ on the real vector bundle $V_0^{\vph_0}$
over $\Si_0^{\si_0}\!\subset\!\Si_0$ lifts to 
a relative $\Pin^{\pm}$-structure~$\wt\fp_0$ on the real vector bundle $\wt{V}_0^{\wt\vph_0}$
over $\wt\Si_0^{\wt\si_0}\!\subset\!\wt\Si_0$.
A tuple~$\fo_{\x}$ of orientations of the fibers of $V_0^{\vph_0}|_{x_r}$ 
as in~\eref{fobfxfn_e} lifts to a tuple~$\wt\fo_{\x}$ of orientations of
the corresponding fibers of~$\wt{V}_0^{\wt\vph_0}$.\\

Let $D_0$ be a real CR-operator on $(V_0,\vph_0)$.
We denote its lift to a real CR-operator on $(\wt{V}_0,\wt\vph_0)$ by~$\wt{D}_0$. 
The orientation 
\BE{Cdegenorient_e3b}
\wt\fo_{\fp_0}\big(\fo_{\x}\big)\equiv\fo_{\wt\fp_0}\big(\wt{V}_0,\wt\vph_0;\wt\fo_{\x}\big)\EE 
of $\wt{D}_0$ and the complex orientation of $V_0|_{\nod^+}$ determine an orientation 
\BE{rk1Csplit_e4b}\fo_{\fp_0}\big(\fo_{\x}\big)\equiv \fo_{\fp_0}\big(V_0,\vph_0;\fo_{\x}\big)\EE
of $D_0$ via the isomorphism~\eref{Cdegses_e2}.
In an analogy with the intrinsic orientation of the CROrient~\ref{CROdegenC_prop} property,
we call \eref{rk1Csplit_e4b} the \sf{intrinsic orientation} of~$D_0$
induced by~$\fp_0$ and~$\fo_{\x}$.\\

Suppose in addition that $(\cU,\wt\fc)$ is a flat family of deformations of~$\cC_0$
as in~\eref{cUsymmdfn_e},
$(V,\vph)$ is a real bundle pair over $(\cU,\wt\fc)$ extending~$(V_0,\vph_0)$,
$s_1^{\R},\ldots,s_N^{\R}$ are sections of $\cU^{\wt\fc}$ over~$\De_{\R}$ with 
\hbox{$s_r^{\R}(0)\!=\!x_r$} for all $S^1_r\!\in\!\pi_0(\Si_0^{\si_0})$,
and  \hbox{$\cD\!\equiv\!\{D_{\t}\}$} is a family of real CR-operators on
$(V_{\t},\vph_{\t})$ as in~\eref{DVvphdfn_e} extending~$D_0$.
The decorated structure on~$\cC_0$, the relative $\Pin^{\pm}$-structure~$\fp_0$ on~$V_0^{\vph_0}$,
and the tuple of orientations~$\fo_{\x}$ of fibers of~$V_0^{\vph_0}$ 
induce a decorated structure on the fiber $(\Si_{\t},\si_{\t})$ of~$\pi$ 
for every $\t\!\in\!\De_{\R}$,
a relative $\Pin^{\pm}$-structure~$\fp_{\t}$ on~$V_{\t}^{\vph_{\t}}$,
and a tuple~$\fo_{\x;\t}$ of orientations of fibers of~$V_{\t}^{\vph_{\t}}$ 
as above the CROrient~\ref{CROdegenC_prop} property and in~\eref{fobfxfn_e2}.
The latter in turn determine an orientation
$$\fo_{\t}\equiv \fo_{\fp_{\t}}\big(V_{\t},\vph_{\t};\fo_{\x;\t}\!\big)$$
of $D_{\t}$ for each $\t\!\in\!\De_{\R}^*$ as above Lemma~\ref{Vvphorient_lmm2}.
These orientations vary continuously with~$\t$ and extend to an orientation
\BE{rk1Clim_e4b}\fo_{\fp_0}'(\fo_{\x})\equiv \fo_{\fp_0}'\big(V_0,\vph_0;\fo_{\x}\big)\EE
of~$D_0$ as above the CROrient~\ref{CROdegenC_prop} property on page~\pageref{CROdegenC_prop}.
In an analogy with the limiting orientation of the CROrient~\ref{CROdegenC_prop} property,
we call \eref{rk1Clim_e4b} the \sf{limiting orientation} of~$D_0$
induced by~$\fo_{\x}$.

\begin{lmm}\label{VvphorientCdeg_lmm}
Suppose $\cC_0$ is a decorated marked symmetric surface as in~\eref{cCdegC_e} 
which contains precisely one conjugate pair $(\nod^+,\nod^-)$ of nodes and no other nodes
and carries precisely one real marked point~$x_r$ on each connected component~$S^1_r$ 
of~$\Si_0^{\si_0}$, 
$(V_0,\vph_0)$ is a real bundle pair over~$\cC_0$,
and $\fp$ is a relative $\Pin^{\pm}$-structure on the real vector bundle $V_0^{\vph_0}$
over $\Si_0^{\si_0}\!\subset\!\Si_0$.
Let $\fo_{\x}$ be a tuple of orientations of~$V_0^{\vph_0}|_{x_r}$ 
as in~\eref{fobfxfn_e} and 
$D_0$ be a real CR-operator on~$(V_0,\vph_0)$.
The intrinsic and limiting orientations, \eref{rk1Csplit_e4b} and~\eref{rk1Clim_e4b},
of~$D_0$ are the~same. 
\end{lmm}

\begin{proof} Let 
$\wt{D}_{0;\la}$ be the lift of a real CR-operator $D_{0;\la}$  on $\la(V_0,\vph_0)$
to~$\la(\wt{V}_0,\wt\vph_0)$.
We~take 
\BE{VvphorientCdeg_e5}D_{0;\pm}=D_0\!\oplus\!(2\!\pm\!1)D_{0;\la}, \qquad
\wt{D}_{0;\pm}=\wt{D}_0\!\oplus\!(2\!\pm\!1)\wt{D}_{0;\la}.\EE
These two decompositions and the one in~\eref{RBPstab_e} induce 
the exact triples of Fredholm operators given by the rows 
in the diagram of Figure~\ref{CROCdeg_fig}.
The exact triple~\eref{Cdegses_e} of Fredholm operators induces 
the exact triples of Fredholm operators given by the columns 
in this diagram.\\

\begin{figure}
$$\xymatrix{& 0\ar[d] & 0\ar[d] & 0\ar[d]  \\
0\ar[r] & D_0\ar[r]\ar[d] & D_{0;\pm}\ar[r]\ar[d] & (2\!\pm\!1)D_{0;\la}\ar[r]\ar[d] & 0 \\
0\ar[r] & \wt{D}_0\ar[r]\ar[d] & \wt{D}_{0;\pm} \ar[r]\ar[d]  
& (2\!\pm\!1)\wt{D}_{0;\la}\ar[r]\ar[d] & 0 \\
0\ar[r] & V_0|_{\nod^+}\ar[r]\ar[d] & V_{0;\pm}|_{\nod^+}\ar[r]\ar[d]
& (2\!\pm\!1)\la_{\C}(V_0|_{\nod^+})\ar[r]\ar[d] & 0 \\
& 0 & 0 & 0}$$ 
\caption{Commutative square of exact rows and columns of Fredholm operators 
for the proof of Lemma~\ref{VvphorientCdeg_lmm}}
\label{CROCdeg_fig}
\end{figure}

By definition, the columns in Figure~\ref{CROCdeg_fig} respect 
the complex orientations on the vector spaces in the bottom row,
the orientations
\BE{VvphorientCdeg_e7}
\wt\fo_{\fp_0}\big(\fo_{\x}\big)\equiv\fo_{\wt\fp_0}\big(\wt{V}_0,\wt\vph_0;\wt\fo_{\x}\big), 
\quad
\wt\fo_{\fp_0}\equiv\fo_{\wt\fp_0}(\wt{V}_{0;\pm},\wt\vph_{\pm}), \quad\hbox{and}\quad
\wt\fo_{0;\la}(\fo_{\x})\equiv\wt\fo\big(\la(\wt{V},\wt\vph);\wt\fo_{\x}\big)\EE
of $\wt{D}_0$, $\wt{D}_{0;\pm}$, and $\wt{D}_{0;\la}$
constructed above Lemma~\ref{Vvphorient_lmm2} and 
Propositions~\ref{genspinorient_prp} and~\ref{gencOpincg_prp}, respectively,
and the associated intrinsic orientations 
$\fo_{\fp_0}(\fo_{\x})$, 
\BE{VvphorientCdeg_e4a}\fo_{\fp_0}\equiv \fo_{\fp_0}\big((V_0)_{\pm},(\vph_0)_{\pm}\big),\EE
and $\fo_0(\fo_{\x})$ as in~\eref{rk1Csplit_e4b}, 
the CROrient~\ref{CROdegenC_prop} property with $\os_0\!=\!\fp_0$, and~\eref{rk1Csplit_e4}.
The bottom row respects the complex orientations of $V_0|_{\nod^+}$,  $V_{0;\pm}$, and
$\la_{\C}(V_0|_{\nod^+})$.
By definition, the middle row respects the three orientations in~\eref{VvphorientCdeg_e7};
the top row respects the limiting orientations $\fo_{\fp_0}'(\fo_{\x})$, 
\BE{VvphorientCdeg_e4b}\fo_{\fp_0}'\equiv \fo_{\fp_0}'\big((V_0)_{\pm},(\vph_0)_{\pm}\big),\EE
and $\fo_0'(\fo_{\x})$ as in~\eref{rk1Clim_e4b}, 
the CROrient~\ref{CROdegenC_prop}\ref{CdegenSpin_it} property with $\os_0\!=\!\fp_0$, 
and~\eref{rk1Clim_e4}.\\

By the CROrient~\ref{CROdegenC_prop} property and Proposition~\ref{genrk1Cdegen_prp},
$$\fo_{\fp_0}=\fo_{\fp_0}' \qquad\hbox{and}\qquad \fo_0(\fo_{\x})=\fo_0'(\fo_{\x}),$$
respectively.
Since the (real) dimension of $V_0|_{\nod^+}$ is even, Lemma~\ref{3by3_lmm} and 
the previous paragraph
thus imply that the orientations~\eref{rk1Csplit_e4b} and~\eref{rk1Clim_e4b}
the~same.
\end{proof}

We next suppose that $\cC_0$ is a decorated symmetric surface with 
one $H3$~node~$\nod$ as above~\eref{wtcCH3deg_e}
and in the top left diagram of Figure~\ref{genH3degen_fig} on page~\pageref{genH3degen_fig}
so that $x_r\!\in\!S^1_r$ for each $r\!\in\![N]$ different from the index~$r_{\bu}$
of the singular topological component of~$\Si_0^{\si_0}$ and $x_{r_{\bu}}\!=\!\nod$
(as in Proposition~\ref{genrk1H3degen_prp}, we allow a ``marked point" to be a node).
We also suppose that $(V_0,\vph_0)$ is a  real bundle pair over~$(\Si_0,\si_0)$.
Let $\wt\cC_0$, $(\wt\Si_0,\wt\si_0)$, $\wt\Si_0^b$, 
$\nod_1\!\in\!S^1_{\bu1}$, and $\nod_2\!\in\!S^1_{\bu2}$ be as in and below~\eref{wtcCH3deg_e} and
$(\wt{V}_0,\wt\vph_0)$ be the lift of $(V_0,\vph_0)$ to a real bundle pair over~$(\wt\Si_0,\wt\si_0)$.
A relative $\Pin^{\pm}$-structure~$\fp_0$ on the real vector bundle $V_0^{\vph_0}$
over $\Si_0^{\si_0}\!\subset\!\Si_0$ lifts to 
a relative $\Pin^{\pm}$-structure~$\wt\fp_0$ on the real vector bundle $\wt{V}_0^{\wt\vph_0}$
over $\wt\Si_0^{\wt\si_0}\!\subset\!\wt\Si_0$.
A tuple~$\fo_{\x}$ of orientations of the fibers of $V_0^{\vph_0}|_{x_r}$ 
as in~\eref{fobfxfn_e} lifts to a tuple~$\wt\fo_{\x}$ of orientations of
the corresponding fibers of~$\wt{V}_0^{\wt\vph_0}$ as in~\eref{wtfoH3gendfn_e}.\\

Let $D_0$ be a real CR-operator on $(V_0,\vph_0)$.
We denote its lift to a real CR-operator on $(\wt{V}_0,\wt\vph_0)$ by~$\wt{D}_0$. 
The orientation 
\BE{H3degenorient_e3b}
\wt\fo_{\fp_0}\big(\fo_{\x}\big)\equiv\fo_{\wt\fp_0}\big(\wt{V}_0,\wt\vph_0;\wt\fo_{\x}\big)\EE 
of $\wt{D}_0$ and the orientation $\fo_{\nod}\!\equiv\!\fo_{x_{r_{\bu}}}$
of $V_0|_{\nod}$ determine an orientation 
\BE{H3split_e4b}\fo_{\fp_0}\big(\fo_{\x}\big)\equiv \fo_{\fp_0}\big(V_0,\vph_0;\fo_{\x}\big)\EE
of $D_0$ via the isomorphism~\eref{Rdegses_e2}.
We call \eref{H3split_e4b} the \sf{intrinsic orientation} of~$D_0$
induced by~$\fp_0$ and~$\fo_{\x}$.\\

Suppose in addition that $(\cU,\wt\fc)$ is a flat family of deformations of~$\cC_0$
as in~\eref{cUsymmdfn_e},
$(V,\vph)$ is a real bundle pair over $(\cU,\wt\fc)$ extending~$(V_0,\vph_0)$,
$s_1^{\R},\ldots,s_N^{\R}$ are sections of $\cU^{\wt\fc}$ over~$\De_{\R}$ 
satisfying~\eref{H3sectcond_e}, 
and  \hbox{$\cD\!\equiv\!\{D_{\t}\}$} is a family of real CR-operators on
$(V_{\t},\vph_{\t})$ as in~\eref{DVvphdfn_e} extending~$D_0$.
The decorated structure on~$\cC_0$ determines an open subspace $\De_{\R}^+\!\subset\!\De_{\R}$
and a decorated structure on the fiber $(\Si_{\t},\si_{\t})$ of~$\pi$ 
for every $\t\!\in\!\De_{\R}$.
The relative $\Pin^{\pm}$-structure~$\fp_0$ on~$V_0^{\vph_0}$
and the tuple of orientations~$\fo_{\x}$ of fibers of~$V_0^{\vph_0}$ induce 
a relative $\Pin^{\pm}$-structure~$\fp_{\t}$ on~$V_{\t}^{\vph_{\t}}$
and a tuple~$\fo_{\x;\t}$ of orientations of fibers of~$V_{\t}^{\vph_{\t}}$ 
as above the CROrient~\ref{CROdegenH3_prop} property and~\eref{H3degenorient_e3},
respectively.
For each $\t\!\in\!\De_{\R}^*$, 
the decorated structure on~$(\Si_{\t},\si_{\t})$,
the relative $\Pin^{\pm}$-structure~$\fp_{\t}$, and
the tuple~$\fo_{\x;\t}$ determine an orientation
\BE{H3degenorient_e4b}\fo_{\t}\equiv \fo_{\fp_{\t}}\big(V_{\t},\vph_{\t};\fo_{\x;\t}\!\big)\EE
of $D_{\t}$ as above Lemma~\ref{Vvphorient_lmm2}.
We denote~by 
\BE{H3lim_e4b}\fo_{\fp_0}^+(\fo_{\x})\equiv \fo_{\fp_0}^+\big(V_0,\vph_0;\fo_{\x}\big)\EE
the orientation of~$D_0$ obtained as the continuous extension of 
the orientations~\eref{H3degenorient_e4b} with $\t\!\in\!\De_{\R}^+$.
We call~\eref{H3lim_e4b} the \sf{limiting orientation} of~$D_0$.\\

Let $W_1(V_0,\vph_0)_{r_{\bu}}$ be as in~\eref{W1Vvphrbudfn_e}.
Define
$$W_1(V_0,\vph_0)\equiv\big|\big\{S_r^1\!\in\!\pi_0(\Si_0^{\si_0})\!:
r\!\neq\!r_{\bu},\,w_1(V_0^{\vph_0})|_{S_r^1}\!\neq\!0\big\}\big| 
+\big|\{r\!\in\![2]\!:w_1(V_0^{\vph_0})|_{S_{\bu r}^1}\!\neq\!0\big\}\big|.$$

\begin{lmm}\label{VvphorientH3deg_lmm}
Suppose $\cC_0$ is a decorated marked symmetric surface as in~\eref{cCdegC_e} 
which contains precisely one~$H3$ node~$\nod$ and no other nodes
and carries precisely one real marked point~$x_r$ on each smooth connected component~$S^1_r$ 
of~$\Si_0^{\si_0}$ and a marked point $x_{r_{\bu}}\!=\!\nod$ on the nodal connected 
component of~$\Si_0^{\si_0}$.
Let $(V_0,\vph_0)$ be a rank~$n$ real bundle pair over~$(\Si_0,\si_0)$,
$\fp$ be a relative $\Pin^{\pm}$-structure on the real vector bundle $V_0^{\vph_0}$
over $\Si_0^{\si_0}\!\subset\!\Si_0$,
$\fo_{\x}$ be a tuple of orientations of $V_0^{\vph_0}|_{x_r}$
as in~\eref{fobfxfn_e}, and $D_0$ be a real CR-operator on $(V_0,\vph_0)$.
The intrinsic and limiting orientations, \eref{H3split_e4b} and~\eref{H3lim_e4b},
of~$D_0$ are the~same if and only~if
\begin{equation*}\begin{split}
&\big(\lr{w_1(V_0^{\vph_0}),[S^1_{\bu1}]_{\Z_2}}\!+\!1\big)
\lr{w_1(V_0^{\vph_0}),[S^1_{\bu2}]_{\Z_2}}
+n\big(W_1(V_0,\vph_0)\!-\!r_{\bu}\big)\!-\!W_1(V_0,\vph_0)_{r_{\bu}} \\
&\hspace{1.5in}=\begin{cases}
2\Z\in\Z_2,&\hbox{if}~\fp\!\in\!\cP^-_{\Si_0}(V_0^{\vph_0});\\
|\pi_0(\Si_0^{\si_0})|\!+\!W_1(V_0,\vph_0)\!+\!2\Z\in\Z_2,
&\hbox{if}~\fp\!\in\!\cP^+_{\Si_0}(V_0^{\vph_0}).
\end{cases}\end{split}\end{equation*}
\end{lmm}

\begin{proof}
Let $D_{0;\la}$, $\wt{D}_{0;\la}$,  $D_{0;\pm}$, and $\wt{D}_{0;\pm}$ be 
as in the proof of Lemma~\ref{VvphorientCdeg_lmm}.
We denote by $\fo_{\nod}^{\pm}$ the canonical orientation~on 
$$\la\big(V_{0;\pm}^{\vph_{0;\pm}}\big)\approx 
\la\big(V_0^{\vph_0}\big)^{\otimes(3\pm1)}\,.$$
The decompositions in~\eref{RBPstab_e} and~\eref{VvphorientCdeg_e5}
induce the exact triples of Fredholm operators given by the rows 
in the diagram of Figure~\ref{CROH3deg_fig}.
The exact triple~\eref{Rdegses_e} of Fredholm operators induces 
the exact triples of Fredholm operators given by the columns 
in this diagram.\\

\begin{figure}
$$\xymatrix{& 0\ar[d] & 0\ar[d] & 0\ar[d]  \\
0\ar[r] & D_0\ar[r]\ar[d] & D_{0;\pm}\ar[r]\ar[d] & (2\!\pm\!1)D_{0;\la}\ar[r]\ar[d] & 0 \\
0\ar[r] & \wt{D}_0\ar[r]\ar[d] & \wt{D}_{0;\pm} \ar[r]\ar[d]  
& (2\!\pm\!1)\wt{D}_{0;\la}\ar[r]\ar[d] & 0 \\
0\ar[r] & V_0^{\vph_0}|_{\nod}\ar[r]\ar[d] & V_{0;\pm}^{\vph_{0;\pm}}|_{\nod}\ar[r]\ar[d]
& (2\!\pm\!1)\la(V_0^{\vph_0}|_{\nod})\ar[r]\ar[d] & 0 \\
& 0 & 0 & 0}$$ 
\caption{Commutative square of exact rows and columns of Fredholm operators 
for the proof of Lemma~\ref{VvphorientH3deg_lmm}}
\label{CROH3deg_fig}
\end{figure}

By definition, the left and middle columns in this figure 
and the exact triple formed by each summand in the right column respect 
the orientations
$\fo_{\nod}$, $\fo_{\nod}^{\pm}$, and $\fo_{\nod}$ of the vector spaces in the bottom row,
the orientations~\eref{VvphorientCdeg_e7} of $\wt{D}_0$, $\wt{D}_{0;\pm}$, 
and $\wt{D}_{0;\la}$, respectively,
and the associated intrinsic orientations~$\fo_{\fp_0}(\fo_{\x})$, $\fo_{\fp_0}$, and $\fo_0(\fo_{\x})$
as in~\eref{H3split_e4b}, \eref{VvphorientCdeg_e4a}, and \eref{rk1H3split_e4}.
The bottom row respects the orientations $\fo_{\nod}$, $\fo_{\nod}^{\pm}$, and~$\fo_{\nod}$.
By definition, the middle row respects the orientations in~\eref{VvphorientCdeg_e7};
the top row respects the limiting orientations 
$\fo_{\fp_0}'(\fo_{\x})$, $\fo_{\fp_0}'$, and $\fo_0'(\fo_{\x})$
as in~\eref{H3lim_e4b}, \eref{VvphorientCdeg_e4b}, and~\eref{cOlimor_e4}.\\

By the CROrient~\ref{CROdegenH3_prop}\ref{HdegenSpin_it} property applied to $(V_0,\vph_0)_{\pm}$, 
$\fo_{\fp_0}\!=\!\fo_{\fp_0}'$ if and only~if
\BE{VvphorientH3deg_e5a}
(n\!+\!1)\big(|\pi_0(\Si_0^{\si_0})|\!-\!r_{\bu}\big)\in2\Z.\EE
By Proposition~\ref{genrk1H3degen_prp} applied to $\la(V_0,\vph_0)$, 
$\fo_0'(\fo_{\x})\!=\!\fo_0(\fo_{\x})$ if and only~if \eref{genrk1H3degen_e0} holds.
By Lemma~\ref{3by3_lmm}, the total number of rows and columns in Figure~\ref{CROH3deg_fig} which 
respect the orientations is of the same parity~as
\BE{VvphorientH3deg_e5} \big(\ind\,(2\!\pm\!1)D_{0;\la}\big)
\big(\!\dim V_0^{\vph_0}|_{\nod}\big)\cong
n\big(\big|\pi_0(\Si_0^{\si_0})|\!+\!W_1(V_0,\vph_0)\big)\mod2.\EE
Combining the statements in this paragraph and the previous one, we obtain
the claim in the $\Pin^-$-case.\\

Since the exact triple formed by the individual summands in the right column in 
Figure~\ref{CROH3deg_fig} respects the orientations 
$\fo_0(\fo_{\x})$, $\wt\fo_{0;\la}(\fo_{\x})$, and $\fo_{\nod}$,
Lemma~\ref{3by3_lmm} implies that the right column in the $\Pin^+$-case respects
the direct sum orientations if and only if the number 
$$\big(\ind\,D_{0;\la}\big)\big(\dim\,\la(V_0^{\vph_0}|_{\nod})\big)\cong
\big|\pi_0(\Si_0^{\si_0})\big|\!+\!W_1(V_0,\vph_0) \mod2$$
is even.
Combining this with the statements in the previous two paragraphs,
we obtained the claim in the $\Pin^+$-case.
\end{proof}

We are now in a position to obtain an analogue of the CROrient~\ref{CROSpinPinSES_prop} 
property for the orientations constructed above Lemma~\ref{Vvphorient_lmm2}.
We first prove Proposition~\ref{Vvphorient_prp4} below under the assumption 
that the complex vector bundle~$V'$ in~\eref{RBPsesdfn_e} is topologically trivial
as a consequence of the already established CROrient~\ref{CROSpinPinSES_prop}\ref{CROsesSpin_it} 
property on page~\pageref{CROSpinPinSES_prop}.
We then reduce the general case to this special case by bubbling off 
the non-trivial part of~$V'$ onto conjugate pairs of~spheres and 
using the already established CROrient~\ref{CROdegenC_prop}\ref{CdegenSpin_it} property
on page~\pageref{CROdegenC_prop} and Lemma~\ref{VvphorientCdeg_lmm}. 

\begin{prp}\label{Vvphorient_prp4}
Suppose $(\Si,\si)$ is a smooth decorated symmetric surface, 
$\ce$ and $\fo_{\x}''$ are as in Lemma~\ref{Vvphorient_lmm1b},
$\os'$ is a relative $\OSpin$-structure on the real vector bundle~$V'^{\vph'}$
over \hbox{$\Si^{\si}\!\subset\!\Si$},  
$\fo_{\x}'$ is the tuple orientations for $(V',\vph')$ as in~\eref{fobfxfn_e} 
obtained by restricting the orientation~$\fo'$ of~$V'^{\vph'}$ determined by~$\os'$,
and $\fp''$ is a relative $\Pin^{\pm}$-structure on~$V''^{\vph''}$.
The homotopy class of isomorphisms~\eref{CROses_e0} respects the orientations 
\BE{Vvphorient4_e3}
\fo_{\fp}(\fo_{\x})\equiv\fo_{\llrr{\os',\fp''}_{\ce_{\R}}}
\!\big(V,\vph;\llrr{\fo_{\x}',\fo_{\x}''}_{\ce_{\R}}\big), ~~
\fo_{\os'}\equiv\fo_{\os'}\big(V',\vph'\big), ~~
\fo_{\fp''}(\fo_{\x}'')\equiv\fo_{\fp''}\big(V'',\vph'';\fo_{\x}''\big)\EE
if and only if 
\BE{Vvphorient4_e4}
(\rk\,V')(\rk\,V''\!+\!1)\binom{|\pi_0(\Si^{\si})|}{2}\in2\Z.\EE
\end{prp}

\begin{proof} Let $\fp\!=\!\llrr{\os',\fp''}_{\ce_{\R}}$  
and $\fo_{\x}\!=\!\llrr{\fo_{\x}',\fo_{\x}''}_{\ce_{\R}}$.\\ 

(1) We first assume that the restriction of~$V'$ to each connected component of~$\Si$ 
is of degree~0.
Let $D'$, $D''$, and $D_{\la}''$ be real CR-operators on $(V',\vph')$, $(V'',\vph'')$,
and $\la(V'',\vph'')$, respectively, and set
$$D=D'\!\oplus\!D'', \qquad D_{\pm}''=D''\!\oplus\!(2\!\pm\!1)D_{\la}''.$$
By the assumption that $V'^{\vph'}$ is orientable and \cite[Lemma~3.2]{RBP},
we can identify the rank~1 real bundle pair $\la(V',\vph')$  over~$(\Si,\si)$
with the trivial rank~1 real bundle
$$(V_1,\vph_1)\equiv \big(\Si\!\times\!\C,\si\!\times\!\fc\big)$$
so that the induced identification of the real line bundles $\la(V'^{\vph'})$ 
and $\Si^{\si}\!\times\!\R$ respects the orientation~$\fo'$ and 
the standard orientation on~$\R$.
We denote by~$D_{\la}$ the real CR-operator on~$\la(V,\vph)$ corresponding to~$D_{\la}''$
under the induced identification~of 
$$\la(V,\vph)\approx \la(V',\vph')\!\otimes\!\la(V'',\vph'')$$
with $\la(V'',\vph'')$ and take
$$D_{\pm}=D\!\oplus\!(2\!\pm\!1)D_{\la}.$$
The decompositions of $D_{\pm}$ and $D_{\pm}''$ above induce 
the exact triples of Fredholm operators given by 
the middle and bottom rows in the first diagram of Figure~\ref{CROses_fig}.
The exact sequence~$\ce$ in~\eref{RBPsesdfn_e} induces the exact triples 
given by the left and middle columns in this diagram.\\

By definition, the middle and bottom rows in the diagram respect the orientations 
$\fo_{\fp}(\fo_{\x})$ of~$D$ and $\fo_{\fp''}(\fo_{\x})$ of~$D''$, 
$\fo_{\fp}(V_{\pm},\vph_{\pm})$ of~$D_{\pm}$ and 
$\fo_{\fp''}(V_{\pm}'',\vph_{\pm}'')$ of~$D_{\pm}$, and 
\BE{Vvphorient4_e6}\fo_{\la}(\fo_{\x})\equiv \fo\big(\la(V,\vph);\fo_{\x}\big)  \qquad\hbox{and}\qquad
\fo_{\la}(\fo_{\x}'')\equiv \fo\big(\la(V'',\vph'');\fo_{\x}''\big)\EE
of $D_{\la}$ and $D_{\la}''$, respectively.
By the assumption on the identification of $\la(V',\vph')$ with $(V_1,\vph_1)$,
the identification in the right column respects the orientations~\eref{Vvphorient4_e6}.
By the second identity in the RelSpinPin~\ref{RelSpinPinSES_prop}\ref{RelDSEquivSum_it} property 
on page~\pageref{RelDSEquivSum_it}
and the CROrient~\ref{CROSpinPinSES_prop}\ref{CROsesSpin_it} property, 
the middle column respects the orientations~$\fo_{\os'}$, $\fo_{\fp}(V_{\pm},\vph_{\pm})$, and 
$\fo_{\fp''}(V_{\pm}'',\vph_{\pm}'')$
if and only if \eref{Vvphorient4_e4} holds.
Along with Lemma~\ref{3by3_lmm}, the last three statements imply that 
the left column in the first diagram in Figure~\ref{CROses_fig}
respects the orientations~\eref{Vvphorient4_e3}
if and only if \eref{Vvphorient4_e4} holds,
provided the restriction of~$V'$ to each connected component of~$\Si$ 
is of degree~0.\\

\begin{figure}
$$\xymatrix{& 0\ar[d] & 0\ar[d]   \\
0\ar[r] & D'\ar[r]^{\id}\ar[d] & D'\ar[r]\ar[d] & 0\ar[d] \\
0\ar[r] & D\ar[r]\ar[d] & D_{\pm} \ar[r]\ar[d]  & (2\!\pm\!1)D_{\la}\ar[r]\ar[d]^= & 0 \\
0\ar[r] & D''\ar[r]\ar[d] & D_{\pm}'' \ar[r]\ar[d]& (2\!\pm\!1)D_{\la}''\ar[r]\ar[d] & 0 \\
& 0 & 0 & 0\\
& 0\ar[d] & 0\ar[d] & 0\ar[d]  \\
0\ar[r] & D_0'\ar[r]\ar[d] & D_{00}'\!\oplus\!D_0'^+\ar[r]\ar[d] 
& \wt{V}_0'|_{S^+}\ar[r]\ar[d] & 0 \\
0\ar[r] & D_0\ar[r]\ar[d] & D_{0;0}\!\oplus\!D_0^+ \ar[r]\ar[d]  
& \wt{V}_0|_{S^+}\ar[r]\ar[d] & 0 \\
0\ar[r] & D_0''\ar[r]\ar[d] & D''_{0;0}\!\oplus\!D_0''^+ 
\ar[r]\ar[d]& \wt{V}_0''|_{S^+}\ar[r]\ar[d] & 0 \\
& 0 & 0 & 0}$$ 
\caption{Commutative squares of exact rows and columns of Fredholm operators 
for the proof of Proposition~\ref{Vvphorient_prp4}}
\label{CROses_fig}
\end{figure}

(2) We now deduce the general case from the special case in~(1) and Lemma~\ref{VvphorientCdeg_lmm}.
Let $(\Si_0,\si_0)$ be a decorated marked symmetric surface consisting of~$(\Si,\si)$
with conjugate pairs of copies of~$S^2$ attached 
at one conjugate pair of points on each elemental component of~$(\Si,\si)$.
We denote~by 
$$\Si_{00},S^+,\Si_0^+\!\subset\!\Si_0$$ 
the copy of~$\Si$, the set of the nodal points carried 
by the distinguished half $\Si_{00}^b\!=\!\Si^b$,
and the union of the additional copies of~$S^2$ attached at~$S^+$.
Let $(\cU,\wt\fc)$ be a family of deformations of $\cC_0\!\equiv\!(\Si_0,\si_0)$ 
over $\De\!\subset\!\C^{2|S^+|}$ 
so that the fiber~$\cC_{\t_1}$ of $(\cU,\wt\fc)$ over some $\t_1\!\in\!\De_{\R}$ is~$(\Si,\si)$ 
and $s_1^{\R},\ldots,s_N^{\R}$ be sections of~$\cU^{\wt\fc}$ over $\De_{\R}$ with
$s_r^{\R}(\t_1)\!=\!x_r$ for each $S^1_r\!\in\!\pi_0(\Si^{\si})$.
Let
\BE{Vvphorient4_e9}0\lra (\wt{V}',\wt\vph')\lra (\wt{V},\wt\vph)\lra (\wt{V}'',\wt\vph'')\lra0\EE
be a short exact sequence of real bundle pairs over~$(\cU,\wt\fc)$  restricting to~$\ce$ over~$\t_1$
so that the restriction of $\wt{V}'$ to each topological component of~$\Si_{00}$ is of degree~0;
such a sequence exists because the restriction of $V'$ to each elemental component 
of~$(\Si,\si)$ is of even degree.
We denote the exact sequence~\eref{Vvphorient4_e9} by~$\wt\ce$
and the restrictions of $(\wt{V}',\wt\vph')$, $(\wt{V}',\wt\vph')$, and
$(\wt{V}'',\wt\vph'')$ to~$\Si_{00}$ by  
$(\wt{V}'_{00},\wt\vph'_{00})$, $(\wt{V}_{00},\wt\vph_{00})$, 
and $(\wt{V}''_{00},\wt\vph''_{00})$, respectively.\\

For each $\t\!\in\!\De_{\R}$,
the tuples $\fo_{\x}'$, $\fo_{\x}$, and~$\fo_{\x}''$ of orientations of fibers of 
$V'^{\vph'}$, $V^{\vph}$, and~$V''^{\vph''}$ induce tuples 
$\fo_{\x;\t}'$, $\fo_{\x;\t}$, and~$\fo_{\x;\t}''$ of orientations of fibers of 
$\wt{V}_{\t}'^{\wt\vph_{\t}'}$, $\wt{V}_{\t}^{\wt\vph_{\t}}$, and~$\wt{V}_{\t}''^{\wt\vph_{\t}''}$
via the vector bundles 
$s_r^{\R*}\wt{V}'^{\wt\vph'}$, $s_r^{\R*}\wt{V}^{\wt\vph}$, and~$s_r^{\R*}\wt{V}''^{\wt\vph''}$.
Let $\wt\os'$ and $\wt\fp''$ be relative $\OSpin$-structures on the real vector bundles
$\wt{V}'^{\wt\vph'}$ and $\wt{V}''^{\wt\vph''}$, respectively, 
over~$\cU^{\wt\fc}\!\subset\!\cU$ restricting to~$\os'$ and~$\fp''$ on~$\Si_{\t_1}$
so that $w_2(\wt\os')$ and $w_2(\wt\fp'')$
vanish on each added copy of~$S^2$.
We~define
\begin{gather}\notag
\wt\fp=\bllrr{\wt\os',\wt\fp''}_{\wt\ce_{\R}},\qquad
\os_{00}'=\wt\os'|_{\Si_{00}},~\fp_{00}=\wt\fp|_{\Si_{00}},~
\fp_{00}''=\wt\fp''|_{\Si_{00}}, \\
\os_{\t}'=\wt\os'|_{\Si_{\t}},~\fp_{\t}=\wt\fp|_{\Si_{\t}},~
\fp_{\t}''=\wt\fp''|_{\Si_{\t}}~~\forall\,\t\!\in\!\De_{\R}\,,\\
\label{Vvphorient4_e6b}
\fo_{\os';\t}=\fo_{\os_{\t}'}\big(\wt{V}_{\t}',\wt\vph_{\t}'\big),~
\fo_{\fp;\t}(\fo_{\x})=\fo_{\fp_{\t}}\!\big(\wt{V}_{\t},\wt\vph_{\t};\fo_{\x;\t}\big), ~
\fo_{\fp_{\t}''}(\fo''_{\x})=
\fo_{\fp_{\t}''}\big(\wt{V}_{\t}'',\wt\vph_{\t}'';\fo_{\x;\t}''\big)
~~\forall\,\t\!\in\!\De_{\R}^*.
\end{gather}

\vspace{.15in}

Let \hbox{$\cD'\!\equiv\!\{D_{\t}'\}$} and \hbox{$\cD''\!\equiv\!\{D_{\t}''\}$}
be families of real CR-operators on $(V_{\t}',\vph_{\t}')$ and $(V_{\t}'',\vph_{\t}'')$,
respectively, as in~\eref{DVvphdfn_e}, and 
\BE{Vvphorient4_e8a}D_{\t}=D_{\t}'\!\oplus\!D_{\t}''~~\forall\t\!\in\!\De_{\R}\,.\EE
We denote by $D_{00}'$, $D_{00}$, and $D_{00}''$ the real CR-operators on 
$(\wt{V}'_{00},\wt\vph'_{00})$, $(\wt{V}'_{00},\wt\vph'_{00})$, and $(\wt{V}''_{00},\wt\vph''_{00})$,
respectively, induced by $D_0'$, $D_0$, and~$D_0''$ and by 
$D_{00}^+$, $D_{00}^+$, and $D_{00}^+$ the real CR-operators on 
$\wt{V}'|_{\Si_0^+}$, $\wt{V}|_{\Si_0^+}$, and $\wt{V}'|_{\Si_0^+}$, 
respectively, induced by $D_0'$, $D_0$, and~$D_0''$.
The orientations~\eref{Vvphorient4_e6b} of $D_{\t}'$, $D_{\t}$, and $D_{\t}''$
extend continuously to orientations
\BE{Vvphorient4_e8}\fo_{\os';0}'\equiv\fo_{\os_0'}\big(\wt{V}'_0,\wt\vph'_0\big),~~
\fo_{\fp;0}'(\fo_{\x})\equiv\fo_{\fp_0}\big(\wt{V}_0,\wt\vph_0;\fo_{\x;0}\big),~~\hbox{and}~~
\fo_{\fp'';0}'(\fo_{\x}'')\equiv\fo_{\fp_0''}\big(\wt{V}_0'',\wt\vph_0'';\fo_{\x;0}''\big)\EE
of $D_0'$, $D_0$, and $D_0''$, respectively;
these are the analogues of the limiting orientations of 
the CROrient~\ref{CROdegenC_prop}\ref{CdegenSpin_it} property.
We denote~by 
\BE{Vvphorient4_e8b}\fo_{\os';0}\equiv\fo_{\os_0'}\big(\wt{V}'_0,\wt\vph'_0\big),~~
\fo_{\fp;0}(\fo_{\x})\equiv\fo_{\fp_0}\big(\wt{V}_0,\wt\vph_0;\fo_{\x;0}\big),~~\hbox{and}~~
\fo_{\fp'';0}(\fo_{\x}'')\equiv\fo_{\fp_0''}\big(\wt{V}_0'',\wt\vph_0'';\fo_{\x;0}''\big)\EE
the corresponding analogues of the intrinsic orientation of 
the CROrient~\ref{CROdegenC_prop}\ref{CdegenSpin_it} property.
By Lemma~\ref{VvphorientCdeg_lmm} applied $|S^+|$ times, 
these orientations are the same as the limiting orientations~\eref{Vvphorient4_e8}.\\

The exact triple~\eref{Cdegses_e}, the decomposition~\eref{D1D2split_e}, and
the isomorphism~\eref{CROdoubl_e0} induce 
the exact triples of Fredholm operators given by 
the rows in the second diagram of Figure~\ref{CROses_fig}.
The restriction of the exact sequence~\eref{Vvphorient4_e9} to~$\Si_0$
induces the exact triples given by the columns in this diagram.
Since the decompositions~\eref{Vvphorient4_e8a} induce an isomorphism
$$\la(\cD)\approx\la(\cD')\!\otimes\!\la(\cD'')$$
of line bundle over~$\De_{\R}$
and the orientations~\eref{Vvphorient4_e8} restrict to the orientations~\eref{Vvphorient4_e3} 
for $\t=\!\t_1$, it is sufficient to show that the left column
respects the limiting orientations~\eref{Vvphorient4_e8b}
if and only if  \eref{Vvphorient4_e4} holds.\\

The right column and the exact triple formed by the second summands
in the middle column respect the complex orientations of all terms.
Since  the restriction of $\wt{V}_{00}'$ to each topological component of~$\Si_{00}$ is of degree~0,
(1)~above implies that the exact triple formed by the first summands
in the middle column respects the orientations 
$$\fo_{\os';00}\equiv\fo_{\os_{00}'}\big(\wt{V}'_{00},\wt\vph'_0\big),~~
\fo_{\fp;00}(\fo_{\x})\equiv\fo_{\fp_{00}}\big(\wt{V}_0,\wt\vph_{00};\fo_{\x;0}\big),~~\hbox{and}~~
\fo_{\fp'';00}(\fo_{\x}'')\equiv\fo_{\fp_{00}''}\big(\wt{V}_{00}'',\wt\vph_{00}'';\fo_{\x;0}''\big)$$
of $D_{00}'$, $D_{00}$, and~$D_{00}''$, respectively, if and only if  \eref{Vvphorient4_e4} holds.
Since the (real) index of $D_0'^+$ is even, 
Lemma~\ref{3by3_lmm} then implies that the middle column respect the direct sum orientations
\BE{Vvphorient4_e10}
\wt\fo_{\os';0}\equiv\wt\fo_{\os_0'}\big(\wt{V}'_0,\wt\vph'_0\big),~~
\wt\fo_{\fp;0}(\fo_{\x})\equiv\wt\fo_{\fp_0}\big(\wt{V}_0,\wt\vph_0;\fo_{\x;0}\big),~~\hbox{and}~~
\wt\fo_{\fp'';0}(\fo_{\x}'')\equiv\wt\fo_{\fp_0''}\big(\wt{V}_0'',\wt\vph_0'';\fo_{\x;0}''\big)\EE
if and only if  \eref{Vvphorient4_e4} holds.\\

The orientations 
\BE{Vvphorient4_e12}\fo_{\os';0}^{\C}\equiv\fo_{\os_0'}^{\C}\big(\wt{V}'_0,\wt\vph'_0\big),~~
\fo_{\fp;0}^{\C}(\fo_{\x})\equiv\fo_{\fp_0}\big(\wt{V}_0,\wt\vph_0;\fo_{\x;0}\big),~~\hbox{and}~~
\fo_{\fp'';0}^{\C}(\fo_{\x}'')\equiv\fo_{\fp_0''}\big(\wt{V}_0'',\wt\vph_0'';\fo_{\x;0}''\big)\EE
of $D_0'$, $D_0$, and $D_0''$, respectively, determined by the orientations~\eref{Vvphorient4_e10} 
and the complex orientations of the vector spaces in the right column via
the rows in the second diagram of Figure~\ref{CROses_fig} are the analogues
of the $\C$-split orientations of Corollary~\ref{CROdegenC_crl}.
By the proof of Lemma~\ref{Vvphorient_lmm3}\ref{CRODisjUn_it}, 
these orientations are the same as the orientations~\eref{Vvphorient4_e8b}.
Along with the sentence after~\eref{Vvphorient4_e8b}, 
this implies that the rows in this diagram respect the orientations~\eref{Vvphorient4_e8},
the orientations~\eref{Vvphorient4_e10}, and the complex orientations of 
$\wt{V}_0'|_{S^+}$, $\wt{V}_0|_{S^+}$, and~$\wt{V}_0''|_{S^+}$.
Since the right column respects these complex orientations,
the middle column respects the orientations~\eref{Vvphorient4_e10} 
if and only if  \eref{Vvphorient4_e4} holds, 
and the (real) dimension of~$\wt{V}'|_{S^+}$ is even,
Lemma~\ref{3by3_lmm} implies that the left column respects  
the limiting orientations~\eref{Vvphorient4_e8} if and only if 
\eref{Vvphorient4_e4} holds.
\end{proof}

\subsection{Properties of twisted orientations}
\label{twisteddetLB_subs}

Suppose $\cC$ is a smooth decorated marked symmetric surface as in~\eref{cCsymdfn_e},
$(V,\vph)$ is a real bundle pair over~$\cC$,
$\fp$ is a relative $\Pin^{\pm}$-structure on the real vector bundle~$V^{\vph}$
over $\Si^{\si}\!\subset\!\Si$, and $\fo_{\x}$ is a tuple of orientations of
fibers of~$V^{\vph}$ as in~\eref{foxcCdfn_e}.
For any real CR-operator~$D$ on the real bundle pair~$(V,\vph)$, 
we denote by $\fo_{\cC;\fp}(V,\vph;\fo_{\x})$ the orientation of
the twisted determinant~\eref{wtlacCDdfn_e} of~$D$ induced by the orientations
$\fo_{\cC}(V^{\vph},\fo_{\x})$ of $\la_{\cC}^{\R}(V,\vph)$
and $\fo_{\fp}(V,\vph;\fo_{\x})$ of~$D$ defined in 
Section~\ref{TargetOrient_subs} and~\ref{detOrient_subs}, respectively.\\

Suppose $\cC$ is a marked symmetric surface with the underlying symmetric surface~$(S^2,\tau)$ 
and $k\!=\!2$ real marked points, 
$(V,\vph)$ is a rank~1 odd-degree real bundle pair over $(S^2,\tau)$,
and $\fo_{x_1}$ is an orientation of~$V_{x_1}^{\vph}$. 
We then denote~by 
\BE{Pin0Orientcrl_e}\fo_{\cC;0}^{\pm}\big(V,\vph;\fo_{x_1}\big)\equiv 
\fo_{\cC;\io_{S^2}(\fp_0^{\pm}(V^{\vph}))}\big(V,\vph;\fo_{x_1}\big)\EE
the orientation on $\wt\la_{\cC}(D)$ determined by the image~\eref{PinNormDfn_e}
of the $\Pin^{\pm}$-structure $\fp_0^{\pm}(V^{\vph})$ of
Examples~\ref{Pin1pMB_eg}, \ref{Pin1mMB_eg}, and~\ref{SpinDfn1to3_eg}
under the second map in~\eref{vsSpinPin_e0} with \hbox{$X\!=\!S^2$}.\\

The first three statements of Corollary~\ref{Vvphorient_crl2} below extend
the three statements of the CROrient~\ref{CROp_prop} property on page~\pageref{CROp_prop}
to real bundle pairs~$(V,\vph)$ over smooth marked symmetric surfaces~$\cC$
without the balancing restriction.
The last statement of this corollary, 
the first and last statements of Corollary~\ref{Vvphorient_crl3}, 
and Corollary~\ref{Vvphorient_crl4}
extend the CROrient~\ref{CROSpinPinStr_prop}\ref{CROPinStr_it}, \ref{CROPin2SpinRed_prop}, 
\ref{CRODisjUn_prop}\ref{pDisjUn_it}, and \ref{CROSpinPinSES_prop}\ref{CROsesPin_it} properties
in the same way.
Corollaries~\ref{VvphorientCdeg_lcrl} and~\ref{VvphorientH3deg_lcrl} extend 
the CROrient~\ref{CROdegenC_prop}\ref{CdegenPin_it} and~\ref{CROdegenH3_prop}\ref{HdegenPin_it} 
properties of Section~\ref{OrientPrp_subs2} 
to real bundle pairs~$(V_0,\vph_0)$ over marked symmetric surfaces~$\cC_0$
without the balancing restriction.
According to the middle statement of Corollary~\ref{Vvphorient_crl3},
the orientations $\fo_{\cC;\fp}(V,\vph;\fo_{\x})$ satisfy the conclusion 
of the CROrient~\ref{CRONormal_prop}\ref{CROnormPin_it} property on 
page~\pageref{CRONormal_prop}.

\begin{crl}\label{Vvphorient_crl2}
Suppose $\cC$ is a smooth decorated marked symmetric surface, 
$(V,\vph)$ is a rank~$n$ real bundle pair over~$\cC$,
$\fp$ is a relative $\Pin^{\pm}$-structure on the real vector bundle $V^{\vph}$ 
over $\Si^{\si}\!\subset\!\Si$,
and $D$ is a real CR-operator on~$(V,\vph)$.
Let $\fo_{\x}$ be a tuple of orientations of~$V^{\vph}$ at points $x_{r}^*$
in $S^1_r\!\subset\!\Si^{\si}$ as in~\eref{foxcCdfn_e}.
\begin{enumerate}[label=(\arabic*),leftmargin=*]

\item\label{rk1flipcrl_it} The orientation $\fo_{\cC;\fp}(V,\vph;\fo_{\x})$
on~$\wt\la_{\cC}(D)$ 
does not depend on the choice of half-surface $\Si_*^b$ of 
an elemental component~$\Si_*$ of $(\Si,\si)$ if and only~if 
\begin{equation*}\begin{split}
\sum_{S_r^1\in\pi_0(\Si_*^{\si})}\!\!\!
\bigg(\!\!\blr{w_1(V^{\vph}),[S^1_r]_{\Z_2}}\binom{k_r\!-\!1}{1}
\!+n\binom{k_r\!-\!1}{2}\!\!\!\bigg)
=\vp_{\fp}(\Si_*)\in \Z_2\,.
\end{split}\end{equation*}

\item\label{kr1intercrl_it} The interchange in the ordering of two consecutive components~$S_r^1$ 
and~$S_{r+1}^1$  of~$\Si^{\si}$ preserves the orientation $\fo_{\cC;\fp}(V,\vph;\fo_{\x})$
if and only~if 
$$n\big(k_rk_{r+1}\!+\!1\big)\!+\!2\Z\neq
\big(\blr{w_1(V^{\vph}),[S^1_r]_{\Z_2}}\!+\!1\big)
\big(\blr{w_1(V^{\vph}),[S^1_{r+1}]_{\Z_2}}\!+\!1\big)\in \Z_2.$$

\item\label{kr1inter2crl_it} 
The interchange of two real marked points $x_{j_i^r(\cC)}$ and $x_{j_{i'}^r(\cC)}$
on the same connected component $S_r^1$ of~$\Si^{\si}$
with \hbox{$2\!\le\!i,i'\!\le\!k_r$} preserves $\fo_{\cC;\fp}(V^{\vph},\fo_{\x})$.
The combination of the interchange of the real points $x_{j_1^r(\cC)}$ and $x_{j_i^r(\cC)}$
with $2\!\le\!i\!\le\!k_r$ and the replacement of the component 
\hbox{$\fo_{x_r^*}\!\equiv\!\fo_{x_{j_1^r(\cC)}^*}$}
in~\eref{foxcCdfn_e} by~$\fo_{x_{j_i^r(\cC)}^*}$ preserves $\fo_{\cC;\fp}(V,\vph;\fo_{\x})$ 
if and only if \eref{Vvphinter2_e} holds.

\item\label{rk1horientcrl_it} 
The reversal of the component orientation~$\fo_{x_r^*}$ in~\eref{foxcCdfn_e}
preserves $\fo_{\cC;\fp}(V,\vph;\fo_{\x})$ if and only if 
$$k_r\!+\!2\Z=\blr{w_1(V^{\vph}),[S^1_r]_{\Z_2}}\!+\!1 \in2\Z.$$

\item\label{rk1etacrl_it} If in addition $\eta\!\in\!H^2(\Si,\Si^{\si};\Z_2)$,
then the orientations $\fo_{\cC;\fp}(V,\vph;\fo_{\x})$ and $\fo_{\cC;\eta\cdot\fp}(V,\vph;\fo_{\x})$ 
are the same if and only~if $\lr{\eta,[\Si^b]_{\Z_2}}\!=\!0$.

\end{enumerate}
\end{crl}

\begin{proof}
The five statements of this corollary follow from 
the four statements of Lemma~\ref{Vvphorient_lmm} and 
the four statements of Lemma~\ref{Vvphorient_lmm2}.
\end{proof}

\begin{crl}\label{Vvphorient_crl3}
Suppose $\cC$, $(V,\vph)$, $\fp$, and $D$ are as in Corollary~\ref{Vvphorient_crl2}.
\begin{enumerate}[label=(\arabic*),leftmargin=*]

\item\label{Pin2SpinRedcrl_it} If $\fo$ is an orientation on~$V^{\vph}$ and
$\fo_{\x}$ is a tuple of orientations of~$V^{\vph}$ at points $x_{r}^*$
in $S^1_r\!\subset\!\Si^{\si}$ as in~\eref{foxcCdfn_e} obtained by restricting~$\fo$,
then the orientation $\fo_{\cC;\fp}(V,\vph;\fo_{\x})$ on~$\wt\la_{\cC}(D)$ 
corresponds to the homotopy class of isomorphisms of $\la(D)$ and $\la_{\cC}^{\R}(V,\vph)$
determined by the orientations $\fo_{\fR_{\fo}^{\pm}(\fp)}(V,\vph)$
and $\la_{\cC}^{\R}(\fo)$ 
if and only if \eref{CROPin2SpinRed_e} holds.

\item\label{CROnormcrl_it} 
If $(\Si,\si)\!=\!(S^2,\tau)$, $k\!=\!2$, $\rk\,V\!=\!1$, $\deg V\!=\!1$,
and $\fo_{x_1}$ is an orientation of~$V_{x_1}^{\vph}$,  then
\eref{Pin0Orientcrl_e} is the orientation of $\wt\la_{\cC}(D)$
corresponding to the homotopy class of the isomorphism~\eref{cOevdfn_e0} with $a\!=\!1$.

\item\label{CRODisjUncrl_it} If $\fo_{\x}$ is a tuple of orientations of~$V^{\vph}$
as in~\eref{foxcCdfn_e} and $\cC\!=\!\cC_1\!\sqcup\!\cC_2$ is a decomposition of~$\cC$ into 
decorated marked symmetric surfaces as above~\eref{V1V2dfn_e4}
of genera~$g_1$ and~$g_2$, respectively, then
$$\fo_{\cC;\fp}(V,\vph)=
\fo_{\cC_1;\fp|_{\cC_1}}\!\big(V|_{\cC_1},\vph|_{\cC_1};\fo_{\x}|_{\cC_1}\!\big)
\!\sqcup 
\fo_{\cC_2;\fp|_{\cC_2}}\!\big(V|_{\cC_2},\vph|_{\cC_2};\fo_{\x}|_{\cC_2}\!\big)$$
if and only~if \eref{CRODisjUn_e0c} with $V_1\!\equiv\!V|_{\cC_1}$ and $V_2\!\equiv\!V|_{\cC_2}$
holds.

\end{enumerate}
\end{crl}

\begin{proof} \ref{Pin2SpinRedcrl_it}
By the assumption on $\fo_{\x}$, the orientations $\la_{\cC}^{\R}(\fo)$ and 
$\fo_{\cC}(V^{\vph},\fo_{\x})$ of  $\la_{\cC}^{\R}(V,\vph)$ are the same.
Along with the definition of the orientation $\fo_{\cC;\fp}(V,\vph;\fo_{\x})$ of~$\wt\la_{\cC}(D)$
above Corollary~\ref{Vvphorient_crl2} and Lemma~\ref{Vvphorient_lmm3}\ref{Pin2SpinRed_it},
this implies the claim.\\

\ref{CROnormcrl_it} By Lemma~\ref{Vvphorient_lmm3}\ref{CROnorm_it},
the orientation~\eref{Pin0Orient_e} of~$D$ is the same as 
the intrinsic orientation $\fo(V,\vph;\fo_{x_1})$ defined above Proposition~\ref{cOpincg_prp}.
By Proposition~\ref{cOpincg_prp}\ref{cOpincg_it2},
$$\fo\big(V,\vph;\fo_{x_1}\big)\equiv\fo_{2,0}\big(V,\vph;\fo_{x_1}\big)
=\fo_{2,0}\big(\fo_{x_1}\big)\,.$$
The claim thus follows from the definitions of the twisted orientation~\eref{Pin0Orientcrl_e}
and the evaluation orientation~$\fo_{2,0}(\fo_{\x})$.\\

\ref{CRODisjUncrl_it} Since the decomposition $\cC\!=\!\cC_1\!\sqcup\!\cC_2$ 
respects the orderings  of the connected components of the fixed loci, 
the isomorphism~\eref{V1V2dfn_e4} respects the orientations
$$\fo_{\cC}\big(V^{\vph},\fo_{\x}\big), \quad
\fo_{\cC_1}\!\big(V_1^{\vph_1},\fo_{\x}|_{\cC_1}\big),\quad\hbox{and}\quad 
\fo_{\cC_2}\big(V_2^{\vph_2},\fo_{\x}|_{\cC_2}\big).$$
The claim thus follows from Lemma~\ref{Vvphorient_lmm3}\ref{CRODisjUn_it}.
\end{proof}

\begin{crl}\label{Vvphorient_crl4}
Suppose $\cC$, $\ce$, and $\fo_{\x}''$ are as in Lemma~\ref{Vvphorient_lmm1b},
$\os'$ is a relative $\OSpin$-structure on the real vector bundle~$V'^{\vph'}$
over $\Si^{\si}\!\subset\!\Si$,  and
$\fo_{\x}'$ is the tuple orientations for $(V',\vph')$ as in~\eref{foxcCdfn_e2} 
obtained by restricting the orientation~$\fo'$ of~$V'^{\vph'}$ determined by~$\os'$.
If $\fp''\!\in\!\cP_{\Si}(V''^{\vph''})$, then
$$\fo_{\cC;\llrr{\os',\fp''}_{\ce_{\R}}}\!\big(V,\vph;\llrr{\fo_{\x}',\fo_{\x}''}_{\ce_{\R}}\big)
=\big(\fo_{\os'}(V',\vph')\la_{\cC}^{\R}(\os')\!\big)_{\ce} 
\fo_{\cC;\fp''}\big(V'',\vph'';\fo_{\x}''\big)$$
if and only~if \eref{CROsesPin_e0} holds.
\end{crl}

\begin{proof}
This claim follows immediately from Proposition~\ref{Vvphorient_prp4} 
and Lemma~\ref{Vvphorient_lmm1b}.
\end{proof}

It remains to describe the behavior of the orientations $\fo_{\cC;\fp}(V,\vph;\fo_{\x})$  
under the flat degenerations of the CROrient~\ref{CROdegenC_prop} and~\ref{CROdegenH3_prop} 
properties in Section~\ref{OrientPrp_subs2}.
We first suppose that $\cC_0$ is a decorated symmetric surface with 
one conjugate pair of nodes~$\nod^{\pm}$ as in~\eref{cCdegC_e},
$(V_0,\vph_0)$ is a real bundle pair over~$(\Si_0,\si_0)$,
$\fp_0$ is a relative $\Pin^{\pm}$-structure on the real vector bundle $V_0^{\vph_0}$
over $\Si_0^{\si_0}\!\subset\!\Si_0$, and 
$\fo_{\x}$ is a tuple of orientations of the fibers of $V_0^{\vph_0}$ at 
points $x_r^*\!\in\!S^1_r$ as in~\eref{foxcCdfn_e}.
Let $D_0$ be a real CR-operator on~$(V_0,\vph_0)$.
We define the associated \sf{intrinsic} and \sf{limiting} orientations,
\BE{genintvslimcrl_e} 
\fo_{\cC_0;\fp_0}\big(\fo_{\x}\big)\equiv\fo_{\cC_0;\fp_0}\big(V_0,\vph_0;\fo_{\x}\big)
\qquad\hbox{and}\qquad 
\fo_{\cC_0;\fp_0}'\big(\fo_{\x}\big)\equiv\fo_{\cC_0;\fp_0}'\big(V_0,\vph_0;\fo_{\x}\big),\EE
of the twisted determinant
\BE{wtlacC0D0dfn_e}\wt\la_{\cC_0}(D_0)\equiv \la_{\cC_0}^{\R}(V_0,\vph_0)^*\!\otimes\!\la(D_0)\EE
of $D_0$ to be the orientations induced by the intrinsic and limiting orientations
of $\la_{\cC_0}^{\R}(V_0,\vph_0)$ defined after Lemma~\ref{Vvphorient_lmm1b}
and the intrinsic and limiting orientations of~$D_0$ 
in~\eref{rk1Csplit_e4b} and~\eref{rk1Clim_e4b}, respectively.
The observation below~\eref{lacC0Rdfn2_e} and Lemma~\ref{VvphorientCdeg_lmm}
immediately give the following conclusion.

\begin{crl}\label{VvphorientCdeg_lcrl}
The intrinsic and limiting orientations~\eref{genintvslimcrl_e} on 
$\wt\la_{\cC_0}(D_0)$ are the~same. 
\end{crl}

Suppose finally that $\cC_0$ is a decorated marked symmetric surface as in~\eref{cCdegC_e} 
which contains precisely one~$H3$ node~$\nod$ and no other nodes,
$(V_0,\vph_0)$ and $\fp_0$ are as above Corollary~\ref{VvphorientCdeg_lcrl},
$\fo_{\x}$ is a tuple of orientations of the fibers of $V_0^{\vph_0}$ at 
points $x_r^*\!\in\!S^1_r$ for $r\!\neq\!r_{\bu}$ and $x_{r_{\bu}}^*\!=\!\nod$
as above~\eref{H3orientlift_e}. 
Let $D_0$ be a real CR-operator on~$(V_0,\vph_0)$.
We define the associated \sf{intrinsic} and \sf{limiting} orientations,
\BE{genintvslimH3_e}\begin{split}
\fo_{\cC_0;\fp_0;1}\big(\fo_{\x}\big)\equiv\fo_{\cC_0;\fp_0;1}\big(V_0,\vph_0;\fo_{\x}\big), &\qquad
\fo_{\cC_0;\fp_0;2}\big(\fo_{\x}\big)\equiv\fo_{\cC_0;\fp_0;2}\big(V_0,\vph_0;\fo_{\x}\big),\\
\hbox{and}\qquad 
\fo_{\cC_0;\fp_0}^+\big(\fo_{\x}\big)&\equiv\fo_{\cC_0;\fp_0}^+\big(V_0,\vph_0;\fo_{\x}\big),
\end{split}\EE
of the twisted determinant~\eref{wtlacC0D0dfn_e}
of $D_0$ to be the orientations induced by the intrinsic and limiting orientations
of $\la_{\cC_0}^{\R}(V_0,\vph_0)$ as below~\eref{H3orientisom_e} and in~\eref{cCOlimor_e4}
and the intrinsic and limiting orientations of~$D_0$ 
in~\eref{H3split_e4b} and~\eref{H3lim_e4b}.
Let $\ep_{\R}(\cC_0;V_0^{\vph_0})$, $\de_{r^*}(V_0,\vph_0)$, and
$\ep_{r^*}(\cC_0;V_0^{\vph_0})$ be as above Lemma~\ref{VvphorientH3deg_lmm1}
and $W_1(V_0,\vph_0)_{r_{\bu}}$ and $W_1(V_0,\vph_0)$ be as above 
Lemma~\ref{VvphorientH3deg_lmm}.
These two lemmas immediately give the following conclusion.

\begin{crl}\label{VvphorientH3deg_lcrl}
Let $r^*\!=\!1,2$.
The intrinsic and limiting orientations, 
$\fo_{\cC_0;\fp_0;r^*}(\fo_{\x})$ and $\fo_{\cC_0;\fp_0}^+(\fo_{\x})$
in~\eref{genintvslimcrl_e} on  $\wt\la_{\cC_0}(D_0)$ are the~same if and only~if
\begin{equation*}\begin{split}
&(n\!+\!1)\bigg(\!\!\big(k_{r_{\bu}}\!-\!1\!\big)\de_{\R}(\cC_0)\!+\!\de_{r^*}(\cC_0)
+\!\!\!\!\!\!\!\sum_{\begin{subarray}{c}S_r^1\in\pi_0(\Si_0^{\si_0})\\ 
r>r_{\bu}\end{subarray}}\!\!\!\!\!\!\!\!\!(k_r\!-\!1)\!\!\bigg)
+\sum_{\begin{subarray}{c}S_r^1\in\pi_0(\Si_0^{\si_0})\\ 
r>r_{\bu}\end{subarray}}\!\!\!\!\!\!\!\!\!(k_r\!-\!1)
-\!W_1(V_0,\vph_0)_{r_{\bu}}\\
&\quad
+\big(k_{r_{\bu}}\!-\!1\!+\!\blr{w_1(V_0^{\vph_0}),[S^1_{\bu1}]_{\Z_2}\!+\![S^1_{\bu2}]_{\Z_2}}\big)
\de_{\R}(\cC_0)
+\big(k_{\bu r^*}\!-\!1\!+\!\lr{w_1(V_0^{\vph_0}),[S^1_{\bu r^*}]_{\Z_2}}\big)j_{r^*}'(\cC_0)\\
&\quad
+\lr{w_1(V_0^{\vph_0}),[S^1_{\bu1}]_{\Z_2}}
\big(k_{\bu2}\!-\!1\!+\!\lr{w_1(V_0^{\vph_0}),[S^1_{\bu2}]_{\Z_2}}\big)
\!+\!\big(k_{\bu1}\!-\!1\big)r^*\\
&\hspace{.2in}=n\big(r^*\!-\!1\big)\!+\!2\Z\!+\!\begin{cases}
n(|\pi_0(\Si_0^{\si_0})|\!+\!W_1(V_0,\vph_0))\in\Z_2,
&\hbox{if}~\fp\!\in\!\cP^-_{\Si_0}(V_0^{\vph_0});\\
(n\!+\!1)(|\pi_0(\Si_0^{\si_0})|\!+\!W_1(V_0,\vph_0))\in\Z_2,
&\hbox{if}~\fp\!\in\!\cP^+_{\Si_0}(V_0^{\vph_0}).
\end{cases}\end{split}\end{equation*}
\end{crl}

\chapter{Real Enumerative Geometry}
\label{Geom_chap}

Enumerative geometry of algebraic varieties is a field of 
mathematics that dates back to the nineteenth century.
The general goal of this subject is to count geometric objects 
that satisfy pre-specified geometric conditions.
The objects are typically complex curves in smooth algebraic manifolds
that are usually required to be of a specified genus,
to represent a given homology class, 
and to meet a given collection of geometric constraints.
The prototypical example of such a count is 
the number~$N_d$ of degree~$d$ rational curves that pass through $3d\!-\!1$ 
points in general position in the complex projective plane~$\P^2$.
As the space of $(3d\!-\!1)$-tuples of points in~$\P^2$ in general position
is path-connected, $N_d$ does not depend on the choice of such a tuple.
It is fairly straightforward to see that $N_1,N_2\!=\!1$ and $N_3\!=\!12$.
The result $N_4\!=\!620$ of~\cite[p378]{Ze} first obtained in the middle 
of the nineteenth century is significantly trickier.
The higher-degree numbers~$N_d$ remained unknown until the early~1990s,
when a recursive formula for these numbers was deduced in~\cite{KM}
from the WDVV relation of string theory. 
It was mathematically confirmed shortly after in~\cite{RT,MS0}
using the interpretation of the numbers~$N_d$ in terms of moduli spaces of 
$J$-holomorphic maps from the Riemann sphere~$\P^1$ as in~\cite{Gr}.\\

While a degree~$d$ polynomial over~$\C$ has precisely $d$ roots (counted with multiplicity),
a degree~$d$ polynomial over~$\R$ can have any number~$d'$ of roots so that  
$0\!\le\!d'\!\le\!d$ and $d'$ is of the same parity~as~$d$.
The number of real curves in a real algebraic variety 
that are of a specified genus, represent a given homology class, 
meet a given collection of geometric constraints likewise
generally depends on the position of the constraints and not just on their type.
It may even come as a surprise that there can be nontrivial lower bounds on
counts of real curves and other geometric objects.
For example, the number of real cubics passing through 8 general points in~$\R\P^2$ can be only
8, 10, or~12 (the complex count); see~\cite{DegKhar}.\\

Invariant {\it signed} counts of real rational $J$-holomorphic curves in compact 
real symplectic fourfolds and sixfolds $(X,\om,\phi)$, 
now known as \sf{Welschinger's invariants},
were first defined in~\cite{Wel4,Wel6,Wel6b} and interpreted 
in terms of moduli spaces of $J$-holomorphic maps from the disk~$\bD^2$ in~\cite{Sol}. 
An adaptation of the interpretation of~\cite{Sol} in terms of real $J$-holomorphic maps
from the Riemann sphere~$\P^1$ later appeared in~\cite{Ge2}.
The two moduli-theoretic perspectives on Welschinger's invariants lead to 
the WDVV-type relations for these curve counts announced in~\cite{Sol2} and
established 12~years later in~\cite{RealWDVV,RealWDVV3}.
In contrast to the intrinsically defined curve signs of~\cite{Wel4}
for real symplectic fourfolds,
the curve signs of \cite{Sol,Ge2} depend on the choice of a relative $\Pin$-structure 
on the fixed locus~$X^{\phi}$ of the involution~$\phi$ on~$X$.
The curve signs of~\cite{Wel6,Wel6b} for real symplectic sixfolds depend on
the choice of an $\OSpin$-structure on~$X^{\phi}$, but in a more intrinsic way 
than in~\cite{Sol,Ge2}.
The main theorems of this chapter, Theorems~\ref{WelSolComp_thm} and~\ref{WelSolComp3_thm}, 
provide precise relations between the signs of~\cite{Wel4,Wel6,Wel6b} and \cite{Sol,Ge2}.
While the systematic approach behind the proofs of these theorems is motivated by 
the approach outlined in \cite[Section~8]{Sol},
the statements of these theorems are different from the sign comparison predictions
in~\cite{Sol,Sol2}.\\

Section~\ref{PinSpinImm_sec} obtains the key topological results concerning $\Pin$-structures
on a surface~$Y$ used in the proof of Theorem~\ref{WelSolComp_thm}.
Theorems~\ref{ShExSeqImm_thm} and~\ref{ShExSeqImm_thm2} in this section relate
the restriction of a $\Pin$-structure~$\fp$ on~$Y$ along an immersed loop~$\al$ in~$Y$
and the number of nodes of~$\al$.
The proof of these theorems uses the {\it classical} perspective of 
Definition~\ref{PinSpin_dfn} on $\Spin$- and $\Pin$-structures.\\

In Section~\ref{g0realmaps_sec}, we describe
invariant counts of real rational curves in real symplectic manifolds
as the degrees of relatively orientable pseudocycles 
from moduli spaces of real pseudoholomorphic maps with signed marked points.
This further re-interpretation of the aforementioned moduli-theoretic perspectives 
on Welschinger's invariants originates in~\cite{RealWDVV} in 
the case of real symplectic fourfolds and was extended in~\cite{RealWDVV3}
to real symplectic sixfolds.
We explain why these degrees give well-defined counts of real rational curves,
which depend only the homology classes of the constraints,
and gather key properties of the resulting invariants in 
Theorems~\ref{RGWs_thm} and~\ref{RGWs_thm2}.
Just as in~\cite{Sol,Ge2}, our description of these invariants involves
orienting determinants of real Cauchy-Riemann operators on the pullbacks
of the real bundle~$(TX,\nd\phi)$ 
by parametrizations of the real curves being counted.
These parametrizations are immersions, if the constraints are chosen generically.
Proposition~\ref{cNsgn_prp} re-interprets the already given definition of real curve counts
in terms of orienting determinants of real Cauchy-Riemann operators on the normal bundles
$(\cN u,\vph)$ to the relevant immersions.
Both definitions of the real curve counts in Section~\ref{g0realmaps_sec} and their comparison
use the properties of orientations
of determinants of real Cauchy-Riemann operators provided by Theorem~\ref{CROrient_thm},
which in turn uses the {\it trivializations} perspective of 
Definition~\ref{PinSpin_dfn3} on $\Spin$- and $\Pin$-structures
and the {\it trivializations} perspective of 
Definition~\ref{RelPinSpin_dfn3} on relative $\Spin$- and $\Pin$-structures.\\

Section~\ref{g0real_sec} recalls Welschinger's definitions of invariant signed counts
of real rational $J$-holomorphic curves in symplectic fourfolds and sixfolds
in~\cite{Wel4,Wel6,Wel6b} and compares them with the definitions described
in Section~\ref{g0realmaps_sec}.  
The two comparison theorems, Theorems~\ref{WelSolComp_thm} and~\ref{WelSolComp3_thm},
are deduced from Proposition~\ref{cNsgn_prp};
the proof of the first theorem also uses Theorem~\ref{ShExSeqImm_thm}.
In particular, 
the proof of Theorem~\ref{WelSolComp_thm} depends on the perspectives of
Definitions~\ref{PinSpin_dfn} and~\ref{PinSpin_dfn3} on $\Spin$- and $\Pin$-structures
being equivalent in a way respecting the properties of these
structures listed in Section~\ref{SpinPinProp_subs};
such an equivalence is provided by Theorem~\ref{SpinStrEquiv_thm}.
Section~\ref{g0real_sec} also contains examples of real curve invariants of 
some simple real symplectic fourfolds and sixfolds;
these examples illustrate some of the implications of the properties 
Theorems~\ref{RGWs_thm} and~\ref{RGWs_thm2} and 
of the comparisons Theorems~\ref{WelSolComp_thm} and~\ref{WelSolComp3_thm}.

\section{Pin-structures and immersions}
\label{PinSpinImm_sec}

Let $Y$ be a smooth manifold.
For $\al\!\in\!\cL(Y)$, define $-\al\!\in\!\cL(Y)$ by
$$-\al(z)=\al(\ov{z})~~\forall\,z\!\in\!S^1\!\subset\!\C.$$
A smooth map $\al\!:S^1\!\lra\!Y$ is an \sf{immersion}\gena{immersion}
if $\nd_z\al\!\neq\!0$ for all \hbox{$z\!\in\!S^1$}.
For such a map, let
$$\cN\al\equiv\frac{\al^*TY}{\Im\,\nd\al}\lra S^1$$
denote its \sf{normal bundle}.
Thus, an immersion~$\al$ determines an exact sequence
\BE{ShExSeqImm_e} 0\lra TS^1\lra \al^*TY\lra \cN\al\lra 0\EE
of real vector bundles over~$S^1$, which we denote by~$\ce(\al)$.
A \sf{node}\gena{node!of immersion} 
of such a~map is a pair $\{z_1,z_2\}$ of distinct points of~$S^1$ 
such that $\al(z_1)\!=\!\al(z_2)$;
a triple self-intersection point of~$\al$ corresponds to 3~nodes in this definition.
We call an immersion~$\al$ \sf{admissible}\gena{admissible curve} 
if 
$$\Im\,\nd_{z_1}\al\neq\Im\,\nd_{z_2}\al$$
for every node $\{z_1,z_2\}$ of~$\al$.
The number of nodes of an admissible immersion~$\al$ is finite;
we denote it by~$\de(\al)$.\nota{deal@$\de(\al)$}
Let $\cL^*(Y)\!\subset\!\cL(Y)$\nota{LYst@$\cL^*(Y)$} denote
the set of admissible immersions into~$Y$.

\subsection{Main statements and examples}
\label{PinSpinImm_subs0}

Suppose $\al\!:S^1\!\lra\!Y$ is an immersion.
By Example~\ref{LBSpinStr_eg}, the standard orientation on $S^1\!\subset\!\C$ determines 
an $\OSpin$-structure $\os_0(TS^1)$ on the first non-trivial vector bundle in 
the short exact sequence~\eref{ShExSeqImm_e}.
By the SpinPin~\ref{SpinPinObs_prop} property in Section~\ref{SpinPinProp_subs}, 
each of the other two bundles in~\eref{ShExSeqImm_e} admits a $\Pin^{\pm}$-structure.
Along with the $H^1(Y;\Z_2)$-equivariance of the second map in~\eref{SpinPinSESdfn_e0}  
in the second input and the SpinPin~\ref{SpinPinStr_prop} property,
this implies that the~map
\BE{ImmBij_e} \cP^{\pm}(\cN\al)\lra\cP^{\pm}(\al^*TY), 
\qquad \fp''\lra \llrr{\os_0(TS^1),\fp''}_{\ce(\al)},\EE
is an $H^1(Y;\Z_2)$-equivariant bijection.
By the SpinPin~\ref{SpinPinStr_prop} property, $\cP^{\pm}(\cN\al)$ consists of two elements.\\

Let $Y$ be a smooth surface (manifold of real dimension~2), 
possibly with boundary, and \hbox{$\al\!:S^1\!\lra\!Y$} be an immersion.
The last non-trivial vector bundle in~\eref{ShExSeqImm_e} is then a line bundle.
The paragraph above Lemma~\ref{PinpCan_lmm} and 
the statement of Lemma~\ref{Pin1_lmm} identify the left-hand side 
of~\eref{ImmBij_e} with $\Z_2$ by specifying a $\Pin^{\pm}$-structure 
$\fp_0^{\pm}(\cN\al)$ corresponding to $0\!\in\!\Z_2$.\\

For a $\Pin^{\pm}$-structure~$\fp$ on~$TY$, we define 
$\ft_{\fp}(\al)\!\in\Z_2$\nota{tp@$\ft_{\fp}(\al)$} 
by 
\BE{tpdfn_e} 
\bllrr{\os_0(TS^1),\fp_{\ft_{\fp}(\al)}^{\pm}(\cN\al)}_{\ce(\al)}=\al^*\fp\,.\EE
If $\cN\al$ is orientable, the identification of $\cP^{\pm}(\cN\al)$ with~$\Z_2$
does not depend on the orientation of~$\al$.
Along with the first statement in
the SpinPin~\ref{SpinPinSES_prop}\ref{DSEquivSum_it} property, this implies~that
$$\ft_{\fp}(-\al)=\ft_{\fp}(\al) ~~\forall\,\fp\!\in\!\cP^{\pm}(TY)
\qquad\hbox{if}~w_1(\al^*TY)=0.$$
If $\cN\al$ is not orientable, the identification of $\cP^+(\cN\al)$ with~$\Z_2$
does not depend on the orientation of~$\al$ and 
the identification of $\cP^-(\cN\al)$ with~$\Z_2$ does.
Along with the first statement in
the SpinPin~\ref{SpinPinSES_prop}\ref{DSEquivSum_it} property, this implies~that 
\BE{tpminga_e}\ft_{\fp}(-\al)=\ft_{\fp}(\al) ~~\forall\,\fp\!\in\!\cP^-(TY), \quad
\ft_{\fp}(-\al)=\ft_{\fp}(\al)\!+\!\blr{w_1(Y),[\al]_{\Z_2}}
~~\forall\,\fp\!\in\!\cP^+(TY).\EE
For $\os\!\in\!\OSp(TY)$, let $\ft_{\os}(\al)\!\equiv\!\ft_{\fp}(\al)$\nota{tos@$\ft_{\os}(\al)$}
be the element of~$\Z_2$ determined by the preimage~$\fp$ of~$\os$ under~\eref{Pin2SpinRed_e}.
By the SpinPin~\ref{Pin2SpinRed_prop} property and~\eref{DSorient2_e},
this element is the same for the two possibilities for the domain of~\eref{Pin2SpinRed_e}.\\

We call a topological surface~$Y$ \sf{closed} if it is compact and without boundary.
For such a surface and $\al,\be\!\in\!H_1(Y;\Z_2)$, let
$$\al\!\cdot\!\be\equiv \blr{(\PD_Y\al)(\PD_Y\be),[Y]_{\Z_2}}\in\Z_2$$
denote the homology intersection number of $\al$ and $\be$.

\begin{thm}\label{ShExSeqImm_thm}
Let $Y$ be a smooth closed surface.
For every $\Pin^-$-structure~$\fp$ on~$TY$, the~map 
$$\mu_{\fp}\!: H_1(Y;\Z_2)\lra \Z_2\nota{mp@$\mu_{\fp}$}, \qquad
\mu_{\fp}([\al]_{\Z_2})=1\!+\!\ft_{\fp}(\al)\!+\!\de(\al)
~~\forall\,\al\!\in\!\cL^*(Y),$$
is well-defined and satisfies
\BE{ShExSeqImm_e1}
\mu_{\fp}\big(\al\!+\!\be\big)
=\mu_{\fp}(\al)\!+\!\mu_{\fp}(\be)\!+\!\al\!\cdot\!\be\!+\!
\lr{w_1(Y),\al}\lr{w_1(Y),\be} \quad\forall\,\al,\be\!\in\!H_1(Y;\Z_2).\EE
Furthermore,
\BE{ShExSeqImm_e1b} \mu_{\eta\cdot\fp}(\al)=\mu_{\fp}(\al)+\lr{\eta,\al}
\qquad\forall~\eta\!\in\!H^1(Y;\Z_2),\,\al\!\in\!H_1(Y;\Z_2).\EE
\end{thm}

\begin{thm}\label{ShExSeqImm_thm2}
Let $Y$ be a smooth closed surface.
For every $\Pin^+$-structure~$\fp$ on~$TY$, the~map 
$$\mu_{\fp}\!: H_1(Y;\Z)\lra \Z_2\nota{mp@$\mu_{\fp}$}, \qquad
\mu_{\fp}([\al]_{\Z})\!=\!1\!+\!\ft_{\fp}(\al)\!+\!\de(\al)
~~\forall\,\al\!\in\!\cL^*(Y),$$
is well-defined and satisfies
\BE{ShExSeqImm_e2}
\mu_{\fp}(0)=0, \quad
\mu_{\fp}\big(\al\!+\!\be\big)=\mu_{\fp}(\al)\!+\!\mu_{\fp}(\be)\!+\!
\al_{\Z_2}\!\cdot\!\be_{\Z_2}
\quad\forall\,\al,\be\!\in\!H_1(Y;\Z),\EE
where $\al_{\Z_2},\be_{\Z_2}\!\in\!H_1(Y;\Z_2)$ are the $\Z_2$-reductions of~$\al,\be$.
Furthermore,
\BE{ShExSeqImm_e2b} \mu_{\eta\cdot\fp}(\al)=\mu_{\fp}(\al)+\lr{\eta,\al_{\Z_2}}
\qquad\forall~\eta\!\in\!H^1(Y;\Z_2),\,\al\!\in\!H_1(Y;\Z).\EE
\end{thm}

\begin{eg}\label{ShExSeqImm_eg0}
Let $\al\!:S^1\!\lra\!\R^2$ be the standard inclusion of the unit circle
and 
$$\al_-\!:S^1 \lra\SO(2)$$ 
be as in the $n\!=\!2$ case of the proof of Lemma~\ref{SO2Spin_lmm}.
The standard orientations on~$S^1$ and~$\R^2$ determine an orientation 
on~$\cN\al$ via the exact sequence~\eref{ShExSeqImm_e};
it is given by the inward pointing radial direction. 
The $\Spin$-structure~on 
$$T\R^2|_{S^1}=\big(\R^2\!\times\!\R^2\big)\big|_{S^1}=S^1\!\times\!\R^2$$ 
induced by the canonical $\Spin$-structures on the oriented line bundles
$TS^1$ and~$\cN\al$ via the exact sequence~\eref{ShExSeqImm_e} is described by the~section
$$S^1\lra\SO\big(T\R^2|_{S^1}\big)=S^1\!\times\!\SO(2), \qquad
\ne^{2\pi\fI t}\lra \big(\ne^{2\pi\fI t},\al_-(t)\big).$$
Since this loop generates $\pi_1(\SO(2))$, this $\OSpin$-structure differs
from the restriction of the unique $\OSpin$-structure~$\os$ on the oriented vector bundle~$T\R^2$
over~$\R^2$. Thus, $\ft_{\os}(\al)\!=\!1$.
This implies that $\mu_{\fp}([\al])\!=\!0$ for any $\Pin$-structure~$\fp$ on a smooth closed surface~$Y$
and a homotopically trivial admissible loop~$\al$ in~$Y$.
\end{eg}

\begin{eg}\label{ShExSeqImm_eg0a}
Let $\os_0$ be the $\OSpin$-structure on the 2-torus $(S^1)^2$
induced by the standard homotopy class of trivializations
$$T\big(\!(S^1)^2\big)=TS^1\!\times\!TS^1=(S^1)^2\!\times\!\R^2$$
of $TY$.
For any parametrization $\al$ of a circle $z\!\times\!S^1$ or $S^1\!\times\!z$ in~$(S^1)^2$,
$$\al^*\os_0=\bllrr{\os_0(TS^1),\os_0(\cN\al)}_{\ce(\al)}$$
and so $\ft_{\os_0}(\al)\!=\!0$.
Thus,
\BE{ShExSeqImm0a_e3}\mu_{\os_0}(\al)=1 \qquad\forall\,
\al\!\in\!H_1\big((S^1)^2;\Z_2\big)\!-\!\{0\}.\EE
\end{eg}

\begin{eg}\label{ShExSeqImm_eg0b}
Let $Y$ be the infinite Mobius band, $S\ga_{\R;1}\!\subset\!Y$ be its unit circle bundle,
and $\fp_0^+$ be the $\Pin^+$-structure on~$Y$ defined in Example~\ref{PinpMB_eg}.
By the SpinPin~\ref{SpinPinStr_prop}\ref{PinStr_it} property, $Y$ admits two $\Pin^+$-structures.
By~\eref{ShExSeqImm0b_e5} and the second statement in~\eref{tpminga_e}, 
$\ft_{\fp_0^+}(\al)\!=\!0$ 
for every smooth parametrization $\al\!:S^1\!\lra\!Y$ of~$S\ga_{\R;1}$.
Since the restriction homomorphism
$$H^1(Y;\Z_2)\lra H^1\big(S\ga_{\R;1};\Z_2\big)$$
is trivial, the SpinPin~\ref{SpinPinStr_prop}\ref{PinStr_it} property and
the $H^1(Y;\Z_2)$-equivariance of~\eref{SpinPinSESdfn_e0b} in the first input
then imply that $\ft_{\fp_1^+}(\al)\!=\!0$ 
for the other $\Pin^+$-structure~$\fp_1^+$ on~$Y$ as~well.
\end{eg}

\begin{eg}\label{RP2pin_eg4}
Let $\fp_0\!\equiv\!\fp_0^-(\R\P^2)$ and $\fp_1\!\equiv\!\fp_1^-(\R\P^2)$ 
be the two $\Pin^-$-structures on~$\R\P^2$ as in Example~\ref{RP2pin_eg}.
By~\eref{RP2pin2_e7} and the first statement in~\eref{tpminga_e}, 
$\ft_{\fp_1^-}(\al)\!=\!0$ 
for every smooth homotopically non-trivial embedding~$\al$ of~$S^1$ into~$\R\P^2$.
Along with the conclusion of Example~\ref{ShExSeqImm_eg0} and~\eref{ShExSeqImm_e1b}, 
this implies~that 
\BE{RP2pin4_e3}
\mu_{\fp_1^-}([\al])=\begin{cases}
0,&\hbox{if}~[\al]\!=\!0\!\in\!H_1(\R\P^2;\Z_2);\\
1,&\hbox{if}~[\al]\!\neq\!0\!\in\!H_1(\R\P^2;\Z_2);
\end{cases}
\quad
\mu_{\fp_0^-}([\al])=0~~\forall\,[\al]\!\in\!H_1\big(\R\P^2;\Z_2\big).\EE
\end{eg}

\begin{rmk}\label{ShExSeqImm_rmk}
Let $Y$ be as in Theorems~\ref{ShExSeqImm_thm} and~\ref{ShExSeqImm_thm2}.
Denote by~$\cQ^-(Y)$ and~$\cQ^+(Y)$ the spaces of~maps 
$$\mu\!: H_1(Y;\Z_2)\lra \Z_2 \qquad\hbox{and}\qquad
\mu\!: H_1(Y;\Z)\lra \Z_2$$
satisfying \eref{ShExSeqImm_e1} and~\eref{ShExSeqImm_e2}, respectively, 
with~$\mu_{\fp}$ replaced by~$\mu$.
The group~$H^1(Y;\Z_2)$ acts on~$\cQ^-(Y)$ and~$\cQ^+(Y)$  
freely and transitively by
\begin{alignat*}{2}
\{\eta\!\cdot\!\mu\}(\al)&=\mu(\al)\!+\!\lr{\eta,\al}
&\quad&\forall~\eta\!\in\!H^1(Y;\Z_2),\,\al\!\in\!H_1(Y;\Z_2), \qquad\hbox{and}\\
\{\eta\!\cdot\!\mu\}(\al)&=\mu(\al)\!+\!\lr{\eta,\al_{\Z_2}}
&\quad&\forall~\eta\!\in\!H^1(Y;\Z_2),\,\al\!\in\!H_1(Y;\Z),
\end{alignat*}
respectively.
By the SpinPin~\ref{SpinPinObs_prop}\ref{PinObs_it} property in Section~\ref{SpinPinProp_subs}
and Wu's relations~\eref{WuRel2_e}, $TY$ admits a $\Pin^-$-structure.
For the same reasons, $TY$ admits a $\Pin^+$-structure if \hbox{$\cQ^+(Y)\!\neq\!\eset$}.
Thus, the assignments 
\BE{ShExSeqImmR_e5}\cP^{\pm}(Y)\lra \cQ^{\pm}(Y), \qquad \fp\lra\mu_{\fp}\,,\EE
of Theorems~\ref{ShExSeqImm_thm} and~\ref{ShExSeqImm_thm2} 
are natural $H^1(Y;\Z_2)$-equivariant bijections.
A \sf{quadratic enhancement} (\sf{of the intersection form}) \sf{on}~Y is 
a~map 
$$q\!: H_1(Y;\Z_2)\lra \Z_4 \quad\hbox{s.t.}\quad
q\big(\al\!+\!\be\big)
=q(\al)\!+\!q(\be)\!+\!\io\big(\al\!\cdot\!\be\big)
~~\forall\,\al,\be\!\in\!H_1(Y;\Z_2),$$
where $\io\!:\Z_2\!\lra\!\Z_4$ is the inclusion.
A natural $H^1(Y;\Z_2)$-equivariant bijection between $\cP^-(Y)$ and
the space of quadratic enhancements on~$Y$ of the same flavor as~\eref{ShExSeqImmR_e5}
is provided by \cite[Theorem~3.2]{KirbyTaylor},
based on the perspective on the Kervaire invariant in~\cite{Brown}. 
\end{rmk}

\subsection{Admissible immersions into surfaces}
\label{ShExSeqImm_subs}

We deduce Theorems~\ref{ShExSeqImm_thm} and~\ref{ShExSeqImm_thm2} from
Proposition~\ref{ShExSeqImm_prp} below.
The latter is in turn a consequence of the four lemmas stated
in the present section and proved in the next.
We recall that Wu's relations \cite[p132]{MiSt} imply that  
\BE{WuRel2_e}  w_2(Y)=w_1^2(Y)\in H^2(Y;\Z_2) \quad\hbox{and}\quad
\lr{w_1(Y),\al}=\al^2\in\Z_2~~\forall\,\al\!\in\!H_1(Y;\Z_2)\EE
for a smooth closed surface~$Y$.

\begin{prp}\label{ShExSeqImm_prp}
Let $Y$ be a smooth closed surface and 
$\al_1,\ldots,\al_m\!:S^1\!\lra\!Y$ be admissible loops.
If
\BE{ShExSeqImm_e0}\sum_{i=1}^m[\al_i]_{\Z_2}=0\in H_1(Y;\Z_2)\EE
and $\fp\!\in\!\cP^-(TY)$, then 
\BE{ShExSeqImm_e0a}
\sum_{i=1}^m\big(\ft_{\fp}(\al_i)\!+\!\de(\al_i)\big)=m+
\sum_{i<j}\!\Big(\al_i\!\cdot\!\al_j\!+\!\blr{w_1(Y),[\al_i]_{\Z_2}}
\blr{w_1(Y),[\al_j]_{\Z_2}}\Big)\in\Z_2.\EE
If \eref{ShExSeqImm_e0} holds with $\Z$-coefficients and $\fp\!\in\!\cP^+(TY)$, then
\BE{ShExSeqImm_e0b}
\sum_{i=1}^m\big(\ft_{\fp}(\al_i)\!+\!\de(\al_i)\big)=m+
\sum_{i<j}\al_i\!\cdot\!\al_j\in\Z_2.\EE
\end{prp}

\begin{proof}[{\bf{\emph{Proof of Theorem~\ref{ShExSeqImm_thm}}}}] 
Let $\fp$ be a $\Pin^-$-structure on~$TY$.
If admissible loops $\al_1,\al_2\!:S^1\!\lra\!Y$ agree in $H_1(Y;\Z_2)$,
then~\eref{ShExSeqImm_e0a} and~\eref{WuRel2_e} give
$$\big(\ft_{\fp}(\al_1)\!+\!\de(\al_1)\big)+\big(\ft_{\fp}(\al_2)\!+\!\de(\al_2)\big)
=\al_1\!\cdot\!\al_1\!+\!\blr{w_1(Y),[\al_1]_{\Z_2}}=0\,.$$
Thus, the~map
$$\mu_{\fp}\!: H_1(Y;\Z_2)\lra \Z_2, \qquad 
\mu_{\fp}\big([\al]_{\Z_2}\big)=\ft_{\fp}(\al)\!+\!\de(\al)\!+\!1,$$
is well-defined ($\al\!:S^1\!\lra\!Y$ is an admissible loop above).
By the $m\!=\!3$ case of~\eref{ShExSeqImm_e0a} and~\eref{WuRel2_e}, 
\begin{equation*}\begin{split}
\mu_{\fp}\big(\al\!+\!\be\big)\!+\!\mu_{\fp}(\al)\!+\!\mu_{\fp}(\be)
&=\al\!\cdot\!\be\!+\!\lr{w_1(Y),\al}\lr{w_1(Y),\be}
+(\al\!+\!\be)\!\cdot\!(\al\!+\!\be)\!+\!\lr{w_1(Y),\al\!+\!\be}\\
&=\al\!\cdot\!\be\!+\!\lr{w_1(Y),\al}\lr{w_1(Y),\be}.
\end{split}\end{equation*}
This establishes~\eref{ShExSeqImm_e1}.
Since the $H^1(Y;\Z_2)$-action  of the SpinPin~\ref{SpinPinStr_prop} property is natural and free,
\BE{ShExSeqImm2_e10} \ft_{\eta\cdot\fp}(\al)
=\ft_{\fp}(\al)\!+\!\blr{\al^*\eta,[S^1]_{\Z_2}}
=\ft_{\fp}(\al)+\blr{\eta,[\al]_{\Z_2}}\EE
for every $\eta\!\in\!H^1(Y;\Z_2)$ and every admissible loop $\al\!:S^1\!\lra\!Y$.
This implies~\eref{ShExSeqImm_e1b}.
\end{proof}

\begin{proof}[{\bf{\emph{Proof of Corollary~\ref{ShExSeqImm_thm2}}}}] 
Let $\fp$ be a $\Pin^+$-structure on~$TY$.
If admissible loops $\al_1,\al_2\!:S^1\!\lra\!Y$ agree in $H_1(Y;\Z)$,
then~\eref{ShExSeqImm_e0b} and~\eref{WuRel2_e}  give
$$\big(\ft_{\fp}(\al_1)\!+\!\de(\al_1)\big)+\big(\ft_{\fp}(-\al_2)\!+\!\de(\al_2)\big)
=\al_1\!\cdot\!\al_1=\blr{w_1(Y),[\al_1]_{\Z_2}}\,.$$
Along with the second statement in~\eref{tpminga_e}, this implies that the~map
$$\mu_{\fp}\!: H_1(Y;\Z)\lra \Z_2, \qquad 
\mu_{\fp}\big([\al]_{\Z_2}\big)=\ft_{\fp}(\al)\!+\!\de(\al)\!+\!1,$$
is well-defined and satisfies
\BE{ShExSeqImm2_e12} \mu_{\fp}(\al)\!+\!\mu_{\fp}(-\al)=\blr{w_1(Y),\al_{\Z_2}}
\qquad\forall\,\al\!\in\!H_1(Y;\Z).\EE
By the $m\!=\!3$ case of~\eref{ShExSeqImm_e0b} and~\eref{WuRel2_e}, 
\begin{equation*}\begin{split}
\mu_{\fp}\big(\al\!+\!\be\big)\!+\!\mu_{\fp}(-\al)\!+\!\mu_{\fp}(-\be)
&=\al\!\cdot\!\be\!+\!(\al\!+\!\be)\!\cdot\!(\al\!+\!\be)
=\al\!\cdot\!\be\!+\!\lr{w_1(Y),\al}+\lr{w_1(Y),\be}.
\end{split}\end{equation*}
Along with~\eref{ShExSeqImm2_e12}, 
this implies the second claim in~\eref{ShExSeqImm_e2}.
The first claim in~\eref{ShExSeqImm_e2} follows immediately from
Example~\ref{ShExSeqImm_eg0}.
The relation~\eref{ShExSeqImm_e2b} follows from~\eref{ShExSeqImm2_e10}.
\end{proof}

We will deduce Proposition~\ref{ShExSeqImm_prp} from the four lemmas stated below;
they are proved in Section~\ref{SESlmmpfs_subs}.
Lemma~\ref{ShExSeqImm_lmm2a} obtains an analogue of this proposition 
for the boundary components of a compact bordered surface. 
Lemma~\ref{ShExSeqImm_lmm4} is used in the proof of Proposition~\ref{ShExSeqImm_prp}
to deduce its case for disjoint embedded loops from Lemma~\ref{ShExSeqImm_lmm2a}.
Lemma~\ref{ShExSeqImm_lmm3a}
and the crucial Lemma~\ref{ShExSeqImm_lmm3} reduce the general case of this proposition
to this special case.

\begin{lmm}\label{ShExSeqImm_lmm2a}
Let $Y$ be a compact surface with boundary components $\al_1,\ldots,\al_m$.
If $\fp$ is a $\Pin^+$-structure on~$TY$, then
\BE{ShExSeqImm1_e2} 
\sum_{i=1}^m\ft_{\fp}(\al_i)=\chi(Y)\!+\!2\Z\in\Z_2\,.\EE
If $\fp$ is a $\Pin^-$-structure on~$TY$, then
\BE{ShExSeqImm1_e} 
\sum_{i=1}^m\ft_{\fp}(\al_i)=m\!+\!2\Z\in\Z_2\,.\EE
\end{lmm}

\begin{lmm}\label{ShExSeqImm_lmm4}
Let $Y$ be a smooth surface and $C_1,\ldots,C_m\!\subset\!Y$ be a collection
of disjoint simple circles such~that
\BE{ShExSeqImm2_e}\sum_{i=1}^m\big[C_i\big]_{\Z_2}=0\in H_1(Y;\Z_2).\EE
\begin{enumerate}[label=(\arabic*),leftmargin=*]

\item\label{SESimm4_it1}  There exists an immersion $f\!:\Si\!\lra\!Y$ from a compact surface 
with boundary components $(\prt\Si)_1,\ldots,(\prt\Si)_m$ so that 
the restriction of~$f$ to each~$(\prt\Si)_i$ is a homeomorphism onto~$C_i$.

\item\label{SESimm4_it2}  If in addition the circles~$C_i$ are oriented and~\eref{ShExSeqImm2_e} 
holds with $\Z$ coefficients, then 
there exists an immersion $f\!:\Si\!\lra\!Y$ from a compact oriented surface 
with boundary components $(\prt\Si)_1,\ldots,(\prt\Si)_m$ so that 
the restriction of~$f$ to each~$(\prt\Si)_i$ is an orientation-preserving homeomorphism onto~$C_i$. 

\end{enumerate}
\end{lmm}

Let $Y$ be a smooth surface and
$\al_1,\ldots,\al_m\!:S^1\!\lra\!Y$ be admissible loops.
By applying a small deformation,
we can deform these loops so that all their
intersection and self-intersection points are simple nodes, 
i.e.
$$\big|\big\{(i,z)\!\in\!\{1,\ldots,m\}\!\times\!S^1\!:\al_i(z)\!=\!p\big\}\big|\le2
~~\forall\,p\!\in\!Y,\quad
\Im\,\nd_{z_1}\al_i\neq\Im\,\nd_{z_2}\al_j~~\forall\,(i,z_1)\!\neq\!(j,z_2).$$
We can then smooth out all intersection and self-intersection points 
in accordance with the orientations
of the loops inside of small coordinate charts as in Figure~\ref{smoothingcol_fig0}.
We call the new collection
$$\al_1',\ldots,\al_{m'}'\!:S^1\lra Y$$
of oriented loops a {\it smoothing} of the collection $\al_1,\ldots,\al_m$.

\begin{figure}
\begin{pspicture}(.3,-1.9)(10,2.2)
\psset{unit=.4cm}
\psline[linewidth=.05](9,-3)(15,-3)\psline[linewidth=.05](9,3)(15,3)
\psline[linewidth=.05](9,-3)(9,3)\psline[linewidth=.05](15,-3)(15,3)
\psline[linewidth=.1,arrowsize=6pt]{->}(7,-5)(17,5)
\psline[linewidth=.1,arrowsize=6pt]{->}(7,5)(17,-5)
\pscircle*(12,0){.3}\rput(12,-.8){$p$}
\rput(7,4){$\al_i$}\rput(7,-4){$\al_j$}
\psline[linewidth=.05](29,-3)(35,-3)\psline[linewidth=.05](29,3)(35,3)
\psline[linewidth=.05](29,-3)(29,3)\psline[linewidth=.05](35,-3)(35,3)
\psline[linewidth=.1,arrowsize=6pt](27,-5)(30.5,-1.5)
\psline[linewidth=.1,arrowsize=6pt]{->}(33.5,1.5)(37,5)
\psline[linewidth=.1,arrowsize=6pt](27,5)(30.5,1.5)
\psline[linewidth=.1,arrowsize=6pt]{->}(33.5,-1.5)(37,-5)
\psarc[linewidth=.1](32,-3){2.12}{45}{135}
\psarc[linewidth=.1](32,3){2.12}{225}{315}
\rput(27,4){$\al_{i'}'$}\rput(27,-4){$\al_{j'}'$}
\end{pspicture}
\caption{Smoothing a node of a collection of oriented curves}
\label{smoothingcol_fig0}
\end{figure}

\begin{lmm}\label{ShExSeqImm_lmm3a}
Let $Y$ be a smooth surface and
$\al_1,\ldots,\al_m\!:S^1\!\lra\!Y$ be admissible loops.
If $\al_1',\ldots,\al_{m'}'$ is a {\it smoothing} of $\al_1,\ldots,\al_m$, then
\BE{ShExSeqImm3_e1a} \sum_{i=1}^{m'}[\al_i']_{\Z}=\sum_{i=1}^m[\al_i]_{\Z}
\in H_1(Y;\Z).\EE
If in addition $Y$ is closed, then
\BE{ShExSeqImm3_e1b} 
m'\equiv m+\sum_{i=1}^m\de(\al_i)+\sum_{i<j}\!\al_i\!\cdot\!\al_j \!\!\mod2\,.\EE
\end{lmm}

\begin{lmm}\label{ShExSeqImm_lmm3}
Let $Y$, $\al_1,\ldots,\al_m$, and $\al_1',\ldots,\al_{m'}'$ 
be as in Lemma~\ref{ShExSeqImm_lmm3a}.
\begin{enumerate}[label=(\alph*),leftmargin=*]

\item If $\fp$ is a $\Pin^+$-structure on~$TY$, then
\BE{ShExSeqImm3_e2}
\sum_{i=1}^{m'}\ft_{\fp}(\al_i')=\sum_{i=1}^m\ft_{\fp}(\al_i)\,.\EE

\item If $\fp$ is a $\Pin^-$-structure on~$TY$, then
\BE{ShExSeqImm3_e3}\begin{split}
\sum_{i=1}^{m'}\ft_{\fp}(\al_i')&+
\sum_{i<j}\!\blr{w_1(Y),[\al_i']_{\Z_2}}\blr{w_1(Y),[\al_j']_{\Z_2}}\\
&\qquad=\sum_{i=1}^m\ft_{\fp}(\al_i)+
\sum_{i<j}\!\blr{w_1(Y),[\al_i]_{\Z_2}}\blr{w_1(Y),[\al_j]_{\Z_2}}\,.
\end{split}\EE
\end{enumerate}
\end{lmm}

\begin{proof}[{\bf{\emph{Proof of Proposition~\ref{ShExSeqImm_prp}}}}]
By Lemmas~\ref{ShExSeqImm_lmm4} and~\ref{ShExSeqImm_lmm3}, 
we can replace the collection $\al_1,\ldots,\al_m$
by its smoothing.
Thus, we can assume that the images $C_1,\ldots,C_m$ of the loops~$\al_1,\ldots,\al_m$
are disjoint circles.
Along with~\eref{ShExSeqImm_e0} and~\eref{WuRel2_e}, this implies~that
\BE{ShExSeqImm_e4} \blr{w_1(Y),[C_i]_{\Z_2}}=0 \qquad\forall\,i\!=\!1,\ldots,m.\EE
Thus, \eref{ShExSeqImm_e0a} and~\eref{ShExSeqImm_e0b} in this situation reduce~to
\BE{ShExSeqImm_e5}
\sum_{i=1}^m\ft_{\fp}(\al_i)\equiv m \mod2.\EE
Since the normal bundle of each $C_i$ is orientable by~\eref{ShExSeqImm_e4},
$\ft_{\fp}(\al_i)$ does not depend on the orientation of~$\al_i$ even for 
\hbox{$\fp\!\in\!\cP^+(TY)$} in this case.\\

\noindent
Let  $f\!:\Si\!\lra\!Y$  be as in Lemma~\ref{ShExSeqImm_lmm4}\ref{SESimm4_it1}.
Thus,
\BE{ShExSeqImm_e7}
\sum_{i=1}^m\ft_{\fp}(\al_i)\equiv \sum_{i=1}^{m'}\ft_{f^*\fp}\big(f|_{(\prt\Si)_i}\big)
\!\!\mod2 \qquad \forall\,\fp\!\in\!\cP^{\pm}(TY).\EE
Along with~\eref{ShExSeqImm1_e} with $(Y,\fp)$ replaced by $(\Si,f^*\fp)$,
this implies \eref{ShExSeqImm_e5} for all $\fp\!\in\!\cP^-(TY)$.\\

\noindent
Suppose \eref{ShExSeqImm_e0} holds with $\Z$-coefficients.
Let  $f\!:\Si\!\lra\!Y$  be as in Lemma~\ref{ShExSeqImm_lmm4}\ref{SESimm4_it2}.
Since $\Si$ is oriented, $\chi(\Si)$ is of the same parity as~$m'$.
Along with~\eref{ShExSeqImm1_e2} with $(Y,\fp)$ replaced by $(\Si,f^*\fp)$,
this implies \eref{ShExSeqImm_e5} for all $\fp\!\in\!\cP^+(TY)$.
\end{proof}

\begin{eg}\label{Pinsum_eg}
Let $Y$ be a smooth closed surface.
By~\eref{ShExSeqImm_e1} and~\eref{WuRel2_e},
\hbox{$\mu_{\fp}(0)\!=\!0$} for every \hbox{$\fp\!\in\!\cP^-(TY)$};
this is consistent with the conclusion of Example~\ref{ShExSeqImm_eg0}.
By~\eref{ShExSeqImm_e2} and~\eref{WuRel2_e},
$$\mu_{\fp}(\al)\!+\!\mu_{\fp}(-\al)=\blr{w_1(Y),\al_{\Z_2}}
\qquad\forall\,\al\!\in\!H_1(Y;\Z),\fp\!\in\!\cP^+(TY)\,;$$
this is consistent with the second statement in~\eref{tpminga_e}.
This relation in particular implies that the normal bundle to an (embedded) circle~$C$
in a smooth closed connected unorientable surface~$Y$ 
representing the unique two-torsion element~$\al$ of $H_1(Y;\Z)$ 
is trivializable if $TY$ admits a $\Pin^+$-structure.
This can be obtained directly from the evenness of~$\chi(Y)$,
as such a class~$\al$ is the sum of the basis elements $\al_1,\ldots,\al_m$
in the proof of Corollary~\ref{RP2pin_crl} below.
\end{eg}

\begin{crl}\label{RP2pin_crl}
Let $Y$ be a closed connected unorientable surface.
If the dimension of $H_1(Y;\Z_2)$ is at most~3, then
$TY$ admits a canonical $\Pin^-$-structure $\fp_0^-(Y)$;
it is preserved by every homeomorphism of~$Y$.
\end{crl}

\begin{proof}
We first recall an observation from the proof of \cite[Lemma~2.2]{XCapsSetup}.
By \cite[Theorem~77.5]{Mu}, $Y$ is the connected sum of $m$~copies of~$\R\P^2$ 
for some $m\!\in\!\Z^+$; the dimension of $H_1(Y;\Z_2)$ then is~$m$.
By \cite[Theorem~77.5]{Mu}, $Y$ can be represented by 
the labeling scheme $\al_1\al_1\!\ldots\!\al_m\al_m$; see Figure~\ref{Tm_fig}.
From the labeling scheme, it is immediate that  
$\al_i\!\cdot\!\al_i\!=\!1$ and $\al_i\!\cdot\!\al_j\!=\!0$ if $i\!\neq\!j$.
Thus, there exists a basis $\al_1,\ldots,\al_m$ for $H_1(Y;\Z_2)$ diagonalizing 
the $\Z_2$-intersection form.\\

\begin{figure}
\begin{pspicture}(-7,-2.3)(10,1)
\psset{unit=.4cm}
\psline[linewidth=.04](6,0)(3,1.73)
\psline[linewidth=.04,arrowsize=5pt]{->}(6,0)(4.5,.86)
\psline[linewidth=.04](3,1.73)(0,0)
\psline[linewidth=.04,arrowsize=5pt]{->}(3,1.73)(1.5,.86)
\psline[linewidth=.04](0,0)(0,-3.46)
\psline[linewidth=.04,arrowsize=5pt]{->}(0,0)(0,-1.73)
\psline[linewidth=.04](6,0)(6,-3.46)
\psline[linewidth=.04,arrowsize=5pt]{->}(6,-3.46)(6,-1.73)
\psline[linewidth=.04](0,-3.46)(3,-5.19)
\psline[linewidth=.04,arrowsize=5pt]{->}(0,-3.46)(1.5,-4.32)
\psline[linewidth=.04](3,-5.19)(6,-3.46)
\psline[linewidth=.04,arrowsize=5pt]{->}(3,-5.19)(4.5,-4.32)
\pscircle*(6,0){.15}\pscircle*(3,1.73){.15}\pscircle*(0,0){.15}
\pscircle*(6,-3.46){.15}\pscircle*(3,-5.19){.15}\pscircle*(0,-3.46){.15}
\rput(5.3,1.3){$\al_1$}\rput(.8,1.3){$\al_1$}\rput(3,-.2){$\al_1'$}
\rput(5.3,-4.66){$\al_3$}\rput(.8,-4.76){$\al_2$}
\rput(-.6,-1.73){$\al_2$}\rput(6.7,-1.73){$\al_3$}
\pscircle*(.75,.43){.12}\pscircle*(3.75,1.3){.12}
\pnode(.75,.43){A}\pnode(3.75,1.3){B}
\ncarc[nodesep=0,arcangleA=-76.1,arcangleB=-53.9,ncurv=.6,linestyle=dashed,linewidth=.03]{-}{A}{B}
\end{pspicture}
\caption{Labeling scheme for $Y\!=\!\R\P^2\!\#\R\P^2\!\#\R\P^2$ and a deformation~$\al_1'$ of 
the loop $\al_1$ used to compute the intersection product on $H_1(Y;\Z_2)$}
\label{Tm_fig}
\end{figure}

\noindent
By the Universal Coefficient Theorem for Cohomology \cite[Theorems~53.5]{Mu2},
the homomorphism
$$\ka\!:H^1(Y;\Z_2)\lra \Hom_{\Z_2}\big(H_1(Y;\Z_2),\Z_2\big), \qquad
\big\{\ka(\eta)\big\}(\al)=\lr{\eta,\al},$$
is an isomorphism.
By the SpinPin~\ref{SpinPinObs_prop} property in Section~\ref{SpinPinProp_subs} and~\eref{WuRel2_e}, 
$TY$ admits a $\Pin^-$-structure~$\fp$.
Along with~\eref{ShExSeqImm_e1b}, these two statements imply that $TY$ admits 
a unique $\Pin^-$-structure $\fp_0^-(Y)$ so~that 
\BE{RP2pi_e3} \mu_{\fp_0^-(Y)}(\al_i)=0 \qquad\forall\,i\!=\!1,\ldots,m.\EE
If $m\!\le\!3$,
a basis $\al_1,\ldots,\al_m$ diagonalizing the $\Z_2$-intersection form is unique
up to the permutations of its elements.
The $\Pin^-$-structure determined by~\eref{RP2pi_e3} then does not depend
on the choice of such a basis and is preserved by every homeomorphism of~$Y$.
\end{proof}

\begin{rmk}\label{RP2pin_rmk} If $m\!=\!4$ and $\al_1,\ldots,\al_4$ are as 
in the proof of Corollary~\ref{RP2pin_crl}, then 
$$\al_1'\equiv\al_2\!+\!\al_3\!+\!\al_4, \quad 
\al_2'\equiv\al_1\!+\!\al_3\!+\!\al_4, \quad 
\al_3'\equiv\al_1\!+\!\al_2\!+\!\al_4, \quad 
\al_4'\equiv\al_1\!+\!\al_2\!+\!\al_3$$
is another basis for $H_1(Y;\Z_2)$ diagonalizing the $\Z_2$-intersection form 
and \hbox{$\mu_{\fp_0^-(Y)}(\al_i')\!=\!1$} for all~$i$.
Every disjoint collections \hbox{$C_1,\ldots\!C_m\!\subset\!Y$} of circles representing 
a diagonalizing basis for $H_1(Y;\Z_2)$ presents~$Y$ as the real blowup of~$S^2$
at $m$~distinct points (with $S^2$ obtained from~$Y$ by contracting the circles).
For any two collections $x_1,\ldots,x_m$ and $x_1',\ldots,x_m'$ of distinct points on~$S^2$,
there exists a homeomorphism from the blowup~$Y$ at the first set to the blowup~$Y'$
at the second set so that the $i$-th ``exceptional divisor" $C_i$ for the first blowup
is taken to~$C_i'$ for all~$i$.
For any two diagonalizing bases $\al_1,\ldots,\al_m$ and $\al_1',\ldots,\al_m'$
for $H_1(Y;\Z_2)$, there thus exists a homeomorphism~$f$ of~$Y$ so that \hbox{$f_*\al_i\!=\!\al_i'$}
for all~$i$.
This implies that Corollary~\ref{RP2pin_crl} does not extend to $m\!\ge\!4$.
\end{rmk}

\subsection{Proofs of Lemmas~\ref{ShExSeqImm_lmm2a}-\ref{ShExSeqImm_lmm3}}
\label{SESlmmpfs_subs}

\noindent
We now establish the four lemmas stated in Section~\ref{ShExSeqImm_subs}.

\begin{proof}[{\bf{\emph{Proof of Lemma~\ref{ShExSeqImm_lmm2a}}}}]
We can assume that $Y$ is connected and $m\!\ge\!1$.
Let $\Si$ be the closed surface obtained from~$Y$ by attaching disks along 
the boundary components of~$Y$.
By the SpinPin~\ref{SpinPinObs_prop} property in Section~\ref{SpinPinProp_subs} and~\eref{WuRel2_e}, 
$T\Si$ admits a $\Pin^-$-structure;
if the Euler characteristic
$$\chi(\Si)=\chi(Y)\!+\!m$$
is even, then $T\Si$ also admits a $\Pin^+$-structure.
By Example~\ref{ShExSeqImm_eg0}, 
\eref{ShExSeqImm1_e2} and \eref{ShExSeqImm1_e} hold if $\fp$ is the restriction~$\fp_0$ of 
a $\Pin^{\pm}$-structure on~$T\Si$.
By Corollary~\ref{X2VB_crl2b}, any other $\Pin^{\pm}$-structure~$\fp$ 
on~$TY$ differs from~$\fp_0$ on an even number of boundary components.
This establishes~\eref{ShExSeqImm1_e2} if $\chi(Y)$ and~$m$ are of the same parity
and~\eref{ShExSeqImm1_e} in all cases.\\

\noindent
Suppose $\chi(Y)$ and~$m$ are not of the same parity, i.e.~$\chi(\Si)$ is odd, 
and thus $\Si$ does not admit a $\Pin^+$-structure by the SpinPin~\ref{SpinPinObs_prop} property.
Let $\Si'\!\subset\!\Si$ be the surface obtained from~$Y$ by attaching disks
along the boundary components $\al_2,\ldots,\al_m$ of~$Y$.
By the SpinPin~\ref{SpinPinObs_prop} property, $T\Si'$ admits a $\Pin^+$-structure~$\fp'$.
If $\ft_{\fp'}(\al_1)\!=\!1$, then  Example~\ref{ShExSeqImm_eg0} implies that 
$\fp'$ extends to a $\Pin^+$-structure on~$T\Si$.
Thus, 
$$\ft_{\fp'}(\al_1)=0\equiv\chi(\Si')\!=\!\chi(Y)\!+\!(m\!-\!1) \mod2. $$
Along with Example~\ref{ShExSeqImm_eg0}, this implies that
\eref{ShExSeqImm1_e2} holds if $\fp$ is the restriction~$\fp_0$ of~$\fp'$.
By Corollary~\ref{X2VB_crl2b}, any other $\Pin^+$-structure~$\fp$ 
on~$TY$ differs from~$\fp_0$ on an even number of boundary components.
This establishes~\eref{ShExSeqImm1_e2} in the remaining case.
\end{proof}

\begin{proof}[{\bf{\emph{Proof of Lemma~\ref{ShExSeqImm_lmm4}}}}]
Let $\Si_1,\ldots,\Si_n$ be the topological components of the surface
obtained by cutting~$Y$ along the circles $C_1,\ldots,C_m$.
Triangulate~$Y$ so that each~$C_i$ is a subcomplex (and no edges of
the same 2-simplex are identified).
Then,
\BE{ShExSeqImm2_e3}\sum_{i=1}^m\!C_i=\prt\sum_{j=1}^k\!\De_j\in C_1(Y;\Z_2)\EE
for some 2-simplices $\De_j\!\subset\!Y$.
Since the right-hand side of~\eref{ShExSeqImm2_e3} contains no 1-simplicies
in the interior of any~$\Si_j$,
$\Si\!\equiv\!\De_1\!\cup\!\ldots\!\cup\!\De_k$
is the union of some of the bordered surfaces $\Si_1,\ldots,\Si_n$;
the boundary of~$\Si$ is \hbox{$C_1\!\cup\!\ldots\!\cup\!C_m$}.\\

\noindent
By the above, the boundary of each surface~$\Si_j$ above is a union of 
some of the circles $C_1,\ldots,C_m$.
Suppose the circle~$C_i$ are oriented and~\eref{ShExSeqImm2_e} holds with $\Z$~coefficients.
Then,
\BE{ShExSeqImm2_e6}\sum_{i=1}^m\!C_i=\prt\sum_{j=1}^k\!a_j\De_j\in C_1(Y;\Z)\EE
for some $a_j\!\in\!\Z^+$ and oriented simplicies $\De_j\!\subset\!Y$ 
so that each unoriented simplex appears at most once.
For each $a\!\in\!\Z^+$, let $\Si_a'$ be the union of the 2-simplices~$\De_j$ 
with $a_j\!=\!a$.
Since the right-hand side of~\eref{ShExSeqImm2_e6} contains no 1-simplicies
in the interior of any~$\Si_j$, each~$\Si_a'$ 
is the union of some of the bordered surfaces $\Si_1,\ldots,\Si_n$.
The orientations of the simplices~$\De_j$ orient~$\Si_a'$.
Define
$$\Si=\bigsqcup_{a=1}^{\i}\{1,\ldots,a\}\!\times\!\Si_a',
\quad
f\!:\Si\lra Y,~~f(b,x)=x~~\forall\,x\!\in\!\Si_a',~b\!=\!1,\ldots,a,~a\!\in\!\Z^+\,.$$
By~\eref{ShExSeqImm2_e6}, the initial orientation of each~$C_i$ agrees 
with its boundary orientation as a component of~$\prt\Si$.
\end{proof}

\begin{proof}[{\bf{\emph{Proof of Lemma~\ref{ShExSeqImm_lmm3a}}}}]
The right-hand sides in~\eref{ShExSeqImm3_e1a} and~\eref{ShExSeqImm3_e1b}
 are invariant under small deformations of the loops~$\al_i$.
Thus, we can assume that all intersection and self-intersection points of
$\al_1,\ldots,\al_m$ are simple nodes.
Smoothing each node then changes the number of loops by one.
This implies~\eref{ShExSeqImm3_e1b};
\eref{ShExSeqImm3_e1a} is immediate from the construction of the smoothing.
\end{proof}

\begin{proof}[{\bf{\emph{Proof of Lemma~\ref{ShExSeqImm_lmm3}}}}]
The right-hand sides in~\eref{ShExSeqImm3_e2} and~\eref{ShExSeqImm3_e3}
are invariant under small deformations of the loops~$\al_i$.
Thus, we can assume that all intersection and self-intersection points of
$\al_1,\ldots,\al_m$ are simple nodes.
We show below that \eref{ShExSeqImm3_e2} and~\eref{ShExSeqImm3_e3}
hold if $\al_1',\ldots,\al_{m'}'$ are obtained from $\al_1,\ldots,\al_m$
by smoothing one node or intersection point~$p$ of $\al_1,\ldots,\al_m$.
This suffices to establish the two identities.\\

\noindent
If $U_0,U_1\!\subset\!S^1$ are open (proper, connected) arcs covering~$S^1$,
we denote by \hbox{$U_{01}\!\subset\!U_0\!\cap\!U_1$} the connected component 
that follows~$U_0\!-\!U_1$ and precedes~$U_1\!-\!U_0$ with respect to the standard orientation of~$S^1$
and by $U_{10}\!\subset\!U_0\!\cap\!U_1$ the other connected component.
Let 
$$\bI_2,\bI_{2;1}\in\O(2), \qquad \wt\bI_2\in\Spin(1)\subset\Pin^{\pm}(1),
 \qquad\hbox{and}\qquad \wt\bI_2,\wh\bI_2\in\Spin(2)\subset\Pin^{\pm}(2)$$
be as in Sections~\ref{SO2Spin_subs} and~\ref{O2Pin_subs}.
Fix a metric on~$TY$.
We establish~\eref{ShExSeqImm_e0a} and~\eref{ShExSeqImm_e0b} using 
the perspective of Definition~\ref{PinSpin_dfn}\ref{PinStrDfn_it}.\\

Let  $\al\!:S^1\!\lra\!Y$ be an immersion.
We trivialize $\al^*TY$ with its metric over~$U_0$ and~$U_1$ 
by taking the first vector field to be the oriented vector field on~$S^1$
of unit length with respect to the pullback metric 
and choosing the second vector field so that the transition function~$g_{01}^{\al}$
from~$U_1$ to~$U_0$ equals~$\bI_2$ on~$U_{01}$.
This implies~that  
$$g_{01}^{\al}\big|_{U_{10}}=\begin{cases}\bI_2,&\hbox{if}~w_1(\al^*TY)\!=\!0;\\
\bI_{2;1},&\hbox{if}~w_1(\al^*TY)\!\neq\!0.
\end{cases}$$
The transition function~$\wt{g}_{01}'$ for $\os_0(TS^1)$ on $TS^1$ equals~$\wt\bI_1$.\\

\noindent
If $\cN\al$ (or equivalently $\al^*TY$) is orientable,
the canonical $\Pin^{\pm}$-structure $\fp_0^{\pm}(\cN\io)$ on~$\cN\al$ 
is specified by the transition function~$\wt{g}_{01}''\!=\!\wt\bI_1$.
By~\eref{DSCorr1_e11}, the transition function 
\BE{ShExSeqImm3_e5}\wt{g}_{01}\!\equiv\!\io_{2;1}\big(\wt{g}_{01}',\wt{g}_{01}''\big)\!:
U_0\!\cap\!U_1\lra\Pin^{\pm}(2)\EE
specifying the $\Pin^{\pm}$-structure $\llrr{\os_0(TS^1),\fp_0^{\pm}(\cN\io)}_{\ce(\al)}$ 
then equals~$\wt\bI_2$. 
The above trivializations for $\al^*TY$ can thus be lifted to
trivializations for the $\Pin^{\pm}$-structure $\al^*\fp$
on $\al^*TY$ so~that the associated transition 
function~$\wt{g}_{01}^{\fp}$ satisfies
\BE{ShExSeqImm3_e7a}\wt{g}_{01}^{\fp}\big|_{U_{01}}=\wt\bI_2, \qquad 
\wt{g}_{01}^{\fp}\big|_{U_{10}}=\wh\bI_2^{\,\ft_{\fp}(\al)},\EE
if $\al^*TY$ is orientable.\\

\noindent
If the vector bundle $\al^*TY$ is not orientable, the transition function~\eref{ShExSeqImm3_e5} 
specifying the $\Pin^{\pm}$-structure $\llrr{\os_0(TS^1),\fp_0^{\pm}(\cN\io)}_{\ce(\al)}$
is described~by
$$\wt{g}_{01}\big|_{U_{01}}=\wt\bI_2, \qquad 
\wt{g}_{01}\big|_{U_{10}}=\wt\bI_{2;1}^{\st},$$
for a certain element $\wt\bI_{2;1}^{\star}\!\in\!\Pin^{\pm}(2)$
in the preimage of~$\bI_{2;1}$ under~\eref{Pindfn_e}.
The above trivializations for $\al^*TY$ can thus be lifted to
trivializations for the $\Pin^{\pm}$-structure $\al^*\fp$
on $\al^*TY$ so~that the associated transition 
function~$\wt{g}_{01}^{\fp}$ satisfies
\BE{ShExSeqImm3_e7b}\wt{g}_{01}^{\fp}\big|_{U_{01}}=\wt\bI_2, \qquad 
\wt{g}_{01}^{\fp}\big|_{U_{10}}=\wt\bI_{2;1}^{\star}\wh\bI_2^{\,\ft_{\fp}(\al)},\EE
if $\al^*TY$ is not orientable.\\

\begin{figure}
\begin{pspicture}(.3,-2.8)(10,3.2)
\psset{unit=.4cm}
\psline[linewidth=.05](9,-3)(15,-3)\psline[linewidth=.05](9,3)(15,3)
\psline[linewidth=.05](9,-3)(9,3)\psline[linewidth=.05](15,-3)(15,3)
\psline[linewidth=.1,arrowsize=6pt]{->}(7,-5)(17,5)
\psline[linewidth=.1,arrowsize=6pt]{->}(7,5)(17,-5)
\pscircle*(12,0){.3}\rput(12,-.8){$p$}
\rput(7,4){$\al$}\rput(7,-4.2){$\al$}\rput(15.7,0){$U'$}
\psline[linewidth=.05]{->}(10,-2.5)(11,-1.5)\psline[linewidth=.05]{<-}(14,-2.5)(13,-1.5)
\psline[linewidth=.05]{->}(10,2.5)(11,1.5)\psline[linewidth=.05]{<-}(14,2.5)(13,1.5)
\psline[linewidth=.05]{->}(11.4,-2.7)(12.6,-2.7)\psline[linewidth=.05]{->}(11.4,2.7)(12.6,2.7)
\psline[linewidth=.05]{->}(9.5,0)(10.8,0)\psline[linewidth=.05]{->}(13.2,0)(14.5,0)
\psarc[linewidth=.08,linestyle=dashed](8,-6){1.41}{135}{270}
\psarc[linewidth=.08,linestyle=dashed](8,6){1.41}{90}{225}
\psarc[linewidth=.08,linestyle=dashed](16,-6){1.41}{-90}{45}
\psarc[linewidth=.08,linestyle=dashed](16,6){1.41}{-45}{90}
\psline[linewidth=.08,linestyle=dashed](8,-7.41)(16,-7.41)
\psline[linewidth=.08,linestyle=dashed](8,7.41)(16,7.41)
\rput(12.2,6.7){$U_a$}\rput(12.2,-6.7){$U_b$}
\rput(9.2,4){$U_{a0}$}\rput(14.8,4){$U_{0a}$}
\rput(9.2,-4){$U_{b0}$}\rput(15,-4){$U_{0b}$}
\psline[linewidth=.05](29,-3)(35,-3)\psline[linewidth=.05](29,3)(35,3)
\psline[linewidth=.05](29,-3)(29,3)\psline[linewidth=.05](35,-3)(35,3)
\psline[linewidth=.1,arrowsize=6pt](27,-5)(30.5,-1.5)
\psline[linewidth=.1,arrowsize=6pt]{->}(33.5,1.5)(37,5)
\psline[linewidth=.1,arrowsize=6pt](27,5)(30.5,1.5)
\psline[linewidth=.1,arrowsize=6pt]{->}(33.5,-1.5)(37,-5)
\psarc[linewidth=.1](32,-3){2.12}{45}{135}
\psarc[linewidth=.1](32,3){2.12}{225}{315}
\rput(27,4){$\al_a$}\rput(27,-4){$\al_b$}\rput(35.7,0){$U'$}
\psline[linewidth=.05]{->}(32,-.88)(33.3,-.88)\psline[linewidth=.05]{->}(33.5,-.88)(34.8,-.88)
\psline[linewidth=.05]{->}(30,-.88)(31.3,-.88)
\psline[linewidth=.05]{->}(32,.88)(33.3,.88)\psline[linewidth=.05]{->}(33.5,.88)(34.8,.88)
\psline[linewidth=.05]{->}(30,.88)(31.3,.88)
\psline[linewidth=.05]{->}(30,-2.5)(31,-1.5)\psline[linewidth=.05]{<-}(34,-2.5)(33,-1.5)
\psline[linewidth=.05]{->}(30,2.5)(31,1.5)\psline[linewidth=.05]{<-}(34,2.5)(33,1.5)
\psline[linewidth=.05]{->}(31.4,-2.7)(32.6,-2.7)\psline[linewidth=.05]{->}(31.4,2.7)(32.6,2.7)
\psline[linewidth=.05]{->}(29.5,0)(30.8,0)\psline[linewidth=.05]{->}(33.2,0)(34.5,0)
\psarc[linewidth=.08,linestyle=dashed](28,-6){1.41}{135}{270}
\psarc[linewidth=.08,linestyle=dashed](28,6){1.41}{90}{225}
\psarc[linewidth=.08,linestyle=dashed](36,-6){1.41}{-90}{45}
\psarc[linewidth=.08,linestyle=dashed](36,6){1.41}{-45}{90}
\psline[linewidth=.08,linestyle=dashed](28,-7.41)(36,-7.41)
\psline[linewidth=.08,linestyle=dashed](28,7.41)(36,7.41)
\rput(32.2,6.7){$U_a$}\rput(32.2,-6.7){$U_b$}
\rput(29.2,4){$U_{a0}$}\rput(34.8,4){$U_{0a}$}
\rput(29.2,-4){$U_{b0}$}\rput(35,-4){$U_{0b}$}
\end{pspicture}
\caption{The space~$S$ and the loops $\al,\al_a,\al_b\!\subset\!S$}
\label{smoothingcol_fig1}
\end{figure}

Let $p\!\in\!Y$ be either a simple node of a loop~$\al$ so that smoothing~$p$ breaks~$\al$
into loops~$\al_a$ and~$\al_b$, as in Figure~\ref{smoothingcol_fig1}, 
or a transverse intersection point of distinct loops~$\al_a$ and~$\al_b$ so that 
smoothing~$p$ combines~$\al_a$ and~$\al_b$ into a single loop~$\al$. 
Since $\al$ is homologous to the sum of~$\al_a$ and~$\al_b$,
\BE{ShExSeqImm3_e10} 
\blr{w_1(Y),[\al_a]_{\Z_2}}+\blr{w_1(Y),[\al_b]_{\Z_2}}
=\blr{w_1(Y),[\al]_{\Z_2}}\,.\EE
Let $U\!\subset\!Y$ be a square coordinate neighborhood of~$p$ 
intersecting the collection $\al_1,\ldots,\al_m$ along the two diagonals
and $U''\!\subset\!U'\!\subset\!U$ be neighborhoods of~$p$ so that 
$\ov{U''}\!\subset\!U'$ and $\ov{U'}\!\subset\!U$.
We smooth out~$p$ inside of~$U''$ and define
\begin{gather*}
S=\big(U'\!\sqcup\!\al^{-1}(Y\!-\!\ov{U''})\big)\!\big/\!\!\sim, \qquad
U'\ni\al(z)\sim z\!\in\!S^1~~\forall\,z\!\in\!\al^{-1}(U'\!-\!\ov{U''}),\\
f\!:S\lra Y, \quad f(x)=\begin{cases}x,&\hbox{if}~x\!\in\!U';\\
\al(z),&\hbox{if}~z\!\in\!\al^{-1}\big(Y\!-\!\ov{U''}\big);
\end{cases} \qquad U_0=S\!-\!\al^{-1}(Y\!-\!U).
\end{gather*}
Let $U_a,U_b\!\subset\!\al^{-1}(Y\!-\!\ov{U'})$ be the connected components 
so that $U_a$ is contained in the domain of the loop~$\al_a$ and 
$U_b$ is contained in the domain of the loop~$\al_b$.
The collection $\{U_a,U_b,U_0\}$ is then an open cover of~$S$.
For $\bu\!=\!\{a,b\}$, we denote by \hbox{$U_{0\bu}\!\subset\!U_0\!\cap\!U_{\bu}$} 
the connected component that follows $U_0\!-\!U_a\!\cup\!U_b$ and precedes $U_{\bu}\!-\!U_0$ 
with respect to the standard orientation of~$S^1$
and by \hbox{$U_{\bu0}\!\subset\!U_0\!\cap\!U_{\bu}$} the other connected component.\\

Let $\xi\!\in\!\Ga(S;f^*TY)$ be such that $\xi|_{U_a}$ is 
the standard oriented unit vector field on
the domain~$S^1$ of~$\al_a$,  $\xi|_{U_b}$ is its analogue for~$\al_b$, and
$\xi|_{U'}$ is as depicted in Figure~\ref{smoothingcol_fig1}.
We trivialize $f^*TY$ over~$U_a$, $U_b$, and~$U_0$ 
by taking the first vector field to be~$\xi$ 
and choosing the second vector field so that the transition functions~$g_{0a}^f$
from~$U_a$ to~$U_0$ and $g_{0b}^f$ from~$U_b$ to~$U_0$ satisfy
$$g_{0a}^f\big|_{U_{0a}}=\bI_2 \qquad\hbox{and}\qquad
g_{0b}^f\big|_{U_{0b}}=\bI_2 \,.$$
This implies~that  
$$g_{0\bu}^f\big|_{U_{\bu0}}=\begin{cases}
\bI_2,&\hbox{if}~w_1(\al_{\bu}^*TY)\!=\!0,\,\bu\!=\!a,b;\\
\bI_{2;1},&\hbox{if}~w_1(\al_{\bu}^*TY)\!\neq\!0,\,\bu\!=\!a,b.
\end{cases}$$
In light of~\eref{ShExSeqImm3_e7a} and~\eref{ShExSeqImm3_e7b},
these trivializations for $f^*TY$ can be lifted to
trivializations for the $\Pin^{\pm}$-structure $f^*\fp$ on $f^*TY$ 
so~that the associated transition functions~$\wt{g}_{0a}^{\fp}$ and~$\wt{g}_{0b}^{\fp}$
satisfy
\BE{ShExSeqImm3_e11}
\wt{g}_{0a}^{\fp}\big|_{U_{0a}},\wt{g}_{0b}^{\fp}\big|_{U_{0b}}=\wt\bI_2, \quad
\wt{g}_{0\bu}^{\fp}\big|_{U_{\bu0}}= \begin{cases}
\wh\bI_2^{\,\ft_{\fp}(\al_{\bu})},&\hbox{if}~w_1(\al_{\bu}^*TY)\!=\!0,\,\bu\!=\!a,b;\\
\wt\bI_{2;1}^{\star}\wh\bI_2^{\,\ft_{\fp}(\al_{\bu})},&\hbox{if}~w_1(\al_{\bu}^*TY)\!\neq\!0,\,
\bu\!=\!a,b.\end{cases}\EE

\vspace{.15in}

Let $\io\!:S^1\!\lra\!S$ be the continuous map so that $\al\!=\!f\!\circ\!\io$.
This map identifies the open subsets 
$$U_a,U_b,U_{0a},U_{a0},U_{0b},U_{b0}\subset S$$
with their preimages.
We denote by \hbox{$U_{0}^{\vdash}\!\subset\!\io^{-1}(U_0)$} 
the connected component containing~$U_{a0}$ and~$U_{0b}$
and by $U_0^{\dashv}\!\subset\!\io^{-1}(U_0)$ the other connected component;
see Figure~\ref{smoothingcol_fig2}.
The collection $\{U_a,U_b,U_0^{\vdash},U_0^{\dashv}\}$ is then 
an open cover of~$S^1$ and
$$U_a\!\cap\!U_0^{\dashv}=U_{0a}, \quad
U_a\!\cap\!U_0^{\vdash}=U_{a0}, \quad
U_b\!\cap\!U_0^{\vdash}=U_{0b}, \quad
U_b\!\cap\!U_0^{\dashv}=U_{b0};$$
all other intersections are empty.\\

\begin{figure}
\begin{pspicture}(-3.2,-1.5)(10,2)
\psset{unit=.4cm}
\rput(6.7,0.2){$U_0^{\vdash}$}\rput(17.3,0.2){$U_0^{\dashv}$}
\psarc[linewidth=.1](8,-2){1.41}{225}{270}\psarc[linewidth=.1](8,2){1.41}{90}{135}
\psarc[linewidth=.1](16,-2){1.41}{-90}{-45}\psarc[linewidth=.1](16,2){1.41}{45}{90}
\psarc[linewidth=.1,arrowsize=6pt]{->}(10,0){4.24}{135}{225}
\psarc[linewidth=.1,arrowsize=6pt]{->}(14,0){4.24}{-45}{45}
\psline[linewidth=.08,linestyle=dashed](8,-3.41)(16,-3.41)
\psline[linewidth=.08,linestyle=dashed](8,3.41)(16,3.41)
\rput(12.2,2.7){$U_a$}\rput(12.2,-2.7){$U_b$}
\rput(7.2,4){$U_{a0}$}\rput(16.8,4){$U_{0a}$}
\rput(7.2,-4){$U_{0b}$}\rput(16.8,-4){$U_{b0}$}
\rput(8.2,2.7){$\wt{g}_{\vdash a}^{\fp}$}\rput(16.4,2.6){$\wt\bI_2$}
\rput(15.8,-2.6){$\wt{g}_{\dashv b}^{\fp}$}\rput(7.8,-2.6){$\wt\bI_2$}
\end{pspicture}
\caption{Open cover of the domain of $\al$ and transition functions~\eref{ShExSeqImm3_e15} 
for~$\al^*\fp$}
\label{smoothingcol_fig2}
\end{figure}

The pullbacks by~$\io$ of the above trivializations of $f^*TY$  induce
trivializations  $\Phi_a,\Phi_b,\Phi_0^{\vdash},\Phi_0^{\dashv}$
of \hbox{$\io^*f^*TY\!=\!\al^*TY$} over $U_a,U_b,U_0^{\vdash},U_0^{\dashv}$, respectively,
so that the first vector field in each trivialization is 
the standard oriented unit vector field on the domain~$S^1$ of~$\al$.
The pullbacks by~$\io$ of the above trivializations of $f^*\fp$ induce 
trivializations  $\wt\Phi_a,\wt\Phi_b,\wt\Phi_0^{\vdash},\wt\Phi_0^{\dashv}$
of $\io^*f^*\fp\!=\!\al^*\fp$ over $U_a,U_b,U_0^{\vdash},U_0^{\dashv}$, respectively.
By~\eref{ShExSeqImm3_e11}, the associated transition functions 
$\wt{g}_{a\dashv}^{\fp}$ from~$U_0^{\dashv}$ to~$U_a$, 
 $\wt{g}_{\vdash a}^{\fp}$ from~$U_a$ to~$U_0^{\vdash}$,  
$\wt{g}_{b\vdash}^{\fp}$ from~$U_0^{\vdash}$ to~$U_b$, and
 $\wt{g}_{\dashv b}^{\fp}$ from~$U_b$ to~$U_0^{\dashv}$ are given~by
\BE{ShExSeqImm3_e15}\begin{aligned}
\wt{g}_{a\dashv}^{\fp}&=\wt\bI_2, &\quad
\wt{g}_{b\vdash}^{\fp}&=\wt\bI_2,\\
\wt{g}_{\vdash a}^{\fp}&= \begin{cases}
\wh\bI_2^{\,\ft_{\fp}(\al_a)},&\hbox{if}~w_1(\al_a^*TY)\!=\!0;\\
\wt\bI_{2;1}^{\star}\wh\bI_2^{\,\ft_{\fp}(\al_a)},&\hbox{if}~w_1(\al_a^*TY)\!\neq\!0;
\end{cases}
&\quad
\wt{g}_{\dashv b}^{\fp}&= \begin{cases}
\wh\bI_2^{\,\ft_{\fp}(\al_b)},&\hbox{if}~w_1(\al_b^*TY)\!=\!0;\\
\wt\bI_{2;1}^{\star}\wh\bI_2^{\,\ft_{\fp}(\al_b)},&\hbox{if}~w_1(\al_b^*TY)\!\neq\!0.
\end{cases}\end{aligned}\EE

\vspace{.15in}

Suppose $w_1(\al_a^*TY)\!=\!0$.
Multiplying the trivializations $\wt\Phi_b$ and $\wt\Phi_0^{\vdash}$ 
by $(\wh\bI_2^{\,\ft_{\fp}(\al_a)})^{-1}$ turns 
the last transition function in~\eref{ShExSeqImm3_e15} into
\BE{ShExSeqImm3_e17a}
\wt{g}_{\dashv b}^{\fp}\wh\bI_2^{\,\ft_{\fp}(\al_a)}= \begin{cases}
\wh\bI_2^{\,\ft_{\fp}(\al_a)+\ft_{\fp}(\al_b)},&\hbox{if}~w_1(\al_b^*TY)\!=\!0;\\
\wt\bI_{2;1}^{\star}\wh\bI_2^{\,\ft_{\fp}(\al_a)+\ft_{\fp}(\al_b)},&\hbox{if}~w_1(\al_b^*TY)\!\neq\!0;
\end{cases}\EE
and makes the remaining three transition functions~$\wt\bI_2$.
This modification does not effect the associated trivializations
of~$\al^*TY$.
Comparing~\eref{ShExSeqImm3_e17a} with~\eref{ShExSeqImm3_e7a} and~\eref{ShExSeqImm3_e7b}, 
we conclude~that 
$$\ft_{\fp}(\al)=\ft_{\fp}(\al_a)+\ft_{\fp}(\al_b)$$
if $w_1(\al_a^*TY)\!=\!0$ or $w_1(\al_b^*TY)\!=\!0$.\\

Suppose $w_1(\al_a^*TY)\!\neq\!0$ and $w_1(\al_b^*TY)\!\neq\!0$.
Multiplying the trivializations $\wt\Phi_b$ and $\wt\Phi_0^{\vdash}$ 
by $(\wt\bI_{2;1}^{\star}\wh\bI_2^{\,\ft_{\fp}(\al_a)})^{-1}$
turns 
the last transition function in~\eref{ShExSeqImm3_e15} into
\BE{ShExSeqImm3_e17b}
\wt{g}_{\dashv b}^{\fp}\wt\bI_{2;1}^{\star}\wh\bI_2^{\,\ft_{\fp}(\al_a)}
=\big(\wt\bI_{2;1}^{\star}\big)^{\!2\,}\wh\bI_2^{\,\ft_{\fp}(\al_a)+\ft_{\fp}(\al_b)}
=\begin{cases}
\wh\bI_2^{\,\ft_{\fp}(\al_a)+\ft_{\fp}(\al_b)},&\hbox{if}~\fp\!\in\!\cP^+(TY);\\
\wh\bI_2^{\,\ft_{\fp}(\al_a)+\ft_{\fp}(\al_b)+1},&\hbox{if}~\fp\!\in\!\cP^-(TY);
\end{cases}\EE
and makes the remaining three transition functions~$\wt\bI_2$.
This modification multiplies the associated trivializations~$\Phi_b$ and $\Phi_0^{\vdash}$ 
of~$\al^*TY$ by~$\bI_{2;1}$ and thus leaves the first vector field unchanged.
Comparing~\eref{ShExSeqImm3_e17b} with~\eref{ShExSeqImm3_e7b}, 
we conclude~that 
$$\ft_{\fp}(\al)=\begin{cases}
\ft_{\fp}(\al_a)\!+\!\ft_{\fp}(\al_b),&\hbox{if}~\fp\!\in\!\cP^+(TY);\\
\ft_{\fp}(\al_a)\!+\!\ft_{\fp}(\al_b)\!+\!1,&\hbox{if}~\fp\!\in\!\cP^-(TY);
\end{cases}$$
if $w_1(\al_a^*TY)\!\neq\!0$ and $w_1(\al_b^*TY)\!\neq\!0$.\\

By the last two paragraphs,
$$\ft_{\fp}(\al)=\begin{cases}
\ft_{\fp}(\al_a)\!+\!\ft_{\fp}(\al_b),&\hbox{if}~\fp\!\in\!\cP^+(TY);\\
\ft_{\fp}(\al_a)\!+\!\ft_{\fp}(\al_b)\!+\!
\blr{w_1(Y),[\al_a]_{\Z_2}}\blr{w_1(Y),[\al_b]_{\Z_2}},&\hbox{if}~\fp\!\in\!\cP^-(TY).
\end{cases}$$
Combining this with~\eref{ShExSeqImm3_e10}, 
we obtain~\eref{ShExSeqImm3_e2} and~\eref{ShExSeqImm3_e3},
with $\al_1',\ldots,\al_{m'}'$ being loops obtained from $\al_1,\ldots,\al_m$
by smoothing one node or intersection point~$p$ of $\al_1,\ldots,\al_m$.
\end{proof}

\begin{rmk}\label{ShExSeqImm3_rmk}
Let $\al$ be as in the paragraph containing~\eref{ShExSeqImm3_e7b}.
By Examples~\ref{Pin1pMB_eg} and~\ref{Pin1mMB_eg}, 
the canonical $\Pin^{\pm}$-structure $\fp_0^{\pm}(\cN\io)$ on~$\cN\al$ 
is specified by the transition function satisfying
\hbox{$\wt{g}_{01}''|_{U_{01}}\!=\!\wt\bI_1$} and 
\hbox{$\wt{g}_{01}''|_{U_{10}}\!=\!\wt\bI_{1;1}$}.
Along with~\eref{wtioPindfn_e}, this implies that 
\hbox{$\wt\bI_{2;1}^{\st}\!=\!\wt\bI_{2;1}$}.
\end{rmk}

\section{Counts of stable real rational maps}
\label{g0realmaps_sec}

We next describe invariant counts of real rational curves in real symplectic manifolds 
and properties of these counts based on a perspective
originating in~\cite{Sol} and developed further in~\cite{Ge2,RealWDVV,RealWDVV3}.
The relevant \sf{real setting} for these invariants is described in Section~\ref{RealSetting_subs}.
We recall cases in which such invariants can be defined in Section~\ref{g0maps_subs}
and summarize their properties in Theorems~\ref{RGWs_thm} and~\ref{RGWs_thm2};
the implications of the two theorems are illustrated 
in Section~\ref{Wel6eg_subs}.
The relevant notation and orientations for 
the Deligne-Mumford spaces $\ov\cM_{k,l}^{\tau}$ of stable real curves
and for  the moduli spaces of stable real maps are specified in Sections~\ref{DM_subs}
and~\ref{MapSpaces_subs}, respectively.
In Section~\ref{MapSignDfn_subs}, 
we define the signs of the constrained maps appearing in the invariant counts
of Theorems~\ref{RGWs_thm} and~\ref{RGWs_thm2}, 
indicate why these counts are indeed invariant in the cases specified
in Section~\ref{g0maps_subs}, and complete the proof of Theorem~\ref{RGWs_thm}.
For generic choices of the constraints, these maps are immersions.
Proposition~\ref{cNsgn_prp}, established in Section~\ref{SignImmer_subs},  
interprets the sign definition of Section~\ref{MapSignDfn_subs} in 
terms of the normal bundles to these maps. 
This proposition is used in Section~\ref{SignImmer_subs} to complete the proof
of Theorem~\ref{RGWs_thm2} and in Section~\ref{WelSolPf_thm} 
to relate the moduli-theoretic definitions of the real curve signs arising 
from~\cite{Sol} with the intrinsic ones introduced in~\cite{Wel4,Wel6,Wel6b}.
Proposition~\ref{cNsgn_prp} is a consequence of Lemma~\ref{cNorient_lmm},
proved in Section~\ref{OrientImmer_subs}, which reformulates 
the definitions of the orientations at stable immersive maps 
in terms of the normal bundles to the immersions.

\subsection{The real setting}
\label{RealSetting_subs}

Let $(X,\om,\phi)$ be a compact connected 
\sf{real symplectic manifold}\gena{real symplectic manifold} 
of (real) dimension~$2n$,
i.e.~$(X,\om)$ is a symplectic manifold and $\phi$ is a diffeomorphism of~$X$
with itself so that $\phi^2\!=\!\id_X$ and $\phi^*\om\!=\!-\om$.
The fixed locus~$X^{\phi}$ of~$\phi$ is then a Lagrangian submanifold of~$(X,\om)$.
Define
$$H_2^{\phi}(X;\Z)=\big\{B\!\in\!H_2(X;\Z)\!:\phi_*(B)\!=\!-B\big\}
\nota{H2X\phi@$H_2^{\phi}(X;\Z)$}.$$

\vspace{.18in}

For a connected component~$Y$ of~$X^{\phi}$, let 
$$\prt_{Y;\Z_2}\!: H_2\big(X,Y;\Z\big)\lra H_1(Y;\Z) \lra H_1(Y;\Z_2)
\nota{dY2@$\prt_{Y;\Z_2}$}$$
be the composition of the boundary homomorphism of the homology relative exact sequence 
for the pair~$(X,Y)$ with the mod~2 reduction of the coefficients.
We denote~by
$$\fd_{X;Y}\!:H_2(X,Y;\Z)  \lra H_2^{\phi}(X;\Z)\nota{dY@$\fd_{X;Y}$}$$
the natural homomorphism
which glues each map $f\!:(\Si,\prt\Si)\!\lra\!(X,Y)$ from an oriented bordered surface
with the map $\phi\!\circ\!f$ from $(\Si,\prt\Si)$ with the opposite orientation; 
see \cite[Sec~1.1]{RealWDVV3}.
Define
\begin{gather}
\label{wtH2phidfn_e}
\wt{H}_2^{\phi}(X,Y)=
\big\{\big(\fd_{X;Y}(\be),\prt_{Y;\Z_2}(\be)\!\big)\!:\be\!\in\!H_2(X,Y;\Z)\big\}
\subset H_2^{\phi}(X;\Z)\!\oplus\!H_1(Y;\Z_2),\nota{H2Xphi2@$\wt{H}_2^{\phi}(X,Y)$}\\
\notag
\fd_Y\!:H_2(X;\Z)\lra \wt{H}_2^{\phi}(X;\Z),\quad \fd_Y(B)=\big(B\!-\!\phi_*B,0\big).
\nota{dY1@$\fd_Y$}
\end{gather}
If $(\Si,\si)$ is a symmetric surface as in Section~\ref{NRS_subs} with fixed locus~$\Si^{\si}$
and \hbox{$u\!:\Si\!\lra\!X$} is a continuous map such that \hbox{$u\!\circ\!\si=\phi\!\circ\!u$}
and $u(\Si^{\si})\!\subset\!Y$, then
$$\big(u_*[\Si]_{\Z},u_*[\Si^{\si}]_{\Z_2}\big)\in \wt{H}_2^{\phi}(X,Y).$$
By \cite[Proposition~4.1]{BHH},
$$ \blr{w_2(X),B}=\blr{w_1(Y),b} \qquad\forall\,(B,b)\in \wt{H}_2^{\phi}(X,Y)\,.$$
We also define
\begin{equation*}\begin{split}
\wh{H}_*(X,Y)
&=\bigoplus_{p\neq n-1}\!\!\!\!H_p(X;\Z)\!\oplus\!H_{n-1}(X\!-\!Y;\Z),
\nota{H2whX@$\wh{H}_*(X,Y)$,$\wh{H}^*(X,Y)$}\\
\wh{H}^*(X,Y)
&=\bigoplus_{p\neq n+1}\!\!\!\!H^p(X;\Z)\oplus\!\!H^{n+1}(X,Y;\Z).
\end{split}\end{equation*}
Let $\wh{H}^{2*}(X,Y)\!\subset\!\wh{H}^*(X,Y)$\nota{H2whX2@$\wh{H}^{2*}(X,Y)$} 
be the subspace of even-degree elements.\\

We denote by $\cJ_{\om}$ the space of $\om$-compatible (or -tamed) 
almost complex structures on~$X$ and by
$\cJ_{\om}^{\phi}\!\subset\!\cJ_{\om}$ the subspace of almost complex structures~$J$
such that \hbox{$\phi^*J\!=\!-J$}.
Let
$$c_1(X,\om)\equiv c_1(TX,J)\in H^2(X;\Z)\nota{c1@$c_1(X,\om)$}$$
be the first Chern class of $TX$ with respect to some $J\!\in\!\cJ_{\om}$;
it is independent of such a choice.
Define
$$\ell_{\om}\!:H_2(X;\Z)\lra\Z, \quad 
\ell_{\om}(B)=\blr{c_1(X,\om),B}\!+\!n\!-\!3\nota{l1@$\ell_{\om}$}.$$
For $\wt{B}\!\equiv\!(B,b)\in\wt{H}_2^{\phi}(X,Y)$, 
we define $\ell_{\om}(\wt{B})\!=\!\ell_{\om}(B)$.\\

For $J\!\in\!\cJ_{\om}$ and $B\!\in\!H_2(X;\Z)$, 
a subset $C\!\subset\!X$ is a \sf{genus~0} (or \sf{rational}) 
\sf{irreducible $J$-holomorphic degree~$B$ curve}\gena{rational curve!irreducible}
if there exists a simple (not multiply covered) \hbox{$J$-holomorphic} map 
\BE{Cudfn_e}u\!:\P^1\lra X \qquad\hbox{s.t.}\quad  C=u(\P^1),~~u_*[\P^1]=B.\EE
If $J\!\in\!\cJ_{\om}^{\phi}$, such a curve~$C$ is \sf{real}\gena{rational curve!real} 
if 
in addition \hbox{$\phi(C)\!=\!C$};
if so, then $B\!\in\!H_2^{\phi}(X;\Z)$.
We say that such a curve is of degree
\BE{wtBdfn_e}\wt{B}\!\equiv\!(B,b)\in \wt{H}_2^{\phi}(X,Y)\EE
if there exists a simple \hbox{$J$-holomorphic} map~$u$ as in~\eref{Cudfn_e} 
so~that 
\BE{Cudfn_e4}u\!\circ\!\tau= \phi\!\circ\!u \qquad\hbox{and}\qquad 
u_*\big([S^1]_{\Z_2}\big)=b\in H_1\big(Y;\Z_2\big),\EE
where $\tau$ is the standard conjugation on~$\P^1$ as in~\eref{tauetadfn_e}
and $S^1\!\subset\!\P^1$ is its fixed locus.\\

Let $k\!\in\!\Z^{\ge0}$ and 
$\bh\!\equiv\!(h_i\!:H_i\!\lra\!X)_{i\in[l]}$ be a tuple of cycles in~$X$,
i.e.~smooth maps from compact oriented manifolds, or more generally of pseudocycles
as in~\cite{pseudo}.
Let $\wt{B}\!\equiv\!(B,b)$ be an element of $\wt{H}_2^{\phi}(X,Y)$ such~that 
\BE{dimkbhcond_e}
\ell_{\om}(\wt{B})=(n\!-\!1)k\!+\!\sum_{i=1}^l\!\big(2n\!-\!\dim H_i\!-\!2\big).\EE
If $J\!\in\!\cJ_{\om}^{\phi}$ is generic, the elements of~$\bh$ are in general position,
and 
\BE{bpRdfn_e}\p\!\equiv\!\big(p_1^{\R},\ldots,p_k^{\R}\big)\in Y^k\EE
is a tuple of points in~$Y$ in general position, then
the space~$\fM_{\p;\bh}(\wt{B};J)$ of real rational irreducible $J$-holomorphic 
degree~$\wt{B}$ curves $C\!\subset\!X$ 
passing through the points $p_1^{\R},\ldots,p_k^{\R}$ and
meeting the images of $h_1,\ldots,h_l$ is a zero-dimensional manifold.\\

Suppose in addition $l^*\!\in\![l]$.
As described in Section~\ref{MapSignDfn_subs},
a relative $\Pin^{\pm}$-structure~$\fp$ on $Y\!\subset\!X$ 
determines a sign $\fs_{\fp;l^*;\bh}(C)$ for each element~$C$ of $\fM_{\p;\bh}(\wt{B};J)$
if 
\BE{BKcond_e0} k\!+\!2\Z\neq \blr{w_2(X),B}\in\Z_2.\EE
A relative $\OSpin$-structure~$\fp$ on~$Y$ determines
a sign~$\fs_{\fp;l^*;\bh}(C)$ without the assumption~\eref{BKcond_e0}.
For a $\Pin^{\pm}$- or~$\OSpin$-structure~$\fp$ on~$Y$, we~set
$$\fs_{k,l^*;\fp}(C)=\fs_{k,l^*;\io_X(\fp)}(C),$$
with $\io_X$ as in~\eref{vsSpinPin_e0}.

\subsection{Invariance and properties}
\label{g0maps_subs}

Let $(X,\om,\phi)$, $Y$, $k,l,l^*$, $\wt{B}$, $\bh$, $\p$, and $J$ be as at the end of 
Section~\ref{RealSetting_subs}.
If in addition $(X,\om,B)$ satisfies certain positivity conditions,
the discrete set $\fM_{\p;\bh}(\wt{B};J)$ is finite.
The~sum
\BE{mapcountdfn_e}N_{\wt{B};k,l^*;\bh}^{\phi;\fp}(Y)\equiv 
\sum_{C\in\fM_{\p;\bh}(\wt{B};J)}\!\!\!\!\!\!\!\!\!\!\!\fs_{\fp;l^*;\bh}(C)\EE
is well-defined and independent of the choices of~$J$, 
cycles~$h_i$ in their respective homology classes in~$\wh{H}_*(X,Y)$,
and the points~$p_i^{\R}$ in~$Y$ if 
\BE{Rcond_e0}\dim h_i\in2\Z, \qquad
\frac12\big(\!\dim\,h_i\big)~\begin{cases}\cong n,&\hbox{if}~i\!\in\![l]\!-\![l^*];\\
\not\cong n,&\hbox{if}~i\!\in\![l^*]\!-\!\{1\};
\end{cases}\quad\hbox{mod}~2,
\EE
and either 
\begin{enumerate}[label=(C\arabic*),leftmargin=*]

\item\label{Rcond_it1} 
$n\!=\!2$, $\fp$ is a relative $\Pin^-$-structure on $Y\!\subset\!X$ with 
$\lr{w_2(\fp),B'}\!=\!0$ for every $(B',b')\!\in\!\wt{H}_2^{\phi}(X,Y)$, or

\item\label{Rcond_it2} $n\!=\!2$, $\fp$ is a relative $\Pin^+$-structure on~$Y\!\subset\!X$ 
with $\lr{w_2(\fp),B'}\!=\!\lr{w_1(X),b'}$ for every $(B',b')\!\in\!\wt{H}_2^{\phi}(X,Y)$, or

\item\label{Rcond_it3} 
$n\!=\!3$, $\ell_{\om}(\wt{B})\!\neq\!0$, $\fp$ is an $\OSpin$-structure on~$Y$, 
and either $k\!\neq\!0$ or $b\!\neq\!0$, or

\item\label{Rcond_it4}  
$n\!=\!3$, $\ell_{\om}(\wt{B})\!\neq\!0$, $\fp$ is an $\OSpin$-structure on~$Y$,
and there exists a finite-order automorphism~$\psi$ of $(X,\om,\phi,Y)$ 
which restricts to an orientation-reversing diffeomorphism of~$Y$ and
acts trivially on~$\wt{H}_2^{\phi}(X,Y)$, or

\item\label{Rcond_it5}  
$n\!\not\in\!2\Z$, $k\!=\!0$, $l^*\!=\!l$, $\ell_{\om}(\wt{B})\!\neq\!0$, $b\!\neq\!0$, 
$\fp$ is an $\OSpin$-structure on~$Y$, and 
$\ell_{\om}(B')\!\in\!4\Z^{\ge0}$ for all spherical classes
$B'\!\in\!H_2^{\phi}(X;\Z)$ with $\om(B')\!>\!0$.

\end{enumerate}
The invariance in the first three cases above is established 
in the proof of Theorem~1.3 in~\cite{Sol};
see also the proofs of Theorems~1.6 and~1.7 in~\cite{Ge2}
and Proposition~5.2 in~\cite{RealWDVV}.
The proof of Proposition~1.3 in~\cite{RealWDVV3} establishes the invariance
in the fourth setting and applies in the third as well.
The last setting is a special case of that of Theorems~1.6 and~1.7 in~\cite{Ge2}.
We describe the underlining reason for the invariance in all cases 
in Section~\ref{InvarPf_subs}.\\

Under the condition~\eref{Rcond_e0}, each $l$-tuple~$\bh$ determines one or two possible  
values of~$l^*$.
Under any of the conditions \ref{Rcond_it1}-\ref{Rcond_it5}, 
the associated sign $\fs_{\fp;l^*;\bh}(C)$ is same for the admissible choices of~$l^*$
and can be defined for~$l\!=\!0$; see Section~\ref{MapSignDfn_subs}.
By Theorem~1.1 in~\cite{pseudo}, every homology class in a smooth manifold can be represented
by a pseudocycle.
In the five cases above, we thus obtain a multilinear symmetric functional 
\BE{RGWdfn_e}\blr{\cdot,\ldots,\cdot}_{\wt{B};Y;k}^{\phi,\fp}\!: 
\wh{H}^{2*}(X,Y)^{\oplus l}\lra\Z\EE
enumerating  real rational irreducible $J$-holomorphic 
degree~$\wt{B}$ curves $C\!\subset\!X$.\\
%

For $y\!\in\!Y$, we denote by~$S(\cN_yY)$ a small sphere 
in the fiber of the normal bundle~$\cN Y$ of~$Y$ in~$X$.
An $\OSpin$-structure~$\fp$ on~$Y$ determines an orientation
on~$S(\cN_yY)$; see Section~2.5 in~\cite{RealWDVV3}.
We denote the resulting homology class in~$H_{n-1}(X\!-\!Y;\Z)$
by~$[S(\cN_yY,\fp)]$ and the Lefschetz-Poincare dual 
of~$[S(\cN_yY,\fp)]$ in~$H^{n-1}(X,X\!-\!Y;\Z)$ by~$\ga_{X,Y}^{\circ}$.
The last class generates the kernel of the surjective homomorphism
$$H^4(X,Y;\Z)\lra H^4(X;\Z)$$
in the cohomology relative exact sequence for the pair $(X,Y)$.

\begin{thm}\label{RGWs_thm}
Suppose $(X,\om,\phi)$ is a compact connected real symplectic manifold of (real) dimension~$2n$,
$Y\!\subset\!X^{\phi}$ is a connected component,
$k\!\in\!\Z^{\ge0}$, $\wt{B}\!\in\!\wt{H}_2^{\phi}(X,Y)$, 
and $\fp$ is a relative $\Pin^{\pm}$- or $\OSpin$-structure on $Y\!\subset\!X^{\phi}$.
If one of the conditions~\ref{Rcond_it1}-\ref{Rcond_it5} above holds, 
then the multilinear functionals~\eref{RGWdfn_e} satisfy the following properties.
\begin{enumerate}[label=(RGW\arabic*),leftmargin=*]

\item\label{dimvan_it} $\lr{\ga_1,\ldots,\ga_l}_{\wt{B};Y;k}^{\phi,\fp}\!=\!0$ 
unless 
\eref{dimkbhcond_e} with $2n\!-\!\dim H_i$ replaced by $\deg\ga_i$ holds.

\item\label{dimvan_it2} If $\phi^*\ga_i\!=\!-(-1)^{(\deg\,\ga_i)/2}\ga_i$
for some $i\!\in\![l]$, then 
$\lr{\ga_1,\ldots,\ga_l}_{\wt{B};Y;k}^{\phi,\fp}\!=\!0$.

\item\label{deg0_it} If $\wt{B}\!=\!0$, 
$\displaystyle\lr{\ga_1,\ldots,\ga_l}_{\wt{B};Y;k}^{\phi,\fp}=
\begin{cases}
\lr{\ga_1,\pt},&\hbox{if}~(k,l)\!=\!(1,1);\\
0,&\hbox{otherwise}.\end{cases}$

\item\label{ins1_it} $\displaystyle\lr{1,\ga_2,\ldots,\ga_l}_{\wt{B};Y;k}^{\phi,\fp}=
\begin{cases}
1,&\hbox{if}~(\wt{B},k,l)\!=\!(0,1,1);\\
0,&\hbox{otherwise}.\end{cases}$

\item\label{div_it} If $\ga_0\!\in\!H^2(X;\Z)$, 
$\displaystyle\lr{\ga_0,\ga_1,\ldots,\ga_l}_{\wt{B};Y;k}^{\phi,\fp}=\lr{\ga_0,\wt{B}}
\lr{\ga_1,\ldots,\ga_l}_{\wt{B};Y;k}^{\phi,\fp}$.

\item\label{sphere_it} If $n\!=\!3$, 
$\displaystyle\lr{\ga_{X,Y}^{\circ},\ga_1,\ldots,\ga_l}_{\wt{B};Y;k}^{\phi,\fp}
=2\lr{\ga_1,\ldots,\ga_l}_{\wt{B};Y;k+1}^{\phi,\fp}$.

\item\label{OrientRev_it} If $n\!=\!3$ and 
$\ov\fp$ is the $\OSpin$-structure on $Y$ corresponding to~$\fp$ as in
the SpinPin~\ref{OSpinRev_prop} property in Section~\ref{SpinPinProp_subs}, then
$\displaystyle \lr{\cdot,\ldots,\cdot}_{\wt{B};Y;k}^{\phi,\ov\fp}
=-(-1)^k\lr{\cdot,\ldots,\cdot}_{\wt{B};Y;k}^{\phi,\fp}$.

\end{enumerate}
\end{thm}

\vspace{.1in}

The vanishing property~\ref{dimvan_it} holds because the dimensions of 
the relevant moduli spaces and the constraints are different unless 
\eref{dimkbhcond_e} with $2n\!-\!\dim H_i$ replaced by $\deg\ga_i$ holds.
This in turn implies the vanishing statement in~\ref{ins1_it}.
Since every degree~0 $J$-holomorphic map is constant,
\ref{dimvan_it} also leads to the vanishing statement in~\ref{deg0_it}.
We discuss the remaining properties of the real GW-invariants~\eref{RGWdfn_e}
stated in Theorem~\ref{RGWs_thm} in Section~\ref{MapSignDfn_subs},
after defining the sign~$\fs_{\fp;l^*;\bh}(C)$.\\

Let $(X',\om',\phi')$ and $(X'',\om'',\phi'')$ be 
compact connected real symplectic manifolds of dimensions~2 and~4, respectively,
$Y'\!\subset\!X'^{\phi'}$ and $Y''\!\subset\!X''^{\phi''}$
be connected components, 
$$(X,\om,\phi)\equiv(X',\om',\phi')\!\times\!(X'',\om'',\phi''), \quad\hbox{and}\quad
Y\equiv Y'\!\times\!Y''.$$
The restrictions $H^2(X',Y')\!\lra\!H^2(X')$ and $H^4(X'',Y'')\!\lra\!H^4(X'')$
are then isomorphisms.
Pairs
$$(\ga_1',\ga_1'')\in H^0(X')\!\oplus\!H^4(X') \qquad\hbox{and}\qquad 
(\ga_2',\ga_2'')\in H^2(X')\!\oplus\!H^2(X')$$
thus determine elements $\ga_1'\!\times\!\ga_1''$ and $\ga_2'\!\times\!\ga_2''$ of
$H^4(X,X'\!\times\!Y'')$ and $H^4(X,Y'\!\times\!X'')$, respectively.
Both restrict to elements of~$H^4(X,Y)$, which we denote in the same way.
Thus, the cross product induces a bilinear~map
$$\times\!:H^{2*}(X';\Z)\!\times\!\wh{H}^{2*}(X'';\Z)\lra  \wh{H}^{2*}(X,Y).$$

\begin{thm}\label{RGWs_thm2} 
Suppose $(X',\om',\phi')$, $(X'',\om'',\phi'')$, $Y',Y''$, $(X,\om,\phi)$, and~$Y$
are as above, $\os'$ and~$\os''$ are OSpin-structures
on~$Y'$ and~$Y''$, respectively, and
$$\ga_1',\ldots,\ga_l'\in H^{2*}(X';\Z), \qquad
\ga_1'',\ldots,\ga_l''\in \wh{H}^{2*}(X'';\Z).$$
Let $k\!\in\!\Z^{\ge0}$,
$$\wt{B}\equiv \big(\wt{B}',\wt{B}''\big)\in 
\wt{H}_2^{\phi'}\!\!(X',Y')\!\oplus\!
\wt{H}_2^{\phi''}\!\!(X'',Y'')\subset \wt{H}_2^{\phi}(X,Y),~~
\os=\llrr{\pi'^*\os',\pi''^*\os''}_{\oplus}\in\OSpin(Y),$$
where $\pi',\pi''\!:Y\!\lra\!Y',Y''$ are the two projections.
\begin{enumerate}[label=(RGW\arabic*),leftmargin=*]

\setcounter{enumi}{7}

\item\label{RGWprod_it0} If $\wt{B}'\!=\!0$, then
$${}\hspace{-.8in}\blr{\ga_1'\!\times\!\ga_1'',\ldots,
\ga_l'\!\times\!\ga_l''}_{\wt{B};Y;k}^{\phi,\os}
=\begin{cases}\!\!\left(\prod\limits_{i=1}^l\lr{\ga_i',\pt}\!\!\right)\!
\blr{\ga_1'',\ldots,\ga_l''}_{\wt{B}'';Y'';1}^{\phi,\os''},
&\hbox{if}~k\!=\!1,\,\deg\ga_i'\!=\!0~\forall\,i\!\in\![l];\\
0,&\hbox{otherwise}.\end{cases}$$

\item\label{RGWprod_it1} If $(\Si,\si)\!=\!(\P^1,\tau)$, 
$\wt{B}'\!=\!([\P^1]_{\Z},[S^1]_{\Z_2})$, 
and $\os'\!=\!\os_0(TS^1)$ is the standard $\OSpin$-structure for either orientation
on~$S^1$, then
$${}\hspace{-.8in}\blr{\ga_1'\!\times\!\ga_1'',\ldots,
\ga_l'\!\times\!\ga_l''}_{\wt{B};Y;k}^{\phi,\os}
=\begin{cases}0,&\hbox{if}~k\!+\!2\big|\!\{i\!\in\![l]\!:\deg\ga_i'\!=\!2\}\!\big|\!<\!3;\\
\blr{\ga_1'',\ldots,\ga_l''}_{\wt{B}'';Y'';3}^{\phi,\os''},&\hbox{if}~k\!=\!3,~
\ga_i'\!=\!1~\forall\,i\!\in\![l];\\
\lr{\ga_1',[\P^1]_{\Z}}
\blr{\ga_1'',\ldots,\ga_l''}_{\wt{B}'';Y'';1}^{\phi,\os''},&\hbox{if}~k\!=\!1,~
\ga_i'\!=\!1~\forall\,i\!\in\![l]\!-\!\{1\}.
\end{cases}$$
\end{enumerate}
\end{thm}

\vspace{.1in}

In the setting of Theorem~\ref{RGWs_thm2}, we can compute the invariants~\eref{RGWdfn_e}
using an almost complex structure $J\!=\!J'\!\times\!J''$ with $J'\!\in\!\cJ_{\om'}^{\phi'}$
and $J''\!\in\!\cJ_{\om''}^{\phi''}$.
If $\wt{B}'\!=\!0$, every real rational irreducible $J$-holomorphic 
degree~$\wt{B}$ curve $C\!\subset\!X$ is then of the form $y'\!\times\!C''$
for some $y'\!\in\!Y'$ and a real rational irreducible $J''$-holomorphic 
degree~$\wt{B}''$ curve $C''\!\subset\!X''$.
This implies the first statement in~\ref{RGWprod_it0}
as well as the second one up to~sign, which is confirmed in Section~\ref{MapSignDfn_subs}.
If $(\Si,\si)\!=\!(\P^1,\tau)$, we can take $J'$ to be the standard complex structure~$J_{\P^1}$ 
on~$\P^1$.
If $\wt{B}'\!=\!([\P^1]_{\Z},[S^1]_{\Z_2})$, every real rational irreducible $J$-holomorphic 
degree~$\wt{B}$ curve $C\!\subset\!X$ is then the graph of a real $J''$-holomorphic 
degree~$\wt{B}$ map from~$\P^1$ to~$X''$.
Since $(\P^1,\tau)$ with $k$~real points and $l'$~conjugate pairs of points
has a positive-dimensional family of automorphisms if $k\!+\!2l'\!<\!3$, 
this establishes the first statement in~\ref{RGWprod_it1}. 
Since any two configurations of $k$~real points and $l'$~conjugate pairs of points
on $(\P^1,\tau)$ are equivalent if $k\!+\!2l'\!=\!3$
and any such configuration has no automorphism,
we also obtain the second and third statements in~\ref{RGWprod_it1} up to~sign,
which is confirmed in Section~\ref{OrientImmer_subs}
as a corollary of Proposition~\ref{cNsgn_prp}.\\

If $n\!=\!2$, (the tangent bundle~of) $Y$ admits a $\Pin^-$-structure;
see the SpinPin~\ref{SpinPinObs_prop}\ref{PinObs_it} property on page~\pageref{SpinPinObs_prop}
and~\eref{WuRel2_e}.
By the RelSpinPin~\ref{RelSpinPinObs_prop}\ref{RelPinObs_it} property 
on page~\pageref{RelSpinPinObs_prop} and~\eref{WuRel2_e},
\hbox{$Y\!\subset\!X$} also admits a relative $\Pin^+$-structure~$\fp$ 
with $w_2(\fp)\!=\!w_2(X)$.
In particular, the invariants~\eref{RGWdfn_e} can be defined for 
every compact real symplectic fourfold.

\subsection{Moduli spaces of real curves}
\label{DM_subs}

Let $\PSL_2^{\tau}\C\!\subset\!\PSL_2\C$ be the subgroup
of automorphisms of~$\P^1$ commuting with the conjugation~$\tau$ on~$\P^1$.
For $k\!\in\!\Z^{\ge0}$, let $[k]\!=\!\{1,\ldots,k\}$ as before.
For $k,l\!\in\!\Z^{\ge0}$ with $k\!+\!2l\!\ge\!3$, 
we denote by~$\cM_{k,l}^{\tau}$ the Deligne-Mumford moduli space of 
equivalence classes of smooth real genus~0 curves
with separating fixed locus, $k$~real marked points, and $l$~conjugate pairs of marked points.
Thus, 
\begin{equation*}\begin{split}
\cM_{k,l}^{\tau}\approx
\big\{\!\big((x_i)_{i\in[k]},(z_i^+,z_i^-)_{i\in[l]}\big)\!:\,
x_i\!\in\!S^1,\,z_i^{\pm}\!\in\!\P^1\!-\!S^1,\,z_i^+\!=\!\tau(z_i^-),&\\
x_i\!\neq\!x_j,\,z_i^+\!\neq\!z_j^+,z_j^-~\forall\,i\!\neq\!j&\big\}\big/\PSL_2^{\tau}\C.
\end{split}\end{equation*}
The elements of~$\cM_{k,l}^{\tau}$ are the equivalence classes of the marked symmetric Riemann
surfaces
\BE{cCdfn0_e}\cC\equiv 
\big(\Si\!=\!\P^1,(x_i)_{i\in[k]},(z_i^+,z_i^-)_{i\in[l]},\si\!=\!\tau\big)\EE
modulo the reparametrizations by~$\PSL_2^{\tau}\C$.
This space is a smooth manifold of dimension $k\!+\!2l\!-\!3$.
Its topological components 
are indexed by the possible distributions of the points~$z_i^+$ between the interiors 
of the two disks cut out by the fixed locus~$S^1$ of 
the standard involution~$\tau$ on~$\P^1$ and by the orderings of 
the real marked points~$x_i$ on~$S^1$.\\

If $k\!+\!2l\!\ge\!4$ and $i\!\in\![k]$, let
\BE{ff1dfn_e}\ff_{k,l;i}^{\R}\!:\cM_{k,l}^{\tau}\lra\cM_{k-1,l}^{\tau}\EE
be the forgetful morphism dropping the $i$-th real marked point.
The associated exact sequence 
$$0\lra \ker\nd\ff_{k,l;i}^{\R}\lra T\cM_{k,l}^{\tau} 
\xlra{\nd\ff_{k,l;i}^{\R}} \ff_{k,l;i}^{\R\,*}T\cM_{k-1,l}^{\tau} \lra0 $$
induces an isomorphism
\BE{cMorientR_e}
\la\big(T\cM_{k,l}^{\tau}\big)\approx 
\ff_{k,l;i}^{\R\,*}\la\big(T\cM_{k-1,l}^{\tau}\big)
\!\otimes\!\big(\!\ker\nd\ff_{k,l;i}^{\R}\big)\,.\EE
If $k\!+\!2l\!\ge\!5$ and $i\!\in\![l]$, we similarly denote~by
\BE{ff2dfn_e}\ff_{k,l;i}\!:\cM_{k,l}^{\tau}\lra\cM_{k,l-1}^{\tau}\EE
the forgetful morphism dropping the $i$-th conjugate pair of marked points.
It induces an isomorphism
\BE{cMorientC_e}
\la\big(T\cM_{k,l}^{\tau}\big)\approx 
\ff_{k,l;i}^{\,*}\la\big(T\cM_{k,l-1}^{\tau}\big)
\!\otimes\!\la\big(\!\ker\nd\ff_{k,l;i}\big)\,.\EE
For each $\cC\!\in\!\cM_{k,l}$ as in~\eref{cCdfn0_e},
$$\ker\nd_{\cC}\ff_{k,l;i}\approx T_{z_i^+}\P^1 $$
is canonically oriented by the complex orientation of the fiber~$\P^1$ at~$z_i^+$.
We denote the resulting orientation of the last factor in~\eref{cMorientC_e} by~$\fo_i^+$.\\

Suppose $l\!\in\!\Z^+$ and $\cC\!\in\!\cM_{k,l}^{\tau}$ is as in~\eref{cCdfn0_e} with $\Si\!=\!\P^1$.
Let $\bD^2_+\!\subset\!\C\!\subset\!\P^1$ be the disk cut out by 
the fixed locus~$S^1$ of~$\tau$ which contains~$z_1^+$. 
We orient \hbox{$S^1\!\subset\!\bD_+^2\!\subset\!\C$} in the standard way.
If $k\!+\!2l\!\ge\!4$ and $i\!\in\![k]$,
this determines an orientation~$\fo_i^{\R}$ of the fiber
$$\ker\nd_{\cC}\ff_{k,l;i}^{\R}\approx T_{x_i}S^1 $$
of the last factor in~\eref{cMorientR_e} over~$\ff_{k,l;i}^{\R}(\cC)$.\\

Let $(x_1,x_{j_2(\cC)},\ldots,x_{j_k(\cC)})$ be the ordering of the real marked points of~$\cC$
starting with~$x_1$ and going in the direction of the standard orientation of~$S^1$.
We denote by $\de_{\R}(\cC)\!\in\!\Z_2$ the sign of the permutation sending
$$\vp_{\cC}\!:\big\{2,\ldots,k\big\}\lra \big\{2,\ldots,k\big\}, \quad
\vp_{\cC}(i)=j_i(\cC)\,.$$
If $k\!=\!0$, we take $\de_{\R}(\cC)\!=\!0$.
For $l^*\!\in\![l]$, let 
\BE{decldfn_e}\de_{l^*}^c(\cC)=
\big|\big\{i\!\in\![l]\!-\![l^*]\!:z_i^+\!\not\in\!\bD_+^2\big\}\big|+2\Z\in\Z_2.\EE
In particular, $\de_{\R}(\cC)\!=\!0$ if $k\!\le\!2$ and $\de_l^c(\cC)\!=\!0$.
The functions~$\de_{\R}$ and~$\de_{l^*}$ are locally constant on~$\cM_{k,l}^{\tau}$.\\

The space $\cM_{1,1}^{\tau}$ is a single point;
we take $\fo_{1,1}\!\equiv\!+1$ to be its orientation as a plus point.
We identify the one-dimensional space 
$\cM_{0,2}^{\tau}$ with $\R^+\!-\!\{1\}$ via the cross ratio 
\BE{cM02ident_e}
\vph\!:\cM_{0,2}^{\tau}\lra\R^+\!-\!\{1\}, ~~
\vph\big([(z_1^+,z_1^-),(z_2^+,z_2^-)]\big)= 
\frac{z_2^+\!-\!z_1^-}{z_2^-\!-\!z_1^-}:\frac{z_2^+\!-\!z_1^+}{z_2^-\!-\!z_1^+}
=\frac{|1\!-\!z_1^+/z_2^-|^2}{|z_1^+\!-\!z_2^+|^2}\,;\EE
see Figure~\ref{M02_fig}.
This identification, which is the {\it opposite} of \cite[(3.1)]{RealEnum} and 
\cite[(1.12)]{RealGWsII}, determines an orientation~$\fo_{0,2}$ on~$\cM_{0,2}^{\tau}$.\\

\begin{figure}
\begin{pspicture}(-3.2,-1.5)(10,2)
\psset{unit=.4cm}
\pscircle[linewidth=.1](-3,0){1}\pscircle[linewidth=.1](-3,2){1}
\pscircle*(-3,1){.2}\pscircle*(-3.7,-.7){.2}\pscircle*(-2.3,-.7){.2}
\pscircle*(-3.7,2.7){.2}\pscircle*(-2.3,2.7){.2}
\rput(-3.6,-1.5){\sm{$1^-$}}\rput(-1.9,-1.5){\sm{$2^+$}}
\rput(-3.6,3.5){\sm{$1^+$}}\rput(-1.9,3.5){\sm{$2^-$}}
\psline[linewidth=.1,linestyle=dotted](-1,1)(1,1)
\pscircle[linewidth=.1](3.5,1){1.5}
\psarc[linewidth=.05](3.5,3.6){3}{240}{300}
\psarc[linewidth=.05,linestyle=dashed](3.5,-1.6){3}{60}{120}
\pscircle*(2.44,-.06){.2}\pscircle*(4.56,-.06){.2}\pscircle*(2.44,2.06){.2}\pscircle*(4.56,2.06){.2}
\rput(2.3,-.9){\sm{$1^-$}}\rput(4.7,-.9){\sm{$2^+$}}
\rput(2.3,2.9){\sm{$1^+$}}\rput(4.7,2.9){\sm{$2^-$}}
\psline[linewidth=.1,linestyle=dotted](6.5,1)(8.5,1)
\pscircle[linewidth=.1](11,1){1}\pscircle[linewidth=.1](13,1){1}
\pscircle*(12,1){.2}\pscircle*(10.3,.3){.2}\pscircle*(13.7,.3){.2}
\pscircle*(10.3,1.7){.2}\pscircle*(13.7,1.7){.2}
\rput(10.16,-.54){\sm{$1^-$}}\rput(13.84,-.54){\sm{$2^-$}}
\rput(10.16,2.54){\sm{$1^+$}}\rput(13.84,2.54){\sm{$2^+$}}
\psline[linewidth=.1,linestyle=dotted](15.5,1)(17.5,1)
\pscircle[linewidth=.1](20.5,1){1.5}
\psarc[linewidth=.05](20.5,3.6){3}{240}{300}
\psarc[linewidth=.05,linestyle=dashed](20.5,-1.6){3}{60}{120}
\pscircle*(19.44,-.06){.2}\pscircle*(21.56,-.06){.2}\pscircle*(19.44,2.06){.2}\pscircle*(21.56,2.06){.2}
\rput(19.3,-.9){\sm{$1^-$}}\rput(21.7,-.9){\sm{$2^-$}}
\rput(19.3,2.9){\sm{$1^+$}}\rput(21.7,2.9){\sm{$2^+$}}
\psline[linewidth=.1,linestyle=dotted](23,1)(25,1)
\pscircle[linewidth=.1](27,0){1}\pscircle[linewidth=.1](27,2){1}
\pscircle*(27,1){.2}\pscircle*(26.3,-.7){.2}\pscircle*(27.7,-.7){.2}
\pscircle*(26.3,2.7){.2}\pscircle*(27.7,2.7){.2}
\rput(26.4,-1.5){\sm{$1^-$}}\rput(28.1,-1.5){\sm{$2^-$}}
\rput(26.4,3.5){\sm{$1^+$}}\rput(28.1,3.5){\sm{$2^+$}}
\psline[linewidth=.1](-3,-3)(27,-3)
\pscircle*(-3,-3){.2}\pscircle*(12,-3){.2}\pscircle*(27,-3){.2}
\rput(-3,-3.8){0}\rput(12,-3.8){1}\rput(27,-3.8){$\i$}
\rput(-3,-3.8){0}\rput(12,-3.8){1}\rput(27,-3.8){$\i$}
\end{pspicture}
\caption{The structure of the Deligne-Mumford compactification $\ov\cM_{0,2}^{\tau}$
of~$\cM_{0,2}^{\tau}$.}
\label{M02_fig}
\end{figure}

We now define an orientation~$\fo_{k,l}$ on $\cM_{k,l}^{\tau}$ for $l\!\in\!\Z^+$ and
$k\!+\!l\!\ge\!3$ inductively.
If $k\!\ge\!1$, we take~$\fo_{k,l}$ to be so that the $i\!=\!k$ case 
of the isomorphism~\eref{cMorientR_e} is compatible with the orientations~$\fo_{k,l}$,
$\fo_{k-1,l}$, and~$\fo_k^{\R}$ on the three line bundles involved.
If $l\!\ge\!2$, we take~$\fo_{k,l}$ to be so that the $i\!=\!l$ case 
of the isomorphism~\eref{cMorientC_e} is compatible with the orientations~$\fo_{k,l}$,
$\fo_{k,l-1}$, and~$\fo_l^+$.
By a direct check, the orientations on~$\cM_{1,2}^{\tau}$ induced
from~$\cM_{0,2}^{\tau}$ via~\eref{cMorientR_e} and~$\cM_{1,1}^{\tau}$ via~\eref{cMorientC_e} are the same.
Since the fibers of $\ff_{k,l;l}|_{\cM_{k,l}^{\tau}}$ are even-dimensional,
it follows that the orientation~$\fo_{k,l}$ on~$\cM_{k,l}^{\tau}$ is well-defined
for all $l\!\in\!\Z^+$ and $k\!\in\!\Z^{\ge0}$ with $k\!+\!2l\!\ge\!3$.
This orientation is as above \cite[Lemma~5.4]{Ge2};
for $k\!=\!0$, this orientation is the opposite of that taken 
in \cite[Section~3]{RealEnum} and \cite[Section~1.4]{RealGWsII}.\\

For $l^*\!\in\![l]$, we denote by $\fo_{k,l;l^*}$ the orientation on 
$\cM_{k,l}^{\tau}$ which equals~$\fo_{k,l}$ at $[\cC]\!\in\!\cM_{k,l}^{\tau}$ 
if and only if $\de_{\R}(\cC)\!=\!\de_{l^*}^c(\cC)$.
If $l\!\ge\!2$ and $\cC$ is as in~\eref{cCdfn0_e} with $z_1^+\!=\!0$ and $z_2^+\!\in\!\R$,
then the natural isomorphism
\BE{cMorientisom_e}
T_{[\cC]}\cM^\tau_{k,l}\approx T_{z^+_2}\R^+\!\oplus\!
\big(T_{z^+_3}\P^1\!\oplus\!\ldots\!\oplus\!T_{z^+_l}\P^1
\!\oplus\! T_{x_1}S^1\!\oplus\!\ldots\!\oplus\!T_{x_k}S^1\big)\EE
has sign $(-1)^{\de_{\R}(\cC)+\de_{l^*}^c(\cC)}$
with respect to the orientations $\fo_{k,l;l^*}$ on $T_{[\cC]}\cM^\tau_{k,l}$,
the opposite~$\ov{\fo_{\R}}$ of the standard orientation on $T_{z^+_2}\R^+$, 
and the standard orientations on~$T_{z^+_i}\P^1$ and~$T_{x_i}S^1$.
If $l\!=\!1$, the same statement holds with $T_{z^+_2}\R^+,T_{z^+_i}\P^1,T_{x_1}S^1$ 
dropped and $\de_{l^*}^c(\cC)\!=\!0$.
The next lemma is straightforward.

\begin{lmm}\label{cMorient_lmm}
The orientations $\fo_{k,l;l^*}$ on $\cM_{k,l}^{\tau}$ with	
$k,l\!\in\!\Z^{\ge0}$ and $l^*\!\in\![l]$ such that \hbox{$k\!+\!2l\!\ge\!3$} and
satisfy the following properties:
\begin{enumerate}[label=($\fo_{\cM}\arabic*$),leftmargin=*]


\item\label{cMorientC_it} 
the isomorphism~\eref{cMorientC_e} with $(l,i)$ replaced by $(l\!+\!1,l^*\!+\!1)$
respects the orientations $\fo_{k,l+1;l^*+1}$, $\fo_{k,l;l^*}$, and~$\fo_{l^*+1}^+$; 

\item $\fo_{k,l;l^*}$  is preserved by the interchange of 
two real points $x_i$ and $x_j$ with $2\!\le\!i,j\!\le\!k$;

\item\label{cM1Rch_it} if $2\!\le\!i\!\le\!k$ and $\cC\!\in\!\cM_{k,l}^{\tau}$,
$\fo_{k,l;l^*}$ is preserved at~$\cC$  by the interchange of 
the real points $x_1$ and $x_{j_i(\cC)}$ with $2\!\le\!i\!\le\!k$
if and only if $(k\!-\!1)(i\!-\!1)\!\in\!2\Z$;

\item\label{cM1cij_it} if $1\!\le\!i,j\!\le\!l$, $\cC\!\in\!\cM_{k,l}^{\tau}$,
and the marked points $z_i^+$ and~$z_j^+$ of~$\cC$ lie on the same connected component
of $\P^1\!-\!S^1$, then $\fo_{k,l;l^*}$ is preserved at~$\cC$  by the interchange of 
conjugate pairs $(z_i^+,z_i^-)$ and $(z_j^+,z_j^-)$;

\item if $l^*\!<\!i\!\le\!l$ (resp.~$1^*\!<\!i\!\le\!l^*$),
$\fo_{k,l;l^*}$ is preserved (resp.~reversed) by
the interchange of the points in a conjugate pair $(z_i^+,z_i^-)$;

\item if $\cC\!\in\!\cM_{k,l}^{\tau}$,
$\fo_{k,l;l^*}$ is preserved at~$\cC$  by the interchange 
of the points in the conjugate pair $(z_1^+,z_1^-)$ if and only~if 
$$k\neq0~~\hbox{and}~~l\!-\!l^*\cong\binom{k}{2}~\tn{mod}~2 
\quad\hbox{or}\quad
k=0~~\hbox{and}~~l\!-\!l^*\cong1~\tn{mod}~2.$$ 

\end{enumerate}
\end{lmm}

\subsection{Moduli spaces of real maps}
\label{MapSpaces_subs}

We denote by~$\fj$ the standard complex structure on~$\P^1$.
Suppose  $(X,\om,\phi)$ is a real symplectic manifold and $k,l\!\in\!\Z^{\ge0}$.
If $\wt{B}\!\equiv\!(B,b)\!\in\!\wt{H}_2^{\phi}(X,Y)$ and $J\!\in\!\cJ_{\om}^{\phi}$,
a \sf{real $J$-holomorphic degree~$\wt{B}$ map} from~$(\P^1,\tau)$ is 
a smooth map \hbox{$u\!:\P^1\!\lra\!X$} so that 
$$u\!\circ\!\tau=\phi\!\circ\!u, \quad u_*[\P^1]=B, \quad u_*\big([S^1]_{\Z}\big)=b, \quad
\bp_Ju|_z\!\equiv\!\frac12\big(\nd_zu\!+\!J\!\circ\!\nd_zu\!\circ\!\fj\big)
=0~~\forall~z\!\in\!\P^1.$$
Such a map is called \sf{simple}\gena{simple map} if $u\!\neq\!u'\!\circ\!h$ for any 
holomorphic map \hbox{$h\!:\P^1\!\lra\!\P^1$} of degree larger than~1.
We denote by $\fM_{k,l}(\wt{B};J)$ the moduli space of equivalence classes
of real degree~$\wt{B}$ $J$-holomorphic maps from~$(\P^1,\tau)$
with $k$~real marked points and $l$~conjugate pairs of marked points
modulo reparametrizations by~$\PSL_2^{\tau}\C$. 
The elements of $\fM_{k,l}(\wt{B};J)$ are the equivalence classes of tuples
\BE{udfn_e} \u=\big(u\!:\P^1\!\lra\!X,(x_i)_{i\in[k]},(z_i^+,z_i^-)_{i\in[l]},\tau\big)\EE
such that 
\BE{cCdfn_e}\cC\equiv \big(\P^1,(x_i)_{i\in[k]},(z_i^+,z_i^-)_{i\in[l]},\tau\big)\EE
is a marked symmetric Riemann surface with
$k$~real marked points and $l$~conjugate pairs of marked points
and $u$ is a real $J$-holomorphic degree~$\wt{B}$ map from~$(\P^1,\tau)$.
A marked map $\u$ as in~\eref{udfn_e} is called \sf{simple}\gena{simple marked map} 
if $u$ is simple.
We denote~by
$$\fM_{k,l}^*(\wt{B};J)\subset\fM_{k,l}(\wt{B};J)$$
the subspace of the equivalence classes of simple maps.
For a generic $J\!\in\!\cJ_{\om}^{\phi}$, this subspace is 
a smooth manifold of dimension~$\ell_{\om}(\wt{B})$.\\

Suppose $B\!\neq\!0$ or $k\!+\!2l\!\ge\!5$.
If $i\!\in\![l]$, let 
\BE{ffMdfn_e}\ff_{k,l;i}\!:\fM_{k,l}^*(\wt{B};J)\lra \fM_{k,l-1}^*(\wt{B};J)\EE
be the forgetful morphism dropping the $i$-th conjugate pair of marked points.
Similarly to~\eref{cMorientC_e}, it induces an isomorphism
\BE{cMorientC2_e}
\la\big(T\fM_{k,l}^*(\wt{B};J)\!\big)\approx 
\ff_{k,l;i}^{\,*}\la\big(T\fM_{k,l-1}^*(\wt{B};J)\!\big)
\!\otimes\!\la\big(\!\ker\nd\ff_{k,l;i}\big)\,.\EE
For each $\u\!\in\!\fM_{k,l}^*(\wt{B};J)$ as in~\eref{udfn_e},
$$\ker\nd_{\u}\ff_{k,l;i}\approx T_{z_i^+}\P^1 $$
is canonically oriented by the complex orientation of the fiber~$\P^1$ at~$z_i^+$.
We again denote the resulting orientation of the last factor in~\eref{cMorientC2_e} by~$\fo_i^+$.\\

For each $i\!\in\![k]$, let 
$$\ev_i^{\R}\!:\fM_{k,l}^*(\wt{B};J)\lra Y,\quad
\ev_i^{\R}\big([u,(x_i)_{i\in[k]},(z_i^+,z_i^-)_{i\in[l]},\tau]\big)=u(x_i),$$
be the evaluation morphism for the $i$-th real marked point.
For each $i\!\in\![l]$, let 
$$\ev_i^+\!:\fM_{k,l}^*(\wt{B};J,\nu)\lra X,\quad
\ev_i^+\big([u,(x_i)_{i\in[k]},(z_i^+,z_i^-)_{i\in[l]},\tau]\big)=u(z_i^+),$$
be the evaluation morphism for the positive point of the $i$-th conjugate pair of marked points.
We define
\BE{fMevdfn_e}\ev\!\equiv\!\prod_{i=1}^k\!\ev_i^{\R}\times\!\prod_{i=1}^l\!\ev_i^+\!:
\fM_{k,l}^*(\wt{B};J)\lra Y^k\!\times\!X^l\EE
to be the \sf{total evaluation} map.\\

For a real $J$-holomorphic map~$u$ from~$(\P^1,\tau)$, let 
\begin{equation*}\begin{split}
D_{J;u}^{\phi}\!:  \Ga(u)
\equiv&\big\{\xi\!\in\!\Ga(\P^1;u^*TX)\!:\,\xi\!\circ\!\tau\!=\!\nd\phi\!\circ\!\xi\big\}\\
&\lra
\Ga^{0,1}(u)\equiv
\big\{\ze\!\in\!\Ga(\P^1;(T^*\P^1,\fj)^{0,1}\!\otimes_{\C}\!u^*(TX,J))\!:\,
\ze\!\circ\!\nd\tau=\nd\phi\!\circ\!\ze\big\}
\end{split}\end{equation*}
be the linearization of the $\bp_J$-operator on the space of 
real maps from~$(\P^1,\tau)$.
This is a real CR-operator on the real bundle pair $u^*(TX,\nd\phi)$ over~$(\P^1,\tau)$
in the terminology of Section~\ref{NRS_subs}.
For a tuple~$\u$ as in~\eref{udfn_e} so that \hbox{$[\u]\!\in\!\fM_{k,l}(\wt{B};J)$},
define
\begin{gather*}
\la_{\u}^{\R}(X)=\bigotimes_{i=1}^k\la\big(T_{u(x_i)}Y\big), \qquad
\la_{\u}^{\C}(X)=\bigotimes_{i=1}^l\la\big(T_{u(z_i^+)}X\big),\\
\wt\la_{\u}\big(D_J^{\phi},X\big)=
\la_{\u}^{\R}(X)^*\!\otimes\!\la_{\u}^{\C}(X)^*\!\otimes\!\la\big(D_{J;u}^{\phi}\big).
\end{gather*}
For a generic $J\!\in\!\cJ_{\om}^{\phi}$, the real CR-operator $D_{J;u}^{\phi}$ is surjective 
for every $[\u]\!\in\!\fM_{k,l}^*(\wt{B};J)$ and thus
$$\la\big(D_{J;u}^{\phi}\big)=\la\big(\!\ker D_{J;u}^{\phi}\big).$$

\vspace{.15in}

Suppose $l\!\in\!\Z^+$.
Given a tuple~$\u$ as in~\eref{udfn_e}, 
we choose $\bD^2_+\!\subset\!\P^1$ to be the half-disk cut out 
by the fixed locus $S^1\!\subset\!\P^1$ of~$\tau$ so that $z_1^+\!\in\!\bD^2_+$.
A relative $\Pin^{\pm}$-structure (resp.~$\OSpin$-structure)~$\fp$ 
on $Y\!\subset\!X$ 
pulls back to a relative $\Pin^{\pm}$-structure
(resp.~$\OSpin$-structure)~$u^*\fp$ on the fixed locus of
the real bundle pair $u^*(TX,\nd\phi)$ over~$(\P^1,\tau)$.
If~\eref{BKcond_e0} holds,
a relative $\Pin^{\pm}$-structure~$\fp$ on $Y\!\subset\!X$ thus determines an orientation
on $\la_{\u}^{\R}(X)^*\!\otimes\!\la(D_{J;u}^{\phi})$
as in Theorem~\ref{CROrient_thm}\ref{CROrientPin_it}.
Along with the orientation of $\la_{\u}^{\C}(X)$ determined 
by the symplectic orientation~$\fo_{\om}$ of~$(X,\om)$,
it thus determines an orientation 
$$\fo_{\fp;\u}^D\equiv \fo_{\cC;\fp}(u^*TX,u^*\nd\phi)$$
of the line $\wt\la_{\u}(D_J^{\phi},X)$ 
for each tuple~$\u$ as in~\eref{udfn_e}  
so that \hbox{$[\u]\!\in\!\fM_{k,l}(\wt{B};J)$};
this orientation varies continuously with~$\u$.
A relative $\OSpin$-structure~$\fp$ on $Y\!\subset\!X$
similarly determines an orientation $\la_{\cC}^{\R}(\os)\!\equiv\!\la_{\cC}^{\R}(\fo)$
on $\la_{\u}^{\R}(X)$
as above the CROrient~\ref{CROPin2SpinRed_prop} property in Section~\ref{OrientPrp_subs1}, 
an orientation on $\la(D_{J;u}^{\phi})$ as in Theorem~\ref{CROrient_thm}\ref{CROrientSpin_it},
and an orientation~$\fo_{\fp;\u}^D$ on $\wt\la_{\u}(D_J^{\phi},X)$ varying continuously
with \hbox{$[\u]\!\in\!\fM_{k,l}(\wt{B};J)$}.\\

If $k\!+\!2l\!\ge\!3$, let 
$$\ff_{k,l}\!:\fM_{k,l}^*(\wt{B};J)\lra \cM_{k,l}^{\tau}$$
be the forgetful morphism dropping the map component~$u$ from each tuple~$\u$ as
in~\eref{udfn_e}.
If $J\!\in\!\cJ_{\om}^{\phi}$ is generic and $\u\!\in\!\fM_{k,l}^*(\wt{B};J)$, 
the associated exact sequence
$$0\lra \ker D_{J;u}^{\phi}\!=\!\ker\nd_{\u}\ff_{k,l}\lra T_{\u}\fM_{k,l}^*(\wt{B};J)
\xlra{\nd_{\u}\ff_{k,l}} T_{\ff_{k,l}(\u)}\cM_{k,l}^{\tau} \lra0 $$
of vector spaces induces an isomorphism
\BE{OrientSubs_e3}\begin{split}
\la_{\u}(\ev)&\equiv
\la_{\u}^{\R}(X)^*\!\otimes\!\la_{\u}^{\C}(X)^*\!\otimes\!
\la_{\u}\big(\fM_{k,l}^*(\wt{B};J)\big)
\approx \wt\la_{\u}\big(D_J^{\phi},X\big)\!\otimes\!
\la_{\ff_{k,l}(\u)}\big(\cM_{k,l}^{\tau}\big).
\end{split}\EE
By the previous paragraph, a relative $\Pin^{\pm}$-structure~$\fp$ 
(resp.~$\OSpin$-structure~$\fp$) on $Y\!\subset\!X$
determines an orientation~$\fo_{\fp;\u}^D$ on $\wt\la_{\u}(D_J^{\phi},X)$ 
if~\eref{BKcond_e0} holds  (resp.~$\fo_{\fp;\u}^D$ on $\wt\la_{\u}\big(D_J^{\phi},X)$).
For each $l^*\!\in\![l]$, $\fo_{\fp;\u}^D$ 
and the orientation~$\fo_{k,l;l^*}$ 
of Lemma~\ref{cMorient_lmm} on the second factor in~\eref{OrientSubs_e3} in turn determine
an orientation~$\fo_{\fp;l^*;\u}$ (resp.~$\fo_{\os;l^*;\u}$) on~$\la_{\u}(\ev)$;
this orientation varies continuously with~$\u$.\\

If $\u'$ is a marked map obtained from $\u\!\in\!\fM_{k,l}^*(\wt{B};J)$ by interchanging 
the points in a conjugate pair~$(z_i^+,z_i^-)$, 
we identify $\la_{\u}^{\C}(X)$ with $\la_{\u'}^{\C}(X)$ via 
the isomorphism $\nd_{u(z_i^+)}\phi$ in the $i$-th factor.
This identification is orientation-preserving with respect to~$\fo_{\om}$
if and only if $n\!\in\!2\Z$.
For $i\!\in\![k]$ and $\u\!\in\!\fM_{k,l}^*(\wt{B};J)$ with the associated marked curve~$\cC$
as in~\eref{cCdfn_e}, let
$$j_i(\u)=j_i(\cC)\in[k]$$
be as in Section~\ref{DM_subs}.
We set $\binom{-1}{2}\!\equiv\!0$ as before.
The next two lemmas follow immediately from Lemma~\ref{cMorient_lmm}
and the CROrient~\ref{CROos_prop}, \ref{CROp_prop},  and~\ref{CROSpinPinStr_prop} properties 
in Section~\ref{OrientPrp_subs1}.

\begin{lmm}\label{orient_lmm}
Suppose $(X,\om,\phi)$ is a real symplectic manifold of dimension~$2n$, 
$Y\!\subset\!X^{\phi}$ is a topological component,
$k,l\!\in\!\Z^{\ge0}$ with \hbox{$k\!+\!2l\!\ge\!3$}, $l^*\!\in\![l]$, 
$\wt{B}\!\in\!\wt{H}_2^{\phi}(X,Y)$, and $J\!\in\!\cJ_{\om}^{\phi}$ is generic.
If~$k$ and $\wt{B}\!\equiv\!(B,b)$ satisfy~\eref{BKcond_e0}, 
then every  relative $\Pin^{\pm}$-structure~$\fp$ on $Y\!\subset\!X$ 
determines a relative orientation~$\fo_{\fp;l^*}$ the map~\eref{fMevdfn_e}  
with the following properties:
\begin{enumerate}[label=($\fo_{\fp}\arabic*$),leftmargin=*]

\setcounter{enumi}{-1}

\item\label{etaMorient_it} 
if $\eta\!\in\!H^2(X,Y;\Z_2)$ and $\u\!\in\!\fM_{k,l}^*(\wt{B};J)$,
the orientations $\fo_{\fp;l^*}$ and $\fo_{\eta\cdot\fp;l^*}$ at~$\u$
are the same if and only if $\lr{\eta,u_*[\bD^2_+]_{\Z_2}}\!=\!0$;

\item\label{fforient_it} the orientations $\fo_{\fp;l^*+1}\fo_{\om}$ and 
$\fo_{l^*+1}^+\fo_{\fp;l^*}$ 
of the composition
$$\fM_{k,l+1}^*(\wt{B};J)\xlra{f_{k,l+1;l^*+1}} \fM_{k,l}^*(\wt{B};J) 
\stackrel{\ev}{\lra} Y^k\!\times\!X^l$$
are the same;

\item\label{RinterchM_it} 
the interchange of two real points $x_i$ and $x_j$ with $2\!\le\!i,j\!\le\!k$ 
preserves~$\fo_{\fp;l^*}$;

\item\label{RinterchM2_it}  
if $2\!\le\!i\!\le\!k$ and  $\u\!\in\!\fM_{k,l}^*(\wt{B};J)$,
$\fo_{\fp;l^*}$ is preserved at~$\u$ by the interchange of 
the real points $x_1$ and $x_{j_i(\u)}$ with $2\!\le\!i\!\le\!k$
if and only if $(k\!-\!1)(i\!-\!1)n\!\in\!2\Z$;

\item\label{fM1cij_it} if $1\!\le\!i,j\!\le\!l$, $\u\!\in\!\fM_{k,l}^*(\wt{B};J)$,
and the marked points $z_i^+$ and~$z_j^+$ of~$\u$ lie on the same connected component
of $\P^1\!-\!S^1$, then $\fo_{\fp;l^*}$ is preserved at~$\u$  by the interchange of 
conjugate pairs $(z_i^+,z_i^-)$ and $(z_j^+,z_j^-)$;

\item\label{fM1c1i_it} if $l^*\!<\!i\!\le\!l$ (resp.~$1^*\!<\!i\!\le\!l^*$) and 
$\u\!\in\!\fM_{k,l}^*(\wt{B};J)$,
$\fo_{\fp;l^*}$ is preserved (resp.~reversed) at~$\u$ by
the interchange of the points in a conjugate pair $(z_i^+,z_i^-)$
if and only if $n\!\in\!2\Z$;

\item\label{orient1pm_it} if $\wt{B}\!=\!(B,b)$ and $\u\!\in\!\fM_{k,l}^*(\wt{B};J)$, 
$\fo_{\fp;l^*}$ is preserved at~$\u$
by the interchange of the points in the conjugate pair  $(z_1^+,z_1^-)$ if and only~if 
$$\hspace{-.3in}
\frac{\lr{c_1(X,\om),B}\!+\!1\!-\!k}{2}+l\!-\!l^*\!+\!n\binom{k\!-\!1}{2}
\!+\!
\lr{w_2(\fp),B}=\begin{cases}0,&\hbox{if}~\fp\!\in\!\cP_X^-(Y);\\
\lr{w_2(X),B},&\hbox{if}~\fp\!\in\!\cP_X^+(Y).
\end{cases}$$

\end{enumerate}
\end{lmm}

\begin{lmm}\label{orient2_lmm}
Suppose $(X,\om,\phi)$, $n$, $Y$, $k,l,l^*$, $\wt{B}$, and $J$ 
are in Lemma~\ref{orient_lmm}.
Every  relative $\OSpin$-structure~$\os$ on $Y\!\subset\!X$ 
determines a relative orientation~$\fo_{\os;l^*}$ of the map~\eref{fMevdfn_e}  
with the following properties:
\begin{enumerate}[label=($\fo_{\os}\arabic*$),leftmargin=*]

\setcounter{enumi}{-1}

\item\label{OSpinPrp_it} the statements \ref{etaMorient_it}, \ref{fforient_it}, \ref{RinterchM_it}, 
\ref{fM1cij_it}, and~\ref{fM1c1i_it} in Lemma~\ref{orient_lmm}
with $\fp$ replaced by~$\os$ hold
and \hbox{$\fo_{\ov\os;l^*}\!=\!-(-1)^k\fo_{\os;l^*}$}; 

\setcounter{enumi}{2}

\item\label{RinterchM2os_it}  
if $2\!\le\!i\!\le\!k$ and  $\u\!\in\!\fM_{k,l}^*(\wt{B};J)$,
$\fo_{\os;l^*}$ is preserved at~$\u$ by the interchange of 
the real points $x_1$ and $x_{j_i(\u)}$ with $2\!\le\!i\!\le\!k$
if and only if $(k\!-\!1)(i\!-\!1)(n\!+\!1)\!\in\!2\Z$;

\setcounter{enumi}{5}

\item\label{orient1os_it} if $\wt{B}\!=\!(B,b)$ and $\u\!\in\!\fM_{k,l}^*(\wt{B};J)$, 
$\fo_{\os;l^*}$ is preserved at~$\u$
by the interchange of the points in the conjugate pair  $(z_1^+,z_1^-)$ if and only~if 
$$\frac{\lr{c_1(X,\om)}}{2}\!+\!k\!-\!1\!+\!l\!-\!l^*\!+\!(n\!+\!1)\binom{k\!-\!1}{2}
\!+\!\lr{w_2(\os),B}=0.$$

\end{enumerate}
\end{lmm}

The next two lemmas concern properties of the orientations $\fo_{\fp;l^*}$ 
and~$\fo_{\os;l^*}$ in special cases.
The two statements of first lemma follow from 
the CROrient~\ref{CROPin2SpinRed_prop}, \ref{CROSpinPinSES_prop}\ref{CROsesSpin_it}, 
and~\ref{CRONormal_prop}\ref{CROnormSpin_it} properties in Section~\ref{OrientPrp_subs1}.

\begin{lmm}\label{orient1b_lmm}
Suppose $(X,\om,\phi)$, $n$, $Y$, $k,l,l^*$, $\wt{B}$, $J$,
and $\fp$ are as in Lemma~\ref{orient_lmm} and 
the pair $(k,\wt{B})$ satisfies~\eref{BKcond_e0}.
\begin{enumerate}[label=($\fo_{\fp}\arabic*$),leftmargin=*]

\setcounter{enumi}{6}

\item\label{deg0Pin_it} If $k,l\!=\!1$ and $\wt{B}\!=\!0$, 
$\fo_{\fp;1}$ is the orientation of~$\la_{\u}(\ev)$ induced by
the diffeomorphism~$\ev_1^{\R}$ and the orientation~$\fo_{\om}$;

\item\label{prodPin_it} If $\fo$ is an orientation on~$TX^{\phi}$,
$\fs$ is the image of~$\fp$ under the isomorphism~\eref{RelPin2SpinRed_e}, 
and $\os\!\equiv\!(\fo,\fs)$, then $\fo_{\fp;l^*}\!=\!\fo_{\os;l^*}$.

\end{enumerate}
\end{lmm}

\vspace{.1in}

Let $(X',\om',\phi')$, $(X'',\om'',\phi'')$, $Y',Y''$, $(X,\om,\phi)$, $Y$, $k,l$,
$\wt{B}''$, and $\wt{B}$ be as in Theorem~\ref{RGWs_thm2}\ref{RGWprod_it0},
$l^*\!\in\![l]$, and $J\!=\!J'\!\times\!J''$ for some generic 
$J'\!\in\!\cJ_{\om'}^{\phi'}$ and $J''\!\in\!\cJ_{\om''}^{\phi''}$.
There is then a natural isomorphism
\BE{deg0spl_e}
\fM_{k,l}^*(\wt{B};J)\approx Y'\!\times\!\fM_{k,l}^*(\wt{B}'';J'').\EE
For $[\u]\!\in\!\fM_{k,l}^*(\wt{B};J)$, let 
$[\u']\!\in\!\fM_{k,l}^*(0;J')$ and $[\u'']\!\in\!\fM_{k,l}^*(\wt{B}'';J'')$
be the projections of~$\u$ to~$X'$ and~$X''$.
Relative $\OSpin$-structure~$\os'$ on $Y'\!\subset\!X'$ determines homotopy 
classes~$\fo'\la_{\u'}^{\R}(\os')$ orientations of 
$\la(T_{u'(\P^1)}Y')$, $\la_{\u'}^{\R}(X')$, and~$\la_{\u'}^{\C}(X')$
and thus a homotopy class~$\fo'\la_{\u'}^{\R}(\os')$ of isomorphisms
$$\la\big(T_{u'(\P^1)}Y'\big)\approx \la_{\u'}^{\R}(X')\!\otimes\!\la_{\u'}^{\C}(X').$$
If the pair $(k,\wt{B})$ satisfies~\eref{BKcond_e0},
a relative $\Pin^{\pm}$-structure~$\fp''$ on $Y''\!\subset\!X''$ determines homotopy 
class~$\fo_{\fp'';l^*}$ of isomorphisms
$$\la_{\u''}\big(T\fM_{k,l}^*(\wt{B}'';J'')\!\big)\approx
\la_{\u''}^{\R}(X'')\!\otimes\!\la_{\u''}^{\C}(X''),$$
respectively.
So does a relative $\OSpin$-structure~$\fp''$, whether or not 
the pair $(k,\wt{B})$ satisfies~\eref{BKcond_e0}.

Along with the isomorphisms~\eref{RBPsesdfn_e4} and~\eref{deg0spl_e}, 
these homotopy classes determine an orientation 
$(\fo'\la_{\u'}^{\R}(\os')\!)\fo_{\fp'';l^*}$ of the line~$\la_{\u}(\ev)$ in~\eref{OrientSubs_e3}
as below~\eref{CROses_e0}.
Let \hbox{$\pi',\pi''\!:X\!\lra\!X',X''$} be the two projections.
The next statement follows from the CROrient~\ref{CROSpinPinSES_prop}\ref{CROsesPin_it}
property in Section~\ref{OrientPrp_subs1}.

\begin{lmm}\label{orient2b_lmm}
With the assumptions as above and $\fp\!=\!\llrr{\pi'^*\os',\pi''^*\fp''}_{\oplus}$,
$$\fo_{\fp;l^*}\big|_{\u}=\big(\fo'\la_{\u'}^{\R}(\os')\!\big)\fo_{\fp'';l^*}\!\big|_{\u''}
\qquad\hbox{iff}\quad \frac{(\dim\,X')(\dim\,X'')}4\binom{k}2\in 2\Z.$$
\end{lmm}

\subsection{Definition of curve signs}
\label{MapSignDfn_subs}

Let $(X,\om,\phi)$, $n$, $Y$, $k,l,l^*$, $\wt{B}$, and $J$
be as in Lemma~\ref{orient_lmm}.
Suppose that either $\fp$ is a relative $\Pin^{\pm}$-structure on $Y\!\subset\!X$
and the pair $(k,B)$ satisfies the condition~\eref{BKcond_e0} or
$\fp$ is a relative $\OSpin$-structure on~$Y$.
For $l'\!\in\![l]$ and a tuple 
\BE{bhdfn_e}\bh \equiv(h_i\!:H_i\!\lra\!X)_{i\in[l']}\EE
of maps, define
$$\cZ_{k,l;\bh}^*(\wt{B};J)=\big\{\big(\u,(y_i)_{i\in[l']}\big)\!\in\!
\fM_{k,l}^*(\wt{B};J)\!\times\!\prod_{i=1}^{l'}\!H_i\!:
\ev_i^+(\u)\!=\!h_i(y_i)\,\forall\,i\!\in\![l']\big\}.$$
We denote by
\BE{JakePseudo_e} \ev_{k,l;\bh}\!: \cZ_{k,l;\bh}^*\big(\wt{B};J\big)\lra Y^k\!\times\!X^{l-l'}\EE
the map induced by~\eref{fMevdfn_e}.
If $J\!\in\!\cJ_{\om}^{\phi}$ is generic and the maps $h_i$ are smooth and generically chosen,
\eref{JakePseudo_e} is a smooth map from a smooth manifold of dimension
$$\dim\cZ_{k,l;\bh}^*\big(\wt{B};J\big)=\ell_{\om}(\wt{B})\!+\!k\!+\!2l
-\sum_{i=1}^{l'}\codim\, h_i\,,$$
where $\codim\,h_i\!\equiv\!\dim X\!-\!\dim H_i$.
Orientations on~$H_i$ determine an orientation~$\fo_{\bh}$ on
\BE{Mhdfn_e} M_{\bh}\equiv \prod_{i=1}^{l'}\!H_i\,.\EE
Along with~$\fo_{\bh}$ and the symplectic orientation~$\fo_{\om}$ of~$X$, 
the relative orientation~$\fo_{\fp;l^*}$ of Lemma~\ref{orient_lmm} or~\ref{orient2_lmm}
determines a relative orientation~$\fo_{\fp;l^*;\bh}$ of the map~\eref{JakePseudo_e}.
The precise definition of~$\fo_{\fp;l^*;\bh}$ generally depends on the orientation
conventions, such those in Section~5.1 of~\cite{JSG}, but there is no ambiguity 
if the dimensions of all~$H_i$ are even.\\

We now restrict to the case $l'\!=\!l$ so that~\eref{dimkbhcond_e} holds.
If $J\!\in\cJ_{\om}^{\phi}$ is generic and the maps $h_i$ are smooth and generically chosen,
then~\eref{JakePseudo_e} is a smooth map between manifolds of the same dimension.
If \hbox{$\wt\u\!\in\!\cZ_{k,l;\bh}^*(\wt{B};J)$} is a regular point of~\eref{JakePseudo_e},
we set~$\fs_{\fp;l^*}(\wt\u;\bh)$ to be $+1$ if the isomorphism
$$\nd_{\wt\u}\ev_{k,l;\bh}\!: 
T_{\wt\u}\cZ_{k,l;\bh}^*\big(\wt{B};J\big)\lra T_{\ev_{k,l;\bh}(\wt\u)}Y^k$$
lies in the homotopy class determined by~$\fo_{\fp;l^*;\bh}$ and  
 to be $-1$ otherwise.\\

Suppose $k\!+\!2l\!\ge\!5$, $i\!\in\![l^*]\!-\![1]$, the codimension of~$h_i$ is~2, 
\begin{gather*}
\bh'=\big(h_1,\ldots,h_{i-1},h_{i+1},\ldots,h_l\big), \quad\hbox{and}\\
\wt\u'\!\equiv\!\big(\ff_{k,l;i}(\u),(y_1,\ldots,y_{i-1},y_{i+1},\ldots,y_l)\!\big)
\in\cZ_{k,l-1;\bh'}^*\big(\wt{B};J\big).
\end{gather*}
By the genericity assumptions above, the homomorphism
\BE{JakePseudo_e10b}
T_{z_i^+}\P^1\!\oplus\!T_{y_i}H_i\lra T_{u(z_i^+)}X\!=\!T_{h_i(y_i)}X, \quad
(v,w)\lra \nd_{z_i^+}u(v)\!+\!\nd_{y_i}h_i(w),\EE
is then an isomorphism.
We set $\fs_i(\wt\u)$ to be $+1$ if this isomorphism 
is orientation-preserving and to be $-1$ if it is orientation-reversing.
By~\ref{fforient_it}, \ref{fM1cij_it}, and~\ref{fM1c1i_it} in Lemma~\ref{orient_lmm}
or the corresponding statements in Lemma~\ref{orient2_lmm}, 
\BE{extrasign_e}\fs_{\fp;l^*;\bh}(\wt\u)=\fs_i(\wt\u)\fs_{\fp;l^*-1;\bh'}(\wt\u')\,;\EE
see the proof of the real divisor relation of \cite[Proposition~5.2]{RealWDVV}.\\

Suppose that $k\!+\!2l\!\ge\!5$, $l'\!\in\![l^*]$, and the codimension of~$h_i$ is~2
for every $i\!\in\![l']$.
Let
$$\bh'=\big(h_{l'+1},\ldots,h_l\big)$$
and $\u'\!\in\!\fM_{k,l-l'+1}^*(\wt{B},J)$ be 
the image of $\wt\u\!\in\!\cZ_{k,l;\bh}^*(\wt{B};J)$ under the composition
$$\cZ_{k,l;\bh}^*(\wt{B};J)\subset 
\fM_{k,l}^*(\wt{B},J)\!\times\!M_{\bh}\lra \fM_{k,l}^*(\wt{B},J)
\xlra{\ff_{k,l-l'+2;2}\circ\ldots\circ\ff_{k,l;l'}}  \fM_{k,l-l'+1}^*(\wt{B},J)$$
and $\u\!=\!\ff_{k,1+l;1}(\u')$.
By~\eref{extrasign_e}, the sign 
\BE{fsfpdfn_e}\fs_{\fp;l^*-l';\bh'}(\u)\equiv 
\fs_1(\wt\u)\!\ldots\!\fs_{l'}(\wt\u)\fs_{\fp;l^*;\bh}(\wt\u)
\in\big\{\pm1\big\}\nota{sphu2@$\fs_{\fp;l;\bh}(\u)$}\EE
is independent of generic choices of $h_2,\ldots,h_{l'}$ and 
$\wt\u\!\in\!\cZ_{k,l;\bh}^*(\wt{B};J)$ with the pair $(\u',y_1)$ fixed.
By Lemma~\ref{orient_lmm}\ref{fM1cij_it} or the corresponding statement in Lemma~\ref{orient2_lmm},
this sign is also independent of the choices of~$h_1$ and $(\u',y_1)$ with~$\u$ fixed
if $l'\!\ge\!2$ and the marked points~$z_1^+$ and~$z_2^+$ of~$\wt\u$ lie
in the same component of~$\P^1\!-\!S^1$.\\

If $\phi_1$ is an  involution on~$H_1$ so~that 
\BE{phi1dfncond_e}\phi\!\circ\!h_1\!=\!h_1\!\circ\!\phi_1\!:H_1\lra X,\EE
then the tuple $\wh\u$ obtained from $\wt\u$ by exchanging the marked points $z_1^+$ and~$z_1^-$
and replacing $y_1$ by $\phi_1(y_1)$ lies in $\cZ_{k,l;\bh}^*(\wt{B};J)$.
If in addition the isomorphism
\BE{phi1NCdfn_e}\nd_{y_1}\phi\!: \frac{T_{h_1(y_1)}X}{\nd_{y_1}h_1(T_{y_1}H_1)}\lra
\frac{T_{\phi(h_1(y_1)\!)}X}{\nd_{\phi_1(y_1)}h_1(T_{\phi_1(y_1)}H_1)}\EE
on the normal bundle of~$h_1$ induced by $\phi$ and~$\phi_1$ is orientation-reversing,
then $\fs_1(\wh\u)\!=\!\fs_1(\wt\u)$.
Along with the conclusions of the previous paragraph, this implies that 
the sign~\eref{fsfpdfn_e} is independent of the choices of~$h_1$ and $(\u',y_1)$ with~$\u$ 
fixed if the relative orientation~$\fo_{\fp;l^*}$ is reversed at~$\wt\u$
by the interchange of the points in the conjugate pair~$(z_1^+,z_1^-)$
(so that $\fs_{\fp;l^*;\bh}(\wt\u)\!=\!\fs_{\fp;l^*;\bh}(\wh\u)$).
By Lemmas~\ref{orient_lmm}\ref{orient1pm_it} and~\ref{orient2_lmm}\ref{orient1os_it}, 
the latter is the case~if and only~if
\BE{sgndfncond_e3}l\!-\!l^*\!+\!n\binom{k\!-\!1}{2}\!+\!\lr{w_2(\fp),B}
=\begin{cases}\frac{\lr{c_1(X,\om),B}-1-k}{2}\!+\!2\Z,
&\hbox{if}~\fp\!\in\!\cP_X^-(Y);\\
\frac{\lr{c_1(X,\om),B}-1-k}{2}\!+\!\lr{w_2(X),B},
&\hbox{if}~\fp\!\in\!\cP_X^+(Y);\\
\frac{\lr{c_1(X,\om),B}}{2}\!+\!\binom{k-1}{2}\!+\!k\!+\!2\Z,
&\hbox{if}~\fp\!\in\!\OSp_X(Y).
\end{cases}\EE
The sign $\fs_{\fp;l^*;\bh'}(u(\P^1)\!)\!\equiv\!\fs_{\fp;l^*;\bh'}(\u)$ of $u(\P^1)$
is then independent of generic choices of the codimension~2 maps~$h_i$  
with $i\!\in\![l']$ and of $l'$ (as long as $k\!+\!2l\!\ge\!3$).\\

Under the assumptions in~\ref{Rcond_it1} and~\ref{Rcond_it2} in Section~\ref{g0maps_subs}, 
\eref{sgndfncond_e3} is equivalent~to
$$l\!-\!l^* \cong \frac{\ell_{\om}(\wt{B})\!-\!k}{2} \qquad\hbox{mod~2}.$$
Under the assumptions in~\ref{Rcond_it3}-\ref{Rcond_it5} in Section~\ref{g0maps_subs}, 
\eref{sgndfncond_e3} is equivalent~to
$$l\!-\!l^* \cong  \frac{\ell_{\om}(\wt{B})}{2}\!+\!k \qquad\hbox{mod~2}.$$
By~\eref{dimkbhcond_e} and~\eref{Rcond_e0}, both congruences hold.
This in particular defines the signs~$\fs_{\fp,l^*;\bh}(C)$ in~\eref{mapcountdfn_e}.\\

The relations~\ref{div_it} and~\ref{OrientRev_it} in Theorem~\ref{RGWs_thm}
and the second statement in~\ref{RGWprod_it0} in Theorem~\ref{RGWs_thm2}
follow from~\eref{extrasign_e}, the last statement of 
Lemma~\ref{orient2_lmm}\ref{OSpinPrp_it}, and Lemma~\ref{orient2b_lmm}, respectively.
For the purposes of establishing the first cases of 
the properties~\ref{deg0_it} and~\ref{ins1_it} in Theorem~\ref{RGWs_thm}, 
we take $k\!=\!1$, $l\!=\!l^*\!=\!1$, $\wt{B}\!=\!0$, and $h_1\!=\!\id_X$.
The claims in these two cases then follow from Lemma~\ref{orient1b_lmm}.
The motivation behind~\ref{sphere_it} is that a $J$-holomorphic curve passing 
though a point $p_{k+1}^{\R}\!\in\!Y$
intersects an infinitesimal sphere $S(\cN_{p_{k+1}^{\R}}Y)$ at two points.
By a straightforward computation, the two resulting real rational irreducible $J$-holomorphic 
curves $C\!\subset\!X$ through~$S(\cN_{p_{k+1}^{\R}}Y)$ contribute to the left-hand side
in~\ref{sphere_it} with the same sign as the nearby curve through~$p_{k+1}^{\R}$
contributes to the right-hand side;
see the proof of Proposition~2.1 in~\cite{RealWDVV3}.\\

Suppose $i\!\in\![l^*]\!-\![1]$ (resp.~$i\!\in\![l]\!-\![l^*]$) and
$h_i$ is an automorphism of~$H_i$  with the sign of~$(-1)^n$ (resp.~$(-1)^{n+1}$)
which satisfies~\eref{phi1dfncond_e} with~1 replaced by~$i$.
For each $\wt\u\!\in\!\cZ_{k,l;\bh}^*(\wt{B};J)$, 
the tuple~$\wh\u$ obtained from~$\wt\u$ by exchanging the marked points~$z_i^+$ and~$z_i^-$
and replacing~$y_i$ by~$\phi_i(y_i)$ then lies in~$\cZ_{k,l;\bh}^*(\wt{B};J)$.
The isomorphism in~\eref{phi1NCdfn_e} with~1 replaced by~$i$ 
is orientation-preserving (resp.~orientation-reversing) in this case.
Along with Lemma~\ref{orient_lmm}\ref{fM1c1i_it} or 
the corresponding statement of Lemma~\ref{orient2_lmm},
this implies~that 
$$\fs_{\fp,l^*;\bh}\big(\wh{u}(\P^1)\!\big)=-\fs_{\fp,l^*;\bh}\big(\wt{u}(\P^1)\!\big)\,.$$
Thus, the signed cardinality of~$\cZ_{k,l;\bh}^*(\wt{B};J)$ vanishes.
Combining this conclusion with~\eref{Rcond_e0},
we obtain Theorem~\ref{RGWs_thm}\ref{dimvan_it2}.

\subsection{Proof of invariance}
\label{InvarPf_subs}

The proof of the invariance of the sums~\eref{mapcountdfn_e} under suitable topological
conditions goes back to~\cite{Sol}.
The argument of~\cite{Sol} uses the moduli space of stable nodal disk maps to~$(X,Y)$,
but can be reformulated in terms of the moduli space $\ov\fM_{k,l}(\wt{B};J)$
of stable real genus~0 $J$-holomorphic degree~$\wt{B}$ maps to~$(X,\phi)$ as done
in~\cite{Ge2}. 
After summarizing the invariance argument in~\cite{Sol,Ge2}, 
which follows the general principle in~\cite{Wel4} of showing directly
that the sums~\eref{mapcountdfn_e} are invariant along 
generic paths $\wt\p\!\equiv\!(\p_t)_{t\in[0,1]}$ of points in~$Y$,
$\wt\bh\!\equiv\!(\bh_t)_{t\in[0,1]}$ of constraints in~$X$ or~$X\!-\!Y$, as appropriate, 
and $\wt{J}\!\equiv\!(J_t)_{t\in[0,1]}$ 
of almost complex structures in~$\cJ_{\om}^{\phi}$,
we recast it in terms of pseudocycles as in~\cite{RealWDVV,RealWDVV3}.
We continue with the notation and setup as at the beginning of Section~\ref{g0maps_subs}
and denote by~$\fp$ either a relative $\Pin^-$-structure as 
in~\ref{Rcond_it1} or~\ref{Rcond_it2} in Section~\ref{g0maps_subs}
or an $\OSpin$-structure as in~\ref{Rcond_it3}-\ref{Rcond_it5}.\\

The total evaluation map~\eref{fMevdfn_e} extends to a continuous map
\BE{wtevdfn_e}\wt\ev\!: \ov\fM_{k,l}\big(\wt{B};\wt{J}\big)\equiv 
\bigcup_{t\in[0,1]}\!\!\!\{t\}\!\times\!\ov\fM_{k,l}\big(\wt{B};J_t\big)
\lra \big([0,1]\!\times\!Y\big)^k\!\times\!\big([0,1]\!\times\!X\big)^l;\EE
the domain of this extension is compact.
For a generic path~$\wt{J}$ in~$\cJ_{\om}^{\phi}$, the subspace 
$$\ov\fM_{k,l}^*(\wt{B};\wt{J})\subset \ov\fM_{k,l}(\wt{B};\wt{J})$$
consisting of simple maps is stratified by smooth manifolds.
The top-dimensional stratum of $\ov\fM_{k,l}^*(\wt{B};\wt{J})$ is 
the subspace 
$$\fM_{k,l}^*\big(\wt{B};\wt{J}\big)\equiv 
\bigcup_{t\in[0,1]}\!\!\!\{t\}\!\times\!\fM_{k,l}^*\big(\wt{B};J_t\big)$$
consisting of maps from~$(\P^1,\tau)$.
If~$k\!\neq\!0$ or $b\!\neq\!0$ (as is the case under the assumptions in 
\ref{Rcond_it1}-\ref{Rcond_it3} and~\ref{Rcond_it5} in Section~\ref{g0maps_subs}),
the codimension~1 strata~$\cS$ consist of maps from a pair of disks 
with a boundary point in common.
The relative orientations~$\fo_{\fp;l^*}$ of Lemmas~\ref{orient_lmm} 
and~\ref{orient2_lmm} extend across some codimension~1 strata
(those with $\ep_{l^*}(\cS)\!\cong\!0,1$ mod~4 below), but not others.
We denote the union of $\fM_{k,l}^*(\wt{B};\wt{J})$ with the codimension~1 strata~$\cS$
over which~$\fo_{\fp;l^*}$ extends by $\fM_{k,l}^{\bigst}(\wt{B};\wt{J})$.\\

For generic choices of~$\wt\p$, $\wt\bh$, and~$\wt{J}$ as above,
the preimage~$\ov\fM_{\wt\p;\wt\bh}(\wt{B};\wt{J})$ of the paths~$\wt\p$ and~$\wt\bh$
under the extended total evaluation map~\eref{wtevdfn_e} 
is a one-dimensional manifold~with
\BE{InvarPf_e5}\prt\ov\fM_{\wt\p;\wt\bh}\big(\wt{B};\wt{J}\big)
=\fM_{\p_1;\bh_1}(\wt{B};J_1)\sqcup \fM_{\p_0;\bh_0}(\wt{B};J_0).\EE
The intersection of~$\fM_{\wt\p;\wt\bh}(\wt{B};\wt{J})$ with a codimension~$\fc$ stratum  
of $\ov\fM_{k,l}^*(\wt{B};\wt{J})$ is a submanifold of~$\fM_{\wt\p;\wt\bh}(\wt{B};\wt{J})$ 
 of codimension~$\fc$ and is in particular
empty for $\fc\!\ge\!2$.
The intersection of~$\fM_{\wt\p;\wt\bh}(\wt{B};\wt{J})$ with $\fM_{k,l}^*(\wt{B};\wt{J})$
inherits an orientation~$\fo_{\fp;l^*;\wt\p;\wt\bh}$ 
from the relative orientation~$\fo_{\fp;l^*}$ of Lemma~\ref{orient_lmm}
or~\ref{orient2_lmm}, as appropriate.
If~$\fo_{\fp;l^*}$ extends over
a codimension~1 stratum~$\cS$ of~$\ov\fM_{k,l}^*(\wt{B};\wt{J})$, 
then~$\fo_{\fp;l^*;\wt\p;\wt\bh}$ extends over 
the intersection of~$\fM_{\wt\p;\wt\bh}(\wt{B};\wt{J})$ with~$\cS$. 
A~key observation originating in~\cite{Sol} is that~$\ov\fM_{\wt\p,\wt\bh}^*(\wt{B};\wt{J})$
does not intersect the codimension~1 strata~$\cS$ over which $\fo_{\fp;l^*}$ 
does not extend.\\

By the positivity assumptions in \ref{Rcond_it1}-\ref{Rcond_it3} and~\ref{Rcond_it5},
the image of the complement of the simple maps under the extended evaluation~map,  
$$\wt\ev\big(\ov\fM_{k,l}(\wt{B};\wt{J})\!-\!\ov\fM_{k,l}^*(\wt{B};\wt{J})\!\big)
\subset \big([0,1]\!\times\!Y\big)^k\!\times\!\big([0,1]\!\times\!X\big)^l,$$
is covered by a smooth map from a manifold of dimension 2 less than
that of~$\fM_{k,l}^*(\wt{B};\wt{J})$; see the proof of Proposition~5.2 in~\cite{RealWDVV}.
Thus, 
$$\ov\fM_{\wt\p;\wt\bh}\big(\wt{B};\wt{J}\big)\subset 
\fM_{k,l}^{\bigst}\big(\wt{B};\wt{J}\big)$$
is a compact oriented one-dimensional manifold so that~\eref{InvarPf_e5} 
respects the orientation~$\fo_{\fp;l^*;\wt\p;\wt\bh}$ on $\ov\fM_{\wt\p;\wt\bh}(\wt{B};\wt{J})$,
the orientation~$\fo_{\fp;l^*;\p_1;\bh_1}$ on
$\fM_{\p_1;\bh_1}(\wt{B};J_1)$, and the opposite 
of the orientation~$\fo_{\fp;l^*;\p_0;\bh_0}$ on
$\fM_{\p_0;\bh_0}(\wt{B};J_0)$.
This implies~that
$$\sum_{C\in\fM_{\p_0;\bh_0}(\wt{B};J_0)}\!\!\!\!\!\!\!\!\!\!\!\!\!\!\fs_{\fp;l^*;\bh_0}(C)
=\sum_{C\in\fM_{\p_1;\bh_1}(\wt{B};J_1)}\!\!\!\!\!\!\!\!\!\!\!\!\!\!\fs_{\fp;l^*;\bh_1}(C)$$
and shows that the functionals~\eref{RGWdfn_e} are well-defined at least on the insertions
from~$H^6(X;\Z)$.
The only case in the setting~\ref{Rcond_it4} not covered by~\ref{Rcond_it3} is $k\!=\!0$ 
and $b\!=\!0$.
The sum~\eref{mapcountdfn_e} then vanishes; 
see the proof of Proposition~1.3(1) in~\cite{RealWDVV3}.\\

A codimension~1 two-disk boundary stratum~$\cS$ of $\ov\fM_{k,l}^*(\wt{B};J)$ or 
$\ov\fM_{k,l}^*(\wt{B};\wt{J})$ is characterized~by the distributions~of
the degree~$\wt{B}$ of the map components~$u$ of its elements~$\u$,
the $k$ real marked points, and the $l$ conjugate pairs of marked points
between the two irreducible components of~the domains.
Supposing that $l\!\in\!\Z^+$, we define 
$$\wt{B}_1(\cS)\!\equiv\!(B_1(\cS),b_1(\cS))\in \wt{H}_2^{\phi}(X,Y), \qquad
K_1(\cS)\subset[k], \qquad\hbox{and}\qquad L_1(\cS)\subset[l]$$
to be
the degree of the restriction of the map~$u$ of the elements~$\u$ of~$\cS$ to
the irreducible component~$\P^1_1$ of the domain carrying the marked points~$z_1^{\pm}$,
the set of real marked points carried by~$\P^1_1$,
and the set of conjugate pairs of marked points carried by~$\P^1_1$, respectively.
We define 
$$\wt{B}_2(\cS)\!\equiv\!(B_2(\cS),b_2(\cS))\in \wt{H}_2^{\phi}(X,Y), \qquad
K_2(\cS)\subset[k], \qquad\hbox{and}\qquad L_2(\cS)\subset[l]$$
to be the analogous objects for the irreducible component~$\P^1_2$ of the domain 
not carrying the marked points~$z_1^{\pm}$.
Let 
$$r(\cS)=\begin{cases}1,&\hbox{if}~1\!\not\in\!K_2(\cS);\\
2,&\hbox{if}~1\!\in\!K_2(\cS);\end{cases} \quad
\ep_{l^*}(\cS)=\blr{c_1(X,\om),B_2(\cS)}
\!-\!\big((n\!-\!1)\big|K_2(\cS)|\!+\!2|L_2(\cS)\!-\![l^*]\big|\big).$$ 

\vspace{.2in}

If \ref{Rcond_it1} or~\ref{Rcond_it2} holds, Corollary~\ref{CROdegenP1H3_crl}
implies that the orientation~$\fo_{\fp;\u}^D$ on the line~$\wt\la_{\u}(D_J^{\phi},X)$ 
extends across~$\cS$ if and only~if
\begin{equation*}\begin{split}
&\frac{\lr{c_1(X,\om),B_2(S)}(\lr{c_1(X,\om),B_2(S)}\!-\!1)}2
\!+\!n\frac{|K_2(\cS)|(|K_2(\cS)|\!-\!1)}2\\
&\hspace{1in}
+\!\big(r(\cS)\!-\!1\big)\big(|K_1(\cS)|\!-\!1\big)\big(|K_2(\cS)|\!-\!1\big)
\!+\!\big|K_2(\cS)\big|\lr{c_1(X,\om),B_2(S)}\in2\Z.
\end{split}\end{equation*}
If \ref{Rcond_it3} or \ref{Rcond_it5} holds, 
the same is implied by the CROrient~\ref{CROos_prop}\ref{osflip_it} property
in Section~\ref{OrientPrp_subs1} and the definition of the orientation 
$\la_{\cC}^{\R}(\os)\!\equiv\!\la_{\cC}^{\R}(\fo)$ on~$\la_{\u}^{\R}(X)$ above 
the CROrient~\ref{CROPin2SpinRed_prop} property.
If $|K_2(\cS)|\!+\!2|L_2(\cS)|\!\ge\!2$, the orientation~$\fo_{k,l;l^*}$ on~$\cM_{k,l}^{\tau}$
defined in Section~\ref{DM_subs} 
extends across the image of~$\cS$ in~$\cM_{k,l}^{\tau}$ under the forgetful morphism
dropping the map component if and only~if
$$\frac{|K_2(\cS)|(|K_2(\cS)|\!-\!1)}2
\!+\!\big(r(\cS)\!-\!1\big)\big(|K_1(\cS)|\!-\!1\big)\big(|K_2(\cS)|\!-\!1\big)
\!+\!\big|K_2(\cS)\big|\!+\!\big|L_2(\cS)\!-\![l^*]\big|\in2\Z.$$
If any of the conditions \ref{Rcond_it1}-\ref{Rcond_it3} or~\ref{Rcond_it5} holds 
and  $|K_2(\cS)|\!+\!2|L_2(\cS)|\!\ge\!2$, 
the relative orientation~$\fo_{\fp;l^*}$ of~\eref{fMevdfn_e} defined
via~\eref{OrientSubs_e3} thus extends across the codimension~1 stratum~$\cS$ of
$\ov\fM_{k,l}^*(\wt{B};J)$ or $\ov\fM_{k,l}^*(\wt{B};\wt{J})$ if and only~if 
\BE{orientcond_e}\ep_{l^*}(\cS)\equiv0,1 \qquad\hbox{mod}~4.\EE
In light of Lemma~\ref{orient_lmm}\ref{fforient_it} or 
the corresponding statement of Lemma~\ref{orient2_lmm}, as appropriate, 
the condition $|K_2(\cS)|\!+\!2|L_2(\cS)|\!\ge\!2$ is not necessary
for the last conclusion.\\

For $J\!\in\!\cJ_{\om}^{\phi}$ generic (resp.~generic path~$\wt{J}$ in $\cJ_{\om}^{\phi}$),
we denote by $\fM_{k,l}^{\bigst}(\wt{B};J)$ (resp.~$\fM_{k,l}^{\bigst}(\wt{B};\wt{J})$)
the union of $\fM_{k,l}^*(\wt{B};J)$ (resp.~$\fM_{k,l}^*(\wt{B};\wt{J})$)
with the codimension~1 strata~$\cS$ of $\ov\fM_{k,l}^*(\wt{B};J)$ 
(resp.~$\ov\fM_{k,l}^*(\wt{B};\wt{J})$) satisfying~\eref{orientcond_e}.
For tuples
$$\bh\equiv\big(h_i\!:H_i\!\lra\!X\big)_{i\in[l]} \qquad\hbox{and}\qquad
\wt\bh\equiv\big(\wt{h}_i\!:\wt{H}_i\!\lra\![0,1]\!\times\!X\big)_{i\in[l]}$$
of maps, define
\begin{equation*}\begin{split}
\cZ_{k,l;\bh}^{\bigst}(\wt{B};J)\equiv\big\{\big(\u,(y_i)_{i\in[l]}\big)\!\in\!
\fM_{k,l}^{\bigst}(\wt{B};J)\!\times\!\prod_{i=1}^l\!H_i\!:
\ev_i^+(\u)\!=\!h_i(y_i)\,\forall\,i\!\in\!l\big\},\\
\cZ_{k,l;\wt\bh}^{\bigst}(\wt{B};\wt{J})\equiv\big\{\big(\u,(y_i)_{i\in[l]}\big)\!\in\!
\fM_{k,l}^{\bigst}(\wt{B};\wt{J})\!\times\!\prod_{i=1}^l\!\wt{H}_i\!:
\wt\ev_i^+(\u)\!=\!\wt{h}_i(y_i)\,\forall\,i\!\in\!l\big\},
\end{split}\end{equation*}
where $\wt\ev_i^+$ is the composition of~\eref{wtevdfn_e} with the projection to the $i$-th
$[0,1]\!\times\!X$ component.
Let 
\BE{stevkdfn_e} \ev_{k,l;\bh}\!:\cZ_{k,l;\bh}^{\bigst}(\wt{B};J)\lra Y^k
\quad\hbox{and}\quad
\wt\ev_{k,l;\wt\bh}\!:\cZ_{k,l;\wt\bh}^{\bigst}(\wt{B};\wt{J})\lra Y^k\EE
be the maps induced by the evaluations at the points $x_1,\ldots,x_k$ of the domains.
We~define 
$$ \ev_{k,l;\bh}\!:\ov\cZ_{k,l;\bh}(\wt{B};J)\lra Y^k
\quad\hbox{and}\quad
\wt\ev_{k,l;\wt\bh}\!:\ov\cZ_{k,l;\wt\bh}(\wt{B};\wt{J})\lra Y^k$$
as above with $\fM_{k,l}^{\bigst}(\wt{B};J)$ replaced by $\ov\fM_{k,l}(\wt{B};J)$
and $\fM_{k,l}^{\bigst}(\wt{B};\wt{J})$ replaced by $\ov\fM_{k,l}(\wt{B};\wt{J})$.\\

Suppose in addition that $\bh$ is a tuple of smooth maps from oriented even-dimensional
manifolds in general position satisfying~\eref{dimkbhcond_e} and
$\wt\bh$ is a tuple of smooth maps from oriented odd-dimensional
bordered manifolds in general position so~that 
$$\prt\wt\bh=\{1\}\!\times\!\bh_1-\{0\}\!\times\!\bh_0$$
for some tuples~$\bh_0,\bh_1$ also satisfying the above conditions for~$\bh$.
The map~$\ev_{k,l;\bh}$ in~\eref{stevkdfn_e} 
is then a smooth map between manifolds of the same dimension;
it inherits a relative orientation~$\fo_{\fp;l^*;\bh}$ 
from the relative orientation~$\fo_{\fp;l^*}$ of Lemma~\ref{orient_lmm} or~\ref{orient2_lmm}, 
as appropriate.
The domain of the map~$\wt\ev_{k,l;\wt\bh}$ in~\eref{stevkdfn_e} is a bordered manifold;
this map inherits a relative orientation~$\fo_{\fp;l^*;\wt\bh}$ from
the relative orientation~$\fo_{\fp;l^*}$ of Lemma~\ref{orient_lmm} or~\ref{orient2_lmm}
so~that 
\BE{prtwtev_e} \prt\big(\wt\ev_{k,l;\wt\bh},\fo_{\fp;l^*;\wt\bh}\big)
=\big(\ev_{k,l;\bh_1},\fo_{\fp;l^*;\bh_1}\big)-\big(\ev_{k,l;\bh_0},\fo_{\fp;l^*;\bh_0}\big)\,.\EE

\vspace{.2in}

From now on, we also assume that $\bh,\bh_0,\bh_1$ are generic tuples of pseudocycles
satisfying~\eref{Rcond_e0}, 
$\wt\bh$ is a generic tuple of pseudocycle equivalences between 
the components of $\{0\}\!\times\!\bh_0$ and $\{1\}\!\times\!\bh_1$ 
as defined in~\cite{pseudo}
with the dimension~$n$ components of~$\wt\bh$ being pseudocycle equivalences 
to~$X\!-\!Y$, and
one of the conditions~\ref{Rcond_it1}-\ref{Rcond_it3} or~\ref{Rcond_it5}
in Section~\ref{g0maps_subs} holds.
We show below that the \sf{limit set} of~the map~$\ev_{k,l;\bh}$ in~\eref{stevkdfn_e},
$$\Om(\ev_{k,l;\bh})\equiv 
\bigcup_{\begin{subarray}{c}K\subset\cZ_{k,l;\bh}^{\bigst}(\wt{B};J)\\
\tn{cmpt}~~~~~~~ \end{subarray}}\hspace{-.3in}
\ev_{k,l;\bh}\big(\cZ_{k,l;\bh}^{\bigst}(\wt{B};J)\!-\!K\big),$$ 
is then covered by a smooth map from a manifold of dimension $nk\!-\!2$,
i.e.~2 less than that of
$\cZ_{k,l;\bh}^{\bigst}(\wt{B};J)$.
This map is thus a codimension~0 pseudocycle; its degree,
$$\deg\ev_{k,l;\bh}
\equiv \sum_{C\in\fM_{\p;\bh}(\wt{B};J)}\!\!\!\!\!\!\!\!\!\!\!\fs_{\fp;l^*;\bh}(C)
\equiv N_{\wt{B};k,l^*;\bh}^{\phi;\fp}(Y),$$
is independent of a generic choice of $\p\!\in\!Y^k$.
We also show that the limit set of the~map~$\wt\ev_{k,l;\wt\bh}$ in~\eref{stevkdfn_e},
$\Om(\wt\ev_{k,l;\wt\bh})$, is
covered by a smooth map from a manifold of dimension $nk\!-\!1$,
i.e.~2 less than that of $\cZ_{k,l;\wt\bh}^{\bigst}(\wt{B};\wt{J})$.
This map is thus a pseudocycle equivalence between $\ev_{k,l;\bh_0}$ and~$\ev_{k,l;\bh_1}$,
and so the above degree is also independent of the choice of pseudocycle representatives
for the even-dimensional elements of~$\wh{H}_*(X,Y)$. 
This shows that the functionals~\eref{RGWdfn_e} are well-defined.\\

By the proof of the second case in Theorem~\ref{RGWs_thm}\ref{ins1_it}, 
we can assume that \hbox{$\dim h_i\!\le\!2n\!-\!2$} for every $i\!\in\![l]$.
By the proof of Theorem~\ref{RGWs_thm}\ref{div_it} in Section~\ref{MapSignDfn_subs}, 
we can also assume that $\codim\,h_1\!=\!2$.
We denote by $\dim\bh$ the sum of the dimensions of the domains~$H_i$
of the components of~$\bh$ and~set
$$\codim\,\bh=2nl-\dim\bh\,.$$
Let \hbox{$h\!:Z\!\lra\!X^l$} be a smooth map from a manifold of dimension
$\dim\bh\!-\!2$ covering~$\Om(h_1\!\times\!\ldots\!\times\!h_l)$,
\begin{gather*}
\ev^+\equiv\prod_{i\in[l]}\!\ev^+_i\!: \ov\fM_{k,l}\big(\wt{B};J\big)\lra X^l,\\
\ev_{k,l;h}\!:\ov\cZ_{k,l;h}(\wt{B};J)\!\equiv\!\big\{
(\u,z)\!\in\!\ov\fM_{k,l}(\wt{B};J)\!\times\!Z\!:\ev^+(\u)\!=\!h(z)\big\}
\lra Y^k.
\end{gather*}
For a stratum~$\cS$ of $\ov\fM_{k,l}^*(\wt{B};J)$, let
$$\cS_{\bh}\equiv \ov\cZ_{k,l;\bh}(\wt{B};J)\cap
\bigg(\!\cS\!\times\!\prod_{i=1}^l\!H_i\!\bigg)$$
and $\fc(\cS)\!\in\!\Z^{\ge0}$ be the number of nodes of the domains of the elements of~$\cS$.\\

Since the space $\ov\fM_{k,l}(\wt{B};J)$ is compact,
\BE{JakePseudo_e5}
\Om\big(\ev_{k,l;\bh}\big|_{\cZ_{k,l;\bh}^{\bigst}(\wt{B};J)}\big)
\subset \ev_{k,l;\bh}\big(\ov\cZ_{k,l;\bh}(\wt{B};J)
\!-\!\cZ_{k,l;\bh}^{\bigst}(\wt{B};J)\!\big)
\cup  \ev_{k,l;h}\big(\ov\cZ_{k,l;h}(\wt{B};J)\!\big)\,.\EE
The subspace
$$\ov\cZ_{k,l;\bh}^{\st}(\wt{B};J)
\!-\!\cZ_{k,l;\bh}^{\bigst}(\wt{B};J)\subset
\ov\fM_{k,l}(\wt{B};J)\!\times\!\!\prod_{i=1}^l\!H_i$$
consists of the spaces~$\cS_{\bh}$ corresponding to the strata~$\cS$ 
of $\ov\fM_{k,l}(\wt{B};J)$ with either $\fc(\cS)\!\ge\!2$ nodes
or $\ep_{l^*}(\cS)\!\cong\!2,3$ mod~4 
and of the subspace $\ov\fM_{k,l}(\wt{B};J)\!-\!\ov\fM_{k,l}^*(\wt{B};J)$
of multiply covered maps meeting the pseudocycles in~$\bh$.
By~\eref{dimkbhcond_e} and~$\cS$ being of codimension~$\fc$ in~$\ov\fM_{k,l}^*(\wt{B};J)$,
the dimension of~$\cS_{\bh}$ is at most $nk\!-\!2$ if $\fc(\cS)\!\ge\!2$.
By~\eref{dimkbhcond_e} and the reasoning in the proof of Proposition~5.2 in~\cite{RealWDVV},
the image under~$\ev_{k,l;\bh}$ of 
the subspace $\ov\fM_{k,l}(\wt{B};J)\!-\!\ov\fM_{k,l}^*(\wt{B};J)$
meeting the pseudocycles in~$\bh$ can also be covered by a smooth map
from a manifold of dimension~$nk\!-\!2$.
The same applies to the last set in~\eref{JakePseudo_e5}.\\

Suppose $\cS$ is a codimension~1 stratum of $\ov\fM_{k,l}^*(\wt{B};J)$ 
with $\ep_{l^*}(\cS)\!\cong\!2,3$ mod~4.
For $r\!=\!1,2$, let 
$$k_r=|K_r(\cS)|, ~~ l_r=|L_r(\cS)|,~~  
L_r^*=L_r(\cS)\!\cap\![l^*], ~~ \wt{B}_r=\wt{B}_r(\cS),~~\bh_r\!=\!(h_i)_{i\in L_r(\cS)}.$$  
In particular,
\BE{JakePseudo2_e8b}\begin{split}
&\hspace{.4in}k_1\!+\!k_2=k, \quad l_1\!+\!l_2=l, \quad 
\codim\,\bh_1\!+\!\codim\,\bh_2=\codim\,\bh;\\
&\blr{c_1(X,\om),\wt{B}_1}\!+\!\blr{c_1(X,\om),\wt{B}_2}
=\blr{c_1(X,\om),\wt{B}}
=(n\!-\!1)k\!+\!2\,\codim\,\bh\!-\!2l\!+\!3\!-\!n.
\end{split}\EE
By the definition~$\ep_{l^*}(\cS)$ and~\eref{Rcond_e0},
\BE{JakePseudo2_e8c}\begin{split}
\blr{c_1(X,\om),\wt{B}_2}\!-\!(n\!-\!1)k_2\!-\!2\big(l_2\!-\!|L_2^*|\big)
&\cong2,3\quad\tn{mod}~4,\\
\blr{c_1(X,\om),\wt{B}_2}\!-\!(n\!-\!1)k_2\!-\!2\big(\codim\,\bh_2\!-\!l_2\big)
&\cong2,3\quad\tn{mod}~4.
\end{split}\EE
Along with \eref{JakePseudo2_e8b} and $\lr{c_1(X,\om),B'}\!\in\!2\Z$ if $n\!=\!3$,
the second equality in~\eref{JakePseudo2_e8c} gives 
\BE{JakePseudo2_e8d}\begin{split}
\blr{c_1(X,\om),\wt{B}_1}\!-\!(n\!-\!1)k_1\!-\!2\big(\codim\,\bh_1\!-\!l_1\big)
&\cong2,3\quad\tn{mod}~4,\\
\blr{c_1(X,\om),\wt{B}_1}\!-\!(n\!-\!1)k_1\!-\!2\big(l_1\!-\!|L_1^*|\big)
&\cong2,3\quad\tn{mod}~4.
\end{split}\EE
By the first statement in~\eref{JakePseudo2_e8c} and
the second in~\eref{JakePseudo2_e8d},
$$\big(B_r,k_r,l_r,|L_r^*|\big)\neq(0,0,0,0),(0,0,1,1) 
\qquad\forall\,r=1,2.$$
Since the image of~$h_i$ with $i\!\in\![l]\!-\![l^*]$ is disjoint from~$Y$,
$\cS_\bh\!=\!\eset$ if 
$$\big(B_r,k_r,l_r,|L_r^*|\big)=(0,0,1,0).$$ 
Thus, we can thus assume either
$\wt{B}_r\!\neq\!0$, or $k_r\!+\!2l_r\!\ge\!3$, or $k_r\!=\!2$ for each $r\!=\!1,2$.\\

By the definition of~$\P_1^1$, $l_1\!\ge\!1$.
Suppose $\wt{B}_2\!=\!0$, $l_2,l_2^*\!=\!0$, and $k_2\!=\!2$.
The restriction of~$\ev_{k,l;\bh}$ to~$\cS_{\bh}$ then factors~as 
$$\cS_{\bh}\lra  \cZ_{k-1,l;\bh}^{\st}(\wt{B};J)\!\times\!\fM_{3,0}(0;J)
\lra Y^{k-2}\!\times\!\De_Y \lra Y^k\,,$$
where $\De_Y\!\subset\!Y^2$ is the diagonal.
Thus,  $\ev_{k,l;\bh}(\cS_{\bh})$ is again contained in a smooth manifold of dimension 
\hbox{$nk\!-\!2$}.\\

Suppose now that either $\wt{B}_r\!\neq\!0$ or $2l_r\!+\!k_r\!\ge\!3$ for each $r\!=\!1,2$.
The restriction of~$\ev_{k,l;\bh}$ to~$\cS_{\bh}$ then factors~as
$$\xymatrix{\cS_{\bh}\ar[r] &
\cZ_{k_1,l_1;\bh_1}^{\st}(\wt{B}_1;J)
\times\cZ_{k_2,l_2;\bh_2}^{\st}(\wt{B}_2;J)
\ar@<-3ex>[d]_{\ev_{k_1,l_1;\bh_1}}\ar@<2ex>[d]^{\ev_{k_2,l_2;\bh_2}}\\
& Y^{k_1}\!\times\!Y^{k_2} \ar[r]& Y^{k}.}$$
Thus, $\ev_{k,l;\bh}(\cS_{\bh})$ is covered by a smooth map from a manifold
of dimension
\BE{JakePseudo2_e15}\begin{split}
&\dim\,\cZ_{k_r,l_r;\bh_r}^{\st}(\wt{B}_r;J)\!+\!\dim Y^{k_{3-r}}
=\ell_{\om}(\wt{B}_r)\!+\!k_r\!\!+\!2l_r\!-\!\codim\,\bh_r\!+\!nk_{3-r}\\
&\hspace{1in}=\blr{c_1(X,\om),\wt{B}_r}\!-\!(n\!-\!1)k_r\!\!+\!2l_r\!-\!\codim\,\bh_r\!+\!nk
\!+\!n\!-\!3
\end{split}\EE
for $r\!=\!1,2$.
By~\eref{JakePseudo2_e8b}, the second statement in~\eref{JakePseudo2_e8c},
and the first in~\eref{JakePseudo2_e8d},
\BE{JakePseudo2_e17}\blr{c_1(X,\om),\wt{B}_r}\!-\!(n\!-\!1)k_r\!\!+\!2l_r\!-\!\codim\,\bh_r
\!+\!n\!-\!3\le -2\EE
for either $r\!=\!1$ or $r\!=\!2$.
Thus, $\ev_{k,l;\bh}(\cS_{\bh})$ is covered by a smooth map from a manifold
of dimension $nk\!-\!2$ in this case as well.\\

The above shows that the limit set of the map~$\ev_{k,l;\bh}$ in~\eref{stevkdfn_e}
is covered by a smooth map from a manifold of dimension~$nk\!-\!2$.
The same reasoning shows that the limit set of the map~$\wt\ev_{k,l;\wt\bh}$ 
in~\eref{stevkdfn_e} is covered by a smooth map from a manifold of dimension~$nk\!-\!1$.
In particular, \eref{JakePseudo2_e8b}, the second statement in~\eref{JakePseudo2_e8c}, 
the first statement in~\eref{JakePseudo2_e8d}, and~\eref{JakePseudo2_e17} 
hold with~$\bh$ replaced by~$\wt\bh$.
Since $\wt{h}_i$ is a pseudocycle equivalence in $X\!-\!Y$ if $n\!=\!3$ and 
$\dim\wt{h}_i\!=\!3$,  the image of~$\wt{h}_i$ with $i\!\in\![l]\!-\![l^*]$ 
is disjoint from~$[0,1]\!\times\!Y$ in this case.
The analogue of~\eref{JakePseudo2_e15} is now 
$$\dim\,\cZ_{k_r,l_r;\wt\bh_r}^{\st}(\wt{B}_r;\wt{J})\!+\!\dim Y^{k_{3-r}}
=\blr{c_1(X,\om),\wt{B}_r}\!-\!(n\!-\!1)k_r\!\!+\!2l_r\!-\!\codim\,\wt\bh_r\!+\!nk
\!+\!n\!-\!2.$$

\subsection{Signs at immersions}
\label{SignImmer_subs}

Suppose $(X,\om,\phi)$, $n$, $Y$, $k,l,l^*$, $\wt{B}$,  $J$, and~$\fp$ 
are as at the beginning of Section~\ref{MapSignDfn_subs} with $l\!\ge\!2$.
Let $\u$ be as in~\eref{udfn_e} so that $[\u]\!\in\!\fM^*_{k,l}$
and $u$ is an immersion.
We denote by~$\cC$  the marked domain of~$\u$ as in~\eref{cCdfn_e}.
The conjugation $\nd\phi$ on~$TX$ descends to a conjugation~$\vph$ on the normal bundle
$$\cN u\equiv\frac{u^*TX}{\Im\,\nd u}\lra T\P^1$$
for the immersion~$u$ so that 
\BE{cNvphidfn_e}0\lra (T\P^1,\tau)\stackrel{\nd u}\lra u^*(TX,\nd\phi)
\stackrel\pi\lra (\cN u,\vph)\lra0\EE
is an exact sequence of real bundle pairs over $(\P^1,\tau)$.
Define
$$\cN_{\u} \equiv \cN_{\u}^{\C}\!\oplus\!\cN_{\u}^{\R}
\equiv \big(\cN u\big|_{z^+_3}\!\oplus\!\ldots\!\oplus\!\cN u\big|_{z^+_l}\big)
\!\oplus\!
\big(\cN u^{\vph}\big|_{x_1}\!\oplus\!\ldots\!\oplus\!\cN u^{\vph}\big|_{x_k}\big).$$

\vspace{.2in}

We denote by $\bp$ the real CR-operator on the rank~1 real bundle $(T\P^1,\nd\tau)$ 
over~$(\P^1,\tau)$ determined by the holomorphic $\bp$-operator on~$T\P^1$.
The linearization $D_{\u}\!\equiv\!D_{J;u}^{\phi}$ of the $\bp_J$-operator on the space of 
real maps from~$(\P^1,\tau)$ descends to a real CR-operator~$D_{\u}''$ 
on the real bundle pair $(\cN u,\vph)$ and induces an exact triple
\BE{sesDNc_e} 0\lra \bp \lra D_{\u}\lra D_{\u}''\lra0\EE
of real CR-operators.
Let
$$\ev_{\u}''^{\C}\!:\ker D_{\u}''\lra \cN_{\u}^{\C} \qquad\hbox{and}\qquad
\ev_{\u}''^{\R}\!:\ker D_{\u}''\lra \cN_{\u}^{\R}$$
be the evaluations at the marked points $z^+_3,\ldots,z^+_l$ and
$x_1,\ldots,x_k$, respectively.\\

Let $\os_0(TS^1)$ be the $\OSpin$-structure on $S^1\!\subset\!\bD^2_+\!\subset\!\P^1$
determined by the standard orientation of~$S^1$ as in Example~\ref{LBSpinStr_eg} and 
$\ce_{\R}$ be the real part of the short exact sequence~\eref{cNvphidfn_e}
as in~\eref{RBPsesdfn_e2}.
By the RelSpinPin~\ref{RelSpinPinSES_prop} property in Section~\ref{RelSpinPinProp_subs}, 
the relative $\Pin^{\pm}$-structure (resp.~$\OSpin$-structure)~$\fp$ on $Y\!\subset\!X$
determines a relative $\Pin^{\pm}$-structure (resp.~$\OSpin$-structure)~$\fp_u''$ on 
the vector bundle $\cN u^{\vph}$ over~$S^1\!\subset\!\P^1$ so~that 
\BE{WelSolComp2_e3a}\bllrr{\os_0(TS^1),\fp_{\u}''}_{\!\ce_{\R}}=
\big\{u|_{S^1}\big\}^*\fp
\in\cP_{\P^1}^{\pm}\big(\{u|_{S^1}\}^*TY\big)
\quad\big(\hbox{resp.}~\OSp_{\P^1}\big(\{u|_{S^1}\}^*TY\big)\big).\EE
By Theorem~\ref{CROrient_thm}, $\fp_{\u}''$ in turn determines a homotopy class 
$$\fo_{\fp;\u}''\equiv \fo_{\cC;\fp_{\u}''}\big(\cN u,\vph\big)$$
of isomorphisms $\la(\ker D_{\u}'')\!\approx\!\la(\cN_{\u}^{\R})$.\\

Let $\bh$ as in~\eref{bhdfn_e} with $l'\!=\!l$
be a tuple of smooth maps from oriented manifolds
in general position so that~\eref{dimkbhcond_e} holds.
Suppose in addition~that 
\BE{wtudfn_e}\wt\u\equiv\big([\u],(y_i)_{i\in[l]}\big) \in \cZ_{k,l;\bh}^*(\wt{B};J)\subset
\fM_{k,l}^*\!\times\!\prod_{i=1}^l\!H_i\EE
with $\u$ as above is a regular point of~\eref{JakePseudo_e} and 
the homomorphisms~\eref{JakePseudo_e10b} with $i\!=\!1,2$ are isomorphisms.
Let
$$g_{\wt\u}\!:\bigoplus_{i=3}^lT_{y_i}H_i\lra \cN_{\u}^{\C}$$
be the direct sum of the compositions of the differentials $\nd_{y_i}h_i$
with the projections to~$\cN u|_{z_i^+}$ and
$$\cK_{\wt\u}\equiv\big\{(\xi,w)\!\in\!(\ker D_{\u}'')\!\oplus\!\bigoplus_{i=3}^lT_{y_i}H_i\!:
\ev_{\u}''^{\C}(\xi)\!=\!g_{\wt\u}(w)\big\}.$$
The homotopy class~$\fo_{\fp;\u}''$  of isomorphisms $\la(\ker D_{\u}'')\!\approx\!\la(\cN_{\u}^{\R})$
and the orientations
of~$H_i$ and~$\cN u$ determine a homotopy class~$\fo_{\fp;\bh;\wt\u}$ of isomorphisms
$\la(\cK_{\wt\u})\!\approx\!\la(\cN_{\u}^{\R})$.
By~\eref{dimkbhcond_e} and the transversality assumptions,
the  homomorphism 
$$\ev_{\wt\u}''^{\R}\!:\cK_{\wt\u}\lra \cN_{\u}^{\R} $$
induced by~$\ev_{\u}''^{\R}$ is an isomorphism.
We set~$\fs_{\fp;\bh}(\u)$\nota{sphu@$\fs_{\fp;\bh}(\u)$} 
to be $+1$ if this isomorphism
lies in~$\fo_{\fp;\bh;\wt\u}$  and $-1$ otherwise.
The next proposition compares this sign with the sign~$\fs_{\fp;l^*;\bh}(\u)$
defined in Section~\ref{MapSignDfn_subs}.
Let $\de_{l^*}^c(\wt\u)\!\equiv\!\de_{l^*}^c(\cC)$ be as in~\eref{decldfn_e}.

\begin{prp}\label{cNsgn_prp}
Suppose $(X,\om,\phi)$ is a real symplectic manifold of dimension~$2n$, 
$Y\!\subset\!X^{\phi}$ is a topological component,
$k,l\!\in\!\Z^{\ge0}$ with \hbox{$l\!\ge\!2$}, $l^*\!\in\![l]$, 
$\wt{B}\!\in\!\wt{H}_2^{\phi}(X,Y)$, $J\!\in\!\cJ_{\om}^{\phi}$ is generic, and
that either~$\fp$ is a relative $\Pin^{\pm}$-structure on $Y\!\subset\!X$
and the pair $(k,B)$ satisfies the condition~\eref{BKcond_e0} or
$\fp$ is a relative $\OSpin$-structure on~$Y$.
Let $\bh$ as in~\eref{bhdfn_e} with $l'\!=\!l$ 
be a tuple of smooth maps from oriented even-dimensional manifolds  
 in general position so that~\eref{dimkbhcond_e} holds and
\hbox{$\wt\u\!\in\!\cZ_{k,l;\bh}^*(\wt{B};J)$} be as in~\eref{wtudfn_e} so that 
the homomorphisms~\eref{JakePseudo_e10b} with $i\!=\!1,2$ are 
orientation-preserving isomorphisms.
The signs~$\fs_{\fp;l^*;\bh}(\u)$ and $\fs_{\fp;\bh}(\u)$ are the same
if and only~if
\BE{cNsgn_e}\de_{l^*}^c(\wt\u)\!+\!(n\!-\!1)\binom{k\!+\!1}2
\!+\!(n\!-\!1)+2\Z =\blr{w_2(X),\wt{B}}\in\Z_2.\EE
\end{prp}

\begin{rmk}\label{cNsgn_rmk}
Let $\bh''\!=\!(\bh_3,\ldots,\bh_l)$.
The sign~$\fs_{\fp;\bh}(\u)$ depends only on the image~$[\u'']$ of~$[\u]$ in
$\fM_{k,l-2}^*(\wt{B};J)$
obtained by forgetting the marked points $z_1^{\pm}$ and $z_2^{\pm}$
and possibly on the choice of the disk $\bD^2_+\!\subset\!\P^1$ bounding the real locus~$S^1$
of~$\tau$ (as determined by the condition $z_1^+\!\in\!\bD^2_+$).
If the sign $\fs_{\fp;l^*-2;\bh'}(u(\P^1)\!)$ is well-defined as in Section~\ref{MapSignDfn_subs}, 
Proposition~\ref{cNsgn_prp} provides a comparison between $\fs_{\fp;l^*-2;\bh'}(u(\P^1)\!)$
and the sign $\fs_{\fp;\bh''}(\u'')\!\equiv\!\fs_{\fp;\bh}(\u)$ obtained from
the normal bundle to the curve~$u(\P^1)$ and a choice of the half of this curve 
bounding~$u(S^1)$.
\end{rmk}

With the assumption as in Proposition~\ref{cNsgn_prp}, let $\hb'\!=\!(h_1,h_2)$ and
$$\wt\u'\equiv\big([\u],(y_1,y_2)\!\big) \in
\cZ_{k,l;\bh'}^*\equiv\cZ_{k,l;\bh'}^*(\wt{B};J)\subset
\fM_{k,l}^*\!\times\!Y_1\!\times\!Y_2 \equiv
\fM_{k,l}^*(\wt{B},J)\!\times\!Y_1\!\times\!Y_2\,.$$
We can assume that $z_1^+\!=\!0$ and $z_2^+\!\in\!\R^+$.
Let
\begin{gather}\label{NcNdfn_e}
V_{\u}'\equiv T_{z^+_3}\P^1\!\oplus\!\ldots\!\oplus\!T_{z^+_l}\P^1
\!\oplus\! T_{x_1}S^1\!\oplus\!\ldots\!\oplus\!T_{x_k}S^1,\\
\notag
\p\!\equiv\!\ev_{k,l;\bh'}(\wt\u') \in X_{k,l-2}\!\equiv\!Y^k\!\times\!X^{l-2}.
\end{gather}
By the short exact sequence above~\eref{OrientSubs_e3}, there is a natural isomorphism 
\BE{stablize_e2} T_{[\u]}\fM^*_{k,l}\approx 
\big(\!\ker D_{\u}\big)\!\oplus\!T_{z^+_2}\R^+\!\oplus\!V_{\u}'\,.\EE

\vspace{.2in}

By the transversality of $u$ to~$H_1$ and~$H_2$, 
$$T_{\wt\u'}\cZ_{k,l;\bh'}^*\!\cap\!\big(\!\ker\bp\big)=\{0\} 
\qquad\hbox{and}\qquad
T_{\wt\u'}\cZ_{k,l;\bh'}^*\!\cap\!T_{z^+_2}\R^+=\{0\}\,.$$
Thus, the isomorphism~\eref{stablize_e2} and 
the exact triple~\eref{sesDNc_e} of real CR-operators induce a short exact sequence
\BE{cNCzdfn_e}0\lra V_{\u}' \lra T_{\wt\u'}\cZ_{k,l;\bh'}^* \lra \ker D_{\u}''\lra0\EE
of vector spaces and a homotopy class of isomorphisms
$\la(T_{\wt\u'}\cZ_{k,l;\bh'}^*)\!\approx\!\la(V_{\u}')\!\otimes\!\la(D_{\u}'')$.
Combining it with the homotopy class~$\fo_{\fp;\u}''$  of isomorphisms 
$\la(\ker D_{\u}'')\!\approx\!\la(\cN_{\u}^{\R})$ above, the orientation of~$\cN_{\u}^{\C}$,
and an identification as in~\eref{RBPsesdfn_e4},
we obtain a homotopy class~$\fo_{\fp;\wt\u'}$ of isomorphisms
\BE{cNorient_e0}\la\big(T_{\wt\u'}\cZ^*_{k,l;\bh'}\big)
\approx\la(V_{\u}')\!\otimes\!\la(D_{\u}'')
\approx\la(V_{\u}')\otimes\la(\cN_{\u})=\la\big(T_{\p}X_{k,l-2}\big)\EE
from $\la(T_{\wt\u'}\cZ^*_{k,l;\bh'})$ to $\la(T_{\p}X_{k,l-2})$.
Another homotopy class~$\fo_{\fp;l^*;\bh'}|_{\wt\u'}$ of such isomorphisms
is obtained as at the beginning of Section~\ref{MapSignDfn_subs}.

\begin{lmm}\label{cNorient_lmm}
The relative orientations $\fo_{\fp;l^*;\bh'}|_{\wt\u'}$ and~$\fo_{\fp;\wt\u'}$
of the map~$\ev_{k,l;\bh'}$ as in~\eref{JakePseudo_e} at~$\wt\u'$ 
are the same if and only~if \eref{cNsgn_e} holds.
\end{lmm}

\begin{proof}[{\bf{\emph{Proof of Proposition~\ref{cNsgn_prp}}}}]
Since the dimensions of $H_3,\ldots,H_l$ are even, $\fo_{\fp;l^*;\bh}|_{\wt\u}$
is the relative orientation of the map~$\ev_{k,l;\bh}$ as in~\eref{JakePseudo_e}
induced by~$\fo_{\fp;l^*;\bh'}|_{\wt\u'}$ by viewing $\cZ_{k,l;\bh}^*(\wt{B};J)$ 
as the fiber product~of
$$(\ev_3^+,\ldots,\ev_l^+)\!: \cZ_{k,l;\bh'}\lra X^{l-2}$$
with $h_3\!\times\!\ldots\!\times\!h_l$.
We denote by~$\fo_{\fp;\bh}|_{\wt\u}$ the relative orientation of~$\ev_{k,l;\bh}$
similarly induced by~$\fo_{\fp;\wt\u'}$.
The signs~$\fs_{\fp;l^*;\bh}(\u)$ and~$\fs_{\fp;\bh}(\u)$ are the same
are the signs of~$\ev_{k,l;\bh}$ with respect to the two relative orientations.
Along with Lemma~\ref{cNorient_lmm},  this implies the claim.
\end{proof}

\begin{proof}[{\bf{\emph{Proof of Theorem~\ref{RGWs_thm2}\ref{RGWprod_it1}}}}]
We continue with the notation and setup in the statement of this theorem
and deduce it from Proposition~\ref{cNsgn_prp}; see also Remark~\ref{cNsgn_rmk}.
It remains to verify the signs in the last two cases. 
Let 
$$\pi',\pi''\!:X\lra X'\!\equiv\!\P^1,X''$$
be the component projections.
We take the almost complex structure~$J$ on~$X$ to be $J_{\P^1}\!\times\!J''$,
with $J_{\P^1}$ denoting the standard complex structure on~$\P^1$ and 
$J''\!\in\!\cJ_{\om''}^{\phi''}$ generic.
For each $i\!\in\![l]$, we take \hbox{$h_i\!=\!h_i'\!\times\!h_i''$} with 
$h_i'\!=\!\id_{\P^1}$ if $k\!=\!3$ or $i\!\ge\!2$ and 
$h_i'$ being the inclusion of a generic point $p_1'\!\in\!\P^1$ 
if $k\!=\!1$ and $i\!=\!1$.
Let $\p$ as in~\eref{bpRdfn_e} be generic with 
\hbox{$p_i^{\R}\!=\!(p_i'^{\R},p_i''^{\R})$}.
Define
$$\bh\equiv(h_i)_{i\in[l]}, \qquad \bh''\equiv(h_i'')_{i\in[l]}, \qquad
\p''\equiv(p_i''^{\R})_{i\in[k]}\in Y''^k\,.$$
Fix $C\!\in\!\fM_{\p;\bh}(\wt{B};J)$.\\

Let $\u$ and~$\wt\u$ be tuples as in~\eref{udfn_e} and~\eref{wtudfn_e}
so that their map component $u\!\equiv\!(\id_{\P^1},u'')$ 
satisfies~\eref{Cudfn_e} and~\eref{Cudfn_e4}.
The map~$u''$ then satisfies~\eref{Cudfn_e} and~\eref{Cudfn_e4} 
with~$u$, $\wt{B}$, and~$C$ replaced by~$u''$, $\wt{B}''$, and
$$C''\equiv \pi''(C)\in \fM_{\p'';\bh''}(\wt{B}'';J'')\,,$$
respectively.
Let $\cC$ as in~\eref{cCdfn_e} be the marked domain of~$\u$.
We denote by~$\u''$ and~$\wt\u''$ the tuples obtained from $\u$ and~$\wt\u$
by replacing the map component~$u$ by~$u''$ and the points $y_i\!\in\!H_i$
by their projections~$y_i''$ to the domains~$H_i''$ of~$h_i''$.
The inclusion $\pi''^*TX''\!\lra\!TX$ induces an isomorphism
$$u''^*(TX'',\nd\phi'')\lra \big(\cN u,\vph)$$
of real bundle pairs over~$(\P^1,\tau)$.
This isomorphism identifies the real CR-operator~$D_{\u''}$ with~$D_{\u}''$,
the $\OSpin$-structures~$u''^*\os''$ and~$\os_{\u}''$, 
the differential~$\nd_{y_i''}h_i''$ with the composition of the restriction
of~$\nd_{y_i}h_i$ to~$T_{y_i''}X''$ with the projection to~$\cN u|_{h_i(y_i)}$,
and~$T_{\p''}Y''^k$ with~$\cN_{\u}^{\R}$.\\

Suppose $k\!=\!3$. The natural isomorphism
\BE{RGWprod_e9}T_{[\cC]}\cM_{k,l}^{\tau}\approx 
T_{z_1^+}\P^1\!\oplus\!\ldots\!\oplus\!T_{z_l^+}\P^1\EE
is then orientation-preserving if $\de_{l^*}(\cC)\!\in\!2\Z$.
Thus, the isomorphism
\begin{gather*}
\begin{split}
T_{[\wt\u'']}\cZ_{k,l;\bh''}^*\big(\wt{B}'';J''\big)&\equiv\!\big\{
\big(\xi,(w_i')_{i\in[l]},(w_i'')_{i\in[l]}\big)\!\in\!(\ker D_{\u''})\!\oplus\!
\bigoplus_{i=1}^lT_{z_i^+}\P^1\!\oplus\!\bigoplus_{i=1}^l\!T_{y_i''}H_i''\!:\\
&\hspace{.2in}
\xi(z_i^+)\!+\!\nd_{z_i^+}u''(w_i')\!=\!\nd_{y_i''}h_i''(w_i'')~\forall\,i\!\in\![l]\big\}
\lra \cK_{\wt\u}\subset(\ker D_{\u}'')\!\oplus\!\bigoplus_{i=1}^l\!T_{y_i}H_i,
\end{split}\\
\big(\xi,(w_i')_{i\in[l]},(w_i'')_{i\in[l]}\big) \lra 
\big(\xi,(w_i',w_i'')_{i\in[l]} \big),
\end{gather*}
identifies the homotopy class $\fo_{\os'';l^*;\bh''}$ of isomorphisms 
$\la(T_{[\wt\u'']}\cZ_{k,l;\bh''}^*)\!\approx\!\la^{\R}_{\u''}(X'')$
with the homotopy class $\fo_{\os;\bh;\wt\u}$ of isomorphisms 
$\la(\cK_{\wt\u})\!\approx\!\la(\cN_{\u}^{\R})$ if and only if 
$\de_{l^*}(\cC)\!\in\!2\Z$.
This isomorphism also intertwines the~maps 
$$\nd_{\wt\u''}\ev_{k,l;\bh''}\!: 
T_{\wt\u''}\cZ_{k,l;\bh''}^*\big(\wt{B}'';J''\big)\lra T_{\p''}Y''^k 
\qquad\hbox{and}\qquad \ev_{\wt\u}''^{\R}\!:\cK_{\wt\u}\lra \cN_{\u}^{\R}.$$
We conclude that 
$$\fs_{\os'';l^*;\bh''}(C'')\equiv \fs_{\os;l^*;\bh''}(\u'')
=(-1)^{\de_{l^*}(\cC)}\fs_{\os;\bh}(\u).$$
By Proposition~\ref{cNsgn_prp}, 
$$\fs_{\os;l^*;\bh}(C)\equiv \fs_{\os;l^*;\bh}(\u)
=(-1)^{\de_{l^*}(\cC)}\fs_{\os;\bh}(\u).$$
Combining the last two statements, we obtain 
$\fs_{\os;l^*;\bh}(C)\!=\!\fs_{\os'';l^*;\bh''}(C'')$.
This establishes the equality in the second case of 
Theorem~\ref{RGWs_thm2}\ref{RGWprod_it1}.\\

Suppose $k\!=\!1$. The right-hand side of the analogue of the isomorphism~\eref{RGWprod_e9}
is $T_{z_2^+}\P^1\!\oplus\!\ldots\!\oplus\!T_{z_l^+}\P^1$ in this case.
The isomorphism $T_{[\wt\u'']}\cZ_{k,l;\bh''}^*\!\approx\!\cK_{\wt\u}$ is as 
in the $k\!=\!3$ case, but with $w_1'\!=\!0$.
The rest of the proof of the $k\!=\!3$ case applies without any changes
and yields  $\fs_{\os;l^*;\bh}(C)\!=\!\fs_{\os'';l^*;\bh''}(C'')$.
This establishes the equality in the third case of 
Theorem~\ref{RGWs_thm2}\ref{RGWprod_it1}.
\end{proof}

\subsection{Proof of Lemma~\ref{cNorient_lmm}}
\label{OrientImmer_subs}

We continue with the notation and setup above Lemma~\ref{cNorient_lmm}.
We orient the individual summands of~$V_{\u}'$ in~\eref{NcNdfn_e}
by the complex orientation of~$\P^1$ and the standard orientation of $S^1\!\subset\!\bD^2_+$.
We denote by~$\fo_{\u}'$ the orientation of~$V_{\u}'$ obtained from these orientations
by ordering the summands~$T_{x_i}S^1$ by the position of the marked points~$x_i$
on~$S^1$ starting from~$x_1$.
Fix an orientation $\fo_{\u}''$ on~$\cN_{\u}$. 
Let~$\fo_{\u}$ be the orientation of~$T_{\p}X_{k,l-2}$ determined by 
the orientations~$\fo_{\u}'$ and~$\fo_{\u}''$ via the last identification 
in~\eref{cNorient_e0}.\\

The orientations~$\fo_{\u}$ on $T_{\p}X_{k,l-2}$ and~$\fo_{\u}''$ on~$\cN_{\u}$ 
determine orientations~$\wh\fo_{\fp;\u}^D$ on~$\ker D_{\u}$ and~$\wh\fo_{\u}''$ on~$\ker D''_{\u}$
via the homotopy classes $\fo_{\fp;\u}^D$ and $\fo_{\fp;\u}''$ of isomorphisms
$$\la\big(D_{\u}\big)\approx\la(T_{\p}X_{k,l-2}) \qquad\hbox{and}\qquad
\la\big(D_{\u}''\big)\approx\la(\cN_{\u}),$$ 
respectively.
Let
$$\wh\fo_{\u}'\equiv \fo_0(T\P^1,\nd\tau;\fo_{S^1})$$
be the orientation of $\ker\bp$ determined by the OSpin-structure~$\os_0(TS^1)$ 
on~$S^1$ as in Theorem~\ref{CROrient_thm}\ref{CROrientSpin_it}.\\

The orientation $\fo_{\u}$ on $T_{\p}X_{k,l-2}$ and the symplectic orientation~$\fo_{\om}$ on~$X$
determine orientations~$\wh\fo_{\fp;l^*}$ on~$T_{[\u]}\fM^*_{k,l}$ 
and~$\wh\fo_{\fp;l^*;\bh'}$ on~$T_{\wt\u'}\cZ^*_{k,l;\bh'}$ via
the homotopy classes~$\fo_{\fp;l^*}|_{[\u]}$ and~$\fo_{\fp;l^*;\bh'}|_{\wt\u'}$ 
of isomorphisms
$$\la\big(T_{[\u]}\fM^*_{k,l}\big)
\approx\la\big(T_{y_1}X\!\oplus\!T_{y_2}X\!\oplus\!T_{\p}X_{k,l-2}\big)
\quad\hbox{and}\quad
\la\big(T_{\wt\u}'\cZ^*_{k,l;\bh'}\big)\approx\la\big(T_{\p}X_{k,l-2}\big),$$
respectively.
We denote by $\wh\fo_{\fp}$ the orientation on $T_{\wt\u'}\cZ^*_{k,l;\bh'}$
similarly determined by the orientations~$\fo_{\u}$  and~$\fo_{\fp;\wt\u'}$.
The relative orientations~$\fo_{\fp;l^*;\bh'}|_{\wt\u'}$ and~$\fo_{\fp;\wt\u'}$
of~\eref{JakePseudo_e} at~$\wt\u'$ are the same if and only~if 
the orientations~$\wh\fo_{\fp;l^*;\bh'}$ and~$\wh\fo_{\fp}$ 
of~$T_{\wt\u'}\cZ^*_{k,l;\bh'}$ are the~same.\\

By the sentence containing~\eref{cMorientisom_e},
the isomorphism in the right column of the first diagram in Figure~\ref{WelSolComp2_fig} 
respects the opposite~$\ov{\fo_{\R}}$ of the standard orientation on $T_{z^+_2}\R^+$,
the orientation~$\fo_{\u}'$ on~$V_{\u}'$, and
the orientation~$\fo_{k,l;l^*}$ on~$T_{\cC}\cM^\tau_{k,l}$ 
if and only if $\de_{l^*}^c(\cC)$ is even. 
By~\eref{WelSolComp2_e3a} and the CROrient~\ref{CROSpinPinSES_prop}\ref{CROsesPin_it} 
property in Section~\ref{OrientPrp_subs1}, 
the left column in this diagram 
respects the orientations~$\wh\fo_{\u}'$,
$\wh\fo_{\fp;\u}^D$,  and~$\wh\fo_{\u}''$ if and only if $(n\!-\!1)\binom{k}{2}$ is even.
By the definition of the relative orientation~$\fo_{\fp;l^*}$ on~\eref{fMevdfn_e} 
via~\eref{OrientSubs_e3},
the middle row in this diagram respects the orientations~$\wh\fo_{\fp;\u}^D$,
$\wh\fo_{\fp;l^*}$, and~$\fo_{k,l;l^*}$.
Along with Lemma~\ref{3by3_lmm} and
$$\dim\ker D_{\u}''=n\!-\!1\!+\!\lr{c_1(X,\om),B}\!-\!2\,,$$
this implies that 
the middle column in the first diagram of Figure~\ref{WelSolComp2_fig} respects
the orientations $\wh\fo_{\u}'\!\oplus\!\ov{\fo_{\R}}\!\oplus\!\fo_{\u}'$, 
$\wh\fo_{\fp;l^*}$, and~$\wh\fo_{\u}''$ if and only~if  
the number
\BE{cNorient_e7}\begin{split}
&\de_{l^*}^c(\cC)\!+\!(n\!-\!1)\!\binom{k}{2}\!+\!\big(\!\dim\ker D_{\u}''\big)(k\!+\!1)\\
&\hspace{.3in}
\cong \de_{l^*}^c(\cC)\!+\!(n\!-\!1)\binom{k\!+\!1}{2}\!+\!(n\!+\!1)
\!+\!(k\!+\!1)\lr{c_1(X,\om),B}
\quad\hbox{mod}~2
\end{split}\EE
is even.\\

\begin{figure}
\begin{gather*}
\xymatrix{&0\ar[d] &0\ar[d] & 0\ar[d]\\
0\ar[r]&\ker\bp\ar[r]\ar[d] & 
\big(\!\ker\bp\big)\!\oplus\!T_{z^+_2}\R^+\!\oplus\!V_{\u}'
\ar[r]\ar[d]&T_{z^+_2}\R^+\!\oplus\!V_{\u}'\ar[r]\ar[d]&0\\
0\ar[r]&\ker D_{\u} \ar[r]\ar[d] & 
T_{[\u]}\fM^*_{k,l} \ar[r]^<<<<<<<<<<<{\nd_{\wt\u}\ff_{k,l}}\ar[d]& 
T_{\cC}\cM^\tau_{k,l}
 \ar[r]\ar[d]& 0\\
0\ar[r]&\ker D_{\u}''\ar@{=}[r]\ar[d] 
&\ker D_{\u}'' \ar[d]\ar[r]& 0\\
&0&0}\\
\xymatrix{&0\ar[d]&0\ar[d] & 0\ar[d]\\
0\ar[r]&V_{\u}'\ar[r]\ar[d] & 
\big(\!\ker\bp\big)\!\oplus\!T_{z^+_2}\R^+\!\oplus\!V_{\u}' \ar[r]\ar[d]&
\big(\!\ker\bp\big)\!\oplus\!T_{z^+_2}\R^+\ \ar[d]^{g_0}\ar[r]&0\\
0\ar[r]& T_{\wt\u'}\cZ^*_{k,l;\bh'}\ar[d] \ar@{^(->}[r] 
& T_{[\u]}\fM^*_{k,l} \ar[r]^g\ar[d]& 
\cN_{y_1}H_1\!\oplus\!\cN_{y_2}H_2 \ar[r]\ar[d]& 0\\
0\ar[r]&\ker D_{\u}''\ar@{=}[r]\ar[d] 
&\ker D_{\u}'' \ar[d]\ar[r]& 0\\
&0&0}
\end{gather*}
\caption{Commutative squares of exact rows and columns for the proof of Lemma~\ref{cNorient_lmm}.}
\label{WelSolComp2_fig}
\end{figure}

For $i\!=\!1,2$, let
$$\pi_i\!: T_{y_i}X\lra \cN_{y_i}H_i\!\equiv\!\frac{T_{y_i}X}{T_{y_i}H_i}$$
be the projection to the normal bundle of $H_i$ at~$y_i$.
We denote by~$\fo_{\cN;\wt\u'}$ the orientation on 
\hbox{$\cN_{y_1}H_1\!\oplus\!\cN_{y_2}H_2$}
induced by the symplectic orientation~$\fo_{\om}$ on~$X$ and the given orientations on~$H_1$
and~$H_2$.
Define 
\begin{alignat*}{2}
g_0\!:\big(\!\ker\bp\big)\!\oplus\!T_{z^+_2}\R^+&\lra\cN_{y_1}H_1\!\oplus\!\cN_{y_2}H_2,
&~~ g_0(\xi,v)&=
\big(\pi_1\big(\nd_{z_1^+}u(\xi(z_1^+)\!)\!\big),
\pi_2\big(\nd_{z_2^+}u(\xi(z_2^+)\!+\!v)\!\big)\!\big),\\
g\!:T_{[\u]}\fM^*_{k,l}&\lra\cN_{y_1}H_1\!\oplus\!\cN_{y_2}H_2,
&\quad g(\xi)&=
\big(\pi_1\big(\nd_{[\u]}\ev_1^+(\xi)\!\big),\pi_2\big(\nd_{[\u]}\ev_2^+(\xi)\!\big)\!\big).
\end{alignat*}
Since the isomorphisms~\eref{JakePseudo_e10b} with \hbox{$i\!=\!1,2$} are orientation-preserving,
the conclusion of Example~\ref{CRORientTP1_eg} implies~that 
$g_0$ is an isomorphism which respects the orientations $\wh\fo_{\u}'\!\oplus\!\ov{\fo_{\R}}$
and~$\fo_{\cN;\wt\u'}$.\\

The middle row in the second diagram of Figure~\ref{WelSolComp2_fig} respects 
the orientations~$\wh\fo_{\fp;l^*;\bh'}$ on $T_{\wt\u'}\cZ^*_{k,l;\bh'}$, 
$\wh\fo_{\fp;l^*}$ on $T_{[\u]}\fM^*_{k,l}$, and 
$\fo_{\cN;\wt\u'}$ on \hbox{$\cN_{y_1}H_1\!\oplus\!\cN_{y_2}H_2$}.
The left column in this diagram respects the orientations~$\fo_{\u}'$ on~$V_{\u}'$,
$\wh\fo_{\fp}$ on~$T_{\wt\u'}\cZ^*_{k,l;\bh'}$, 
and~$\wh\fo_{\u}''$ on~$\ker D_{\u}''$.
By the last sentence of the previous paragraph,
the isomorphism in the right column of this diagram respects the orientations
$\wh\fo_{\u}'\!\oplus\!\ov{\fo_{\R}}$ and~$\fo_{\cN;\wt\u'}$.
Combining this with the above conclusion concerning the middle column in the two diagrams
and Lemma~\ref{3by3_lmm}, we conclude that $\wh\fo_{\fp;l^*;\bh'}\!=\!\wh\fo_{\fp}$
if and only~if the right-hand side of~\eref{cNorient_e7} is an even number.
This establishes the claim.

\section{Counts of real rational curves}
\label{g0real_sec}

We recall Welschinger's definitions of invariant counts
of real rational irreducible $J$-holomorphic curves in compact 
real symplectic fourfolds and sixfolds~$(X,\om,\phi)$ in Section~\ref{g0comp_sec}.
We then compare these invariants with the invariant counts~\eref{RGWdfn_e}
of stable real $J$-holomorphic maps from~$(\P^1,\tau)$;
see Theorems~\ref{WelSolComp_thm} and~\ref{WelSolComp3_thm}.
Basic examples of such signed curve and map counts appear in Sections~\ref{Wel4eg_subs}
and~\ref{Wel6eg_subs}.
In Section~\ref{WelSolPf_thm}, we deduce Theorems~\ref{WelSolComp_thm} and~\ref{WelSolComp3_thm}
from Proposition~\ref{cNsgn_prp}.

\subsection{Definitions and main theorem}
\label{g0comp_sec}

Suppose first that $(X,\om,\phi)$ is a compact connected real symplectic fourfold 
and $Y\!\subset\!X^{\phi}$ is a connected component.
For $B\!\in\!H_2(X;\Z)$, we denote~by
$$B^2\equiv\blr{\PD_X(B),B}\in\Z$$
the homology self-intersection number of~$B$.
For $\wt{B}\!\in\!\wt{H}_2^{\phi}(X,Y)$ as in~\eref{wtBdfn_e}, let 
$\wt{B}^2\!=\!B^2$.
Since $\phi^*\om\!=\!-\om$, $\fd_Y(B)^2\!\in\!2\Z$ for every $B\!\in\!H_2(X;\Z)$.\\

Let $J\!\in\!\cJ_{\om}^{\phi}$, 
$C\!\subset\!X$ be a real rational irreducible $J$-holomorphic curve,
and $u$ be as in~\eref{Cudfn_e} and~\eref{Cudfn_e4}.
A point $x\!\in\!C$ is a \sf{simple node} if
$$\big|u^{-1}(x)\big|= 2 \qquad\hbox{and}\qquad
\bigoplus_{z\in u^{-1}(x)}\!\!\!\!\!\Im\,\nd_zu=T_xX\,.$$
A (simple) node $x$ of a real $J$-holomorphic curve~$C$ in~$X$ can be of 3 types: 
\begin{enumerate}[leftmargin=.3in]

\item[(E)]\label{E_it} $x\!\in\!X^{\phi}$ is an isolated point of 
$C^{\phi}\!\equiv\!C\!\cap\!X^{\phi}$,

\item[(H)]\label{H_it} $x\!\in\!X^{\phi}$ is a non-isolated point of~$C^{\phi}$,

\item[($\C$)]\label{C_it} $x\!\in\!X\!-\!X^{\phi}$ is a non-real point of~$C^{\phi}$.

\end{enumerate}
We denote by $\de_E(C)$ and $\de_H(C)$ the numbers of simple nodes of~$C$ of types~$E$
and~$H$, respectively.
The nodes of type~$\C$ come in pairs $\{x,\phi(x)\}$;
we denote the number of such pairs by $\de_{\C}(C)$.
If $u$ is an immersion with transverse intersection points, 
i.e.~$x\!\in\!C$ is a simple node whenever $|u^{-1}(x)|\!\neq\!1$, then
\BE{nodecond_e}\de_E(C)+\de_H(C)+2\de_{\C}(C)=\frac12 
\big(B^2\!-\!\lr{c_1(X,\om),B}\!+\!2\big).\EE

\vspace{.1in}

Suppose in addition that
$B\!\in\!H_2^{\phi}(X;\Z)$ and $l\!\in\!\Z^{\ge0}$ 
are such~that
\BE{dimcond_e} k\equiv \ell_{\om}(B)\!-\!2l\in\Z^{\ge0}\,.\EE
For a generic $J\!\in\!\cJ_{\om}^{\phi}$, 
the set $\fM_{\p}(B;J)$ of rational irreducible real $J$-holomorphic 
degree~$B$ curves $C\!\subset\!X$ passing through a tuple
\BE{bfpdfn_e}\p\!\equiv\!\big((p_1^{\R},\ldots,p_k^{\R}),(p_1^+,\ldots,p_l^+)\big)\in
Y^k\!\times\!(X\!-\!X^{\phi})^l\EE
of $k$ points in~$Y$ and
$l$ points in $X\!-\!X^{\phi}$ in general position is then finite.
Furthermore, every such curve is the image of an immersion~$u$ as in~\eref{Cudfn_e}
with transverse intersection points.
By~\eref{dimcond_e}, either $k\!>\!0$ or $\lr{c_1(X,\om),B}$ is odd.
This implies that for every \hbox{$C\!\in\!\fM_{\p}(B;J)$} there exists 
a  \hbox{$J$-holomorphic} immersion~$u$ as in~\eref{Cudfn_e} satisfying the first condition
in~\eref{Cudfn_e4}.
If $k\!>\!0$,
the set $\fM_{\p}(B;J)$ thus decomposes into the subsets $\fM_{\p}(\wt{B};J)$ of  
real rational irreducible $J$-holomorphic degree~$\wt{B}$ curves with
$$\wt{B}\!=\!(B,b) \in \wt{H}_2^{\phi}(X,Y)\subset 
H_2^{\phi}(X;\Z)\!\oplus\!H_1(Y;\Z_2);$$
if $k\!=\!0$, the above decomposition needs to be taken over all connected components
of~$X^{\phi}$. 
According to \cite[Theorem~0.1]{Wel4}, the~sum
\BE{Welcountdfn_e}N_{\wt{B},l}^{\phi}(Y)\equiv 
\sum_{C\in\fM_{\p}(\wt{B};J)}\!\!\!\!\!\!\!\!\!{(-1)}^{\de_E(C)}\EE
is independent of generic choices of $J\!\in\!\cJ_{\om}^{\phi}$ and 
$\p\!\in\!Y^k\!\times\!X^l$.\\

As noted in Section~\ref{g0maps_subs} and justified in Section~\ref{MapSignDfn_subs}, 
a relative $\Pin^-$-structure~$\fp$ 
on $Y\!\subset\!X$ satisfying~\ref{Rcond_it1} in Section~\ref{g0maps_subs}
determines a sign~$\fs_{\fp}(C)\!\in\!\{\pm1\}$
for every $C\!\in\!\fM_{\p}(\wt{B};J)$.
Let
\BE{SolInvdfn_e}N_{\wt{B},l}^{\phi,\fp}(Y)\equiv
\blr{\underset{l}{\underbrace{\pt,\ldots,\pt}}}_{\wt{B};Y;k}^{\phi,\fp} 
\equiv \sum_{C\in\fM_{\p}(\wt{B};J)}\!\!\!\!\!\!\!\!\!\fs_{\fp}(C)\EE
be the associated invariant curve count~\eref{RGWdfn_e} with point insertions only.
If $\fp$ is a $\Pin^-$-structure on~$Y$, let 
\BE{mapcount2dfn_e}N_{\wt{B},l}^{\phi,\fp}(Y)\equiv N_{\wt{B},l}^{\phi,\io_X(\fp)}(Y)\,.\EE

\begin{thm}\label{WelSolComp_thm}
Let  $(X,\om,\phi)$ be a compact real symplectic fourfold 
and $Y\!\subset\!X^{\phi}$ be a connected component.
There exists a collection 
$$\big(\ws_{\fp}\!:\wt{H}_2^{\phi}(X,Y)\!\lra\!\Z_2\big)_{\fp\in\cP^-(Y)}$$
of group homomorphisms so that 
\begin{gather}
\label{WelSolComp_e0a}
\ws_{\fp}\big(\fd_Y(B')\!\big)=\fd_Y(B')^2/2\!+\!2\Z, \qquad
\ws_{\eta\cdot\fp}(B,b)=\ws_{\fp}(B,b)\!+\!\lr{\eta,b},\\
\label{WelSolComp_e0b} 
\fs_{\fp}(C)=(-1)^{\!\ws_{\fp}(B,b)+l-1+\de_E(C)}
\end{gather}
for all $\fp\!\in\!\cP^-(Y)$, $\eta\!\in\!H^1(Y;\Z_2)$, $B'\!\in\!H_2(X;\Z)$,
$(B,b)\!\in\!\wt{H}_2^{\phi}(X,Y)$ with $B\!\neq\!0$, 
$k,l\!\in\!\Z^{\ge0}$ satisfying~\eref{dimcond_e},
$\p\!\in\!Y^k\!\times\!X^l$ and $J\!\in\!\cJ_{\om}^{\phi}$ generic, 
and $C\!\in\!\fM_{\p}(\!(B,b);J)$.
\end{thm}

Let $\mu_\fp$ be as in Theorem~\ref{ShExSeqImm_thm}.
We show in Section~\ref{WelSolPf_thm} that~\eref{WelSolComp_e0b} holds with
\BE{wsfpdfn_e}\ws_{\fp}\!: \wt{H}_2^{\phi}(X,Y)\lra \Z_2, \quad
\ws_{\fp}(B,b)=\frac{\lr{c_1(X,\om),B}^2\!+\!B^2}{2}+\mu_{\fp}(b)\,.\EE
By~\eref{ShExSeqImm_e1}, this map is a group homomorphism 
which satisfies the first property in~\eref{WelSolComp_e0a}.
By~\eref{ShExSeqImm_e1b}, it also satisfies the second property in~\eref{WelSolComp_e0a}.\\

By~\eref{Welcountdfn_e}, \eref{SolInvdfn_e}, and \eref{WelSolComp_e0b},
\BE{WelSolComp_e0c}
N_{\wt{B},l}^{\phi,\fp}(Y)=(-1)^{\!\ws_{\fp}(\wt{B})+l-1}N_{\wt{B},l}^{\phi}(Y)\,.\EE
Thus, Theorem~\ref{WelSolComp_thm} establishes a comparison between the counts
of real genus~0 $J$-holomorphic curves in~$(X,\om,\phi)$
with the intrinsic signs of~\cite{Wel4} and 
the counts of such curves with signs dependent on a $\Pin^-$-structure~$\fp$ on~$Y$
as in~\cite{Sol,Ge2,RealWDVV}.\\

By the RelSpinPin~\ref{RelSpinPinStr_prop} property, a relative $\Pin^-$-structure~$\fp$ 
on $Y\!\subset\!X$ satisfying~\ref{Rcond_it1} 
is of the~form 
$$\fp=\eta\!\cdot\!\io_X(\fp') \quad\hbox{with}~~
\eta\!\in\!H^2(X,Y;\Z_2),\,\fp'\!\in\!\cP^-(Y)
~\hbox{s.t.}~\lr{\eta|_X,B'}=0~\forall\,
(B',b')\!\in\!\wt{H}_2^{\phi}(X,Y).$$
If $C$ is an element of $\fM_{\p}(\wt{B};J)$, 
$u$ is as in~\eref{Cudfn_e} and~\eref{Cudfn_e4}, 
and $\bD^2_+\!\subset\!\P^1$ is a disk cut out by~$S^1$, then
\BE{WelSolComp_e0e}\fs_{\fp}(C)=(-1)^{\lr{\eta,u_*[\bD^2_+]_{\Z_2}}}\fs_{\fp'}(C)\,.\EE
This identity is a direct consequence of the CROrient~\ref{CROSpinPinStr_prop}\ref{CROPinStr_it}
property on page~\pageref{CROSpinPinStr_prop}. 
By the condition on~$\eta$ above, the exponent in~\eref{WelSolComp_e0e}
does not depend on the choice of~$\bD^2_+$.
The identity~\eref{WelSolComp_e0e} translates into a signed relation 
between the  counts 
$N_{\wt{B},l}^{\phi,\fp}(Y)$ and $N_{\wt{B},l}^{\phi,\fp'}(Y)$
of rational curves with a refined notion of the degree~$\wt{B}$ incorporating an element
of the quotient
$$H_2(X,Y;\Z)\big/\big\{B'\!+\!\phi_*B'\!: B'\!\in\!H_2(X,Y;\Z)\big\}$$
to encode the relative homology class of $u_*[\bD^2_+]_{\Z}$.\\

Suppose next that $(X,\om,\phi)$ is a compact connected real symplectic sixfold, 
$Y\!\subset\!X^{\phi}$ is a connected component, and $\os$ is an $\OSpin$-structure
on~$Y$. 
Let $l\!\in\!\Z^{\ge0}$ and $\wt{B}\!\in\!(B,b)\!\in\!\wt{H}_2^{\phi}(X,Y)$ be so~that 
\BE{dimcond3_e} k\equiv \ell_{\om}(B)/2\!-\!2l\in\Z^{\ge0}\EE
and the conditions in~\ref{Rcond_it3} in Section~\ref{g0maps_subs} are satisfied.
For a generic $J\!\in\!\cJ_{\om}^{\phi}$, 
the set $\fM_{\p}(\wt{B};J)$ of rational irreducible real $J$-holomorphic 
degree~$\wt{B}$ curves $C\!\subset\!X$ passing through a generic tuple~$\p$
of points as in~\eref{bfpdfn_e} is then finite.
Each element $C\!\in\!\fM_{\p}(\wt{B};J)$ is an embedded curve and 
is assigned an intrinsic sign $\fs_{\os}^*(C)\!\in\!\{\pm1\}$ in~\cite{Wel6,Wel6b};
see the proof of Theorem~\ref{WelSolComp3_thm} in Section~\ref{WelSolPf_thm}.
According to \cite[Theorem~0.1]{Wel6} in the projective case and 
\cite[Theorem~4.1]{Wel6b} in general, the~sum
\BE{Welcountdfn3_e}
W_{\wt{B},l}^{\phi,\os}(Y)\equiv 
\sum_{C\in\fM_{\p}(\wt{B};J)}\!\!\!\!\!\!\!\!\!\fs_{\os}^*(C)\EE
is independent of generic choices of $J\!\in\!\cJ_{\om}^{\phi}$ and 
$\p\!\in\!Y^k\!\times\!X^l$.
By the proof of Proposition~1.3 in~\cite{RealWDVV3} or Theorem~\ref{WelSolComp3_thm} below,
this is also the case if  the conditions in~\ref{Rcond_it4}, instead of~\ref{Rcond_it3},
are satisfied.\\

As noted in Section~\ref{g0maps_subs} and justified in Section~\ref{MapSignDfn_subs}, 
an $\OSpin$-structure~$\os$ on~$Y$ determines a sign~$\fs_{\os}(C)\!\in\!\{\pm1\}$
for every $C\!\in\!\fM_{\p}(\wt{B};J)$.
Let
$$N_{\wt{B},l}^{\os,\phi}(Y)\equiv
\blr{\underset{l}{\underbrace{\pt,\ldots,\pt}}}_{\wt{B};Y;k}^{\phi,\io_Y(\os)} 
\equiv \sum_{C\in\fM_{\p}(\wt{B};J)}\!\!\!\!\!\!\!\!\!\fs_{\os}(C)$$
be the associated invariant curve count~\eref{RGWdfn_e} with point insertions only.
The next theorem relates this count to the one in~\eref{Welcountdfn3_e}.

\begin{thm}\label{WelSolComp3_thm}
Let $(X,\om,\phi)$ be a compact real symplectic sixfold,
$Y\!\subset\!X^{\phi}$ be a connected component, 
and $\os$ be an $\OSpin$-structure on~$Y$. 
Then,
\BE{WelSolComp3_e0b} 
\fs_{\os}(C)=(-1)^{\binom{k+1}{2}+1}\fs_{\os}^*(C)\EE
for all $k,l\!\in\!\Z^{\ge0}$ and $\wt{B}\!\in\!\wt{H}_2^{\phi}(X,Y)$
satisfying~\eref{dimcond3_e} and~\ref{Rcond_it3} in Section~\ref{g0maps_subs},
$\p\!\in\!Y^k\!\times\!X^l$ and $J\!\in\!\cJ_{\om}^{\phi}$ generic,
and $C\!\in\!\fM_{\p}(\wt{B};J)$.
\end{thm}

\subsection{Basic examples: fourfolds}
\label{Wel4eg_subs}

For a compact real symplectic manifold $(X,\om,\phi)$ 
with connected real locus~$X^{\phi}$, we define
$$\wt{H}_2^{\phi}(X)\equiv \wt{H}_2^{\phi}\big(X,X^{\phi}\big)\,.$$
We also omit $Y\!=\!X^{\phi}$ from the notation in~\eref{Welcountdfn_e}, 
\eref{mapcount2dfn_e}, and~\eref{Welcountdfn3_e}.\\

Suppose in addition that the (real) dimension of~$X$ is~4,
$\fp$~is a $\Pin^-$-structure on~$X^{\phi}$,
$\wt{B}\!\in\!\wt{H}_2^{\phi}(X)$ is as in~\eref{wtBdfn_e},
and $k,l\!\in\!\Z^{\ge0}$ satisfy~\eref{dimcond_e}.
If $(X,\om,J)$ is projective and Fano, 
then the number~$N_B$ of (complex) rational 
irreducible $J$-holomorphic degree~$B$ curves passing through 
$k\!+\!2l$ generic points in~$X$ is finite and independent of the choice of
the points.
Furthermore,
\BE{RvsCGWs_e}
\big|N_{\wt{B},l}^{\phi}\big|,
\big|N_{\wt{B},l}^{\phi,\fp}\big|\le N_B \qquad\hbox{and}\qquad
N_{\wt{B},l}^{\phi},N_{\wt{B},l}^{\phi,\fp}\cong N_B\mod2.\EE

\begin{eg}\label{RP2pin_eg5}
The complex projective plane $\P^2$ with the Fubini-Study symplectic form and
the standard conjugation
$$\tau_2\!:\P^2\lra \P^2, \qquad 
\tau_2\big([Z_0,Z_1,Z_2]\big)=\big[\ov{Z_0},\ov{Z_1},\ov{Z_2}\big],$$ 
is a compact real symplectic fourfold.
The fixed locus of~$\tau_2$ is the real projective plane~$\R\P^2$.
The group $\wt{H}_2^{\tau_2}(\P^2)$ is freely generated by the pair $(L,L^{\R})$
consisting of the homology class~$L$ 
of a linearly embedded $\P^1\!\subset\!\P^2$ and the nonzero element~$L^{\R}$
of~$H_1(\R\P^2;\Z_2)$.
We thus identify $\wt{H}_2^{\tau_2}(\P^2)$ with~$\Z$.
There is a unique line (degree~1 curve) passing through 2~general points in~$\P^2$;
there is also a unique conic (degree~2 curve) passing through 5~general points in~$\P^2$.
Since both of these curves are rational and smooth (without nodes),
\BE{RP2pin5_e3} N_1,N_2=1, \qquad 
N_{1,0}^{\tau_2},N_{1,1}^{\tau_2}=1, \qquad\hbox{and}\qquad 
N_{2,0}^{\tau_2},N_{2,1}^{\tau_2},N_{2,2}^{\tau_2}=1\,.\EE
Let $\fp_0^-\!\equiv\!\fp_0^-(\R\P^2)$ and $\fp_1^-\!\equiv\!\fp_1^-(\R\P^2)$ 
be the two $\Pin^-$-structures on~$\R\P^2$ as in Example~\ref{RP2pin_eg}.
By~\eref{wsfpdfn_e} and~\eref{RP2pin4_e3},
$$\ws_{\fp_0^-}(d)=d,\quad \ws_{\fp_1^-}(d)=0 \qquad\forall\,d\!\in\!\Z.$$
Along with~\eref{WelSolComp_e0c} and~\eref{RP2pin5_e3}, this implies that 
\begin{alignat*}{5}
N_{1,0}^{\tau_2,\fp_0^-}&=1,&\quad N_{1,1}^{\tau_2,\fp_0^-}&=-1, &\quad
N_{2,0}^{\tau_2,\fp_0^-}&=-1,&\quad N_{2,1}^{\tau_2,\fp_0^-}&=1, &\quad
N_{2,2}^{\tau_2,\fp_0^-}&=-1;\\
N_{1,0}^{\tau_2,\fp_1^-}&=-1,&\quad N_{1,1}^{\tau_2,\fp_1^-}&=1, &\quad
N_{2,0}^{\tau_2,\fp_1^-}&=-1,&\quad N_{2,1}^{\tau_2,\fp_1^-}&=1, &\quad
N_{2,2}^{\tau_2,\fp_1^-}&=-1.
\end{alignat*}
\end{eg}

\vspace{.2in}

Let $\om_{\FS}$ be the Fubini-Study symplectic form on~$\P^1$ and 
$$\om=\pi_1^*\om_{\FS}\!+\!\pi_2^*\om_{\FS}\,,$$
where $\pi_1,\pi_2\!: \P^1\!\times\!\P^1\lra \P^1$
are the two component projections.
The group $H_2(\P^1\!\times\!\P^1;\Z)$ is freely generated by 
the homology classes~$L_1$ and~$L_2$ of the slices 
$\P^1\!\times\!q_2$ and $q_1\!\times\!\P^1$, respectively,
with $q_1,q_2\!\in\!\P^1$.
We thus identify this homology group with~$\Z^2$.
By symmetry, 
$$N_{(a,b)}^{\P^1\!\times\P^1}=N_{(b,a)}^{\P^1\!\times\P^1} 
\qquad\forall\,a,b\!\in\!\Z^{\ge0}\,.$$
Since a degree $(a,0)$-curve in $\P^1\!\times\!\P^1$ is a degree~$a$ cover of 
a horizontal section,
\BE{P1P1wel_e3c} N_{(a,0)}^{\P^1\!\times\P^1}=0 \quad\forall\,a\!\ge\!2.\EE
An irreducible degree $(1,b)$-curve in $\P^1\!\times\!\P^1$ 
corresponds to the graph of a ratio of two degree~$b$ polynomials on~$\C$.
Since every such rational function is determined by its values at $2b\!+\!1$ points,
\BE{P1P1wel_e3}N_{(1,b)}^{\P^1\!\times\P^1}=1 \qquad\forall\,b\!\in\!\Z^{\ge0}.\EE

\vspace{.15in}

The involutions 
$$\tau_{1,1},\tau_{1,1}'\!:\P^1\!\times\!\P^1\lra\P^1\!\times\!\P^1, \quad
\tau_{1,1}(q_1,q_2)=\big(\tau(q_1),\tau(q_2)\!\big),~~
\tau_{1,1}'(q_1,q_2)=\big(\tau(q_2),\tau(q_1)\!\big),$$
on $\P^1\!\times\!\P^1$ are anti-symplectic with respect to~$\om$.
Their fixed loci are~$\R\P^1\!\times\!\R\P^1$ and
$$\gr(\tau)\equiv\big\{(q,\tau(q)\!)\!:q\!\in\!\P^1\big\}\approx S^2,$$ 
respectively.

\begin{eg}\label{P1P1pin_eg5}
The group $\wt{H}_2^{\tau_{1,1}}\!(\P^1\!\times\!\P^1)$ is freely generated by 
the pairs $(L_1,L_1^{\R})$ and $(L_2,L_2^{\R})$, where
$$L_1^{\R},L_2^{\R}\in H_1\big(\R\P^1\!\times\!\R\P^1;\Z_2\big)$$
are the homology classes of the slices $\R\P^1\!\times\!q_2^{\R}$  
and $q_1^{\R}\!\times\!\R\P^1$, respectively, with \hbox{$q_1^{\R},q_2^{\R}\!\in\!\R\P^1$}.
We thus identify $\wt{H}_2^{\tau_{1,1}}\!(\P^1\!\times\!\P^1)$ with~$\Z^2$.
By symmetry, 
\BE{P1P1wel_e4a} N_{(a,b),l}^{\tau_{1,1}}=N_{(b,a),l}^{\tau_{1,1}} 
\qquad\forall\,a,b,l\!\in\!\Z^{\ge0}\,.\EE
By~\eref{P1P1wel_e3c} and~\eref{P1P1wel_e3}, 
\BE{P1P1wel_e4} 
N_{(a,0),l}^{\tau_{1,1}}=0 \quad\forall\,a\!\ge\!2,\,l\!\in\!\Z^{\ge0}
\quad\hbox{and}\quad
N_{(1,b),0}^{\tau_{1,1}},N_{(1,b),1}^{\tau_{1,1}},\ldots,N_{(1,b),b}^{\tau_{1,1}}=1 
~~\forall\,b\!\in\!\Z^{\ge0},\EE
respectively.
Let $\fp_0^-$ be the $\Pin^-$-structure on $(\R\P^1)^2\!=\!(S^1)^2$ 
corresponding to the canonical $\OSpin$-structure~$\os_0$ of Example~\ref{ShExSeqImm_eg0a}
via the bijection~\eref{Pin2SpinRed_e}.
By~\eref{wsfpdfn_e} and~\eref{ShExSeqImm0a_e3},
$$\ws_{\fp_0^-}(a,b)=a\!+\!b \qquad\forall\,a,b\!\in\!\Z.$$
Along with~\eref{WelSolComp_e0c}, \eref{P1P1wel_e4a}, and~\eref{P1P1wel_e4}, 
this implies~that 
\begin{gather}\notag
N_{(a,b),l}^{\tau_{1,1},\fp_0^-}=N_{(b,a),l}^{\tau_{1,1},\fp_0^-}~\forall\,a,b,l\!\in\!\Z^{\ge0}, 
\quad
N_{(a,0),l}^{\tau_{1,1},\fp_0^-}=0 \quad\forall\,a\!\ge\!2,\,l\!\in\!\Z^{\ge0},\\
\label{P1P1pin5_e7a}
N_{(1,b),l}^{\tau_{1,1},\fp_0^-}=(-1)^{b+l}~\forall\,b,l\!\in\!\Z^{\ge0},\,l\!\le\!b.
\end{gather}
\end{eg}

\begin{eg}\label{P1P1pin2_eg5} 
The group $\wt{H}_2^{\tau_{1,1}'}\!(\P^1\!\times\!\P^1)$ is freely generated by 
$(L_1\!+\!L_2,0)$.
We thus identify $\wt{H}_2^{\tau_{1,1}'}\!(\P^1\!\times\!\P^1)$ with~$\Z$.
By~\eref{P1P1wel_e3}, 
\BE{P1P1wel2_e4}  N_{1,0}^{\tau_{1,1}'},N_{1,1}^{\tau_{1,1}'}=1.\EE
Let $\fp_0^-$ be the unique $\Pin^-$-structure on~$S^2$. 
By~\eref{wsfpdfn_e} and the conclusion of Example~\ref{ShExSeqImm_eg0},
$$\ws_{\fp_0^-}(d)=d \qquad\forall\,d\!\in\!\Z.$$
Along with~\eref{WelSolComp_e0c} and~\eref{P1P1wel2_e4}, 
this implies~that 
$$N_{1,0}^{\tau_{1,1}';\fp_0^-}=1, \qquad N_{1,1}^{\tau_{1,1}';\fp_0^-}=-1.$$
\end{eg}

\subsection{Basic examples: sixfolds}
\label{Wel6eg_subs}

We next obtain counts of real rational curves in $\P^3$ with the standard conjugation
and in $(\P^1)^3$ with two different conjugations.
We identify an element of~$H_*(X;\Z)$ \hbox{(resp.~$H_*(X\!-\!Y;\Z)$)}
and its Poincare dual in~$H^*(X;\Z)$ (resp.~$H^*(X,Y;\Z)$).

\subsubsection{The projective space $\P^3$.}\label{RP3_ssubs}
The complex projective space $\P^3$ with the Fubini-Study symplectic form and
the standard conjugation
$$\tau_3\!:\P^3\lra \P^3, \qquad 
\tau_3\big([Z_0,Z_1,Z_2,Z_3]\big)=\big[\ov{Z_0},\ov{Z_1},\ov{Z_2},\ov{Z_3}\big],$$ 
is a compact real symplectic sixfold.
The fixed locus of~$\tau_3$ is the real projective plane~$\R\P^3$.
We denote by \hbox{$L\!\in\!H_2(\P^3;\Z)$} the homology class 
of a linearly embedded $\P^1\!\subset\!\P^3$  
and by \hbox{$L^{\R}\!\in\!H_1(\R\P^3;\Z_2)$} the nonzero element.
The groups~$H_2(\P^3;\Z)$ and~$\wt{H}_2^{\tau_3}(\P^3)$ are freely generated by~$L$
and~$(L,L^{\R})$, respectively; we thus identify them with~$\Z$.\\

The involution
$$\psi\!:\P^3\lra\P^3, \quad
\psi\big([Z_0,Z_1,Z_2,Z_3]\big)=[Z_0,Z_1,Z_2,-Z_3],$$
satisfies the conditions in~\ref{Rcond_it4} in Section~\ref{RealSetting_subs}
with $Y\!=\!\R\P^3$; the same is the case of any automorphisms of~$\P^3$
exchanging a pair of its homogeneous coordinates.
By the $n\!=\!3$ case of~\eref{RPnEES_e3} and~\eref{SpinPinStr_e}, 
\BE{RP3spin5_e2} \psi^*\os \neq \ov\os \qquad\forall~\os\!\in\!\OSp(\R\P^3),\EE 
i.e.~$\psi$ reverses both the orientation and the Spin component of 
every $\OSpin$-structure on~$\R\P^3$.\\

As explained in Section~2.6 of~\cite{RealWDVV3}, the lines $\P^1\!\subset\!\P^3$ 
contained in~$\P^3\!-\!\R\P^3$ determine two distinct homology classes $L_+^c,L_-^c$ 
in~$\P^3\!-\!\R\P^3$ so~that
\BE{RP3spin5_e3} L_+^c\!-\!L_-^c=\big[S(\cN_y\R\P^3,\os)\big]\EE 
for an $\OSpin$-structure~$\os$ on~$\R\P^3$.
If $H_1,H_2\!\subset\!\P^3$ are two projective hyperplanes intersecting transversally in
$\P^3\!-\!\R\P^3$ so that $H_1\!\cap\!H_2$ represents~$L_+^c$, 
then $H_1\!\cap\!\tau_3(H_2)$ represents~$L_-^c$. 
By~\eref{RP3spin5_e3}, $L_+^c$ and $L_-^c$ freely generate $H_2(\P^3\!-\!\R\P^3;\Z)$ and
are interchanged by~$\psi$.\\

For $d,a,b\!\in\!\Z^{\ge0}$ with $a\!+\!2b\!=\!4d$,
let $N_{d;a,b}\!\in\!\Z^{\ge0}$ be the number of (complex) 
rational irreducible holomorphic degree~$d$ curves passing through 
$a$~generic lines and~$b$ generic points in~$\P^3$.
For $d,a,b\!\in\!\Z^{\ge0}$ and an $\OSpin$-structure~$\os$ on~$\R\P^3$, let
$$N^{\tau_3,\os}_{d;a,b}=
\blr{\underset{b}{\underbrace{{\pt,\ldots,\pt}}},
\underset{a}{\underbrace{{(L_+^c\!+\!L_-^c)/2,\ldots,(L_+^c\!+\!L_-^c)/2}}}
}_{d;\R\P^3;2d-a-2b}^{\tau_3,\os}\,.$$
By Theorem~\ref{RGWs_thm}\ref{sphere_it} and~\eref{RP3spin5_e3}, 
$N^{\tau_3,\os}_{d;a,b}\!\in\!\Z$.
Similarly to~\eref{RvsCGWs_e},
$$N_{d;0,2d}\!-\!\big|W^{\tau_3,\os}_{d;b}\big|
\in2\Z^{\ge0}~~\hbox{if}~d\!\ge\!b, \qquad
N_{d;2a,2d-a}\!-\!\big|N^{\tau_3,\os}_{d;a,b}\big|
\in2\Z^{\ge0} ~~\hbox{if}~2d\!\ge\!a\!+\!2b.$$
By~\eref{RP3spin5_e2} and 
the proofs of Theorems~\ref{RGWs_thm}\ref{OrientRev_it} and~\ref{WelSolComp3_thm},
\BE{RP3spin5_e7}
N^{\tau_3,\os}_{d;a,b}=0~\hbox{if}~d\!-\!a\!\in\!2\Z,~~
W^{\tau_3,\os}_{d;b}=0~\hbox{if}~d\!\in\!2\Z, \quad
N^{\tau_3,\psi^*\os}_{d;a,b}=N^{\tau_3,\os}_{d;a,b}, ~~
W^{\tau_3,\psi^*\os}_{d;b}=W^{\tau_3,\os}_{d;b}\EE
for every automorphism~$\psi$ of $(\P^3,\tau_3)$.
Since there is a unique line passing through 2~general points in~$\P^3$,
$N_{1;0,1}\!=\!1$.
Along with the definition of the sign~$\fs_{\os}^*(C)$ in~\cite{Wel6,Wel6b}
and the $k\!=\!2,0$ cases of~\eref{WelSolComp3_e0b}, this implies~that 
$$W^{\tau_3,\os}_{1;0}=W^{\tau_3,\os}_{1;1}=\pm1, \qquad
-N^{\tau_3,\os}_{1;0,0}=N^{\tau_3,\os}_{1;0,1}=-W^{\tau_3,\os}_{1;1}=\mp1.$$
The first equality in the second equation above is also implied by 
the extended WDVV equations for $(\P^3,\om,\tau_3)$
predicted in~\cite{Adam} and established in~\cite{RealWDVV3}.\\

Let $\os_0\!\equiv\!\os_0(\R\P^3)$ and $\os_1\!\equiv\!\os_1^-(\R\P^3)$ 
be the two $\OSpin$-structures on~$\R\P^3$ defined in Example~\ref{RP3spin_eg}.
By the third statement in~\eref{RP3spin5_e7}, the numbers 
$N^{\tau_3,\os_0}_{d;a,b}$ and $N^{\tau_3,\os_1}_{d;a,b}$ are invariant
under linear reparametrizations of~$(\P^3,\tau_3)$;
this implication is a special case of Theorem~1.6 in~\cite{RealGWsIII}.
It also follows from
the  extended WDVV equations for $(\P^3,\om,\tau_3)$, which 
determine all numbers $N^{\tau_3,\os}_{d;a,b}$ from the single input~$N^{\tau_3,\os}_{1;0,1}$;
see Section~1.4 in~\cite{Adam} and Section~7 in~\cite{RealWDVVapp}.\\

The holomorphic normal bundle~$\cN_{\P^3}\P^1$ of any real holomorphic line 
$\P^1\!\subset\!\P^3$ is isomorphic to two copies of a rank~1 real bundle pair~$(V,\vph)$ 
over~$(S^2,\tau)$ of degree~1. 
In particular, $\cN_{\P^3}\P^1$ is \sf{balanced} in the sense of
\cite[Section~2.2]{Wel6} and \cite[Section~1.2.2]{Wel6b}.
As noted in Remark~\ref{SpinDfn1to3_rmk}, 
the base OSpin-structure on~$T(\R\P^3)|_{\R\P^1}$ in \cite{Wel6,Wel6b} is
the OSpin-structure $\os_0|_{\R\P^1}$.
Along with the definition of the sign~$\fs_{\os}^*(C)$ in \cite[Section~2.2]{Wel6}
and \cite[Section~3.2]{Wel6b}, this gives
$$W^{\tau_3,\os_0}_{1;0}=W^{\tau_3,\os_0}_{1;1}=1 \qquad\hbox{and}\qquad
W^{\tau_3,\os_1}_{1;0}=W^{\tau_3,\os_1}_{1;1}=-1.$$
Along with the $k\!=\!2,0$ cases of~\eref{WelSolComp3_e0b}, this implies~that 
\BE{RP3spin5_e11}-N^{\tau_3,\os_0}_{1;0,0}=N^{\tau_3,\os_0}_{1;0,1}=-1
\qquad\hbox{and}\qquad
-N^{\tau_3,\os_1}_{1;0,0}=N^{\tau_3,\os_1}_{1;0,1}=1.\EE
By the paragraph above Theorem~1.5 in~\cite{RealGWsII},
the value of~$+1$ for $N^{\tau_3,\os_0}_{1;0,1}$ obtained in Example~6.3 in~\cite{Teh}
corresponds to the orientation of~\eref{OrientSubs_e3} with the orientation 
on the last factor on the right-hand side as in~\cite{RealEnum,RealGWsII}.
As noted in Section~\ref{DM_subs}, the last orientation is the opposite
of the orientation~$\fo_{0,l}$ used in the present manuscript.
Thus, the last equality in~\eref{RP3spin5_e11} agrees with 
the $m,l\!=\!1$, $t_1\!=\!3$ statement of Example~6.3 in~\cite{Teh}.

\subsubsection{The sixfold $(\P^1)^3$.}
\label{P1P1P1_ssubs}

We now take $X\!=\!(\P^1)^3$ with the symplectic form
$$\om=\pi_1^*\om_{\FS}\!+\!\pi_2^*\om_{\FS}\!+\!\pi_3^*\om_{\FS}\,,$$
where $\pi_1,\pi_2,\pi_3\!:(\P^1)^3\!\lra\!\P^1$ are the component projections.
The group $H_2(\!(\P^1)^3;\Z)$ is freely generated by 
the homology classes~$L_1,L_2,L_3$ of the~slices 
\BE{P1P1P1lndfn_e}\P^1\!\times\!q_2\times\!q_3, \quad q_1\times\!\P^1\!\times\!q_3,
\quad\hbox{and}\quad q_1\!\times\!q_2\!\times\!\P^1,\EE
respectively, with $q_1,q_2,q_3\!\in\!\P^1$.\\

For $\bfd\!\equiv\!(d_1,d_2,d_3)$ with each $d_i\!\in\!\Z$, let
$$|\bfd|=d_1\!+\!d_2\!+\!d_3, \qquad
B(\bfd)=d_1L_1\!+\!d_2L_2\!+\!d_3L_3\in 
H_2\big(\!(\P^1)^3;\Z\big).$$
If in addition $\bfa\!\equiv\!(a_1,a_2,a_3)$ with each $a_i\!\in\!\Z^{\ge0}$
and $b\!\in\!\Z^{\ge0}$ are so that 
\BE{P1Cdimcond_e}(a_1\!+\!a_2\!+\!a_3)\!+\!2b=2(d_1\!+\!d_2\!+\!d_3),\EE
we define $N_{\bfd;\bfa,b}\!\in\!\Z^{\ge0}$ to be the number of (complex) 
rational irreducible holomorphic degree~$B(\bfd)$ curves passing through 
$a_1,a_2,a_3$ generic slices~\eref{P1P1P1lndfn_e} representing~$L_1,L_2,L_3$,
respectively, and $b$ generic points in~$(\P^1)^3$.\\

The involutions $\phi\!\equiv\!\tau\!\times\!\tau_{1,1}$ and
$\phi'\!\equiv\!\tau\!\times\!\tau_{1,1}'$
on $(\P^1)^3$ are anti-symplectic with respect to~$\om$.
Their fixed loci are~$(\R\P^1)^3$ and $\R\P^1\!\times\!\gr(\tau)$,
respectively.
In both cases, we take~$Y$ to be the fixed locus.
The involution
$$\psi\!:(\P^1)^3\lra(\P^1)^3, \quad
\psi\big([Z_1,W_1],[Z_2,W_2],[Z_3,W_3]\big)=
\big([Z_1,-W_1],[Z_2,W_2],[Z_3,W_3]\big),$$
satisfies the conditions in~\ref{Rcond_it4} in Section~\ref{RealSetting_subs}.\\

We take the $\OSpin$-structures on the two fixed loci to~be 
$$\os\equiv\bllrr{\pi_1^*\os_0(TS^1),
\llrr{\pi_2^*\os_0(TS^1),\pi_3^*\os_0(TS^1)}_{\oplus}}_{\!\oplus} ~~\hbox{and}~~
\os'\equiv\bllrr{\pi_1^*\os_0(TS^1),\pi_{23}^{\,*}\os_0(TS^2)}_{\!\oplus}\,,$$
where $\os_0(TS^1)$ is as in Theorem~\ref{RGWs_thm2}\ref{RGWprod_it1}, 
$\os_0(TS^2)$ is the unique $\OSpin$-structure for either orientation on~$S^2$, and
$$\pi_1,\pi_2,\pi_3\!:(\R\P^1)^3\lra\R\P^1\!\approx\!S^1 \qquad\hbox{and}\qquad
\pi_1,\pi_{23}\!:\R\P^1\!\times\!\gr(\tau)\lra \R\P^1,\gr(\tau)\!\approx\!S^2$$
are the component projections.
Since $\psi^*\os\!=\!\ov\os$ and $\psi^*\os'\!=\!\ov\os'$,
Theorem~\ref{RGWs_thm}\ref{OrientRev_it} implies~that 
\BE{P1P1P1van_e}\begin{split}
&\blr{\mu_1,\ldots,\mu_l}_{\wt{B};(\R\P^1)^3;k}^{\phi,\os}=0
\qquad\hbox{and}\qquad
 \blr{\mu_1,\ldots,\mu_l}_{\wt{B};\R\P^1\!\times\!\gr(\tau);k}^{\phi',\os'}=0\\
&\hspace{.5in}\forall~k\!\in\!2\Z,~
\mu_1,\ldots,\mu_l\!\in\!\wh{H}^{2*}\!\big(\!(\P^1)^3,Y\big)~\hbox{s.t.}~
\psi^*\mu_i=\mu_i~\forall\,i.
\end{split}\EE

\begin{eg}\label{P1P1P1_eg1}
The group $\wt{H}_2^{\phi}(\!(\P^1)^3)$ is freely generated by 
the~pairs 
$$\wt{L}_1\equiv\big(L_1,L_1^{\R}\big), \qquad  
\wt{L}_2\equiv\big(L_2,L_2^{\R}\big), 
\qquad\hbox{and}\qquad \wt{L}_3\equiv\big(L_3,L_3^{\R}\big),$$ 
where $L_1^{\R},L_2^{\R},L_3^{\R}\!\in\!H_2(\!(\R\P^1)^3;\Z_2)$ 
are the homology classes of the slices 
$$\R\P^1\!\times\!q_2^{\R}\times\!q_3^{\R}, \qquad
q_1^{\R}\times\!\R\P^1\!\times\!q_3^{\R}, \qquad\hbox{and}\qquad
q_1^{\R}\!\times\!q_2^{\R}\!\times\!\R\P^1,$$ 
respectively, with $q_1^{\R},q_2^{\R},q_3^{\R}\!\in\!\R\P^1$.
The space of the slices~\eref{P1P1P1lndfn_e} of each kind contained in $(\P^1)^3\!-\!(\R\P^1)^3$ 
is path-connected and thus determines a homology class~$L_i^c$
in $(\P^1)^3\!-\!(\R\P^1)^3$.
For $\bfd\!\equiv\!(d_1,d_2,d_3)$ and $\bfa\!\equiv\!(a_1,a_2,a_3)$ with 
$d_i,a_i\!\in\!\Z^{\ge0}$ and $b\!\in\!\Z^{\ge0}$, define
\begin{gather*}
\wt{B}(\bfd)=d_1\wt{L}_1\!+\!d_2\wt{L}_2\!+\!d_3\wt{L}_3, \quad
k_{\bfd}(\bfa,b)=|\bfd|\!-\!(a_1\!+\!a_2\!+\!a_3)\!-\!2b, \quad
W^{\phi,\os}_{\bfd;b}=W^{\phi,\os}_{\wt{B}(\bfd),b}\,,\\
N^{\phi,\os}_{\bfd;\bfa,b}=
\blr{\underset{b}{\underbrace{{\pt,\ldots,\pt}}},
\underset{a_1}{\underbrace{{L_1^c,\ldots,L_1^c}}},
\underset{a_2}{\underbrace{{L_2^c,\ldots,L_2^c}}},
\underset{a_3}{\underbrace{{L_3^c,\ldots,L_3^c}}}
}_{\wt{B}(\bfd);(\R\P^1)^3;k_{\bfd}(\bfa,b)}^{\phi,\os}\,.
\end{gather*}
Similarly to~\eref{RvsCGWs_e},
$$N_{\bfd;0,|\bfd|}\!-\!\big|W^{\phi,\os}_{\bfd;b}\big|
\in2\Z^{\ge0}~~\hbox{if}~|\bfd|\ge2b, \quad
N_{\bfd;2\bfa,2b+k_{\bfd}(\bfa,b)}\!-\!\big|N^{\phi,\os}_{\bfd;\bfa,b}\big|
\in2\Z^{\ge0}~~\hbox{if}~k_{\bfd}(\bfa,b)\ge0.$$
If $h$ is the diffeomorphism interchanging a pair of components of $(\P^1)^3$,
$h^*\os\!=\!\ov\os$.
Combining this with Theorem~\ref{RGWs_thm}\ref{OrientRev_it} and \eref{P1P1P1van_e},
we obtain
\BE{P1P1P1_e3a}
N^{\phi,\os}_{(d_1,d_2,d_3);(a_1,a_2,a_3),b}=
N^{\phi,\os}_{(d_2,d_1,d_3);(a_2,a_1,a_3),b}=
N^{\phi,\os}_{(d_1,d_3,d_2);(a_1,a_3,a_2),b}\,.\EE
By Theorem~\ref{RGWs_thm2}\ref{RGWprod_it0},
\BE{P1P1P1_e5a}
N^{\phi,\os}_{(0,d_2,d_3);(a_1,a_2,a_3),b}=
\begin{cases}N^{\tau_{1,1},\fp_0^-}_{(d_2,d_3),a_1},
&\hbox{if}~a_1\!=\!d_2\!+\!d_3\!-\!1,~a_2,a_3,b\!=\!0;\\
0,&\hbox{otherwise}.
\end{cases}\EE
By Theorems~\ref{RGWs_thm2}\ref{RGWprod_it1} and~\ref{RGWs_thm}\ref{div_it},
$$N^{\phi,\os}_{(1,d_2,d_3);(a_1,a_2,a_3),b}=
\begin{cases}0,&\hbox{if}~d_2\!+\!d_3\!>\!0,~
a_1\!>\!d_2\!+\!d_3\!+\!a_2\!+\!a_3\!-\!2;\\
N^{\tau_{1,1},\fp_0^-}_{(d_2,d_3),a_1},
&\hbox{if}~a_1\!=\!d_2\!+\!d_3\!-\!2,~a_2,a_3,b\!=\!0;\\
N^{\tau_{1,1},\fp_0^-}_{(d_2,d_3),a_1+1},
&\hbox{if}~a_1\!=\!d_2\!+\!d_3\!-\!2,~a_2,a_3\!=\!0,~b\!=\!1;\\
d_2N^{\tau_{1,1},\fp_0^-}_{(d_2,d_3),a_1},
&\hbox{if}~a_1\!=\!d_2\!+\!d_3\!-\!1,~a_2,b\!=\!0,~a_3\!=\!1.
\end{cases}$$
The holomorphic normal bundle~$\cN$ of any real slice $L_i\!\subset\!(\P^1)^3$
is isomorphic to two copies of a rank~1 real bundle pair~$(V,\vph)$ 
over~$(S^2,\tau)$ of degree~0 and is thus balanced.
The base OSpin-structure on~$T(\!(\P^1)^3)|_{L_i^{\R}}$ in \cite{Wel6,Wel6b} is
the OSpin-structure $\os|_{L_i^{\R}}$ or its conjugate $\ov\os|_{L_i^{\R}}$.
Along with the definition of the sign~$\fs_{\os}^*(C)$ in \cite{Wel6,Wel6b}
and the $k\!=\!1$ case of~\eref{WelSolComp3_e0b}, this gives
\BE{P1P1P1_e7a}W_{(1,0,0);0}^{\phi,\os}=1, \qquad 
N^{\phi,\os}_{(1,0,0);(0,0,0),0}=1\,.\EE
The same applies to the numbers $W_{(0,1,0);0}^{\phi,\os},W_{(0,0,1);0}^{\phi,\os}$
and $N^{\phi,\os}_{(0,1,0);(0,0,0),0},N^{\phi,\os}_{(0,0,1);(0,0,0),0}$.
These statements agree with~\eref{P1P1P1_e5a}, \eref{P1P1P1_e3a}, and 
the $b,l\!=\!0$ case of~\eref{P1P1pin5_e7a}.
\end{eg}

\begin{eg}\label{P1P1P1_eg2}
Let $L_1,L_2,L_3$ be as in Example~\ref{P1P1P1_eg2} and 
$L_1^{\R}\!\in\!H_1(\R\P^1\!\times\!\gr(\tau);\Z_2)$ be the homology class of 
the slice $\R\P^1\!\times\!q_2\times\!\tau(q_2)$ with $q_2\!\in\!\R\P^1$.
The pairs 
$$\wt{L}_1\equiv\big(L_1,L_1^{\R}\big) \qquad\hbox{and}\qquad  
\wt{L}_{23}\equiv\big(L_2\!+\!L_3,0\big)$$ 
then freely generate the group $\wt{H}_2^{\phi'}\!(\!(\P^1)^3)$.
The spaces of the slices
$$\P^1\!\times\!q_2\times\!q_3~~\hbox{with}~~q_2,q_3\!\in\!\P^1
\quad\hbox{and}\quad 
q_1\!\times\!\big\{(q_2,h(q_3)\!)\!:q_2\!\in\!\P^1\big\}
~~\hbox{with}~~q_1\!\in\!\P^1,~h\!\in\!\Aut(\P^1),$$
contained in $(\P^1)^3\!-\!\R\P^1\!\times\!\gr(\tau)$
are path-connected and thus determine homology classes~$L_1^c$ and~$L_{23}^c$
in $(\P^1)^3\!-\!\R\P^1\!\times\!\gr(\tau)$.
For $\bfd\!\equiv\!(d_1,d_{23})$ and $\bfa\!\equiv\!(a_1,a_{23})$ with 
$d_i,a_i\!\in\!\Z^{\ge0}$ and $b\!\in\!\Z^{\ge0}$, define
\begin{gather*}
\wt{B}(\bfd)=d_1\wt{L}_1\!+\!d_{23}\wt{L}_{23}, \quad
k_{\bfd}(\bfa,b)=(d_1\!+\!2d_{23})\!-\!(a_1\!+\!a_{23})\!-\!2b,, \quad
W^{\phi',\os'}_{\bfd;b}=W^{\phi',\os'}_{\wt{B}(\bfd),b}\,,\\
N^{\phi',\os'}_{\bfd;\bfa,b}=
\blr{\underset{b}{\underbrace{{\pt,\ldots,\pt}}},
\underset{a_1}{\underbrace{{L_1^c,\ldots,L_1^c}}},
\underset{a_{23}}{\underbrace{{L_{23}^c,\ldots,L_{23}^c}}}
}_{\wt{B}(\bfd);\R\P^1\times\gr(\tau);k_{\bfd}(\bfa,b)}^{\phi',\os'}\,.
\end{gather*}
Similarly to~\eref{RvsCGWs_e},
\begin{gather*}
N_{(d_1,d_{23},d_{23};2a_1,2a_2,2a_3,2b+k_{\bfd}(\bfa,b)}\!-\!
2^{-a_{23}}\big|N^{\phi',\os'}_{(d_1,d_{23});(a_1,a_{23}),b}\big|
\in2\Z^{\ge0}\\ 
\qquad\hbox{if}\qquad
a_2,a_3\!\in\!\Z^{\ge0},~~a_2\!+\!a_3\!=\!a_{23},~~
k_{\bfd}(\bfa,b)\ge0.
\end{gather*}
By Theorem~\ref{RGWs_thm2}\ref{RGWprod_it0},
$$N^{\phi',\os'}_{(0,d_{23});(a_1,a_{23}),b}=
\begin{cases}N^{\tau_{1,1}',\fp_0^-}_{d_{23},a_1},
&\hbox{if}~a_1\!=\!2d_{23}\!-\!1,~a_{23},b\!=\!0;\\
0,&\hbox{otherwise}.
\end{cases}$$
By Theorems~\ref{RGWs_thm2}\ref{RGWprod_it1} and~\ref{RGWs_thm}\ref{div_it},
$$N^{\phi',\os'}_{(1,d_{23});(a_1,a_{23}),b}=
\begin{cases}0,&\hbox{if}~d_{23}\!>\!0,~
a_1\!>\!2d_{23}\!+\!a_{23}\!-\!2;\\
N^{\tau_{1,1}',\fp_0^-}_{d_{23},a_1},
&\hbox{if}~a_1\!=\!2d_{23}\!-\!2,~a_{23},b\!=\!0;\\
N^{\tau_{1,1}',\fp_0^-}_{d_{23},a_1+1},
&\hbox{if}~a_1\!=\!2d_{23}\!-\!2,~a_{23}\!=\!0,~b\!=\!1;\\
2d_{23}N^{\tau_{1,1}',\fp_0^-}_{d_{23},a_1},
&\hbox{if}~a_1\!=\!2d_{23}\!-\!1,~a_{23}\!=\!1,~b\!=\!0.
\end{cases}$$
\end{eg}
Similarly to~\eref{P1P1P1_e7a}
$$W_{(1,0);0}^{\phi',\os'}=1, \qquad 
N^{\phi',\os'}_{(1,0);(0,0),0}=1\,.$$

\subsection{Proofs of Theorems~\ref{WelSolComp_thm} and~\ref{WelSolComp3_thm}}
\label{WelSolPf_thm}

We now deduce \eref{WelSolComp_e0b} and~\eref{WelSolComp3_e0b}
from Proposition~\ref{cNsgn_prp}; see also Remark~\ref{cNsgn_rmk}.
We take $(X,\om,\phi)$, $Y$, $k,l$, $\wt{B}\!\equiv\!(B,b)$, $J$, 
$$\p\!\equiv\!\big((p_1^{\R},\ldots,p_k^{\R}),(p_1^+,\ldots,p_l^+)\big)\in
Y^k\!\times\!(X\!-\!X^{\phi})^l,$$
$C\!\in\!\fM_{\p}(\wt{B};J)$, and $\fp$ or $\os$ as in 
Theorems~\ref{WelSolComp_thm} and~\ref{WelSolComp3_thm}.
Let $h_i$ be the inclusion of~$p_i^+$ into~$X$, $\bh\!\equiv\!(h_1,\ldots,h_l)$,
$\u$ and~$\wt\u$ be tuples as in~\eref{udfn_e} and~\eref{wtudfn_e}
so that their map component~$u$  satisfies~\eref{Cudfn_e} and~\eref{Cudfn_e4},
and $\cC$ as in~\eref{cCdfn_e} be the marked domain of~$\u$.
In particular, $u$ is an immersion.
We define $(\cN u,\vph)$, $\cN_{\u}$, $\cN_{\u}^{\R}$,
$D_{\u}''$, $\ev_{\u}''^{\C}$, $\ev_{\u}''^{\R}$,
$\fp_{\u}''$, $\fo_{\fp;\u}''$, and~$\fs_{\fp;\bh}(\u)$ as in Section~\ref{SignImmer_subs},
with the lower limit of $l\!=\!3$ replaced by $l\!=\!1$ in the definition of~$\cN_{\u}$.
In the case of Theorem~\ref{WelSolComp3_thm}, we denote the last three objects by
$\os_{\u}''$, $\fo_{\os;\u}''$, and~$\fs_{\os;\bh}(\u)$, respectively.
By the definition in Section~\ref{SignImmer_subs}, 
$\fs_{\fp;\bh}(\u)\!=\!+1$ if and only~if the isomorphism
\BE{WelSolComp_e4}\big(\ev_{\u}''^{\C},\ev_{\u}''^{\R}\big)\!:\ker D_{\u}''\lra \cN_{\u}\EE
lies in the homotopy class of isomorphisms $\la(\ker D_{\u}'')\!\approx\!\la(\cN_{\u})$
determined by~$\fo_{\fp;\u}''$ and the symplectic orientation of~$X$.

\begin{proof}[{\bf{\emph{Proof of \eref{WelSolComp_e0b}}}}]
Since~$\p$ is generic, the loop $u_{\R}\!\equiv\!u|_{S^1}$ in~$Y$ is admissible 
in the sense of Theorem~\ref{ShExSeqImm_thm}.
Thus, $\fp_{\u}''\!=\!\fp_0^-(\cN u_{\R})$ if and only~if
\BE{WelSolComp_e5} 1\!+\!\mu_{\fp}(b)\!+\!\de_H(C)=0\in\Z_2.\EE
Along with the CROrient~\ref{CROSpinPinStr_prop}\ref{CROPinStr_it} property,
this implies that the homotopy class~$\fo_{\fp;\u}''$ of isomorphisms 
$\la(\ker D_{\u}'')\!\approx\!\la(\cN_{\u}^{\R})$ is the same as 
the ``canonical" homotopy class~$\fo_{\cC;0}(\cN u,\vph)$ if  and only~if
\eref{WelSolComp_e5} holds.
By Corollary~\ref{NormEv_crl}, the isomorphism~\eref{WelSolComp_e4} lies 
in the homotopy class of such isomorphisms
determined by~$\fo_{\cC;0}(\cN u,\vph)$ and the symplectic orientation of~$X$ if  and only~if
$\de_0^c(\cC)\!\in\!2\Z$.
Thus,
\BE{WelSolComp_e3}\fs_{\fp;\bh}(\u)=(-1)^{\de_0^c(\cC)+1+\mu_{\fp}(b)+\de_H(C)}\,.\EE

\vspace{.15in}

By the definition of the sign $\fs_{\fp}(C)\!\equiv\!\fs_{\fp;0;\bh}(\u)$ 
in Section~\ref{MapSignDfn_subs} and Proposition~\ref{cNsgn_prp},
$$\fs_{\fp}(C)=
(-1)^{\de_0^c(\cC)+\binom{k+1}{2}+1+\lr{w_2(X),B}}\fs_{\fp;\bh}(\u)\,.$$
Combining this with \eref{WelSolComp_e3}, we obtain
$$\fs_{\fp}(C)=(-1)^{\frac{k(k+1)}{2}+\lr{w_2(X),B}+\mu_{\fp}(b)+\de_H(C)}\,.$$
Along with~\eref{nodecond_e} and~\eref{dimcond_e}, 
this establishes~\eref{WelSolComp_e0b}.
\end{proof}

\begin{proof}[{\bf{\emph{Proof of \eref{WelSolComp3_e0b}}}}]
By~\eref{dimcond3_e}, the degree of~$\cN u$ is $2(k\!+\!2l\!-\!1)$ in this case.
By \cite[Proposition~4.1]{BHH},
the real bundle pair $(\cN u,\vph)$ is thus isomorphic to~$2(V,\psi)$ 
for any real line bundle pair $(V,\psi)$ of degree~$k\!+\!2l\!-\!1$.
This isomorphism can be chosen so that it identifies the orientation of~$\cN^{\vph}$
determined by~$\os_{\u}''$ with the canonical orientation on~$2V^{\psi}$ 
as in Example~\ref{CanSpin_eg}\ref{CanSpin_it1}.\\

By the SpinPin~\ref{SpinPinStr_prop}\ref{SpinStr_it} property on page~\pageref{SpinPinStr_prop},  
there are two OSpin-structures compatible  on $\cN u^{\vph}\!=\!2V^{\psi}$ with its canonical orientation.
If $k\!\not\in\!2\Z$, we denote by $\os_0(\cN u^{\vph})$ the canonical OSpin-structure 
on~$\cN u^{\vph}$ as in Example~\ref{CanSpin_eg}\ref{CanSpin_it2}.
This is the preferred OSpin-structure~$\os^{\st}\!(\cN u^{\vph})$ in the first case  
in Section~1.2.2 in~\cite{Wel6b}. 
If $k\!\in\!2\Z$, we denote by $\os_0(\cN u^{\vph})$ the OSpin-structure 
$\os_0(2V^{\psi},\fo_{V^{\psi}}^-)$ as in Examples~\ref{Pin2Spin_eg} and~\ref{SpinDfn1to3_eg};
under the bijection~\eref{SpinPinCorr_e},  $\os_0(2V^{\psi},\fo_{V^{\psi}}^-)$
corresponds to the canonical $\Pin^-$-structure~$\fp_0^-(V^{\psi})$ 
on~$V^{\psi}$ with respect to the orientation of~$S^1$ 
as the boundary of~$\bD^2_+\!\subset\!\P^1$.
As noted in Remark~\ref{SpinDfn1to3_rmk}, the preferred OSpin-structure~$\os^{\st}\!(\cN u^{\vph})$ 
in the second case in \cite[Section~1.2.2]{Wel6b} is~$\os_1(2V^{\psi},\fo_{V^{\psi}}^-)$.\\

Since the degree of the line bundle $V$ is at least $-1$, 
a real CR-operator~$D_{(V,\psi)}$ on the real bundle pair~$(V,\psi)$ is surjective
and the evaluation homomorphism
$$\ker D_{(V,\psi)}\lra V_{\u}\equiv V_{\u}^{\C}\!\oplus\!V_{\u}^{\R}
\equiv \big(V_{z^+_1}\!\oplus\!\ldots\!\oplus\!V_{z^+_l}\big)
\!\oplus\!
\big(V^{\psi}\big|_{x_1}\!\oplus\!\ldots\!\oplus\!V^{\psi}\big|_{x_k}\big)$$
is an isomorphism.
By the surjectivity of~$D_{(V,\psi)}$, $\la(\ker D_{(V,\psi)})\!=\!\la(D_{(V,\psi)})$.
By Proposition~\ref{tauspinorient_prp} and the orienting construction 
described at the beginning of Section~\ref{tauspinorient_subs},
the natural isomorphism
$$\ker D_{2(V,\psi)}\!=\!\ker\big(\!D_{(V,\psi)}\!\oplus\!D_{(V,\psi)}\big)
\approx\big(\!\ker D_{(V,\psi)}\big)\!\oplus\!\big(\!\ker D_{(V,\psi)}\big)$$
is orientation-preserving with respect to the orientation $\fo_0(\cN u^{\vph})$
on the left-hand side induced
by the OSpin-structure $\os_0(\cN u^{\vph})$ on~$\cN u^{\vph}$
and the canonical orientation on the right-hand side.
Thus, the sign of the evaluation isomorphism
\BE{WelSolComp3_e9}\ker D_{2(V,\psi)}\lra\cN_{\u}\EE
with respect to the orientation $\fo_0(\cN u^{\vph})$ on the left-hand side
and the orientation on the right-hand side induced by the orientations of~$\cN u$
and~$\cN u^{\vph}$ is~$(-1)^{\binom{k}{2}}$.\\

By the conclusions of the last two paragraphs and 
the CROrient~\ref{CROSpinPinStr_prop}\ref{CROSpinStr_it} property on page~\pageref{CROSpinPinStr_prop}, 
the sign of the evaluation isomorphism~\eref{WelSolComp3_e9} with respect 
to the orientation on the left-hand side induced
by the OSpin-structure~$\os^{\st}\!(\cN u^{\vph})$ 
and the orientation on the right-hand side induced by the orientations of~$\cN u$
and~$\cN u^{\vph}$ is~$(-1)^{\binom{k}{2}+k+1}$.
Along with the definition of the sign~$\fs_{\os}^*(C)$ in \cite[Section~3.2]{Wel6b}
and the sentence containing~\eref{WelSolComp_e4} above, this gives
\BE{WelSolComp3_e5}  \fs_{\os}^*(C)=(-1)^{\binom{k+1}{2}+1}\fs_{\os;\bh}(\u);\EE
By the definition of the sign $\fs_{\os}(C)\!\equiv\!\fs_{\os;l;\bh}(\u)$ 
in Section~\ref{MapSignDfn_subs} and Proposition~\ref{cNsgn_prp},
\hbox{$\fs_{\os}(C)\!=\!\fs_{\os;\bh}(\u)$}.
Combining this with~\eref{WelSolComp3_e5}, we obtain~\eref{WelSolComp3_e0b}.
\end{proof}

\appendix

\chapter*{Appendices}
\addcontentsline{toc}{chapter}{Appendices}
\renewcommand{\thesection}{\Alph{section}}

\section{\v{C}ech cohomology}
\label{CechH_app}

This appendix contains a detailed review of \v{C}ech cohomology, including for 
sheaves of non-abelian groups, 
describes its connections with singular cohomology and principal bundles,
and classifies oriented vector bundles over bordered surfaces.
We carefully specify the assumptions required in each statement.\\

We generally follow the perspective of \cite[Chapter~5]{Warner}.
In particular, a \sf{sheaf}~$\cS$ over a topological space~$Y$ 
is a topological space along with a projection map $\pi\!:\cS\!\lra\!Y$
so that $\pi$ is a local homeomorphism.
For a \sf{sheaf}~$\cS$ \sf{of~modules} over a ring~$R$ as in \cite{Warner} 
and in Section~\ref{Hisom_subs} below,
$\cS_y\!\equiv\!\pi^{-1}(y)$ is a module over~$R$ for every $y\!\in\!Y$ and
the module operations are continuous with  respect to the topology of~$\cS$.
For a \sf{sheaf}~$\cS$ \sf{of~groups} (not necessarily abelian),
 as in Sections~\ref{Nonabelian_subs}-\ref{CechPB_subs} below,
$\cS_y$ is a group for every $y\!\in\!Y$ and
the group operations are continuous with  respect to the topology of~$\cS$.
For a collection $\{U_{\al}\}_{\al\in\cA}$ of subsets of a space~$Y$
and $\al_0,\al_1,\ldots,\al_p\!\in\!\cA$, we~set
$$U_{\al_0\al_1\ldots\al_p}=U_{\al_0}\!\cap\!U_{\al_1}\!\cap\!\ldots\!\cap\!U_{\al_p}
\subset Y\,.$$

\subsection{Identification with singular cohomology}
\label{Hisom_subs}

For a sufficiently nice topological space~$Y$ and a module~$M$ over a ring~$R$,
the \v{C}ech cohomology group $\wcH^p(Y;M)$ of~$Y$
with coefficients in the sheaf $Y\!\times\!M$ 
of germs of locally constant functions on~$Y$ with values in~$M$
is well-known to be canonically isomorphic to the singular cohomology group
$H^p(Y;M)$ of~$Y$ with coefficients in~$M$.
Proposition~\ref{CechSing_prp} below makes this precise in the $M\!=\!\Z_2$ case 
relevant to our purposes,
making use of the locally $H^k$-simple notion of Definition~\ref{locH1simp_dfn}.
The statement and proof of this proposition apply to an arbitrary module~$M$
over a ring~$R$.
The $p\!=\!1$ case of the isomorphism of Proposition~\ref{CechSing_prp}
is described explicitly at the end of this section.

\begin{prp}\label{CechSing_prp}  
Let $k\!\in\!\Z^{\ge0}$.
For every paracompact locally $H^k$-simple space~$Y$, there exist
canonical isomorphisms
\BE{CechSing_e0}\Phi_Y\!:H^p(Y;\Z_2)\stackrel{\approx}{\lra} \wcH^p(Y;\Z_2),
\qquad p\!=\!0,1\ldots,k\,.\EE
If $Y$ is another paracompact locally $H^k$-simple space and $f\!:Y\!\lra\!Y'$ is 
a continuous map, then the~diagram
\BE{CechSing_e0b}\begin{split}
\xymatrix{ H^p(Y';\Z_2) \ar[rr]^{\Phi_{Y'}}\ar[d]_{f^*}&& \wcH^p(Y';\Z_2) \ar[d]^{f^*} \\
H^p(Y;\Z_2) \ar[rr]^{\Phi_Y}&& \wcH^p(Y;\Z_2) }
\end{split}\EE
commutes for every $p\!\le\!k$.
\end{prp}

\begin{proof}
Let $p\!\in\!\Z^{\ge0}$ and $Y$ be a topological space. 
Denote by $\cS_Y^p\!\lra\!Y$ the sheaf of germs of $\Z_2$-valued singular $p$-cochains on~$Y$ 
as in \cite[5.31]{Warner},  by 
\BE{CechSing_e1}\nd_p\!: \cS_Y^p\lra \cS_Y^{p+1}\EE
the homomorphism induced by the usual differential in the singular cohomology theory,
and by
$$\nd_{p;Y}\!:\Ga\big(Y;\cS_Y^p\big)\lra\Ga\big(Y;\cS_Y^{p+1}\big)$$
the resulting homomorphism between the spaces of global sections.
Let $\cZ^p_Y\!\subset\!\cS_Y^p$ be the kernel of 
the sheaf homomorphism~\eref{CechSing_e1} so that 
\BE{CechSing_e2}\{0\}\lra \cZ^p_Y\lra \cS_Y^p \stackrel{d_p}{\lra} \cZ_Y^{p+1}\EE
is an exact sequence of sheaves.
If $Y$ is locally path-connected,  $\cZ^0_Y\!=\!Y\!\times\!\Z_2$.\\ 

From now on, we assume that $Y$ is paracompact.
By the exactness of~\eref{CechSing_e2}, 
\BE{CechSing_e2b} \ker \nd_{p;Y}=\Ga\big(Y;\cZ^p\big)=\wcH^0\big(Y;\cZ^p\big).\EE
By \cite[p193]{Warner}, each sheaf $\cS_Y^p$ is fine.
By \cite[p202]{Warner}, this implies~that 
\BE{CechSing_e3} \wch{H}^q(Y;\cS_Y^p)=0
\qquad\forall~p\!\in\!\Z^{\ge0},~q\!\in\!\Z^+\,.\EE
Each $\Z_2$-valued singular $p$-cochain~$\vp$ on~$Y$ determines 
a section $(\rho_{y,Y}(\vp))_{y\in Y}$ of $\cS_Y^p$~over~$Y$.
By \cite[5.32]{Warner}, the induced homomorphism
\BE{CechSing_e5}
H^p(Y;\Z_2)\lra H^p\big(\Ga(Y;\cS_Y^*),\nd_{*;Y}\big), \qquad
[\vp]\lra \big[(\rho_{y,Y}(\vp))_{y\in Y}\big],\EE
is an isomorphism.
Combining the $p\!=\!0$ cases of this isomorphism and 
of the identification~\eref{CechSing_e2b}, 
we obtain an isomorphism~\eref{CechSing_e0} for $p\!=\!0$.\\

Suppose $Y$ is locally $H^k$-simple and $p\!\in\!\Z^+$ with $p\!\le\!k$.
The~sequence
\BE{CechSing_e7}
\{0\}\lra \cZ^{p-q-1}_Y\lra \cS_Y^{p-q-1}\lra \cZ^{p-q}_Y \lra\{0\}\EE
of sheaves is then exact for every $q\!\in\!\Z^{\ge0}$ with $q\!<\!p$.
From the exactness of the associated long sequence in \v{C}ech cohomology, \eref{CechSing_e2b},
and~\eref{CechSing_e3}, 
we obtain isomorphisms
\begin{gather*}
\wch\de_Y\!:  H^p\big(\Ga(Y;\cS_Y^*),d_{*;Y}\big)
\equiv\frac{\ker \nd_{p;Y}}{\Im\,\nd_{p-1;Y}} 
=\frac{\wcH^0(Y;\cZ^p)}{\nd_{p-1}(\wcH^0(Y;\cS^{p-1}))}
\stackrel{\approx}{\lra}\wcH^1\big(Y;\cZ_Y^{p-1}\big),\\
\wch\de_Y\!: \wcH^q(Y;\cZ_Y^{p-q}) \stackrel{\approx}{\lra} 
\wcH^{q+1}\big(Y;\cZ_Y^{p-q-1}\big)
~~\forall~q\!\in\!\Z^+,~q\!<\!p.
\end{gather*}
Putting these isomorphisms together, we obtain an isomorphism
\BE{CechSing_e11}\wch\de_Y^p\!: H^p\big(\Ga(Y;\cS_Y^*),d_{*;Y}\big)\lra 
\wcH^p\big(Y;\cZ_Y^0\big)\!=\!\wcH^p\big(Y;\Z_2\big).\EE
Combining~\eref{CechSing_e5} with this isomorphism, we obtain 
an isomorphism as in~\eref{CechSing_e0}  with $p\!>\!0$.\\

Suppose $Y$ is another paracompact locally $H^k$-simple space.
A continuous map \hbox{$f\!:Y\!\lra\!Y'$} induces commutative diagrams
$$\xymatrix{H^p(Y';\Z_2)\ar[r]^>>>>>{\approx}\ar[d]_{f^*}& 
H^p\big(\Ga(Y';\cS_{Y'}^*),d_{*;Y'}\big)\ar[d]^{f^*} &
\!\!\!\!\!\!\!\!\!\{0\}\ar[r]& \cZ^{p-q-1}_{Y'}\ar[r]\ar[d]^{f^*}& 
\cS_{Y'}^{p-q-1}\ar[r]\ar[d]^{f^*}& \cZ^{p-q}_{Y'}\ar[d]^{f^*} \ar[r]&\{0\}\\
H^p(Y;\Z_2)\ar[r]^>>>>>{\approx}& H^p\big(\Ga(Y;\cS_Y^*),d_{*;Y}\big)&
\!\!\!\!\!\!\!\!\!
\{0\}\ar[r]& \cZ^{p-q-1}_Y\ar[r]& \cS_Y^{p-q-1}\ar[r]& \cZ^{p-q}_Y \ar[r]&\{0\}}$$
for all $p,q\!\in\!\Z^{\ge0}$ with $q\!<\!p\!\le\!k$.
Combining the $p\!=\!0$ cases of the first diagram above and
of the identifications~\eref{CechSing_e2b} for~$Y$ and for~$Y'$, 
we obtain a commutative diagram~\eref{CechSing_e0b} for $p\!=\!0$.   
The second commutative diagram above induces a commutative diagram
$$\xymatrix{ 
H^p\big(\Ga(Y';\cS_{Y'}^*),\nd_{*;Y'}\big)\ar[r]^>>>>>>{\approx}\ar[d]_{f^*}& 
\wcH^p\big(Y';\Z_2\big)\ar[d]^{f^*}\\
H^p\big(\Ga(Y;\cS_Y^*),\nd_{*;Y}\big)\ar[r]^>>>>>>{\approx}&  \wcH^p\big(Y;\Z_2\big)}$$
with the horizontal isomorphisms as in~\eref{CechSing_e11}.
Combining this with the first commutative diagram in this paragraph, 
we obtain~\eref{CechSing_e0b} with $p\!>\!0$.
\end{proof}   

Let $Y$ be a paracompact locally $H^1$-simple space. 
We now describe the $p\!=\!1$ case of the isomorphism~\eref{CechSing_e0} explicitly. 
Suppose $\vp$ is a $\Z_2$-valued singular 1-cocycle on~$Y$.
Since $Y$ is locally $H^1$-simple,
there exist an open cover $\{U_{\al}\}_{\al\in\cA}$ of~$Y$ and
a $\Z_2$-valued singular 0-cochain~$\mu_{\al}$ on~$U_{\al}$ for each $\al\!\in\!\cA$
so~that 
$$\nd_{0;U_{\al}}\mu_{\al}=\vp\big|_{U_{\al}}\qquad\forall\,\al\!\in\!\cA.$$
We define a \v{C}ech 1-cocycle $\eta$ on~$Y$ by
\BE{PhiYetadfn_e}\eta_{\al\be}=\mu_{\be}\big|_{U_{\al\be}}-\mu_{\al}\big|_{U_{\al\be}}
\in \cS_Y^1\big(U_{\al\be}\big)\qquad\forall\,\al,\be\!\in\!\cA.\EE
Since $\nd_{0;U_{\al\be}}\eta_{\al\be}\!=\!0$ and $Y$ is locally path-connected,
$\eta_{\al\be}$ is a locally constant function on $U_{\al\be}$.
Thus, $\eta$ takes values in the sheaf of germs of $\Z_2$-valued continuous functions
on~$Y$ and so defines an element $[\eta]$ of $\wch{H}^1(Y;\Z_2)$.
This is the image of~$[\vp]$ under the $p\!=\!1$ case of 
the isomorphism~$\Phi_Y$  in~\eref{CechSing_e0}.\\

Suppose $Y$ is a CW complex and $\vp$ is a $\Z_2$-valued singular 1-cocycle on~$Y$
as above.
For each vertex $\al\!\in\!Y_0$ of~$Y$, let $U_{\al}\!\subset\!Y$ denote 
the (\sf{open}) \sf{star} of~$\al$, i.e.~the union of all open cells~$\mr{e}$ of~$Y$ 
so that~$\al$ is contained in the closed cell~$e$.
In particular, $U_{\al}$ is an open neighborhood of~$\al$ and the collection
$\{U_{\al}\}_{\al\in Y_0}$ covers the 1-skeleton $Y_1\!\!\subset\!Y$.
We~take 
$$\{U_{\al}\}_{\al\in\cA}\equiv\{U_{\al}\}_{\al\in Y_0}\sqcup\big\{Y\!-\!Y_1\};$$
this is an open cover of~$Y$.
By adding extra vertices to~$Y_1$, 
we can ensure that no closed 1-cell is a cycle.
This implies that every closed 1-cell~$e$ of~$Y$ runs between 
distinct vertices~$\al$ and~$\be$ with
\BE{CWcorr_e0}
e\subset U_{\al}\!\cup\!U_{\be}, \qquad e\cap\!U_{\al\be}=\mr{e},\qquad
 e\cap U_{\ga}=\eset~~\forall\,\ga\!\in\!\cA\!-\!\{\al,\be\}.\EE
For every $\al\!\in\!\cA$, there then exists
a $\Z_2$-valued singular 0-cochain~$\mu_{\al}$ on~$U_{\al}$  so~that 
$$\nd_{0;U_{\al}}\mu_{\al}=\vp\big|_{U_{\al}},\quad
\mu_{\al}(\al)=0\qquad\forall~\al\!\in\!\cA.$$
Every closed 1-cell~$e$ of~$Y$  is cobordant to 
the difference of a singular 1-simplex $e_{x\be}$ running from a point $x\!\in\!\mr{e}$
to~$\be$ and a singular 1-simplex $e_{x\al}$ running from~$x$ to~$\al$.
By~\eref{CWcorr_e0}, $e_{x\al}\!\subset\!U_{\al}$ and $e_{x\be}\!\subset\!U_{\be}$.
Since $\vp$ is a cocycle, it follows~that
\begin{equation*}\begin{split}
\vp(e)=\vp\big(e_{x\be}\!-\!e_{x\al}\big)
&=\vp\big(e_{x\be}\big)-\vp\big(e_{x\al}\big)
=\big\{\nd_{0;U_{\be}}\mu_{\be}\big\}\big(e_{x\be}\big)
-\big\{\nd_{0;U_{\al}}\mu_{\al}\big\}\big(e_{x\al}\big)\\
&=\big(\mu_{\be}(\be)\!-\!\mu_{\be}(x)\big)-
\big(\mu_{\al}(\al)\!-\!\mu_{\al}(x)\big)
=\mu_{\al}(x)\!-\!\mu_{\be}(x).
\end{split}\end{equation*}
Along with~\eref{PhiYetadfn_e}, this implies that
the \v{C}ech cohomology class $[\eta]\!\equiv\!\Phi_Y([\vp])$ 
corresponding to~$[\vp]$ under the isomorphism~\eref{CechSing_e0}
is represented by a collection $\{\eta_{\al\be}\}_{\al,\be\in\cA}$
associated with an open cover $\{U_{\al}\}_{\al\in\cA}$ of~$Y$ such~that 
$$ \mr{e}\subset U_{\al\be}, \qquad \eta_{\al\be}\big|_{\mr{e}}=\vp(e)\in\Z_2$$
for all $\al,\be\!\in\!Y_0$ and every closed 1-cell~$e$ with vertices~$\al$ and~$\be$.

\subsection{Sheaves of groups}
\label{Nonabelian_subs}

\noindent
\v{C}ech cohomology groups $\wcH^p$ are normally defined for sheaves or presheaves 
of (abelian) modules over a ring. 
The sets~$\wcH^0$ and $\wcH^1$ can be defined for sheaves or presheaves 
of non-abelian groups as well.
The first set is still a group, while the second is a \sf{pointed set},
i.e.~it has a distinguished element.
A short exact sequence of such sheaves gives rise to 
an exact sequence of the corresponding \v{C}ech pointed sets,
provided the kernel sheaf~$\cR$ lies in the center~$Z(\cS)$ of  
the ambient sheaf~$\cS$; see Proposition~\ref{Snake_prp}.
The main examples of interest are the sheaves~$\cS$ of germs of 
continuous functions over a topological space~$Y$ 
with values in a Lie group~$G$, as in Section~\ref{CechPBs_subs}.\\

\noindent
We denote the center of a group~$G$ by~$Z(G)$.
We call a collection 
$$\big(\big(\de_p\!:C^p\!\lra\!C^{p+1}\big)_{p=0,1,2},*\!:C^0\!\times\!C^1\!\lra\!C^1\big)$$
consisting of maps~$\de_p$ between groups~$C^p$ with the identity element~$\1_p$
and a left action~$*$ a \sf{short cochain complex}~if
\begin{gather}\label{CCprp_e11}
\de_p\1_p=\1_{p+1}, \quad \de_{p+1}\!\circ\!\de_p=\1_{p+2},\quad
\de_1(f\!*\!g)=\de_1g ~~\forall\,f\!\in\!C^0,\,
g\!\in\!\de_1^{-1}\big(Z(C^2)\big),\\
\label{CCprp_e12}
\de_0(f\!\cdot\!f')=f\!*\!(\de_0f'),~~f\!*\!g=(\de_0f)g
\quad
\forall\,f,f'\!\in\!C^0,\,g\!\in\!Z(C^1),\\ 
\label{CCprp_e13}
\de_p\big(g\!\cdot\!g'\big)=\big(\de_pg\big)\big(\de_pg'\big) \quad
\forall\,g\!\in\!C^p,\,g'\!\in\!Z(C^p),\,p\!=\!1,2.
\end{gather}
By the second condition in~\eref{CCprp_e12}, 
\BE{CCprp_e15}
C^0\!*\!\{g\}=\big(\Im\,\de_0\big)\!\cdot\!\{g\} \quad\forall\,\,g\!\in\!Z(C^1).\EE
By both conditions in~\eref{CCprp_e12}, 
\BE{H0dfn_e}H^0(C^*)\equiv
H^0\big((C^p,\de_p)_{p=0,1,2},*\big) \equiv\ker\de_0\equiv\de_0^{-1}(\1_1)\EE
is a subgroup of~$C^0$.
By the last property in~\eref{CCprp_e11}, $*$ restricts to an action on 
$\ker\de_1\!\equiv\!\de_1^{-1}(\1_2)$.
We can thus define
\BE{H1dfn_e} H^1(C^*)\equiv
H^1\big((C^p,\de_p)_{p=0,1,2},*\big)=\ker\de_1\big/C^0\,;\EE
this is a pointed set with the distinguished element given by the image of
$\Im\,\de_0\!\ni\!\1_1$ in~$H^1(C^*)$.\\

\noindent
By~\eref{CCprp_e12} and the $p\!=\!1$ case of~\eref{CCprp_e13},
$\de_0$ and $\de_1$ are group homomorphisms if the group $C^1$ is abelian
and $*$~is the usual action of the 1-coboundaries on the 1-cochains 
via the group operation.
In this case, \eref{H1dfn_e} agrees with the usual definition and
the last condition in~\eref{CCprp_e11} is automatically satisfied.
If in addition the group~$C^2$ is also abelian, 
as happens for the kernel complex~$B^*$ in Lemma~\ref{Snake_lmm} below, 
then the map~$\de_2$ is a group homomorphism as~well and 
$$H^2(C^*)\equiv
H^2\big((C^p,\de_p)_{p=0,1,2},*\big)=\ker\de_2\big/\Im\,\de_1$$
is a well-defined abelian group.\\

\noindent
A \sf{morphism} of short cochain complexes
$$\io\!\equiv\!(\io_p)_{p=0,1,2,3}\!:
\big((B^p,\de_p)_{p=0,1,2},*\big)\lra 
\big((C^p,\de_p)_{p=0,1,2},*\big)$$
is a collection of group homomorphisms $\io_p\!:B^p\!\lra\!C^p$ that 
commute with the maps~$\de_p$ and the actions~$*$.
Such a homomorphism induces morphisms
$$\io_*\!:H^p(B^*)\lra H^p(C^*),\quad p=0,1,$$
of pointed sets, i.e.~$\io_*$ takes the distinguished element of the domain
to the distinguished element of the target; the map~$\io_0$ is a group homomorphism.
The \sf{kernel} of such a morphism is the preimage of the distinguished element
of the target.
The next lemma is an analogue of the Snake Lemma \cite[Proposition~5.17]{Warner}
for short cochain complexes of groups.

\begin{lmm}\label{Snake_lmm}
For every short exact sequence
$$\{\1\}\lra \big((B^p,\de_p)_{p=0,1,2},*\big) \stackrel{\io}{\lra}
\big((C^p,\de_p)_{p=0,1,2},*\big) \stackrel{\fj}{\lra}
\big((D^p,\de_p)_{p=0,1,2},*\big)\lra\{\1\}$$
of short cochain complexes of groups
such that $\io_p(B^p)\!\subset\!Z(C^p)$ for $p\!=\!1,2$,
there exist morphisms 
\BE{Snake_e0a}\prt_p\!:H^p(D^*)\lra H^{p+1}(B^*), \quad p=0,1,\EE
of pointed sets such that the sequence
\BE{Snake_e0}\begin{split}
\{\1\}&\lra  H^0(B^*)\stackrel{\io_*}{\lra} H^0(C^*)\stackrel{\fj_*}{\lra} H^0(D^*) 
\stackrel{\prt_0}{\lra}\\ 
&\stackrel{\prt_0}{\lra}
 H^1(B^*)\stackrel{\io_*}{\lra} H^1(C^*)\stackrel{\fj_*}{\lra} H^1(D^*) 
\stackrel{\prt_1}{\lra} H^2(B^*)
\end{split}\EE
of morphisms of pointed sets is exact.
The maps~$\prt_p$ are natural with respect to morphisms of short exact sequences of 
short cochain complexes of groups.
\end{lmm}

\begin{proof} We proceed as in the abelian case. 
Given $d_p\!\in\!\ker\de_p\!\subset\!D^p$, let $c_p\!\in\!C^p$ be such that 
$\fj_p(c_p)\!=\!d_p$.
Since
$$\fj_{p+1}\big(\de_p(c_p)\big)=\de_p\big(\fj_p(c_p)\big)
=\de_p(d_p)=\1_{p+1}\in D^{p+1},$$
there exists a unique $b_{p+1}\!\in\!B^{p+1}$ such that $\io_{p+1}(b_{p+1})\!=\!\de_p(c_p)$.
By the second condition in~\eref{CCprp_e11}, $b_{p+1}\!\in\!\ker\de_{p+1}$.
We~set 
$$\prt_p\big([d_p]\big)=\big[b_{p+1}\big]\in H^{p+1}(B^*).$$
By the first condition in~\eref{CCprp_e12}, $[b_1]$ is independent 
of the choice of $c_0\!\in\!C^0$ such that \hbox{$\fj_0(c_0)\!=\!d_0$}.
By the $p\!=\!1$ case of~\eref{CCprp_e13} and 
the assumption that  $\io_p(B^p)\!\subset\!Z(C^p)$ for $p\!=\!1,2$,
$[b_2]$ is independent 
of the choice of $c_1\!\in\!C^1$ such that \hbox{$\fj_1(c_1)\!=\!d_1$}.
By the last condition in~\eref{CCprp_e11} and the assumption that 
$\io_2(B^2)\!\subset\!Z(C^2)$, $[b_2]$ does not depend
on the choice of representative~$d_1$ for~$[d_1]$.
Thus, the maps~\eref{Snake_e0a} are well-defined.
By the first condition in~\eref{CCprp_e11}, $\prt_p([\1_p])\!=\![\1_{p+1}]$,
i.e.~$\prt_p$ is a morphism of~pointed sets.
By the construction, the maps~$\prt_p$ are natural with respect to morphisms of 
 exact sequences of short cochain complexes.\\   

It is immediate that \eref{Snake_e0} is exact at $H^0(B^*)$ and $H^0(C^*)$ and that 
$$\fj_*\!\circ\!\io_*\!=\![\1_1]\!:H^1(B^*)\!\lra\!H^1(D^*), \quad
\prt_p\!\circ\!\fj_*\!=\![\1_{p+1}]\!:H^p(C^*)\!\lra\!H^{p+1}(B^*)\,.$$
The exactness of~\eref{Snake_e0} at~$H^1(B^*)$ is immediate from~\eref{CCprp_e15} 
with $g\!=\!\1_1\!\in\!Z(C^1)$.
The exactness  at $H^1(C^*)$ follows from~\eref{CCprp_e15} 
with $g\!=\!\1_1\!\in\!Z(D^1)$,
the second condition in~\eref{CCprp_e12}, and  the assumption that 
$\io_1(B^1)\!\subset\!Z(C^1)$.
The exactness at $H^0(D^*)$ follows from~\eref{CCprp_e15} 
with $g\!=\!\1_1\!\in\!Z(B^1)$, 
both conditions in~\eref{CCprp_e12}, and
the assumption that $\io_1(B^1)\!\subset\!Z(C^1)$.
The exactness at $H^1(D^*)$ follows from the assumption that 
$\io_p(B^p)\!\subset\!Z(C^p)$ for $p\!=\!1,2$ and
the $p\!=\!1$ case of~\eref{CCprp_e13}. 
\end{proof}

We next review the definitions and key properties of the group~$\wcH^0$
and pointed set~$\wcH^1$ for a sheaf~$\cS$ of groups over a topological space~$Y$.
We denote by $Z(\cS)\!\subset\!\cS$ the subsheaf consisting of 
the centers~$Z(\cS_y)$ of the groups~$\cS_y$ with $y\!\in\!Y$
and by $\1_y\!\in\!\cS_y$ the identity element of~$\cS_y$.\\

\noindent
Let $\un{U}\!\equiv\!\{U_{\al}\}_{\al\in\cA}$ be an open cover of~$Y$.
As in the abelian case, the~set 
$$\wch{C}^p(\un{U};\cS)\equiv
\prod_{\al_0,\al_1,\ldots,\al_p\in\cA} \hspace{-.26in}
\Ga\big(U_{\al_0\al_1\ldots\al_p};\cS\big)$$
of \sf{\v{C}ech $p$-cochains} is a group 
under pointwise multiplication of sections:
\begin{gather*}
\cdot: \wch{C}^p(\un{U};\cS)\times\wch{C}^p(\un{U};\cS)
\lra\wch{C}^p(\un{U};\cS),\\
\{h\cdot h'\}_{\al_0\al_1\ldots\al_p}(y)=h_{\al_0\al_1\ldots\al_p}(y)\cdot
h'_{\al_0\al_1\ldots\al_k}(y)~~\forall\,\al_0,\al_1,\ldots,\al_p\!\in\!\cA,\,
y\!\in\!U_{\al_0\al_1\ldots\al_p}\,.
\end{gather*}
The identity element $\1_p\!\in\!\wch{C}^p(\un{U};\cS)$ is given~by
$$(\1_p)_{\al_0\al_1\ldots\al_p}(y)=\1_y\qquad  
\forall\,\al_0,\al_1,\ldots,\al_p\!\in\!\cA,\,
y\!\in\!U_{\al_0\al_1\ldots\al_p}\,.$$

\vspace{.1in}

Define the boundary maps by
\begin{gather*} 
\de_0\!: \wch{C}^0(\un{U};\cS)\lra\wch{C}^1(\un{U};\cS),  \quad
(\de_0 f)_{\al_0\al_1}=f_{\al_0}\big|_{U_{\al_0\al_1}}
\cdot f_{\al_1}^{-1}\big|_{U_{\al_0\al_1}}\,,\\
\de_1\!: \wch{C}^1(\un{U};\cS)\lra\wch{C}^2(\un{U};\cS), \quad
 (\de_1 g)_{\al_0\al_1\al_2}=
g_{\al_1\al_2}\big|_{U_{\al_0\al_1\al_2}}
\!\!\cdot g_{\al_0\al_2}^{-1}\big|_{U_{\al_0\al_1\al_2}}
\!\!\cdot g_{\al_0\al_1}\big|_{U_{\al_0\al_1\al_2}}\,.
\end{gather*} 
We also define a left action of $\wch{C}^0(\un{U};\cS)$ on $\wch{C}^1(\un{U};\cS)$ by
\begin{gather*}
*\!: \wch{C}^0(\un{U};\cS)\times\wch{C}^1(\un{U};\cS)\lra 
\wch{C}^1(\un{U};\cS),\\
\{f\!*\!g\}_{\al_0\al_1}=f_{\al_0}\big|_{U_{\al_0\al_1}}
\!\cdot g_{\al_0\al_1}\!\cdot  f_{\al_1}^{-1}\big|_{U_{\al_0\al_1}}
\!\!\in\Ga(U_{\al_0\al_1};\cS).
\end{gather*}

\vspace{.1in}

We now construct a short cochain complex.
Let
$$C^p\big(\un{U};\cS\big)=\begin{cases}\wch{C}^p(\un{U};\cS),&\hbox{if}~p\!=\!0,1,2;\\
\Abel(\wch{C}^p(\un{U};\cS))&\hbox{if}~p\!=\!3.
\end{cases}$$
For $p\!=\!0,1$, we take 
$$\de_p\!:C^p\big(\un{U};\cS\big)\lra C^{p+1}\big(\un{U};\cS\big)$$ 
to be  the maps defined above.
We take~$\de_2$ to be the composition of the~map
\begin{gather*}
\de_2\!:\wch{C}^2(\un{U};\cS)\lra\wch{C}^3(\un{U};\cS), \\
\big(\de_2h\big)_{\al_0\al_1\al_2\al_3}=
h_{\al_1\al_2\al_3}\big|_{U_{\al_0\al_1\al_2\al_3}}
h_{\al_0\al_2\al_3}^{\,-1}\big|_{U_{\al_0\al_1\al_2\al_3}}
h_{\al_0\al_1\al_3}\big|_{U_{\al_0\al_1\al_2\al_3}}
h_{\al_0\al_1\al_2}^{\,-1}\big|_{U_{\al_0\al_1\al_2\al_3}},
\end{gather*}
with the projection $\wch{C}^3(\un{U};\cS)\!\lra\!C^3\big(\un{U};\cS\big)$.
The tuple 
$$\big(\big(\de_p\!:C^p(\un{U};\cS)\!\lra\!C^{p+1}(\un{U};\cS)\big)_{p=0,1,2},
*\!:C^0(\un{U};\cS)\!\times\!C^1(\un{U};\cS)\!\lra\!C^1(\un{U};\cS)\big)$$
is then a short cochain complex of groups.
We denote the associated group~\eref{H0dfn_e} and the pointed set~\eref{H1dfn_e}
by $\wch{H}^0(\un{U};\cS)$ and $\wch{H}^1(\un{U};\cS)$, respectively.\\

Let $\un{U}'\!\equiv\!\{U_{\al}'\}_{\al\in\cA'}$ be an open cover of~$Y$
refining~$\un{U}$, i.e.~there exists a map \hbox{$\mu\!:\cA'\!\lra\!\cA$} such that 
$U_{\al}'\!\subset\!U_{\mu(\al)}$ for every $\al\!\in\!\cA'$.
Such a \sf{refining map} induces group homomorphisms
\begin{gather}\label{mupdfn_e}
\mu_p^*\!:\wch{C}^p(\un{U};\cS)\lra \wch{C}^p(\un{U}';\cS), \\
\notag
\big(\mu_p^*h\big)_{\al_0\ldots\al_p}=h_{\mu(\al_0)\ldots\mu(\al_p)}
\big|_{U_{\al_0\ldots\al_p}'} \quad \forall~h\!\in\!\wch{C}^p(\un{U};\cS),~
\al_0,\ldots,\al_p\!\in\!\cA'. 
\end{gather}
These homomorphisms commute with $\de_0$, $\de_1$, and the action of 
$\wch{C}^0(\cdot;\cS)$ on~$\wch{C}^1(\cdot;\cS)$.
Thus, $\mu$ induces~maps
\BE{nonabelCech_e5}R_{\un{U}',\un{U}}^0\!: \wcH^0(\un{U};\cS)\lra \wcH^0(\un{U}';\cS)
\quad\hbox{and}\quad
R_{\un{U}',\un{U}}^1\!: \wcH^1(\un{U};\cS)\lra \wcH^1(\un{U}';\cS)\EE
of pointed sets;
the first map above is a group homomorphism.\\

\noindent
If $\mu'\!:\cA'\!\lra\!\cA$ is another refining map, then
${U}_{\al}'\!\subset\!{U}_{\mu(\al)\mu'(\al)}$ for every $\al\!\in\!\cA'$ and thus
$$\mu_0^*\big|_{\ker\de_0}=\mu_0'^*\big|_{\ker\de_0}\!:
\ker\de_0\lra \ker\de_0\subset \wch{C}^0(\un{U}';\cS).$$
For $g\!\in\!\wch{C}^1(\un{U};\cS)$, define
$$h_1g\in\wch{C}^0(\un{U}';\cS) \qquad\hbox{by}\qquad
(h_1g)_{\al}=g_{\mu'(\al)\mu(\al)}|_{U_{\al}'}.$$
If $g\!\in\!\ker\de_1\!\subset\!\wch{C}^1(\un{U};\cS)$, then
\begin{equation*}\begin{split}
g_{\mu(\al_0)\mu(\al_1)}\big|_{{U}_{\mu'(\al_1)\mu(\al_0)\mu(\al_1)}}
\cdot g_{\mu'(\al_1)\mu(\al_1)}^{-1}\big|_{{U}_{\mu'(\al_1)\mu(\al_0)\mu(\al_1)}}
&=g_{\mu'(\al_1)\mu(\al_0)}^{-1}\big|_{{U}_{\mu'(\al_1)\mu(\al_0)\mu(\al_1)}},\\
g_{\mu'(\al_0)\mu(\al_0)}\big|_{{U}_{\mu'(\al_1)\mu'(\al_0)\mu(\al_0)}}
\cdot g_{\mu'(\al_1)\mu(\al_0)}^{-1}\big|_{{U}_{\mu'(\al_1)\mu'(\al_0)\mu(\al_0)}}
&=g_{\mu'(\al_1)\mu'(\al_0)}^{-1}\big|_{{U}_{\mu'(\al_1)\mu'(\al_0)\mu(\al_0)}}
\end{split}\end{equation*}
for all $\al_0,\al_1\!\in\!\cA$.
From this, we find that 
$$\mu_1'^*g=(h_1g)*(\mu_1^*g) ~~~\forall~g\!\in\!\ker\de_1\subset \wch{C}^1(\un{U};\cS).$$

\vspace{.1in}

By the previous paragraph,
the pointed maps~\eref{nonabelCech_e5} are independent of the choice of 
refining map \hbox{$\mu\!:\cA'\!\lra\!\cA$}.
We can therefore define the group $\wcH^0(Y;\cS)$ and 
the pointed set $\wcH^1(Y;\cS)$
as the direct limits of the groups $\wcH^0(\un{U};\cS)$ and 
of the pointed sets $\wcH^1(\un{U};\cS)$, respectively,
over open covers of~$Y$.
The~map
\BE{chH0_e}\Ga(Y;\cS)\lra \wcH^0(Y;\cS), \qquad
f\lra \big(f|_{U_{\al}}\big)_{\al\in\cA}\,,\EE
is a group isomorphism.\\

If $\cS$ is a sheaf of abelian groups, 
as happens for the kernel sheaf~$\cR$ in Proposition~\ref{Snake_prp} below, 
the definitions of $\wcH^0(Y;\cS)$ and $\wcH^1(Y;\cS)$ above reduce to
the ones in \cite[Section~5.33]{Warner}.
Furthermore, 
$$\wcH^2(\un{U};\cS)\equiv
\frac{\ker(\de_2\!:\wch{C}^2(\un{U};\cS)\!\lra\!\wch{C}^3(\un{U};\cS))}
{\Im(\de_1\!:\wch{C}^1(\un{U};\cS)\!\lra\!\wch{C}^2(\un{U};\cS))}$$
is a well-defined abelian group for every open cover~$\un{U}$ of~$Y$. 
The group homomorphisms
$$R_{\un{U}',\un{U}}^2\!: \wcH^2(\un{U};\cS)\lra \wcH^2(\un{U}';\cS)$$
induced by refining maps still depend only on the covers $\un{U}$ and~$\un{U}'$.
The abelian group  $\wcH^2(Y;\cS)$ is again the direct limit of the groups 
$\wcH^2(\un{U};\cS)$ over all open covers~$\un{U}$ of~$Y$.\\

A homomorphism $\io\!:\cR\!\lra\!\cS$ of sheaves of groups over~$Y$ induces 
maps
$$\io_*\!:\Ga(Y;\cR)\lra \Ga(Y;\cS),\quad
\io_*\!:\wcH^0(Y;\cR)\lra \wcH^0(Y;\cS), 
\quad \io_*\!:\wcH^1(Y;\cR)\lra \wcH^1(Y;\cS)$$
between pointed spaces.
The first two maps are group homomorphisms which commute with 
the identifications~\eref{chH0_e}.

\begin{prp}\label{Snake_prp}
Let $Y$ be a paracompact space.
For every short exact sequence
\BE{Snake2_e0a}\{\1\}\lra \cR\stackrel{\io}{\lra} \cS\stackrel{\fj}{\lra} \cT\lra\{\1\}\EE
of sheaves of groups over~$Y$ such that $\io(\cR)\!\subset\!Z(\cS)$,
there exist morphisms 
\BE{Snake2_e0b}\wch\de_p\!: \wcH^p(Y;\cT)\lra \wcH^{p+1}(Y;\cR), \quad p=0,1,\EE
of pointed sets such that the sequence
\BE{Snake2_e0}\begin{split}
\{\1\}&\lra  \wcH^0(Y;\cR)\stackrel{\io_*}{\lra}\wcH^0(Y;\cS)\stackrel{\fj_*}{\lra} \wcH^0(Y;\cT)
\stackrel{\wch\de_0}{\lra}\\ 
&\stackrel{\wch\de_0}{\lra}
 \wcH^1(Y;\cR)\stackrel{\io_*}{\lra} \wcH^1(Y;\cS)\stackrel{\fj_*}{\lra} \wcH^1(Y;\cT) 
\stackrel{\wch\de_1}{\lra} \wcH^2(Y;\cR)
\end{split}\EE
of morphisms of pointed sets is exact.
The maps~$\wch\de_p$ are natural with respect to morphisms of short exact sequences 
of sheaves of groups over~$Y$. 
\end{prp}

\begin{proof}
Let $\un{U}\!\equiv\!\{U_{\al}\}_{\al\in\cA}$ be an open cover of~$Y$,
$$B^p(\un{U})=C^p\big(\un{U};\cR\big), \quad
C^p(\un{U})=C^p\big(\un{U};\cS\big), \quad
D^p(\un{U})=C^p\big(\un{U};\cT\big).$$
Since $\io(\cR)\!\subset\!Z(\cS)$, $\io_*(B^p(\un{U}))\!\subset\!Z(C^p(\un{U}))$
for all~$p$.
By the exactness of~\eref{Snake2_e0a},  the sequence 
$$\{\1\}\lra B^p(\un{U})\stackrel{\io_*}{\lra}
C^p(\un{U})\stackrel{\fj_*}{\lra} D^p(\un{U})$$
of groups is exact.
For $p\!=\!0,1,2,3$,
we denote by $\ov{D}^p(\un{U})\!\subset\!D^p(\un{U})$ the image of~$\fj_*$.
For $p\!=\!0,1$, let
$\ov{H}^p(\un{U};\cT)$ be the \v{C}ech pointed sets determined
by the short cochain complex~$\ov{D}^*(\un{U})$.\\

The sequence
\BE{Snake2_e4}\{\1\}\lra B^p(\un{U})\stackrel{\io_*}{\lra}
C^p(\un{U})\stackrel{\fj_*}{\lra} \ov{D}^p(\un{U})\lra\{\1\}\EE
of short cochain complexes is exact.
By Lemma~\ref{Snake_lmm}, there thus exist morphisms 
\BE{Snake2_e5b}\wch\de_p\!: \ov{H}^p(\un{U};\cT)\lra \wcH^{p+1}(\un{U};\cR), \quad p=0,1,\EE
of pointed sets such that the sequence
\BE{Snake2_e5}\begin{split}
\{\1\}&\lra  \wcH^0(\un{U};\cR)\stackrel{\io_*}{\lra}\wcH^0(\un{U};\cS)\stackrel{\fj_*}{\lra} 
\ov{H}^0(\un{U};\cT)
\stackrel{\wch\de_0}{\lra}\\ 
&\stackrel{\wch\de_0}{\lra}
 \wcH^1(\un{U};\cR)\stackrel{\io_*}{\lra} \wcH^1(\un{U};\cS)\stackrel{\fj_*}{\lra} 
\wcH^1(\un{U};\cT) 
\stackrel{\wch\de_1}{\lra} \wcH^2(\un{U};\cR)
\end{split}\EE
of morphisms of pointed sets is exact.\\

Let $\un{U}'\!\equiv\!\{{U}_{\al}'\}_{\al\in\cA'}$ be an open cover of~$Y$
refining~$\un{U}$ and \hbox{$\mu\!:\cA'\!\lra\!\cA$}
be a refining map.
By the naturality of the morphisms~\eref{Snake2_e5b}, 
the group homomorphisms~\eref{mupdfn_e} induce commutative diagrams
$$\xymatrix{ \ov{H}^p(\un{U};\cT) \ar[d]_{R_{\un{U}',\un{U}}^p}
\ar[rr]^{\wch\de_p}&& \wcH^{p+1}(\un{U};\cR) \ar[d]^{R_{\un{U}',\un{U}}^{p+1}}\\
 \ov{H}^p(\un{U}';\cT)\ar[rr]^{\wch\de_p}&& \wcH^{p+1}(\un{U}';\cR)}$$
of pointed sets.
Taking the direct limit of the morphisms~\eref{Snake2_e5b} over all open covers
of~$Y$, we thus obtain 
morphisms 
\BE{Snake2_e9b}\wch\de_p\!: \ov{H}^p(Y;\cT)\lra \wcH^{p+1}(Y;\cR), \quad p=0,1,\EE
of pointed sets such that the sequence~\eref{Snake2_e0} with 
$\wcH^*(Y;\cT)$ replaced by $\ov{H}^{p+1}(Y;\cT)$ is exact.\\

The inclusions $\fI_p\!:\ov{D}^p(\un{U})\!\lra\!D^p(\un{U})$ of short cochain complexes
commute with the refining homomorphisms~\eref{mupdfn_e} 
and induce morphisms
\BE{Snake2_e11}\fI_*\!:\ov{H}^p(\un{U};\cT)\lra\wcH^p(\un{U};\cT)
\qquad\hbox{and}\qquad
\fI_*:\ov{H}^p(Y;\cT)\lra\wcH^p(Y;\cT)\EE
of pointed sets.
By the paracompactness of~$Y$ and the reasoning in~\cite[p204]{Warner},
for every open cover $\un{U}\!\equiv\!\{{U}_{\al}\}_{\al\in\cA}$ of~$Y$ 
and every element~$d_p$  of $D^p(\un{U})$ 
there exist an open cover $\un{U}'\!\equiv\!\{{U}_{\al}'\}_{\al\in\cA'}$
refining~$\un{U}$, a refining map \hbox{$\mu\!:\cA'\!\lra\!\cA$}, and 
an element~$d_p'$  of $\ov{D}^p(\un{U}')$ such that
$\fI_p(d_p')\!=\!\mu_p^*(d_p)$.
This implies that the second map in~\eref{Snake2_e11} is a bijection.
Composing~\eref{Snake2_e9b} with this bijection, we obtain a morphism
as in~\eref{Snake2_e0b} so that the sequence~\eref{Snake2_e0} is exact.\\

A morphism of short exact sequences of sheaves of groups over~$Y$ as in~\eref{Snake2_e0a}
induces morphisms of the corresponding exact sequences of short cochain 
complexes as in~\eref{Snake2_e4} and of the inclusions~$\fI_p$ above. 
Thus, it also induces morphisms of the corresponding maps as in~\eref{Snake2_e5b} and 
as on the left-hand side of~\eref{Snake2_e11}.
These morphisms commute with the associated maps~\eref{nonabelCech_e5}
and thus induce morphisms of the maps as in~\eref{Snake2_e0b}.
This establishes the last claim.
\end{proof}

\subsection{Sheaves determined by Lie groups}
\label{CechPBs_subs}

For a Lie group~$G$ and a topological space~$Y$,
we denote by $\cS_Y(G)$ the sheaf of germs of continuous $G$-valued functions on~$Y$ and let
$$\wcH^p(Y;G)=\wcH^p\big(Y;\cS_Y(G)\big) \qquad\forall~p\!=\!0,1.$$
If $G$ is abelian, we use the same notation for all $p\!\in\!\Z$.
We begin this section
by applying Proposition~\ref{Snake_prp} to short exact sequences of sheaves
arising from short exact sequences
\BE{LGses_e} \{\1\}\lra K \stackrel{\io}{\lra} G 
\stackrel{\fj}{\lra} Q\lra\{\1\} \EE
of Lie groups.
For certain kinds of exact sequences~\eref{LGses_e},
the topological condition on~$Y$ of Definition~\ref{locLGsimp_dfn} appearing
in the resulting statement of Corollary~\ref{SnakeLG_crl} reduces
to the locally \hbox{$H^1$-simple} notion of Definition~\ref{locH1simp_dfn}. 
For such exact sequences of Lie groups and topological spaces, 
we combine Proposition~\ref{CechSing_prp} and Corollary~\ref{SnakeLG_crl} to obtain 
an exact sequence mixing \v{C}ech and singular cohomology; see Proposition~\ref{SnakeLG_prp}.\\

A homomorphism $\io\!:K\!\lra\!G$ of Lie groups induces a homomorphism
$$\io\!: \cS_Y(K)\lra \cS_Y(G)$$
of sheaves over every topological space and thus morphisms
$$\io_*\!:\wcH^p(Y;K)\lra \wcH^p(Y;G)$$
of pointed sets for $p\!=\!0,1$;
the $p\!=\!0$ case of~$\io_*$ is a group homomorphism.\\

A continuous map \hbox{$f\!:Y\!\lra\!Y'$} induces group homomorphisms
$$f^*\!: \wch{C}^p\big(\un{U};\cS_{Y'}(G)\big)\lra \wch{C}^p\big(f^{-1}(\un{U});\cS_Y(G)\big),
\qquad p\!\in\!\Z,$$
for every open cover $\un{U}$ of~$Y'$
that commute with the \v{C}ech coboundaries and group actions for 
the sheaves~$\cS_{Y'}(G)$ and~$\cS_Y(G)$
constructed in Section~\ref{Nonabelian_subs}
and with the refining homomorphisms as in~\eref{mupdfn_e}.
Thus, $f$ induces morphisms
$$f^*\!:\wcH^p(Y';G)\lra \wcH^p(Y;G)$$
of pointed sets for $p\!=\!0,1$;
the $p\!=\!0$ case of~$f^*$ is a group homomorphism.
If $G$ is abelian, then $f$ induces such a morphism for every $p\!\in\!\Z$ and
this morphism is a group homomorphism.
If in addition $\io$ is a homomorphism of Lie groups as above, then 
the diagram
$$\xymatrix{\wcH^p(Y';K)\ar[rr]^{\io_*}\ar[d]_{f^*}&&\wcH^p(Y';G)\ar[d]^{f^*}\\
\wcH^p(Y;K)\ar[rr]^{\io_*}&&\wcH^p(Y;G)}$$
commutes.
 
\begin{dfn}\label{locLGsimp_dfn}
Let~\eref{LGses_e} be a short exact sequence of Lie groups.
A topological space~$Y$ is \sf{locally simple with respect to~\eref{LGses_e}} 
if it is locally path-connected and  
for every neighborhood $U\!\subset\!Y$ of a point $y\!\in\!Y$
and a continuous map $f_U\!:U\!\lra\!Q$ 
there exist a neighborhood $U'\!\subset\!U$ of~$y$ and a continuous map $f_U'\!:U'\!\lra\!G$ 
such that $f_U|_{U'}\!=\!\fj\!\circ\!f_U'$.
\end{dfn}

For any topological space~$Y$, a short exact sequence~\eref{LGses_e} of Lie groups 
induces an exact sequence 
$$\{\1\}\lra \cS_Y(K)\stackrel{\io}{\lra} \cS_Y(G)\stackrel{\fj}{\lra} \cS_Y(Q)$$ 
of sheaves over~$Y$. 
The last map above is surjective if and only if  
$Y$ is locally simple with respect to~\eref{LGses_e}.
If the restriction of~$\fj$ to the identity component~$G_0$
of~$G$ is a double cover of~$Q_0$ and $\pi_1(Q_0)$ is (possibly infinite) cyclic, 
then the condition of Definition~\ref{locLGsimp_dfn} is equivalent 
to $Y$ being locally $H^1$-simple.
This follows from the lifting property for covering projections 
\cite[Lemma~79.1]{Mu},
Hurewicz isomorphism for~$\pi_1$ \cite[Proposition~7.5.2]{Spanier}, 
and the Universal Coefficient Theorem for Cohomology \cite[Theorem~53.3]{Mu2}.

\begin{crl}\label{SnakeLG_crl}
Let $Y$ be a paracompact space and \eref{LGses_e} be a short exact sequence of Lie groups
such that $\io(K)\!\subset\!Z(G)$.
If $Y$ is locally simple with respect to~\eref{LGses_e},  
then there exist morphisms 
\BE{Snake3_e0b}\wch\de_p\!: \wcH^p(Y;Q)\lra \wcH^{p+1}(Y;K), \quad p=0,1,\EE
of pointed sets such that the sequence
\BE{Snake3_e0}\begin{split}
\{\1\}&\lra  \wcH^0(Y;K)\stackrel{\io_*}{\lra}\wcH^0(Y;G)\stackrel{\fj_*}{\lra} \wcH^0(Y;Q)
\stackrel{\wch\de_0}{\lra}\\ 
&\stackrel{\wch\de_0}{\lra}
 \wcH^1(Y;K)\stackrel{\io_*}{\lra} \wcH^1(Y;G)\stackrel{\fj_*}{\lra} \wcH^1(Y;Q) 
\stackrel{\wch\de_1}{\lra} \wcH^2(Y;K)
\end{split}\EE
of morphisms of pointed sets is exact.
The maps~$\wch\de_p$ are natural with respect to morphisms 
of short exact sequences of Lie groups and with respect to
continuous maps between paracompact spaces
that are locally simple with respect to~\eref{LGses_e}.
\end{crl}

\begin{proof}
Since $Y$ is locally simple with respect to~\eref{LGses_e},
the sequence
\BE{SnakeLG_e3}
\{\1\}\lra \cS_Y(K)\stackrel{\io}{\lra} \cS_Y(G)\stackrel{\fj}{\lra} \cS_Y(Q)\lra\{\1\}\EE
of sheaves over~$Y$ is exact.
Since $\io(K)\!\subset\!Z(G)$, $\io(\cS_Y(K))\!\subset\!Z(\cS_Y(G))$.
The existence of morphisms~\eref{Snake3_e0b} so that the sequence~\eref{Snake3_e0} is exact
thus follows from the first statement of Proposition~\ref{Snake_prp}.\\

A morphism of short exact sequences of Lie groups as in~\eref{LGses_e} satisfying the conditions
at the beginning of the statement of the proposition induces a morphism
of the corresponding short exact sequences of sheaves as in~\eref{SnakeLG_e3}.
Thus, the naturality of~\eref{Snake3_e0b} with respect to morphisms of short exact 
sequences of Lie groups follows 
from the second statement of Proposition~\ref{Snake_prp}.\\

A continuous map $f\!:Y\!\lra\!Y'$ between paracompact spaces that 
are locally simple with respect to~\eref{LGses_e} induces a morphism
of the corresponding exact sequences of
short cochain complexes as in~\eref{Snake2_e4} 
and of the inclusions~$\fI_p$ as in the proof of Proposition~\ref{Snake_prp}.
Thus, it also induces morphisms of the corresponding maps as in~\eref{Snake2_e5b} and 
as on the left-hand side of~\eref{Snake2_e11}.
These morphisms commute with the associated maps~\eref{nonabelCech_e5}
and thus induce morphisms of the maps as in~\eref{Snake3_e0b}.
This establishes the naturality of~\eref{Snake3_e0b} with respect to continuous
maps.
\end{proof}

\begin{prp}\label{SnakeLG_prp}
Let $Y$ be a paracompact locally $H^1$-simple space and 
\BE{SnakeLG_e0a} \{\1\}\lra \Z_2 \stackrel{\io}{\lra} G 
\stackrel{\fj}{\lra} Q\lra\{\1\} \EE
be an exact sequence of Lie groups such that $\io(\Z_2)\!\subset\!Z(G)$
and $\pi_1(Q_0)$ is cyclic.
Then there exist morphisms 
\BE{SnakeLG_e0b}
\wch\de_0\!: \wcH^0(Y;Q)\lra H^1(Y;\Z_2)
\quad\hbox{and}\quad
 \wch\de_1\!: \wcH^1(Y;Q)\lra \wcH^2(Y;\Z_2),\EE
of pointed sets such that the sequence
\BE{SnakeLG_e0}\begin{split}
\{\1\}&\lra  H^0(Y;\Z_2)\stackrel{\io_*}{\lra}\wcH^0(Y;G)\stackrel{\fj_*}{\lra} \wcH^0(Y;Q)
\stackrel{\wch\de_0}{\lra}\\ 
&\stackrel{\wch\de_0}{\lra}
 H^1(Y;\Z_2)\stackrel{\io_*}{\lra} \wcH^1(Y;G)\stackrel{\fj_*}{\lra} \wcH^1(Y;Q) 
\stackrel{\wch\de_1}{\lra} \wcH^2(Y;\Z_2)
\end{split}\EE
of morphisms of pointed sets is exact.
If in addition $Y$  is locally $H^2$-simple, then the same statement 
with $\wcH^2(Y;\Z_2)$ replaced by $H^2(Y;\Z_2)$ also holds.
The maps~$\wch\de_0$ and~$\wch\de_1$ are natural with respect to morphisms 
of exact sequences of Lie groups as in~\eref{SnakeLG_prp} and 
with respect to continuous maps between paracompact locally $H^1$-simple spaces.
\end{prp}

\begin{proof}
Since $Y$ is locally $H^1$-simple, it is locally simple with respect to
the exact sequence~\eref{SnakeLG_e0a} in the sense of Definition~\ref{locLGsimp_dfn}. 
Thus, this proposition with all $H^p(Y;\Z_2)$ replaced by $\wcH^p(Y;\Z_2)$  
is a specialization of Corollary~\ref{SnakeLG_crl}.
By Proposition~\ref{CechSing_prp}, we can then replace $\wcH^0(Y;\Z_2)$ by $H^0(Y;\Z_2)$
and $\wcH^1(Y;\Z_2)$  by $H^1(Y;\Z_2)$.
If in addition $Y$ is locally $H^2$-simple, then $\wcH^2(Y;\Z_2)$ 
can also be replaced by~$H^2(Y;\Z_2)$.
\end{proof}

\subsection{Relation with principal bundles}
\label{CechPB_subs}

Let $G$ be a Lie group and $Y$ be a topological space.
We recall below the standard identification of
the set $\Prin_Y(G)$
of equivalence (isomorphism) classes of principal \hbox{$G$-bundles}
over~$Y$ with the pointed set~$\wcH^1(Y;G)$.
This identification is key for applying Proposition~\ref{SnakeLG_prp} 
to principal \hbox{$G$-bundles}, including to study $\Spin$- and $\Pin^{\pm}$-structures
in the classical perspective of Definition~\ref{PinSpin_dfn}.\\

Suppose $\pi_P\!:P\!\lra\!Y$ is a principal $G$-bundle.
Let $\un{U}\!\equiv\!(U_{\al})_{\al\in\cA}$ be an open cover of~$Y$ so that 
the principal $G$-bundle $P|_{U_{\al}}$ is trivializable for every $\al\!\in\!\cA$. 
Thus, for every $\al\!\in\!\cA$ there exists a homeomorphism
\begin{gather*}
\Phi_{\al}\!:P\big|_{U_{\al}}\lra U_{\al}\!\times\!G 
\qquad\hbox{s.t.}\\ 
\pi_{\al;1}\!\circ\!\Phi_{\al}=\pi_P,\quad
\pi_{\al;2}\big(\Phi_{\al}(pu)\!\big)=\big(\pi_{\al;2}\big(\Phi_{\al}(p)\!\big)\!\big)\!\cdot\!u
~~\forall\,p\!\in\!P|_{U_{\al}},\,u\!\in\!G,
\end{gather*}
where $\pi_{\al;1},\pi_{\al;2}\!:U_{\al}\!\times\!G\!\lra\!U_{\al},G$
are the two projection maps.
Thus, for all $\al,\be\!\in\!\cA$ there exists a continuous map
$$g_{\al\be}\!:U_{\al\be}\lra G \qquad\hbox{s.t.}\quad
\pi_{\al;2}\big(\Phi_{\al}(p)\big)=
g_{\al\be}\big(\pi_P(p)\!\big)\!\cdot\!
\big(\pi_{\be;2}\big(\Phi_{\be}(p)\!\big)\!\big)
\quad\forall\,p\!\in\!P\big|_{U_{\al\be}}\,.$$
These continuous maps satisfy
$$g_{\be\ga}\big|_{U_{\al\be\ga}}\!\cdot\!
g_{\al\ga}^{-1}\big|_{U_{\al\be\ga}}\!\cdot\!g_{\al\be}\big|_{U_{\al\be\ga}}=\1
\qquad\forall\,\al,\be,\ga\!\in\!\cA.$$
Therefore, $g_P\!\equiv\!(g_{\al\be})_{\al,\be\in\cA}$ lies in 
$\ker\de_1\!\subset\!\wch{C}^1(\un{U};\cS_Y(G))$ and thus defines an element
$$[g_P]\in \wcH^1(Y;G).$$
We show below that $[g_P]$ depends only on the isomorphism class of~$P$.\\

Suppose $\un{U}'\!\equiv\!\{U_{\al}'\}_{\al\in\cA'}$ is a refinement of~$\un{U}$.
If $\mu\!:\cA'\!\lra\!\cA$ is a refining map, then 
$$\Phi_{\al}\!\equiv\!\Phi_{\mu(\al)}\big|_{P|_{U_{\al}'}}\!:
P\big|_{U_{\al}'}\lra U_{\al}'\!\times\!G $$
is a trivialization of the principal $G$-bundle $P|_{U_{\al}'}$ for every $\al\!\in\!\cA'$.
The corresponding transition data~is 
$$g_P'\equiv 
\big\{g_{\mu(\al)\mu(\be)}|_{U_{\al\be}'}\!:U_{\al\be}'\!\lra\!G\big\}_{\al,\be\in\cA'}
=\mu_1^*g_P.$$
Since 
$$[g_P]=[\mu_1^*g_P]\in \wcH^1(Y;G),$$
it is thus sufficient to consider trivializations of isomorphic vector bundles
over a common cover (otherwise we can simply take the intersections of 
open sets in the two covers).\\

Suppose $\Psi\!:P\!\lra\!P'$ is an isomorphism of principal $G$-bundles over~$Y$
and the principal $G$-bundle $P'|_{U_{\al}}$ 
is trivializable for every $\al\!\in\!\cA$. 
Thus, for every $\al\!\in\!\cA$ there exists a homeomorphism
\begin{gather*}
\Phi_{\al}'\!:P'\big|_{U_{\al}}\lra U_{\al}\!\times\!G 
\qquad\hbox{s.t.}\\ 
\pi_{\al;1}\!\circ\!\Phi_{\al}'=\pi_{P'},\quad
\pi_{\al;2}\big(\Phi_{\al}'(p'u)\!\big)=\big(\pi_{\al;2}\big(\Phi_{\al}'(p')\!\big)\!\big)\!\cdot\!u
~~\forall\,p'\!\in\!P|_{U_{\al}},\,u\!\in\!G.
\end{gather*}
For every $\al\!\in\!\cA$, there then exists a continuous map
$$f_{\al}\!:U_{\al}\lra G \qquad\hbox{s.t.}\quad
\pi_{\al;2}\big(\Phi_{\al}'\big(\Psi(p)\!\big)\!\big)=
f_{\al}\big(\pi_P(p)\big)\!\cdot\!
\big(\pi_{\al;2}\big(\Phi_{\al}(p)\!\big)\!\big)
\quad\forall\,p\!\in\!P\big|_{U_{\al}}\,.$$
The transition data $g_{P'}\!\equiv\!(g_{\al\be}')_{\al,\be\in\cA}$ 
determined by the collection $\{\Phi_{\al}'\}_{\al\in\cA}$ of trivializations of~$P'$
then satisfies
$$g_{\al\be}'=f_{\al}\big|_{U_{\al\be}}\!\cdot\!g_{\al\be}\!\cdot\!f_{\be}^{-1}\big|_{U_{\al\be}}
\qquad\forall\,\al,\be\!\in\!\cA.$$
Thus, $g_{P'}\!=\!f\!*\!g_P$, where $f\!\equiv\!(f_{\al})_{\al\in\cA}$, and 
$$\big[g_{P'}\big]=\big[g_P\big]\in \wcH^1(Y;G).$$
We conclude that the element $[g_P]\!\in\!\wcH^1(Y;G)$ constructed above
depends only on the isomorphism class of the principal $G$-bundle~$P$
over~$Y$.\\

Conversely, suppose $[g]\!\in\!\wcH^1(Y;G)$. 
Let $\un{U}\!\equiv\!(U_{\al})_{\al\in\cA}$ be an open cover of~$Y$
and $g\!\equiv\!\{g_{\al\be}\}_{\al,\be\in\cA}$ be an element 
of $\ker\de_1\!\subset\!\wch{C}^1(\un{U};\cS_Y(G))$ representing~$[g]$.
Define
\begin{gather*}
\pi_{P_g}\!:
P_g=\Big(\bigsqcup_{\al\in\cA}\{\al\}\!\!\times\!U_{\al}\!\times\!G\Big)\!\!\Big/\!\!\!\sim_g
\lra Y,\\
\big(\al,y,g_{\al\be}(y)u\big)\sim_g\big(\be,y,u\big)
\quad\forall\,\al,\be\!\in\!\cA,\,(y,u)\!\in\!U_{\be}\!\times\!G\,.
\end{gather*}
This is a principal $G$-bundle over~$Y$ with trivializations
$$\Phi_{\al}\!:P_g\big|_{U_{\al}}\lra U_{\al}\!\times\!G, \quad
\Phi_{\al}\big([\al,y,u)]\big)=(y,u),$$
for $\al\!\in\!\cA$ and the associated transition data~$g$.
Thus, 
\BE{gPg_e}\big[g_{P_g}\big]=[g]\in \wcH^1(Y;G).\EE
We show below that the isomorphism class~$[P_g]$ of~$P_g$
depends only on~$[g]$.\\

Suppose $\un{U}'\!\equiv\!\{{U}_{\al}'\}_{\al\in\cA'}$ is a refinement of
$\un{U}$ and $\mu\!:\cA'\!\lra\!\cA$ is a refining map.
The~map
\begin{gather*}
\Psi\!: P_{\mu^*g}\!\equiv\!
\Big(\bigsqcup_{\al\in\cA'}\!\{\al\}\!\!\times\!U_{\al}'\!\times\!G\Big)
\!\!\Big/\!\!\!\sim_{\mu^*g}
\lra P_g\!\equiv\!
\Big(\bigsqcup_{\al\in\cA}\!\{\al\}\!\!\times\!U_{\al}\!\times\!G\Big)
\!\!\Big/\!\!\!\sim_g,\\
\Psi\big([\al,y,u]\big)= \big[\mu(\al),y,u\big],
\end{gather*}
is then an isomorphism of principal $G$-bundles.
Thus, it is sufficient to show that if
$$g,g'\in\ker\de_1\subset\wch{C}^1\big(\un{U};\cS_Y(G)\big) \qquad\hbox{and}\qquad
[g]=[g']\in \wcH^1\big(\un{U};\cS_Y(G)\big),$$
then the principal $G$-bundles $P_g$ and $P_{g'}$ are isomorphic.
By definition, there exists 
$$f\!\equiv\!(f_{\al})_{\al\in\cA}\in\wch{C}^0\big(\un{U};\cS_Y(G)\big) \qquad\hbox{s.t.}\qquad
g'=f\!*\!g.$$
The~map
\begin{gather*}
\Psi\!: P_g\!=\!
\Big(\bigsqcup_{\al\in\cA}\!\{\al\}\!\!\times\!U_{\al}\!\times\!G\Big)\!\!\Big/\!\!\!\sim_g
\lra P_{g'}\!=\!
\Big(\bigsqcup_{\al\in\cA}\!\{\al\}\!\!\times\!U_{\al}\!\times\!G\Big)\!\!\Big/\!\!\!\sim_{g'},\\
\Psi\big([\al,y,u]\big)=\big[\al,y,f_{\al}(y)\!\cdot\!u\big],
\end{gather*}
is then an isomorphism of principal $G$-bundles.\\

Let $P$ be a principal $G$-bundle over~$Y$,
$\{\Phi_{\al}\}_{\al\in\cA}$ be a collection of trivializations of~$P$
over an open cover $\un{U}\!\equiv\!(U_{\al})_{\al\in\cA}$,
and $g_P\!\equiv\!(g_{\al\be})_{\al,\be\in\cA}$ be the corresponding transition data.
The~map
$$\Psi\!: P\lra P_{g_P}\!\equiv\!
\Big(\bigsqcup_{\al\in\cA}\!\{\al\}\!\!\times\!U_{\al}\!\times\!G\Big)
\!\!\Big/\!\!\!\sim_g, \quad
\Psi(p)=\big[\al,\Phi_{\al}(p)\big]~~\forall\,p\!\in\!P|_{U_{\al}},\,\al\!\in\!\cA,$$
is then an isomorphism of principal $G$-bundles.
Along with~\eref{gPg_e}, this implies that the~maps
\BE{PrinH1_e}\begin{aligned}
\Prin_Y(G)&\lra \wcH^1(Y;G), &\qquad [P]&\lra [g_P],\\ 
\wcH^1(Y;G)&\lra \Prin_Y(G), &\qquad [g]&\lra[P_g],
\end{aligned}\EE
are mutual inverses that identify $\Prin_Y(G)$ with $\wcH^1(Y;G)$.\\

If $f\!:Y\!\lra\!Y'$ is a continuous map and $P\!\lra\!Y'$ is a principal $G$-bundle,
then
$$\big[g_{f^*P}\big]=f^*\big[g_P\big]\in\wcH^1(Y;G).$$
Thus, the identifications~\eref{PrinH1_e} are natural with respect to continuous maps.

\begin{crl}\label{w1_crl}
Let $Y$ be a paracompact locally $H^1$-simple space and 
$\Phi_Y$ be as in~\eref{CechSing_e0}.
For every real line bundle~$V$ over~$Y$,
$$ \wcH^1\big(Y;\Z_2\big) \ni\Phi_Y\big(w_1(V)\big)
=\big[g_{\O(V)}\big]\in \wcH^1\big(Y;\O(1)\big)$$
under the canonical identification of the groups $\Z_2$ and $\O(1)$.
\end{crl}

\begin{proof}
By the Universal Coefficient Theorem for Cohomology \cite[Theorem~53.5]{Mu2},
the homomorphism
$$\ka\!:H^1(Y;\Z_2)\lra\Hom\big(\pi_1(Y),H^1(S^1;\Z_2)\big), \qquad
\big\{\ka(\eta)\big\}\big(f\!:S^1\!\lra\!Y\big)=f^*\eta,$$
is injective.
By the naturality of $w_1$, $\Phi_Y$, and~\eref{PrinH1_e}, it is thus sufficient
to show~that 
\BE{w1crl_e3}\wcH^1\big(\R\P^1;\Z_2\big) \ni\Phi_{\R\P^1}\big(w_1(f^*V)\big)
=\big[g_{\O(f^*V)}\big]\in \wcH^1\big(\R\P^1;\O(1)\big)\EE
for every continuous map $f\!:\R\P^1\!\lra\!Y$.
Since every line bundle over the interval $[0,1]$ is trivializable,
the line bundle $f^*V$ is isomorphic to either the trivial line bundle~$\tau_1$ or 
the real tautological line bundle~$\ga_{\R;1}$.
Both sides of~\eref{w1crl_e3} vanish in the first case.
Since \eref{PrinH1_e} is a bijection, this implies that the right-hand side of~\eref{w1crl_e3}
does not vanish in the second case.
The left-hand side of~\eref{w1crl_e3} does not vanish in this case
by the {\it Normalization Axiom} for Stiefel-Whitney classes \cite[p38]{MiSt}.
\end{proof}

\subsection{Orientable vector bundle over surfaces}
\label{OVB_subs}

\noindent
We now combine the description of complex line bundles in terms of \v{C}ech cohomology
and the identification of some \v{C}ech cohomology groups with the singular ones
to characterize orientable vector bundles over surfaces and their trivializations.

\begin{lmm}\label{CLB_lmm}
Let $Y$ be a paracompact locally contractible space.
The homomorphism
$$c_1\!: \LB_{\C}(Y)\lra H^2(Y;\Z), \qquad L\lra c_1(L),$$
from the group of isomorphism classes of complex line bundles is an isomorphism.
\end{lmm}

\begin{proof} By Section~\ref{CechPB_subs}, there is a natural bijection 
$$\LB_{\C}(Y)\lra \wch{H}^1(Y;\C^*);$$
it is a group isomorphism in this case.
By the proof of Proposition~\ref{CechSing_prp}, there are natural isomorphisms
$$\wch{H}^p(Y;\Z)\approx H^p(Y;\Z) \qquad\forall\,p\!\in\!\Z.$$
By the reasoning in \cite[Section~5.10]{Warner}, $\cS_Y(\C)$ is a fine sheaf.
Along with \cite[p202]{Warner}, this implies~that
$$\wch{H}^p\big(Y;\cS_Y(\C)\big)=\{0\} \qquad\forall\,p\!\in\!\Z^+.$$
Since $Y$ is locally contractible, it is locally simple with respect
to the short exact sequence
$$\{0\}\lra \Z \lra \C \stackrel{\exp}{\lra} \C^*\lra \{0\}$$
of abelian Lie groups in the sense of Definition~\ref{locLGsimp_dfn}.
Thus, we obtain a commutative diagram
\BE{CLB_e5}\begin{split}
\xymatrix{ \{0\}\!=\!\wcH^1(\C\P^{\i};\C)\ar[r]\ar[d]^{f^*}& 
\LB_{\C}(\C\P^{\i})\ar[r]^>>>>>>{\wch\de_1}\ar[d]^{f^*}& 
H^2\big(\C\P^{\i};\Z\big) \ar[d]^{f^*} \ar[r]& \wcH^2(\C\P^{\i};\C)\!=\!\{0\}\ar[d]^{f^*}\\
\{0\}\!=\!\wcH^1(Y;\C)\ar[r]& \LB_{\C}(Y) 
\ar[r]^>>>>>>{\wch\de_1}& H^2(Y;\Z) \ar[r]& \wcH^2(Y;\C)\!=\!\{0\}}
\end{split}\EE
of group homomorphisms for every continuous map $f\!:Y\!\lra\!\C\P^{\i}$.\\

\noindent
Let $\ga_{\C}\lra\!\C\P^{\i}$ be the complex tautological line bundle.  
By \cite[Theorem~14.5]{MiSt}, $H^2(\C\P^{\i};\Z)$ is freely generated by~$c_1(\ga_{\C})$.
Along with the exactness of the top row in~\eref{CLB_e5}, 
this implies that $\wch\de_1\!=\!\pm\!c_1$ in this row.
By \cite[Theorem~14.6]{MiSt}, for every complex line bundle $L$ over~$Y$
there exists a continuous map $f\!:Y\!\lra\!\C\P^{\i}$ such that $L\!=\!f^*\ga_{\C}$.
Along with the commutativity of~\eref{CLB_e5}, these statements imply that
$\wch\de_1\!=\!\pm c_1$ in the bottom row in~\eref{CLB_e5} as well. 
The claim now follows from the exactness of this row.
\end{proof}

\begin{rmk}\label{CLB_rmk}
The statement and proof of Lemma~\ref{CLB_lmm} apply to 
any paracompact space~$Y$ satisfying the $k\!=\!2$ case of Definition~\ref{locH1simp_dfn}
with $H^p(\cdot;\Z_2)$ replaced by $H^p(\cdot;\Z)$.
\end{rmk}

\begin{crl}\label{X2VB_crl0}
Let $Y$ be a CW complex with cells of dimension at most~2.
If \hbox{$H^2(Y;\Z)\!=\!\{0\}$}, then every orientable vector bundle~$V$ over~$Y$ 
is trivializable.
\end{crl}

\begin{proof}
Let $n\!=\!\rk\,V$. 
If $n\!=\!1$, then $V$ is an orientable line bundle and is thus trivializable.
Suppose $n\!\ge\!2$.
Since the cells of $Y$ are of dimension at most~2, there exists a rank~2 orientable vector bundle
$L$ over~$Y$ such~that
\BE{X2VB0_e3}V\approx L\oplus\big(Y\!\times\!\R^{n-2}\big).\EE
The real vector bundle~$L$ admits a complex structure~$\fI$.
It can be obtained by fixing an orientation and a metric on~$L$ and defining 
$\fI v\!\in\!L$ for $v\!\in\!L$ nonzero to 
be the vector which is orthogonal to~$v$ and has the same length as~$v$ so that
$v,\fI v$ form an oriented basis for a fiber of~$L$.
By Lemma~\ref{CLB_lmm}, $(L,\fI)$ is trivializable as a complex line bundle.
Along with~\eref{X2VB0_e3}, this establishes the claim.
\end{proof}

\begin{crl}\label{X2VB_crl1}
Let $\Si$ be a surface, possibly with boundary, and $n\!\ge\!3$.
The~map
$$\OVB_n(\Si)\lra H^2(\Si;\Z_2), \qquad V\lra w_2(V),$$
from the set of isomorphism classes of rank~$n$ oriented vector bundles over~$\Si$
is a bijection.
\end{crl}

\begin{proof} We can assume that $\Si$ is connected.
If $\Si$ is not compact or has boundary, then 
$$H^2(Y;\Z),H^2(\Si;\Z_2)=\{0\}.$$ 
By Corollary~\ref{X2VB_crl0}, we can thus assume that $\Si$ is closed
and so $H^2(\Si;\Z_2)\!\approx\!\Z_2$.\\

Let $C\!\subset\!\Si$ be an embedded loop separating~$\Si$ into two surfaces,~$\Si_1$ and~$\Si_2$,
with boundary~$C$.
By Corollary~\ref{X2VB_crl0},
a rank~$n$ oriented vector bundle~$V$ over $\Si$ is isomorphic 
to the vector bundle obtained by gluing $\Si_1\!\times\!\R^n$ and $\Si_2\!\times\!\R^n$
along $C\!\times\!\R^n$ by a clutching map $\vph\!:C\!\lra\!\SO(n)$.
Since $n\!\ge\!3$, $\pi_1(\SO(n))\!\approx\!\Z_2$.
It thus remains to show that there exists an orientable vector bundle~$V$ over~$\Si$
with $w_2(V)\!\neq\!0$.\\

Let $\ga_{\C;1}\!\lra\!\C\P^1$ be the complex tautological line bundle.
Since $w_2(\ga_{\C;1})$ is the image of $c_1(\ga_{\C;1})$ under the reduction homomorphism
$$H^2\big(\C\P^1;\Z\big)\lra H^2\big(\C\P^1;\Z_2\big),$$
$w_2(\ga_{\C;1})\!\neq\!0$ by the proof of Lemma~\ref{CLB_lmm}.
If $f\!:\!\Si\!\lra\!\C\P^1$ is a degree~1 map with respect to the $\Z_2$-coefficients,
then
$$\blr{w_2(f^*\ga_{\C;1}),[\Si]_{\Z_2}}
=\blr{w_2(\ga_{\C;1}),f_*[\Si]_{\Z_2}}
=\blr{w_2(\ga_{\C;1}),[\C\P^1]_{\Z_2}}\neq0.$$
Thus, $w_2$ of the orientable vector bundle 
$$V\equiv f^*\ga_{\C;1}\oplus\big(\Si\!\times\!\R^{n-2}\big)\lra\Si$$
is nonzero.
\end{proof}

\begin{crl}\label{X2VB_crl1b}
Suppose $\wt\Si$ is a compact surface with two boundary components and 
$\wh\Si$ is a closed surface obtained from~$\wt\Si$ by identifying these components
with each other.
Let $n\!\in\!\Z^+$ and $\wh{V}\!\lra\!\wh\Si$ be the orientable vector bundle obtained
from $\wt\Si\!\times\!\R^n$ by identifying its restrictions to~$\prt\wt\Si$ 
via a clutching map \hbox{$\vph\!:S^1\!\lra\!\SO(n)$}.
If $\wh\Si$ is connected and $n\!\ge\!3$, then
$\vph$ is homotopically trivial if and only if $w_2(\wh{V})\!=\!0$. 
\end{crl}

\begin{proof}
By Corollary~\ref{X2VB_crl0}, 
every rank~$n$ orientable vector bundle over~$\wt\Si$ is trivializable.
Thus, every rank~$n$ orientable vector bundle~$\wh{V}$ over~$\wh\Si$ is obtained
from $\wt\Si\!\times\!\R^n$ by identifying its restrictions to
the two components of~$\prt\wt\Si$ 
via a clutching map \hbox{$\vph\!:S^1\!\lra\!\SO(n)$}.
Since 
$$\pi_1\big(\SO(n)\big)\approx\Z_2 \qquad\hbox{and}\qquad 
H^2(\wh\Si;\Z_2)\approx\Z_2,$$ 
the claim thus follows from Corollary~\ref{X2VB_crl1}.
\end{proof}

\begin{crl}\label{X2VB_crl2}
Let $\Si$ be a compact connected surface with boundary components $C,C_1,\ldots,C_k$
and $V$ be an orientable vector bundle over~$\Si$.
If $\rk\,V\!\ge\!3$, then every trivialization of~$V$ over \hbox{$C_1\!\cup\!\ldots\!\cup\!C_k$}
extends to a trivialization~$\Psi$ of~$V$ over~$\Si$ and 
the homotopy class of the restriction of~$\Psi$ to $V|_C$ is determined
by the homotopy class of its restriction to~$V|_{C_1\cup\ldots\cup C_k}$.
\end{crl}

\begin{proof}
Let $n\!=\!\rk\,V$ and choose an orientation on~$V$.
Denote by $\wh\Si$ the connected surface with one boundary component~$C$ obtained
from~$\Si$ by attaching the 2-disks~$D^2_i$ along the boundary components~$C_i$.
Let $\wh{V}$ be the oriented vector bundle over $\wh\Si$ obtained by identifying
each $D^2_i\!\times\!\R^n$ with~$V$ over~$C_i$ via the chosen trivialization~$\phi_i$.
By Corollary~\ref{X2VB_crl0}, the oriented vector bundle~$\wh{V}$ admits a trivialization~$\Psi$.
Since there is a unique homotopy class of trivializations of $\wh{V}|_{D^2_i}$,
the restriction of~$\Psi$ to $V|_{C_i}$ is homotopic to~$\phi_i$ and thus can deformed
to be the~same.\\

\noindent
Suppose $\Psi,\Psi'$ are trivializations of $V\!\lra\!\Si$ restricting to 
the same trivializations~$\phi_i$ of $V|_{C_i}$ for every $i\!=\!1,\ldots,k$.
Denote by $\wh\Si$ (resp.~$\wt\Si$) the closed (resp.~compact) surface 
obtained from two copies of~$\Si$, $\Si$ and~$\Si'$, by identifying
them along the boundary components corresponding to~$C,C_1,\ldots,C_k$
(resp.~$C_1,\ldots,C_k$).
Thus, $\wt\Si$ has two boundary components, each of which corresponds to~$C$,
and $\wh\Si$ can be obtained from $\wt\Si$ by identifying these two boundary components.
Let 
$$\wt{q}\!:\wt\Si\lra\Si \qquad\hbox{and}\qquad \wh{q}\!:\wh\Si\lra\Si$$
be the natural projections.
The trivializations~$\Psi$ and~$\Psi'$ induce a trivialization~$\wt\Psi$ 
of~$\wt{q}^*V$ over~$\wt\Si$
which restricts to $\Psi$ and~$\Psi'$ over $\Si,\Si'\!\!\subset\!\wh\Si,\wt\Si$.
The bundle~$\wh{q}^*V$ over~$\wh\Si$ is obtained from $\wt{q}^*V$
by identifying the copies of $V|_C$ via the clutching map \hbox{$\vph\!:S^1\!\lra\!\SO(n)$}
determined by the difference between the trivializations of $V|_C$ induced by~$\Psi$
and~$\Psi'$. 
Since 
$$w_2\big(\wh{q}^*V\big)=\wh{q}^*w_2(V)=0\in H^2(\wh\Si;\Z_2),$$
Corollary~\ref{X2VB_crl1b} implies that $\vph$ is homotopically trivial.
Thus,
$\Psi$ and $\Psi'$ determine the same homotopy class of trivializations of~$V|_C$.
\end{proof}

\noindent
For an oriented vector bundle $V\!\lra\!Y$, let $\Triv(V)$ denote 
the set of homotopy classes of trivializations of~$V$.
For an oriented vector bundle~$V$ over a surface~$\Si$, we define the~map
\BE{veSidfn_e}\ve_V\!:\Triv\big(V|_{\prt\Si}\big)\lra \Z_2\EE
by setting $\ve_V(\phi)\!=\!0$ for the trivializations~$\phi$ of $V|_{\prt\Si}$ that
extend to trivializations of~$V$ over~$\Si$ and 
$\ve_V(\phi)\!=\!1$ for the trivializations~$\phi$ that do~not.

\begin{crl}\label{X2VB_crl2b}
Let $\Si$ be a compact connected surface with $\prt\Si\!\neq\!\eset$
and $V$ be an oriented vector bundle over~$\Si$.
If $\rk\,V\!\ge\!3$, then the map~\eref{veSidfn_e} is surjective and  
changing the homotopy class of a trivialization~$\phi$ over precisely one
component of~$\prt\Si$ changes the value~$\ve_V(\phi)$.
\end{crl}

\begin{proof}
This follows from $\pi_1(\SO(n))\!\approx\!\Z_2$ and Corollary~\ref{X2VB_crl2}.
\end{proof}

\section{Lie group covers}
\label{LG_app}

\noindent
This appendix reviews basic statements concerning covers of Lie groups 
by Lie groups that are Lie group homomorphisms.
Lemma~\ref{CovExt_lmm} describes the structure of connected Lie group covers.
Lemma~\ref{CovExt_lmm2} and
Proposition~\ref{CovExt_prp} do the same for covers of disconnected Lie groups
with connected restrictions to the identity component of the base.
We conclude with examples involving the groups $\Spin(n)$
and $\Pin^{\pm}(n)$ defined in Sections~\ref{SO2Spin_subs} and~\ref{O2Pin_subs}, respectively.

\subsection{Terminology and summary}
\label{LGsumm_subs}

\noindent
We call a covering projection $q\!:\wt{G}\!\lra\!G$ a \sf{Lie group covering} 
if $\wt{G}$ and $G$ are Lie groups and $q$ is a Lie group homomorphism.
We call such a cover \sf{connected} if $\wt{G}$ is connected;
this implies that so is~$G$.
Lie group coverings 
$$q\!:\wt{G}\lra G \qquad\hbox{and}\qquad  q'\!:\wt{G}'\lra G$$
are equivalent if there exists a Lie group isomorphism $\rho\!:\wt{G}\!\lra\!\wt{G}'$
such that \hbox{$q\!=\!q'\!\circ\!\rho$}.
For a connected Lie group~$G$, we denote by $\Cov(G)$ 
the set of equivalence classes of connected Lie group coverings of~$G$
and by $\pi_1(G)$ its fundamental group based at the identity~$\1$.
For any group~$H$, we denote by $\SG(H)$ the set of subgroups of~$H$. 
The next lemma is established in Section~\ref{connLG_subs}.

\begin{lmm}\label{CovExt_lmm} 
\begin{enumerate}[label=(\alph*),leftmargin=*]

\item\label{CovExtConn_it1} Let $G$ be a connected Lie group. The map 
$$\Cov(G)\lra \SG\big(\pi_1(G)\!\big), \qquad 
\big[q\!:\wt{G}\!\lra\!G\big]\lra q_*\pi_1\big(\wt{G}\big),$$
is a bijection.
For every $[q]\!\in\!\Cov(G)$ as above,
$q^{-1}(\1)$ is contained in the center of~$\wt{G}$.

\item\label{CovExtConn_it2}
Let $q\!:\wt{G}\!\lra\!G$ and $q'\!:\wt{G}'\!\lra\!G'$ be connected Lie group coverings.
A Lie homomorphism \hbox{$\io\!:G\!\lra\!G'$} lifts to a Lie group homomorphism
$\wt\io\!:\wt{G}\!\lra\!\wt{G}'$ if and only~if 
\BE{CovExt1_e}\io_*\big(q_*\pi_1(\wt{G})\!\big)\subset q_*'\pi_1(\wt{G}').\EE
If such a lift~$\wt\io$ exists, it is unique.

\end{enumerate}
\end{lmm}

\vspace{.1in}

For a group $G$, we denote by $\Aut(G)$ the group of automorphisms of~$G$.
For each $k\!\in\!G$, let
$$\fc(k)\!:G\lra G \qquad \fc(k)g= k\!\cdot\!g\!\cdot\!k^{-1},$$
denote the \sf{conjugation homomorphism} of $k$ on~$G$;
it preserves every normal subgroup~$G_0$ of~$G$.\\

If $K$ and $G_0$ are groups and 
\BE{ladfn_e}\fc\!:K\lra\Aut(G_0)\EE
is a homomorphism,
the \sf{semi-direct product of~$G_0$ and~$K$ with respect to~$\fc$} is the group
$$G\equiv G_0\!\rtimes_{\fc}\!K=G_0\!\times\!K, \quad
(g_1,k_1)\!\cdot\!(g_2,k_2)=\big(g_1\!\cdot\!\fc(k_1)g_2,k_1\!\cdot\!k_2\big).$$
The subsets $\{\1\}\!\times\!K$ and $G_0\!\times\!\{\1\}$ of~$G$ are then
a subgroup isomorphic to~$K$ and a normal subgroup isomorphic to~$G_0$, respectively.
The conjugation action of $\{\1\}\!\times\!K$ on $G_0\!\times\!\{\1\}$ under 
the natural isomorphisms is given by~$\fc$.
Furthermore, the~map
$$K\lra  G\big/\big(G_0\!\times\!\{\1\}\big), \qquad
k\lra [(\1,k)],$$
is a group isomorphism.\\

If $G_0$ is a connected Lie group and $K$ is (at most) countable, 
then $G$ above is a Lie group and $G_0\!\times\!\{\1\}$ is its identity component.
Conversely, if $G_0$ is the identity component of a Lie group~$G$, 
$K$ is a subgroup of~$G$ such~that the~map
\BE{Kcond_e}K\lra G/G_0, \qquad k\lra [k],\EE
is a group isomorphism, and \eref{ladfn_e} is the homomorphism induced by the conjugation
of~$G_0$ by the elements of~$K$, then the~map
$$G_0\!\rtimes_{\fc}\!K \lra G, \qquad (g,k)\lra gk,$$
is a Lie group isomorphism.\\

An \sf{extension} of a group $K$ by another group $H$ is a short exact sequence
\BE{GrExtdn_e}\{\wt\1\}\lra H \lra \wt{K} \stackrel{\fj}{\lra} K\lra \{\wt\1\}\EE
of groups.
An \sf{extension} of a group homomorphism \hbox{$\io\!:K\!\lra\!K'$} by a group
homomorphism \hbox{$\wt\io_0\!:H\!\lra\!H'$} is a commutative diagram
\BE{HomExtdfn_e}\begin{split}
\xymatrix{ \{\wt\1\}\ar[r]& H\ar[r]\ar[d]_{\wt\io_0}& 
\wt{K}\ar[r]^{\fj}\ar[d]^{\wt\io}& K\ar[r]\ar[d]^{\io}& \{\wt\1\}\\
\{\wt\1\}\ar[r]& H'\ar[r]& 
\wt{K}'\ar[r]^{\fj'}& K'\ar[r]& \{\wt\1\}}
\end{split}\EE
of extensions of~$K$ by~$H$ and of~$K'$ by~$H'$.
An extension of~$K$ by~$H$ as in~\eref{GrExtdn_e} is isomorphic to another such
extension with $(\wt{K},\fj)$ replaced $(\wt{K}',\fj')$
if there exists a commutative diagram as in~\eref{HomExtdfn_e} with $K'\!=\!H$,
$\io\!=\!\id$, and $\wt\io_0\!=\!\id$.
We denote by $\Ext(K,H)$ the set of equivalence classes of extensions of~$K$ by~$H$.\\

If $H$ in~\eref{GrExtdn_e} is abelian, the conjugation action of~$\wt{K}$ on $H\!\subset\!\wt{K}$
descends to an action of~$K$ on~$H$ by group isomorphisms.
If $H$ is abelian and \hbox{$\fc\!:K\!\lra\!\Aut(H)$} is any homomorphism,
denote~by
$$\Ext_{\fc}(K,H)\subset \Ext(K,H)$$
the subset of equivalence classes of extensions as in~\eref{GrExtdn_e} 
so~that 
\BE{laExtcond_e} \wt{k}\!\cdot\!\wt{g}\!\cdot\!\wt{k}^{-1}=
\fc\big(\fj(\wt{k})\big)\wt{g}
\qquad\forall~\wt{g}\!\in\!H\!\subset\!\wt{K},\,\wt{k}\!\in\!\wt{K}\,;\EE
this condition is well-defined on the equivalence classes.\\

\noindent
Suppose $G$ is a Lie group, $G_0$ is its identity component, and 
$q_0\!:\wt{G}_0\!\lra\!G_0$ is a connected Lie group covering.
In particular, $G_0$ is a normal subgroup of~$G$.
If $k\!\in\!G$ and  the conjugation homomorphism~$\fc(k)$ of~$k$ on~$G_0$
preserves the subgroup \hbox{$q_{0*}(\pi_1(\wt{G}_0))$} of $\pi_1(G_0)$, 
then it lifts uniquely to a Lie group automorphism
$$\wt\fc_{q_0}(k)\!:\big(\wt{G}_0,\1\big)\lra\big(\wt{G}_0,\1\big).$$
We denote by $\Cov_{q_0}(G)$ the set of equivalence classes of
Lie group coverings~$q$  of~$G$ that restrict to~$q_0$ over~$G_0$.
The next two statements are established in Section~\ref{disconnLG_subs}.

\begin{lmm}\label{CovExt_lmm2}
Suppose $G$ is a Lie group, $G_0$ is its identity component, and 
\hbox{$q_0\!:\wt{G}_0\!\lra\!G_0$} is a connected Lie group covering.
If $\Cov_{q_0}\!(G)\!\neq\!\eset$, then
the subgroup $q_{0*}(\pi_1(\wt{G}_0))$ of $\pi_1(G_0)$ is preserved by
the conjugation homomorphism~$\fc(k)$  for every $k\!\in\!G$. 
\end{lmm}

\begin{prp}\label{CovExt_prp}
Suppose $G$ is a Lie group, $G_0$ is its identity component, 
\hbox{$q_0\!:\wt{G}_0\!\lra\!G_0$} is a connected Lie group covering
such that the subgroup $q_{0*}(\pi_1(\wt{G}_0))$ of $\pi_1(G_0)$ is preserved 
by the conjugation homomorphism~$\fc(k)$ for every $k\!\in\!G$, and 
$K$ is a subgroup of $G$ so that the map~\eref{Kcond_e} is an isomorphism.
\begin{enumerate}[label=(\alph*),leftmargin=*]

\item\label{CovExt_it1} The map
\BE{CovExtDC_e}
\big[q\!:\wt{G}\!\lra\!G\big]\lra
\big[\{\1\}\!\lra\!q_0^{-1}(\1)\!\lra\!q^{-1}(K)\!\stackrel{q}{\lra}\!K\!\lra\!\{\1\}\big] \EE
is a bijection from $\Cov_{q_0}\!(G)$ to $\Ext_{\wt\fc_{q_0}}\!(K,q_0^{-1}(\1)\!)$.

\item\label{CovExt_it2} If $\io$ is a Lie group automorphism of~$G$ preserving 
the subgroup $q_{0*}(\pi_1(\wt{G}_0))$ of $\pi_1(G_0)$
and $\wt\io_0$ is 
the lift of $\io_0\!\equiv\!\io|_{G_0}$ to a Lie group automorphism of~$\wt{G}_0$,
then an extension~$\wt\io$ of
$$\io\!:K\lra\io(K)  \qquad\hbox{by}\qquad
\wt\io_0\!:q_0^{-1}(\1)\lra q_0^{-1}(\1)$$
induces a Lie group isomorphism~$\wt\io$, lifting~$\io$ and extending~$\wt\io_0$,
between the Lie group coverings determined by
the associated extensions of~$K$ and of~$\io(K)$ by~$q_0^{-1}(\1)$.
\end{enumerate}
\end{prp}

\subsection{Proof of Lemma~\ref{CovExt_lmm}}
\label{connLG_subs}

By \cite[Theorem~82.1]{Mu}, for every subgroup $H\!\subset\!\pi_1(G)$ there exists 
a connected covering projection \hbox{$q\!:\wt{G}\!\lra\!G$} such that 
$q_*\pi(\wt{G})\!=\!H$.
Since $\pi_1(G)$ is abelian and countably generated, the index of~$H$ in~$G$
is countable and thus $\wt{G}$ is a second countable topological space.
The smooth structure on~$G$ then lifts to a smooth structure on~$\wt{G}$
so that $q$ becomes a local diffeomorphism and all deck transformations of~$q$
are smooth.\\

Let $\wt\1\!\in\!\wt{G}$ be any preimage of the identity $\1\!\in\!G$.
Since the images of $\pi_1(\wt{G})$ and $\pi_1(\wt{G}\!\times\!\wt{G})$ 
under the homomorphisms determined by the continuous~maps
\begin{alignat*}{2}
\wt{G} \lra& G \lra G, &\qquad  \wt{g}&\lra q\big(\wt{g}\big)^{-1},\\
\wt{G}\!\times\!\wt{G} \lra &G\!\times\!G \lra G,
&\qquad \big(\wt{g}_1,\wt{g}_2\big)&\lra q\big(\wt{g}_1\big)\!\cdot\!q\big(\wt{g}_2\big),
\end{alignat*}
are $q_*\pi_1(\wt{G})$, the lifting property for covering projections  \cite[Lemma~79.1]{Mu} 
implies that the inverse~$^{-1}$ and product~$\cdot$ operations on~$G$ lift to continuous maps 
on~$\wt{G}$ and $\wt{G}\!\times\!\wt{G}$ so that the diagrams
$$\xymatrix{(\wt{G},\wt\1)\ar[d]_q\ar[r]^{^{-1}}& (\wt{G},\wt\1) \ar[d]^q&&
(\wt{G}\!\times\!\wt{G},\wt\1\!\times\!\wt\1)\ar[d]_{q\times q}
\ar[r]^>>>>>>{\cdot}& (\wt{G},\wt\1) \ar[d]^q\\ 
(G,\1)\ar[r]^{^{-1}}& (G,\1)&&
(G\!\times\!G,\1\!\times\!\1)\ar[r]^>>>>>>{\cdot}& (G,\1)}$$
commute.
Since $q$ is a local diffeomorphism, the two lifts are smooth.
By the next paragraph, they determine a group structure on~$\wt{G}$ so
that $\wt\1$ is the identity element. 
By construction, $q$ commutes with the group operations~$\cdot$.\\

Since the maps
$$\big(\wt{G},\wt\1\big)\lra\big(\wt{G},\wt\1\big), 
\qquad \wt{g}\lra\wt{g}\!\cdot\!\wt\1,\wt\1\!\cdot\!\wt{g},$$
lift the identity on $(G,\1)$, they are the identity on $(\wt{G},\wt\1)$.
Since the maps
$$\big(\wt{G},\wt\1\big)\lra\big(\wt{G},\wt\1\big), 
\qquad \wt{g}\lra\wt{g}\!\cdot\!\wt{g}^{-1},\wt{g}^{-1}\!\cdot\!\wt{g},$$
lift the constant map on $(G,\1)$, they are the constant map on $(\wt{G},\wt\1)$.
Since the maps 
\BE{LGlift_e7}
\big(\wt{G}\!\times\!\wt{G}\!\times\!\wt{G},\wt\1\!\times\!\wt\1\!\times\!\wt\1\big)
\lra\big(\wt{G},\wt\1\big), 
\qquad \big(\wt{g}_1,\wt{g}_2,\wt{g}_3\big)\lra
\wt{g}_1\!\cdot\!\big(\wt{g}_2\!\cdot\!\wt{g}_3\big),
\big(\wt{g}_1\!\cdot\!\wt{g}_2\!\big)\cdot\!\wt{g}_3,\EE
lift the maps 
$$\big(G\!\times\!G\!\times\!G,\1\!\times\!\1\!\times\!\1\big)\lra(G,\1), 
\qquad \big(g_1,g_2,g_3\big)\lra
g_1\!\cdot\!\big(g_2\!\cdot\!g_3\big),\big(g_1\!\cdot\!g_2\big)\!\cdot\!g_3,$$
which are the same, the two maps in~\eref{LGlift_e7} are also the same.
Thus, the lifted map~$\cdot$ on~$\wt{G}\!\times\!\wt{G}$ constructed in the previous
paragraph defines a group structure on~$\wt{G}$ so that $\wt\1$ is the identity element
and the lifted map~$^{-1}$ is the inverse operation.\\

Let $\wt{h}\!\in\!q^{-1}(\1)$.
Since the~map
$$\big(\wt{G},\wt\1\big)\lra\big(\wt{G},\wt\1\big), 
\qquad \wt{g}\lra \wt{h}\!\cdot\!\wt{g}\!\cdot\!\wt{h}^{-1},$$
lifts the identity on $(G,\1)$, it is the identity on $(\wt{G},\wt\1)$.
Thus, $q^{-1}(\1)$ is contained in the center of the Lie group~$\wt{G}$
and is in particular abelian.\\

Let $q$, $q'$, and $\io$ be as in Lemma~\ref{CovExt_lmm}\ref{CovExtConn_it2}.
By the lifting property for covering projections \cite[Lemma~79.1]{Mu},
$\io$ lifts to a continuous map~$\wt\io$ so that the diagram
$$\xymatrix{ \big(\wt{G},\wt\1\big) \ar[rr]^{\wt\io}\ar[d]_q &&
\big(\wt{G}',\wt\1'\big) \ar[d]^{q'}\\
(G,\1) \ar[rr]^{\io} && (G',\1') }$$
commutes if and only if~\eref{CovExt1_e} holds;
if such a lift~$\wt\io$ exists, it is unique.
If it exists, $\wt\io$ is a smooth~map because 
$q$ and~$q'$ are local diffeomorphisms and $\io$ is a smooth map.
Since the maps 
\BE{LGlift_e9}
\big(\wt{G}\!\times\!\wt{G},\wt\1\!\times\!\wt\1\big)
\lra\big(\wt{G}',\wt\1'\big), 
\qquad \big(\wt{g}_1,\wt{g}_2\big)\lra
\wt\io\big(\wt{g}_1\!\cdot\!\wt{g}_2\big),
\wt\io\big(\wt{g}_1\big)\!\cdot\!\wt\io\big(\wt{g}_2\big),\EE
lift the maps 
$$\big(G\!\times\!G,\1\!\times\!\1\big)\lra(G',\1'), 
\qquad \big(g_1,g_2\big)\lra \io\big(g_1\!\cdot\!g_2\big),\io(g_1)\!\cdot\!\io(g_2),$$
which are the same, the two maps in~\eref{LGlift_e9} are also the same.
Thus, the map~$\wt\io$ is a Lie group homomorphism.
This establishes Lemma~\ref{CovExt_lmm}\ref{CovExtConn_it2}.\\

Suppose $q\!:\wt{G}\!\lra\!G$ is a connected cover of a connected Lie group with
the identity $\1\!\in\!G$ as above and $\wt\1,\wt\1'\!\in\!q^{-1}(\1)$.
Since $\pi_1(G)$ is abelian, $q$ is a regular covering and thus $\wt\1'\!=\!\rho(\wt\1)$
for some deck transformation~$\rho$ of~$q$. 
If $\cdot$ and~$\cdot'$ are Lie group structures on~$\wt{G}$
with the identity elements~$\wt\1$ and~$\wt\1'$, respectively,
so that~$q$ is a group homomorphism with respect to both, then
$$\rho\!:\big(\wt{G},\1\big)\lra \big(\wt{G},\1'\big)$$
is a Lie group isomorphism with respect to~$\cdot$ and~$\cdot'$ by the previous paragraph.
Thus, the Lie group structure on~$G$ (uniquely) determines a Lie group structure on~$\wt{G}$ 
so that $q$ is a Lie group homomorphism.
This establishes Lemma~\ref{CovExt_lmm}\ref{CovExtConn_it1}.

\subsection{Disconnected Lie groups}
\label{disconnLG_subs}

\noindent
We next establish Lemma~\ref{CovExt_lmm2} and Proposition~\ref{CovExt_prp}
and then give some examples.

\begin{proof}[{\bf{\emph{Proof of Lemma~\ref{CovExt_lmm2}}}}]
If $q\!:\wt{G}\!\lra\!G$ is a Lie group covering which extends~$q_0$, then 
\BE{conjlift_e} q_0\big(\wt{k}\!\cdot\!\wt{g}\!\cdot\!\wt{k}^{-1}\big)=
 \fc\big(q(\wt{k})\!\big)q_0(\wt{g})
\qquad\forall~\wt{k}\!\in\!\wt{G},\,\wt{g}\!\in\!\wt{G}_0\,.\EE
Thus, $\fc(k)$ lifts to a Lie group automorphism of~$\wt{G}_0$ for every $k\!\in\!G$.
By Lemma~\ref{CovExt_lmm}\ref{CovExtConn_it2}, this implies that $\fc(k)$
preserves the subgroup $q_{0*}(\pi_1(\wt{G}_0))$ of $\pi_1(G_0)$.
\end{proof}

\begin{proof}[{\bf{\emph{Proof of Proposition~\ref{CovExt_prp}}}}] 
Suppose $q\!:\wt{G}\!\lra\!G$ is a Lie group covering which extends~$q_0$.
The subset \hbox{$q^{-1}(K)\!\subset\!\wt{G}$} is then a subgroup so~that
\BE{qladiag_e}\{\1\}\lra q_0^{-1}(\1)\lra q^{-1}(K)\stackrel{q}{\lra} K\lra \{\1\}\EE
is an extension of $K$ by $q_0^{-1}(\1)$.
Since $q_0^{-1}(\1)$ is abelian by Lemma~\ref{CovExt_lmm}\ref{CovExtConn_it1}, 
the conjugation action of $q^{-1}(K)$ 
on $q_0^{-1}(\1)$ descends to an action of~$K$.
By~\eref{conjlift_e}, the latter action is~$\wt\fc_{q_0}$.
If $q'\!:\wt{G}'\!\lra\!G$ is another Lie group covering which extends~$q_0$
and 
$\rho$ is an equivalence from~$q$ and~$q'$, 
then~$\rho$ induces an equivalence
\BE{ExtExtdfn_e3}\begin{split}
\xymatrix{ \{\1\}\ar[r]& q_0^{-1}(\1)\ar[r]\ar[d]_{\id}& 
q^{-1}(K)\ar[r]^{\fj}\ar[d]^{\rho}& K\ar[r]\ar[d]^{\id}& \{\1\}\\
\{\1\}\ar[r]& q_0^{-1}(\1)\ar[r]& q'^{-1}(K)\ar[r]^{q'}& K\ar[r]& \{\1\}}
\end{split}\EE
between the extensions of $K$ by $q_0^{-1}(\1)$
determined by~$q$ and~$q'$.
Thus, the map~\eref{CovExtDC_e} is well-defined.\\

Suppose conversely that 
\BE{wtKcond_e}\{\1\}\lra q_0^{-1}(\1)\lra \wt{K}\stackrel{\fj}{\lra} K\lra \{\1\}\EE
is an extension of $K$ by $q_0^{-1}(\1)$ such~that 
\BE{wtKcond_e2} \wt{k}\!\cdot\!\wt{h}\!\cdot\!\wt{k}^{-1}=
\wt\fc_{q_0}\big(\fj(\wt{k})\big)\wt{h}
\qquad\forall~\wt{h}\!\in\!q_0^{-1}(\1)\!\subset\!\wt{K},\wt{G}_0,~
\wt{k}\!\in\!\wt{K}\,.\EE
The quotient 
$$q\!:\wt{G}\!\equiv\!\big(\wt{G}_0\rtimes_{\wt\fc_{q_0}\circ\fj}\!\wt{K}\big)\!\big/\!\!\sim, \qquad
\big(\wt{g},\wt{k}\big)\sim\big(\wt{g}\!\cdot\!\wt{h},\wt{h}^{-1}\!\cdot\!\wt{k}\big)
~~\forall\,\wt{h}\!\in\!q_0^{-1}(\1),$$
is then a Lie group so that the~maps 
\BE{wtG0G_e}\wt{G}_0\lra \wt{G}, ~~ \wt{g}\lra \big[\wt{g},\1\big],
\qquad\hbox{and}\qquad
\io\!:\wt{K}\lra\wt{G}, ~~\io(\wt{k})=\big[\1,\wt{k}\big],\EE
are a Lie group isomorphism onto the identity component of~$\wt{G}$
and an injective group homomorphism, respectively.
The~map
$$q\!:\wt{G} \lra G, \qquad
q\big([\wt{g},\wt{k}]\big)=q_0(\wt{g})\!\cdot\!\fj(\wt{k}),$$
is a Lie group covering 
so that its composition with the first map in~\eref{wtG0G_e} is~$q_0$.
Furthermore, the image of~$\io$ is $q^{-1}(K)$ and the diagram
$$\xymatrix{ \{\1\}\ar[r]& q_0^{-1}(\1)\ar[r]\ar[d]_{\id}& 
\wt{K}\ar[r]^{\fj}\ar[d]^{\io}& K\ar[r]\ar[d]^{\id}& \{\1\}\\
\{\1\}\ar[r]& q_0^{-1}(\1)\ar[r]& q^{-1}(K)\ar[r]^q& K\ar[r]& \{\1\}}$$
is an equivalence between extensions of $K$ by $q_0^{-1}(\1)$.
Thus, the map~\eref{CovExtDC_e} is surjective.\\

Let $q'\!:\wt{G}'\!\lra\!G'$ be another Lie group covering and 
$$\xymatrix{\big(\wt{G},\1\big)\ar[d]_q 
&\ar@{_{(}->}[l]  \big(\wt{G}_0,\1\big) \ar[rr]^{\wt\io_0}\ar[d]_{q_0} &&
\big(\wt{G}_0',\1\big) \ar[d]^{q_0'} \ar@{^{(}->}[r] 
& \big(\wt{G}',\1\big) \ar[d]^{q'}\\
(G,\1)  \ar@/_1.5pc/[rrrr]^{\io} &\ar@{_{(}->}[l]  (G_0,\1) \ar[rr]^{\io_0} && (G_0',\1) 
\ar@{^{(}->}[r]  & (G',\1)}$$
be a commutative diagram of Lie group homomorphisms.
Since $\io$ is a group homomorphism,
\BE{discGext_e} \io_0\big(\fc(k)g\big)=\fc\big(\io(k)\big)\io_0(g)
\qquad\forall\,g\!\in\!G_0,~k\!\in\!K.\EE
If $\wt\io\!:\wt{G}\!\lra\!\wt{G}'$ is a Lie group homomorphism lifting~$\io$ and
extending~$\wt\io_0$, then
the diagram
\BE{discGext_e0}\begin{split} 
\xymatrix{ \{\1\}\ar[r]& q_0^{-1}(\1)\ar[r]\ar[d]_{\wt\io_0}& 
 q^{-1}(K)\ar[r]^q\ar[d]^{\wt\io}& K\ar[r]\ar[d]^{\io}& \{\1\}\\
\{\1\}\ar[r]& q_0'^{-1}(\1)\ar[r]& q'^{-1}(\io(K))\ar[r]^>>>>>{q'}& 
\io(K)\ar[r]& \{\1\}}
\end{split}\EE
is an extension of $\io$ by $\wt\io_0$.
Suppose conversely that \eref{discGext_e0} is an extension of $\io$ by $\wt\io_0$.
By~\eref{discGext_e} and Lemma~\ref{CovExt_lmm}\ref{CovExtConn_it2}, 
$$\wt\io_0\big(\fc(\wt{k})\wt{g}\big)=\fc\big(\wt\io(\wt{k})\big)\wt\io_0(\wt{g})
\qquad\forall\,\wt{g}\!\in\!\wt{G}_0,~\wt{k}\!\in\!\wt{K}.$$
This in turn implies that the~map
\BE{wtio0ext_e}\wt\io\!:\wt{G}\lra\wt{G}', \qquad
\wt\io\big(\wt{g}\!\cdot\!\wt{k}\big)=\wt\io_0(\wt{g})\!\cdot\!\wt\io(\wt{k})
\qquad\forall\,\wt{g}\!\in\!\wt{G}_0,~k\!\in\!q^{-1}(K),\EE
is a Lie group homomorphism lifting~$\io$ and extending~$\wt\io_0$.\\

We now apply the conclusion of the previous paragraph with $G\!=\!G'$ and $q_0\!=\!q_0'$.
If~$\io$ is as in Proposition~\ref{CovExt_prp}\ref{CovExt_it2}, then
$\io_0\!\equiv\!\io|_{G_0}$ lifts to a Lie group isomorphism
$$\wt\io_0\!:\big(\wt{G}_0,\1\big)\lra \big(\wt{G}_0,\1\big)$$
by Lemma~\ref{CovExt_lmm}\ref{CovExtConn_it2}.
By the previous paragraph, an extension
\BE{discGext_e5}\begin{split} 
\xymatrix{ \{\1\}\ar[r]& q_0^{-1}(\1)\ar[r]\ar[d]_{\wt\io_0}& 
\wt{K}\ar[r]^{\fj}\ar[d]^{\wt\io}& K\ar[r]\ar[d]^{\io}& \{\1\}\\
\{\1\}\ar[r]& q_0^{-1}(\1)\ar[r]& 
\wt{K}'\ar[r]^{\fj'}& \io(K)\ar[r]& \{\1\}}
\end{split}\EE
of $\io|_K$ by $\wt\io_0|_{q_0^{-1}(\1)}$ induces an isomorphism~\eref{wtio0ext_e},
lifting~$\io$ and extending~$\wt\io_0$,
between the Lie group coverings~$q$ and~$q'$ associated with the first and second 
lines in~\eref{discGext_e5}.
This establishes Proposition~\ref{CovExt_prp}\ref{CovExt_it2}.\\

In particular, an equivalence between 
the extensions of~$K$ by~$q_0^{-1}(\1)$ determined by two covers~$q$ and~$q'$
of~$G$ extending~$q_0$ as in~\eref{CovExtDC_e}
induces an equivalence between~$q$ and~$q'$. 
Thus, the map~\eref{CovExtDC_e} is injective.
This concludes the proof of Proposition~\ref{CovExt_prp}\ref{CovExt_it1}.
\end{proof}

\begin{eg}\label{PinPin_eg}
The groups $\Pin^{\pm}(n)$ defined directly in Section~\ref{O2Pin_subs}
are the quotients of certain semi-direct products
$\Spin(n)\!\rtimes\!\Z_2^{\,2}$  and $\Spin(n)\!\rtimes\!\Z_4$
by natural $\Z_2$-actions. 
The subgroup \hbox{$q_0^{-1}(\1)\!\approx\!\Z_2$} of $\Spin(n)$ appearing in~\eref{wtKcond_e}
is the subgroup generated by the element~$\wh\bI_n$ in the notation of Section~\ref{SO2Spin_subs}.
The subgroup \hbox{$K\!\approx\!\Z_2$} of~$\O(n)$ can be taken to be 
the subgroup generated by any order~2 element of $\O(n)\!-\!\SO(n)$,
i.e.~an order~2 element of $\O(n)$ with an odd number of $(-1)$-eigenvalues.
Any two distinct order~2 elements of~$\O(n)$ with precisely one $(-1)$-eigenvalue
are contained in a subgroup~$G$ isomorphic to~$\O(2)$ and are interchanged by 
an automorphism~$\io$ of~$\O(n)$ which preserves~$G$.
By Proposition~\ref{CovExt_prp}\ref{CovExt_it2}, $\io$ lifts to an isomorphism~$\wt\io$  
between the Lie group coverings $\Pin_{K_1}^{\pm}(n),\Pin_{K_2}^{\pm}(n)$  determined by the extensions
$$\xymatrix{ \{\1\}\ar[r]& \Z_2\ar[r]\ar[d]_{\id}& 
\Pin^{\pm}(1)\ar[r]^q\ar[d]^{\id}&  K_1\ar[r]\ar[d]& \{\1\}\\
\{\1\}\ar[r]& \Z_2\ar[r]&  \Pin^{\pm}(1)\ar[r]^>>>>>{q'}&  K_2\ar[r]& \{\1\}}$$
of the subgroups $K_1,K_2$ generated by the two elements.
The restriction of~$\wt\io$ to $\Pin_{K_1}^{\pm}(n)|_G$ is an isomorphism onto 
$\Pin_{K_2}^{\pm}(n)|_G$. 
By Example~\ref{Pin1Pin2_eg}, this implies~that  
$$\Pin_{K_2}^{\pm}(n)\big|_G=\Pin_{K_1}^{\pm}(n)\big|_G$$
and so $\Pin_{K_2}^{\pm}(n)\!=\!\Pin_{K_1}^{\pm}(n)$.
Thus, the criterion for distinguishing between $\Pin^-(n)$ and~$\Pin^+(n)$ 
above Definition~\ref{PinSpin_dfn} and Example~\ref{Pin1Pin2_eg}
does not depend on the choice of order~2 element of~$\O(n)$ with precisely one $(-1)$-eigenvalue
used to generate the subgroup~$K$.
However, the first statement in~\eref{wtInmsq_e2} implies that 
the subgroup of $\Pin^+(n)$ generated by the preimages of 
an order~2 element of~$\O(n)$ with precisely three $(-1)$-eigenvalues is~$\Z_4$.
Thus, the criterion for distinguishing between $\Pin^+(n)$ and $\Pin^-(n)$ depends
on the choice of the conjugacy class of order~2 elements of~$\O(n)$.
In the notation of Example~\ref{Pin1Pin2_eg}, 
\begin{alignat*}{2}
&\Pin^+(2)=\big(\R/2\pi\Z\!\rtimes_{\fc_+}\!\Z_2^{\,2}\big)\!\big/\!\!\sim,
&\qquad 
&\Pin^-(2)=\big(\R/2\pi\Z\!\rtimes_{\fc_-}\!\Z_4\big)\!\big/\!\!\sim,\\
&\fc_+(a,b)(\th\!+\!2\pi\Z)=(-1)^b\th\!+\!2\pi\Z,
&\qquad 
&\fc_-(a)(\th\!+\!2\pi\Z)=(-1)^a\th\!+\!2\pi\Z,\\
&\big(\th\!+\!2\pi\Z,(a,b)\big)\sim\big(\th\!+\!\pi\!+\!2\pi\Z,(a\!+\!1,b)\big),
&\qquad 
&\big(\th\!+\!2\pi\Z,a\big)\sim\big(\th\!+\!\pi\!+\!2\pi\Z,a\!+\!2\big).
\end{alignat*}
\end{eg}

\begin{eg}\label{SpinLift_eg}
In the case of the lift~$\wt\io_{n;m}$ of $\io_{n;m}$ in~\eref{wtioPin_e}, 
\eref{discGext_e} becomes
\begin{alignat*}{2}
\io_{n;m}'(A)&=\fc(\bI_{n;1})\io_{n;m}'(A) &\qquad &\forall~A\!\in\!\SO(m), \\
\io_{n;n-m}''(\fc(\bI_{n-m;1})A)&=\fc(\bI_{n;1})\io_{n;n-m}''(A)
&\qquad &\forall~A\!\in\!\SO(n\!-\!m).
\end{alignat*}
These two conditions are equivalent to the two equations in~\eref{wtiocomm_e}. 
In the situation of Remark~\ref{PinPin_rmk},   
the diagram~\eref{discGext_e0} would become
$$\xymatrix{ \{\1\}\ar[r]& \Z_2^2\ar[r]\ar[d]_{\wt\io_0}& 
 \wt{K}_m\!\times\!\wt{K}_{n-m}\ar[r]^q\ar[d]^{\wt\io}& K_m\!\times\!K_{n-m}
\ar[r]\ar[d]^{\io}& \{\1\}\\
\{\1\}\ar[r]& \Z_2\ar[r]& \wt{K}_n\ar[r]^>>>>>{q'}&  K_n\ar[r]& \{\1\}}$$
with $K_n$ and $\wt{K}_n$ denoting the group generated by $\bI_{n;1}$ and 
its preimage in $\Pin^{\pm}(n)$, respectively.
However, there is no group homomorphism~$\wt\io$ making this diagram commute.
\end{eg}

\begin{eg}\label{noliftLG_eg}
The connected Lie group double cover~$q_0$ of $\Spin(2)\!\approx\!S^1$
does not extend to a Lie group double cover~$q$~of
\begin{gather*}
\Pin^-(2)\equiv \R/2\pi\Z\!\times\!\Z_2,\\
\cdot\!:
\Pin^-(2)\!\times\!\Pin^-(2)\lra\Pin^-(2), \quad 
(\th_1,k_1)\!\cdot\!(\th_2,k_2)=
\big(\th_1\!+\!(-1)^{k_1}\th_2\!+\!k_1k_2\pi,k_1\!+\!k_2\big).
\end{gather*}
Such an extension would be of the form 
$$q\!:\R/2\pi\Z\!\times\!\Z_2\lra \R/2\pi\Z\!\times\!\Z_2, \quad q(\th,k)=q(2\th,k),
\quad
\wt{g}_0\!\equiv\!(0,1)\!\cdot\!(0,1)\in\big\{(\pi/2,0),(3\pi/2,0)\big\}.$$
The conjugation homomorphism on the double cover of $\Spin(2)$ by~$(0,1)$ is 
the lift of the conjugation homomorphism on~$\Spin(2)$ by~$(0,1)$
and is thus given~by
$$S^1\lra S^1, \qquad \th\!+\!2\pi\Z\lra-\th\!+\!2\pi\Z.$$
Since the lifted homomorphism does not fix $\wt{g}_0$, $\wt{g}_0$ cannot be the square of $(0,1)$
in a Lie group double cover of~$\Pin^-(2)$.
Since every element of $\Pin^-(2)\!-\!\Spin(2)$ is of order precisely~4,
$\Pin^-(2)$ contains no subgroup projecting isomorphically to $\Pin^-(2)/\Spin(2)$
and thus Proposition~\ref{CovExt_prp} does not apply in this case.
\end{eg}

\printindex[gen]
\addcontentsline{toc}{chapter}{Index of Terms}

\printindex[not]
\addcontentsline{toc}{chapter}{Index of Notation}

\end{document}